%% file: thesis-arxiv.tex
\documentclass[
    crop=false,
    a4paper,
    10pt,
]{report}

\usepackage{textcomp} 
\usepackage{booktabs}
\usepackage{amssymb,amsthm,amsmath}

\usepackage{custom-style}
\usepackage[plain]{theoremtemplate-style}
\usepackage{standalone}
\usepackage[margin=2cm]{geometry}

\renewcommand{\todo}[2][]{ }

\providecommand{\datapath}{.}

\renewenvironment{thebibliography}[1]{
  \begin{oldthebibliography}{#1}
  	 \footnotesize
    \setlength{\itemsep}{0.5em}
    \setlength{\parskip}{0em}
}
{
  \end{oldthebibliography}
}
\newcommand{\includebibliography}{%
  \cleardoublepage
  \phantomsection
  \addcontentsline{toc}{chapter}{\bibname}
}

\title{Homotopy Comomentum Maps in Multisymplectic Geometry}
\author{Antonio Michele Miti}

\begin{document}
\begin{minipage}{\textwidth}
\maketitle

\noindent\fbox{%
    \parbox{\textwidth}{%
        \textbf{Disclaimer:}
			
			This preprint is a re-edition of the doctoral thesis bearing the same title.
			\\
			The version accepted for obtaining the PhD title can be found on \href{https://limo.libis.be/primo-explore/fulldisplay?docid=LIRIAS3408626&context=L&vid=Lirias&search_scope=Lirias&tab=default_tab&lang=en_US}{Lirias (KU Leuven)} and \href{https://tesionline.unicatt.it/handle/10280/94031}{Docta (Universit\'a Cattolica del Sacro Cuore)}.

    }%
}
\vspace{2em}

\noindent\fbox{%
    \parbox{\textwidth}{%
		\textbf{Notes for the arXiv version.}
\\
Chapter 4  contains new results not yet published.
\\
Chapters 1, 2 and the appendices are intended as a review on $L_\infty$-algebras and their applications in multisymplectic geometry. They are functional to chapter 4.
\\
Chapters \ref{Chap:LeonidPaper} and \ref{Chap:MauroPaper} contain results accepted for publication and already available on arXiv. See \href{https://arxiv.org/abs/1805.01696}{arXiv:1805.01696}, \href{https://arxiv.org/abs/1906.08790}{arXiv:1906.08790}, and \href{https://arxiv.org/abs/1910.13400}{arXiv:1910.13400}.
    }%
}

\vspace{2em}
\begin{abstract}
\Momaps are a higher generalization of the notion of moment map introduced to extend the concept of Hamiltonian actions to the framework of multisymplectic geometry.
Loosely speaking, higher means passing from considering symplectic $2$-form to consider differential forms in higher degrees.
The goal of this thesis is to provide new explicit constructions and concrete examples related to group actions on multisymplectic manifolds admitting \momaps.
\\
The first result is a complete classification of compact group actions on multisymplectic spheres. The existence of a \momap pertaining to the latter depends on the dimension of the sphere and the transitivity of the group action. Several concrete examples of such actions are also provided.
\\
The second novel result is the explicit construction of the higher analogue of the embedding of the Poisson algebra of a given symplectic manifold
into the corresponding twisted Lie algebroid. 
It is also proved a compatibility condition for such embedding for gauge-related multisymplectic manifolds in presence of a compatible Hamiltonian group action. The latter construction could play a role in determining the multisymplectic analogue of the geometric quantization procedure.
\\
Finally, a concrete construction of a \momap for the action of the group of volume-preserving diffeomorphisms on the multisymplectic 3-dimensional Euclidean space is proposed.
This map can be naturally related to hydrodynamics. For instance, it transgresses to the standard hydrodynamical co-momentum map of Arnol'd, Marsden and Weinstein and others.
A slight generalization of this construction to a special class of Riemannian manifolds is also provided.  
The explicitly constructed \momap can be also related to knot theory 
by virtue of the aforementioned hydrodynamical interpretation.
Namely, it allows for a reinterpretation of (higher-order) linking numbers in terms of multisymplectic conserved quantities. 
As an aside, it also paves the road for a semiclassical interpretation of the HOMFLYPT polynomial in the language of geometric quantization.
\end{abstract}
\end{minipage}

\tableofcontents

\input{chapters/introduction/introduction}

\part{Background}
\input{chapters/Linfinity/Linfinity}

\renewcommand*{\datapath}{chapters/multisymplectic/image}

\input{chapters/multisymplectic/multisymplectic}

\part{Foreground}
\renewcommand*{\datapath}{chapters/compactactionspheres/image}
\input{chapters/compactactionspheres/compactactionspheres} 
\renewcommand*{\datapath}{chapters/higherrogersembedding/image}

\input{chapters/higherrogersembedding/higherrogersembedding}
\renewcommand*{\datapath}{chapters/hydromomaps/image}
\input{chapters/hydromomaps/hydromomaps}

\part{Appendices}
\appendix

\input{chapters/gradedmultilinear/gradedmultilinear}

\input{chapters/multibracketscoderivations/multibracketscoderivations}

\input{chapters/permutators/permutators}

\input{chapters/gradedprelie/gradedprelie}



\includebibliography
\bibliographystyle{hep} 
\bibliography{mypapers,websites,biblio-tidy}


\end{document}

%% file: chapters/introduction/introduction.tex
\chapter*{Introduction and overview}\label{Chap:Introduction}
\chaptermark{Introduction and overview}
\addcontentsline{toc}{chapter}{Introduction}
In this doctoral thesis, we are interested in developing the theory of symmetries on multisymplectic manifolds, especially when they admit a so-called \emph{\momap}.
\\
\emph{Multisymplectic structures} (also called \emph{``$n$-plectic''}) are the rather straightforward generalization of symplectic ones when closed non-degenerate $(n+1)$-forms replace $2$-forms.

Historically, the interest in multisymplectic manifolds, \ie smooth manifolds equipped with an $n$-plectic structure,  has been motivated by the need for understanding the geometrical foundations of first-order classical field theories.
The key point is that, just as one can associate a symplectic manifold to an ordinary classical mechanical system, it is possible to associate a multisymplectic manifold to any classical field system.
Hence, from the mechanical point of view, the passage from considering $2$-forms to higher forms is equivalent to leap from considering a single point-like particle constrained to some manifold to considering a continuous medium, like a filament, a brane or a fluid.
\\
The initial thrust of the theory\footnote{
	Although the first instances (in local coordinates) of such structure could be traced back to the classical work of De Donder and Weyl in the 1930’s \cite{DeDonder35,Weyl35},
	what we mean here is the modern -global- formulation.
	The latter was initiated by Kijowski and Szczyrba \cite{Kijowski1973,KS} in the 1970’s and definitely established in the 1990’s with the works of Cari\~nena-Crampin-Ibort \cite{Carinena1991b} and Gotay \cite{Gotay91}.
}, 
essentially due to its close relationship with the finite-dimensional geometrical description of classical fields, lost its momentum when it was realized that an adequate notion of "observables" was not available. 
The centrality of this concept, due to its fundamental role in the construction of most quantization schemes, prompted the research to focus on different generalized notions of symplectic manifolds, 
predominantly in the direction of infinite-dimensional smooth spaces.

One of the inherent difficulties was that most of the candidates aimed at representing the appropriate "algebra of observables", understood not only as "physically" measurable quantities but also as generators of the time evolution, suffered from the lack of identification of a suitable related algebraic structure.
Essentially, one lacked an appropriate Lie algebra structure. Early approaches tried to cure this defect by dividing out "non-physical terms", \eg by considering quantities "modulo divergences".
\\
A breakthrough happened around the year 2010 when Chris Rogers and John Baez \cite{Rogers2010} realized that a more apt structure to look for was an $L_\infty$ (strongly homotopy) algebra rather than to expect a bonafide Lie algebra.
Namely, Rogers proved that the algebraic structure encoding the observables on a multisymplectic manifold is the one of an $L_{\infty}$-algebra, that is, a graded vector space endowed with skew-symmetric multilinear brackets satisfying the Jacobi identity up to coherent homotopies.
This idea has led to a new surge of interest in the multisymplectic framework and, roughly five years later, a suitable notion of "moment map" made its appearance.
In \cite{Callies2016}, Callies, Fregier, Rogers, and Zambon, gave the definition of \emph{\momap} as a natural generalization of the notion of the ordinary (symplectic) comomentum map unifying several other previous attempts to encode momenta in multisymplectic geometry.
In a nutshell, a \momap is an $L_\infty$-morphism associated to certain infinitesimal actions which preserve the multisymplectic form of a target manifold. 

Therefore, the upshot is as follows: Multisymplectic manifolds, $L_\infty$ observables, and \momaps are higher
\footnote{
	We will mostly understand the "higher" term appearing here in a naive way, that is as "going higher" in the degrees of the considered differential forms.
	A more cogent interpretation of the use of this term can be provided in the language of higher categories and homotopy theory.
}
 generalizations of symplectic manifolds, Poisson algebras, and co-moment maps.

Being the latter concept particularly subtle and technical, there are not so many meaningful examples worked out in full details. 
In this thesis,  we try to address this problem by giving new insights and delivering new concrete constructions related to \momaps trying to further develop the understanding of symmetries in the context of multisymplectic geometry.
For instance, in chapter \ref{Chap:MauroPaper}, we exhibit how one can explicitly construct a \momap pertaining to the infinitesimal action of the volume-preserving diffeomorphisms group of Euclidean space, and other Riemannian manifolds with similar cohomology, upon resorting to Hodge theory.
This construction allows for a physical interpretation in terms of ideal fluids and singular vortices. In addition, the latter can be put in relation with knot theory.

Though mechanics provided the original motivation for the foundation of multisymplectic geometry, we stress nevertheless that mathematical-physics is not the only source where to find instances of this class of structures.
For example, any orientable $n$-dimensional manifold can be considered $(n-1)$-plectic when equipped with a volume form.
Following this purely mathematical premise, in chapter \ref{Chap:LeonidPaper} we focus our attention on multisymplectic actions of compact groups and thus deriving existence results and explicit constructions for \momaps related to actions on spheres.

A natural issue that arises when dealing with both symplectic and multisymplectic structures is to investigate what relationship exists between gauge-related multisymplectic manifolds, \ie manifolds endowed with multisymplectic forms lying in the same de Rham cohomology class. 
In chapter \ref{Chap:MarcoPaper}, we will focus on the two L$_\infty$-algebras of observables associated to a pair of gauge-related multisymplectic manifolds.
To date, no canonical correspondence is known between two gauge-related observables algebras. 
However, we will be able to exhibit a compatibility relation between those observables that are momenta of corresponding \momaps. Although this construction is essentially algebraic in nature, it admits also a geometric interpretation when applied to the particular case of pre-quantizable symplectic forms. 
This provides some evidence that this construction may be related to the higher analogue of geometric quantization for integral multisymplectic forms.

%

\section*{Structure of the thesis}
\addcontentsline{toc}{section}{Structure of the thesis}
\sectionmark{Structure of the thesis}
This thesis consists of an introduction, five chapters and four appendices.
The core of the thesis is split in two parts. The first part (chapters \ref{Chap:Linfinity} and \ref{Chap:MultiSymplecticGeometry})  consists of background material and the second part (chapters \ref{Chap:LeonidPaper}, \ref{Chap:MarcoPaper}, \ref{Chap:MauroPaper}) contains mostly original results.

\begin{figure}[h!]
	\centering
	\resizebox{\columnwidth}{!}{%
		\begin{tikzpicture}[>=stealth,every node/.style={shape=rectangle,draw,rounded corners},]
		    \node [dashed] (a1) {App. \ref{App:GradedMultilinearAlgebra}};
		    \node [dashed] (a2)[below right =of a1] {App. \ref{App:RNAlgebras}};
		    \node [dashed] (a3)[right =of a1] {App. \ref{App:UnshuffleAtors}};
		    \node [dashed] (a4)[right =of a3] {App. \ref{App:PreLie}};
		    \node (c1)[below left =of a1] {Chapter \ref{Chap:Linfinity}};
		    \node (c2)[below =of c1] {Chapter \ref{Chap:MultiSymplecticGeometry}};
		    \node (c3)[below left =of c2] {Chapter \ref{Chap:LeonidPaper}};
		    \node (c4)[below =of c2] {Chapter \ref{Chap:MarcoPaper}};
		    \node (c5)[below right =of c2] {Chapter \ref{Chap:MauroPaper}};
		    \draw[->,dashed] (a1) -- (a2);
		    \draw[->,dashed] (a3) -- (a2);
		    \draw[->,dashed] (a4) -- (a2);
		    \draw[->,dashed] (a2) -- (c1);
		    \draw[->] (c1) -- (c2);
		    \draw[->] (c2) -- (c3);
		    \draw[->] (c2) -- (c4);
		    \draw[->] (c2) -- (c5);
		\end{tikzpicture}
	}
	\caption{Structure of the chapters.}
\end{figure}
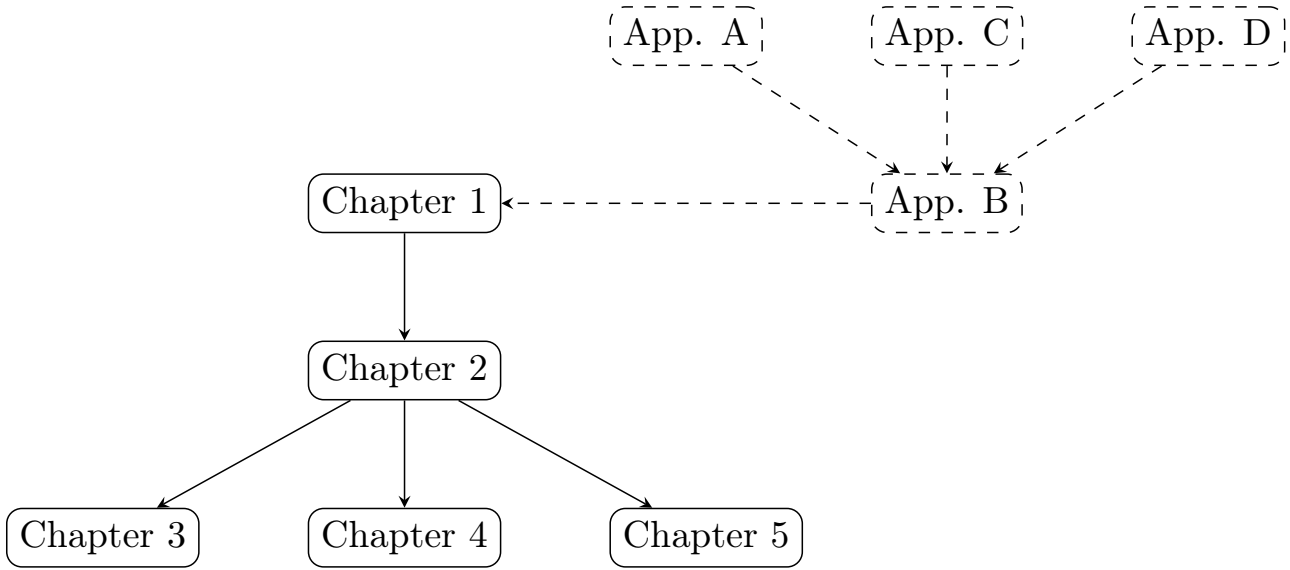

\medskip
Chapter \ref{Chap:Linfinity} is devoted to recap the conventions of multilinear algebra adopted in the thesis (which are more extensively explained in appendices \ref{App:GradedMultilinearAlgebra} and \ref{App:RNAlgebras} ) and to introduce the reader to the language of $L_\infty$-algebras.

Chapter \ref{Chap:MultiSymplecticGeometry} is a short survey on the three cornerstone concepts upon which this work is based.
Namely, the notion of \emph{multisymplectic manifold}, the definition of the pertaining \emph{$L_\infty$-algebra of observables} and the characterization of symmetries admitting a \emph{\momap}.

\medskip
In chapter \ref{Chap:LeonidPaper} we focus on multisymplectic actions of compact groups. 
We observe how, in this case, the cohomological obstructions to the existence of a \momap, already discussed in the literature in term of the equivariant cohomology, can be expressed in terms of the de Rham cohomology on the product of the acting Lie group with the base manifold.
We profit from this by giving a complete classification of proper actions of compact groups on spheres admitting \momaps, providing at the same time several explicit constructions.

In chapter \ref{Chap:MarcoPaper} we deal with the commutativity of a certain diagram involving multisymplectic structures in the same cohomology class and \momaps.
First we discuss a motivation for this purely algebraic problem which comes from the context of geometric quantization of symplectic manifolds.
Then, we show explicitly how the observables $L_\infty$-algebra can be embedded  into the $L_\infty$-algebra corresponding to another geometric object called \emph{Vinogradov algebroid}.
At last we prove how the diagram given by the aforementioned embedding with respect to two gauge related multisymplectic manifold can be closed in presence of a group action admitting a \momap.

Chapter \ref{Chap:MauroPaper} revolves around the explicit construction of 
a \momap for the action  of divergence-free vector fields on $\RR^3$ rapidly vanishing at infinity. 
This map will be called "hydrodynamical" for two main reasons: the acting group can be interpreted as the configuration (displacement) group of an incompressible fluid and the transgression of the corresponding \momap yields the standard hydrodynamical co-momentum map of  Arnol'd, Marsden, Weinstein and others.
As an application of the above \momap, a reinterpretation of the (Massey) higher order linking numbers in terms of conserved quantities within the 
multisymplectic framework is provided and  knot theoretic analogues of first integrals in involution are determined.

\medskip
The chapters in appendix are intended as a survey to prepare the algebraic framework used to present our results. In particular, the language of \RN algebras (see section \ref{Section:MultibracketsAlgebra}) will be extensively used in chapters \ref{Chap:Linfinity},\ref{Chap:MultiSymplecticGeometry} and \ref{Chap:MarcoPaper}.

Appendix \ref{App:GradedMultilinearAlgebra} describes our conventions involving graded vector spaces starting from the notion of \emph{graded objects}. We also recall the notion of tensor algebras and coalgebras and lifts to coalgebra morphisms.

In appendix \ref{App:RNAlgebras} we discuss the algebraic structure of the graded vector space of multilinear maps. 
In particular, we focus on the space of symmetric homogeneous multilinear maps, endow it with a non-associative algebra structure, called \emph{\RN product}, and check its pre-Lie property. 
The ensuing right pre-Lie algebra is then proved to be isomorphic to the algebra of graded skew-symmetric multilinear maps and to the algebra of coderivation of a corresponding symmetric tensor coalgebra. 

Appendix \ref{App:UnshuffleAtors} focuses on the combinatorial aspects involved in the definition of the \RN product and, consequently, in dealing with $L_\infty$-algebras in the "multibrackets presentation".
Namely, it describes \emph{unshuffle permutations}, states some properties involving the operators giving the action of this permutations on a given graded vector spaces, and uses these ideas to yield an explicit expression for the associator pertaining to the non-associative algebra of multilinear maps.

At last, appendix \ref{App:PreLie} is a short collection of some basic definitions and properties on graded pre-Lie algebras found scattered in the literature.

\section*{Results}
\addcontentsline{toc}{section}{Results}
\sectionmark{Results}
This thesis aims at investigating the notion of \momap and trying to work out in full detail concrete examples of interests in a different area of mathematics.
Below we summarize the main results obtained.

The first thing we want to mention cannot be qualified as a new result
but we believe it is noteworthy since it rests on an approach not conventionally used in the multisymplectic literature.
\\
In the presentation of the background material (chapters \ref{Chap:Linfinity} and \ref{Chap:MultiSymplecticGeometry}) and in chapter \ref{Chap:MarcoPaper}, we will heavily rely on the notion of \RN product (which is in particular discussed in appendix \ref{App:RNAlgebras}).
This choice proves to be convenient because it allows to express explicit constructions concerning multilinear maps, specifically the multibrackets of certain $L_\infty$-algebras, without having to deal with the combinatorial details implied by the constructions themselves.
In a certain sense, this places us at an intermediate level of abstraction between the original presentation of Lada and Stasheff \cite{LadaStasheff} and more abstract presentations, for example via operads \cite{Loday}, in which the central matter is no longer how certain multilinear maps act on vectors but it is given by the maps themselves, together with their composition rules.
In other words, this promotes the philosophy to focus on multibrackets, interpreted as morphisms in the category of graded vector spaces, without referring to their internal structure.

The core of the appendix can be subsumed by the following theorem
\begin{bigthm}[\emph{Thm \ref{Theorem:RecapGerstenhaber}}]
	Let $V$ be a graded vector space, denote by $M^{\sym}(V)$ and $M^{\skew}(V)$ the graded vector spaces of symmetric and skew-symmetric multilinear maps from $V$ into itself.
	One has the following sequence of isomorphisms in the category of  graded right pre-Lie algebras
	  \begin{displaymath}
		 \begin{tikzcd}
		 	(M^{\skew}(V[-1]),\skewgerst) \ar[r,"\Dec"',"\cong"]
		 	&  	(M^{\sym}(V),\symgerst) \ar[r,"\widetilde{L}_{\sym}"',"\cong"]
		 	& (\coDer({S(V)}),\symgerst) 
		 \end{tikzcd} 
	 \end{displaymath}
	where:
	\begin{itemize}
		\item[-] $\cs$ and $\ca$ denote respectively the symmetric and skew-symmetric \RN product;
		\item[-] $\coDer(\overline{S(V)})$ denote the space of coderivations on the reduced symmetric Tensor coalgebra $\overline{S(V)}$;
		\item[-] $\Dec$ is the \emph{d\'ecalage} operator pertaining to multilinear maps;
		\item[-] $\widetilde{L}_{\sym}$ denotes are the lift operators to the symmetric tensor coalgebra.
	\end{itemize}	 
\end{bigthm}
	Being the above three algebras isomorphic, it is legitimate to regard them as a single object \ie as \emph{the} \RN algebra pertaining to the graded vector space $V$.
	Accordingly, one could talk of three different "presentations" of the \RN algebra; respectively in terms of skew-symmetric multibrackets, symmetric multibrackets and coderivations.
In proposition \ref{Prop:SymmetricGerstenhaberAssociators} we also provided explicit formulas regarding the failure of the associativity in both the symmetric and skew-symmetric "\RN algebras".
\\
In this framework, it is immediate to view $L_\infty$-algebras structures on the graded vector space $V$ as $2$-nilpotent degree $1$ elements in the \RN algebra.

Next, we proceed to state the main results of the present work, expounded from chapter \ref{Chap:LeonidPaper} onwards.

\subsection*{Chapter \ref{Chap:LeonidPaper}}
The contents of this chapter are based on a paper co-authored by the author of the present thesis with Leonid Ryvkin \cite{Miti2019}. 

In this chapter, we expand the panorama of examples of \momaps by giving new insights on multisymplectic actions of compact groups and thus deriving existence results and explicit constructions for \momaps related to actions on spheres.
\\
The first novel result presented is the solution of the existence problem for \momaps (\Hcmm) pertaining to compact effective group actions on high-dimensional spheres seen as multisymplectic manifolds (with respect to the standard volume form).
\begin{bigthm}[\emph{Prop. \ref{prop:intransitive} and Thm \ref{thm:surprise}}]\label{thm:mainresult} 
	Let $G$ be a compact Lie group acting multisymplectically and effectively on the $n$-dimensional sphere $S^n$ equipped with the standard volume form.
	\\
	Then the action admits a \momap if and only if $n$ is even or the action is not transitive. 
\end{bigthm}
Independently from the proof, which was mainly based on results in algebraic topology, we also exhibited non-trivial classes of examples. 
\begin{itemize}
	\item The action of $SO(n)$ on $S^n$ is not transitive, hence it admits a \momap for all $n$. 
	We give an explicit construction for such a \momap in Subsection \ref{subsecson} that extends the construction given in \cite{Callies2016} only up to the $5$-dimensional sphere.
	\item The action of $SO(n+1)$ on $S^n$ admits a \momap for even $n$ only. 
	For the cases where such a \momap exists, giving explicit formulas seems to be a non-trivial task. 
	We give explicit formulas for the first two components $f_1$ and $f_2$ in terms of the standard basis of $\mathfrak {so}(n)$ in Subsection \ref{subsectra}, leaving an explicit description of the higher components as an open question.
	The core idea will be to focus on the particular cohomology of the acting group rather than working on the analytical problem of finding the primitives required for the construction of the components of a \momap. 
	\item For $SO(4)$ acting on $S^3$ no \momap exists. However, this problem can be fixed by centrally extending the Lie algebra $\mathfrak {so}(n)$ to a suitable $L_\infty$-algebra (\cf  \cite{Callies2016,Mammadova2019}).
	For instance, the Lie algebra  $\mathfrak{so}(2n)$ (giving the action of $SO(2n)$ on $S^{2n-1}$ preserving the $(2n$-$1)$-plectic volume form) could be extended to an $L_\infty$-algebra concentrated in degrees from $0$ to $(2-2n)$.
	\item Apart from the previous constructions, which all pertain to the multisymplectic structure given by the volume, we also discuss the existence of another natural \momap that can be associated to the action of the exceptional group $G_2$ on $S^7$.
\end{itemize}

\subsection*{Chapter \ref{Chap:MarcoPaper}}
The contents of this chapter are based on a preprint co-authored by the author of the present thesis and Marco Zambon \cite{Miti2020}. 

In chapter \ref{Chap:MarcoPaper}, we study the higher (multisymplectic) analogue of the standard embedding of the observables Poisson algebra pertaining to a symplectic manifold into the space of sections of the corresponding Lie algebroid.
Namely, we produce the following explicit construction for the embedding of the observables $L_\infty$-algebra of a given $n$-plectic manifold into the $L_\infty$-algebra pertaining to the corresponding \emph{Vinogradov} (or \emph{higher Courant}) algebroid:
\begin{bigthm}[\emph{Thm. \ref{thm:iso} and Cor. \ref{cor:Psi}}]\label{bigthm:rogerEmb}
	Let be $n \leq 4$.
	Consider an $n$-plectic manifold $(M,\omega)$  and take the corresponding \emph{Vinogradov (higher Courant) algebroid} $E^n:=TM\oplus\Lambda^{n-1}T^\ast M$ twisted by $\omega$.
	Denote by $L_\infty(M,\omega)$ the observables $L_\infty$-algebra associated to the former and by $L_\infty(E^{n},\omega)$ the $L_\infty$-algebra associated to the latter.
	\\
	There is a $L_\infty$-algebra embedding $\Psi$ defined by the following diagram
	\begin{displaymath}
		\begin{tikzcd}[column sep = huge, row sep = large]
			L_\infty(M,\omega) \ar[r,dashed,"\Psi"] \ar[d,sloped,"\sim"] &
			L_\infty(E^n,\omega)
			\\
			(\mathcal{A},\pi) \ar[r,"\sim","\Phi"'] & (\mathcal{A},\mu) \ar[u,hook]
		\end{tikzcd}
	\end{displaymath}
	Where:
	\begin{itemize}
		\item $\mathcal{A}$ is a graded vector subspace of $L_{\infty}(E^n,\omega)$ concentrated in degrees $0\leq k \leq 1-n$. It consists of Hamiltonian pairs (pairs composed of Hamiltonian forms together with the corresponding Hamiltonian field) in degree $0$ and differential $(n-1+k)$-forms in degree $-k$.
		\item $\pi$ and $\mu$ denote respectively the restriction of the $L_\infty$-algebra structures of $L_\infty(M,\omega)$ and $L_\infty(E^{n},\omega)$ to the graded vector space $\mathcal{A}$.
		\item $\Phi$ is a $L_\infty$-isomorphism given in components by
			\begin{displaymath}
				\Phi_n := \left(\frac{2^{n-1}}{(n-1)!} B_{n-1}\right)~\pairing_-^{\ca (n-1)}~: \mathcal{A}^{\wedge k } \to \mathcal{A}
			\end{displaymath}
			where $B_k$ denotes the $k$-th \emph{Bernoulli number}, $\pairing_-$ is the skew-symmetric pairing operator on $E^n$, and superscript $\ca(k)$ denotes the \RN product of the given operator with itself $k$ times.
	\end{itemize}
\end{bigthm}
The proof relies on the observation that the two aforementioned $L_\infty$-structures $\pi$ and $\mu$ on $\mathcal{A}$ can be reconstructed, via the \RN product, out of the same set of just four graded skew-symmetric multilinear maps.
We establish the existence of the embedding for $n\le 4$ but we expect the proof to extend to the case of arbitrary $n$.
\\
Nevertheless, this result is a generalization of a similar construction performed by Rogers \cite{Rogers2013} in the $2$-plectic case involving the ordinary twisted Courant algebroid.

Given the embedding $\Psi$, we investigated its compatibility with \Momaps and gauge transformations:
\begin{bigthm}[\emph{Thm. \ref{thm:comm}}]
	Consider two gauge related $n$-plectic forms $\omega$ and $\widetilde{\omega}:= \omega + \d B$ on the smooth manifold $M$.
	Consider a smooth action of the Lie group $G$ on $M$ admitting \momap with respect to both $\omega$ and $\widetilde{\omega}$. Denote by $f:\mathfrak{g}\to L_\infty(M,\omega)$ and $\tilde{f}:\mathfrak{g}\to L_\infty(M,\widetilde{\omega})$ the two \momaps.
	Denote by $(E^n,\omega)$ the Vinogradov Algebroid twisted by $\omega$.
	\\
	The following diagram commutes in the category of $L_\infty$-algebras:
	\begin{displaymath}
	\begin{tikzcd}
		&
		L_{\infty}(M,\omega) \ar[r,"\Psi"]
		&
		L_{\infty}(E^n,\omega) \ar[dd,"\tau_B"]
		\\[-1em]
		\mathfrak{g}\ar[ru,"f"]
		 \ar[dr,"\widetilde{f}"']
		\\[-1em]
		&
		L_{\infty}(M,\widetilde{\omega}) \ar[r,"\Psi"]
		&
		L_{\infty}(E^n ,\widetilde{\omega})
	\end{tikzcd}	
	\end{displaymath}
	where $\Psi$ is the embedding introduced in theorem \ref{bigthm:rogerEmb} and the rightmost vertical arrow is the gauge transformation isomorphism 
	$L_{\infty}(E^n,\omega) \cong L_{\infty}(E^n,\widetilde{\omega}) $.
\end{bigthm}
Carrying out this construction in the $1$-plectic case, it is possible to interpret the morphism $\Psi$ in term of prequantization.
We suppose that our construction could play a similar role in the prequantization of higher dimensional systems.


\subsection*{Chapter \ref{Chap:MauroPaper}}
The contents of this chapter are based on two papers co-authored by the author of thesis and Mauro Spera \cite{Miti2018,Miti2019a}. 

In this chapter, we exhibit a \momap pertaining to the action of the infinitesimal action of divergence free vector fields, \ie the infinite-dimensional Lie algebra $\sDiff_0$, on the three dimensional Euclidean space, that is a $2$-plectic manifold with respect to the standard volume form $\nu$:
\begin{bigthm}[\emph{Thm. \ref{Thm:HydrodynamicalComoment}}]			The infinitesimal action of $v:\mathfrak{g}\to \mathfrak{X}(\mathbb{R}^3)$, concretely given by the inclusion of divergence free fields in the set of all vector fields, 
			admits a \momap $(f)$ with components $f_j: \Lambda^j {\mathfrak g} \to  \Omega^{2-j} (\R^3)$ given by
			\begin{displaymath}
				\begin{aligned}
					f_1 =&~ \flat\circ {\rm curl}^{-1}
					\\
					f_2 =&~ \Delta^{-1} \circ\delta\circ \mu_2		
				\end{aligned}
			\end{displaymath}
		where $\mu_2(p):= f_{1} (\partial p) +  \iota(v_p) \omega$ (a term introduced in remark \ref{Rem:TermMuByMauro}), $\delta$ is the de Rham induced by the Hodge structure, and the inverse of the vector calculus operators involved, $\curl$  and $\Delta$,  have to be thought of as their corresponding Green operators (hence they are not unique).
\end{bigthm}
This object has an interesting interpretation in the context of fluid dynamics since it transgresses to the standard hydrodynamical co-momentum map of Arnol'd, Marsden and Weinstein and others.
In theorem \ref{Thm:RiemannGeneralization}, we also show how this construction could be then generalized to a suitable class of Riemannian manifolds.
Furthermore, we discuss a covariant phase space interpretation of the coadjoint orbits associated to the Euler evolution for perfect fluids and, in particular, of Brylinski's manifold of smooth knots.

The last observation prepares the ground for an application of the aforementioned  \momap in the context of knot theory.
Namely, we provide a reinterpretation of the (Massey) higher order linking numbers in terms of conserved quantities within the multisymplectic framework thus determining knot theoretic analogues of first integrals in involution.
\begin{bigthm}
[\emph{Prop. \ref{Prop:MasseyMess}, Thm. \ref{thm:VorticityFormExact} and \ref{Thm:MasseyMess}}]
	Consider a $n$-link $L=\coprod_{i=1}^n L_i $ with $L_i:S^1\to \R^3$ the parametrization of the $i$-th component.
	\begin{itemize}
		\item 	The velocity $1$-form $v_L$ (definition \ref{Def:Velocity1Form}) pertaining to $L$ is an Hamiltonian form of $L_\infty(\R^3,\nu)$ and it lie in the image of the  \momap $f$ given in the previous theorem (in other terms, $v_L$ is a "momentum" with respect to the action of divergence free fields).
	\end{itemize}
	\smallskip 
	Consider the case that the $n$-link $L$ above satisfies the property that the cohomology classes of all Massey $2$-forms $\Omega_{i j}:= v_i \wedge v_l$, for any pairs of components $i$ and $j$ of $L$, vanish in the cohomology of the link $H(S^3\setminus L)$.
	In other words, all mutual first order linking numbers are zero.
	\begin{itemize}
		\item The primitive $v_{i j}$ corresponding to the Massey form $\Omega_{i j}$ are Hamiltonian with Hamiltonian vector field given by $\xi_{i j}= \flat \cdot \ast \cdot \Omega_{i j}$. Hence they are in particular momenta with respect to the \momap given above (\ie $f_1(\xi_{i j})= v_{i j}$).
		\item $v_{i j}$ are strictly conserved quantities (see definition \ref{Def:conservedQuantities}) along the flow given by the velocity $1$-form $v_L$ associated to the link.
		\item The  commutation rule (with respect to the binary brackets $\lbrace\cdot,\cdot\rbrace$ of $L_\infty(\R^3,\nu)$) 
		$$\lbrace v_{i j}, v_{k \ell}\rbrace=0$$
		holds for any given primitives of two arbitrary Massey $2$-forms.
	\end{itemize}
\end{bigthm}
The second part of the previous statement can be easily expressed for any higher-order linking number, of order greater than $2$.
The derivation of the previous results essentially relies on the vortex theoretic approach to $n$-links.
Considering \emph{fluid configurations with linked singular vortices} is the cornerstone idea for building our bridge between multisymplectic geometry and knot theory.
\\
By further pursuit of this path, it was also possible  to give a semiclassical interpretation of the HOMFLYPT polynomial, building on the Liu-Ricca hydrodynamical approach to the latter and on the Besana-Spera symplectic approach to framing.

\section*{Previous works}
\addcontentsline{toc}{section}{Previous works}
\sectionmark{Previous works}
As explained in the previous section, most of the original results presented here also appeared in the following four preprints \cite{Miti2018}, \cite{Miti2019}, \cite{Miti2019a} and \cite{Miti2020}.
To date, three of them have been accepted for publication.

This project owes a lot to the recent development in multisymplectic geometry published in the last ten years.
\\
In large part, our research is built upon the works written by Chris Rogers \cite{Rogers2010,Rogers2011,Rogers2013}, Marco Zambon \cite{Zambon2012,Fregier2015,Callies2016} and Leonid Ryvkin  \cite{zbMATH06448534,Ryvkin2016,Ryvkin2016a,Ryvkin2018}, together with their collaborators John Baez, Martin Callies, Yael Fregier, Camille Laurent-Gengoux, and Tilmann Wurzbacher.

All our considerations regarding the possible application of multisymplectic geometry in knot theory are instead motived by a series of articles authored by Alberto Besana, Vittorio Penna and Mauro Spera \cite{BeSpe06,Pe-Spe89,Pe-Spe92,Pe-Spe00,Pe-Spe02,Spe06}.

The contents of the appendix draw inspiration from numerous lecture notes and sources available online.
Above all, especially in what pertains the identification of the algebraic structure of the graded space of skew-symmetric multibrackets, we should pinpoint the course on deformation theory by Marco Manetti \cite{Manetti-website}.

We also want to mention that the present work has been submitted shortly after the publication of the thesis, bearing a similar title, authored by Leyli Mammadova \cite{Mammadova2020a}.
Both theses share a comparable background material  but differ, though, in the different flavour used to present it.
Her novel results mainly concern two further generalizations of the notion of \momap given by the so-called\emph{ weak homotopy comoment map} and \emph{$L_2$ moment maps}.

\section*{Conventions}\label{Sec:conventions}
\addcontentsline{toc}{section}{Conventions}
\sectionmark{Conventions}
%
Throughout the thesis we will essentially work within two categories only.
On the algebraic side, we will mostly work in the category of graded (non-necessarily finite dimensional) vector spaces and, on the geometric side,  we will mainly deal with finite dimensional smooth real manifolds.
%
	Accordingly, the composition of maps\footnote{In this text, function composition is always meant as precomposition $(f\circ g)(x) = f(g(x))$; in other terms, we apply transformations on the left. The same convention holds for the composition of morphism in any category as in \cite[S 1.8]{MacLane1978}.}, mostly indicated with the symbol $\circ$, will be often denoted by simple juxtaposition in those situations where it cannot be confused with the point-wise product of functions.
%
%
When dealing with multilinear maps, we will introduce the Gerstenhaber, symmetric \RN, and skew-symmetric \RN products, denoted respectively by $\gerst, \symgerst$ and $\skewgerst$.
These can be thought of as suitable notion of composition between two multilinear maps.

Our conventions in graded multilinear algebra will be extensively recalled in section \ref{Sec:ConventionsMultiLinearAlgebras} and further justified in appendix \ref{App:GradedMultilinearAlgebra}. 
We only mention here that, for any given graded vector space $V$,
 we will seamlessly identify 
 linear maps from $V^{\otimes n},V^{\odot n}$ or $V^{\wedge n}$ into $W$ 
 with maps $\bigtimes^k V \to W$ of arity $k$ enjoying respectively the extra property of being multilinear, multilinear graded symmetric and multilinear graded skew-symmetric.
\\
%
In what pertains to homological algebra, given any cochain complex $C=(C^\bullet,d)$ we denote by $Z^k(C)=ker(d^{(k)})$ the subgroup of cocycles and by $B^k(C)=d~C^{k-1}$ the subgroup of coboundaries. 
In the case of chain complexes, we employ the same notation with lower indices.

On the geometric side, all our objects will be smooth unless differently specified.
We will take for granted the notions of smooth manifold, smooth map, diffeomorphism, smooth bundle, vector field, differential form, and  de Rham calculus.
In particular, all the manifolds considered will be smooth, finite-dimensional, Hausdorff, and second-countable.
In section \ref{Sec:MultiCartan}, we will state our conventions regarding the Cartan calculus of multi-vector fields. 
We only stress here that the contraction operator $\iota$ with decomposable multi-vector fields $p=\xi_1\wedge\dots\wedge \xi_n$ will be given by $\iota_p = \iota_{\xi_n}\dots\iota_{\xi_1}$.

The vast majority of this work deals with the notion of symmetry in the sense of group actions on smooth manifolds.
We take mostly for granted the notions of Lie group, Lie algebra, $G$-principal bundle, smooth Lie group action, (infinitesimal) Lie algebra action, and equivariant map.
We will denote by $Ad$ the adjoint action of a Lie group on itself, by $ad$ the representation of the Lie group on its Lie algebra, and by $ad^\ast$ the coadjoint representation of the group on the dual of its Lie algebra.
\emph{Left (smooth) actions} will be employed throughout, unless differently specified (see for instance \ref{rem:RightActionMess}). 
Hence, by $\theta:G\action M$ we mean a group homomorphism  $\vartheta: G \to \Diff(M)$ into the diffeomorphism group, understanding the group structure on $\Diff(M)$ again as precomposition, that is smooth in the appropriate sense. (A right action is, on the other hand, given by an antihomomorphism.)
The corresponding Lie algebra action will be the morphism of Lie algebras $\mathfrak{g}\to \mathfrak{X}(M)$ given by fundamental vector fields according to the prescription of equation \eqref{eq:LeftFundVF}.

In what concerns the combinatorics involved in this work, most of the constructions that we will encounter will be expressed in terms of unshuffles.
We denote by $S_n$ the group of permutations of $n$ elements and by $\ush{i_1,\dots, i_\ell}$ the subgroup of $(i_1,\dots,i_\ell)$-\emph{unshuffle permutatations}. Recall that $\sigma \in S_{n}$ is a $(i,n-i)$-unshuffle if $\sigma_{k}<\sigma_{k+1}$ for any $k\neq i$. Further details are given in appendix \ref{App:UnshuffleAtors}.

%
In the appendix, we will make elementary use of some basic concepts in category theory like functor, natural transformation, limits, colimits, and monoidal structure.

\ifstandalone
	\bibliographystyle{../../hep} 
	\bibliography{../../mypapers,../../websites,../../biblio-tidy}
\fi

\cleardoublepage


%% file: chapters/Linfinity/Linfinity.tex
\chapter{$L_\infty$-algebras}\label{Chap:Linfinity}
$L_{\infty}$-algebras, also known as \emph{strongly homotopy Lie algebras} or SH Lie algebras, are a generalization of \emph{Lie algebras} where one requires that the Jacobi equation is satisfied only up a controlling term.
\\
According to Stasheff \cite{Stasheff2019}, the idea of considering "Jacobi equations up to homotopy" began to sprout in conjunction with the developments in homotopy theory occurring during the mid-20th century.
However, the precise mathematical formalization of "strongly homotopy Lie algebras"  only appears in 1993, in a greatly influential paper by Lada and Stasheff \cite{LadaStasheff}. 
The authors seemed prompted to crystallize this concept after its progressive and ubiquitous appearance, roughly sparked in the eighties,  in the many different branches of theoretical physics. 
In fact, some first examples had already been found in supergravity, string theory and quantization (see \cite{nlab:l-infinity-algebra} for an updated list of applications in physics).
\\
Besides their physical applications, $L_{\infty}$-structures took a crucial role in deformation theory exemplified by the Deligne's leading principle that \emph{deformation of any given algebraic structure or geometric structure is governed by a strong homotopy algebra} \cite{Deligne}.
\\
In the last decade, 
they also began to assume a prominent role in the context of the geometric approach to \emph{classical} (in the sense of "local" and "prequantum") \emph{field theory}, \ie the classical mechanics of system with a continuum of degrees of freedom.
Namely, Baez and Rogers noticed the existence of an infinite-dimensional $L_\infty$-algebra behaving like the analogue for continuum media of the observables (Poisson) algebra of ordinary mechanical systems.
This observation leads to the introduction of the so-called \emph{(Rogers) $L_\infty$-algebra of observables} to any multisymplectic manifold \cite{Rogers2010}. 
\\
The latter concept will be of pivotal importance in the following chapters. 
In particular, most of the results involving the multisymplectic analogue of \emph{moment maps} can only be correctly phrased in the language of $L_\infty$-structures.

\medskip
In this chapter, we present some background material on $L_{\infty}$-algebras.
The contents are not new in their substance (which can largely be found in the seminal articles \cite[\S 3]{LadaStasheff}\cite{LadaMarkl} as well as in more recent surveys like \cite[Ch. 2, \S 1]{Schatz2009}, \cite[\S 6]{Doubek2007},\cite{Ryvkin2016a} and \cite{Reinhold2019})
but are somehow "new" in the form of their presentation.
Namely, our discussion relays heavily on the notion of \emph{\RN} products between graded multilinear maps.
\\
The main motivation for this unusual choice is that it allows to perform computations on "multibrackets" in a relatively more agile way with respect to dealing with expressions evaluated on arbitrary lists of objects.
Furthermore, this notation allows to keep track more easily of many prefactors involved in the combinatorics underlying the theory of $L_\infty$-algebras without requiring to pass to others layers of abstraction (like coderivations or operads) that are not perfectly suited for the kind of explicit constructions that we desired to deliver.

Specifically, in section \ref{Sec:ConventionsMultiLinearAlgebras} we explicitly state the conventions in graded multi-linear algebra, including the definition of \RN algebras, that will be employed throughout this thesis.
This section is meant as a summary of the contents of appendices \ref{App:GradedMultilinearAlgebra} and \ref{App:RNAlgebras} which are self-contained surveys on the \RN products starting from the basics in graded linear algebra.
\\
In section \ref{Sec:LinfinityAlgebras}, we discuss three possible equivalent perspectives on $L_\infty$-algebras.
Namely as a graded vector space with a family of skew-symmetric multibrackets satisfying \emph{higher Jacobi} equations,
		as a graded vector space with a family of symmetric multibrackets satisfying (symmetric)  \emph{higher Jacobi} equations
	and as a cofree graded cocommutative coalgebra endowed with a degree $1$ codifferential.
We also give explicit expressions for morphisms and their composition.
Finally, in subsection \ref{SubSection:studycase} we discuss some examples related to the specific $L_\infty$-structures studied in the rest of the thesis.

\section{Conventions in Graded multi-linear algebra}\label{Sec:ConventionsMultiLinearAlgebras}
Most of the objects of this thesis sit in the category of \emph{$\ZZ$-graded vector spaces}.
From now on, we will often drop the $\ZZ$ and simply talk about \emph{graded vector space}.
In this section we summarize the conventions that will be adopted in the body of this text. The rationale of these notations, especially in what pertains the so-called "Koszul convention", are more extensively explained in appendices \ref{App:GradedMultilinearAlgebra} and \ref{App:RNAlgebras}.
\\
We think of graded vector space as a functor between the set $\ZZ$, seen as a discrete category, and the category of (not necessarily finite-dimensional) vectors spaces $\Vect$ over the field field $\mathbb{R}$ in characteristic $0$ (it will always be the field  of real numbers in what pertains this text).
We write a graded vector space as
\begin{displaymath}
	V= (k \mapsto V^k)
	~.
\end{displaymath}
In practical terms, a graded vector space can be thought of as a family of vector spaces parametrized by $\mathbb{Z}$.
We denote by $V^\oplus = \bigoplus_{n\in\ZZ} V^n$ the ordinary vector space obtained as a direct sum of all the components of the graded vector space $V$ (see lemma \ref{lemma:totalfunctor} and remark \ref{rem:totaloplus} in appendix).
\\
We call $V^k$ the \emph{$k$-th component of $V$} and an element $v\in V^k$ is called \emph{homogeneous element in degree $k$}.
We denote by $|\cdot|$ the \emph{grading map} giving for any homogeneous element $v\in V^k$, for any $k\in \ZZ$, its degree $k$.
We say that the graded vector space is \emph{concentrated in degree $n$} if $V^k=0$ for $k\neq n$. Ordinary vector spaces will be regarded as graded vector spaces concentrated in degree $0$.
\\
Morphisms between graded vector spaces $\varphi:V\to W$ are natural transformations between the corresponding functors, hence, they are defined by a collection of linear maps maps 
\begin{displaymath}
	\varphi = \lbrace\varphi^n \in \Hom_{\Vect}(V^n,W^n) \rbrace_{n \in \ZZ}
	~.
\end{displaymath}
We will refer to elements in $\Hom_{\Vect^\ZZ}(V,W)$ as \emph{graded-morphisms} or \emph{degree preserving (linear) maps}.
Epimorphisms, monomorphisms and isomorphisms in $\Vect^{\ZZ}$ are given by collections of epimorphisms, monomorphisms and isomorphisms in $\Vect$, hence by surjective, injective and invertible linear maps respectively.

The category $\Vect^{\ZZ}$ is closed, the hom space between two given graded vector spaces is a graded vector space given by
\begin{equation}\label{eq:gvecthomspace}
	\Hom_{\Vect^{\ZZ}}(V,W) = \left(k \mapsto \Hom_{\Vect}(V^k,W^k)\right)
	~.
\end{equation}	
We will often neglect this "internal grading" understanding it as the ordinary vector space obtained by direct sum on all the components, note that one has
\begin{displaymath}
	\Hom_{\Vect^{\ZZ}}(V,W)^{\oplus} = \Hom_{\Vect}(V^{\oplus},W^{\oplus})
	~.
\end{displaymath}
The category $\Vect^{\ZZ}$ inherit from $\Vect$ the property to be Abelian and thence complete and cocomplete. In particular limits in $\Vect^{\ZZ}$ are constructed out of collections of limits on the components. 
For instance, the \emph{direct sum} of graded vector spaces is defined as the graded vector space given as follows:
\begin{displaymath}
	V\oplus W := (k \mapsto V^k\oplus W^k)
~.
\end{displaymath}
%

The action of the \emph{shift endofunctor} $[k]:\Vect^{\ZZ}\to\Vect^{\ZZ}$, for any $k\in\ZZ$, is given on components as
\begin{displaymath}
	(V[k])^i = V^{k+i}
	~.
\end{displaymath} 
An homogeneous vector of $v\in V$ in degree $|v|=n$ can be seen as an homogeneous vector of $V[k]$ in degree $(n-k)$, the latter will be denoted as
\begin{displaymath}
	v_{[k]} \in (V[k])^{(n-k)}
	~
\end{displaymath}
\ie $|v_{[k]}| = |v|-k$.
By definition, one has the following identification of functors $[k][\ell]=[\ell][k]=[k+\ell]$ (see subsection \ref{sec:degreeshifts}).
\\
Beside the shift functor, we will sometimes make use of the following two functors altering the grading of graded vector spaces.
We call \emph{Truncation up to degree $n$} the functor
\begin{equation}\label{eq:truncationFunctor}
	\morphism{\trunc_{n}~}
	{\Vect^\ZZ}
	{\Vect^\ZZ}
	{\left(k\mapsto V^k\right)}
	{\left(k\mapsto \begin{cases} V^k & \text{ if } k<n\\ 0 & \text{ if } k\geq n \end{cases} \right)}
	~,
\end{equation}
and call \emph{degree reversing} the functor
\begin{equation}\label{eq:degreeReversingFunctor}
	\morphism{\setminus~}
	{\Vect^\ZZ}
	{\Vect^\ZZ}
	{\left(k\mapsto V^k\right)}
	{\left(k\mapsto (\setminus V)^k = V^{-k}\right)}
	~.
\end{equation}
In both cases, the action on morphisms is the obvious one.

We call \emph{homogeneous map in degree $n$} between the graded vector spaces $V$ and $W$ any graded morphism 
	\begin{displaymath}
		\morphism{f}
		{V}
		{W[n]}
		{v}
		{(f(v))_{[n]}}
		~,
	\end{displaymath}
we denote by $|f|=n$ the degree of the homogeneous map $f$.
Homogeneous maps in degree $0$ are the graded morphisms defined above.
The action of the shift functor $[\ell]$ on an homogeneous map $f:V\to W[|f|]$ (\cf  remark \ref{rem:shiftHomogMapsTrickySign}) is given by
		\begin{displaymath}
			\morphism{f[\ell]}
			{V[\ell]}
			{W[|f|][\ell]}
			{v_{[\ell]}}
			{(f(v))_{[|f|][\ell]}}
			~.
		\end{displaymath}
\\
By the closedness property of the category $\Vect^{\ZZ}$, also homogeneous maps in a fixed degree $n$ forms a graded vector space.
Neglecting again this internal grading, we define the graded vector space of homogeneous map in any degree as
\begin{equation}\label{eq:homogmapshomspace}
	\underline{\Hom}_{\Vect^{\ZZ}}(V,W):=
	\left(
		k \mapsto \left(\Hom_{Vect^{\ZZ}}(V,W[k])\right)^{\oplus}
	\right)
	~.
\end{equation}
All the algebraic structures considered in this thesis will be graded linear, hence they could be thought of as graded vector spaces endowed with extra structure hence as subcategories of $\Vect^{\ZZ}$. Therefore, from now on, we will omit the subscript $\Vect^{\ZZ}$ when denoting the hom-spaces defined in equations \eqref{eq:homogmapshomspace} and \eqref{eq:gvecthomspace}.
Furthermore, we will often refer to element in $\underline{\Hom}^k(V,W)$ simply as \emph{linear maps} (in degree $k$).
We stress that in our convention only homogeneous maps of degree $0$ are "morphisms", hence when we will write a diagram in the category $\Vect^{\ZZ}$ of graded vector space, arrows have always to be interpreted as degree $0$ linear maps.
We will also employ the following decorated notation for particular subclasses of morphisms:
\begin{center}
\begin{tabular}{ c l }
 $\rightarrowtail$ & monomorphism, \ie degree-wise injective, \\ 
 $\twoheadrightarrow$ & epimorphism, \ie degree-wise surjective, \\  
 $\hookrightarrow$ & inclusion monomorphism, \\
 $\xrightarrow{\sim}$ & isomorphisms, \ie degree-wise invertible.
\end{tabular}
\end{center}

The category of graded vector spaces inherits the symmetric monoidal structure from the category of ordinary vector spaces.
Namely, we define the \emph{tensor product of two graded vector spaces} $V$ and $W$ as the graded vector space
\begin{displaymath}
	V\otimes W :=
	\left(
		k \mapsto \bigoplus_{i+j=k} V_i\otimes W_j
	\right)
~.
\end{displaymath}
The action of the tensor product functor on any given graded morphisms $f_i:V_i\to W_i$, \ie degree $0$ homogeneous maps, is a graded morphism given by
		\begin{displaymath}
			\morphism{f_1\otimes f_2}
			{V_1\otimes V_2}
			{W_1 \otimes W_2}
			{u_1\otimes u_2}
			{f_1(u_1)\otimes f_2(u_2)}
		\end{displaymath}
The respective \emph{symmetric braiding} is given, for any graded vector space $V$ and $W$, by the graded linear isomorphism $B_{V,W}$, defined on homogeneous elements by
\begin{displaymath}
	\morphism{B_{V,W}}
	{V\otimes W}
	{W \otimes V}
	{x\otimes y}
	{(-)^{|x||y|} y \otimes x}
~.
\end{displaymath}
The choice of that sign prefactor in the definition of the braiding takes the name of \emph{Koszul convention}.
\\
Out of this convention (see section \ref{sec:bloodyKoszulConvention}), one can introduce the \emph{d\'ecalage} isomorphism, defined on any $n$-tuple of graded vector spaces $(V_1,\dots, V_n)$ as:
	\begin{equation}\label{eq:deca}
		\isomorphism{\dec}
		{V_1[1]\otimes\dots\otimes V_n[1]}
		{(V_1\otimes\dots \otimes V_n)[n]}
		{v_{1[1]}\otimes\dots v_{n[1]}}
		{(-)^{\sum_{i=1}^{n}(n-i)|v_i|}(v_1\otimes\dots\otimes v_n)_{[n]}}
		~.
	\end{equation}

For any homogeneous maps $f \in \underline{\Hom}^{|f|}(V,W), f' \in \underline{\Hom}^{|f'|}(V',W')$ we define, 
with a slight abuse of notation (see notation \ref{not:abuseotimes}),
$f\otimes f' \in \underline{\Hom}^{|f|+|f'|}(V\otimes V', W \otimes W')$ as the homogeneous map acting on homogeneous elements as
		\begin{equation}\label{eq:KoszulConvTensorProducts}
			\morphism{f~{\otimes}~ f'}
		 {V\otimes V'}
		 {(W\otimes W')[|f|+|f'|]}
		 {v \otimes v'}
		 {(-)^{|v||f'|}(f(v)\otimes f'(v'))_{[|f|+|f'|]}}
		 ~.
		\end{equation}
		Accordingly, there is also a sign rule for the composition of tensor products of homogeneous maps:
		\begin{equation}\label{Eq:TensorHomogeneousMaps}
			(f'\otimes g') \circ (f \otimes g) = (-)^{|g'||f|}(f'\circ f)\otimes (g'\circ g)
			~.
		\end{equation}

Since $\Vect^{\ZZ}$ is a symmetric monoidal category, for any graded vector space $V$ and for any positive integer $n$ there are two canonical representation of the symmetric group $S_n$ over $V^{\otimes n}=\otimes^n V$.
We denote by $B_\sigma$ and $P_\sigma$ the even, respectively odd, representation of the permutation $\sigma\in S_n$ on a given $V^{\otimes n}$, namely
	\begin{displaymath}
		\morphism{B_\sigma}
		{V^{\otimes n}}
		{V^{\otimes n}}
		{x_1\otimes \dots \otimes x_n}
		{\epsilon(\sigma;x_1,\dots,x_n) x_{\sigma_1}\otimes \dots \otimes x_{\sigma_n}}
	\end{displaymath}
and $P_\sigma = (-)^\sigma B_{\sigma}$.
The coefficient $\epsilon(\sigma;x_1,\dots,x_n)$ is the so-called \emph{Koszul sign} and it is defined as the sign of the sub-permutation of $\sigma$ involving only elements in odd-degree (see remark \ref{rem:aboutKoszulSign} and \ref{rem:practicalComputKoszulSign} for further details).
We will usually omit the dependence on the list of graded vectors since it should appear clear from the context. Namely we will often abuse the notation $\epsilon(\sigma;x_1,\dots,x_n)$ by writing $\epsilon(\sigma)$, and writing
 $\chi(\sigma):= \epsilon(\sigma)(-)^{\sigma}$ for the sign involved in the definition of $P_\sigma$, also called \emph{odd Koszul sign}.
 \\
Given any subset $I\subset S_n$ of the group of permutations of $n$-elements, we denote by $B_I$ and $P_I$ the operator giving the sum on all the elements of $S$:
	\begin{displaymath}
		B_I = \sum_{\sigma \in I} B_\sigma~,\qquad P_I=\sum_{\sigma \in I} P_\sigma~
		~,
	\end{displaymath}
we will often refer to them as the even and odd \emph{permutator of $I$}.
Observe that the d\'ecalage isomorphism determine a natural transformation between the even and odd representation of $S_n$, namely the following diagram commutes in the graded vector spaces category (see remark \ref{rem:decVsBraiding} and proposition \ref{Prop:DecalageOfPermutation}):
			\begin{displaymath}
				\begin{tikzcd}
				(V[1])^{\otimes n} \ar[r,"\dec"] \ar[d,"B_{\sigma}"'] 
				&
				(V^{\otimes n})[n] \ar[d,"(P_{\sigma}){[n]}"] 
				\\
				(V[1])^{\otimes n} \ar[r,"\dec"] 
				&
				(V^{\otimes n})[n]
				\end{tikzcd}
				~.
			\end{displaymath}	
We denote by $\odot^n V = V^{\odot n}$ and $\Lambda^n V = V^{\wedge n}$ the spaces of coinvariant elements with respect to the even and odd representation of $S_n$. One has the following splitting sequence in the category of graded vector spaces
			\begin{equation}\label{eq:symskewsplitting}
				\begin{tikzcd}
					0 \ar[r, shift left =.5ex] 
					&[-5ex]
					\Lambda^n V \ar[r, shift left =.75ex,hookrightarrow,"N_a"] 
					&
					\bigotimes^n V 
					\ar[r,two heads, shift left =.75ex,two heads,"\pi_s"] 
					\ar[l,two heads, shift left =.75ex,two heads,"\pi_a"]
					&
					\bigodot^n V 
					\ar[r]			
					\ar[l, shift left =.75ex,hookrightarrow,"N_s"]
					&[-5ex] 0
				\end{tikzcd}~.
			\end{equation}	
We refer to proposition \ref{lemma:splitsequencebrutta} in appendix (specifically to equation \eqref{eq:SymSkewOperatorsDef-appendix}) for the explicit definition of all these mappings. We only mention here that, on the symmetric side, one has
	\begin{equation}\label{eq:SymSkewOperatorsDef}
		\begin{aligned}[c]
			\pi_s:&~ x_1\otimes\dots\otimes x_n \mapsto x_1\odot\dots\odot x_n ~,
			\\
			{N_s}:&~ x_1\odot\dots\odot x_n \mapsto 
			\left(\sum_{\sigma\in S_n}x_{\sigma_1}\otimes\dots \otimes x_{\sigma_n}\right) ~,
			\end{aligned}
	\end{equation}
	and their composition yields the \emph{(graded) symmetrizator operator}
	\begin{displaymath}
				\symAtor_{(n)}:= \frac{1}{n!}~N_s\circ \pi_s \equiv \left(\sum_{\sigma \in S_n} \frac{1}{n!} B_\sigma \right)~.
	\end{displaymath}
The d\'ecalage isomorphism restrict compatibly with these space of coinvariants, namely the following diagram commutes (see lemma \ref{Lemma:DecalageRestrictToSymTens}):
\begin{equation}\label{eq:decalageRestriction}
		\begin{tikzcd}[column sep = 6em,row sep =small]
			(V[1])^{\odot n} \ar[r,"dec\big|_{V[1]^{\odot n}}"] \ar[d,hook]&
			(V^{\wedge n})[n]\ar[d,hook]
			\\
			(V[1])^{\otimes n} \ar[r,"\dec"] &
			(V^{\otimes n}) [n]
			\\
			(V[1])^{\wedge n} \ar[r,"dec\big|_{V[1]^{\wedge n}}"'] \ar[u,hook]&
			(V^{\odot n})[n] \ar[u,hook]
		\end{tikzcd}
	\end{equation}
	\subsection{Tensor algebras and coalgebras}
	Given a graded vector space $V$ we denote the  tensor, symmetric tensor, and skew-symmetric tensor spaces\footnote{In the literature these spaces are usually called \emph{reduced tensor spaces} (see \eg \cite{Manetti-website-coalgebras,Reinhold2019}) in order to distinguish them from the (augmented) tensor spaces. The latter are defined by a similar summation with index $n$ starting from $0$ and by implicitly assuming that ${T(V)}^0=S(V)^0=\Lambda(V)^0 = \RR$.
		We do not need this case in the body of the thesis, the (augmented) case is only mentioned in appendix \ref{App:GradedMultilinearAlgebra}.} 
	 respectively as the graded vector spaces:
	\begin{displaymath}
		\overline{T(V)} := \bigoplus_{n\geq 1} V^{\otimes n}
		~,\qquad
		\overline{S(V)} := \bigoplus_{n\geq 1} V^{\odot n}
		~,\qquad
		\overline{\Lambda(V)} := \bigoplus_{n\geq 1} V^{\wedge n}
		~.
	\end{displaymath}
	At times, we will also adopt the notation $\overline{S(V)}=S^{\geq 1}(V)$ to stress the fact that the direct summutation runs on non-negative degrees.
	Observe that, according to this definition, $\overline{T(V)}^k = V^{\otimes k}$ if and only if $V$ is a graded vector space concentrated in degree $1$.
	These spaces carry a natural structure of (not unital) graded associative algebra given by the action of $\otimes,\odot$ and $\wedge$ on vectors, and a (not counital) coassociative coalgebra structures given by deconcatenation. 
	In particular $\overline{S(V)}$ and $\overline{\Lambda(V)}$ are graded (co)commutative and graded (co)anticommutative respectively.
	We call \emph{decomposable element} in the Tensor algebra any element that can be expressed as tensor product of given vectors.

	In what follows we will essentially only concerned with the coalgebra structure on $\overline{S(V)}$, we refer to appendix \ref{App:GradedMultilinearAlgebra}, and references therein, for a more exhaustive discussion of the subject.
	For the rest of this subsection we will only state some basic notations\footnote{We point out that that the notation employed here slightly depart to the more decorated notation  used in the appendix. See notation \ref{not:droppingthesuperc}.}.
	We will denote, with a slight abuse of language, the \emph{deconcatenation} and \emph{unshuffled deconcatenation} with the same symbol $\Delta$ (see definitions \ref{def:Decon} and \ref{def:unshuffleDecon} for the explicit expressions). 
	We recall that a \emph{coalgebra morphism} from $\overline{T(V)}$ to $\overline{T(W)}$ is a graded vector space morphism $f:\overline{T(V)}\to \overline{T(W)}$, hence a degree $0$ homogeneous map, compatible with the coproduct.
	We denote by
	\begin{displaymath}
		\Hom_{coAlg}(\overline{T(V)},\overline{T(W)}):=
		\left\lbrace\left.
			F \in \Hom(\overline{T(V)},\overline{T(W)}) ~ \right\vert ~
			\Delta\circ F = F\otimes F \circ  \Delta
		\right \rbrace
	\end{displaymath}
	the graded set of all coalgebra morphisms. We stress that it does not respect the linear structure of $\Hom(\overline{T(V)},\overline{T(W)}$ thus it is only a graded subset.
	\\
	Fixed a coalgebra morphism $F\in\Hom_{coAlg}(\overline{T(V)},\overline{T(W)})$, we recall that a \emph{degree $k$ $F$-coderivation} (see section \ref{sec:abtractFcoderivations}) of the coalgebra $\overline{T(V)}$ is a degree $k$ homogeneous linear map $Q\colon \overline{T(V)} \to \overline{T(W)}$ such that $\Delta\circ Q=(Q \otimes F + F\otimes Q)\circ \Delta$.
To explain the terminology, notice that this equation is what one obtains dualizing the property of being a derivation of an algebra.	
We denote the graded vector space of coderivations along $F$ as
	\begin{displaymath}
		\coDer(\overline{T(V)},\overline{T(W)};F) ~.
	\end{displaymath}
	We denote by $\coDer(\overline{T(V)})$ the special case where $W=V$ and $F$ is given by the identity.
	Elements in $\coDer(\overline{T(V)})$ will be simply called \emph{coderivations}.
	The same definitions applies, mutatis mutandis, also to the symmetric coalgebra case.

	Tensor (co)algebras enjoy the special properties to be (co)-free objects with respect to certain subcategories of graded (co)-algebras.
	This feature can be expressed by universal properties that can in turn be translated into the following lemmas (stated for the symmetric case but analogue results hold for $\overline{T(V)}$ and $\overline{\Lambda(V)}$):
\begin{lemma}[Lift to a coalgebra morphism (\emph{Prop. \ref{Prop:UniversalPropertyCocommutativeGradedCoalgebras} and Thm. \ref{Theorem:HomCoAlgISOSymmetric}})]\label{lem:liftToMorph}
   There is a bijection between morphisms of coalgebras $F:\overline{S(V)}\to \overline{S(W)}$ and morphisms of graded vector spaces $f\colon: \overline{S(V)} \to W$.
   Namely the following graded sets are isomorphic
   \begin{displaymath}
   	\begin{tikzcd}[ column sep=1em,row sep=-1ex]
   		\Hom_{coAlg}(\overline{S(V)},S(W)) \ar[r,equal,"\sim"] & 
   		\Hom(\overline{S(V)},W) \ar[r,equal,"\sim"] &
   		\bigoplus_{n\geq 1}\Hom(V^{\odot n},W)
   		\\
   		F \ar[r,mapsto]
   		&
   		\pr_W \circ F \ar[r,mapsto]
   		&
   		(\pr_W \circ F\eval_{V},\dots,\pr_W \circ F \eval_{V^{\odot n}},\dots)
   		\\
   		L_{\sym}(f) &
   		f= \bigoplus_{i\geq 1} f_i \ar[l,mapsto] &
   		(f_1,\dots,f_n,\dots) \ar[l,mapsto]
   	\end{tikzcd}
   \end{displaymath}
   where the direct function is called \emph{corestriction} and is obtained by postcomposition with the standard projection $\pr_W: \overline{T(W)}\twoheadrightarrow W$,
   and the inverse is called \emph{(unique) lift to a coalgebra morphism} and is given explicitly by 
\begin{equation}\label{Eq:LiftMorphismMapping-noappendix}
	   	\hspace{-.025\textwidth}
			\morphism{L_{\sym}}
			{{\Hom}({\overline{S(V)}},W)}
			{\Hom_{\text{coAlg}}({\overline{S(V)}},{S(W)})}
			{f}
			{\displaystyle
			\sum_{n>0}\sum_{s=1}^n 	
			\pi_s \circ \left[\mkern-60mu\sum_{\mkern75mu \substack{i_1+\dots+i_s = n \\ 0<i_1\leq i_2 \leq \dots \leq i_s}}\mkern-55mu 
			(f_{i_1}\otimes \dots \otimes f_{i_s}) \circ B^{<}_{i_1,\dots,i_s} \right]\circ N_s~}
\end{equation}	
		$N_s,\pi_s$ are the operator defined by equation \eqref{eq:SymSkewOperatorsDef} and $B^<_{i_1,\dots,i_n}\equiv B_{\ush{i_1,\dots,i_n}^<}$ denotes the sum on all the ordered $(i_1,\dots,i_n)$-unshuffles (ordered unshuffleator, see definition \ref{Def:OrderedUnshuffleator}), \ie the sum 
runs through all $(k_1,\cdots,k_\ell)$-unshuffles $\sigma$ satisfying the extra condition
	\begin{displaymath}
		\sigma(k_1+\dots+k_{j-1}+1)<\sigma(k_1+\dots+k_{j}+1) 
		\quad \text{if}~k_{j-1}=k_j ~.
	\end{displaymath}		
\end{lemma}
Given a coalgebra morphism $F:\overline{S(V)}\to \overline{S(W)}$, we call its \emph{$k$-th} corestriction the symmetric multilinear operator
\begin{displaymath}
	f_k= \pr_W \circ F \eval_{V^{\odot k}}~: V^{\odot k} \to W~.
\end{displaymath}
\begin{notation}\label{rem:operatorSln}
	The linear operator between square brackets in equation \eqref{Eq:LiftMorphismMapping-noappendix} will appear again in our constructions. We single out the following definition for future reference.
	\\
	Given a $f\in\Hom(\overline{T(V)},W)$, eventually expressible as a collection $(f_i,\dots,f_k,\dots)$ with $f_k\in \Hom(V^{\otimes k},W)$, 
	we denote by $\mathcal{S}_{\ell,m} (f)$, 	for any $1\leq\ell\leq m $,  the graded multilinear map $V^{\otimes m} \to W^{\otimes \ell}$ defined as
	\begin{equation}\label{Eq:Soperator}
		\mathcal{S}_{\ell,m} (f) =
		\left( 
			\sum_{\substack{k_{1}+\cdots+k_{\ell}=m\\1\leq k_{1}\leq\cdots\leq k_{\ell}}}
			(f_{k_1}\otimes\cdots\otimes f_{k_\ell})\circ 
			B_{k_1,\ldots,k_\ell}^<		
		\right)
		~.
	\end{equation}
	With this notation one has:
	\begin{displaymath}
			L_{\sym}(f) = 	\sum_{m>0}\sum_{\ell=1}^m 	
			\pi_s \circ \Sop_{\ell,m}(f)\circ N_s
			~.
	\end{displaymath}
	\note{Forse la scelta di questo simbolo non è felice e crea confusione con i simmetrizzatori $\symAtor$}
\end{notation}
Equation \eqref{Eq:LiftMorphismMapping-noappendix} permits to lift uniquely any graded morphisms, \ie degree $0$ homogeneous maps, to a coalgebra morphisms.
A similar construction applies also to homogeneous maps in arbitrary degree, yielding this time a unique lift to a coderivation:
\begin{lemma}[Lift to a coderivation {\cite[Lemma 2.4]{LadaMarkl}}]\label{lem:liftToCoder}
	\note{Dalla appendice vedo che questo lemma si può esprimere per ogni $f:\overline{S(V)}\to W$, ma le espressioni sono ancora piu' complicate}
  Given any graded morphism $f_1:V\to W$, denote by $f:\overline{S(V)}\to W$ its trivial extension to the whole symmetric tensor space, there exists a bijection between degree $k$  coderivations $Q:\overline{S(V)}\to \overline{S(W)}$ and degree $k$ linear maps $q\colon \overline{S(V)}\to W$.
   Namely the following graded vector spaces are isomorphic
   \begin{displaymath}
   	\hspace{-.075\textwidth}
   	\begin{tikzcd}[ column sep=2em,row sep=-1ex]
   		\coDer(\overline{S(V)},\overline{S(W)};L_{\sym}(f)) \ar[r,equal,"\sim"] & 
   		\underline{\Hom}(\overline{S(V)},W) \ar[r,equal,"\sim"] &
   		\bigoplus_{n\geq 1}\underline{\Hom}(V^{\odot n},W)
   		\\
   		Q \ar[r,mapsto]
   		&
   		\pr_W \circ Q \ar[r,mapsto]
   		&
   		(\pr_W \circ Q\eval_{V},\dots,\pr_W \circ Q \eval_{V^{\odot n}},\dots)
   		\\
   		\widetilde{L}_{\sym}(q) &
   		q= \bigoplus_{i\geq 1} q_i \ar[l,mapsto] &
   		(q_1,\dots,q_n,\dots) \ar[l,mapsto]
   	\end{tikzcd}
   \end{displaymath}
   where the direct function is the \emph{corestriction} as in lemma \ref{lem:liftToCoder} and inverse is called \emph{(unique) lift to a $F$-coderivation}, with $F=L_{\sym}(f)$, and is given explicitly by 
\begin{equation}\label{Eq:LiftCoDerMapping-noappendix}
	   	\hspace{-.025\textwidth}
			\morphism{\widetilde{L}_{\sym}}
			{\underline{\Hom}^k({\overline{S(V)}},W)}
			{\coDer^k({\overline{S(V)}},{\overline{S(W)}};F)}
			{q}
			{\displaystyle
			\sum_{n>0}\sum_{s=1}^n 	
			\pi_s \circ \left[\left(
			q_{n-s+1}\otimes f^{\otimes(s-1)}\right) \circ B_{n-s+1,s-1}
			\right]\circ N_s~}
		\end{equation}	
		where $N_s,\pi_s$ are the operators defined by equation \eqref{eq:SymSkewOperatorsDef} and $B_{n-s+1,s-1}\equiv B_{\ush{n-s+1,s-1}}$ denotes the sum on all $(n-s+1,s-1)$-unshuffles.
\end{lemma}
	Note that the term between square brackets in equation \eqref{Eq:LiftCoDerMapping-noappendix} is a graded linear operator $V^{\otimes n}\to W^{\otimes s}$. In the case that $V=W$ and $f_1=\Unit$, one has that $F=\Unit$ and $f^{\otimes n} = \mathbb{1}_n$; the corresponding term play a role in the definition of the \emph{\RN product}.
	\\
	More explicitly, for any given degree $k$ homogeneous map $m\in \underline{\Hom}^k(\overline{S(V)},V)$, the action of its lift on homogeneous elements is given by:	
	\begin{equation}\label{eq:coder}
	\widetilde{L}_{\sym}(m)~(x_1,\dots,x_n):= \sum_{i=1}^n 
	\mkern-50mu \sum_{\mkern70mu\sigma \in \ush{i,n-i}}\mkern-50mu
	\epsilon(\sigma)~ m_{i}(x_{\sigma_1},\dots,x_{\sigma_{i}})\odot x_{\sigma_{i+1}}\odot\dots\odot x_{\sigma_n}
	~,
\end{equation}
where $\ush{i,j}$ denotes the subgroup of $(i,j)$-unshuffles in the permutation group $S_{i+j}$.

	The use of the term "lift" follows from the commutativity of the following diagrams in the category of graded vector spaces:
	\begin{displaymath}
	\begin{tikzcd}[column sep =huge]
		\overline{S(V)} \ar[r,"F=L_{\sym}(f)"] \ar[dr,"f"']& \overline{S(W)} \ar[d,"\pr_W"]\\
		&W	
	\end{tikzcd}
	\quad,\qquad
	\begin{tikzcd}[column sep =huge]
		\overline{S(V)} \ar[r,"Q=\widetilde{L}_{\sym}(q)"] \ar[dr,"q"']
		& \overline{S(W)}[|q|] \ar[d,"\pr_W{[|q|]}"]\\
		&W[|q|]	
	\end{tikzcd}
	\quad .
	\end{displaymath}
	Observe that for degree $0$ linear maps, \ie graded morphisms, both the lift to a coalgebra morphism and the lift to a coderivation are well-defined.
	This specific case requires to use the decorated notation $L$ and $\widetilde{L}$ to distinguish between the two possible lifts.

\subsection{Differential graded vector spaces}\label{sec:HomologicalAlgebrasConventions}
	We call a \emph{differential graded vector space} any pair $(V,\d)$ composed of a graded vector space $V\in \Vect^\ZZ$ together with a $2$-nilpotent homogeneous map in degree $1$ from $V$ into itself, \ie $\d \in \underline{\End}^1(V)$ and $\d\cdot \d = 0$.
	
	Most of the time we will make use of the language of homological algebra. 
	In those terms, a differential graded vector space can be simply seen as a \emph{cochain complex} \ie as a sequence of vector spaces $V^k$ together with operators $d^{(k)}\equiv d\eval_{V^k}: V^k\to V^{k+1}$ such that $d^{(k+1)}\circ d^{(k)}=0$. 
	Diagrammatically, the differential graded vector space $(V,\d)$ will be depicted as 
	\begin{displaymath}
		\begin{tikzcd}
			\dots \ar[r]
			&
			V^{k-1} \ar[r,"\d"]
			&
			V^{k} \ar[r,"\d"]
			&
			V^{k+1} \ar[r]
			&
			\dots
		\end{tikzcd}
		~.
	\end{displaymath}
	Accordingly, operator $\d$ will be called \emph{coboundary operator} and  homogeneous elements $v\in V^k$ will be called \emph{$k$-cochains}.
	We will refer to elements in $(\ker(d))^n\subset V^n$ as \emph{$n$-cocycles} and elements in $(\Im(\d))^n \subset V^n$ as \emph{$n$-coboundaries}.
		We will also employ the shorthand notation $B(V):=\Im(d)$ and $Z(V):=\ker(d)$
	for the graded vector spaces of coboundaries and cocycles. Clearly $B(V)\subset Z(V)$.
	We will call \emph{the cohomology of $(V,\d)$} the graded vector spaces
	\begin{displaymath}
		H(V,\d)= \dfrac{\ker(d)}{\Im(d)} :=
		\left( k \mapsto \dfrac{\ker(d^{(n)})}{\Im(d^{(n-1)})}\right)
	\end{displaymath}
	and its $k$-component $H^k(V,\d)$ will be called \emph{$k$-th cohomology group of $(V,\d)$}.

	Recall that in homological algebra there is the dual notion of \emph{chain complex} defined by a $2$-nilpotent boundary operator $\partial$ with decrease the grading. Therefore we are tacitly adopting the \emph{cohomological convention} when dealing with differential graded vector space.
	One can simply get a "homological" differential graded vector space applying the reverse ordering functor $\setminus:\Vect^\ZZ\to \Vect^\ZZ$ to the ("cohomological") differential graded vector space $(V,\d)$.

	We will call a \emph{chain map} between two cochain complexes $(V_1,\d_1)$ and $(V_2,\d_2)$ any graded morphism $f:V_1\to V_2$ which commutes with the coboundary operators, that is, the following diagram commutes in the category of (ordinary) vector spaces:
	\begin{displaymath}
		\begin{tikzcd}
			\dots \ar[r]
			&
			V^{k-1}_1 \ar[r,"\d_1"]\ar[d,"f^{k-1}"]
			&
			V^{k}_1 \ar[r,"\d_1"]\ar[d,"f^k"]
			&
			V^{k+1}_1 \ar[r]\ar[d,"f^{k+1}"]
			&
			\dots
			\\			
			\dots \ar[r]
			&
			V^{k-1}_2 \ar[r,"\d_2"]
			&
			V^{k}_2 \ar[r,"\d_2"]
			&
			V^{k+1}_2 \ar[r]
			&
			\dots
		\end{tikzcd}
		~.
	\end{displaymath}
	More synthetically, the latter diagram can be seen as a commutative square in the category $\GVect$ of graded vector spaces:
	\begin{displaymath}
		\begin{tikzcd}
			V_1 \ar[r,"\d_1"] \ar[d,"f"'] & V_1 \ar[d,"f"]
			\\
			V_2 \ar[r,"\d_2"] & V_2
		\end{tikzcd}
		~.
	\end{displaymath}
	We will call \emph{chain homotopy} from a chain map $h:V_1 \to V_2$ to a chain map $f:V_1 \to V_2$ a degree $-1$ homogeneous map $\varphi\in \underline{\Hom}^1(V_1,V_2)$ such that 
\begin{equation}\label{eq:chainhomotopy}
	\d_2 \circ \varphi + \varphi \circ \d_1 = f-h
	~.
\end{equation}
Diagrammatically, we will denote a chain homotopy as a $2$-morphism:
\begin{displaymath}
	\begin{tikzcd}[column sep = huge]
		 V_1 \ar[r,blue,"h", bend left=20, ""{name=U, below}]
		\ar[r,"f"', bend right=20, ""{name=D}]
		&[1em] V_2
		\ar[Rightarrow, purple,"\varphi", from=U, to=D]
\end{tikzcd}		
\end{displaymath}
	
	Consider now a differential graded vector space $C=(C,\d)$, if $C$ is a (counital, cocommutative) coalgebra and $d\in \coDer^1(C)\subset \End^1(C)$ is a degree $1$ coderivation, we will talk of a \emph{graded differential coalgebra} and call $\d$ a \emph{codifferential} rather than a coboundary operator.
	When the above algebra $C$ is also cocommutative and \emph{cofree}, \ie $C\cong \overline{S(V)}$ for a certain graded vector space $V$, one talk about a \emph{cofree graded differential coalgebra}.	This last notion will play a central role in the coalgebraic presentation of a $L_\infty$-algebra.

	We will call a \emph{(cochain) bicomplex} any $(\ZZ\times\ZZ)$-graded (or \emph{bi-graded}) vector space 
	\begin{displaymath}
		V=\left( (k,\ell)\mapsto V^{k,\ell}\right)
	\end{displaymath}
	together with a pair of $2$-nilpotent homogeneous map $\d_v \in \underline{\End}^{(1,0)}(V)$ and $\d_h \in \underline{\End}^{(0,1)}(V)$ such that 
	\begin{displaymath}
		\d_v\cdot \d_h + \d_h \circ \d_v=0
		~.
	\end{displaymath}		
	Namely, this is given by the following commuting diagram
	 \begin{displaymath}
	 	\begin{tikzcd}[column sep = huge]
	 		& \vdots & \vdots & \\
	 	\cdots \ar[r,"(-)^{i+1}\d_h"]& V^{i+1,j} \ar[u,"\d_v"]\ar[r,"(-)^{i+1}\d_h"]& V^{i+1,j+1} \ar[r,"(-)^{i+1}\d_h"]\ar[u,"\d_v"]& \cdots \\
	 	\cdots \ar[r,"(-)^i\d_h"]& V^{i,j} \ar[u,"\d_v"]\ar[r,"(-)^i\d_h"]& V^{i,j+1}\ar[r,"(-)^i\d_h"] \ar[u,"\d_v"]& \cdots \\
		 & \vdots \ar[u,"\d_v"]& \vdots \ar[u,"\d_v"]
	 	\end{tikzcd}
	 \end{displaymath}
	where each row and column are separately cochains complexes.
	(In other words, we are adopting the "anticommuting square" convention with the wording of \cite{nlab:double_complex}).
	\\
	Given a bicomplex $V^{\bullet,\bullet}=(V,\d_h,\d_v)$ we call \emph{total complex} of $V^{\bullet,\bullet}$ the $\ZZ$-graded vector space
	\begin{displaymath}
		\tot(V) = \left( k \mapsto \bigoplus_{i+j=k}V^{i,j}\right)
	\end{displaymath}
	together with the coboundary operator $\d_{tot} = \d_h + \d_v$.
	\\
	Given two differential graded vector spaces $(V,\d_v)$ and $(H,\d_h)$, 
	we call \emph{tensor product complex}, denoted as $(V,\d_v)\otimes (H,\d_h)$, 
	the total complex corresponding to the bigraded vector space
	\begin{displaymath}
		(V\otimes W)^{\bullet,\bullet} :=
		\left(
			(k,\ell) \mapsto V^k\otimes W^\ell
		\right)		
	\end{displaymath}
together with the coboundary operators $\d_v\otimes \Unit_H$ and $\Unit_V\otimes \d_h$. 
Observe that the \emph{Koszul convention} implies that the action of $\d_{tot}$ on an homogeneous element $x\otimes y\in V\otimes W$ reads as follows:
\begin{displaymath}
	\d_{tot}(x\otimes y) = (\d_v(x))\otimes y + (-)^{|x|} x\otimes (\d_h(y)) 
	~.
\end{displaymath}
Iterating this procedure, we will call the \emph{$n$-iterated tensor product} of a cochain complex $(V,\d)$ the cochain complex $(V^{\otimes n},\d_{\otimes n})$ where
	\begin{displaymath}
		V^{\otimes n} = \left(
			k \mapsto \bigoplus_{i_1+\dots+i_n = k} V^{i_1}\otimes V^{i_2}\otimes \dots V^{i_n}
			\right)
	\end{displaymath}
and
\begin{equation}\label{eq:totalcoboundaryoperator}
	\d_{\otimes n} = \sum_{i=1}^n \Unit_{i-1}\otimes \d \otimes \Unit_{n-i}
	~.
\end{equation}

\subsection{Multibrackets and Coderivations}\label{sec:RNstuff}
 Consider a graded vector space $V$, for any $n\geq 0$ and $k \in \ZZ$ we denote by 
	\begin{displaymath}
		M_{n,k}(V,W):= \underline{\Hom}^k (V^{\otimes n},W)	
	\end{displaymath}	
	the graded vector space of degree $k$ (homogeneous) $n$-multilinear maps. 
	We take for granted, and kept implied, the universal property of multilinear maps, hence we will be free to understand elements of $M_{n,k}(V,W)$ as $n$-ary functions $V\times\dots \times V \to W[k]$ with the extra property of being separately linear in each entry. 
	Accordingly, homogeneous linear maps from $V^{\otimes n}$ to $W$ will be said \emph{of arity $n$}, and we will often denote the image of a multilinear map $\mu_n$ on $x_1\otimes\dots\otimes x_n$ as $\mu_n(x_1,\dots,x_n)$	separating elements by commas and omitting the symbol $\otimes$.
	 
	The same applies to graded symmetric and graded skew-symmetric multilinear maps, we denote by 
	\begin{displaymath}
		M_{n,k}^{\sym}(V,W):= \underline{\Hom}^k (V^{\odot n},W)
		~,\quad
		M_{n,k}(V,W)^{\skew}:= \underline{\Hom}^k (V^{\wedge n},W)
		~,
	\end{displaymath}
	respectively the spaces of degree $k$  symmetric and skew-symmetric  $n$-multilinear homogeneous maps on $V$ {with values in $W$}. 
	It follows from the splitting sequence \eqref{eq:symskewsplitting} that $M_{n,k}(V,W)=M_{n,k}^{\sym}(V,W)\oplus M_{n,k}^{\skew}(V,W)$.
	 When $W=V$ we will lighten the notation omitting the second entry.
 
 Considering all the possible arities and degrees collectively, and neglecting the "internal" grading of the graded vector space $\underline{\Hom}(V^{\otimes n},W[k])$, one obtains the $(\NN_0\times \ZZ)$-graded (bi-graded) vector space $M_{\bullet,\bullet}(V,W)$. The same reasoning applies also to the subspace of (skew)-symmetric multilinear maps.
 \\
	The d\'ecalage isomorphism, {introduced in eq. \eqref{eq:deca}}, induces by precomposition an isomorphism of graded vector spaces
	\begin{displaymath}
		\isomorphism{\Dec}
		{\underline{\Hom}^k(V^{\otimes n}, W)}
		{\underline{\Hom}^{k+n-1}(V[1]^{\otimes n}, W[1])}
		{\mu}
		{\Dec(\mu):= \mu[n]\circ \dec}
	\end{displaymath}
	with
	\begin{equation}\label{eq:ExplicitDECA}
		\Dec(\mu)(x_{1\,[1]},\dots,x_{n\,[1]})
		=
		(-)^{\sum_{i=1}^{n}(n-i)|x_1|}
		\left( \mu(x_1,\dots,x_n)\right)_{[n]}
		~.
	\end{equation}
	Slightly more in general, we will also intend the d\'ecalage of a map $f\in \underline{Hom}(V^{\otimes n},W^{\otimes m})$ as $\Dec(f):= \dec^{-1}[|f|+n-m]\circ f[n] \circ \dec$.
	It is possible to read the operator $\Dec$ as a genuine graded isomorphism (not bi-graded) by appropriately contracting indices $n$ and $k$ to give a certain $\ZZ$-grading.
	Conventionally, we introduce the \emph{graded vector spaces of graded symmetric and graded skew-symmetric multilinear maps} from $V$ to $W$ as 
 \begin{equation}\label{eq:RNspaces}
	\begin{aligned}
  		M^{sym}(V):=&~ 
  		\left( k \mapsto \bigoplus_{n} M^{\sym}_{n, k}(V)\right)
  		~, 
  		\\
  		M^{skew}(V) :=&~ 
  		\left( k \mapsto \bigoplus_{n+i=k+1} M^{\skew}_{n, i}(V)\right) 
  		~.
  	\end{aligned}
 \end{equation}
 According to this choice, the d\'ecalage operator defined above restricts to a well-defined isomorphism of graded vector spaces (see remark \ref{rem:DecasGradedMorph})
 \begin{displaymath}
 		Dec: M^{skew}(V) \xrightarrow{\quad \sim \quad} M^{sym}(V[1])
			~.
 \end{displaymath}
 
	We finally notice that in equation \eqref{eq:RNspaces} is implied the choice of two possible $\ZZ$-grading on  $M(V)$. Namely, given an homogeneous multilinear map $\mu \in M^{n,k}(V,W)$, we introduce the two gradings 
	\begin{displaymath}
		|\mu|=k ~,\qquad  ||\mu||=k+n-1 ~,
\end{displaymath}	
	the first is the degree in the sense of homogeneous maps and the second is the grading as an element of 
	\begin{displaymath}
	\left(
 				k \mapsto \bigoplus_{n+m=k+1}M_{n,k}(V,W)
 			\right) ~.
\end{displaymath}

\subsubsection{\RN product for multibrackets}\label{sec:RNProdMB}	 
 	The graded vector spaces of graded symmetric and graded skew-symmetric multilinear maps constitute graded right pre-Lie algebras (see appendix \ref{App:PreLie}) $(M^{sym}(V),\cs)$ and $(M^{skew}(V),\ca)$ when endowed with the so-called \emph{\RN} product (see \cite{Nijenhuis1967} for the original definition on ordinary vector spaces and \cite[Thm. 3.3]{Lecomte1992} or \cite[\S 1.2]{Delgado2018b} for the graded case - note that our definition differs for a sign-; the formula for the graded symmetric case can be also found in \cite[\S 1.1]{Bandiera2016}).
 Denoting by $B_{i_1,\dots, i_\ell}$ and $P_{i_1,\dots,i_\ell}$ the operators summing on all even and odd actions of permutations in $\ush{i_1,\dots,i_\ell}$ (see "conventions" section on page \pageref{Sec:conventions} or appendix \ref{App:UnshuffleAtors}), i.e
	\begin{displaymath}
		\begin{aligned}
		B_{i_1,\dots,i_\ell} ~\left( x_1\otimes x_2 \otimes \dots \right)
		=& \mkern-30mu\sum_{\mkern50mu\sigma \in \ush{i_1,\dots,i_\ell}} \mkern-30mu
		\epsilon(\sigma) x_{\sigma_1}\otimes x_{\sigma_2} \otimes \dots	
		\\
		P_{i_1,\dots,i_\ell} ~\left( x_1\otimes x_2 \otimes \dots \right) 
		=& 
		\mkern-30mu \sum_{\mkern50mu\sigma \in \ush{i_1,\dots,i_\ell}} \mkern-30mu
		\epsilon(\sigma)(-)^\sigma x_{\sigma_1}\otimes x_{\sigma_2} \otimes \dots	
		~,
		\end{aligned}
	\end{displaymath}	  
	the \RN product can be succinctly written as
	\begin{equation}\label{Eq:RNProducts-myway}
		\begin{aligned}
			\mu_n \cs \mu_m	
			=& 
			\mu_n \circ (\mu_m \otimes \mathbb{1}_{n-1}) \circ B_{m,n-1}	
			\\
			\mu_n \ca \mu_m	
			=& 
			(-)^{|\mu_m|(n-1)}
			~\mu_n \circ (\mu_m \otimes \mathbb{1}_{n-1}) \circ P_{m,n-1}
			~.
		\end{aligned}
	\end{equation}
 	More explicitly, evaluating on homogeneous element $x_i\in V$, the products read as follows:
 	
	\begin{equation}\label{Eq:RNProducts-explicit}
		\mathclap{
		\begin{aligned}
		 \mu_n \cs \mu_m &(x_1,\dots,x_{m+k-1}) =
		 \\
		 =&~
		 \sum_{\sigma \in \ush{m,n-1}}
		 \mkern-20mu		 
		  \epsilon(\sigma) 
		 \mu_n\Big(\mu_m(x_{\sigma_1},\dots,x_{\sigma_m}),x_{\sigma_{m+1}}\dots,x_{\sigma_{m+k-1}}	\Big)
		 \\[2em]
		 \mu_n \ca \mu_m &(x_1,\dots,x_{m+k-1}) =
		 \\
		 =&~ 
		 (-)^{|\mu_m|(n-1)}\mkern-30mu
		 \sum_{\sigma \in \ush{m,n-1}}\mkern-20mu
		  (-)^\sigma \epsilon(\sigma) 
		 \mu_n\Big(\mu_m(x_{\sigma_1},\dots,x_{\sigma_m}),x_{\sigma_{m+1}}\dots,x_{\sigma_{m+k-1}}\Big)
		\end{aligned}
		}
	\end{equation}
 where $\epsilon(\sigma) $ is the Koszul sign.
 \\
 These products are not associative and non-associativity is measured by the \emph{associators} multilinear operators\footnote{For the skew-symmetric case, replace $\cs$ with $\ca$.}
 
 \begin{displaymath}
 	 \alpha(\cs ;\mu_\ell,\mu_m,\mu_n) := (\mu_\ell \cs \mu_m) \cs \mu_n - \mu_\ell \cs (\mu_m \cs \mu_n)
 	 ~.
 \end{displaymath}
 Explicitly, they are given by the following equation (see proposition \ref{Prop:SymmetricGerstenhaberAssociators}):
	\note{
   	\begin{displaymath}
  		\begin{aligned}
	  		\alpha&({\cs};\mu_\ell,\mu_m,\mu_n) =
	  			\mu_\ell \circ (\mu_m \otimes \mu_n) \circ B_{\ell,m,n}
	  		\\
  			\alpha&({\ca};\mu_\ell,\mu_m,\mu_n) = 
  				\pm
	  			\mu_\ell \circ (\mu_m \otimes \mu_n) \circ P_{\ell,m,n}
  		\end{aligned}
\end{displaymath} 
} 
	\begin{align}\label{Eq:explicitassociators}
\alpha&({\cs};\mu_\ell,\mu_m,\mu_n)(x_1,\dots,x_{m+n+\ell-2}) = 
  	 \\
  	 =&
  	 \mkern-35mu\sum_{\qquad\sigma \in \ush{m,n,\ell-2}}\mkern-40mu 
  	 \mathscr{s}(\cs;\sigma)~
  	  \mu_\ell\Big(\mu_m(x_{\sigma_1},\dots,x_{\sigma_m}),\mu_n(x_{\sigma_{m+1}},\dots,x_{\sigma_{m+n}}),x_{\sigma_{m+n+1}},\dots,x_{\sigma_{m+n+\ell-2}}\Big) \notag 	
	\end{align}
%
  where $\mathscr{s}(\cs;\sigma)$ and $\mathscr{s}(\ca;\sigma)$ are sign prefactors given explicitly by
  \begin{displaymath}
  	\begin{aligned}
  	\mathscr{s}(\cs;\sigma)&=(-)^{|\mu_n|(|x_1|+\dots+|x_m|)}\epsilon(\sigma)
  	\\
 	 \mathscr{s}(\ca;\sigma)&= 
 	 (-)^{{|\mu_n|(m+\ell)} +{(|\mu_m|(\ell-1))} +m(n+1)}
 	  (-)^{\sigma}\mathscr{s}(\cs;\sigma)
  	\end{aligned}
  \end{displaymath}

	In our conventions, the operator $\Dec$, giving the  d\'ecalage of multilinear maps, is compatible with the \RN products (see theorem \ref{Thm:ManettiFactorizationOnGerstenhaberAlgebras}).
	Namely it induces an isomorphism in the category of graded right pre-Lie algebras:
		\begin{equation}\label{eq:Dec}
			Dec: (M^{skew}(V), \ca) \xrightarrow{\quad \sim \quad} (M^{sym}(V[1]),\cs)
			~.
		\end{equation}

	\begin{remark}[Comparing our conventions with the literature]\label{rem:comparesignmesswithliterature}
		It has to be pointed out how our conventions differ from those found in the literature.
		\begin{itemize}
			\item Notice that our definition of d\'ecalage of multilinear maps $\Dec$ does not get a sign coming from the degree of homogeneous map in input as it appears in the foundational paper \cite[Eq. 3]{LadaStasheff}, or in many other references in our bibliography (e.g \cite[\S 1]{Fiorenza2006}\cite[Rem 1.7]{Bandiera2016}). 
		This discrepancy can be motivated from a different convention in defining the action of the shifted functor on homogeneous map in odd degree or  from a different choice of isomorphism $W[|\mu_n|][n]\cong W[1][|\mu_n|+n-1]$ in the diagram defining $\Dec$:
		\begin{displaymath}
			\begin{tikzcd}
				V^{\otimes n}[n] \ar[r,"\mu_n{[n]}"] & W[|\mu_n|][n] \ar[r,equal] & W[1][|\mu_n|+n-1]
				\\
				(V[1])^{\otimes n} \ar[u,"\dec"] \ar[urr,dashed,"\Dec(\mu_n)"']
			\end{tikzcd}
		\end{displaymath}

			\item Our definition of the skew-symmetric \RN product, see equation \eqref{Eq:RNProducts-myway}, differs by a sign with respect to the corresponding definition given in \cite[Thm. 3.3]{Lecomte1992}. 
			This is a byproduct of a different sign convention in the definition of the Gerstenhaber product (see remark \ref{Rem:signsProblemwithGerstenhaberProducts}) together with our different convention regarding the sign of the d\'ecalage of multilinear maps.
			In appendix \ref{App:RNAlgebras} we introduced this sign in order to guarantee the property of preservation of the \RN products under D\'ecalage, expressed in equation \eqref{eq:Dec} (to be read as a diagram in the category of graded algebras). 
		\end{itemize}
	\end{remark}

\subsubsection{\RN product for coderivations}\label{sec:RNProdCoder}
Let be $V$ a $\ZZ$-graded vector space and consider the graded vector space of homogeneous maps from $V$ into itself $\underline{End}(V):=\underline{\Hom}(V,V)$.
The latter forms an associative algebra with respect to the obvious composition of homogeneous maps
\begin{displaymath}
	\begin{tikzcd}
		V \ar[rr,bend left = 30,"(g\cdot f)"]\ar[r,"f"] 
		&
		V[|f|]\ar[r,"{g[|f|]}"]
		& 
		V[|g|+|f|]
	\end{tikzcd}	
\end{displaymath}
and therefore it forms a graded Lie algebra with respect to the (graded) commutator
\begin{displaymath}
	[f,g]_\circ = f\circ g - (-)^{|f||g|} g \circ f
	~.
\end{displaymath}

Consider then the graded vector subspace of coderivations $\coDer(\overline{S(V)}) \subset \underline{End}(\overline{S(V)},\overline{S(V)})$.
The composition   of two   coderivations $Q_1$ and $Q_2$ is a linear map $Q_1\circ Q_2\colon \overline{S(V)} \to \overline{S(V)}$, which in general fails to be a coderivation.
However the graded commutator 
\begin{displaymath}
	[Q_1,Q_2]:= Q_1\circ Q_2 - (-)^{|Q_1||Q_2|} Q_2 \circ Q_1
\end{displaymath} is a coderivation of degree $|Q_1|+|Q_2|$. 
In other words, the space of coderivations is graded Lie subalgebra of $(\underline{End}(V),[\cdot,\cdot])$ but it is not an associative subalgebra with respect to the composition of graded linear maps.
It is then customary to introduce the \emph{\RN product} (or "composition") on $\coDer(\overline{S(V)})$ defined as
\begin{displaymath}
	\morphism{\cs}
	{\coDer(\overline{S(V)})\otimes\coDer(\overline{S(V)})}
	{\coDer(\overline{S(V)})}
	{Q\otimes Q'}
	{\widetilde{L}_{\sym}\left(\pr_V \circ Q \circ Q' \right)}
	~,
\end{displaymath}
where $\pr_V: \overline{S(V)}\twoheadrightarrow V$ denotes the standard projection on $V^{\odot 1}$ and $\widetilde{L}_{\sym}$ is the "lift to a coderivation" operator introduced in lemma \ref{lem:liftToCoder}.
The composition $\cs$, which is not associative, makes the space of coderivations into a graded right pre-Lie algebra, in particular, it induces the same Lie brackets inherited from $\underline{\End}(V)$ since
\begin{displaymath}
	[Q,Q']_{\cs}= 
	\widetilde{L}_{\sym}\left(\pr_V \left( Q \circ Q'-(-)^{|Q||Q'|}~Q' \circ Q \right)\right)
	= \cancel{\widetilde{L}_{\sym}\circ\pr_V}\circ [Q,Q']_{\circ}
	~.
\end{displaymath}
The naming comes from the following sequence of isomorphisms in the category of graded pre-Lie algebra structures:
\begin{equation}\label{eq:NocciolodiAppendice}
	\begin{tikzcd}
		(M^{\skew}(V),\ca)	\ar[r,"\Dec"]&
		(M^{\sym}(V[1]),\cs) 
		\ar[r,"\widetilde{L}_{\sym}"]&
		\coDer(S(V[1]),\ca)		 
	\end{tikzcd}
	~,
\end{equation}
implied by the natural identification 
\begin{displaymath}
	M^{\sym}(V) = \bigoplus_{n\geq 1} \underline{\Hom}(V^{\odot n},V) =
		\underline{\Hom}(\overline{S(V)},V)
\end{displaymath}
given by the commutation rule of categorical (co)limits with the $\Hom$ functor.
\begin{remark}
	Notice that here, and in appendix \ref{App:RNAlgebras}, we decided to present this topic defining the products on $M^{\sym}(V)$ and $\coDer(\overline{S(V)})$ independently and proving they are isomorphic as a second step.
	Often in the literature, see for example \cite{Bandiera2016,Manetti-website-coalgebras,Miti2020}, one chooses to go the opposite way.
	Namely, one starts from the definition of $\cs$ on $\coDer(\overline{S(V)})$ and then pullback the product to $M^{\sym}(V)$, along the lift, and to $M^{\skew}(V[-1])$ along the d\'ecalage.
		For instance, the  \emph{Nijenhuis-Richardson product} of two given maps $a,b\in \underline{\Hom}(\overline{S(V)},V)$ can be explicitly given by
\begin{equation}\label{eq:compsymm}
 a\cs b:= \pr_V(C_a\circ C_b) ~,
\end{equation}
where $C_a$ denotes the lift of $a$ to a coderivation on $\overline{S(V)}$.
It is not too difficult to see that this is explicitly obtained by summing insertions of $b_j$ in $a_i$, where $a_i:=a \vert_{V^{\odot i}}$ (see equation \eqref{Eq:RNProducts-explicit}).
\end{remark}
Let us stress again that $\widetilde{L}_{\sym}(a\cs b) \neq \widetilde{L}_{\sym}(a)\circ \widetilde{L}_{\sym}(b)$ (the latter is not even a coderivation in general). However, as shown above, the lift to a coalgebra coderivations preserve the commutator bracket, \ie $$[\widetilde{L}_{\sym}(a),\widetilde{L}_{\sym}(b)] = \widetilde{L}_{\sym}([a,b]_{\cs})~.$$

In chapter \ref{Chap:MarcoPaper}, we will make use of the following construction for producing a coalgebra isomorphism starting from a degree $0$ coderivation:  
\begin{lemma}[Exponential of a coderivation]\label{lem:coderExpo}
	Given a $m$-nilpotent coderivation $Q\in \coDer^0(\overline{S(V)})$, \ie $Q^n=0$ for any $n\geq m$, then the corresponding exponential operator
	\begin{displaymath}
		e^Q := \sum_{n\geq 0} \dfrac{Q^k}{k!}~ \in \End(\overline{S(V)})
	\end{displaymath}
	is an endomorphism of coalgebras.
\end{lemma}
\begin{remark}\label{rem:codermor}
	We notice that if $C$ is a degree $0$ coderivation of $\overline{S(V)}$, then $e^C$ is a morphism of coalgebras, provided it converges.
	The nilpotency guarantees that the summation involved is actually made up of a finite number of summands.

	A slightly more general statement can be provided as follows.
	Endow $\overline{S(V)}$ with the filtration $\mathcal{F}_0\subset\mathcal{F}_1\subset \mathcal{F}_2\subset\dots$ where $\mathcal{F}_0=\{0\}$ and $\mathcal{F}_k:=\oplus_{i=1}^k V^{\odot i}$ for any $k\geq 1$.
	\\
	Whenever $C$ maps $\mathcal{F}_k$ to $\mathcal{F}_{k-1}$ for all $k\geq 1$,
	it follows that $(e^C-\Id)$ has the same property. 
	Therefore $e^C\eval_{V^{\odot n}}$ is a finite sum for all 
$n$, and $e^C$ converges.
	\\
	Let now $C$ be the coderivation obtained lifting the graded morphism $p\in \underline{\Hom}^0(\overline{S(V)},V)$,	$C$ satisfies the above property whenever we have $p\eval_{S^1V}=0$. This follows from eq. \eqref{eq:coder}.
\end{remark}

Observe that out of given a multilinear map $f:\overline{S(V)}\to V$ one can construct two corresponding coalgebra morphisms: 
$L_{\sym}(f)\in \Hom_{coAlg}(V,V)$, given by the lift to a coalgebra morphism (Lemma \ref{lem:liftToMorph}), 
and $\exp(\widetilde{L}_{\sym}(f))\in \Iso_{coAlg}(V,V)$ given by lemma \ref{lem:coderExpo}.
Clearly the two do not coincide, for instance their first components results:
	\begin{displaymath}
		L_{\sym}(f)\eval_{V} = f_1 \quad,\qquad
		\exp(\widetilde{L}_{\sym}(f))\eval_{V}= \exp(f_1)
		~.
	\end{displaymath}

\section{Lie infinity structures}\label{Sec:LinfinityAlgebras}
Lie infinity structures, from now on $L_\infty$, are generalizations of \emph{differential graded Lie algebras} (DGLA), therefore of cochain complexes and (graded) Lie algebras at the same time.
The key ideas are two (we refer to \cite[\S 2]{Ryvkin2016a} for a concise introductory exposition):
\begin{enumerate}
	\item start from a DGLA $L$ and weaken the Jacobi equation condition requiring it to be satisfied only up to a chain homotopy;
	\item require that the failure of the ordinary Jacobi equation identity is controlled by a skew-symmetric $3$-multilinear map from $L$ to itself satisfying a similar higher "weaker" Jacobi equation thus allowing for a possible infinite sequence of higher multibrackets with arity greater than $3$.
\end{enumerate}

Given a graded vector space $V$, there are several alternative ways of presenting what does it mean to endow it with $L_\infty$-structure. 
Historically, the first precise definition it is due to Lada and Stasheff:
	\begin{definition}[$L_\infty$-algebra \emph{(Lada, Stasheff) \cite{LadaStasheff}}]
	\label{Def:LInfinityStasheff}
		We call \emph{$L_\infty$-algebra} the pair 
		\begin{displaymath}
			\Big( L, \lbrace \mu_k \rbrace_{k\in \mathbb{N}} \Big)
		\end{displaymath}
		given by a $\mathbb{Z}$-graded vector space $L$ together with
		a family, parametrized by integers $k\geq 1$, of homogeneous graded skew-$k$-multilinear maps 
		\begin{displaymath}
			\mu_k : \wedge^k L \rightarrow L[2-k]
		\end{displaymath}
		(usually called \emph{multi-brackets}) 	satisfying \emph{"Higher Jacobi"} relations
		\begin{equation}\label{Eq:HigherJacobiStasheff}
			0 = 
			\mkern-30mu\sum_{\substack{i+j=m+1\\ \sigma \in \ush{i,m-i}}}\mkern-30mu
			(-)^{i(j+1)} (-)^\sigma \epsilon (\sigma; x) 
			~\mu_j \Big( \mu_i(x_{\sigma_1},\dots, x_{\sigma_i}), x_{\sigma_{i+1}},\dots, x_{\sigma_m}\Big)
		\end{equation}
		$\forall m\geq 1$ and $x_i$ homogeneous elements in $L$.	Recall that $\ush{i,m-i}$ denotes the subgroup of unshuffle permutations of $m$ elements, namely $\sigma \in \ush{i,m-i}$ if $\sigma(j)<\sigma(j+1)$ for every $j\neq i$ (see appendix \ref{App:UnshuffleAtors}).
	\end{definition}

	\begin{remark}[Homological and Cohomological convention]
		It is been noted already by Lada and Stasheff that are two possible convention in the previous definition. 
		Namely that the unary operator $\mu_1$ lowers degrees (homological convention) or it raises degrees (cohomological convention).
		As explained in section \ref{sec:HomologicalAlgebrasConventions}, in the thesis we are adopting the cohomological convention.
		\\
		Passing from a notation to the other can be easily achieved by the reverse-grading functor $\setminus$. 
		Therefore, in the homological case, the degree of the $k$-multibracket $\mu_k$ changes from $(2-k)$ to $(k-2)$.
	\end{remark}

	\begin{example}[Cochain complexes]\label{ex:cochaincomplexLinfinity}
		A $L_\infty$-algebra with multibrackets $\mu_k=0$ for any $k\geq 2$ is a cochain complex.
		Indeed, it is given by a pair $(L,\mu_1)$ together with the only non-trivial higher Jacobi equation (\eqref{Eq:HigherJacobiStasheff} with $m=1$) reading as 
		\begin{displaymath}
			\mu_1(\mu_1(x))=0 \qquad \forall x \in L
		~.
		\end{displaymath}
		Essentially, any $L_\infty$-algebra $(L,{\mu_k}_{k\geq 1})$ has an underlying cochain complex, or differential graded vector spaces, obtained by neglecting all multibrackets of arity greater than $1$.
	\end{example}

	\begin{example}[Differential graded Lie algebras]\label{Ex:dglaAsLinfinity}
		A $L_\infty$-algebra $(L,{\mu_k}_{k\geq 1})$ with multibrackets $\mu_k=0$ for any $k\geq 3$ is a differential graded Lie algebra (DGLA, see definition \ref{def:DGLA}).
		Namely, denoting by $\d$ and $[\cdot,\cdot]$ the unary and binary operators of $L$, equation \eqref{Eq:HigherJacobiStasheff} with $m=1$ yields the $2$-nilpotency condition of $\d$, and the other two non-trivial higher Jacobi equations, given by $m=2,3$, reads as follows:
		\begin{displaymath}
			\begin{aligned}
				\d [x_1,x_2] &= [\d x_1, x_1] -(-)^{|x_1|} [x_1, \d x_2] ~~,
				\\
				0 &=[[x_1,x_2],x_3] - (-)^{|x_3||x_2|}[[x_1,x_3],x_2]+(-)^{|x_1|(|x_2|+|x_3|)}[[x_2,x_3],x_1] ~~.
			\end{aligned}		
		\end{displaymath}
	The first one expresses that $\d$ is graded derivation with respect to the algebraic product $[\cdot,\cdot]$ and the second one is the graded Jacobi identity for $[\cdot,\cdot]$.
	When $L$ is concentrated in degree $0$, the only non-trivial multibrackets is given by $[\cdot,\cdot]$, hence any $L_\infty$-algebra concentrated in degree $0$ is simply a \emph{Lie algebra}.
	\end{example}

	We will make us of the following nomenclature:
\begin{definition}\label{Def:groundedLinfinity}
	\begin{itemize}
		\item A $L_\infty$-algebra is called an \emph{Abelian $L_\infty$-algebra} if all $k$-ary brackets with $k\geq 2$ are trivial, \ie is a plain chain complex.\cite[\S 1.0.3.]{Fiorenza2014a}.
		\item A $L_\infty$-algebra is called a \emph{grounded $L_\infty$-algebra} if $k$-ary brackets, for any $k\geq 2$, are trivial when evaluated on elements in non-zero degree, \ie
		$\mu_k(x_1,\dots,x_k)=0$ for any $\sum_{i=1}^k|x_i| \neq 0$.
		\cite[Def 2.33]{Ryvkin2016a}. (It is called \emph{property (P)} in \cite{Callies2016}.)
		\item A $L_\infty$-algebra is called a \emph{$L_n$-algebra} (or \emph{Lie-$n$}-algebra) if the underlying vector space is concentrated in degrees from $-n$ to $0$. 
		(In the homological convention it would be concentrated in degrees $(0,\dots,n)$.)
		The corresponding $L_\infty$-structure consists of $n$+$1$ multibrackets $\{\mu_1,\dots,\mu_{n+1}\}$.
	\end{itemize}
\end{definition}
In this thesis, we will mainly concerned with $L_\infty$-algebras of the last two kinds.
\begin{remark}[Curved $L_\infty$-algebras]
	It is often found in literature the notion of \emph{curved $L_\infty$-algebra}.
	This is obtained from definition \ref{Def:LInfinityStasheff} by additionally allowing for an element $\mu_0$, "$0$-ary bracket", in degree $2$ (or $-2$ in the homological convention) and allowing indexes $i,j$ and $m$ in equation \eqref{Eq:HigherJacobiStasheff} to be zero. Specifically, when $m=0$ this would mean that $\mu_1(\mu_0)=0$ hence $\mu_0\in Z^2(L,\mu_1)$ is a cochain in the chain complex $(L,\mu_1)$
\end{remark}
Taking advantage of the \RN formalism introduced in section \ref{sec:RNstuff}, it is possible to encode $L_\infty$-structures in a particularly succinct way.
\begin{remark}[Reading definition \ref{Def:LInfinityStasheff} in NR-algebraic terms]
	By its very definition, a $k$-multibracket $\mu_k$ of a $L_\infty$-algebra $(L,\{\mu_k\}_{k\geq 1})$ is an element of $M^{\skew}_{k,2-k}(L)$.
According to the definition of the \RN product (see section \ref{Section:MultibracketsAlgebra}), 
		"higher Jacobi equations" \eqref{Eq:HigherJacobiStasheff} can be synthetically recast as the vanishing of the multilinear operators $J_m \in M^{\skew}_{m,3-m}(L)$ for any $m\geq 1$,
		explicitly given by
		\begin{equation}\label{Eq:JacobiatorVit}
			\begin{split}
			J_m :&= \sum_{k=1}^m (-)^{k(m-k)} \mu_{m-k+1} \circ (\mu_k\otimes \Unit_{m-k})\circ P_{k,m-k} =
			\\
			&=
			\sum_{k=1}^m \mu_{m-k+1}\skewgerst\mu_k
			~.
			\end{split}
		\end{equation}

		With the grading defined on the skew-symmetric \RN algebra,
		see equation \eqref{eq:RNspaces},
		every multilinear maps $\mu_k$ of a $L_\infty$-structure is a degree $1$ element in $M^{\skew}(V)$ and, in particular, their sum is homogeneous.
		Denoting by $\mu = \sum_{k\geq 1} \mu_k \in (M^{\skew}(L))^1$ the direct sum of all multibrackets, which completely encodes the $L_\infty$-structure given by $\{\mu_k\}_{k\geq 1}$,	
		 follows that
		\begin{displaymath}
			\sum_{m\geq 1} J_m = \mu \skewgerst \mu = \frac{1}{2} [ \mu,\mu]_{\skewgerst}	
			~,
		\end{displaymath}
		in other words, the $L_\infty$-structure $\mu$ is a \emph{Maurer-Cartan element} in the graded Lie algebra $(M^{\skew}(L),[\cdot,\cdot]_{\skewgerst})$ (see definition \ref{def:MCelements}).
\end{remark}

	Summing up, we get the following reformulation of definition \ref{Def:LInfinityStasheff}:
		\begin{definition}[${L_\infty}$-algebras]\label{Def:LInfinityTony}
				We call \emph{$L_\infty$-algebra} a pair $(V, \mu)$ where $V$ is a graded vector space, and $\mu\in (M^{\skew}(V))^1=\bigoplus_{k\geq 1} M^{\skew}_{k,2-k}(V)$ is the direct sum of skew-symmetric multibrackets satisfying the \emph{higher Jacobi equations} \ie
				\begin{equation}\label{Eq:HigherJacobiTony}
					J_n :=~ \mu \skewgerst \mu \eval_{L^{\otimes n}}\equiv \sum_{k=1}^{n-1} \mu_k \skewgerst\mu_{n-k} =~0
					\qquad \forall n\geq 2 ~,
				\end{equation}
	where $\mu_k$ denotes the projection of $\mu$ into $M^{\skew}_{k,2-k}(V) \subset (M^{\skew}(V))^1$.
		\end{definition}

	\begin{remark}[Understanding higher Jacobi equations as homotopies]
		The higher Jacobi equations (equations \eqref{Eq:HigherJacobiStasheff} or \eqref{Eq:HigherJacobiTony}) can be made slightly more expressive introducing the so-called \emph{Jacobiator}\footnote{Note that we are slightly departing from the more common notation, see for example \cite{Vitagliano2013}, which reserves the name "Jacobiator" to the operator $J_m$ (equation \eqref{Eq:JacobiatorVit}) instead of $j_m$ (equation \eqref{Eq:Jacobiator}.} 
		multilinear operator $j_m\in M^{\skew}_{m,3-m}(L)$ defined as follows
		\begin{equation}\label{Eq:Jacobiator}
			j_m 
			= \sum_{k=2}^{m-1} \mu_{m-k+1} \skewgerst \mu_k
			= J_m - \left(  \mu_m \skewgerst \mu_1 +\mu_1 \skewgerst \mu_m \right)
			~.
		\end{equation}
		Observe that, when $m=3$ one gets
		\begin{align*}
			j_3 (x_1,x_2,x_3)=&~+ 
			\mu_2(\mu_2(x_1,x_2),x_3) +
			\\
			&~- (-)^{|x_3||x_2|}\,\mu_2(\mu_2(x_1,x_2),x_3)+
			\\
			&~+(-)^{|x_1|(|x_2|+|x_3|)}\,\mu_2(\mu_2(x_2,x_3),x_1) 
			=
			\\
			=&
			(-)^{|x_1||x_3|}\left(
			 (-)^{|x_1||x_3|}\mu_2(\mu_2(x_1,x_2),x_3) + \cyc
			\right)
			~~,
		\end{align*}	
		where $\cyc$ denotes sum over cyclic permutations.
		Hence equation \eqref{Eq:Jacobiator} recovers the usual definition of Jacobiator for (graded) Lie algebras.

		In example \ref{ex:cochaincomplexLinfinity} has been shown that the unary bracket of a $L_\infty$-algebra determines a cochain complex boundary operator $\mu_1$ on $L$, let us denote $\mu_1=\partial$.
		Employing the notations introduced in section \ref{sec:HomologicalAlgebrasConventions}, we introduce a coboundary operator on $V^{\otimes m}$ given by
		\begin{displaymath}
				\partial_{{\otimes m}} := (-)^{m-1} \d_{\otimes m} =
				(-)^{m-1}\sum_{k=1}^m \Unit_{k-1}\otimes ~\partial~ \otimes \Unit_{m-k}
			~.		
		\end{displaymath}	
		(differing from the total coboundary operator given in equation \eqref{eq:totalcoboundaryoperator} by an overall sign).
		Observe then that, for any given graded skew-symmetric $n$-multilinear map $f$, the following equality holds
		\begin{displaymath}
			\begin{aligned}
				f \circ 	\partial_{{\otimes n}} &(x_1,\dots,x_n) =
				\\
				=& (-)^{n-1}\left(
				f(\partial x_1,\dots,x_n) + \dots + (-)^{|x_1|+\dots+|x_{n-1}|} f(x_1,\dots,\partial x_n) \right)=
				\\
				=&
				(-)^{n-1}\sum_{i=1}^{n}(-)^{|x_i|(|x_1|+\dots+|x_{i-1}|)} f(\partial x_i,x_1,\dots \widehat{x_i},\dots, x_n)				
				=
				\\
				=&(-)^{n-1} f \circ \left(\partial \otimes \Unit_{n-1}\right) \circ P_{1,n-1} ~ (x_1,\dots, x_n) 
				=
				\\
				=&
				f \ca \mu_1 (x_1,\dots, x_n)
				~.
			\end{aligned}	
		\end{displaymath}
		Therefore, one can conclude that the condition $J_k=0$, \ie when higher Jacobi equations hold, is equivalent to say that
		\begin{displaymath}
			\mu_m \circ \partial_{{\otimes m}} + \partial \circ \mu_m	= -j_m
			~, 
		\end{displaymath}
		hence the $m$-ary multibracket $\mu_m$ is a chain-homotopy between the higher Jacobiator $j_m$ and $0$: 
\begin{displaymath}
	\begin{tikzcd}[column sep = huge]
		 (L^{\otimes m}, \partial_{\otimes m}) \ar[r,blue,"j_m", bend left=20, ""{name=U, below}]
		\ar[r,"0"', bend right=20, ""{name=D}]
		&[1em] (L,\partial)[3-m]
		\ar[Rightarrow, purple,"\mu_m", from=U, to=D]
\end{tikzcd}	
~.
\end{displaymath}
		This justify the slogan: "\emph{$L_\infty$-algebra is the notion that one obtains from a Lie algebra when one requires the Jacobi identity to be satisfied only up to a higher coherent chain homotopy}"
		\footnote{
			Notice that the term "coherent homotopy" has a precise meaning  in the context of homotopy theory. 
			What we are doing here is rather providing a basic justification of the reason why the term "homotopy" appears in conjunction with these structures.
			Notably, the first name attributed to this algebraic structure has been "strongly homotopy Lie algebra"\cite{LadaStasheff}.
		}
		\cite{nlab:l-infinity-algebra}\cite{Rogers2010}\cite{Shahbazi2016}.
		Observe at last that the first two higher Jacobi identities implies that $j_1$ and $j_2$ are automatically zero. 
		Summing up, the first Jacobi equations $J_m=0$ reads as follows:
		\begin{itemize}
			\item when $m=1$ means that $\partial \circ \partial = 0$, \ie $\partial=\mu_1$ is a coboundary operator and $(L,\partial)$ is a cochain complex;
			\item when $m=2$ means that $\mu_2 \circ \partial_{\otimes 2} + \partial \circ \mu_2=0$, hence $\mu_2$ is a chain map from $(L^{\otimes 2},\partial_{\otimes_2})$ to $(L,\partial)$;
			\item when $m=3$ means that $\mu_3$ is a chain homotopy from the Jacobiator $j_2$ to the zero map,
			\begin{displaymath}
	\begin{tikzcd}[column sep = huge]
		 (L^{\otimes 3}, \partial_{\otimes 3}) \ar[r,blue,"j_3", bend left=20, ""{name=U, below}]
		\ar[r,"0"', bend right=20, ""{name=D}]
		&[1em] (L,\partial)
		\ar[Rightarrow, purple,"\mu_3", from=U, to=D]
	\end{tikzcd}	
~.
\end{displaymath} 
			In particular $\mu_3$ is a degree $-1$ homogeneous map between the cochain complexes $(L^{\otimes 3},\partial_{\otimes 3})$ and $ (L,\partial)$.
		\end{itemize}
	\end{remark}

	\subsection{Coalgebraic approach to $L_\infty$-structures}
		In a nutshell, the gist of the previous discussion is subsumed by stating that the set of all possible $L_\infty$-structure on a given graded vector space $L$ is given by Maurer-Cartan elements:
	\begin{displaymath}
		\mathbb{L}_{\infty}(L) := \text{MC}(M^{\skew}(L),[\cdot,\cdot]_{\ca}) 
		\equiv
		\lbrace
			\mu \in M^{\skew} ~ \left\vert
			\quad
				||\mu||=1
			~,\quad
				\mu\ca\mu = 0
		\right.\rbrace
	~.
	\end{displaymath}
	(Compare with definition \ref{def:MCelements}).	
	This claim, joined with diagram \eqref{eq:NocciolodiAppendice} 
	(which basically subsumes the contents of appendix \ref{App:RNAlgebras}), implies that one gets two other completely equivalent presentations of a $L_\infty$-algebra structure over the graded vector space $L$
	\begin{equation}\label{eq:VisioneGlobaleLinfinito}
		\begin{aligned}
		\mathbb{L}_{\infty}(L) 
		\cong&~
		\lbrace
			\nu \in M^{\sym}(L[1]) ~ \left\vert
			\quad
				|\nu|=1
			~,\quad
				\nu\cs\nu = 0
		\right.\rbrace
		\cong 
		\\
		\cong&~
		\lbrace
			Q\in \coDer(S(L))	) ~ \left\vert
			\quad
				|Q|=1
			~,\quad
				Q\circ Q = 0
				\right.			
		\rbrace
	~.
		\end{aligned}
	\end{equation}	
	\begin{remark}
		Observe that in the last term of equation \eqref{eq:VisioneGlobaleLinfinito} appears the composition of two coderivation, which in principle is not a coderivation, instead that the \RN product introduced in section \ref{sec:RNProdCoder}.
	This is due by the commutativity of the following diagram in the category of graded vector space
	\begin{displaymath}
		\begin{tikzcd}[column sep = huge]
			\overline{\overline{S(V)}} \ar[r,"\widetilde{L}_{\sym}(\mu)"] \ar[dr,"\mu"']
			\ar[rr,bend left=20,"\widetilde{L}_{\sym}\big(\mu\cdot \widetilde{L}_{\sym}(\mu)\big) = \widetilde{L}_{\sym}(\mu\symgerst\mu)=0"]
			&[2em]
			\overline{S(V)}[1] \ar[r,"{\widetilde{L}_{\sym}(\mu)[1]}"]\ar[dr,"{\mu[1]}"']\ar[d]
			&[2em]
			\overline{S(V)}[2]\ar[d]
			\\
			&
			V[1]
			&
			V[2]
		\end{tikzcd}
	\end{displaymath}
	noting that the commutativity of the uppermost triangle in the above diagram makes sense only if the curved arrow is the zero arrow.
	\end{remark}	
	
	In this spirit, one could introduce the following two definition:
	\begin{definition}[Shifted Lie infinity  algebra (${L_\infty[1]}$-structures)
	(see \emph{\cite[Def. 5]{Kajiura2006b} or \cite[Def. 4]{Vitagliano2013}})]\label{Def:LInfinityShifted}
		We call \emph{shifted $L_\infty$-algebra} a pair $(W, \{\ell_k\}_{k\geq 1})$ where $W$ is a graded vector space, and $\ell_k\in M^{\sym}_{k,1}$ are symmetric multibrackets satisfying the \emph{higher Jacobi equations} \ie 
		\begin{displaymath}
			J_n :=~ \sum_{k=1}^{n-1} \ell_k \symgerst\ell_{n-k} =~0
			\qquad \forall n\geq 2 ~.
		\end{displaymath}
	\end{definition}
	\begin{definition}[Chevalley-Eilenberg complex for a ${L_\infty[1]}$-structure]\label{def:CELinfty1struct}
		We call \emph{Chevalley-Eilenberg complex} pertaining to the ${L_\infty[1]}$-algebra $(W, \nu)$, where $\nu=\sum_{k\geq 1} \nu_k$, the cofree graded codifferential coalgebra
		\begin{displaymath}
			\CE(W, \{\ell_k\}_{k\geq 1}) := ( \overline{S(W)}, \widetilde{L}_{\sym}(\nu) )
		\end{displaymath}
		given by the (cofree cocommutative) symmetric tensor coalgebra $\overline{S(W)}$ together with the lift of $\nu$ to a coderivation.
	\end{definition}
	\begin{remark}
		We point out that our naming here diverges from the literature. Often the name of "Chevalley-Eilenberg complex of a $L_\infty$-algebra $L$ is reserved to the dual of definition \ref{def:CELinfty1struct}, hence by the graded vector space $\overline{S(W)}^\ast \cong \Hom(\overline{S(W)},\RR)$ endowed with a certain differential.
		A more general notion can be given in terms of $L_\infty$-algebras representations, see \eg \cite{Reinhold2019}.
	\end{remark}
	The key point, as stated in \cite[Thm. 2.3]{LadaMarkl}, is the existence of a one-to-one correspondence between $L_\infty$-algebra structures on $V$ and degree $1$ nilpotent $2$-coderivations on $\overline{S(V[1])}$.
	When $(W, \nu)$ is a differential graded Lie algebra, $\CE(W,\nu)$ is sometimes called \emph{Quillen construction} (see \cite[\S 2]{Fiorenza2006}).
	\begin{reminder}[(Ordinary) Chevalley-Eilenberg (co)-chain complex]\label{Rem:CEconventions}
		Given an (ungraded) Lie algebra $\mathfrak{g}=(\mathfrak{g},[\cdot,\cdot])$ it is customary to call its \emph{Chevalley-Eilenberg chain complex} the pair $(\wedge^\bullet \mathfrak{g},\partial)$, depicted by
	\begin{displaymath}
		\begin{tikzcd}
			\cdots \ar[r] &
			\wedge^n\mathfrak{g} \ar[r, "\partial_{\CE}"] &
			\wedge^{n-1}\mathfrak{g} \ar[r] & \cdots
		\end{tikzcd}
		~,
	\end{displaymath}	
		 where the first data, in slight contradiction with our previous notation, is the graded vector space
		$\wedge^\bullet \mathfrak{g} = (k \mapsto \Lambda^k(\mathfrak{g}))$
	and the \emph{ CE boundary operator} $\partial: \wedge^\bullet \mathfrak{g} \to \wedge^{\bullet-1} \mathfrak{g}$		
	is given explicitly, for any homogeneous element $\xi_i \in \mathfrak{g}$, by the following equation:
\begin{equation} \label{eq:CE_boun}
	\partial_{\CE} (\xi_1 \wedge \xi_2 \wedge \dots \wedge \xi_k) := 
	\mkern-15mu\sum_{1\leq i< j \leq n}\mkern-15mu
	 (-1)^{i+j}\, [\xi_i, \xi_j] \wedge \xi_1 \wedge \dots
\wedge {\hat \xi}_i \wedge \dots \wedge {\hat \xi}_j \wedge \dots \wedge \xi_n,
\end{equation}
with $\partial_0 = 0$ and $\widehat{\cdot}$ denoting deletion as usual.
	\\	
	Dually, one can define a cochain complex structure on $\wedge^{\bullet}\mathfrak{g}^*$ introducing the \emph{Chevalley-Eilenberg coboundary (differential)} $\delta_{\CE}: \wedge^k\mathfrak{g}^* \to \wedge^{k+1} \mathfrak{g}^*$,  whose action on an element $\phi \in \wedge^\bullet \mathfrak g^*$ is given by $\delta_{CE} \phi := \phi \circ \partial$.
	\\
	Henceforth, by \emph{$k$-th CE homology group} and \emph{$k$-th CE cohomology group} we will intend the following two vector spaces:
	\begin{displaymath}
		H_k(\mathfrak{g}):= H_k(\wedge^\bullet \mathfrak{g},\partial_{\CE})
		~, \qquad
		H^k(\mathfrak{g}^*)= H^k(\wedge^\bullet \mathfrak{g}^\ast,\delta_{\CE})
	~.
	\end{displaymath}
	More in general, one can define the \emph{Chevalley-Eilenberg cohomology over a representation $\rho$}, \ie a Lie algebra morphism $\rho: \mathfrak{g}\to \End(V)$ for a certain vector space $V$, as the cochain complex
	\begin{displaymath}
		\begin{tikzcd}
			\cdots \ar[r] &
			\Hom_{\Vect}(\wedge^n\mathfrak{g},V) \ar[r, "\delta_{\CE}"] &
						\Hom_{\Vect}(\wedge^{n+1}\mathfrak{g},V) \ar[r] & \cdots
		\end{tikzcd}
	\end{displaymath}
	with coboundary operator defined, for any $\omega: \mathfrak{g}^{\wedge n}\to V$ and for any elements $x_i \in \mathfrak{g}$, by
		\begin{displaymath}
			\begin{aligned}
			(\delta \omega) (x_1,\dots,x_{n+1})
			=&
			+\sum_{i=1}^{n+1}(-)^{i+1}\rho(x_i)\cdot \omega(x_1,\dots,\hat{x_i},\dots,x_{n+1})+
			\\
			&
			+\mkern-25mu\sum_{1\leq j < k \leq n+1}\mkern-25mu
			(-)^{j+k} \omega([x_j,x_k],x_1,\dots,\hat{x_j},\dots,\hat{x_k},\dots,x_{n+1})
			~.
			\end{aligned} 
		\end{displaymath}
	
	\end{reminder}
	\begin{remark}
	Observe that the naming in definition \ref{def:CELinfty1struct} is completely compatible with the previous construction. The Lie algebra $\mathfrak{g}$ is in particular a $L_\infty$-algebra structure with only non-trivial multibrackets given by $\mu_2=[\cdot,\cdot]$. Applying the d\'ecalage and the lift one gets the following diagram in the graded vector space category
	\begin{displaymath}
		\begin{tikzcd}
			&
			\overline{S(\mathfrak{g}[1])} \ar[r,"\widetilde{L}_{\sym}(\Dec(\mu_2))"] \ar[d,two heads]
			&[7em]
			(\overline{S(\mathfrak{g}[1])})[1] \ar[dd,two heads]
			\\
			(\mathfrak{g}^{\wedge 2})[2] \ar[r,"\dec"]
			\ar[drr,bend right=10,"{\mu_2[2]}"']
			&
			(\mathfrak{g}[1]^{\odot 2}) \ar[dr,"\Dec(\mu_2)"']
			&
			\\
			&
			&
			\mathfrak{g}[1][1]
		\end{tikzcd}
	\end{displaymath}	
	where $\CE(\mathfrak{g})=(S(\mathfrak{g}[1]),\widetilde{L}_{\sym}(\Dec(\mu_2))$ is precisely the Chevalley-Eilenberg complex as defined in \ref{def:CELinfty1struct}.
	Explicitly, one can see that $\CE(\mathfrak{g}) = \setminus (\Lambda^\bullet \mathfrak{g},-\partial)$, where $\setminus$ denotes the degree-reversing functor introduced in equation \eqref{eq:degreeReversingFunctor}. 
	In fact, according to diagram \eqref{eq:decalageRestriction}, one has
	\begin{displaymath}
	\CE(\mathfrak{g})^{-k} = (\mathfrak{g}[1])^{\odot k} \cong \mathfrak{g}^{\wedge k}[k]
	~,
	\end{displaymath}
	 and, from lemma \ref{lem:liftToCoder}, one has that
	\begin{displaymath}
		\mathclap{
		\morphism{\widetilde{L}_{\sym}(\mu_2)}
		{\CE(\mathfrak{g})^{-k}}
		{\CE(\mathfrak{g})^{-k+1}}
		{x_1\wedge\dots \wedge x_n}
		{\displaystyle \mkern-30mu\sum_{\sigma\in \ush{2,n-2}}\mkern-30mu(-)^\sigma
			[x_i, x_j] \wedge x_1 \wedge \dots
\wedge {\hat x}_i \wedge \dots \wedge {\hat x}_j \wedge \dots \wedge x_n}
~.}
	\end{displaymath}
		The latter can be read as the fact that the lift $\mu_2$ is equivalent to $\partial_{\CE}$, defined in reminder \ref{Rem:CEconventions}, modulo a sign. 
		Namely, $\widetilde{L}_{\sym}(\mu_2) \equiv -\partial_{\CE}$.
	\end{remark}		

	Spelling out, an $L_\infty$-algebra structure on $L$ is equivalently given by
	\begin{itemize}
		\item	a family $\lbrace \mu_k \rbrace_{k\geq 1}$ of graded skew-symmetric multibrackets $\mu_k \in M^{\skew}_{k,2-k}(L)$ on $L$ satisfying the higher Jacobi equation $J_m =\sum_{k=1}^m \mu_k\ca\mu_{m-k+1}=0$;
		(that is a $2$-nilpotent degree $1$ element in the \RN algebra $(M^{\skew}(V),\ca)$);
		\item a family $\lbrace \nu_k \rbrace_{k\geq 1}$ of graded symmetric multibrackets $\nu_k \in M^{\sym}_{k,1}(L[1])$ satisfying the higher Jacobi equation $J_m =\sum_{k=1}^m \nu_k\cs\nu_{m-k+1}=0$;
		(that is a $2$-nilpotent degree $1$ element in the \RN algebra $(M^{\sym}(V),\cs)$);
		\item a $2$-nilpotent, degree one coderivation $Q \in \coDer(S(L[1]))$, \ie a codifferential on the (cofree) graded coalgebra $S(V[1])$. 
		(In the terms of \cite{nlab:l-infinity-algebra}, an $L_{\infty}$-algebra is \emph{cofree cocommutative differential coalgebra}, \ie a dg-coalgebra whose underlying coalgebra is isomorphic to the symmetric tensor coalgebra for a given graded vector space.)
\end{itemize}		
	Although the characterization in terms of $2$-nilpotent coderivations may appear particularly convoluted, it has several advantages.
	Specifically, it transparently provides the notion of commutator and linear combination of $L_\infty$-structures (beware that such operations do not do not yield $L_\infty$-structures in general), the notion of (tangent) differential graded Lie algebra $(\overline{S(V)},[\nu,\cdot]_{\cs},[\cdot,\cdot]_{\cs})$ governing the deformations of a given $L_\infty[1]$-algebra $(V,\nu)$, and the notion of $L_\infty$-morphism, as we will see in the following subsection.

\subsection{$L_\infty$-morphisms}\label{subsec:LinftyMorphi}
	In the coalgebraic framework there is an obvious notion of morphism between $L_\infty$-algebras:
	\begin{definition}[$L_\infty{[1]}$-morphisms (coalgebraic framework)]\label{Def:LinftyMorphism-coalg}
		Given two $L_\infty[1]$-algebras $(V,\mu)$ and $(W,\nu)$, 
		an \emph{$L_\infty[1]$-morphism} between $(V,\mu)$ and $(W,\nu)$ is a morphism of graded differential cocommutative coalgebra $F: ({S([1])},\widetilde{L}_{\sym}(\mu))\to ({S(W[1])},\widetilde{L}_{\sym}(\nu))$, \ie a graded coalgebra morphism such that the following diagram commutes in the graded vector space category
		\begin{displaymath}
			\begin{tikzcd}
				{S(V[1])} \ar[r,"F"] \ar[d,"\widetilde{L}_{\sym}(\mu)"']& 
				{S(W[1])} \ar[d,"\widetilde{L}_{\sym}(\nu)"]
				\\
				{S(V[1])[1]} \ar[r,"{F[1]}"'] & 
				{S(W[1])[1]}
			\end{tikzcd}	
			~.
		\end{displaymath}			
	\end{definition}
	According the universal property of cocommutative graded coalgebra given by lemma \ref{lem:liftToMorph}, coalgebra morphisms 
	${\overline{S(V)}}\to {\overline{S(W)}}$
	 are in one-to-one correspondence with graded morphisms ${\overline{S(V)}} \to W$, which are collections of symmetric multilinear maps from $V$ to $W$.
	Therefore, any $L_\infty$-morphism can be translated into a collection of
	(graded) symmetric maps that are compatible with the multi-brackets:
	\begin{definition}[${L_\infty[1]}$-morphisms (graded symmetric framework)]\label{Def:LinftyMorphism-sym}
		Given two $L_\infty[1]$-algebras $(V,\mu)$ and $(W,\mu')$, a \emph{$L_\infty[1]$-morphism} between $(V,\mu)$ and $(W,\mu')$ is given by a collection of degree $0$, graded symmetric, multilinear maps
		\begin{displaymath}
			(f)= \left\{ f_k: (V[1])^{\odot k} \to W[1]\right\}_{k\geq 1}
		\end{displaymath}
		such that the auxiliary graded symmetric, multilinear operator $K_k^f: (V[1])^{\odot k} \to W[1]$
		defined as
		\begin{displaymath}
			K_m^f :=
			 \sum_{\ell=1}^m 
			\left( 
				f_{m-\ell+1} \symgerst \mu_\ell 
			- \mu_\ell' \circ \Sop_{\ell,m}(f) 
			\right)
			~,
		\end{displaymath}
		where $\Sop_{\ell,m}(f)$ is the operator defined in Notation \ref{rem:operatorSln},
		vanishes for all $k \geq 1$.
		Namely, we denote by $\Sop_{\ell,m} (f)$ the multilinear map 
		\begin{equation}\label{Eq:SoperatorSymCase}
		\Sop_{\ell,m} (f) =
		\left( 
			\sum_{\substack{k_{1}+\cdots+k_{\ell}=m\\1\leq k_{1}\leq\cdots\leq k_{\ell}}}
			(f_{k_1}\otimes\cdots\otimes f_{k_\ell})\circ 
			B_{k_1,\ldots,k_\ell}^<		
		\right)~: V[1]^{\otimes m} \to W[1]^{\otimes \ell}
		~,
		\end{equation}
		where the unshuffleator $B_{k_1,\ldots,k_\ell}^<$ runs through all $(k_1,\cdots,k_\ell)$-unshuffles $\sigma$ satisfying the extra condition
	\begin{displaymath}
		\sigma(k_1+\dots+k_{j-1}+1)<\sigma(k_1+\dots+k_{j}+1) 
		\quad \text{if}~k_{j-1}=k_j ~.
	\end{displaymath}
	\end{definition}

	After d\'ecalage, you get a similar expression for the formulation with skew-symmetric multibrackets:
	\begin{definition}[$L_\infty$-morphisms (graded skew-symmetric framework)]
\label{Def:LinftyMorphism-skew}
		An \emph{$L_\infty$-morphism} between $(V,\mu)$ and $(W,\mu')$ is given by a collection of homogeneous, graded skew-symmetric, multilinear maps
		\begin{displaymath}
			(f)= \left\{ f_k: V^{\wedge k} \to W[1-k]\right\}_{k\geq 1}
		\end{displaymath}
		such that the auxiliary graded symmetric, multilinear operator $\bar{K}_m^f: L^{\wedge m} \to L'[2-m]$
		defined as
		\begin{displaymath}
			\bar{K}_m^f :=
			 \sum_{\ell=1}^m 
			\left( 
				f_{m-\ell+1} \skewgerst \mu_\ell 
			- \mu_\ell' \circ \bar{\Sop}_{\ell,m}(f) 
			\right)
			~,
		\end{displaymath}
		vanishes for all $k \geq 1$.
		We denote by $\bar{\Sop}_{\ell,m} (f)$ the multilinear map 
		\begin{equation}\label{Eq:SoperatorSkewCase}
		\bar{\Sop}_{\ell,m} (f) =
		\left( 
			\sum_{\substack{k_{1}+\cdots+k_{\ell}=m\\1\leq k_{1}\leq\cdots\leq k_{\ell}}}
			(-)^{\sum_{i=1}^{\ell-1}(|f_{k_i}|)(\ell-i)}
			(f_{k_1}\otimes\cdots\otimes f_{k_\ell})\circ P_{k_1,\ldots,k_\ell}^<		
		\right)~,
		\end{equation}
		where the odd unshuffleator $P_{k_1,\ldots,k_\ell}^<$ runs through all ordered $(k_1,\cdots,k_\ell)$-unshuffles $\sigma$.
	\end{definition}

	\begin{remark}[Understanding $K_m^f$]
	Although the term may appear, at first glance, difficult to read, it is actually constructed according to a simple combinatorial pattern.
		\\
		The first summand	is a combination of all the possible ways to compose a bracket $\mu$ with a single component of $(f)$ to give a map of arity $m$.
		\\
		The second term is a combination of all the ways one can compose several components of the morphism $(f)$ with a primed multibracket $\mu'$ to get a map with arity $m$. 
		\\
		In particular, the term $\Sop_{\ell,m}(f)$ is a combination of all the possible ways one can tensor multiply different components of $(f)$ to give a graded (skew-)symmetric map
		$$
			\Sop_{\ell,m} : L^{\otimes (k_1 + \dots +k_\ell)}
			=
			L^{\otimes m} \to L'^{\otimes \ell} ~.
		$$
		The same applies to $\bar{\Sop}$ and $\bar{K}$ with extra signs coming from the different degrees and from the sum on odd permutations.
	\end{remark}

	\begin{remark}[$L_\infty$-morphisms in the literature]\label{rem:CorrectAttributions}
		The notion of $L_\infty$-morphism is basically old as the notion of $L_\infty$-algebra itself.
		Initially, see \cite{LadaStasheff}\cite{LadaMarkl}, only "strict" (see definition \ref{Def:strictLinfi}) $L_\infty$-morphisms and (weak) $L_\infty$-morphisms valued in a DGLA were considered.
		The explicit expression of definition \ref{Def:LinftyMorphism-sym} can be tracked in \cite[Def. 2.7]{Kajiura2006}\cite[Def. 6]{Kajiura2006b}\cite[Def. 5]{Vitagliano2013}.
		The explicit formulation in the skew-symmetric framework (\ref{Def:LinftyMorphism-skew}) has been already worked out in	 \cite[Definition 4.7.1]{Allocca2010}, see also \cite[ex. 2.20]{Ryvkin2016a}.	
	\end{remark}

	A crucial point is that the previous three definitions are equivalent:
	\begin{lemma}
		Given a $L_\infty$-morphism $(f):(V,\mu)\to(W,\nu)$ as defined in \ref{Def:LinftyMorphism-skew}, $\Dec(f)$ yields a $L_\infty[1]$-morphism from $(V[1],\Dec(\mu))$ to $(W[1],\Dec(\nu))$ in the sense of definition \ref{Def:LinftyMorphism-sym} and $L_{\sym}(Dec(f))$ yields a morphism between the corresponding Chevalley-Eilenberg complexes as given by definition \ref{Def:LinftyMorphism-coalg}.
	\end{lemma}
	\begin{proof}
		Take a $F$ as in definition \ref{Def:LinftyMorphism-coalg}, by the universal property of cocommutative coalgebras, one has that $F= L_{\sym}(f)$ for a certain collection of graded symmetric of multilinear map $(f)$.
		Being
		\begin{displaymath}
			F\circ \widetilde{L}_{\sym}(\mu) - \widetilde{L}_{\sym}(\nu) \circ F
		\end{displaymath}
		 a coderivation, the vanishing condition can be read as the following equation of multilinear maps (obtained by postcomoposition with the standard projection $\pr_W: \overline{S(W)}\twoheadrightarrow W$)
		\begin{equation}\label{Eq:intermadiateLinfinityMorphCond}
			f\circ \widetilde{L}_{\sym}(\mu) - \nu \circ F
			= 0			
			 ~.
		\end{equation}	
		In other terms, expression \eqref{Eq:intermadiateLinfinityMorphCond} lift to the above coderivation. 
		The first summand can be written as $f\symgerst \mu$ extending in the obvious way the operator $\symgerst$ to
		\begin{displaymath}
			\symgerst: M(V,W)\otimes M(V,V) \to M(V,W)
		\end{displaymath}
		as the concatenation of multilinear operators with matching domains and codomain.
		Restricting equation \eqref{Eq:intermadiateLinfinityMorphCond} to $V[1]^{\otimes m}$ one gets the equation $K^f_m = 0$ observing that term \eqref{Eq:SoperatorSymCase} is exactly the same appearing in the construction of the symmetric lift given in remark \ref{rem:operatorSln}.
		Pictorially, the whole situation is subsumed by the following commutative diagram in the category of graded vector spaces
		\begin{displaymath}
			\begin{tikzcd}
				&
				&
				W[1][1]
				\\
				&
				\overline{S(V[1])} \ar[r,"F"] \ar[dl,"q"'] \ar[d,"Q"] \ar[ur,"f"]
				&
				\overline{S(W[1])}	\ar[rd,"q'"] \ar[d,"Q'"] \ar[u,two heads]
				&
				\\
				V[1][1]
				&
				\left(\overline{S(V[1])}\right)[1] \ar[l,two heads] \ar[r,"F"]
				&			
				\left(\overline{S(W[1])}\right)[1] \ar[r,two heads]
				&
				W[1][1]
				\\
			\end{tikzcd}
		\end{displaymath}
		Note that $Q,Q'$ are in particular coderivations and $F$ is a coalgebra morphism, respectively they are the unique lift of $q,q'$ to coderivation and  of $f$ to a coalgebra morphism. This prove equivalence between definitions \ref{Def:LinftyMorphism-sym} and \ref{Def:LinftyMorphism-coalg}.
			To prove the last equivalence observe that
			\begin{displaymath}
				\Dec\left(
				\bar{K}^f_m
			\right)
			=\sum_{\ell=1}^m 
			\left( 
				\Dec(f_{m-\ell+1}) \symgerst \Dec(\mu_\ell) 
			- \Dec(\nu_\ell) \circ \Dec(\bar{\Sop}_{\ell,m}(f)) 
			\right)
			\end{displaymath}
			which equates ${K}^f_m$ upon noticing (see \ref{Cor:DecalageofTensorsProducts} in appendix) that
			\begin{displaymath}
				\Dec(f_{k_1}\otimes \dots \otimes f_{k_\ell}) =
				(-)^{\sum_{i=1}^\ell |f_{k_i}|(\ell-i)}\Dec(f_{{k_1}})\otimes \dots \otimes \Dec(f_{{k_\ell}})						
			\end{displaymath}
			whenever $|f_{k_i}|= (1-k_i \mod{2})$.
		(The last step can also be tracked in \cite[Lemma 2.18]{Ryvkin2016a}).
	\end{proof}

	\begin{remark}[Compare our notation with the literature]
		Let  us briefly discuss how our definitions of operators $K_m^f$, given in \ref{Def:LinftyMorphism-sym} and \ref{Def:LinftyMorphism-skew}, compare with the equations present in the literature mentioned in remark \ref{rem:CorrectAttributions}.
		\\
		Plugging in homogeneous elements, equation $K_m^f(x_1,\dots,x_m)=0$ in definition \ref{Def:LinftyMorphism-sym} reads as	
		
		 \begin{displaymath}
		 	\mathclap{
		 	\begin{aligned}
		 		&\sum_{\ell=1}^m \mkern-50mu
		 		\sum_{\mkern70mu \sigma \in \ush{\ell, m-\ell}}\mkern-40mu 
		 		\mkern-18mu
		 		\epsilon(\sigma) f_{m-\ell+1	}
		 		(\mu_\ell(x_{\sigma_1},\dots,x_{\sigma_{\ell}}),~ x_{\sigma_{\ell+1}},\dots x_{\sigma_\ell})
		 		=
		 		\\
		 		&=
		 		\mkern-35mu
		 		\sum_{\substack{1\leq \ell \leq m\\k_{1}+\cdots+k_{\ell}=m\\1\leq k_{1}\leq\cdots\leq k_{\ell}}}
		 		\mkern-30mu
			(-)^{\sum_{i=1}^{\ell-1}(|f_{k_i}|)(\ell-i)}
				\mkern-70mu
				\sum_{\mkern70mu\sigma \in \ush{k_1,\dots,k_\ell}^<} 
		 		\mkern-70mu
		 		\epsilon(\sigma)
		 		f_{k_1}(x_{\sigma_1},\dots, x_{\sigma_{k_1}})\otimes \dots \otimes
		 		f_{k_\ell}(x_{\sigma_{(m+1-k_{\ell})}},\dots, x_{\sigma_{m}})
		 	\end{aligned}
		 	}
		\end{displaymath}		 	
		 which coincides with \cite[Def. 1.4]{Vitagliano2013} by substituting $\ell$ with $i$ and $m-1+1$ with $j$.
		 	\\
		 Performing the same computation for the operator 	$K_m^f$ in definition \ref{Def:LinftyMorphism-skew} yields
		 \begin{displaymath}
		 	\mathclap{
		 	\begin{aligned}
		 		&\sum_{\ell=1}^m \mkern-50mu
		 		\sum_{\mkern70mu \sigma \in \ush{\ell, m-\ell}}\mkern-40mu 
		 		\mkern-18mu
		 		(-)^{|\mu_k|(m-\ell)}	 		
		 		\chi(\sigma) f_{m-\ell+1	}
		 		(\mu_\ell(x_{\sigma_1},\dots,x_{\sigma_{\ell}}),~ x_{\sigma_{\ell+1}},\dots x_{\sigma_\ell})
		 		=
		 		\\
		 		&=
		 		\mkern-30mu
		 		\sum_{\substack{1\leq\ell\leq m\\k_{1}+\cdots+k_{\ell}=m\\1\leq k_{1}\leq\cdots\leq k_{\ell}}}
				\mkern-32mu
				(-)^{\sum_{i=1}^{\ell-1}(|f_{k_i}|)(\ell-i)}
				\mkern-40mu
				\sum_{\sigma \in \ush{k_1,\dots,k_\ell}^<} 
		 		\mkern-40mu
		 		\chi(\sigma)(-)^\beta~
		 		f_{k_1}(x_{\sigma_1},\dots, x_{\sigma_{k_1}})\otimes \dots \otimes
		 		f_{k_\ell}(x_{\sigma_{(m+1-k_\ell)}},\dots, x_{\sigma_{m}})
		 	\end{aligned}
		 	}
		\end{displaymath}				 
		where $(-)^\beta$ is the sign coming from the Koszul convention when evaluating $f_{k_1}\otimes\dots \otimes f_{k_\ell}$ on homogeneous elements.
		Practically, the latter sign prefactor is computed as the Koszul sign of the following permutation of graded elements
		\begin{displaymath}
			\left(
			\substack{
				{f_{k_1},\dots,f_{k_\ell},x_{\sigma_1},\dots x_{\sigma_{m}}}
				\\
				{f_{k_1},x_{\sigma_1}\dots,x_{\sigma_{k_1}},f_{k_2},x_{\sigma_{k_1+1}}\dots,x_{\sigma_{k_1+k_2}},\dots}
		}\right)
		\end{displaymath}
		(a completely explicit expression can be found in \cite[Rem. 4.7.1]{Allocca2010}).
		Recalling that $|\mu_\ell|= 2-\ell$, $|f_{k_i}| = 1- k_i$  and noting that
		\begin{displaymath}
			\sum_{i=1}^\ell |f_{k_i}|(\ell-i) = 
			\sum_{i=1}^\ell\left( \ell -i - k_i(\ell -i)\right) =
			\dfrac{\ell(\ell-1)}{2} + \sum_{i=1}^\ell (k_i(\ell -i))
			~,
		\end{displaymath}
		one recovers the explicit definition of $L_\infty$-morphisms as stated, for example, in \cite[Def. 4.7.1]{Allocca2010} and \cite[Lem. 2.18]{Ryvkin2016a}.	
	\end{remark}

\begin{remark}[Unfolding the notion of $L_\infty$-morphisms]\label{rem:morphismsunfolded}
	Observe that in the following two "extreme" cases $\Sop$ and $\overline{\Sop}$ coincide, as multilinear operators, to
	\begin{displaymath}
		\begin{aligned}
			\Sop_{1,m}(f) &= \overline{\Sop}_{m,m}(f) = f_m
			\\
			\Sop_{m,m}(f) &= \overline{\Sop}_{m,m}(f) =\underbrace{f_1\otimes \ldots \otimes f_1}_\text{$m$-times}
			~.
		\end{aligned}
\end{displaymath}
	(That is because when $\ell=m$ then necessarily $k_1=k_2=\dots=k_m=1$ and $|f_1|=0$ in the skew-symmetric setting.)
	\\	
	Therefore the operator $K_m^f$ can be rewritten as
	\begin{displaymath}
		K_m^f := f_1 \circ \mu_m - \mu'_m \circ f_1^{\otimes m} + \kappa^m_f 
	\end{displaymath}
	where
	\begin{displaymath}
		\kappa_m^f :=
			 \sum_{\ell=1}^{m-1} 
			\left( 
			f_{m-\ell+1} \cs \mu_\ell 
			- \mu_\ell' \circ \Sop_{\ell,m}(f) 
			\right)
			~.
	\end{displaymath}
	The $L_\infty$-morphism condition, \ie the vanishing of operator $K_m^f$, can therefore be seen as the fact that the multibracket $\mu_m$, as a map from $L^{\otimes m} \to L[1]$, is not preserved by $f_1$.
	In other words, the following diagram does not commute:
	\begin{displaymath}
		\begin{tikzcd}
				V^{\otimes m} \ar[d,"\mu"] \ar[r,"f_1^{\otimes m}"]\ar[dr,phantom,"\cancel{\circlearrowleft}"]& 
				W^{\otimes m} \ar[d,"\mu'"]
				\\
				V[1] \ar[r,"f_1{[1]}"'] & W[1]
		\end{tikzcd}
		~.
	\end{displaymath}		
	Non-commutativity, however, is "mild" in the sense that its failure is explicitly controlled by the term $\kappa_m^f$ introduced above.
	The same discussion applies, mutatis mutandis, to the skew-symmetric case.
\end{remark}
The previous remark suggests to introduce the following definition:
\begin{definition}[Strict $L_\infty$-morphisms]\label{Def:strictLinfi}
	A $L_\infty$-morphism $(f)$ is called "\emph{strict}" when $f_i = 0$ for all $i\geq 2$.
	Hence, in the terms of remark \ref{rem:morphismsunfolded}, each multibracket is preserved by $f_1$.
\end{definition}

	\begin{notation}
		From now on, we will take for granted the equivalence between the three previous presentations of a $L_{\infty}$-algebra structure.
		We will denote by $(f)$ the $L_\infty$-morphism from $V\to W$ when seen as a collection of multilinear maps $f_k:\underline{\Hom}(V^{\otimes k},W)$, hence with a slight abuse of notation when considering the multibrackets symmetric or skew-symmetric.
		By $F\in \Hom_{\text{coAlg}}(\overline{S(V)},\overline{S(V)})$ we denote the lift of the symmetric multibrackets $(f_1,\dots,)$ to a coalgebra morphism, which is in particular also a cochain map.
	\end{notation}

\subsection{Composition of $L_\infty$-morphisms}
Once the concept of morphism between $L_\infty$-algebras has been clarified, it is crucial to specify how such morphisms act under composition to completely understand $L_\infty$-algebras as a category.
\\
In the coalgebraic formulation there is a natural definition of composition of $L_\infty$-morphisms simply given by the composition of the two corresponding differential coalgebra morphisms.
This easily translates to the following definition in the multibrackets setting:
	\begin{definition}[$L_\infty$-morphisms composition (multibrackets framework)]\label{Def:CompositionFormula}
		If $(f):L\to L^{\prime}$ and $(g):L^{\prime}\to
		L^{\prime\prime}$ are morphisms of $L_{\infty}$-algebras, the \emph{composition
		$(g)\circ (f):L\to L^{\prime\prime}$} is a $L_\infty$-morphism defined as 
		$(g \circ f)_k:=\{(g\MBComp f)_{k}\}_{k\in\mathbb{N}}$, where
			\begin{displaymath}
				(g \circ f )_m :=
				\pr_{L^{\prime\prime}} \circ \widetilde{L}_{\sym}(g) \circ \widetilde{L}_{\sym}(f) \eval_{L^{\otimes m}}
				=
				\sum_{\ell=1}^m g_\ell 	\circ \Sop_{\ell,m}(f)
			\end{displaymath}
			and $\Sop_{\ell,m}(f)$ is defined in equations \eqref{Eq:SoperatorSymCase} ,\eqref{Eq:SoperatorSkewCase} in the symmetric and skew-symmetric framework respectively.
	\end{definition}
	(One can retrieve this formula decalaging the analogue formula in \cite[\S 1]{Vitagliano2013}. See also \cite[Cor. 1.3.3 and Prop 1.5.3]{Manetti-website-coalgebras}.)
	\\
	Pictorially, the composition of two $L_\infty$-morphisms is given by the commutativity of the following diagram in the graded vector space category: 
		\begin{displaymath}
		\begin{tikzcd}[column sep=huge]
			\overline{S(L[1])} \ar[r,"L_{\sym}(f)"] \ar[dr,"f"']
			\ar[rr,bend left=30,"L_{\sym}\big(g\circ L_{\sym}(f)\big):=L_{\sym}(g\circ f)"]
			&
			\overline{S(L'[1])}\ar[r,"{L_{\sym}(g)}"]\ar[dr,"g"']\ar[d]
			&
			\overline{S(L''[1])}\ar[d]
			\\
			&
			L'[1]
			&
			L''[1]
		\end{tikzcd}
	\end{displaymath}
	Recall that the identity coalgebra morphism $\Unit: \overline{S(V)} \to \overline{S(V)}$ is precisely given by the lift of the projection $\pr_V:\overline{S(V)}\twoheadrightarrow V$ to a coalgebra morphism.
	Hence, it corresponds to the following definition in the graded coalgebra setting:
		\begin{definition}[Identity $L_\infty$-morphism]
		Consider an $L_\infty$-algebra $L$, the \emph{identity morphism} pertaining to $L$ is the strict $L_\infty$-morphism $(\mathrm{id}):L\longrightarrow L$ given by:
		\begin{displaymath}
			\begin{aligned}
				\mathrm{id}_{1}= \mathrm{id}_L  &  :L\longrightarrow L \\
				\mathrm{id}_{k}=0  &  :L^{\otimes k}\longrightarrow L
				\qquad \forall ~k \geq 1
			\end{aligned}
		\end{displaymath}
	\end{definition}
In the same spirit one gets the following notion of $L_\infty$-isomorphism:
	\begin{definition}[$L_\infty$-isomorphisms]
		Consider two $L_\infty$-algebras $L$ and $L'$.
		An $L_\infty$-morphism $(f):L\to L^{\prime}$ is an \emph{$L_\infty$-isomorphism} if the corresponding lift $F:\overline{S(L[1])}\to \overline{S(L^{\prime}[1])}$ is an invertible coalgebra morphism, \ie if $f_1:L_1\to L_1^{\prime}$ is a chain complex isomorphism.
	\end{definition}
	A characterization in the multibrackets frameworks is given by the following lemma
	\begin{lemma}\label{Lem:InvertibleTensorColalgebraMorphisms}
		A coalgebra morphism $F: \overline{S(V)} \to \overline{S(W)}$ is invertible if and only if $f_1: V \to W $ is an invertible graded morphism.
	\end{lemma}
	\begin{proof}
		Being $F$ invertible, there exists $G: \overline{S(W)} \to \overline{S(V)}$ such that $F\circ G = \Unit$ and $G\circ F = \Unit$.
		These two compositions can therefore be characterized as lifts of the standard projection $p:\overline{S(X)} \to X$ with $X$ equal to $W$ and $V$ respectively, hence $f_1 \circ g_1 = \id$ and $g_1 = f_1^{-1}$.
		
		On the opposite, given a $F$ such that $f_1$ is invertible, one can construct an inverse $G$ of $F$ iteratively as the lift of the following components:
		\begin{displaymath}
			\begin{cases}
				g_1 = f_1^{-1}
				\\
				g_2 = -f_1^{-1} \circ f_2 \circ \left(f_1^{-1}\otimes f_1^{-1}\right)
				\\
				\vdots
				\\
				g_m = -f_1^{-1} \circ f_m \circ \left(f_1^{-1}\right)^{\otimes m}
				\displaystyle
				-
					\sum_{\ell=2}^{m-1} g_\ell \circ  \Sop_{\ell, m}(f)
				\circ \left(f_1^{-1}\right)^{\otimes m}
				\qquad ~.
			\end{cases}
		\end{displaymath}
	\end{proof}
	Recalling that every $L_\infty$-algebra has an underlying chain complex, it is natural to introduce the notion of \emph{quasi-isomorphism}
	\begin{definition}[$L_\infty$-quasi-isomorphism]\label{Def:LinfintyQuasiIso}
		A $L_\infty$-morphism $(f):(L,\mu)\to (L',\mu')$ is called a \emph{quasi-isomorphim} if $f_1: (L,\mu_1)\to (L',\mu_1')$ is a quasi isomorphism of chain complexes, \ie $f_1$ induces an isomorphism of the corresponding cohomology groups.	
	\end{definition}
	
\begin{remark}[$L_\infty$-algebras cohomology]
	Observe that, for any given $L_\infty$-algebra $(L,\mu)$, one has two naturally associated cohomologies:
	the cohomology of the underlying chain complex $(L,\mu_1)$, obtained by neglecting all multibrackets of arity greater than one (see example \ref{ex:cochaincomplexLinfinity}), and the cohomology of the Chevalley-Eilenberg complex $(S(L[1]),\widetilde{L}_{\sym}(\Dec(\mu)))$ (see definition \ref{def:CELinfty1struct}).
	Although the former one has been used in definition \ref{Def:LinfintyQuasiIso} to characterize quasi-isomorphism, the latter is more apt to be intended as the proper cohomology of $(L,\mu)$.
	More refined notions of \emph{$L_\infty$-algebra cohomology} are expressed in terms of representations, we point out \cite{Reinhold2019} for an extensive introduction on this topic.
\end{remark}

\begin{remark}[$L_\infty$-homotopies]
	 We have shown that understanding an $L_\infty$-algebra as a (commutative cofree) differential coalgebra  readily rewards a natural notion of $L_\infty$-morphism as a graded morphism respecting at the same time the codifferential and the coalgebra structure.
	 \\
	 Similarly, it might seem just as natural to introduce the notion of
	 \emph{$L_\infty$-homotopy} between two $L_\infty$-morphisms, seen as graded differential coalgebra morphisms $F,G: (C,Q)\to (C',Q')$, as a certain chain homotopy $H$, in the sense of equation \eqref{eq:chainhomotopy}, together with suitable coalgebraic conditions.
	 For instance, it could be given by a degree $-1$ \emph{(F,G)-coderivation} $H$ from $C$ to $C'$ , \ie $\Delta' \circ H = ( F \otimes H + H \otimes G)\circ \Delta)$, 
	 together with properties
	 \begin{align*}
 		H Q + Q' H =&~ F-G
 		~,
 		\\
 		(H\otimes (G Q) + (F Q)\otimes H ))\Delta =&~ 0
 		~.
	 \end{align*}
	However, such a kind of definition is unsatisfactory from the more conceptual point of view of \emph{homotopy theory} (see \cite[Remark 4]{Dotsenko2016}).
	In this language, the first step would be to identify a \emph{model structure} for the category of $L_\infty$-algebras. An account on many different models, and the discussion of their equivalence, can be found in \cite{Pridham2010}.
	\\
	An useful and explicit model of \emph{weak equivalence} (homotopies or 2-morphisms) has been given by Dolgushev in \cite{Dolgushev2007}.
	Namely, given two $L_\infty$-algebras $L$ and $L'$, there is an auxiliary \emph{pro-nilpotent} $L_\infty$-algebra $\mathcal{U}(L,L')$, with $k$-multibrackets denoted by $\lbrace\dots\rbrace_k$, encoding $L_\infty$-morphisms as \emph{Maurer-Cartan} elements, \ie
	\begin{displaymath}
		\Hom_{Lie\infty}(L,L') = 
		MC(\mathcal{U}(L,L')) :=
		\left\lbrace 
			f \in \mathcal{U} ~\left\vert\quad \sum_{n=1}^\infty \dfrac{1}{n!}\lbrace f,\dots, f \rbrace_n
		\right\rbrace\right.
		~.
	\end{displaymath}
	Note that the pro-nilpotency condition guarantees convergence of the above sum.
	Hence, two $L_{\infty}$-morphisms are \emph{equivalent} if there corresponding Maurer Cartan elements in $\mathcal{U}$ are equivalent.
	(See also \cite[Appendix A]{Fregier2015}).
\end{remark}

	\begin{remark}[The category of $L_\infty$-algebras]
		The discussion up to this point can be summarized by saying that $L_{\infty}$-algebras build up a subcategory of differential graded vector spaces, the insertion is given by "forgetting" all multibrackets of arity greater than two.
		More pedantically, we have constructed three $L_\infty$-algebras categories:
		\begin{enumerate}
			\item[(1)] Objects are graded vector spaces endowed with graded skew-symmetric multibrackets (definition \ref{Def:LInfinityStasheff} or \ref{Def:LInfinityTony}), together with morphisms given by definition \ref{Def:LinftyMorphism-skew};
			\item[(2)] Objects are graded vector spaces endowed with graded symmetric multibrackets (definition \ref{Def:LInfinityShifted}), together with morphisms given by definition \ref{Def:LinftyMorphism-sym};
			\item[(3)] Objects are cofree cocommutative differential coalgebras, morphisms are differential coalgebra morphisms.
\end{enumerate}	
		The three categories are isomorphic	\footnote{Isomorphism of categories is an extremely stronger notion than \emph{equivalence of category}. That is why we talk about different presentations of the same concept rather than three equivalent categories.}.
		The invertible functor $(1)\to (2)$ acts as the shift endofunctor on graded vector spaces and as the d\'ecalage $\Dec$ on multibrackets and the components of any given $L_\infty$-morphism.
		The invertible functor $(2)\to (3)$ acts as the symmetric tensor coalgebra functor on graded vector spaces, as the lift to coderivations $\widetilde{L}_{\sym}$ on multibrackets, and as the lift to coalgebra morphism on the components of a $L_\infty[1]$-morphism.
	\end{remark}

	\begin{remark}[Other formalization of $L_\infty$-algebra structures]
		It is worth to mention, without claiming to be exhaustive, 
		three other extremely fruitful frameworks to understand homotopy algebras structures.
		Both encode $L_\infty$-algebras as a special case, generalizing the concept in quite different directions.
		\begin{itemize}
			\item
				\emph{A $L_\infty$-algebra is a $L_\infty$-algebroids over a point (\ie $0$-dimensional)} \cite[Appendix A]{Sati2012a}. 
		 The latter can be encoded in terms of Graded manifolds. Namely a \emph{$L_\infty$}-algebroid is $Q$-manifold that is a $\mathbb{Z}$-graded manifold endowed with a homogeneous vector field $Q$ homogeneous in degree one and such that $Q^2=0$ (\emph{homological vector field}). 
		 See \cite[\S 1]{Jurco2018} for an introduction to this approach with application to prequantum field theories.

			\item \emph{The category of $L_\infty$-algebras is the vertical categorification of ordinary Lie algebras}. This categorical construction has been extensively discussed in \cite{Baez2003} for the case of $L_2$-algebras. The general idea can be found in \cite{nlab:l-infinity-algebra}.
			
			\item
				\emph{A $L_\infty$-algebra is an algebra over the homotopy Lie operad in the category of chain complexes} \cite{Markl1998}.				
				The theory of operads is a vast topic; a full account is contained in \cite{Loday}, a motivational introduction can be found in \cite{Vallette2014}.
				For a contained exposition mainly geared toward the definition of $L_\infty$-algebras as as an algebra over the $S$ operad see \cite[\S 2]{Kimura1995}.
		\end{itemize}
		While we will not use this machinery in what follows, it is interesting to note that working with the algebraic framework of multilinear maps with  composition (\RN products or Gerstenhaber product in the non-symmetric case), as it will be done in chapter \ref{Chap:MarcoPaper}, can be seen as a little step toward the operadic approach.
		The point is recognizing the role of the \emph{endomorphism operad} \cite[\S 2.2.]{Vallette2014} as the appropriate setting for studying the composition of multi-linear maps.
		
	\end{remark}

	\subsection{Examples}\label{SubSection:studycase}
		In this subsection, we focus on some particular classes of $L_\infty$-structures that will appear recurrently in the next chapters.
		
		In section \ref{Sec:hcmm}, we will recall the notion of \emph{\momap}. 
		This is a particular $L_\infty$-morphism from an ordinary Lie algebra:
\begin{example}[$L_\infty$-morphisms defined on Lie algebras]\label{Rem:LftyMorphasChainMap}
	Recall that in example \ref{Ex:dglaAsLinfinity} we have shown how a differential graded Lie algebra can be seen as a $L_\infty$-algebra where all multibrackets in arity greater than two are trivial.
	\\
	Consider an ordinary Lie algebra $\mathfrak{g}$, \ie a differential graded Lie algebra concentrated in degree $0$, and an arbitrary $L_\infty$-algebra $(L,\mu)$.
	An $L_\infty$-morphism $(f):\mathfrak{g}\to (L,\mu)$, reads, in the skew-symmetric presentation, as a collection of skew-symmetric multilinear maps $f_k:\wedge^k \mathfrak{g}\to L^{1-k}$ such that 
	\begin{displaymath}
		f_{m}\ca [\cdot,\cdot] = \sum_{\ell=1}^m \mu_\ell \circ \overline{\Sop}_{\ell,m}(f)	
		\qquad \forall m\geq 2
		~.
	\end{displaymath}		
	Recognizing that
	\begin{displaymath}
		\begin{aligned}
		(f_{m}\ca [\cdot,\cdot])(x_1,\dots,x_{m+1})=&~
		\mkern-40mu\sum_{\mkern30mu\sigma \in \ush{2,m-1}} \mkern-40mu
		(-)^\sigma f_m([x_{\sigma_1},x_{\sigma_2}],x_{\sigma_3},\dots x_{\sigma_{m+1}} )
		=
		\\
		=&~ - (f_{m}\circ \partial_{\CE}) (x_1,\dots,x_{m+1})
		\end{aligned}
	\end{displaymath}
	the previous condition can be read as the coboundary condition
	\begin{displaymath}
		\delta_{CE} f_m = - \sum_{\ell=1}^m \mu_\ell \circ \overline{\Sop}_{\ell,m}(f)	
		\qquad \forall m\geq 2
	\end{displaymath}
	in the Chevalley-Eilenberg cochain complex over the trivial representation $\rho: \mathfrak{g}\to 0$ (see reminder \ref{Rem:CEconventions}).
	\\
	When $L$ is in particular a cochain complex, \ie the only non-trivial multibracket is $\mu_1$, the previous situation is expressed by the commutation of the following diagram (in the category of ordinary vector spaces):
	\begin{displaymath}
		\begin{tikzcd}
			\cdots \ar[r,"\partial"] &
			\Lambda^k\mathfrak{g} \ar[r,"\partial"] \ar[d,"f_k"]&
			\Lambda^{k-1}\mathfrak{g} \ar[r,"\partial"]\ar[d,"f_{k-1}"] &
			\cdots \ar[r,"\partial"] &
			\mathfrak{g} \ar[r,"\partial"]\ar[d,"f_1"] &
			0
			\\
			\cdots \ar[r,"\mu_1"] &
			L^{-k} \ar[r,"\mu_1"] &
			L^{-k+1} \ar[r,"\mu_1"] &
			\cdots \ar[r,"\mu_1"] &
			L^0 \ar[r,"\mu_1"] &
			\cdots
		\end{tikzcd}
	\end{displaymath}
	hence it can be seen as a chain map between the Chevalley-Eilenberg complex, with a suitable degree re-parametrization, and $(L,\mu_1)$, namely
	\begin{displaymath}
		(f): \left(S( \mathfrak{g}[1])\right)[-1] \to (L,\mu_1)
		~.
	\end{displaymath}
\end{example}			
		
		Most of the subsequent work revolves around the notion of \emph{(higher) observables $L_\infty$-algebras} (see section \ref{Section:RogersObservables}). 
		These $L_\infty$-algebras enjoy the convenient properties to be \emph{grounded}, \ie all multibrackets with arity greater than one are non-trivial only on the "ground degree" zero (definition \ref{Def:groundedLinfinity}), and 
		concentrated in finitely many degrees (negative degrees in the cohomological convention and positive degrees in the homological convention).
		\\		
		In the case of a grounded $L_\infty$-algebra the axioms defining multibrackets and morphisms are considerably simpler:
\begin{example}[Axioms of grounded $L_\infty$-algebras {\cite[Rem. 6.2]{Ryvkin2018}}]\label{Rem:GroundedEasyAxioms}
	The higher Jacobi equations for a grounded $L_\infty$-algebra $(L,\{\mu_i\}_{i\geq 1})$ (in the skew-symmetric multibrackets presentation) read, by degree reasons, as
	\begin{displaymath}
		0 =\mu_k \ca \mu_2 - \mu_1 \circ \mu_{k+1} \qquad (\forall k\geq 1)
		~.
	\end{displaymath}
	In the spirit of example \ref{Rem:LftyMorphasChainMap}, the previous equation can be recast as
	\begin{displaymath}
		\delta (\mu_k ) = \d \circ \mu_{k+1} \qquad (\forall k\geq 1)
		~,
	\end{displaymath}
	denoting $\mu_1 = \d$ and 
	\begin{displaymath}
		\big(\delta (\mu_k)\big)(x_1,...,x_{k+1})
		=\sum_{i<j}(-1)^{i+j}\mu_k(\mu_2(x_i,x_j), x_1 ,\dots,\widehat{x_i},\dots,\widehat{x_j},\dots,x_{k+1})~.
\end{displaymath}		
	(See \cite[Rem. 6.2]{Ryvkin2018} or \cite[\S 4]{Reinhold2019} for the definition of the Chevalley-Eilenberg coboundary $\delta$ in the context of $L_\infty$-algebras).
\end{example}		
\begin{example}[Morphisms into a grounded $L_\infty$-algebras {(\cite[Lem. 2.36]{Ryvkin2016a})}]\label{Rem:GroundedEasyMorphisms}
Also morphisms into a ground $L_n$ algebra (see definition \ref{Def:groundedLinfinity}) take a particular form.
Consider $(f):(L,\mu)\to (L',\mu')$ with $L'$ grounded, according to definition \ref{Def:LinftyMorphism-skew}, the following equation must hold for $m\geq 1$:
\begin{displaymath}
	\begin{aligned}
		0 =&
		\left(\sum_{\ell=1}^m f_{m-\ell+1}\ca \mu_\ell\right)
		- \mu_1'\circ \Sop_{1,m}(f) - \mu_m' \circ \Sop_{m, m}(f) =
		\\
		=&
		\left(\sum_{\ell=1}^m f_{m-\ell+1}\ca \mu_\ell\right)
		- \mu_1\circ f_m - \mu_m' \circ f_1^{\otimes m}		
	\end{aligned}
\end{displaymath}
When studying infinitesimal Lie algebra actions, we will be particularly interested in $L_\infty$-morphisms from a Lie algebra into a grounded $L_\infty$-algebra.
In this case the previous equation reads as follows:
\begin{displaymath}
	\begin{aligned}
		0 =&
		f_{m-1}\ca \mu_2
		- \mu_1 \circ f_m - \mu_m' \circ f_1^{\otimes m}			 =
		\\
		=&
		- \delta_{CE} \left( f_{m-1} \right)
		- \mu_1 \circ f_m - \mu_m' \circ f_1^{\otimes m}		
		~.	
	\end{aligned}	
\end{displaymath}
\end{example}	
Note also that, from the $L_\infty$-perspective, $L_n$-algebras are considerably easier to handle. 
Being concentrated in degrees $\lbrace k \in \mathbb{Z}~\vert -(n-1) \leq k \leq q\rbrace$ means that a $L_n$-algebra involves only finitely many non-trivial brackets $\mu_1,\dots,\mu_{n+1}$ satisfying only finitely many higher Jacobi equations $J_1,\dots, J_{n+1}$ (see for instance \cite[Lem. 2.35]{Ryvkin2016a}).

In the following we will also make use of the following construction which
allows to "pushforward" $L_\infty$-structures on a given graded vector space $V$ along any collection of multibrackets $(p_i: V^{\odot i} \to V)_{i\geq 1}$.
\begin{remark}\label{rem:exp}
Let $(V,\mu)$ be an $L_{\infty}[1]$-algebra, and denote the corresponding codifferential by $Q:=\widetilde{L}_{\sym}(\mu)$. 
A degree $0$ linear map
$p\colon \overline{S(V)}\to V$ gives rise to a degree $0$ coderivation $C_p:=\widetilde{L}_{\sym}(p)$ of $\overline{S(V)}$. 
In turn, by Remark \ref{rem:codermor}, $C_p$ determines an isomorphism of coalgebras
$e^{C_p}$.
\\
Consider now $$Q':=e^{C_p}\circ Q\circ e^{-C_p} ~.$$
One checks easily that 
\begin{itemize}
\item $Q'$ is also a codifferential on $\overline{S(V)}$, and thus corresponds to a new $L_{\infty}[1]$-algebra structure $\mu'$ on $V$,
\item $e^{C_p}$ intertwines $Q$ and $Q'$, hence it corresponds to
an $L_{\infty}[1]$-isomorphism $f$ from $(V,\mu)$ to $(V,\mu')$.
\end{itemize}
Explicitly, one has
 \begin{align*}
	m'&=~\pr_V(Q') 
	~=
	\\
	&=~	
	\pr_V\left(Q+[C_p,Q]+\frac{1}{2!}[C_p, [C_p,Q]]+\dots \right) 
	~=
	\\
 	&=~ m+[p,m]_{\cs}+\frac{1}{2!}[p, [p,m]_{\cs}]_{\cs}+\dots
\end{align*}
and
$$f=\pr_V(e^{C_p})=pr_V+p+\frac{1}{2!}(p\cs p)+\dots,$$
where $[\;,\;]_{\cs}$ denotes the graded commutator with respect to the symmetric \RN product $\cs$.
Under d\'ecalage one gets a corresponding result in the skew-symmetric framework.
\end{remark}

	In chapter \ref{Chap:MarcoPaper} we will discuss the so-called \emph{Vinogradov $L_\infty$-algebra}. 
	That is a (not grounded) $L_\infty$-algebra realized from a differential graded Lie algebra
	applying the following construction:
\begin{theorem}[Getzler {\cite[Thm. 3]{Getzler1991}}]\label{Thm:Getzler}
	Given a differential graded Lie algebra $(L,d,[\cdot,\cdot])$
	\begin{displaymath}
		\begin{tikzcd}
			\cdots \ar[r] & L_{-1}\ar[r,"\d"] & L_0 \ar[r,"\d"] & L_1 \ar[r]& \cdots
		\end{tikzcd}
	\end{displaymath}
	the truncation of its underlying chain complex in negative degrees 
	\begin{displaymath}
	 \mathbb{L} : = 	\trunc_{0} L =
	 \begin{cases}
	 	L_i~, & i< 0 
	 	\\
	 	0~, & i \geq 0
	 \end{cases}
	\end{displaymath}
	constitutes a $L_\infty[1]$-algebra with the following symmetric multibrackets
	\begin{displaymath}
		\lbrace a \rbrace_1 = 
		\begin{cases}
			\d a ~, & |a| < 0 \\
			0 ~, & |a| \geq 0		
		\end{cases}
	\end{displaymath}
	and
	\begin{displaymath}
		\begin{aligned}
		\lbrace a_0,\dots, a_n\rbrace_{n+1} =&
		b_n  \left( \underbrace{[\cdot,\cdot]\symgerst \dots ([\cdot,\cdot]}_{n\text{ times}} \symgerst D ) \right)	(a_0,\dots,a_{n})
		=
		\\
		=&
		b_n\left(\sum_{\sigma \in S_{n+1}} \epsilon(\sigma) [[\dots,[[D a_{\sigma_0},a_{\sigma_1}],a_{\sigma_2}]\dots],a_{\sigma_n}]\right)
		\end{aligned}
	\end{displaymath}			
	where 
	\begin{displaymath}
			D := \d - \{\cdot\}_1 = \d \circ \pi_g ~, 
	\end{displaymath}		
	$\pi_g: L \twoheadrightarrow L_{-1}$ is the projection on the "ground" degree,
	and
	\begin{displaymath}
		b_n = (-)^n \dfrac{B_n}{n!}
	\end{displaymath}
	is a numerical constant containing the $n$-th Bernuolli number $B_n$ (see  \cite{Weisstein} for a survey on the definition and property of Bernoulli numbers).
\end{theorem}
\begin{remark}[Getzler theorem in the literature]
	Observe that the previous result has been originally stated in the \emph{homological convention} \cite{Getzler1991}.
	In the same preprint, it is also explained how it can be obtained as a corollary of \cite[Theorem 3.1]{Fiorenza2006}.
	Namely, Fiorenza and Manetti has shown how to naturally endow the mapping cone of a given DGLA morphism $\chi:(L,\d,[\cdot,\cdot]) \to (M,\d',[\cdot,\cdot]')$, \ie the differential graded vector space $(C(\chi),\d_C)$ where
	\begin{displaymath}
		C(\chi) = L\oplus M[1]
		\quad,\qquad
			\left(
		\morphism{\d_c}
		{L^i\oplus M^{i-1}}
		{L^{i+1}\oplus M^{i}}
		{\pair{x}{y}}
		{\pair{\d x}{\chi(x)- \d' y}}
		\right)
		~,
	\end{displaymath}
	with a $L_\infty$-structure via \emph{homotopy transfer} (see \cite[\S 4,5]{Fiorenza2006}). 
	In \cite[Rem. 5.7]{Fiorenza2006}, they also discuss the, somewhat, "miraculous" appearance of the Bernoulli numbers.
\end{remark}
Although we will not make explicit use of it in the following, we want to mention a procedure for constructing algebras which, in fashion similar to theorem \ref{Thm:Getzler}, realizes higher multibrackets by iteration on multibrackets in lowest arity.
	\begin{example}[Derived brackets]
		A huge class of examples of $L_\infty$-algebras can be produced with  the so-called \emph{derived bracket construction} due to Voronov \cite{Voronov2005}.
		Namely, it is possible to show how to construct explicitly a $L_{\infty}$-structure out of  any quadruple $(L,\mathcal{a},P,\Delta)$, called \emph{V-data}, composed of
		\begin{itemize}
			\item a graded Lie algebra $L=(L,[\cdot,\cdot])$;
			\item an Abelian Lie subalgebra $\mathcal{a}\hookrightarrow L$;
			\item a projection $P:L\twoheadrightarrow \mathcal{a}$ whose kernel is a Lie subalgebra of $L$;
			\item an homogeneous, degree $1$, element $\Delta \in \ker(P)_1$ such that $[\Delta,\Delta]=0$.
		\end{itemize}
		See \cite[Thm 2]{Voronov2004} or \cite[\S 1]{Fregier2015b} for a condensed survey.
	\end{example}

Observe that given any two $L_\infty$-algebras one can construct a third $L_\infty$-algebra given by their direct sum:
	\begin{definition}[Direct sum of $L_\infty$-algebras {\cite[ex. 6.5]{Doubek2007}}]
		Given any two $L_\infty$-algebras $V=(V,\mu)$ and $W=(W,\nu)$, we call \emph{direct sum of $V$ and $W$} the $L_\infty$-algebra $V\oplus W$ 
		with underlying vector space given by the direct sum $V\oplus W$ of the underlying vector spaces of $V$ and $W$ and with multibrackets given, for any $k\geq 1$, by
	\begin{displaymath}
		\morphism{\mu_k^\oplus:=\mu_k\oplus \mu_k'}
		{(V\oplus W)^{\wedge k}}
		{V\oplus W}
		{\pair{v_1}{w_1}\otimes \dots \pair{v_k}{w_k}}
		{\pair{\mu_k(v_1,\dots,v_k}{\mu'_k(w_1,\dots,w_k)}}
		~.
	\end{displaymath}		
\end{definition}
An important results show how any $L_\infty$-algebra can be equivalently expressed as the direct sum of two $L_\infty$-algebras of a "simpler" type:
	\begin{theorem}[Decomposition theorem {\cite[Lemma 4.9]{Kontsevich2003}}]
		Any $L_\infty$-algebra is $L_\infty$-isomorphic to the direct sum of a \emph{minimal $L_\infty$-algebra} (the unary bracket is zero) with a \emph{linear contractible $L_\infty$-algebra} (all brackets are zero except the unary one and all cohomology groups vanish).
		\end{theorem}

\ifstandalone
	\bibliographystyle{../../hep} 
	\bibliography{../../mypapers,../../websites,../../biblio-tidy}
\fi

\cleardoublepage


%% file: chapters/multisymplectic/multisymplectic.tex
\chapter{Multisymplectic manifolds and symmetries}\label{Chap:MultiSymplecticGeometry}
\emph{Multisymplectic structures} (also called \emph{``$n$-plectic''}) are a rather straightforward generalization of symplectic ones where closed non-degenerate $(n+1)$-forms replace $2$-forms.
\\
Historically, the interest in multisymplectic manifolds, \ie smooth manifolds equipped with an $n$-plectic structure,  has been motivated by the need for understanding the geometrical foundations of first-order classical field theories.
The key point is that, just as one can associate a symplectic manifold to an ordinary classical mechanical system (\eg a single point-like particle constrained to some manifold), it is possible to associate a multisymplectic manifold to any classical field system (\eg a continuous medium like a filament or a fluid).
It is important to stress that mechanical systems are not the only source of inspiration for instances of this class of structures. For example, any orientable $n$-dimensional manifold can be considered $(n-1)$-plectic when equipped with a volume form and semisimple Lie groups have a natural interpretation as $2$-plectic manifolds.
\\
As proposed by Rogers in \cite{Rogers2010} (see also cite{Zambon2012}), this generalization can be expanded by introducing a higher analogue of the Poisson algebra of smooth functions (also known as ``observable algebra'')  to the multisymplectic case.
Namely, he proved that observables on a multisymplectic manifold can be algebraically encoded by an $L_{\infty}$-algebra; that is a graded vector space endowed with skew-symmetric multilinear brackets satisfying the Jacobi identity up to coherent homotopies.
\\
The latter concept allowed for a natural extension of the notion of moment map, called \emph{\momap}, originally defined in \cite{Callies2016}, associated to an infinitesimal action of a Lie group on a manifold preserving the multisymplectic form.

This chapter delivers a self-contained survey to the theory of multisymplectic manifolds, as defined by Catrijn, Ibort, and de Le\'on \cite{CatIbort}, adopting the algebraic framework proposed originally by Rogers in order to study symmetries admitting \momaps.
Most results can be found in the literature \cite{Rogers2010,Callies2016,Fregier2015,zbMATH06448534} but we present some of the proofs for a clearer and self-contained exposition.
For general background on symplectic geometry and (co)momentum maps we quote, among others 
\cite{Abraham1978,GS84,Arn-Khe}.

\section{Reminder on multi-Cartan calculus}\label{Sec:MultiCartan}
From a purely algebraic point of view, one can define the \emph{Cartan calculus} on every pair $(\mathfrak{g}, A)$ composed by a Lie algebra $\mathfrak{g}$ and a graded associative algebra $A$ as a triple of operations
	\begin{displaymath}
		\d \in \Der^1(A) ~,
		\qquad
		\mathcal{L} : \mathfrak{g} \to \Der^{0}(A) ~,
		\qquad
		\iota: \mathfrak{g} \to \Der^{-1}(A) ~;
	\end{displaymath}
satisfying the following commutation rules for any $x,y\in\mathfrak{g}$:
	\begin{eqnarray}
		[\d , \d] &=& 0
		\\
		{[\mathcal{L}_x , \d]} &=& 0
		\\
		{[\iota_x , \d]} &=& \mathcal{L}_x \label{Eq:CartanMagic}
		\\
		{[\iota_x , \iota_y]} &=& 0
		\\
		{[\mathcal{L}_x , \mathcal{L}_y]} &=& \mathcal{L}_{[x,y]_\mathfrak{g}}
		\\
		{[\mathcal{L}_x , \iota_y]} &=& \iota_{[x,y]_{\mathfrak{g}}}
	\end{eqnarray}
where $[\cdot,\cdot]_{\mathfrak{g}}$ denotes the Lie bracket on $\mathfrak{g}$ and $[\cdot,\cdot]$ denotes the graded commutator between endomorphisms on $A$. (See \cite[\S 3]{Meinrenken2004}).

Specializing this to the geometric setting of smooth manifolds, \ie taking as the Lie algebra the $C^\infty(M)$-module of vector fields and considering the graded commutative $C^\infty(M)$-algebra of differential forms 
	\begin{displaymath}
		\mathfrak{g} = \mathfrak{X}(M) = \Gamma(TM) ~,
		\qquad
		A^k = \Omega^k(M) = \Gamma(\Lambda^k T^\ast M) ~, 	
	\end{displaymath}
we get the standard Cartan calculus on smooth manifold where $\d$ is the \emph{de Rham} differential, $\mathcal{L}_v$ is the \emph{Lie derivative} along the vector field $v$ and $\iota_v$ is the \emph{insertion (contraction or inner composition)} of the vector field $v$ into a given form.

When working on manifolds with a fixed ("preferred") differential form in degree greater than two, computations involving several vector fields come immediately into play. 
Hence, we are led to generalize the previous apparatus replacing $\mathfrak{g}$ with its corresponding \emph{Gerstenhaber algebra} (see remark \ref{Rem:GerstehaberofLieAlg}) or, in the geometric setting, to consider the algebra of \emph{multi-vector fields}.

\begin{reminder}[Gerstenhaber Algebra structure on $CE(\mathfrak{g})$]\label{Rem:GerstehaberofLieAlg}
	Given a Lie algebra $(\mathfrak{g},[\cdot,\cdot])$, consider the Chevalley-Eilenberg chain complex $CE(\mathfrak{g})$ given by $\Lambda^{\geq 1} (\mathfrak{g}[1])$ together with the boundary operator $\partial$ given in equation \eqref{eq:CE_boun} (see reminder \ref{Rem:CEconventions}).
	
	One can introduce on $CE(\mathfrak{g})$ a skew-symmetric bilinear operator of degree $-1$ called \emph{Schouten bracket} defined on decomposable elements by
		\begin{equation}\label{Eq:Schouten}
			\begin{aligned}
			[x_1 & \wedge\dots\wedge x_m , y_1 \wedge\dots \wedge y_n ]
			=\\
			&=
			\mkern-10mu\sum_{\substack{1\leq i \leq m \\ 1 \leq j \leq n}}\mkern-10mu
			(-)^{i+j} [x_i,y_j] ~\wedge x_1\wedge\dots\wedge \widehat{x_i} \wedge \dots \wedge x_m
			 \wedge y_1\wedge\dots\wedge \widehat{y_j} \wedge \dots \wedge y_n 
			~.
			\end{aligned}
		\end{equation}
	and then extended by linearity on the whole graded vector space. 
	This operation makes $CE(\mathfrak{g})$ into a \emph{Gerstenhaber algebra}, in particular the following properties are satisfied
	\begin{displaymath}
		[x,[y,z]] = [[x,y],z] + (-)^{(|x|-1)(|y|-1)}[y,[x,z]]
		\tag{Graded Jacobi}
	\end{displaymath}
	\begin{displaymath}
		[x,y \wedge z] = [x,y] \wedge z + (-)^{(|x|-1)|y|} y \wedge [x,z]
		\tag{Graded Poisson}~.
	\end{displaymath}
	Observe that the Schouten bracket and the Chevalley-Eilenberg operator (equation \eqref{eq:CE_boun}) satisfy the following compatibility relation (see for example \cite[Lemma 3.12]{Ryvkin2016}) for any $p,q \in CE(\mathfrak{g})$
	\begin{displaymath}
		\partial (p \wedge q) = 
		(\partial p) \wedge q 
		+ (-)^{|p|} p \wedge (\partial q)
		+ (-)^{|p|} [ p,q]
	~.
	\end{displaymath}
\end{reminder}

\begin{definition}[(Gerstenhaber) Algebra of multi-vector fields]
	We call \emph{algebra of multi-vector fields} the Gerstenhaber algebra associated to the Lie algebra of vector fields $\mathfrak{X}(M)$. 
	Explicitly, one has
	\begin{displaymath}
		\mathfrak{X}^k(M) = \Gamma(\Lambda T^n M) \cong \Lambda_{C^{\infty}(M)} \mathfrak{X}(M)
	\end{displaymath}
	together with $\partial$ and $[\cdot,\cdot]$ as defined in equations \eqref{eq:CE_boun} and \eqref{Eq:Schouten} respectively.
\end{definition}
As before, we call \emph{decomposable multi-vector field} any element $p\in \mathfrak{X}^{\bullet}(M)$ which can be expressed as a wedge product of a fixed number of tangent vector fields.

Hence, multi-vector calculus, is determined by the specification of the following three operations on the algebra  $\Omega(M)$ (\ie differential forms together with the wedge product):

	\begin{displaymath}
		\mathclap{
		\d \in \Der^1(\Omega(M)) ~,
		\quad
		\mathcal{L} : \mathfrak{X}^k(M) \to \Der^{1-k}(\Omega(M)) ~,
		\quad
		\iota: \mathfrak{X}^k(M) \to \Der^{-k}(\Omega(M)) ~;
		}
	\end{displaymath}
Operator $\d$ is not affected by this transition from vector fields to multi-vector fields, and is still given by the de Rham operator. 
One has only to make clear how to construct $\iota_p$ and $\mathcal{L}_p$ along a given multi-vector field $p$. 
It is possible to define such operators iterating the definition of the corresponding operator on ordinary vector fields.
\begin{definition}[(Multi-)contraction (interior product)]
	Given a differential form $\alpha \in \Omega(M)$, for any given decomposable vector field $p=v_1\wedge\dots\wedge v_n$ we define the \emph{multi-contraction} operation of $\alpha$ by $p$ as the consecutive contraction of $\alpha$ with each component of $p$, \ie
	\begin{displaymath}
		\iota_p\, \alpha = 
		\iota(v_1\wedge\dots\wedge v_n) \alpha = 
		\iota_{v_n}\dots \iota_{v_1} \alpha
		~.
	\end{displaymath}
	The corresponding interior product 
	$$\iota_\blank^k: \mathfrak{X}^k(M)\otimes \Omega^\ell(M) \to \Omega^{\ell-k}(M)$$ 
	is extended to a well-defined derivation $\iota_{v}\in \Der^{-k}(\Omega(M))$ by $C^{\infty}(M)$-linearity.
\end{definition}
\begin{remark}
	There is an obvious matter of choice in the definition of $\iota$ given by the order in which the components $v_i$ of $p$ are inserted into the form. 
	Clearly, all of this conventions only differs by a sign.
	Here, we are following the convention used, for example, in \cite{Rogers2010,Callies2016,Ryvkin2018}.
	
	Another natural choice would be to use the opposite order (see e.g \cite{Delgado2018,Delgado2018b}). The two conventions differs by a sign given by the following equation:
	\begin{displaymath}
		\iota_{x_1}\dots\iota_{x_n} = (-)^{\bar{\sigma}_n} \iota_{x_n}\dots \iota_{x_1}
	\end{displaymath}
		where $\bar{\sigma}_n$ is the permutation reversing the order of the list of $n$ elements, therefore
	\begin{displaymath}
		(-)^{\bar{\sigma}_n} =(-)^{\sum_{i=1}^n (n-i)} = (-)^{\frac{n (n-1)}{2}}
		~.
	\end{displaymath}
\end{remark}
There is a sign emerging from the previous convention that will appear recurrently when dealing with multisymplectic construction. We single out a compact definition and a couple of simple properties for later reference:
\begin{definition}[Total Koszul sign]\label{Def:SigmaSign}
	We call \emph{($k$-th) total Koszul sign}, the coefficient
	\begin{equation}\label{Eq:SigmaSign}
		\varsigma(k) = -(-)^{\frac{k(k+1)}{2}}
	\end{equation}
	corresponding to minus the Koszul sign of the permutation reversing $k+1$ degree $1$ elements
	\footnote{The first total Koszul signs read $+,+,-,-,+,+,-,-,\dots$ for $k$ equal to $1,2,\dots$.}.
\end{definition}
\begin{lemma}\label{lem:varsigmasignprops}
	The total Koszul signs satisfies the following properties
	\begin{equation}
		\varsigma(k+1)=-\varsigma(k-1) ~,\qquad \varsigma(k-1) = (-)^k \varsigma(k)
		~.
	\end{equation}
\end{lemma}

The notion of Lie derivative can be extended to the multi-vector fields context simply by enforcing the analogue of the Cartan's "magic" rule (equation \eqref{Eq:CartanMagic} in the above list).
\begin{definition}[(Multi-)Lie Derivative]
	Given a differential form $\alpha \in \Omega(M)$, for any given multi-vector field $p\in\mathfrak{X}(M)$ we define the \emph{multi-Lie derivative} of $\alpha$ along $p$ as the differential form given by the graded commutator
	\begin{displaymath}
		\mathcal{L}_p(\alpha) := \d \iota_p \alpha - (-)^{|p|} \iota_p \d \alpha
		~.
	\end{displaymath}
\end{definition}

These operators satisfy a "graded" analogue of the six commutation rules defining the ordinary Cartan calculus:
\begin{proposition}[Multi-Cartan commutation rules, {\cite[A.3]{Forger2003}}]
	For any multi-vector field $v \in \mathfrak{X}(M)$, the graded derivations $\mathcal{L}_v,\iota_v,\d$ in degrees $|\mathcal{L}_v| =1 -|v|,~ |\iota_v| = -|v|,~ |\d|=1$ satisfy the following commutation rules:
		\begin{eqnarray}
		[\d , \d] &=& 0
		\\
		{[\d, \mathcal{L}_v ]} &=& 0
		\\
		{[\d, \iota_x]} &=& \mathcal{L}_x \label{Eq:MultiCartanMagic}
		\\
		{[\iota_x , \iota_y]} &=& 0
		\\
		{[\iota_y,\mathcal{L}_x]} &=& -\iota_{[x,y]_{sc}}
		\\
		{[\mathcal{L}_x , \mathcal{L}_y]} &=& (-)^{(|x|-1)(|y|-1)}\mathcal{L}_{[x,y]_{sc}}
	\end{eqnarray}
	where $[\cdot,\cdot]_{sc}$ are understood as graded commutators pertaining to the (associative) composition of graded linear operators and $[\cdot,\cdot]$ denotes the Schouten bracket of multi-vector fields.
\end{proposition}
\begin{remark}\label{rem:Lieder}
From the previous commutation rules follows the following explicit expression of the Lie derivative along the wedge of two ordinary graded vector fields:
\begin{displaymath}
	\mathcal{L}_{x\wedge y} = [\d, \iota_{x\wedge y} ] = 
	[\d, \iota_y \iota_x] = (-)^{|y|}\iota_y [d, \iota_x] + [\d,\iota_y] \iota_x =
	\mathcal{L}_y \iota_x +(-)^{|y|}\iota_y \mathcal{L}_x
	~.
\end{displaymath}
Thence, by multiple iteration of the Cartan's formula \eqref{Eq:CartanMagic}, one gets the following explicit expression for the Lie derivative along a decomposable multi-vector field
\begin{displaymath}
	\mathcal{L}_{x_1 \wedge\dots\wedge x_n}
			=(-)^n \biggr[
		\iota(\partial(v_{1}\wedge\dots\wedge v_{n}))+\sum_{k=1}^{m} (-1)^{k} 
		\iota( v_{1} \wedge \dots
		\wedge \widehat{v}_{k} \wedge \dots \wedge {v}_{n})\mathcal{L}_{v_i}\biggr]
		~,
\end{displaymath}
where $\widehat{(\cdot)}$ denotes deletion and $\partial$ is the Chevalley-Eilenberg coboundary operator pertaining to the Lie algebra $\mathfrak{X}(M)$ (see equation \eqref{eq:CE_boun}).
\end{remark} 
This leads to the following recurring formula:
\begin{lemma}[Multi-Cartan magic formula {(see \cite[Lem 3.4]{Madsen2013} or \cite[Lem 2.18]{Ryvkin2016})}]\label{lemma:multicartan}
	Given any $m$ vector fields $x_i \in \mathfrak{X}(M)$, the following equation holds:
	\begin{equation}
	\begin{aligned}
		(-1)^m \d~ \iota(x_1\wedge\dots\wedge x_m) =&~\phantom{+} \iota(x_1\wedge\dots\wedge x_m) \textrm{d} ~+\\
		&+ \iota(\partial \, x_1\wedge\dots\wedge x_m) ~+\\
		&+ \sum_{k=1}^{m} (-1)^k  \iota( x_1\wedge\dots\wedge  \widehat{x_k}\, \wedge\dots\wedge x_m) \mathcal{L}_{x_k}
		~.
	\end{aligned}
	\end{equation}
\end{lemma}

\section{Multisymplectic manifolds}\label{Sec:MultiSymGeo}
In this work, we will adopt the notion of \emph{multisymplectic manifold} as introduced by Catrijn, Ibort, and de Le\'on  in \cite{CatIbort}
\begin{definition}[Multisymplectic manifold \cite{CatIbort} ($n$-plectic manifold)]\label{Def:MultisymplecticManifold}
	We call a \emph{pre-multisymplectic manifold} of degree $(n+1)$, a pair $(M,\omega)$ composed of a smooth manifold $M$ together with a closed differential $(n+1)$-form $\omega\in\Omega^{n+1}(M)$, called \emph{pre-multisymplectic form}.
	\\
	A pre-multisymplectic manifold is called \emph{multisymplectic} if $\omega$ is not degenerate, \ie if the flat map of $\omega$,
	\begin{equation}\label{Eq:OmegaFlat}
		\morphism{\omega^\flat}
		{TM}
		{\Lambda^{n}T^*M}
		{(x,u)}
		{(x,\iota_u \omega_x)}
		~,
\end{equation}		
	is an injective bundle morphism over $M$.
	\\
	For a fixed degree $(n+1)$ of the form $\omega$, such manifolds are also called \emph{"$n$-plectic"}. 
\end{definition}
Multisymplectic manifolds form a category with the following notion of morphism:
\begin{definition}[Multisymplecto-morphism ($(n)$-plecto-morphism)]
	Given two $n$-plectic manifolds $(M,\omega)$ and $(M',\omega')$ we call multisymplecto-morphism ($n$-plecto-morphism) $\phi: (M,\omega) \to (M',\omega')$ any diffeomorphism $\phi:M\to M'$ such that
	\begin{displaymath}
		\phi^* \omega' = \omega 
		~.
	\end{displaymath}
\end{definition}

\begin{remark}[On other notions of higher symplectic structures]
	We stress that there are numerous candidates proposed as the "higher analogue" of symplectic structures scattered across the literature; most of them are not entirely equivalent to the notion of multisymplectic structure adopted in this text.
	\\
	For instance, we mention  the notion of \emph{Poly-symplectic structures}\cite{Gunther1987a} which are given by vector-valued symplectic forms in $\Omega^{2}(M,V)$.
	Other recurrent names are the so-called \emph{$k$-symplectic structures}\cite{Awane1992} and the \emph{$k$-almost cotangent structures}\cite{DeLeon1988} that are now fully understood to be equivalent to the polysymplectic ones (see \cite{Blacker2019} for a review).
	\\
	There is also a notion of \emph{higher symplectic} geometry\cite{nlab:higher_symplectic_geometry} in the sense of categorification\cite{Baez2010} based on the idea to consider as the underlying space $M$ some generalized notion of smooth spaces, for instance, \emph{Lie $\infty$ algebroids}. 
\end{remark}

There are two extreme cases of multisymplectic manifolds given by the following two examples:
\begin{example}[Symplectic forms and volume forms]\label{Ex:VolumesAreMultiSymp}
	\quad
	\begin{itemize}
		\item Symplectic manifolds are $n$-plectic manifolds with $n=1$.
		\item Orientable manifolds of dimension $(n+1)$ are $n$-plectic with respect to a choice of volume form.
	\end{itemize}
	Observe that in both of these cases the flat bundle map $\omega^\flat$ is also bijective, by purely dimensional reasons, since:
	\begin{displaymath}
		\text{dim}(\Lambda^{n}T^\ast M ) = \binom{n+1}{n} = \binom{n+1}{1} = \text{dim}(\Lambda^{1}T^\ast M )
		~.
	\end{displaymath}		
\end{example}
In subsection \ref{Sec:MSExamples} we will recall some other examples thoroughly studied in the literature.

\subsection{Linear case}
	When considering a "flat" smooth manifold $M$, \ie globally diffeomorphic to the tangent space $T_p M$ at any point $p\in M$, giving a pre-$n$-plectic form on $M$ is tantamount to give a $n$-form on its linear model $V\cong T_p M$. 
	This justify the following definition:
	\begin{definition}[Multisymplectic vector space]
		We call \emph{$k$-plectic vector space} any $\R$-vector space $V$ endowed with a form $\omega\in \Lambda^{k+1} V^*$ that is not degenerate, \ie $ \omega(v,\cdot) = 0 \Leftrightarrow v=0$.
		\\
		Two multisymplectic $k$-plectic vector spaces $(V,\omega), (V',\omega')$ are said \emph{isomorphic} if there exists a linear isomorphism map $L: V \to V'$ pull-backing $\omega'$ to $\omega$ (\emph{linear multisymplectomorphism}), i.e such that $L^\ast \omega' = \omega$.	
	\end{definition}

\begin{remark}
	Clearly, the previous notion of multisymplectic manifold can be recovered from the notion of multisymplectic vector space.
	A \emph{multisymplectic bundle}, is a vector bundle $E\to M$ equipped with a multisymplectic (linear) structure of fixed order (say $k$) on each fibre with smooth dependence on the point of the base manifold.
	If $E$ is indeed the tangent bundle of $M$, the above structure is encoded by a differential form $\omega \in \Omega^k(M)$. If $\omega$ is non-degenerate in the above sense it will be called a \emph{almost multisymplectic structure}. 
	The special subclass of almost multisymplectic structures that are also \emph{integrable}, in the sense that  $\omega$ is closed, are the multisymplectic manifolds given in definition \ref{Def:MultisymplecticManifold}.	
\end{remark}

	\begin{example}["Canonical" linear multisymplectic form]\label{Ex:CanonicalLinearForm}
		Given any vector space $V$, for any $1\leq k \leq \text{dim}(V)$, the vector space $V\oplus \Lambda^k V^\ast$ is naturally $k$-plectic when equipped with the canonical multisymplectic form (see \cite[Prop. 2.2]{CatIbort}) 
		\begin{displaymath}
			\morphism{\Omega}
			{(V\oplus \Lambda^k V^\ast)^{\otimes (k+1)}}
			{V\oplus \Lambda^k V^\ast}
			{\pair{v_1}{\alpha_1}\otimes\dots\otimes\pair{v_{k+1}}{\alpha_{k+1}}}
			{\displaystyle\sum_{i=1}^{k+1}(-)^i \alpha_i(v_1,\dots,\widehat{v_i},\dots,v_{k+1})}
			~.
		\end{displaymath}
	\end{example}

	In the $1$-plectic case, the remarkable phenomenon occurs that every symplectic vector space of the same dimension "looks the same".
	Namely, every $1$-plectic vector space is isomorphic to the space $L\oplus L^*$, for a given Lagrangian subspace $L$, together with the canonical symplectic form defined in example \ref{Ex:CanonicalLinearForm}.
	In other words, the multisymplectic vector space of example \ref{Ex:CanonicalLinearForm}, when $k=1$, gives a canonical model for any $1$-plectic vector space in a given dimension.
	\\
	On the other hand, the situation is much wilder for arbitrary $k$-plectic vector spaces.
	On a given $n$ dimensional space, there may be several non-isomorphic families of $k$-plectic structures or even none. 
	The following theorem elucidates the picture:
	\begin{theorem}[Martinet-Capdevielle-Westwick-Djocovi\'c-Ryvkin {\cite[Thm. 4.2]{Ryvkin2018}}]\label{Thm:LeoClassification}
		Let $\Sigma_n^k$ be the number of non-equivalent isomorphism classes of $(k+1)$-plectic forms on $n$-dimensional vector spaces. 
		The following properties hold:
		\begin{itemize}
			\item $\Sigma_n^n=1$ for all $n$ (volume forms).
			\item $\Sigma^1_n$ as well as $\Sigma^{n{-}1}_n$ are zero for $n>1$.
			\item $\Sigma^2_n$ is 0 for $n$ odd and one for $n$ even (symplectic forms).
			\item $\Sigma_{n}^{n-2}=\lfloor \frac{n}{2}\rfloor-1$, when $(n \mod 4)\neq 2$ (for $n\geq 4$) and
			$\Sigma_{n}^{n-2}= \frac{n}{2}$, when $(n \mod 4)=2$ (for $n\geq 4$).
			\item $\Sigma_6^3=3$, $\Sigma_7^3=8$  ,  $\Sigma_8^3=21$,  $\Sigma_7^4=15$ and $\Sigma_8^5=31$ .
		 \item $\Sigma_n^k=\infty$ in all other cases.
		\end{itemize}
		\noindent For dimensions up to $n=9$ the numbers look as follows, where the rows range from 0-forms to $n$-forms:
		\begin{align*}
			\mathclap{
			\begin{array}{c|ccccccccccccccccccccc}
				n=0&&&&&&&&&&&-\\[1.5em]
				n=1&&&&&&&&&&-&&1\\[1.5em]
				n=2&&&&&&&&&-&&0&&1\\[1.5em]
				n=3&&&&&&&&-&&0&&0&&1\\[1.5em]
				n=4&&&&&&&-&&0&&1&&0&&1\\[1.5em]
				n=5&&&&&&-&&0&&0&&1&&0&&1\\[1.5em]
				n=6&&&&&-&&0&&1&&3&& 3&&0&&1\\[1.5em]
				n=7&&&&-&&0&&0&& 8&& 15&&2&&0&&1\\[1.5em]
				n=8&&&-&&0&&1&& 21&&\infty&&31&&3&&0&&1\\[1.5em]
				n=9&&-&&0&&0&&\infty&&\infty&&\infty&&\infty&&3&&0&&1\\[1.5em]
			\end{array}
			}
		\end{align*}
	\end{theorem}

\begin{remark}[Normal forms in multisymplectic geometry]
	The \emph{Darboux theorem} is a celebrated result in symplectic geometry stating that,
	for any given point $p$ of a symplectic manifold $(M,\omega)$ of dimension $2n$,
	there exists a local coordinate chart $(U_p,(x_1,\dots,x_{2n})$,
	called \emph{Darboux coordinates},
	such that $\omega\vert_{U_p}$ can be expressed in the normal form
	\begin{displaymath}
		\omega\eval_{U_p} = \sum_{i=1}^n \d x^i \wedge \d x^{n+i}
		~.
	\end{displaymath}
	The existence of such coordinate charts is a consequence of two phenomena which appear to be false in the general $k$-plectic case:
	\begin{enumerate}
		\item\label{Item:canonical} All $1$-plectic vector spaces with the same dimension are isomorphic to the canonical symplectic vector space given in example \ref{Ex:CanonicalLinearForm}.
		\item\label{Item:flatness} All symplectic manifolds are locally diffeomorphic to  the linear symplectic manifold $(T_p M,\omega_p)$, \ie for any $p\in M$ there exists a suitable neighbourhood $U_p$ together with a diffeomorphism $\phi:U_p \to T_p M$ such that $\phi(p)=0$ and $\phi^\ast \omega_p = \omega$, hence, in such a coordinate chart, $\omega\vert_{U_p}$ is given by the constant-coefficient differential form $\omega_p$.
	\end{enumerate}
	The general failure of \ref{Item:canonical} is completely understood by theorem \ref{Thm:LeoClassification}, \ie the existence of several non-equivalent isomorphism classes. 
	
	The failure of \ref{Item:flatness} is slightly more subtle and requires to define a notion of flatness at any point $p$, \ie the existence of a local chart $\varphi:U_p \to T_p M$ such that $\varphi(p)=0$ and $\varphi^\ast \omega_p = \omega$.
	Note that, if the latter holds, it is not guaranteed that local models for every point $p\in M$ will lie in the same isomorphism class.

	Other than $1$-plectic structure, also $(\text{dim}(M)-1)$-plectic structures (\ie volume forms) and multi-cotangent bundle structures admits a normal form. Several other examples are studied in \cite[\S 4]{Ryvkin2018}.
\end{remark}

\subsection{Examples}\label{Sec:MSExamples}
In addition to the extreme cases stated in example \ref{Ex:VolumesAreMultiSymp}, several other examples of multisymplectic structure can be found in the literature. 
Below we mention some examples mainly taken from \cite{CatIbort,Callies2016,Ryvkin2018}.

\begin{example}[Multicontangent bundles {\cite[\S 6]{CatIbort}}]\label{Ex:Multicotangent}
	Consider a smooth manifold $Q$, the corresponding \emph{Multicotangent bundle} $M = \Lambda^n T^\ast Q$ is naturally $n$-plectic.
	\\
	We denote by $\pi_Q$ the fibration $\Lambda^n T^\ast Q \twoheadrightarrow Q$ and by $\pi_M$ the fibration $\Lambda^n T^\ast M \twoheadrightarrow M$.
	On $M$ one can introduce a differential form $\theta\in\Omega^n(M)$, called \emph{tautological $n$-form} defined as the unique section of 
	$\Gamma(\Lambda^n T^\ast M, M)$ such that the diagram of smooth maps 
	\begin{center}
	\begin{tikzcd}
		M \equiv \Lambda^n T^\ast Q \ar[r,"\theta"] 
		&
		\Lambda^n T^\ast M \ar[d,two heads,"\pi_M"]
		\\
		Q \ar[r,"\alpha"] \ar[u,"\alpha"]
		&
		\Lambda^n T^\ast Q
	\end{tikzcd}
	\end{center}
	commutes for any section $\alpha \in \Gamma(\Lambda^n T^\ast Q, Q)$.
	In terms of pull-backs of differential forms, it means that $\alpha^\ast \theta = \alpha$.
	Explicitly, commutativity means that for any $q\in Q$, $\eta\in \Lambda^n T^\ast_q Q$ and vectors $u_1,\dots, u_2 \in T_{(q,\eta)}M$ the following equation holds
	\begin{displaymath}
		(q, \eta((\pi_Q)_\ast u_1, \dots, (\pi_Q)_\ast u_n ))
		=
		(q, \theta_{(q,\eta)}(u_1,\dots u_n) )
		~,
	\end{displaymath}
	where $(\pi_Q)_\ast$ denotes the push forward along the projection $T_{(q,\eta)}\pi_Q$.
	This can be read as the condition that
	\begin{displaymath}
		\theta_{\eta} (u_1,\dots,u_n) = \eta((\pi_Q)_\ast u_1, \dots, (\pi_Q)_\ast u_n )
	\end{displaymath}
	or	
	\begin{displaymath}
		\theta\eval_{(q,\eta)} = (\pi_Q)^\ast_{(q,\eta)} \eta
	\end{displaymath}
	justifying the name of \emph{"tautological form"}.
	Choosing local coordinates $(q^1,\dots,q^k)$ on $Q$ and denoting by $p_I$, with $I=(i_1,\dots,i_n)$ a multi-index with $1\leq i_1<i_2<\dots<i_n\leq k$, the canonical conjugated coordinates on the fibres, one gets that the tautological form takes the following expression:
	\begin{displaymath}
		\mathclap{
		\theta\eval_{(q^i,p_I)} 
		= 
		\sum_{1\leq i_1<i_2<\dots<i_n\leq k}
		p_{i_1,\dots,i_n} \d q^{i_1}\wedge\dots \d q^{i_n}
		= \sum_{I} p_I \d q^I	
		~.
		}
	\end{displaymath}		
	\\
	The canonical $n$-plectic form on $\Lambda^n T^\ast Q$, also known as \emph{multicanonical}\cite{Forger2005}, is given by the exterior derivative of $\theta$:
	\begin{displaymath}
		\omega := -\d \theta ~.
	\end{displaymath}
	In the above coordinates it reads as
	\begin{displaymath}
		\omega\eval_{(q^i,p_I)} 
		= \sum_{I} -\d p_I \wedge \d q^I 
		~.
	\end{displaymath}		
	Non-degeneracy of $\omega$ can be read as the condition that, for any $v= \sum_i x^i \partial_{i} + \sum_I z_I \partial_{p^I}$, the form
	\begin{displaymath}
		\iota_v \,\omega =
		\mkern-15mu
		\sum_{\substack{I=(i_1,\dots,i_n)\\1\leq i_1<\dots<i_n\leq k}}
		\mkern-15mu
		\left(
			z_I \d q^I
			+
			\d p^I \wedge \left(
				\mkern-50mu\sum_{\mkern70mu j\in\lbrace i_1,\dots,i_n\rbrace}
				\mkern-50mu 
				x^j\, \d q^{i_1}\wedge\dots \wedge \widehat{\d q^j} \wedge \dots \wedge \d q^{i_n}
			\right)
		\right)
	\end{displaymath}
	vanishes only if $v=0$. The latter condition is met since the right-hand side is a composition of linearly independent $n$-forms. 

		This construction is the ``higher analogue'' of the canonical symplectic structure naturally defined on any cotangent bundle. 
		Note, however, that this is not yet the ``higher analogue'' of a \emph{phase space}, as explained in the following example.		
\end{example}

\begin{example}[Multiphase space {\cite{Carinena1991b,Gimmsy1,CatIbort}}]
%
%
	Consider now the case that $Q$ is fibred over another smooth manifold $\Sigma$, \ie is the total space of a smooth bundle $\pi:Q\to \Sigma$.
	Denote by $V(\pi)$ the vertical sub-bundle $V(\pi)\hookrightarrow TQ$.
	Recall that the fibres of $V(\pi)$ are given by $V_y(\pi) = ker(\pi_\ast)_y$ for any $y \in Q$.
	Consider then the subbundle of $r$-semibasic $n$-forms of $Q$, concretely defined by the following fibres
	\begin{displaymath}
		\Lambda^n_r Q =
		\left\lbrace
			\alpha \in \Lambda^n T^\ast Q ~
			\vert \quad
			\iota_{v_{1}}\dots \iota_{v_r} \alpha = 0,\quad \forall v_i \in V_{\pi(\alpha)}(\pi)
		\right\rbrace
		~.
	\end{displaymath}	
	For any $r\leq n$ one has the following diagram in the category of smooth manifolds
		\begin{displaymath}
			\begin{tikzcd}[column sep= small,row sep=small]
				\Lambda^n_r Q \ar[rr,hookrightarrow,"i"] \ar[dr]& & \Lambda^n T^\ast Q \ar[dl]\\
				& Q \ar[d,"\; \pi"] & \\
				& \Sigma &
				~.
			\end{tikzcd}			
		\end{displaymath}
	One can then pullback the canonical multisymplectic form $\omega$ given in example \ref{Ex:Multicotangent} from $\Lambda^n T^\ast Q$ to $\Lambda^n_r Q$ to obtain the pre-multisymplectic manifold
	\begin{displaymath}
		(\Lambda^n_{r} Q,~ i^\ast\omega)
	\end{displaymath}
	that can be proved to be $n$-plectic (see for example \cite[Prop./Def. 2.14]{Ryvkin2018}).
	
 A crucial result is that, in the particular case $r=2$, the smooth manifold $\Lambda^n_2 Q$ is diffeomorphic to the twisted affine dual of the first jet bundle \cite{Saunders1989} of $\pi:Q\to\Sigma$ (see \cite[Prop. 2.1]{Gimmsy1} or  \cite[Example 2.4]{Baez2010}).
	The choice of such naming comes from the geometric mechanics of classical field theories with \emph{configuration bundle} given by $\pi$, for details see \cite{Carinena1991b,Gimmsy1,Ryvkin2018}
\end{example}

\begin{example}[Semisimple Lie groups \cite{CatIbort}]
	Any Lie group $G$ with (finite-dimensional) Lie algebra $\mathfrak{g}$ is canonically pre-$2$-plectic, if $\mathfrak{g}$ is semisimple the canonically pre-$2$-plectic form is non-degenerate.
	\\	
	The \emph{canonical $3$-form} for $G$ is the form $\omega \in \Omega^3(M)$ given, 	for any $g\in G$ and $u_i\in T_g G$, by
	\begin{displaymath}
		\omega\eval_g (u_1,u_2,u_3) = \left\langle \theta^L\eval_g(u_1),\left[\theta^L\eval_g(u_2),\theta^L\eval_g(u_3)\right]\right\rangle	
	\end{displaymath}			
	where:
	\begin{itemize}
		\item $\pairing$ denotes the \emph{Killing form} of $\mathfrak{g}$ defined as the bilinear operator
	\begin{displaymath}
		\morphism{\langle\cdot,\cdot\rangle}
		{\mathfrak{g}\otimes \mathfrak{g}}
		{\R}
		{(x,y)}
		{\text{trace}(\text{ad}_x \cdot \text{ad}_y )}
		~.
	\end{displaymath}
	From the properties of the trace and from the $Ad$-equivariance of the Lie bracket one has that $\pairing$ is symmetric and $Ad$-invariant.
		\item $ad$ denotes the adjoint representation of $\mathfrak{g}$ explicitly given by the Lie algebra morphism
	\begin{displaymath}
		\morphism{\text{ad}}
		{\mathfrak{g}}
		{\End(\mathfrak{g})}
		{x}
		{(\text{ad}_x :~y\mapsto [x,y])}
	\end{displaymath}
	(where $[\cdot,\cdot]$ on $\End(\mathfrak{g})$ is given by the commutator of linear operators);
		\item $\theta^L$ denotes the \emph{Maurer-Cartan} form of $G$, that is the left invariant $\mathfrak{g}$-valued $1$-form $\theta^L\in \Omega^1(M,\mathfrak{g})$ defined as the tangent map along the left multiplication
		\begin{displaymath}
			\morphism{L_g}
			{G}
			{G}
			{h}
			{g\cdot h}
			~.
		\end{displaymath}		 
		Namely:
		\begin{displaymath}
			\theta^L\eval_g = T_g(L_{g^{-1}}): T_g G \to T_e G \cong \mathfrak{g}
			~.
		\end{displaymath}
	\end{itemize}
	The closedness of $\omega$ follows from being both left and right invariant on the Lie group $G$.
	If the group is semisimple then $\mathfrak{g}=[\mathfrak{g},\mathfrak{g}]$; therefore the pairing $\pairing$ and $\omega$ are non-degenerate (see \cite[Ex. 3.6]{Ryvkin2018} for the complete argument).
\end{example}

\begin{example}[Cosymplectic manifolds \cite{CatIbort}]
	A \emph{cosymplectic} manifold is a triple $(M,\Phi,\eta)$ consisting of an orientable $(2n+1)$-dimensional manifold $M$ together with closed forms $\Phi\in\Omega^1(M)$ and $\eta\in \Omega^2(M)$ such that $\Phi\wedge\eta^{n}$ is a volume form.
	The form $\omega= \Phi\wedge\eta^n$ is a closed non-degenerate $(2n+1)$-form, hence $2n$-plectic over $M$. 
	Non-degeneracy can be ascertained easily by writing $(M,\phi,\eta)$ in Darboux-like coordinates (see \cite{Cappelletti2013} for a complete survey on the subject).
\end{example}

\begin{example}[(Subclasses) of K\"ahler manifolds]
	Any \emph{quaternionic almost K\"ahler manifold} is $3$-plectic with multisymplectic form given by the \emph{fundamental $4$-form} $\Omega$  \cite{Swann1991}. 
	Recall that \emph{Hyper-K\"ahler} manifolds are a subcase of the latter.
	In \cite{Madsen2012} is proposed a construction recovering all homogeneous strictly nearly K\"ahler $6$-manifolds as $2$-plectic.
\end{example}

\begin{example}[Sum and products of multisymplectic manifolds {\cite{Ryvkin2018}}]\label{Ex:SumProductMultiSymp}
	Given a finite number of multisymplectic manifolds $(M_i;\omega_i)$, one can canonically endow the product manifold $\times_i M_i$ wit a multisymplectic structure given by the wedge product
	\begin{displaymath}
		\omega_{prod} = \bigwedge_{i} \pi_i^\ast \omega_i ~,
	\end{displaymath}
	where $\pi_j : \times_i M_i \to M $ is the standard projection.
	If all the multiplicand manifolds are $k$-plectic, one can also define a $k$-plectic "sum" structure given by
	\begin{displaymath}
		\omega_{sum}=\sum_{i} \pi_i^\ast \omega_i
		~.
	\end{displaymath}
\end{example}	

\begin{example}[$G_2$-structures {\cite[ex 3.7]{Ryvkin2018}}]
	Given a seven-dimensional manifold $M$, a closed $G_2$-structure on $M$ consists of a closed differential 3-form $\omega$ admitting, for all $p\in M$,  a basis $(e^1,...,e^7)$ of $T^*_pM$ such that 
	\begin{displaymath}
	\omega_p=  e^{123}+e^{145}-e^{167}+e^{246}+e^{257}+e^{347}-e^{356} ~,
	\end{displaymath}
	where $e^{ijk}$ denotes $e^i\wedge e^j\wedge e^k$. 
	Being the addenda defining $\omega$ linear independent, $\omega^\flat$ has a $7$-dimensional range, hence $\omega$ is a $2$-plectic structure.
\end{example}

The following is an example involving a $\infty$-dimensional smooth manifold:
\begin{example}[Affine connections on a Principal bundle {\cite[\S 10]{Callies2016}}]
	Let be $G$ a Lie group with a finite-dimensional Lie algebra $\mathfrak{g}$,
	consider a smooth $G$-principal bundle 
	$P=(G\hookrightarrow P \twoheadrightarrow M)$ over a $(n+1)$-dimensional compact smooth manifold $M$.
	The action $G\action P$ on the principal bundle is considered on the right (\cf subsection \ref{subsec:prequantizationReminder}).
	\\
	Denote by $\mathcal{A}$ the space of all \emph{principal $G$-connections} on $P$. 
	Recall that $\mathcal{A}$ is an infinite-dimensional affine space modelled on the vector space 
	$\vec{\mathcal{A}}= (\Omega^1_{hor}(P)\otimes \mathfrak{g})^G$, where superscript $G$ denotes invariance with respect to the aforementioned right action.
	Therefore, $\mathcal{A}$ is smooth and "flat" in the sense that $T_a \mathcal{A}\cong \vec{A}$ for any affine connection $a\in \mathcal{A}$ (See \cite{Kobayashi1996}).
	\\
	Any $G$-invariant polynomial $q \in S^{n+1}(\mathfrak{g}^*)$ define a constant $(n+1)$-differential form, hence pre-$n$-plectic form, on the space of affine connections $\mathcal{A}$.
	If the chosen polynomial is in particular also non-degenerate, \ie the  mapping
	\begin{displaymath}
		\lambdamorphism{\mathfrak{g}}
		{S^n(\mathfrak{g}^*)}
		{x}
		{\iota_x q}
	\end{displaymath}
	in injective, the corresponding canonical pre-$n$-plectic structure on $\mathcal{A}$ is in particular non-degenerate \cite[Prop. 10.3]{Callies2016}.
\end{example}

\subsection{Special classes of differential forms and vector fields} 
Exactly as it happens in symplectic geometry, fixing a (pre-)$n$-plectic structure $\omega$ on $M$ 
provides a criterion for identifying special classes of vector fields and differential forms. 

	\begin{definition}[Multisymplectic and Hamiltonian vector fields]\label{def:Hamiltonianvfields}
		Given a pre-$n$-plectic manifold $(M,\omega)$, a vector field $v\in \mathfrak{X}(M)$ is said \emph{(multi)-symplectic} (resp. \emph{Hamiltonian}) if $\iota_v \omega$ is a closed (resp. exact) form.
		We denote the corresponding vector spaces as
		\begin{displaymath}
			\begin{aligned}
				\mathfrak{X}_{\msy}(M,\omega) =&~
				\lbrace v \in \mathfrak{X}(M) ~|~ \iota_v(\omega) \in Z^n(M) \rbrace
				\\
				\mathfrak{X}_{\ham}(M,\omega) =&~
				\lbrace v \in \mathfrak{X}(M) ~|~ \iota_v(\omega) \in B^n(M) \rbrace
				~.
			\end{aligned}
		\end{displaymath}			
	\end{definition}
	It follows from the Cartan rule that the flow of a multisymplectic vector field preserves the pre-$n$-plectic form (strictly, in the sense of definition \ref{Def:conservedQuantities}). 
	Hence the elements of $\mathfrak{X}_{\msy}(M,\omega)$ are the infinitesimal counterpart of diffeomorphisms preserving the multisymplectic structure.

	\begin{lemma}
		Let be $(M,\omega)$ a pre-$n$-plectic manifold and consider the vector spaces introduced in definition \ref{def:Hamiltonianvfields}.
		One has the following inclusion of Lie algebras
		\begin{displaymath}
			\begin{tikzcd}
				\mathfrak{X}_{\ham}(M,\omega) \ar[r,hook] &
				\mathfrak{X}_{\msy}(M,\omega) \ar[r,hook] &
				(\mathfrak{X}(M),[\cdot,\cdot])
			\end{tikzcd}
			~.
		\end{displaymath}
		In the case that $H^{n}_{d R}(M)=0$ the first inclusion is an isomorphism $\mathfrak{X}_{\msy}(M)\cong\mathfrak{X}_{\ham}(M)$.	
		\\
		Furthermore, $\mathfrak{X}_{\ham}(M)$ is a Lie algebra ideal of $\mathfrak{X}_{\msy}(M)$, \ie $[\mathfrak{X}_{\ham}(M),\mathfrak{X}_{\msy}(M)]\subset \mathfrak{X}_{\ham}(M)$.	
	\end{lemma}	
	 Accordingly, we can select a special class of differential forms:
\begin{definition}[Hamiltonian $(n-1)$-forms]\label{Def:Hamiltonianform}
	Given a pre-$n$-plectic manifold $(M,\omega)$, a differential form $\alpha \in \Omega^{n-1}(M)$ is called \emph{Hamiltonian} if is a primitive of $\iota_{\vHam_\alpha} \omega \in B^{n}(M)$ for a certain Hamiltonian vector field $\vHam_\alpha$, called the \emph{Hamiltonian vector field of $\alpha$}.
	\\
	Hamiltonian $(n-1)$-forms constitute a subspace of $\Omega^{n-1}(M)$ denoted as
	\begin{displaymath}
\Omega^{n{-}1}_{Ham}(M,\omega)=\left\{\alpha\in \Omega^{n-1}(M)~|~\exists\, \vHam_\alpha\in\mathfrak X(M) 
			~:~ d\alpha=-\iota_{\vHam_\alpha}\omega \right\}
			~.	
	\end{displaymath}
	The equation $d\alpha +\iota_{\vHam_\alpha}\omega=0$ it is also known as the \emph{Hamilton-DeDonder-Weyl (HDDW) equation} (see \cite{Ryvkin2018}).
\end{definition}
\begin{remark}[Sign conventions]
	The appearance of the minus sign in the HDDW equation is purely conventional.
	Here we are adopting a convention mostly found in the multisymplectic geometry literature (\eg \cite{Callies2016,Ryvkin2018}).
	The opposite choice can be found for instance in \cite{CannasdaSilva2001}.
\end{remark}

\begin{remark}\label{Rem:ClosedformsTrivialHamiltonian}
	Observe that, except for the case $n$ equal to $1$ or $\dim(M)-1$ where $\omega^\flat$ is bijective, not all the $(n-1)$-forms are Hamiltonian.
	\\
	Note also that if $\omega$ is non-degenerate, \ie $n$-plectic rather than pre-$n$-plectic, if there exists an Hamiltonian vector field $v_\alpha$ pertaining to $\alpha$ it ought to be unique. Hence there is a well-defined surjective map
	\begin{displaymath}
		\epimorphism{\pi_{\ham}}
		{\Omega^{n-1}_{\ham}(M,\omega)}
		{\mathfrak{X}_{\ham}(M,\omega)}
		{\alpha}
		{\vHam_{\alpha}}
		~.
	\end{displaymath}
	Finally, observe that closed $(n-1)$-forms are clearly Hamiltonian with trivial Hamiltonian vector field, hence there is a short exact sequence of vector spaces:
	\begin{displaymath}
		\begin{tikzcd}[column sep = small]
			0 \ar[r]&
			\Omega^{n-1}_{cl}(M) \ar[r,hook] &[1em]
			\Omega_{\ham}^{n-1}(M,\omega) \ar[r,two heads,"\pi_{\ham}"] &[1.5em]
			\mathfrak{X}_{\ham}(M,\omega) \ar[r] & 
			0
		\end{tikzcd}
		~.
	\end{displaymath}	 
\end{remark}

\begin{remark}[Hamiltonian pairs and Hamilton-DeDonder-Weyl equation]\label{Rem:HamiltonianPairs}
	Slightly more in general (see for instance \cite{Delgado2018b,Ryvkin2018}), given a pre-$n$-plectic manifold $(M,\omega)$, one could define a \emph{Hamiltonian pair} as the pair of multivector fields and differential forms satisfying the so-called \emph{Hamilton-DeDonder-Weyl} equation pertaining to $\omega$.
	Namely, one can define
	\begin{displaymath}
		Ham^{n-k}(M) = \left\lbrace \left.
			\pair{\vHam_\alpha}{\alpha}\in \mathfrak{X}^k(M)\oplus\Omega^{n-k}(M)
			~\right\vert~ \iota_{\vHam_\alpha}\omega + \d \alpha = 0 
		\right\rbrace
		~,
	\end{displaymath}
	for any $ 0\leq k\leq n-1$.
	Note that one can make sense of the space of Hamiltonian $(n-k)$-forms as a pullback in the category of ordinary vector spaces
	\begin{displaymath}
		\begin{tikzcd}
			Ham^{n-k}(M) \arrow[dr, phantom, "\scalebox{1.5}{$\lrcorner$}" , very near start, color=black] \ar[r,dashed]\ar[d,dashed]
			& 
			\Omega^{n-k}(M) \ar[d,"\d"]\\
			\mathfrak{X}^k(M) \ar[r,"\iota(\blank)\omega"] & 
			\Omega^{n-k+1}(M)
		\end{tikzcd}
		~.
	\end{displaymath}
	In the non-degenerate case, one has that $Ham^{n-1}(M)\cong \Omega_{\ham}^{n-1}(M)$.
\end{remark}

\begin{notation}[On the ubiquitous role of the multicontraction operator]\label{Rem:SignedMultiContraction}
	Once one fixes a "preferred" differential form $\omega$ over a given smooth manifold $M$, and being the former a $C^\infty(M)$-multilinear map on vector fields, it is natural to study the properties of the pair $(M,\omega)$ by probing the chosen form $\omega$ on arbitrary vector fields.
	Hence, the contraction (insertion) operation against the multisymplectic form will be a key ingredient in most constructions involved in multisymplectic geometry.
\\	
In this spirit, let us anticipate some notations involving multicontraction that will have a recurring role in the following, although they may appear trivial and redundant at this stage.
	\\
	\smallskip
	Given a differential form $\omega \in \Omega^n(M)$ we define the \emph{signed $k$-multicontraction of $\omega$} as the following linear operator
	\begin{equation}\label{Eq:signedMulticontraction}
		\morphism{\iota_{\mathfrak{X}}^k \omega}
		{\mathfrak{X}^k(M)}
		{\Omega^{n-k}(M)}
		{x_1\wedge\dots\wedge x_k}
		{\varsigma(k)\iota(x_1\wedge\dots\wedge x_k) \omega}
		~.
	\end{equation}
	Observe that this maps can be arranged as a graded morphism between two chain complexes
	\begin{equation}\label{Eq:multicontractionaschainmap}
	\hspace{-.025\textwidth}
		\begin{tikzcd}[column sep = small]
			0 &
			\mathfrak{X}(M) \ar[l,"\partial"']\ar[d,"\iota_{\mathfrak{X}}^1\omega"]&
			\mathfrak{X}^2(M) \ar[l,"\partial"']\ar[d,"\iota_{\mathfrak{X}}^2\omega"] &
			\cdots \ar[l,"\partial"']&
			\mathfrak{X}^k(M) \ar[l,"\partial"']\ar[d,"\iota_{\mathfrak{X}}^k\omega"]&
			\cdots \ar[l,"\partial"']&
			\mathfrak{X}^{n+1}(M) \ar[l,"\partial"']\ar[d,"\iota_{\mathfrak{X}}^{n+1}\omega"]&
			\cdots \ar[l,"\partial"']
			\\
			\cdots &
			\Omega^n(M) \ar[l,"\d"] &
			\Omega^{n-1}(M) \ar[l,"\d"] &
			\cdots \ar[l,"\d"]&
			\Omega^{n+1-k}(M) \ar[l,"\d"]&
			\cdots \ar[l,"\d"]&
			\Omega^{0}(M) \ar[l]&
			0 \ar[l,"\d"]
		\end{tikzcd}
	\end{equation}	
	namely, following \cite[\S 3.2]{Delgado2018b}, one can observe that $\iota^k_{\mathfrak{X}}\omega$ can be read as the $k$-th component of a homogeneous map of degree $n$ given by the following graded vector space morphism
	\begin{displaymath}
		\iota_{\mathfrak{X}}\omega ~: \setminus \mathfrak{X}(M) \to \Omega(M)[n+1]
		~,
	\end{displaymath}
where $\setminus \mathfrak{X}(M)$ denotes the graded vector space $\mathfrak{X}(M)$ taken with reverse ordering(\ie $(\setminus \mathfrak{X}(M))^{-k}=\mathfrak{X}(M)^k$, see equation \eqref{eq:degreeReversingFunctor}). 
This graded map yields a chain-complexes morphism if certain other conditions are met (see lemma \ref{lem:signedmulticontractionLinfinitymorphism}).
	\\
	\smallskip
	Collectively, this can be abstracted as the image of the following operator
	\begin{displaymath}
		\morphism{\iota_\mathfrak{X}}
		{\Omega(M)}
		{\underline{\Hom}(\setminus \mathfrak{X}(M),\Omega(M))}
		{\omega}
		{\left\lbrace
			\iota^{-k}_\mathfrak{X}\omega: \mathfrak{X}(M)^{-k} \to \Omega^{|\omega|+k}(M)
		\right\rbrace_{k\leq -1}}
	\end{displaymath}
	Starting from section \ref{Sec:HCMM}, we will mainly focus on certain graded linear subspaces of the space of multivectors fields. 
	More generally, for any given linear map $\vAct:\mathfrak{g}\to \mathfrak{X}^1(M)$ we will denote the signed multi-contraction along the image of $\vAct$ as the precomposition 
	\begin{equation}\label{eq:restrictedmulticontraction}
		\iota_{\mathfrak{g}}^k \omega := \left(\iota_{\mathfrak{X}(M)}^k \omega \right)
		\circ
		\vAct^{\otimes k}
		~.
	\end{equation}		
	In chapter \ref{Chap:MarcoPaper} we will make frequent use of a (skew)-symmetrized version of the ordinary contraction operator defined, for any $n\in \mathbb{Z}$, as
	\begin{displaymath}
		\morphism{\pairing_\pm}
		{\left(\mathfrak{X}(M)\oplus \Omega(M){[n]}\right)^{\otimes 2}}
		{\Omega(M){[n-1]}}
		{\pair{x_1}{\alpha_1}\otimes \pair{x_2}{\alpha_2}}
		{\frac{1}{2}\left(\iota_{x_1}\alpha_2 \pm \iota_{x_2}\alpha_1\right)}
		~.
	\end{displaymath}
\end{notation}
	
A first justification for introducing the notation given in remark \ref{Rem:SignedMultiContraction} comes from the next lemma:
\begin{lemma}[{KKS $L_\infty$-cocycle \cite[Prop. 3.8]{Fiorenza2014a}}]\label{lem:signedmulticontractionLinfinitymorphism}
	Call $\mathcal{B}_{(n)}$ the following shifted truncation of the de Rham complex of the pre-$n$-plectic manifold $(M,\omega)$:
	\begin{displaymath}
		\mathcal{B}_{(n)} := \trunc_n(\Omega(M))[n]\oplus \d \Omega^{n-1}(M)
		~;
	\end{displaymath}
	the latter is given component-wise by the following diagram in the category of vector spaces:
	\begin{displaymath}
		\begin{tikzcd}[column sep = small, row sep = tiny]
			0 \ar[r] & 
			\mathcal{B}_{(n)}^{-n}\ar[d,equal] \ar[r,"\d"] &\dots \ar[r,"\d"] &
			\mathcal{B}_{(n)}^{-k} \ar[r,"\d"] \ar[d,equal]& \dots &			
			 \mathcal{B}_{(n)}^{-1}\ar[d,equal]  \ar[r,"\d"] & 
			\mathcal{B}_{(n)}^{0}\ar[d,equal]  \ar[r] & 0
			 \\
			 & C^\infty(M) && \Omega^{n-k}(M) & & \Omega^{n-1}(M) & \d \Omega^{n-1}(M)
		\end{tikzcd}
		~.
	\end{displaymath}		
	The $k$-multilinear maps $\iota_{\mathfrak{X}_{\ham}}^k\omega$ defined in equation \eqref{Eq:signedMulticontraction}, with $1\leq k \leq n+1$, taken together with the restriction to Hamiltonian vector fields (see eq. \eqref{eq:restrictedmulticontraction}),
	define the components $(1\leq k \leq n+1)$ of a $L_\infty$-morphism
	\begin{displaymath}
		(\iota_{\mathfrak{X}_{\ham}}\omega): \mathfrak{X}^1_{\ham}(M) \to \mathcal{B}_{(n)}
	\end{displaymath}
	between the Lie algebra of Hamiltonian vector fields and the Abelian $L_\infty$-algebra $\mathcal{B}_{(n)}$.
\end{lemma}
\begin{proof}
	Observe first that the diagram \eqref{Eq:multicontractionaschainmap} can be read as a bona fide chain map when restricted to multisymplectic vector fields since lemma \ref{lemma:multicartan} guarantees,
	for any $p=v_1\wedge\dots \wedge v_k \in \Lambda^k \mathfrak{X}(M)$,
	that
	\begin{displaymath}
			\left( \d \circ (\iota^k_{\mathfrak{X}}\omega) - (\iota_{\mathfrak{X}}^{k-1}\omega ) \circ \partial
			\right)(p) =
			(-)^k\varsigma(k) \left(\iota_p \cancel{\d \omega}
				+\sum_{i=1}^k \iota_{x_k}\dots \widehat{\iota_{x_i}}\dots \iota_{x_1} \mathcal{L}_{x_i}\omega			
			 \right)
	\end{displaymath}	 
	and the right-hand side vanishes for $p\in\Lambda^k \mathfrak{X}_{\msy}(M,\omega)$.
	Hence $\iota_{\mathfrak{X}}\omega$ defines a chain map $\iota_{\mathfrak{X}}\omega ~: \setminus \mathfrak{X}(M) \to \Omega(M)[n+1]$ from the order reversed complex of multisymplectic multivector fields to the $[n+1]$ shift of the de Rham complex (which puts $(n+1)$-forms in degree zero).
	\\
	Restricting furthermore to $\mathfrak{X}_{\ham}(M,\omega)\subset \mathfrak{X}_{\msy}(M,\omega)$ this can be refined to a chain map
	\begin{displaymath}
		\iota_{\mathfrak{X}_{\ham}} : \setminus \mathfrak{X}_{\ham}(M,\omega) \to \mathcal{B}_{(n)}[1]
		~.		
	\end{displaymath}
	An overall shift of the previous map yields the following chain complex morphism
	\begin{displaymath}
		\iota_{\mathfrak{X}_{\ham}}[1] : \left(S^{\geq 1}\left(\mathfrak{X}^1_{\ham}(M,\omega)[1]\right)\right)[-1] \to \mathcal{B}_{(n)}
	\end{displaymath}
	giving, according to remark \ref{Rem:LftyMorphasChainMap}, the sought $L_\infty$-morphism.
\end{proof}
\begin{remark}
	This construction has been proposed in \cite{Fiorenza2014a} as a $n$-plectic analogue of the Kirillov-Konstant-Souriau $2$-cocycle of symplectic mechanics.
	In \cite[Rem 6.2]{Callies2016} the chain complex $\mathcal{B}_{(n)}$ is denoted as $B^n A$ to remind of an analogue construction involved when studying classifying spaces on smooth bundles.
\end{remark}

\paragraph{Conserved quantities}
We will often adopt the following nomenclature borrowed from \cite{Ryvkin2016} providing different nuances of conservation along a given vector field.

\begin{definition}[Conserved quantities {\cite{Ryvkin2016}}]\label{Def:conservedQuantities}
	Let be $v\in \mathfrak{X}^1(M)$ a vector field over the smooth manifold $M$.
	A differential form $\alpha\in \Omega^k(M)$ is called:
	\begin{itemize}
		\item \emph{locally conserved} along $v$ if $\mathcal{L_v \alpha}$ is a closed form;
		\item \emph{globally conserved} along $v$ if $\mathcal{L_v \alpha}$ is an exact form;
		\item \emph{strictly conserved} along $v$ if $\mathcal{L_v \alpha}=0$.
	\end{itemize}
		We denote by $C_{loc}(v),C(v),C_{str}(v)$ the graded vector space of locally, globally, strictly preserved forms respectively.
\end{definition}
This definition can be easily translated to group actions implying the conservation along the fundamental vector fields.
The following inclusion relations between conserved quantities follow from the Cartan formula:
\begin{lemma}[{\cite[\S 1.3]{Ryvkin2016}}]
	Given a vector field $v\in \mathfrak{X}(M)$, the following diagram holds in the category of graded vector spaces
	\begin{displaymath}
		\begin{tikzcd}
			\d C_{loc}(v) \ar[r,hook] & C_{str}(v) \ar[r,hook] & C(v) \ar[r,hook] & C_{loc}(v) \ar[r,hook] & \Omega(M)
			\\
			& & \Omega_{cl}(M) \ar[u,hook]
		\end{tikzcd}
		~.
	\end{displaymath}	
	Furthermore, one has that $C_{str}(v)$ is a graded subalgebra of $(\Omega(M),\wedge)$ and $C(v),C_{loc}(v)$ are graded modules over $(C_{str}(v) \cap \Omega_{cl}(M))$.
\end{lemma}
In particular, one has that any Hamiltonian $(n-1)$-form is globally conserved along its corresponding Hamiltonian vector field:
\begin{equation}
	\mathcal{L}_{v_H} H = \d \iota_{v_H} H + \cancel{\iota_{v_H}\d H}
	\qquad \forall H \in \Omega^{n-1}_{\ham}(M,\omega)
\end{equation}
(or strictly conserved in the $1$-plectic case).

\section{Higher observables}\label{Section:RogersObservables}
There is an algebraic counterpart to the geometric theory of smooth manifolds represented by the study of the unital commutative algebra of smooth functions on a given manifold, also known as \emph{smooth observables} \cite{Nestruev2010}.
One could informally think about this relationship as a suitable extension of the Gelfand duality \cite{nlab:gelfand_duality}, establishing an equivalence between the category of compact topological spaces and the category of commutative $C^\ast$-algebras, from ordinary topology to differential geometry.

When the studied manifold is in particular symplectic, its dual -algebraic- object is given by its corresponding \emph{Poisson algebra}. 
This is the unital commutative algebra of smooth functions over $M$ together with a Lie bracket $\{\cdot,\cdot\}$, called the Poisson bracket, that is compatible with the pointwise product in the sense of being a derivation in each entry, \ie
	\begin{displaymath}
		\lbrace f, g\cdot h \rbrace = g\cdot\lbrace f,  h \rbrace + \lbrace f, g \rbrace\cdot h
		\qquad \forall f,g,h \in C^{\infty}(M)
		~.
	\end{displaymath}
Explicitly, the Poisson bracket corresponding to $(M,\omega)$ takes the following expression:
\begin{displaymath}
	\morphism{\lbrace\cdot,\cdot\rbrace}
	{C^{\infty}(M)\otimes C^{\infty}(M)}
	{C^{\infty}(M)}
	{(\alpha,\beta)}
	{\iota_{\vHam_\beta}\iota_{\vHam_\alpha}\omega},
\end{displaymath}
	where $\vHam_\alpha$ denotes the Hamiltonian vector field pertaining to $\alpha$.
	
Note that there is a neat interpretation of these concepts when one reads a symplectic manifold in the geometric mechanics setting, \ie as the space of states of a certain mechanical system with finite degrees of freedom.
Namely, one can see the Poisson algebra of $M$ as the set of physical observables of the system, also known as "classical observables".
The latter has to be interpreted in the sense of measurable quantities yielding a single value (real number) on every state.
	\begin{figure}[h!]
	 	 \begin{center}
	 	   \includegraphics[width=0.22\textwidth]{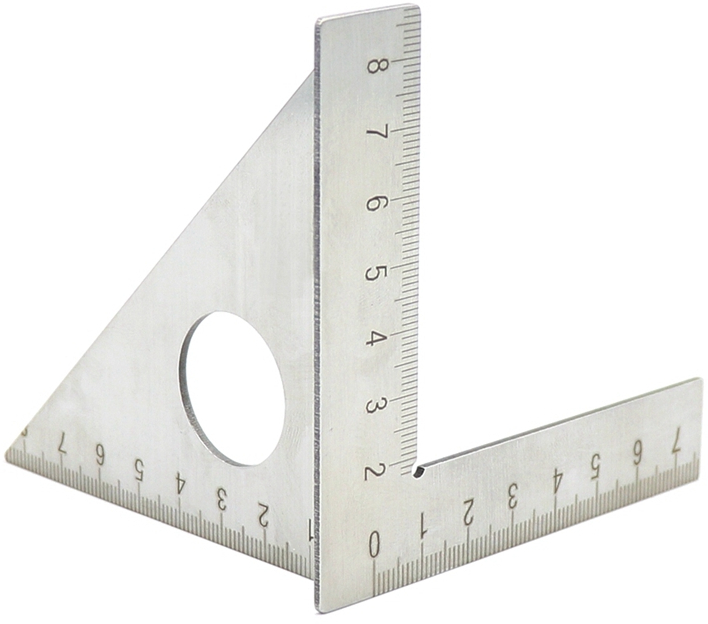}
		  \end{center}
		  \caption{A (3-dimensional) rigid ruler.}
		  \label{Fig:Regolo}		
	\end{figure}
An example is given by the act of reading the coordinates of the position of a certain point-particle 
counting the ticks on the axes of a rigid ruler (see figure \ref{Fig:Regolo}) displaced in the physical space.

If the system is in particular conservative, one can single out a particular observable $H$, called "the Hamiltonian (function)", to be interpreted as the "energy" of the system. 
The flow along $v_H$ yields the time evolution of the system and taking the Poisson bracket of the Hamiltonian with a certain observable encodes how the value of the measurable quantity change along the motion of the system.

	\medskip
	Let us now discuss a possible candidate to be the corresponding algebraic counterpart of a generic multisymplectic manifold.
	When trying to generalize the construction of the Poisson bracket from the observables $C^{\infty}(M)$ of a symplectic manifold to an arbitrary $n$-plectic manifold, one could follow two paths:
	\begin{itemize}
		\item passing to multivector fields and thus considering Hamiltonian pairs composed of a $k$-vector field and a $(n+1-k)$-form satisfying HDDW equation (\eg \cite[\S 4]{Herman2017} );
		\item Focus on Hamiltonian $1$-forms, as suggested by Baez and Rogers in
		\cite{Baez2010,Rogers2010}.	
	\end{itemize}		
	We will adhere to the second path. 
	In this case the definition of Poisson bracket can be translated verbatim to $\Omega_{\ham}^{n-1}(M,\omega)$ as follows:
	\begin{definition}[Binary multi-bracket for Hamiltonian $(n-1)$-forms]\label{Def:BinaryBracketofHamiltonianForms}
		Let be $(M,\omega)$ an $n$-plectic manifold, we call \emph{binary multibrackets} of Hamiltonian $(n-1)$-forms the following bilinear operator
		\begin{displaymath}
			\morphism{\lbrace\cdot,\cdot\rbrace}
			{\Omega^{n-1}_{\ham}(M,\omega)\otimes\Omega^{n-1}_{\ham}(M,\omega)}
			{\Omega^{n-1}_{\ham}(M,\omega)}
			{(\alpha,\beta)}
			{\iota_{\vHam_\beta}\iota_{\vHam_\alpha}\omega}
		\end{displaymath}
		Note that $\lbrace\cdot,\cdot \rbrace = (\iota^2_{\mathfrak{X}}\omega) \circ \pi_{\ham}$ \ie is given by the precomposition of the signed multicontraction map given in equation \eqref{Eq:signedMulticontraction} with $\pi_{\ham}:\Omega^{n-1}_{\ham}(M,\omega) \twoheadrightarrow \mathfrak{X}_{\ham}(M,\omega)$, the standard projection from Hamiltonian forms to their corresponding Hamiltonian vector fields.
	\end{definition}
	When $n=1$, one clearly recovers the Lie algebra structure on the "classical observables" smooth functions of $M$. 
	When $n\geq 2$, the binary bracket shares some properties of the Poisson one but, crucially, it fails to give a Lie algebra structure and a derivation\footnote{In particular, we do not have a canonical associative algebra structure on $\Omega_{\ham}^{n-1}(M,\omega)$. Clearly Hamiltonian forms are not closed under the wedge product.}.
	
	\begin{lemma}\label{Lem:BinBrackofHamFormsisHamiltonian}
		Let be $(M,\omega)$ an $n$-plectic manifold, the binary bracket given in definition \ref{Def:BinaryBracketofHamiltonianForms} is a well-defined skew-symmetric bilinear maps valued in $\Omega^{n-1}_{\ham}(M,\omega)$ satisfying the following properties:
		\begin{itemize}
			\item For any $\alpha,\beta \in \Omega^{n-1}_{\ham}(M,\omega)$, one has 
			$$\d ~\lbrace \alpha, \beta\rbrace = -\iota([\vHam_\alpha,\vHam_\beta])\omega~,$$
			\ie the following diagram commutes
			\begin{displaymath}
				\begin{tikzcd}
					\big(\Omega^{n-1}_{\ham}(M,\omega)\big)^{\otimes 2}
					\ar[r,"{\lbrace \cdot , \cdot \rbrace}"]
					\ar[d,"\pi_{\ham}^{\otimes 2}"'] &
					\Omega^{n-1}_{\ham}(M,\omega) \ar[d,"\pi_{\ham}"]
					\\
					\big( \mathfrak{X}_{\ham}(M,\omega)\big)^{\otimes 2} \ar[r,"{[\cdot,\cdot]}"']  &
					\mathfrak{X}_{\ham}(M,\omega)
				\end{tikzcd}
			\end{displaymath}
			\item It satisfies the Jacobi identity up to an exact term, \ie, for any $\alpha_1,\alpha_2,\alpha_3 \in \Omega^{n-1}_{\ham}(M,\omega)$, one has
			\begin{displaymath}
				\lbrace\lbrace \alpha_1,\alpha_2 \rbrace,\alpha_3 \rbrace + \cyc
				=
				\d~\iota(\vHam_{\alpha_1}\wedge \vHam_{\alpha_2} \wedge \vHam_{\alpha_3}) \omega~.
			\end{displaymath}
		\end{itemize}
	\end{lemma}
\begin{proof}
	Both claims follow by specializing corollary \ref{lemma:multicartan} to the cases $m=2,3$ and noting that
	\begin{displaymath}
		\begin{split}
		\lbrace\lbrace \alpha_1,\alpha_2 \rbrace,\alpha_3 \rbrace + \cyc
		=&~
		\left(\iota_{\vHam_3}\iota_{[\vHam_1,\vHam_2]} - \iota_{\vHam_2}\iota_{[\vHam_1,\vHam_3]} + \iota_{\vHam_1}\iota_{[\vHam_2,\vHam_3]}\right)\omega 
		=\\
		=&~ 
		\iota(-\partial \vHam_1\wedge\vHam_2\wedge\vHam_3)\omega
		~.
		\end{split}
	\end{displaymath}
\end{proof}

\begin{remark}[Lie algebra of Hamiltonian forms modulo boundaries]\label{Rem:LieAlgebraHamModulo}
	Modding out exact terms, \ie considering the vector space
	\begin{displaymath}
		\mathfrak{g}= \dfrac{\Omega^{n-1}_{\ham}(M,\omega)}{\d \Omega^{n-1}(M)}~,
	\end{displaymath}
	which is isomorphic to $\mathfrak{X}_{\ham}(M,\omega)$ when $\omega$ is not degenerate, the binary bracket $\lbrace\cdot,\cdot\rbrace$ induces on $\mathfrak{g}$ a well-defined Lie bracket since closed forms are Hamiltonian with trivial Hamiltonian vector fields (compare with the standard Lie algebra structure considered in Hydrodynamics, see section \ref{Thm:HydroBracket}).
\end{remark}
The upshot of lemma \ref{Lem:BinBrackofHamFormsisHamiltonian} is that $(\Omega_{\ham}^{n-1}(M,\omega),\lbrace\cdot,\cdot\rbrace)$ fails to be a Lie algebra (except for the case $n=1$) since $\lbrace\cdot,\cdot \rbrace$ does not satisfy the Jacobi identity in general.
However, the failure is somewhat mild, namely, it is controlled by the exterior derivative of the differential $(n-2)$-form $(\iota^3_{\mathfrak{X}}\omega) \circ \pi_{\ham}$.

In the earliest development of multisymplectic geometry, which was basically motivated by the study of mechanical systems with a continuous set of degrees of freedom (classical field theories), it has been suggested to cure this problem by dividing out inconsistencies, \ie, by considering everything up to exact terms ("Total divergences" in the physics lingo). See for example \cite{Carinena1991b} or section \ref{Sec:IdroPoisson} for the hydrodynamics case.
\\
With the advent of homological methods in mathematical physics (Kontsievich, Stasheff, Vinogradov) this point of view has become obsolete.
Indeed, before modding out exact terms, one should be naturally led to check whether the ternary bracket $(\iota^3_{\mathfrak{X}}\omega) \circ \pi_{\ham}$, controlling the failing of the Jacobi identity pertaining to $\lbrace\cdot,\cdot\brace$, satisfies an higher analogue of the Jacobi equation.
More specifically, one should look for a suitable completion of $\lbrace\cdot,\cdot\rbrace$ and $(\iota^3_{\mathfrak{X}}\omega) \circ \pi_{\ham}$ to a family of multibrackets defining a $L_\infty$-algebra structure. 
Such completion has been proved to be possible; the following explicit construction is due by Rogers:
\begin{definition}[$L_\infty$-algebra of observables \emph{(\cite[Thm. 5.2]{Rogers2010}, \cf also \cite{Barnich1998})}]\label{Def:RogersAlgebra}
	Given a $n$-plectic manifold $(M,\omega)$, we call \emph{$L_{\infty}$-algebra of observables} (higher observables) associated to $(M,\omega)$
	the $L_\infty$-algebra $L_{\infty}(M,\omega)=(L,\{[\cdot,\cdots,\cdot]_k \}_{k\geq 1})$ where:
	\begin{itemize}
		\item 
	the underlying graded vector space $L$ is given by
	\begin{equation}\label{eq:Lspace}
		L^i=\begin{cases}
			\Omega_{\ham}^{n-1}(M,\omega) 
			& \quad~\text{if } i=0
			\\
			\Omega^{n-1+i}(M) 
		 	&  \quad~\text{if } 1-n \leq i\leq -1
		 	\\			
			0 & \quad ~\text{otherwise} 
			~;
		\end{cases}
	\end{equation}
	\item 
	the $n+1$ non-trivial multibrackets 
	$\lbrace [\cdot,\cdots,\cdot]_k : L^{\wedge k} \to L ~| 1\leq k \leq n+1\rbrace$
	are defined, for any given $\alpha_i\in L$, as
	\begin{displaymath}
		[\alpha]_1 = 
		\begin{cases}
			0 & \quad\text{if~} |\alpha| = 0	
			\\
			\d \alpha & \quad \text{if~} |\alpha| \leq -1
		\end{cases}
	\end{displaymath}
	and, for $ 2 \leq k \leq n+1$, as
	\begin{displaymath}
		[\alpha_1,\dots,\alpha_k]_k = 
		\begin{cases}
			\varsigma(k) \iota(\vHam_{\alpha_1}\wedge\dots\wedge\vHam_{\alpha_k})~\omega
			& \quad\text{if~} |\alpha_i|=0 \text{ for } 1\leq i \leq k
			\\
			0 & \quad\text{otherwise}		
		\end{cases}
		~.
	\end{displaymath}
	\end{itemize}
		In the above equation, $\vHam_{\alpha_k}=\pi_{\ham}(\alpha_k)$ denotes the Hamiltonian vector field associated to $\alpha_k\in \Omega^{n-1}_{\ham}(M,\omega)$ (see definition \ref{Def:Hamiltonianform}) and $\varsigma(k) := - (-1)^{\frac{k(k+1)}{2}}$ is the total Koszul sign (see definition \ref{Def:SigmaSign}).
\end{definition} 
\begin{remark}
	We emphasize that here, as in \cite{Callies2016,Ryvkin2016}, is adopted an index convention different from what originally employed by \cite{Rogers2010}. 
Namely, we adopt the "cohomological convention" with non-trivial components concentrated in negative degrees.
\end{remark}

\begin{remark}[The chain complex underlying $L_\infty(M,\omega)$]\label{Rem:RogersChainComplex}
	The $L_\infty$-algebra of observables consists of a cochain complex	 $(L,d)$
	\begin{displaymath}
		\begin{tikzcd}[column sep= small,row sep=small]
			0 \ar[r]
			& L^{1-n}\ar[symbol=\coloneqq,d] \ar[r]
			& \cdots\ar[r]
			& L^{-k}\ar[symbol=\coloneqq,d]\ar[r]
			& \cdots\ar[r]
			& L^{-1}\ar[symbol=\coloneqq,d]\ar[r]
			& L^0\ar[symbol=\coloneqq,d]\ar[r]
			&[-2em] 0\\
			& \Omega^0(M)\ar["d",r]
			& \cdots\ar["d",r]
			& {\Omega^{n-1-k}(M)} 	\ar["d",r]&\cdots\ar["d",r]
			& \Omega^{n-2}(M)\ar["d",r]
			& \Omega^{n-1}_{\textrm{Ham}}(M,\omega)&
		\end{tikzcd},
	\end{displaymath}
	which is a truncation of the de-Rham complex with a degree shift putting (Hamiltonian) $(n-1)$-form in degree zero, \ie
	\begin{displaymath}
		L= \trunc_{n-1}(\Omega(M))[n-1]\oplus \Omega^{n-1}_{\ham}(M,\omega)
		~,
	\end{displaymath}
	endowed with $n$ (skew-symmetric) multibrackets $(2 \leq k \leq n+1)$ obtained by the signed multicontraction
		\begin{equation}\label{Eq:RogersKBrackets}
		\begin{tikzcd}[column sep= small,row sep=0ex]
				[\cdot,\dots,\cdot]_k:= (\iota^k_{\mathfrak{X}}\omega) \circ \pi_{\ham} \colon& \Lambda^k\left(\Omega^{n-1}_{\textrm{Ham}}\right) 	\arrow[r]& 				\Omega^{n+1-k} \\
				& \alpha_1\wedge\dots\wedge\alpha_k 	\ar[r, mapsto]& 	\varsigma(k)\iota_{\vHam_{\alpha_k}}\dots\iota_{\vHam_{\alpha_1}}\omega 
		\end{tikzcd}
		~.		
	\end{equation}
\end{remark}

\begin{remark}[Reinterpretation in terms of $\mathcal{B}_{(n)}$]
	Observe that extending the chain complex underlying $L_{\infty}(M,\omega)$ to the right with the space of Hamiltonian vector fields, one gets a cochain sub-complex of the complex $\mathcal{B}_{(n)}$ introduced in proposition \ref{lem:signedmulticontractionLinfinitymorphism}.
	Namely, one has the inclusion $L[1]\oplus \mathfrak{X}_{\ham}(M)\hookrightarrow \mathcal{B}_{(n)}$ with components given by the vertical arrows in the following commutative diagram:
	\begin{displaymath}
		\begin{tikzcd}[column sep = small]
			C^{\infty}(M) \ar[r,"\d"] \ar[d,equal]
			&
			\dots \ar[r,"\d"]
			&
			\Omega^{n-2}(M) \ar[r,"\d"] \ar[d,equal]
			&[1.5em]
			\Omega^{n-1}_{\ham}(M,\omega) \ar[r,"\pi_{\ham}"] \ar[d,hook]
			&[1.5em]
			\mathfrak{X}_{\ham}(M,\omega) \ar[d,"-\iota_{\mathfrak{X}_{\ham}}^1\omega"]
			\ar[r]
			&
			0
			\\
			C^{\infty}(M) \ar[r,"\d"]
			&
			\dots \ar[r,"\d"]
			&
			\Omega^{n-2}(M) \ar[r,"\d"]
			&
			\Omega^{n-1}(M) \ar[r,"\d"]
			&
			\d \Omega^{n-1}(M) \ar[r]
			&
			0
		\end{tikzcd}
		~.
	\end{displaymath}	
	The two rightmost squares express respectively the fact that all closed $(n-2)$-forms are Hamiltonian with trivial Hamiltonian vector field, see remark \ref{Rem:ClosedformsTrivialHamiltonian}, and the very definition of Hamiltonian forms, see remark \ref{Rem:HamiltonianPairs}.
\end{remark}
\begin{remark}[Precursors]
	It should be noted that an early precursor of the (Chris) Rogers' construction can be found in the work of Claude Roger 	\cite[Thm. 6.1]{Roger2012}.
	Briefly, he showed that for any compact, connected, orientable $(p+1)$-dimensional manifold $(p\geq 1)$ there is an associated $L_p$-algebra (see also \cite[Cor. 6.11]{Zambon2012} for a dual version).
	Roger's work was in the context of \emph{unimodular vector fields} (also known as \emph{divergence-free} or \emph{solenoidal} fields), a framework that is akin to the construction that we are going to study in chapter \ref{Chap:MauroPaper}.
\end{remark}

\begin{lemma}[Theorem 5.2 in \cite{Rogers2010}] The $L_\infty$-algebra of observables $L_\infty(M,\omega)$ is a grounded $L_n$-algebra (see \ref{Def:groundedLinfinity}).
\end{lemma}
\begin{proof} 
	By its very definition ,$L$ is concentrated in degrees $(1-n),\dots,0$ and $[\cdots]_k = (\iota^k_{\mathfrak{X}(M)}\omega) \circ \pi_{\ham}$ is non-zero only when evaluated on a degree $0$ element \ie on $k$ Hamiltonian $(n-1)$-forms. 
	According to remark \ref{Rem:GroundedEasyAxioms}, one has only to prove the following equality
	\begin{displaymath}
		[\cdot]_1 \ca\, [\cdots]_{k+1} = [\cdots]_k \ca\, [\cdot,\cdot]_2 ~,
	\end{displaymath}
	where $\ca$ denotes the skew-symmetric \RN product (see equation \eqref{Eq:RNProducts-explicit}), for any $k\geq 1$.
	Note in particular that on the left-hand side the \RN products reads simply as the composition $\circ$. 
	On given Hamiltonian forms $\alpha_1,\dots, \alpha_{k+1}$ this can be read as
\begin{displaymath}
	\mathclap{
	\begin{aligned}
		\varsigma(k&+1) \d \iota(\vHam_{\alpha_1}\wedge\dots\wedge\vHam_{\alpha_{k+1}} ) \omega
		=
		\\ 
		=&~
		\varsigma(k)\sum_{i<j}(-1)^{i+j}
		\iota(
			\vHam_{[\alpha_i,\alpha_j]_2} \wedge \vHam_{\alpha_1} \wedge \dots \wedge
			\widehat{\vHam_{\alpha_i}} \wedge\dots\wedge
			\widehat{\vHam_{\alpha_j}} \wedge\dots\wedge
			\vHam_{\alpha_{k+1}}
		)\omega
		=
		\\
		=&~ \varsigma(k) \iota(\partial \vHam_{\alpha_1}\wedge\dots\wedge\vHam_{\alpha_{k+1}}) \omega
			~,
	\end{aligned}
	}
\end{displaymath}
where, in the first equation, lemma \ref{Lem:BinBrackofHamFormsisHamiltonian} has been used .
One can easily read the previous equation as the multi-Cartan magic formula of Lemma \ref{lemma:multicartan}, restricted to the case of being $\omega$ a closed form and the involved vector fields being Hamiltonian.
\end{proof} 

\begin{remark}[Homological and cohomological convention for $L_{\infty}(M,\omega)$]
	We emphasize that in definition \ref{Def:RogersAlgebra} we presented the $L_\infty$-algebra of observables with the so-called \emph{cohomological convention}; namely $[\cdot]_1$ acts by raising the degree and the underlying cochain complex $L$ is concentrated in degrees $(1-n),\dots,n$.
	\\
	One can give an equivalent definition in the \emph{homological convention} of $L_\infty(M,\omega$
		as the chain complex $L_\bullet$
	\begin{center}
		\begin{tikzcd}[column sep= small,row sep=small]
				0 \ar[r]& L_{n-1}\ar[symbol=\coloneqq,d] \ar[r]&
		\cdots\ar[r]&L_{k-2}\ar[symbol=\coloneqq,d]\ar[r]&\cdots\ar[r]&
		L_1\ar[symbol=\coloneqq,d]\ar[r]&L_0\ar[symbol=\coloneqq,d]\ar[r]&0\\
		& \Omega^0\ar["d",r]&\cdots\ar["d",r]&{\Omega^{n+1-k}}
		\ar["d",r]&\cdots\ar["d",r]&\Omega^{n-2}\ar["d",r]
		&\Omega^{n-1}_{\textrm{Ham}}&
		\end{tikzcd},
	\end{center}				
	(which is a truncation of the de-Rham complex with inverted grading)
	endowed with $n$ (skew-symmetric) multibrackets $(2 \leq k \leq n+1)$
	\begin{equation}
		\begin{tikzcd}[column sep= small,row sep=0ex]
				[\cdot,\dots,\cdot]_k \colon& \Lambda^k\left(\Omega^{n-1}_{\textrm{Ham}}\right) 	\arrow[r]& 				\Omega^{n+1-k} \\
				& \sigma_1\wedge\dots\wedge\sigma_k 	\ar[r, mapsto]& 	\varsigma(k)
				~\iota_{v_{\sigma_k}}\dots\iota_{v_{\sigma_1}}\omega 
		\end{tikzcd}		
	\end{equation}		
	Clearly the two definition are not unrelated, one can retrieve one from the other by applying the reverse-ordering functor $\setminus \blank$.
\end{remark}

\begin{remark}[Local existence of the $L_\infty$-algebra of observables]
	The introduction of the algebra of observables, as given in definition \ref{Def:RogersAlgebra}, could appear rather artificial at first sight.
	However, there is a motivation, proposed in \cite[\S 5]{Rogers2010}, for the local existence of such an object.
	Observe that if $(M,\omega)$ is contractible (\eg if $M$ is an open neighbourhood in a given bigger multisymplectic manifold) the cohomology groups of the chain complex $L$ (\cf  remark	\ref{Rem:RogersChainComplex}) reads as follows
	\begin{displaymath}
		H^{i}(L) =
		\begin{cases}
			\dfrac{\Omega_{\ham}^{n-1}(M,\omega)}{\d \Omega^{n-2}(M,\omega)}:= \mathfrak{g}
			& \text{ if } i=0
			\\[1em]
			0 & \text{ if } 1 \leq i \leq 1-n
			\\[.5em]
			\R & \text{ if } i = 1-n
		\end{cases}
		~.
	\end{displaymath}
	Hence, the complex (augmentation of $L$)
	\begin{displaymath}
		\begin{tikzcd}[column sep = small]
			\widetilde{L}:\quad &
			0 \ar[r] &
			\R \ar[r,hook] &
			C^{\infty}(M)\ar[r,"\d"] &
			\cdots \ar[r,"\d"] &
			\Omega^{n-1}_{\ham}(M,\omega)
		\end{tikzcd}
	\end{displaymath}
	 is a resolution (in the sense of \cite[\S 2.1]{Barnich1998}) of the Lie algebra $(\mathfrak{g},\lbrace \cdot,\cdot\rbrace)$, introduced in remark \ref{Rem:LieAlgebraHamModulo}. In particular all cohomology groups are zero except for $H^0(\widetilde{L})=\mathfrak{g}$.
	 \\
	 The key point is that, from a general homological algebra result \cite[Thm. 7]{Barnich1998}, one can naturally endow a resolution of a Lie algebra with a $L_\infty$-algebra structure obtaining by suitable iteration of the complex coboundary operator and the Lie bracket.
	Hence, this proves global existence of a $L_\infty$-algebra structure with the same cohomology as $L$ for any contractible multisymplectic manifold or local existence for a general multisymplectic manifold.
	Definition \ref{Def:RogersAlgebra} has been originally proposed by Rogers as a suitable globalization of such local construction.
	\\
	Since in the non-degenerate case $\mathfrak{g}$ is isomorphic to the Lie algebra of Hamiltonian vector fields, one can think about $L_\infty(M,\omega)$ as a $L_\infty$-extension of $\mathfrak{X}_{\ham}(M,\omega)$.
\end{remark}

\begin{remark}[Standard projection to Hamiltonian vector fields]\label{Rem:ProjectiontoHamVfields}
	By construction, the restriction of $[\cdot,\cdot]_2$ to Hamiltonian forms, coincides with the bracket $\lbrace\cdot,\cdot\rbrace$ given in definition \ref{Def:BinaryBracketofHamiltonianForms}.
	Lemma \ref{Lem:BinBrackofHamFormsisHamiltonian} implies then that the surjection $\pi_{\ham}$ mapping Hamiltonian forms into their corresponding Hamiltonian vector fields can be lifted to a strict $L_\infty$-morphism 
	\begin{displaymath}
		\begin{tikzcd}
			L_\infty(M,\omega) \ar[r,dashed] \ar[d]& \mathfrak{X}^1_{\ham}(M,\omega)
			\\
			\Omega^{n-1}_{\ham}(M,\omega) \ar[ur,two heads,"\pi_{\ham}"',sloped]
		\end{tikzcd}
		~,
	\end{displaymath}
	defined by the projection 
	\begin{displaymath}
		\morphism{\pi_{\ham}}
		{L_\infty(M,\omega)}
		{\mathfrak{X}^1_{\ham}(M,\omega)}
		{\alpha}
		{
			\begin{cases}
				\vHam_{\alpha} & \text{ if } |\alpha|=0
				\\
				0 & \text{ otherwise}
			\end{cases}	
		}
		~,
	\end{displaymath}		
	denoted with the same symbol with a slight abuse of notation\footnote{Note that any Lie algebra can be seen as an $L_\infty$-algebra concentrated in degree $0$, therefore any $L_\infty$-morphism $L\to\mathfrak{g}$ is simply given by a linear map $L_0 \to \mathfrak{g}$ preserving the binary bracket.}.
\end{remark}

\begin{remark}[Degenerate case]\label{Rem:DegenerateCase}
	Definition \ref{Def:RogersAlgebra} translates verbatim to the pre-$n$-plectic, \ie degenerate, 
	case replacing $\Omega^{n-1}_{\ham}(M,\omega)$ with the vector space of Hamiltonian pairs $Ham^{n-1}(M,\omega)$ as given in remark \ref{Rem:HamiltonianPairs} (cfr \cite[Thm 5.2]{Rogers2010},\cite[\S 3]{Ryvkin2016a},\cite[Thm 4.7]{Callies2016},\cite[Prop. 3.2]{Fiorenza2014a}).
	Namely, one can define a graded vector space
	
	\begin{displaymath}
	 \mathclap{
		Ham_{\infty}(M,\omega)^i =
		\begin{cases}
			\left\lbrace\left. \pair{\vHam}{\alpha}\in\mathfrak{X}^1(M)\oplus\Omega^{n-1}(M) ~\right|~ \d\alpha = -\iota_{\vHam}\omega \right\rbrace
			& \text{ if } i=0
			\\[1em]
			\Omega^{n-1+i}(M) & \text{ if } 1-n \leq i \leq 1
		\end{cases}
	 }
	\end{displaymath}
	with multibrackets defined as in equation \eqref{Eq:RogersKBrackets} interpreting the map $\pi_{\ham}$ as the projection on the first term of the direct sum:
	\begin{displaymath}
		\pi_{\ham}: Ham_{\infty}(M,\omega) \twoheadrightarrow \mathfrak{X}_{\ham}(M,\omega)
	\end{displaymath}
	(which can be again seen as a strict $L_\infty$-morphism).
	\\
	When $\omega$ is non-degenerate, there is an $L_\infty$-isomorphism $Ham_{\infty}(M,\omega)\cong L_\infty(M,\omega)$ (see \cite[Prop. 4.8]{Callies2016} and section \ref{Chap:MarcoPaper}).
	\note{aggiungere riferimento preciso al lemma nel capitolo sul lavoro con Marco}
\end{remark}

\begin{remark}[Dg Leibniz algebra of higher observables]
	As pointed out by Rogers in \cite[\S 6]{Rogers2010}, there is another natural algebraic structure that one can construct out of $\Omega_{\ham}^{n-1}(M,\omega)$.
	Recall that, in the symplectic case, \ie $n$-plectic case with $n=1$, the Poisson bracket between two give smooth functions $f,g\in C^{\infty}(M)$ can be seen as
	\begin{displaymath}
		\lbrace f,g \rbrace = \mathcal{L}_{\vHam_f} g 
	\end{displaymath}
	hence, $\lbrace f,\cdot \rbrace$ is a degree $0$ derivation making $\Omega^0_{\ham}(M)=C^{\infty}(M)$ into a Poisson algebra.
	\\
	When $n\geq 2$ one should not expect that $L_\infty(M,\omega)$ would acts like a Poisson algebra. 
	In such case, the Cartan magic formula would rather imply that
	\begin{displaymath}
		\mathcal{L}_{\vHam_\alpha} \beta = \lbrace \alpha, \beta \rbrace + \d \iota_{\vHam_\alpha}\beta~,
	\end{displaymath}
	thus suggesting to introduce the following non-symmetric binary bracket
		\begin{displaymath}
			\morphism{\llbracket \cdot,\cdot\rrbracket }
			{\Omega^{n-1}_{\ham}(M,\omega)\otimes\Omega^{n-1}_{\ham}(M,\omega)}
			{\Omega^{n-1}_{\ham}(M,\omega)}
			{(\alpha,\beta)}
			{\mathcal{L}_{\vHam_\alpha}\beta}
		\end{displaymath}
	Notice that $\llbracket\cdot,\cdot\rrbracket$ equals $\lbrace\cdot,\cdot\rbrace$ modulo boundary terms but with the different flavour of measuring the rate of change ("conservation", in the spirit of definition \ref{Def:conservedQuantities}) of the second Hamiltonian observable along the flow given by the first.
	The bracket $\llbracket\cdot,\cdot \rrbracket$ can be easily extended to a binary bracket on $L$ as follows
	\begin{displaymath}
		\llbracket x, y \rrbracket =
		\begin{cases}
			\mathcal{L}_{\vHam_x} \beta & \text{ if } |x|=0
			\\
			0 & \text{ otherwise}
		\end{cases}
		\qquad \forall \alpha,\beta \in L
	\end{displaymath}
	and it has been proved in \cite[Prop. 6.3]{Rogers2010} that the triple $(L,\delta,\llbracket\cdot,\cdot\rrbracket )$, with $\delta=[\cdot]_1$, forms a \emph{differential graded Leibniz algebra} (or \emph{dg Loday algebra}). 
	Namely the following two compatibility relations hold for any $x,y,z\in L$:
	\begin{displaymath}
		\begin{aligned}
			\delta \llbracket x,y\rrbracket
			=&
			\llbracket\delta x,y\rrbracket +
			(-)^{|x|}\llbracket x,\delta y\rrbracket
			\\
			\llbracket x,\llbracket y, z\rrbracket\rrbracket
			=&
			\llbracket\llbracket x,y \rrbracket, z\rrbracket +
			(-)^{|x||y|}\llbracket y,\llbracket x, z\rrbracket\rrbracket
			~.
		\end{aligned}
	\end{displaymath}

\end{remark}

	\begin{remark}[Observables $L_\infty$-algebra as a homotopy pullback]\label{Rem:FiorenzaPullback}
		It has been noticed in \cite[Thm. 3.12]{Fiorenza2014a} that the $L_\infty$-algebras of observables can be seen as a \emph{homotopy pullback} along the "signed multicontraction" $L_\infty$-morphism $(\iota_{\mathfrak{X}_{\ham}}\omega): \mathfrak{X}^1_{\ham}(M) \to \mathcal{B}_{(n)}$ introduced in lemma \ref{lem:signedmulticontractionLinfinitymorphism}.
		Namely it sits in the following pullback diagram
		\begin{displaymath}
			\begin{tikzcd}
				L_{\infty}(M,\omega)\arrow[dr, phantom, "\scalebox{1.5}{$\lrcorner$}" , very near start, color=black] \ar[r]\ar[d,"\pi_{\ham}"']
			& 
			0 \ar[d]\\
			\mathfrak{X}^1(M) \ar[r,"(\iota_{\mathfrak{X}_{\ham}})"] & \mathcal{B}	
			\end{tikzcd},
		\end{displaymath}
		commuting modulo homotopies, \ie $2$-morphisms in the category of $L_\infty$-algebras.
		\note{
			"$Ham_{\infty}(M,\omega)$ is the \emph{homotopy fibre} of the \emph{structure maps} $(\iota_{\mathfrak{X}}\circ \pi_{\ham}) = (\iota_{\mathfrak{g}})$ with $\mathfrak{g}=\mathfrak{X}_{\ham}$
			
			The notation $\mathcal{B}^n$ is to remind an analogous construction for the "classifying space of vector bundles"
		}
		\note{
			sarebbe interessante dare esplicitamente questa omotopia (sospetto che coinvolga il pairing). Servirebbe una definizione esplicita di homotopia

			In section 3.2 of dolgushev,hoffung baez has been shown the mapping space for $L_\infty$-algebras ie.e the simplicial set $\text{Map}_{\bullet}(L,L')$ where $0$-simplcies are $L_\infty$-morphisms and 1-simplicies are homotopies between morphisms.
			In this term the homotopy-pullback diagram means that $[\iota_{\mathfrak{X}}\circ \pi_{\ham}]=[0]\in \pi_0\text{Map}_{\bullet}(\mathfrak{g},\mathcal{B})$.
			
				}
		Without going into the precise definition of homotopy\cite{Buijs2013}, or equivalence \cite{Dolgushev2007}\cite[Appendix]{Fregier2015}, between two $L_\infty$-morphism, one could notice that this weaker notion of commutativity is another reverberation of the multi-Cartan magic rule.
		For instance, in the $1$-plectic case the lower-left corner of the above diagram restricts to
		\begin{displaymath}
			\begin{tikzcd}
			C^{\infty}(M) \ar[dr,dashed,"(f)"] \ar[d,"\pi_{\ham}"'] & \\
			\mathfrak{X}_{\ham} \ar[r,"(\iota_{\mathfrak{X}_{\ham}})"'] &
			C^{\infty}(M)[1]\oplus \d C^{\infty}(M)
			\end{tikzcd}
		\end{displaymath}
		where we denoted with $(f)$ the $L_2$-morphism resulting from the composition.
		This consists of two components.
		The first one, $f_1= (\iota^1_{\mathfrak{X}(M)}\omega) \circ \pi_{\ham}$, can be seen a chain map homotopically equivalent to $0$. 
		Namely the identity map is a chain homotopy $\id: f_1 \Rightarrow 0$ since, according to the very definition of Hamiltonian vector fields, the following diagram commutes:
		\begin{displaymath}
		\begin{tikzcd}
			& C^{\infty}(M) \ar[dl,"-\id",sloped]\ar[d,"f_1"]
			\\
			C^{\infty}(M) \ar[r,"\d"]&
			\d C^{\infty}(M)
		\end{tikzcd}
		~.
		\end{displaymath}
	On the other hand, the second component $f_2:C^{\infty}(M)\otimes C^{\infty}(M) \to C^{\infty}(M)$ can be seen itself as a chain homotopy from $0$ to the chain map $g$ defined by the commutation of the right square in the diagram below
	\begin{displaymath}
		\begin{tikzcd}
			&
			C^{\infty}(M)\otimes C^{\infty}(M) \ar[r,"\pi\otimes\pi"]\ar[d,dashed,"g"]\ar[dl,"f_2"',dashed]&
			\mathfrak{X}_{\ham}\otimes \mathfrak{X}_{\ham}\ar[d,"{-[\cdot,\cdot]}"] \\
			C^{\infty}(M) \ar[r,"\d"]&
			\d C^{\infty}(M) & 
			\mathfrak{X}_{\ham} \ar[l,"\iota_{\mathfrak{X}_{\ham}}^i\omega"]
		\end{tikzcd}
	\end{displaymath}
	since, as follows from the first point in lemma \ref{Lem:BinBrackofHamFormsisHamiltonian}, the whole diagram commutes.
	\end{remark}

	Many non-trivial explicit examples of the algebra of observables $L_\infty(M,\omega)$ can be found in \cite{Rogers2011}\cite{Callies2016}\cite{Ryvkin2018}.
	More will be studied in chapters \ref{Chap:LeonidPaper} and \ref{Chap:MauroPaper}.
	We only mention here an example directly following from example \ref{Ex:SumProductMultiSymp}:
\begin{example}[Higher observables for sums and products]
	Notice that the observables $L_\infty$-algebra construction behaves compatibly with the sum and product constructions of multisymplectic manifolds given in example \ref{Ex:SumProductMultiSymp}.
	Namely, when $(M,\omega)$ and $(\tilde{M},\tilde{\omega})$ are two multisymplectic structure of the same order, there exists a strict $L_\infty$-morphisms \cite[Ex. 6.5]{Ryvkin2018}
	\begin{displaymath}
		\lambdamorphism{L_\infty(M,\omega)\oplus L_\infty(\tilde{M},\tilde{\omega})}
		{L_\infty(M\times \tilde{M},\pi_M^\ast\omega+\pi_{\tilde{M}}^\ast\tilde{\omega})}
		{(\alpha,\beta)}
		{\pi_M^\ast \alpha + \pi_{\tilde{M}}^\ast \beta}
		~,
	\end{displaymath}
	and, for any $\omega,\tilde{\omega}$, possibly in different order, there exists a non-strict $L_\infty$-morphism \cite[Thm. 4.2]{Shahbazi2016}	
	\begin{displaymath}
		L_\infty(M,\omega)\oplus L_\infty(\tilde{M},\tilde{\omega})
		\to
		L_\infty(M\times \tilde{M}, \pi_M^\ast \omega \wedge \pi_{\tilde{M}}^\ast \tilde{\omega})
		~.
	\end{displaymath}
\end{example}

\section{Symmetries and \Momaps}\label{Sec:HCMM}
When one fixes a form $\omega$ on a manifold $M$ it is natural to highlight the group actions preserving this extra structure, also known as "symmetries".
\\
Consider a multisymplectic manifold $(M,\omega)$, we introduce the following nomenclature pertaining to global and local symmetries:
\begin{definition}[Multisymplectic actions]
	A smooth action $\vartheta:G\action M$ of a Lie group $G$ on $M$ is called \emph{multisymplectic} if it preserves the multisymplectic form, \ie 
	\begin{displaymath}
		\hat{\vartheta}_g^* \omega = \omega \qquad \forall g\in G
	\end{displaymath}
	where $\hat{\vartheta}_g=\vartheta(\cdot,g)$ is the diffeomorphism giving the action of $G$ on $M$.
	\\
	A Lie algebra morphism 	$\vAct: \mathfrak{g} \rightarrow \mathfrak{X} (M)$	is called an \emph{infinitesimal multisymplectic action} if it acts by multisymplectic vector fields, \ie 
	\begin{displaymath}
		\mathcal{L}_{\vAct_\xi} \omega = 0 \quad \forall \xi \in \mathfrak{g}
		~.
	\end{displaymath}		
\end{definition}
Clearly, these two concepts are not unrelated:
\begin{lemma}
	Consider a multisymplectic manifold $(M,\omega)$.
	If the action $\vartheta:G\action M$ is multisymplectic then the corresponding infinitesimal action $\vAct: \mathfrak{g}\to \mathfrak{X}(M)$, given by the fundamental vector fields \ie
	\begin{equation}\label{eq:LeftFundVF}
		v_\xi(m) = \left.\dfrac{d}{dt}\right\vert_0	\vartheta(m,\exp(- t\xi)) \qquad \forall m \in M , \xi \in \mathfrak{g}
		~,
	\end{equation}		
	is multisymplectic in the sense of definition \ref{def:Hamiltonianvfields}.
	\\
	If the Lie group $G$ is in particular connected, also the converse is true.
\end{lemma}

A multisymplectic action acts infinitesimally by multisymplectic vector fields, we give a special name when the fundamental vector fields are Hamiltonian:
	\begin{definition}[Weakly Hamiltonian actions]\label{Def:WeakHamActions}
		The Lie algebra morphism 	$\vAct: \mathfrak{g} \rightarrow \mathfrak{X} (M)$
		is an \emph{infinitesimal weakly Hamiltonian action} if $\mathfrak{g}$ acts by Hamiltonian vector fields, \ie $\text{Im}(v) \subseteq \mathfrak{X}_{\ham}(M,\omega)$.
		A smooth action $\vartheta:G\action M$ of a Lie group $G$ on $M$ is called \emph{weakly Hamiltonian} if the corresponding infinitesimal action via fundamental vector fields (see equation \eqref{eq:LeftFundVF}) is Hamiltonian.
	\end{definition}

\begin{remark}[Lift condition for weakly Hamiltonian actions {\cite[\S 4.1]{Ryvkin2016a}}]\label{Rem:weakHamasLift}
	The property for the action $\vartheta: G \action M$ to be weakly Hamiltonian can be easily read as the existence of a lift $f_1$ of the corresponding infinitesimal action to the vector space of Hamiltonian $(n-1)$-form, \ie the following diagram commutes in the category of vector spaces:
	\begin{displaymath}
		\begin{tikzcd}[column sep = huge]
			&
			\Omega^{n-1}_{\ham}(M,\omega)
			\ar[d,two heads,"\pi_{\ham}"]
			\\
			&
			\mathfrak{X}_{\ham}(M,\omega)\ar[d,hook]
			\\
			\mathfrak{g} \ar[r,"\vAct"] \ar[ruu,"f_1"]
			&
			\mathfrak{X}(M)		
		\end{tikzcd}
		~.
	\end{displaymath}
	\\
	Such condition can be expressed in a cohomological flavour by considering the following exact sequence of vector spaces :
	\begin{displaymath}
		\begin{tikzcd}[row sep = small,column sep = small]
			0 \ar[r]
			&
			\Omega^{n-1}_{cl}(M) \ar[r,hook]
			&
			\Omega^{n-1}_{\ham}(M,\omega) \ar[r,"\pi_{\ham}"]
			&[1.5em]
			\mathfrak{X}_{\ham}(M,\omega) \ar[r,"\gamma"]
			&
			H_{dR}^n(M)
			\ar[r]
			&
			0
			\\[-1em]
			& & &
			v \ar[r,mapsto]
			&
			{[\iota_{v}\omega]}
		\end{tikzcd}
		~,
	\end{displaymath}
	meaning that $\vartheta$ acts by Hamiltonian vector fields whenever the class $[\gamma \circ \vAct_\xi]=0$ for any $\xi \in \mathfrak{g}$.
	\\
	In particular, if $H^n_{dR}(M)=0$ or the acting Lie algebra satisfy the equation $[\mathfrak{g},\mathfrak{g}]= \mathfrak{g}$ (\eg $\mathfrak{g}$ is the Lie algebra of a semisimple Lie group) any multisymplectic action is also weakly Hamiltonian. 
\end{remark}

\subsection{\Momaps}\label{Sec:hcmm}
When studying weakly Hamiltonian actions on symplectic manifolds, the auxiliary concept of \emph{moment map} takes an exceptionally important role.
The latter is in particular instrumental in many celebrated results in symplectic geometry like the Hamiltonian formulation of the Noether theorem, the classification of toric manifolds and the Kostant-Kirillov-Souriau coadjoint orbits method (see \cite{CannasdaSilva2001} for a complete review).
Most of these results have also a non-trivial application in the geometrical approach to mechanics (see, for instance, \cite{Abraham1978}).
\begin{reminder}[Moment maps in symplectic geometry]\label{Rem:SymplecticMomaps}
	Let be $(M,\omega)$ a symplectic (\ie $1$-plectic) manifold.
	Consider a symplectic smooth action $\vartheta: G \action M$ preserving the $2$-form $\omega$.
	Denote by $\vAct: G \to \mathfrak{X}_{\msy}(M,\omega)$ the corresponding infinitesimal action by fundamental vector fields.
	We define the following:
	\begin{itemize}
		\item A \emph{weak moment map pertaining to $\vartheta$} is a smooth map $\hat{\mu}:M\to \mathfrak{g}^\ast$ such that
		\begin{displaymath}
			\d \langle \hat{\mu}(x), \xi \rangle = -\iota_{v_\xi} \omega_x
			\qquad \forall x\in M, \xi \in \mathfrak{g}~.
		\end{displaymath}
		\item A \emph{(strong) moment map pertaining to $\vartheta$} is a $Ad^\ast$-equivariant weak moment map, \ie the following diagram commutes in the category of smooth manifolds
		\begin{displaymath}
			\begin{tikzcd}
				M \ar[r,"\hat{\mu}"]\ar[d,"\vartheta_g"] & \mathfrak{g}^\ast \ar[d,"Ad^\ast_g"]
				\\
				M \ar[r,"\hat{\mu}"] & \mathfrak{g}^\ast
			\end{tikzcd}
			~.
		\end{displaymath} 
	\end{itemize}
	Dually, see for instance \cite[\S 22.1]{CannasdaSilva2001}, one can define the following:
	\begin{itemize}
		\item
			A \emph{weak comoment map pertaining to $\vartheta$} is a linear map $\check{\mu}:\mathfrak{g} \to C^{\infty}(M)$ such that
		\begin{displaymath}
			\d \check{\mu}(\xi) = -\iota_{v_{\xi}}\omega
			\qquad
			\forall \xi \in \mathfrak{g}
		\end{displaymath}
		
		\item 
			A \emph{(strong) comoment map pertaining to $\vartheta$} is a weak comoment map that is also a Lie algebra morphism, \ie the following diagram commutes in the category of vector spaces
		\begin{displaymath}
			\begin{tikzcd}
				\mathfrak{g}^{\otimes 2} \ar[r,"\check{\mu}^{\otimes 2}"]\ar[d,"{[\cdot,\cdot]}"'] &
				C^{\infty}(M)^{\otimes 2} \ar[d,"{[\cdot,\cdot]}"]
				\\
				\mathfrak{g} \ar[r,"\check{\mu}"] & C^{\infty}(M)
			\end{tikzcd}
			~.
		\end{displaymath}	
	\end{itemize}
	The duality, hence the "co" in the namings of the previous list, comes from the fact that $\check{\mu}$ is the dual  of $\hat{\mu}$ in the following sense
	\begin{equation}\label{eq:dualitymomaps}
		\hat{\mu}(\xi)\eval_p = \langle \hat{\mu}_p, \xi \rangle
		\qquad \forall \xi \in \mathfrak{g}, p\in M
		~.
	\end{equation}
	In other words
	\begin{displaymath}
			\hat{\mu} 
			\in 
			\Hom_{\text{smooth}}\big(M,~\Hom_{\text{vect}}(\mathfrak{g},\mathbb{R})\big)
			~,\quad
			\check{\mu} 
			\in 
			\Hom_{\text{vect}}\big(\mathfrak{g},~\Hom_{\text{smooth}}(M,\mathbb{R})\big)
	\end{displaymath}		
	come respectively from the \emph{currying} \cite{nlab:currying} of the same "evaluation" vector bundle map
	\begin{displaymath}
		\mu: M \times \mathfrak{g} \to \mathbb{R}_M
	\end{displaymath}	
	with respect to the first or second entry.
	\\
	The upshot is that a (co-)moment can exist only if the action is \emph{weakly Hamiltonian} in the sense of definition \ref{Def:WeakHamActions}.
	In that case, a weak comoment map is precisely a choice of Hamiltonian form for any fundamental vector field, hence the following diagram commutes in the category of vector spaces
	\begin{displaymath}
		\begin{tikzcd}
			& C^{\infty}(M)\ar[d,"\pi_{ham}"]\\
			\mathfrak{g} \ar[r,"\vAct"'] \ar[ur,dashed,"\check{\mu}"]& \mathfrak{X}
		\end{tikzcd}
		~.
	\end{displaymath}
	If the comoment map is "strong", the action is said to be \emph{(strongly) Hamiltonian} and the previous diagram commutes in the category of Lie algebras.
\end{reminder}

	The term "moment" comes from the following crucial examples that is prototypical in geometric mechanics:
\begin{example}[Linear and angular momenta]\label{Ex:SymplecticMechanicalMomenta}
	Given an action $\vartheta:G\action Q$, its lift $\vartheta^L:G\action M=T^\ast Q$ acts via symplectic vector fields with respect to the \emph{canonical symplectic form}.
	In particular this action preserves the tautological $1$-form $\theta$ and it can be shown to be Hamiltonian with comoment map given by
	\begin{displaymath}
		\morphism{\check{\mu}}
		{\mathfrak{g}}
		{C^{\infty}(M)}
		{\xi}
		{-\iota_{v_{\xi}}\theta}
		~.
	\end{displaymath}
	If one takes $Q=\R^3$, $Q$ and $T^* Q\cong \R^6$ can be interpreted as the \emph{configuration space} and the \emph{phase space} of a point-particle in the physical space.
	The comoment map with respect to the action of the translation or the rotation group on $Q$ gives, respectively, the linear and angular momenta of the point-particle freely moving in space \cite[\S 22.4]{CannasdaSilva2001}.
\end{example}

\medskip
In the multisymplectic context, the generalization of  the (co)momentum maps of the symplectic case leads to the more refined concept of "\momap".
\begin{definition}[\Momaps \cite{Callies2016}]\label{Def:HCMM}
	Let $v:\mathfrak g\to \mathfrak X(M)$ be a multisymplectic Lie algebra action, \ie it preserves the symplectic form $\omega \in \Omega^{n+1}(M)$.
	We call a \emph{\momap} pertaining to $v$ any $L_\infty$-morphism 
	$$
		(f)=
		\Big\{f_k:\Lambda\mathfrak{g}^k \to (L_{\infty}(M,\omega))^{1-k}\subseteq \Omega^{n-k}
		\Big\}_{k=1,...,n}
	$$
	from $\mathfrak g$ to $L_\infty(M,\omega)$ 
	satisfying the HDDW equation on the first component
	$$
		df_1(\xi)=-\iota_{v_\xi}\omega \qquad \forall \xi\in\mathfrak{g}
		~.
	$$ 
\end{definition}
\begin{notation}
	A group action $G\action (M,\omega)$ is called \emph{Hamiltonian}, if the corresponding infinitesimal Lie algebra action admits a \momap.
	We will often refer to the \momap of a certain Lie group action understanding it as the \momap corresponding to the infinitesimal action for the Lie algebra of the given group.
\end{notation}

\begin{remark}
More conceptually, a \momap is an $L_\infty$-morphism $(f):\mathfrak{g}\to L_\infty(M,\omega)$ lifting the action $v:\mathfrak{g}\to \mathfrak{X}(M)$, 
\ie making the following diagram commutative in the $L_\infty$-algebras category:
\begin{center}
\begin{tikzcd}[column sep = large]
	& L_\infty (M,\omega)\ar[d,"\pi_{\ham}"]\\
	\mathfrak{g} \ar[ur,"(f)",dashed]\ar[r,"v"'] & \mathfrak{X}(M)
\end{tikzcd}
\end{center}
where the vertical arrow $\pi_{\ham}$ is the trivial $L_\infty$-extension 
of the linear function mapping any Hamiltonian form to the unique corresponding Hamiltonian vector field given in remark \ref{Rem:ProjectiontoHamVfields}.
	\\
	Hence, comparing with remark \ref{Rem:weakHamasLift}, one can see that the existence of a \momap is a stronger condition in the sense that it requires a lift of $\vAct$ not only in the plain vector space category but in the $L_\infty$-algebra category
		\begin{displaymath}
		\begin{tikzcd}[ column sep = huge]
			&[3em]
			L_{\infty}(M,\omega)
			\ar[d,two heads,"\pi_{\ham}"]
			\\
			\mathfrak{g} \ar[rd,"\vAct"'] \ar[r,"f_1"]\ar[ru,"(f)"]		
			&
			\Omega^{n-1}_{\ham}(M,\omega)
			\ar[d,two heads,"\pi_{\ham}"]			
			\\
			&
			\mathfrak{X}_{\ham}(M,\omega)
		\end{tikzcd}
		~.
	\end{displaymath}

\end{remark}
In the following we will often make use of an explicit version of definition \ref{Def:HCMM} subsumed by the following lemma:
\begin{lemma}[\cite{Callies2016}]\label{Lem:ExplicitHCCM}
	A \momap $(f)$ for the infinitesimal	multisymplectic action	of ${\mathfrak g}$ on $M$ 
	is given explicitly by a sequence of linear maps
	\begin{displaymath}
		(f)  = \Big\{ f_i: \,\,\, \Lambda^i {\mathfrak g} \to \Omega^{n-i}(M)
		\quad \vert \quad 
		0\leq i \leq n+1  \Big\}	
	\end{displaymath}
	fulfilling a set of equations:
	\begin{equation}\label{eq:fk_hcmm}
	-f_{k-1} (\partial p) = d f_k (p) + \varsigma(k) \iota(v_p) \omega
	\end{equation}
	together with the condition
	\begin{displaymath}
		f_0 = f_{n+1} = 0 
	\end{displaymath}
	for all $p \in \Lambda^k\mathfrak{g}$ and $k=1,\dots n+1$.
	(Recall that $\partial$ denotes the Chevalley-Eilenberg boundary operator defined in equation \eqref{eq:CE_boun}
	and $\varsigma(k)$ is the sign coefficient given in equation \eqref{Eq:SigmaSign}.)
\end{lemma}
\begin{proof}
	Equation \eqref{eq:fk_hcmm} is a simple application of remark \ref{Rem:GroundedEasyAxioms} to the grounded $L_\infty$-algebra of higher observables.
\end{proof}

\begin{remark}\label{Rem:TermMuByMauro}
	\note{
		Usare $\mu_k$ può creare confusione con il capitolo higher rogers visto che è il simbolo che uso spesso per i multibrackets
	}
Formula \eqref{eq:fk_hcmm} can be read as follows: when $k=1$ it tells us that $\vAct$ acts via Hamiltonian vector fields and $f_1$ is a linear map choosing an Hamiltonian form $f_1(x)$ pertaining to the Hamiltonian vector field $v_x$ for any $x\in \mathfrak{g}$. 
Hence, $f_1$ is the choice of a primitive for the contraction of $\omega$ with any fundamental vector field of the action.
\\
When $k\geq 2$, equation \eqref{eq:fk_hcmm} can be read as the condition that the auxiliary closed differential form 
$$
\mu_k := f_{k-1} (\partial p) +  \varsigma(k) \iota(v_p) \omega
$$
must actually be {\it exact}, with potential $-f_k(p)$.
Closure of $\mu_k$ is again a consequence of lemma \ref{lemma:multicartan}  together with $\d \omega = 0$.
	Namely one has:
	\begin{equation}\label{eq:MauroMukForm}
		\begin{aligned}
			\d \mu_k 
			&=
			\d (f_{k-1} (\partial p) +  \varsigma(k) \iota(v_p) \omega)
			=
			\\
			&=
			\varsigma(k) (-1)^k \iota(v_{\partial p})\omega -
\varsigma(k-1) \iota(v_{\partial p})\omega
			=\\
			&=
			[-\varsigma(k+1) -\varsigma(k-1)] \iota(v_{\partial p})\omega = 0
			~.
		\end{aligned}
	\end{equation}
\end{remark}
The previous remark leads to the following sufficient condition to the existence of a \momap:
\begin{theorem}[{\cite[Thm. 9.6]{Callies2016}}]
	Consider a Lie algebra $\mathfrak{g}$ acting on a multisymplectic manifold $(M,\omega)$.
	The action $\vAct:\mathfrak{g}\to \mathfrak{X}(M)$ admits a \momap if it acts by Hamiltonian vector fields, \ie if there exists a function $\phi: \mathfrak{g}\to \Omega^{n-1}_{\ham}(M,\omega)$ such that $\d \phi = -\iota_{v} \omega$, and $H_{dR}^k(M)=0$ all $1 \leq k \leq n-1$.
\end{theorem}
\begin{remark}[{\cite[Rem. 9.7]{Callies2016}}]
	The previous theorem applies also when the infinitesimal Lie algebra action $\vAct$ comes from an infinite-dimensional Lie group $G$ that is locally exponential.
\end{remark}

\begin{remark}
	The reason why the term "homotopy" appears in the definition of the multisymplectic analogue of an ordinary comoment map is due to the fact that it can be interpreted as a homotopy between cochain-complexes.
	Observe that the components of $(f)$ can be arranged in the following (non-commutative) diagram inside the category of vector spaces:
	\begin{wideeq}
		\begin{tikzcd}[ampersand replacement=\&]
		\bigwedge^{n+2}\mathfrak{g} 
		\ar[r,"\partial"]
		\ar[d,blue,"\iota_{\mathfrak{g}}^{n+2}\omega",pos=0.4]
		\&
		\bigwedge^{n+1}\mathfrak{g} 
		\ar[r,"\partial"]
		\ar[dl,purple,"f_{n+1}",sloped,pos=0.4]
		\ar[d,blue,"\iota_{\mathfrak{g}}^{n+1}\omega",pos=0.4]
		\& 
		\bigwedge^{n}\mathfrak{g} 
		\ar[r] 
		\ar[d,blue,"\iota_{\mathfrak{g}}^{n}\omega",pos=0.4]	
		\ar[dl,purple,"f_{n}",sloped,pos=0.4]
		\&[-3em] 
		\cdots \ar[r]
		\&[-3em] 
		\bigwedge^{2}\mathfrak{g}
		\ar[r,"\partial"]
		\ar[d,blue,"\iota_{\mathfrak{g}}^{2}\omega",pos=0.4]
		\& 
		\bigwedge^{1}\mathfrak{g} 
		\ar[r,"\partial"]
		\ar[d,blue,"\iota_{\mathfrak{g}}^{1}\omega",pos=0.4]
		\ar[dl,purple,"f_{1}",sloped,pos=0.4]
		\&[-2em] 
		0 \ar[d,blue,hook]
		\\
		0 \ar[r,hook]
		\&
		C^{\infty}(M) \ar[r,"\d"] 
		\&
		\Omega^{1}(M) \ar[r] 
		\& 
		\cdots \ar[r]
		\&
		\Omega^{n-1}(M) \ar[r,"\d"]
		\&
		\Omega^{n}(M) \ar[r,"\d"]
		\&
		\cdots
		\end{tikzcd}
	\end{wideeq}
	The top line can be interpreted as the reversed Chevalley-Eilenberg chain complex, i.e $ \setminus CE(\mathfrak{g}) = S^{\bullet}(\mathfrak{g}[1])$ and the line below as a suitable shifted truncation of the de Rham complex, $\trunc_{n}\Omega(M)[n+1] \supset L$, including the chain complex underlying $L_\infty(M,\omega)$.
	Hence, equation \eqref{eq:fk_hcmm} express the condition that $(f)$ is a chain homotopy between the chain map $\iota_\mathfrak{g}\omega$ (see remark \ref{Rem:SignedMultiContraction}) and the zero map, \ie
	\begin{displaymath}
\begin{tikzcd}[column sep = huge]
\setminus CE(\mathfrak{g}) \arrow[r,blue,"\iota_{\mathfrak{g}}", bend left=20, ""{name=U, below}]
\arrow[r,"0"', bend right=20, ""{name=D}]
& \qquad\qquad\trunc_{n}\Omega(M)[n+1]
\arrow[Rightarrow, purple,"(f)", from=U, to=D]
\end{tikzcd}		
	\end{displaymath}
\end{remark}
\begin{remark}
	The universal property of the homotopical pullback introduced in remark \ref{Rem:FiorenzaPullback} implies that if a weakly Hamiltonian action $\vAct:\mathfrak{g}\to \mathfrak{X}_{\ham}(M)$ satisfies the property that 
	$(\iota_\mathfrak{g}):= ((\iota_{\mathfrak{X}})\circ \vAct )  $ is null-homotopic, 
	then there exists a $(f)$ (unique modulo homotopy) such that the following diagram commutes up to homotopies
		\begin{displaymath}
			\begin{tikzcd}
				\mathfrak{g}\ar[ddr,bend right=30,"\vAct"']\ar[rrd,bend left=30]\ar[dr,"\exists ! (f)"] 
				&&\\
				&L_{\infty}(M,\omega)\arrow[dr, phantom, "\scalebox{1.5}{$\lrcorner$}" , very near start, color=black] \ar[r]\ar[d,"\pi_{\ham}"']
				& 
				0 \ar[d]\\
				&\mathfrak{X}^1(M) \ar[r,"(\iota_{\mathfrak{X}})"] & \mathcal{B}	
			\end{tikzcd}
			~.
		\end{displaymath}
		The latter means that the two outer triangles commute and there exists two homotopies $(\iota_{\mathfrak{X}_{\ham}})\Rightarrow 0$ and $(\iota_{\mathfrak{g}})\Rightarrow 0$.
\end{remark}

\begin{remark}
	As noted in \cite[Rem. 5.2]{Callies2016}, one can generalize the standard characterization of the image of a comoment map as a Poisson sub-algebra in symplectic geometry also in the multisymplectic case.
	\\
	In the latter case, the image of $(f)$ has to be understood as the cochain complex $I \hookrightarrow L$ given by the following components
	\begin{displaymath}
		I^k =
		\begin{cases}
			\text{Im}(f_n)\subset  C^{\infty}(M) 
			& \text{ if } k=1-n
			\\
			\text{Im}(f_{1-k})\oplus \d\text{Im}(f_{2-k}) \subset \Omega^{n+k-1}(M)
			& \text{ if } 2-n \leq k \leq 0
		\end{cases}
		~.
	\end{displaymath}
\end{remark}

\begin{remark}[Relation with other notions of multisymplectic moment maps]\label{Rem:OtherNotionsofComoments}
	It is important to remark that definition \ref{Def:HCMM} is not the only notion of "multisymplectic" moment map proposed in the literature.
	We mention among other the so-called \emph{covariant multimoment map}\cite{Carinena1991b} (see also \cite{Gimmsy1}), the \emph{multimoment map}\cite{Madsen2012} and the \emph{weak moment map}\cite{Herman2017}.
	The relationship of these notions with the definition of \momap can be read in \cite[\S 12]{Callies2016} and \cite{Mammadova2020}.
\end{remark}	

A big part of chapters \ref{Chap:MauroPaper} and \ref{Chap:LeonidPaper} will be devoted to give new explicit constructions of \momaps. 
Several other examples can be found in \cite{Callies2016}\cite{Ryvkin2018}.

\subsection{Equivariant \momaps}
	In symplectic geometry, the $G$-equivariance condition took a central role in discerning between the weaker and stronger notion of moment map.
	When passing to the dual notion, this condition is equivalent to the Lie algebra morphism property of the comoment map.
	
	In the $n$-plectic case, the latter condition's higher analogue translates to the requirement for the \momap to be a morphism in the $L_\infty$-algebra category; this it is not equivalent to the $G$-equivariance condition anymore.
	Therefore, it is still possible to recognize among \momaps those which happen to be equivariant in the sense that all components are equivariant with respect to the coadjoint action.
	Namely, we have the following definition:

\begin{definition}[Equivariant \Momaps]\label{Def:EquivariantMomap}
	A \momap pertaining to the multisymplectic action $G\action M$ is called \emph{$G$-equivariant} if
	\begin{displaymath}
		{\mathcal L}_{v_\xi} f_k({p}) = f_k ([\xi, {p}])
		\qquad \forall \xi \in \mathfrak{g} \; ,~ p \in \Lambda^k
		\mathfrak{g},
	\end{displaymath}
where $[\xi, p]$ is the adjoint action of $\mathfrak{g}$ on $\Lambda^k \mathfrak{g}$.
Explicitly, the adjoint action reads on decomposable elements as follows: 
	\begin{equation}\label{eq:adjointactionwedge}
		[\xi, x_1\wedge\dots\wedge x_k] =
		\sum_{l=0}^k
		(-1)^{k-l} [v,x_l] 
		\wedge x_k \wedge \dots \wedge \hat{x_l} \wedge \dots \wedge 	x_1
		~.
	\end{equation}
\end{definition}
An example of equivariant \momap is given by the following example generalizing example \ref{Ex:SymplecticMechanicalMomenta}:
\begin{example}[Fields mechanical momentum]
	Consider a smooth action $\vartheta^Q : G\action Q$ on a given a smooth manifold $Q$. 
	There is a natural lift for this action to $M=\Lambda^n Q$ given by the bundle map
	\begin{displaymath}
		\morphism{\vartheta^M_g}
		{\Lambda^n T^\ast Q}
		{\Lambda^n T^\ast Q}
		{(q,\alpha)}
		{(\vartheta^Q_g(q),((T_q~\vartheta^Q_g)^{-1})^\ast \alpha )}
	\end{displaymath}
	over $\vartheta^Q_g$ for any $g\in G$.
	\\
	This action preserves the tautological $n$-form $\theta$. Being then $\theta$ a $G$-invariant primitive of $\omega = \d \theta$, this action admits a \momap given by \cite[lem 8.1]{Callies2016}](see also lemma \ref{lem:extexact} in the following chapter).
\end{example}

\subsection{Induced \momaps and isotropy subgroups}
In this subsection we will discuss how \momaps behave under restriction to subgroups and invariant submanifolds. 
This will be useful for constructing an $SO(n)$-comoment for $S^n$ in chapter \ref{Chap:LeonidPaper}. 
Let $G$ be a Lie group with Lie algebra $\mathfrak{g}$, acting on a pre-$n$-plectic manifold $(M,\omega)$ with \momap $(f)\colon \mathfrak{g}\to L_{\infty}(M,\omega)$. 
One can obtain new actions either restricting to a Lie subgroup of $G$ or restricting to an invariant submanifold of $(M,\omega)$. 
\begin{proposition}[Lemma 3.1 in \cite{Shahbazi2016}]\label{prop:restrict}
	Let $H\subset G$ be a Lie subgroup, and denote by $j\colon \mathfrak{h}\hookrightarrow \mathfrak{g}$ the corresponding Lie algebra inclusion. 
	The restricted action of $H$ on $(M,\omega)$ has \momap $(f \circ j): \mathfrak{h}\to L_{\infty}(M,\omega)$, given in components by $f_i\circ j:\Lambda^i\mathfrak{h}\to \Omega^{n-i}(M)$ for $i=1,\dots,n$.
\end{proposition}
\begin{center}
	\begin{tikzcd}
		H \ar[symbol=\action]{r} \ar[hook]{d}& M \ar[equal]{d}& & \mathfrak{h} \ar[hook,"j"']{d} \ar[dashed]{dr} & \\
		G \ar[symbol=\action]{r}& M & & \mathfrak{g}\ar["(f)"']{r} & L_{\infty}(M,\omega) 
	\end{tikzcd}
\end{center}
\begin{proposition}[Lemma 3.2 in \cite{Shahbazi2016}]
	\label{prop:NenMGinvariant} 
	Let $N\overset{i}{\hookrightarrow} M$ be a $G$-invariant submanifold of $M$.
	\\
	Then the action $G\action \left(N,i^{\ast}\omega\right)$ has \momap $(i^{\ast} \circ f): \mathfrak{g}\to L_{\infty}\left(N,i^{\ast}\omega\right)$, 
	given in components by $i^*\circ f_i:\Lambda^i\mathfrak{h}\to \Omega^{n-i}(N)$ for $i=1,\dots,n$.
\end{proposition}
\begin{center}
	\begin{tikzcd}
		G \ar[symbol=\action]{r} \ar[equal]{d}& N \ar[hook,"i"]{d}& & & L_{\infty}(N, i^\ast\omega) \\
		G \ar[symbol=\action]{r}& M & & \mathfrak{g}\ar["(f)"']{r}\ar[dashed]{ur} & L_{\infty}(M, \omega) \ar["i^\ast"']{u} 
	\end{tikzcd}
\end{center}

Moreover, we can produce a new \momap by considering a different multisymplectic form obtained contracting $\omega$ with cycles in the Lie algebra homology of $\mathfrak{g}$, \ie elements of the vector space
\begin{equation}
	Z_k(\mathfrak{g}) = \{ p \in \Lambda^k \mathfrak{g} \; \vert: \partial p = 0 \}
	\quad.
\end{equation}
\begin{proposition}[Proposition 3.8 in \cite{Ryvkin2016}]\label{prop:indmomap}
	Let $p\in Z_{k}(\mathfrak{g})$, for some $k \geq 1$,
	denote by $G_{p}$ the corresponding isotropy group for the adjoint action of $G$ on $\Lambda^{k}\mathfrak{g}$, and by $\mathfrak{g}_p=\{x\in \mathfrak{g}: [x,p]=0\}$ its Lie algebra.
	\\
	If $G_p^0$ is the connected component of the identity in $G_p$, then the action $G_p^0\action \left(M,\iota(v_p)\omega\right)$ admits a comoment
	$(f^p) \colon \mathfrak{g}_p \to L_{\infty}(M,\iota(v_p)\omega)$ with components $(i=1,\dots,n-k)$:
	\begin{displaymath}
		\morphism{f^p_i}
		{\Lambda^i\mathfrak{g}_p}
		{\Omega^{n-k-i}(M)}
		{q}
		{\varsigma(k)~f_{i+k}(q\wedge p)}
		~.
	\end{displaymath}
\end{proposition}
\begin{remark} In the context of multisymplectic geometry, ``weak comoment maps'' (\cf \cite{Herman2018}), and ``multimoments'' (\cf \cite{Madsen2013} and remark \ref{Rem:OtherNotionsofComoments}),  the subspace $ Z_{k}(\mathfrak{g})$ is often referred to as ``the $k$-th Lie kernel''.
\end{remark}
\begin{proposition}[Remark 3.9 in \cite{Ryvkin2016}]\label{prop:indmomap_equiv} 
	If the \momap $(f) \colon \mathfrak{g} \to L_{\infty}(M,\omega)$ is also $G$-equivariant, then the map $(f^p)$ defined in proposition \ref{prop:indmomap} is $G_p^0$-equivariant.\\
	Another equivariant \momap
	for the action of $G_p^0$ on $(M,\iota(v_p)\omega)$ is given in components $(i=1,\dots,n-k)$ by:
	\begin{displaymath}
		\morphism{f^p_i}
		{\Lambda^i\mathfrak{g}_p}
		{\Omega^{n-k-i}(M)}
		{q}
		{(-1)^{k}\iota(v_q)(f_{k}(p))}
	\end{displaymath}
	which, in general, may differ from the one given in Proposition \ref{prop:indmomap}. 
\end{proposition}
\noindent

\vspace{1em}
Provided certain conditions are met, it is possible to induce a comoment on an invariant submanifold of $M$ even if the obvious pullback vanishes (\eg~when  $\omega$ is a top dimensional form).
In chapter \ref{Chap:LeonidPaper} we will make use of the following corollary subsuming the contents of the previous propositions:
\begin{corollary}\label{cor:inducedmachinery}
	Let $G\action(M,\omega)$ be a multisymplectic group action.
	If there exists:
	\begin{itemize}[topsep=1pt,itemsep=0pt, partopsep=1pt]
		\item another multisymplectic manifold $(N,\eta)$ containing $M$ as a $G$-invariant embedding $j:M\hookrightarrow N$;
		\item a Lie group $H \supset G$ containing $G$ as a Lie subgroup;
		\item a multisymplectic action $H \action (N,\eta)$ with equivariant \momap $s:\mathfrak{h}\to L_\infty(N,\eta)$;
		\item an element $p\in Z_k(\mathfrak{h})$ in the Lie kernel of $\mathfrak{h}$ such that $G \subset H_p$ and $\omega = j^\ast \iota_p \eta$;
	\end{itemize}
	then the action $G\action(M,\omega)$ 	admits an equivariant \momap,	given in components $(i=1,\dots,n-k)$ by:
	\begin{displaymath}
		\morphism{f_i}
		{\Lambda^i\mathfrak{g}}
		{\Omega^{n-k-i}(M)}
		{q}
		{\displaystyle(-1)^{k}j^\ast\left(\iota(v_q)(s_{k}(p))\right)}
		~.
	\end{displaymath}
\end{corollary}
\begin{proof}
	Starting from the given comoment $(s)$ it is possible to construct another comoment $(s^p)$ resorting to proposition \ref{prop:indmomap_equiv}.
	The sought \momap descends from $(s^p)$ via the consecutive application of propositions \ref{prop:restrict} and \ref{prop:NenMGinvariant}
	\begin{center}
	\begin{tikzcd}[column sep=large]
		H \ar[symbol=\action]{r} & (N,\eta)
		& & \mathfrak{h} \ar["(s)"]{r} & L_{\infty}(N,\eta)
		\\
		H_p \ar[symbol=\action]{r} \ar[hook]{u}& (N,\iota_{v_p}\eta)  
		& & \mathfrak{h}_p\ar[hook]{u} \ar[dashed,"(s^p)"]{r}[swap]{Prop. \ref{prop:indmomap_equiv}} & L_{\infty}(N,\iota_{v_p}\eta) \ar[equal]{d} 
		\\
		G \ar[symbol=\action]{r} \ar[hook]{u}& (N,\iota_{v_p}\eta)   \ar[equal]{u}
		& & \mathfrak{h}\ar[hook]{u} \ar[dashed,"Prop. \ref{prop:restrict}"']{r} & L_{\infty}(N,\iota_{v_p}\eta) \ar["j^\ast"]{d} 
		\\
		G \ar[symbol=\action]{r} \ar[equal]{u}& (M,\omega =j^\ast \iota_{v_p}\eta) \ar[hook]{u}  
		& & \mathfrak{g}\ar[equal]{u} \ar[dashed,"Prop. \ref{prop:NenMGinvariant}"']{r}& L_{\infty}(M,\omega)  
	\end{tikzcd}
	\end{center}
	together with the observation that if the starting \momap $(s)$ is equivariant then the induced maps are such.
\end{proof}

\subsection{Conserved quantities along \momaps}
Recall that one of the primary motivation for introducing the notion of moment map in symplectic geometry was to characterize the generators of symmetries of Hamiltonian systems in terms of conservation laws.

Let us take as a generalization of an Hamiltonian system in multisymplectic geometry the triple composed of a smooth manifold $M$ together with a $n$-plectic form $\omega$ and a fixed Hamiltonian form $H\in \Omega^{n-1}_{\ham}(M,\omega)$. 
One could wonder if the momenta, \ie the image of the \momap with respect to a certain action, are preserved along the flow determined by $H$. 
The answer to this question is subsumed by the following theorem
\begin{theorem}[{\cite{Ryvkin2016}}]
	Let be $(M,\omega)$ be a multisymplectic manifold and
	consider a multisymplectic action $\vAct:\mathfrak{g}\to (M,\omega)$	admitting a \momap $(f):\mathfrak{g}\to L_{\infty}(M,\omega)$.
	Let be $H$ an Hamiltonian $(n-1)$-form that is locally (resp. globally or strictly) preserved along all the fundamental vector fields of $\vAct$.
	The image of the component $f_k$ of the \momap are preserved in the following sense (compare with definition \ref{Def:conservedQuantities}):
	
	\smallskip
	\begin{adjustbox}{center}
		\begin{tabular}{| l | p{9em}| p{9em} | p{9em} | }
 			\hline
 			& $H$ locally preserved along $\vAct$& $H$ globally preserved along $\vAct$& $H$ strictly preserved along $\vAct$ \\
			\hline
 			$f_k(Z_k(\mathfrak{g}))$ & locally conserved & locally conserved & globally conserved \\ 
 			$f_k(B_k(\mathfrak{g}))$ & globally conserved & globally conserved & globally conserved \\ 
 			\hline
		\end{tabular}
 \end{adjustbox}
	\smallskip
	\\ 
	where $Z_k(\mathfrak{g})$ and $B_k(\mathfrak{g})$ denote respectively $k$-cycles and $k$-boundaries in the Chevalley-Eilenberg cohomology of $\mathfrak{g}$.
\end{theorem}
Note that conservation, in the sense of definition \ref{Def:conservedQuantities}, of momenta (\ie elements in the image of $(f)$) is only assured in the images of cocycles of the Chevalley-Eilenberg chain complex.
This can be seen as the multisymplectic analogue of Noether's theorem  as explained in \cite{Herman2017} employing the languages of \emph{weak homotopy moment map}.

\section{Gauge transformations}\label{Section:GaugeTransformations}
In chapter \ref{Chap:MarcoPaper}, we will discuss the possible relationship existing between the $L_\infty$-algebras of observables pertaining to two multisymplectic forms differing by an exact form.
Let us first fix the following nomenclature:
\begin{definition}[$B$-related closed form]
	Given two closed form $\tilde{\omega},~\omega \in \Omega^{n+1}(M)$, we call them \emph{$B$-related} or \emph{gauge related} if there exists $B\in\Omega^n(M)$ such that
	\begin{displaymath}
		\tilde{\omega}= \omega + \d B
\end{displaymath}	
\ie if the two closed forms sit in the same de Rham cohomology class with difference given by $\d B$.
\end{definition}

Consider two $B$-related $n$-plectic forms $\omega$ and $\tilde{\omega}$ on $M$ and a Lie algebra action 
$v:\mathfrak{g} \to \mathfrak{X}(M)$ which is multisymplectic with respect to $\omega$, one has the following
\begin{lemma}
	The infinitesimal action $v$ is multisymplectic with respect to both $\omega$ and $\tilde{\omega}$ if and only if $B$ is locally conserved along $v$.
\end{lemma}
\begin{proof}
	Being the action multisymplectic with respect to both $\omega$ and $\tilde{\omega}$ means that
	\begin{displaymath}
		0 =
		\mathcal{L}_{v_\xi} \tilde{\omega} =
		\mathcal{L}_{v_\xi} \omega - \mathcal{L}_{v_\xi} \dd B 
	\end{displaymath}
	therefore $\dd \mathcal{L}_{v_\xi}B = 0$ for all $\xi \in \mathfrak{g}$.
\end{proof}

Consider now the case that $v$ admits a \momap $(f):\mathfrak{g}\to L_\infty(M,\omega)$ with respect to $\omega$,
\begin{lemma}
	A necessary condition for the existence of a \momap pertaining to $v$ with respect to both the symplectic forms $\omega$ and $\tilde{\omega}$ is that $B$ is globally conserved.
\end{lemma}
\begin{proof}
	Let $\tilde{f}:\mathfrak{g}\to L_\infty(M,\tilde{\omega})$ be an \Hcmm with respect to $\tilde{\omega}$, one has
	\begin{displaymath}
		\dd \tilde{f}_1(\xi) =
		-\iota_{v_\xi}\tilde{\omega}=
		-\iota_{v_\xi}\omega + \iota_{v_\xi}\dd B= 0
	\qquad \forall \xi \in \mathfrak{g}
	~,	
	\end{displaymath}
	therefore, $\iota_{v_\xi}\dd B$ must be exact. 
	A primitive $h$ of the latter must satisfy the following equation
	\begin{displaymath}
		\dd h = \iota_{v_\xi} \dd B = \mathcal{L}_{v_\xi} B - \dd \iota_{v_\xi}B
		~,
	\end{displaymath}
	which means that $\mathcal{L}_{v_\xi} B $ has to be exact for all $\xi \in \mathfrak{g}$.
\end{proof}
\begin{remark}
	Observe that this is a simpler instance of the notion of equivalence between \momap proposed in \cite{Fregier2015}.
\end{remark}
	When considering gauge-related multisymplectic structures,   the following lemma holds,  see   \cite[Beginning of \S7.2]{Fregier2015}.
\begin{lemma}[Gauge transformation of \Hcmm]\label{lem:momaps}
Let the infinitesimal action $v:\mathfrak{g}\to \mathfrak{X}(M)$ 
	{preserve the $n$-symplectic form $\omega$ and admit a homotopy comoment map $(f):\mathfrak{g}\to L_\infty(M,\omega)$. 
	Suppose that $B\in \Omega^n(M)$ is strictly conserved. 
	Then $\widetilde{\omega}=\omega + dB$, which we assume to be $n$-plectic, is also preserved and}
  admits a homotopy comoment map $(\widetilde{f}):\mathfrak{g}\to L_\infty(M,\widetilde{\omega})$, with components 
	\begin{displaymath}
		\widetilde{f}_k = (f_k +\mathsf{b}_k)
		~: {\wedge}^k \mathfrak{g}\to L^{1-k}\subseteq \Omega^{n-k}(M)
		~.
	\end{displaymath}
	where 
	\begin{displaymath}
		\morphism{\mathsf{b}_k:=(-)^k \iota^k_{\mathfrak{g}}B}
		{\Lambda^k\mathfrak{g}}
		{\Omega^{n-k}}
		{x_1\wedge\dots\wedge x_k}
		{- \varsigma(k+1) \iota(x_1\wedge\dots\wedge x_k) B}
		~.
	\end{displaymath}	
\end{lemma}
\begin{proof}

	Fix $p\in \Lambda^k\mathfrak g$, then:
	\begin{displaymath}
		\begin{aligned}
			\dd \tilde{f}_k(p) ~+&~ \tilde{f}_{k-1}(\partial p ) =
			\\
			\equal{}&
			\dd f_k(p) + f_{k-1}(\partial p )	 
			-\varsigma(k)\dd \iota_{v_p} B
			-\varsigma(k-1) \iota_{\partial v_p} B
			=\\
			\equal{Lem. \ref{Lem:ExplicitHCCM}}&
			- \varsigma(k)\iota_{v_p}\omega
			-\varsigma(k)\dd \iota_{v_p} B
			-\varsigma(k-1) \iota_{\partial v_p} B
		\end{aligned}
	\end{displaymath}
	Plugging in the (multi)Cartan formula (Lemma \ref{lemma:multicartan}) and using the hypothesis that $\mathcal{L}_{v_\xi} B = 0$ we conclude that
	\begin{displaymath}
		\begin{aligned}
			\dd \tilde{f}_k(p) +& \tilde{f}_{k-1}(\partial p )
			=
			\\
			=&~
			- \varsigma(k)\iota_{v_p}\omega
			-[\varsigma(k+1)(-)^k + \varsigma(k)] \iota_{\partial p} B
			-\varsigma(k+1)(-)^k \iota_{v_p}\dd B
			=\\
			=&~
			- \varsigma(k)\iota_{v_p}\omega
			+\varsigma(k) \iota_{v_p}\dd B
			=\\
			=&~
			- \varsigma(k)\iota_{v_p}\tilde{\omega}
		\end{aligned}
		~.
	\end{displaymath}
	The equality $\varsigma(k+1)(-)^k =- \varsigma(k)$ follows from the definition of $\varsigma(k)$ (see lemma \ref{lem:varsigmasignprops}).
\end{proof}
 This last lemma will have a key role in \ref{Chap:MarcoPaper} (see section \ref{Sec:DiagramGaugeTransf}) where the operator $\mathsf{b}_k$ will be recast in terms of the \emph{pairing operator} mentioned in remark \ref{Rem:SignedMultiContraction}.

\section{Transgression}
We already mentioned several times how multisymplectic geometry originated from the long-standing problem of identifying the apt geometric formulation to encode field-like mechanical systems.
Other approaches to the geometric mechanics of infinite-dimensional systems (in sense of being parametrized by a continuous infinity of degrees of freedom) involve formal methods on infinite-dimensional manifolds.
\footnote{We use the term "formal" because the problem of formalizing the concept of smoothness in infinite-dimensional spaces is particularly difficult. 
Above all, there is none completely general and uniquely recognized definition of smooth infinite-dimensional space (see \cite{Stacey2008} for a rundown on some of these possible notions).}
Further inputs coming from physics allow a slightly finer characterization of the infinite-dimensional space to be considered.
\\
The key idea is that the configuration space of such systems consists of smooth sections of a certain smooth fibre bundle called "configuration bundle" (the points of the base are the parameters and the fibres consist of all the values admissible by the field).
Under certain condition, it is also possible to equip this infinite-dimensional configuration space with a canonical symplectic structure (see \cite{Miti2015} and references therein for a slightly more complete account).
In other terms, configurations of a field theory are encoded by a (weakly) symplectic mapping space.
\\
In many practical cases, the field configurations are mappings into a multisymplectic manifold (see e.g \cite{Forger2005}).
The notion of \emph{transgression on mapping space} makes possible to relate the geometry of (finite-dimensional) multisymplectic manifold to the geometry of infinite-dimensional symplectic manifolds.

	Consider a smooth $m$-dimensional manifold $M$ together with a smooth, compact $s$-dimensional manifold $\Sigma$. 
	Let us denote by $M^{\Sigma}$ the mapping space of smooth functions from $\Sigma$ to $M$:
	\begin{displaymath}
			M^\Sigma := \hom_{smooth}(\Sigma, M)	
			~.
	\end{displaymath}	
	The space $M^\Sigma$ carries a Fr\'echet manifold structure (see \cite{Kri-Mich})
	and the tangent space on a given smooth function $\gamma\in M^\Sigma$ is given by the space of sections of the pullback bundle of $TM$ along $\gamma$, \ie :
	\begin{displaymath}
		T_\gamma( M^\Sigma) = \Gamma(\gamma^\ast (TM))
		~.
\end{displaymath}	
	In particular, there is a well-defined notion of differential form over $M^\Sigma$ giving on any submanifold $\gamma$ a $C^{\infty}(\Sigma)$-multilinear $\R$-valued form on $T_\gamma (M^\Sigma)$.
	\note{check $C^{\infty}(\Sigma)$-multilinear}
	One can induce a differential form on $M^\Sigma$ from any differential form on $M$:
	\begin{definition}[Transgression]
		We call \emph{transgression} the graded map
		\begin{displaymath}
			(-)^\ell: \Omega(M) \to \Omega(M^\Sigma)[-s]
		\end{displaymath}	
		defined on any $\alpha\in \Omega^n(M)$, $\gamma \in M^{\Sigma}$ and $v_1,\dots v_{n-s}\in T_{\gamma}(M^\Sigma)$, by
		\begin{displaymath}
			\alpha^\ell \eval_{\gamma}(v_1,\dots,v_{n-s})= \int_{\Sigma} \gamma^\ast \left( \iota_{v_1\wedge\dots\wedge v_{n-s}}\alpha\eval_\gamma \right)
			~.
		\end{displaymath}
		Alternatively, one can see the transgression map as the composition of the pullback along the  	\emph{evaluation map}
	\begin{displaymath}
		\morphism{ev}
		{\Sigma\times M^\Sigma}
		{M}
		{(x,\gamma)}
		{\gamma(x)}
		~.
	\end{displaymath}
		with the operation of integration along the fibres (of the trivial fibration $\Sigma\times M^\Sigma \to M^\Sigma$):
		\begin{displaymath}
			\int_\Sigma ~: \Omega({\Sigma\times M^\Sigma}) \to \Omega(M^\Sigma) ~.
		\end{displaymath}
		See, for instance, \cite[\S 3.7]{Bry}.
	\end{definition}	

	\begin{example}[Transgression on loop spaces]\label{Ex:TransgressionOnLoops}
		Let $M$ be a smooth oriented manifold and $LM = C^\infty(S^{1},M)	$ be the \emph{Free loop space} (see \cite{Bry} and Remark \ref{Rem:BryLoopSpaces}).
		The transgression on the loop space is given by a	 degree $-1$ chain map
				\[
						\begin{tikzcd}[column sep= small,row sep=0ex]
					    \ell \colon& \Omega^{\bullet}(M) 	\arrow[r]& \Omega^{\bullet -1}(LM) \\
			  		  & \alpha (\textvisiblespace)\arrow[r, mapsto]& 	\left.\alpha^{\ell} \right\vert_{\gamma} =
			    		\left.{\displaystyle\int^{2\pi}_{0} \iota_{\dot{\gamma}}\alpha(\textvisiblespace)} \right\vert_{\gamma(s)} ~ ds
						\end{tikzcd}	
					\]		
	\end{example}				
	Since the transgression commutes with the de Rham differential one has that $n$-plectic forms on $M$ transgress to pre-$(n-s)$-plectic forms on $M^\Sigma$.
		\note{questo risultatino sarebbe stato bello metterlo, ma salto per motivi di tempo}

	A particularly interesting property of \momaps, already noticed in \cite[\S 6]{Fregier2015}  and \cite[\S 11]{Callies2016}, is that	$G$-actions on pre-$2$-plectic manifolds $(M,\omega)$ can be transgressed to ordinary moment maps on the associated pre-symplectic loop space $LM:=M^{S^1}$.
	\\
	A similar result holds in the case of more general mapping spaces.
	Observe first  that any smooth action $\vartheta:G\action M$ can be lifted to an action $\vartheta^\Sigma:G\action M^\Sigma$ given point-wise by
		\begin{displaymath}
			\morphism{\vartheta^{\Sigma}}
			{G\times M^\Sigma}
			{M^\Sigma}
			{(g,\gamma)}
			{\left(x\mapsto \vartheta(g,\gamma(x))\right)}
			~.
		\end{displaymath}
		The key results is that, when $\vartheta:G\action (M,\omega)$ admits a \momap, the same holds for $\vartheta^\Sigma: G \action (M^\Sigma,\omega^\ell)$
\begin{theorem}[{\cite[Corollary 6.3]{Callies2016}}]
	Consider a pre-$n$-plectic manifold $(M,\omega)$ and let
	$	f:\mathfrak{g}\to L_\infty(M,\omega) $	 
	be a \momap with respect to the action $\vartheta:G\action (M,\omega)$.
	\\
	Then the $L_\infty$-morphism $f^\ell: \mathfrak{g}\to L_{\infty}(M,\omega^\ell)$ with components given by
	\begin{displaymath}
		(f^\ell)_k = (f_k)^\ell : \Lambda^k \mathfrak{g} \to \Omega^{n-s-k}(M^\Sigma)	
	\end{displaymath}		
	is a \momap with respect to the action $\vartheta^{\Sigma}:G\action (M^\Sigma,\omega^\ell)$.
\end{theorem}
This last result provides further hints on how the higher symplectic geometry on $M$ can interact with the ordinary geometry on $M^{\Sigma}$.
		We will return to the topic of transgression on loop spaces in chapter \ref{Chap:MauroPaper}.

\ifstandalone
	\bibliographystyle{../../hep} 
	\bibliography{../../mypapers,../../websites,../../biblio-tidy}
\fi

\cleardoublepage


%% file: chapters/compactactionspheres/compactactionspheres.tex
\chapter{Action of compact groups on multisymplectic manifolds}\label{Chap:LeonidPaper}
%
%
\note{
	\textbf{Abstract:}
	We investigate the existence of \momaps (comoments) for high-dimensional spheres seen as multisymplectic manifolds.
	Especially, we solve the existence problem for compact effective group actions on spheres and provide explicit constructions for such comoments in interesting particular cases.
}
In this chapter, based on \cite{Miti2019}, we focus our attention to \momap pertaining to multisymplectic actions of compact groups.
More specifically, 
we solve the existence problem for compact effective group actions on oriented high-dimensional spheres seen as multisymplectic manifolds and provide explicit constructions of some relevant cases.
The main result can be subsumed by the following important theorem:
\begin{theorem}\label{thm:Leomainresult} (Proposition \ref{prop:intransitive} and Theorem \ref{thm:surprise})
	Let $G$ be a compact Lie group acting multisymplectically and effectively on the $n$-dimensional sphere $S^n$ equipped with the standard volume form, then the action admits a \momap if and only if $n$ is even or the action is not transitive. 
\end{theorem}

The outline of the chapter is as follows:
in the first section we briefly survey the theory of cohomological obstructions to the existence of \momaps as introduced in \cite{Callies2016} and further developed in \cite{Fregier2015,zbMATH06448534}. 
We include proofs for some known results in order to achieve a complete and self-contained exposition. 
The main novelty in this section is an intrinsic proof of Theorem \ref{cpteqcom}, which does not depend on the choice of a model for equivariant cohomology. 
\\
We then prove theorem \ref{thm:Leomainresult} for the non-transitive case in Section \ref{secintrans} and the transitive case in Section \ref{sectra}. 
In addition to proving the abstract theorem, we give explicit constructions for important classes of group actions and highlight interesting phenomena. 

\begin{remark}[Right-actions convention]\label{rem:RightActionMess}
	In this chapter, all groups considered will be connected unless stated otherwise, all actions will be smooth and on the right. A right smooth action $\vartheta:G\action M$ will be encoded by the smooth map $\vartheta: G\times M \to M$, \ie with the acting group $G$ put as the first multiplicand of the scalar product (in place of the more natural choice  $\vartheta: M\times G \to M$), in order to agree with the sign conventions usually found in the literature.
	The corresponding infinitesimal action (via fundamental vector fields) will be given by the Lie algebra morphism
$\vAct: \mathfrak{g}\to \mathfrak{X}(M)$ defined as
	\begin{displaymath}
		\vAct_\xi(m) = \left.\dfrac{d}{dt}\right\vert_0	\vartheta(m,\exp( t\xi)) \qquad \forall m \in M , \xi \in \mathfrak{g}
		~;
	\end{displaymath}			
	(\cf equation \eqref{eq:LeftFundVF}, notice the opposite sign).
	\\
	Recall that preferring right action rather than left action is mostly a matter of taste.
	Furthermore, for any given left action $\lambda:G\action M$, the smooth action $\rho:G\action M$, defined as $\rho_g = \lambda_{g^{-1}}$ for any $g\in G$, is a right action.	
\end{remark}

\section{Cohomological obstructions to \Hcmm for compact groups}
In this section, we give a short introduction to the cohomological obstructions for \momaps focusing on their geometric description.
Most results can be found in the literature  \cite{Rogers2010,Callies2016,Fregier2015,zbMATH06448534}
but we present some of the proofs for a clearer and self-contained exposition.
Our main contribution is an intrinsic proof of Theorem \ref{cpteqcom} which does not depend on a choice of model for the equivariant cohomology.

Consider an infinitesimal action of $\mathfrak{g}$ on the pre-$n$-plectic manifold $(M,\omega)$ preserving the form $\omega$.
As shown in \cite{Fregier2015, zbMATH06448534}, \momaps for this action can be interpreted as a primitive of a certain cocycle in a cochain complex.

\begin{definition}\label{def:comomentbicomplex} 
The bi-complex naturally associated to the action of $\mathfrak{g}$ on $M$ is defined by
	\begin{displaymath}
		(C_\mathfrak{g}^{\bullet,\bullet} = \Lambda^{\geq 1} 
		\mathfrak{g}^*\otimes \Omega^\bullet(M), \delta_\text{CE},d),	
	\end{displaymath}
	where $d$ denotes the de Rham differential and $\delta_{CE}:\Lambda^k\mathfrak g^*\to \Lambda^{k+1}\mathfrak g^*$ the Lie algebra cohomology differential (see reminder \ref{Rem:CEconventions}), defined on generators by
	\begin{displaymath}
	\delta_{CE}(f)(\xi_1,...,\xi_k)=\sum_{i<j}(-1)^{i+j}f([\xi_i,\xi_j],\xi_1,\dots,\hat{\xi}_i,\dots,\hat{\xi}_j,\dots,\xi_{k})
	~.
	\end{displaymath}
	The  corresponding total complex is given by
	\begin{displaymath}
		(C_\mathfrak{g}^{\bullet}, d_{\tot} = 
		\delta_\text{CE}\otimes \id + \id\otimes d),
	\end{displaymath}
where, according to the Koszul sign convention, $d_{\tot}$ acts on an element of $\Lambda^k \mathfrak{g}^*\otimes \Omega^\bullet(M)$ as $\delta_\text{CE} + (-1)^k d$.
\end{definition}

\begin{theorem}[Proposition 2.5 in \cite{Fregier2015}, Lemma 3.3 in \cite{zbMATH06448534}]\label{Thm:CoboundaryinCGRANDE}
Let $(M,\omega)$ be a pre-$n$-plectic manifold and $v:\mathfrak g\to \mathfrak X(M)$ be an infinitesimal multisymplectic action. 
The primitives of the natural cocycle
	\begin{equation}\label{eq:omegatildeobstruction}
		\tilde{\omega} = \sum_{k=1}^{n+1} (-1)^{k-1} \iota^k_\mathfrak{g} \omega \in C_\mathfrak{g}^{n+1},
	\end{equation}
	where, as already introduced in remark \ref{Rem:SignedMultiContraction},
	\begin{displaymath}
		\morphism{\iota^k_\mathfrak{g}}
		{\Omega^\bullet(M)}
		{\Lambda^k \mathfrak{g}^\ast \otimes \Omega^{\bullet-k}(M)}
		{\omega}
		{\left(p \mapsto \iota_{v_p} \omega  \right)}
		~,
	\end{displaymath}
\noindent are in one-to-one correspondence with comoments of $v$.
In particular, a \comoment exists if and only if~$[\tilde{\omega}]=0\in H^{n+1}(C_\mathfrak g^\bullet,d_ {tot})$.
\end{theorem}
\note{
	Controllare che la vecchia notazione $\omega_k$ non viene alla fine utilizzata: 
	\begin{align*}
		\iota^k_\mathfrak{g} \colon \Omega^\bullet(M)
		&\to \Lambda^k \mathfrak{g}^\ast \otimes \Omega^{\bullet-k}(M)
		\\ \omega&\mapsto \omega_k = 
		\left(p \mapsto \iota_{v_p} \omega  \right) ,
	\end{align*}
}

\begin{proof}  
Being linear maps, the components $f_k$ can be regarded as elements of a vector space
	\begin{displaymath}
		f_k \in \Lambda^k \mathfrak{g}^\ast \otimes \Omega^{n-k}(M)\cong \text{Hom}_{\text{Vect}}(\Lambda^k \mathfrak{g}, \Omega^{n-k}(M))
	\end{displaymath}
	satisfying equation (\ref{eq:fk_hcmm})
	or, equivalently, as vectors $\tilde{f}_k =\varsigma(k) f_k$ satisfying
	\begin{equation}\label{eq:fk_hcmm_tilde}
		\tilde{f}_{k-1}(\partial p ) + (-1)^k d \tilde{f}_k ( p) = (-1)^{k-1}\iota(v_p) \omega .
	\end{equation}
	The last equation is obtained multiplying equation (\ref{eq:fk_hcmm}) by the sign factor 
	$\varsigma{(k-1)}$ and noting that $\varsigma{(k-1)}\varsigma{(k)} = (-1)^k$.
	\\
	Upon considering the direct sum of these vectors
	\begin{displaymath}
		\tilde{f} = \left(\sum_{k=1}^n \tilde{f}_k \right) \in 
		\bigoplus_{k=1}^n \left(\Lambda^k \mathfrak{g}^\ast \otimes \Omega^{n-k}(M)\right)
	\end{displaymath}
	equation (\ref{eq:fk_hcmm_tilde}) can be recast as:
	\begin{equation}\label{eq:fk_hcmm_tilde_complex_1}
		\left[\delta_{\text{CE}}\otimes \id + \id\otimes d \right] \tilde f =
		\sum_{k=1}^{n+1} (-1)^{k-1} \iota^k_\mathfrak{g} \omega
	\end{equation}
	where $\iota^k_\mathfrak{g}$ is the operator defined above.
	Note that we are implicitly using the Koszul convention, therefore the action of $\id\otimes d$ on a homogeneous element $f_k \in \Lambda^k \mathfrak{g}^*\otimes \Omega^\bullet(M)$ yields $(\id\otimes d) f_k = (-1)^k (\id \otimes d f_k) $.
	\\
	If we set
		\begin{displaymath}
			\tilde{\omega} = \sum_{k=1}^{n+1} (-1)^{k-1} \iota^k_\mathfrak{g} \omega \in C_\mathfrak{g}^{n+1},
		\end{displaymath}
	equation (\ref{eq:fk_hcmm_tilde_complex_1}) corresponds to
	\begin{equation}\label{eq:fk_hcmm_tilde_complex_2}
		d_{\tot} \tilde{f} = \tilde{\omega}
	\end{equation}
	which is exactly the condition of $\tilde{f}$ being a primitive of $\tilde{\omega}$.
	\\
	It follows from Lemma \ref{lemma:multicartan} that $d_{tot}\tilde \omega=0$ for all actions preserving $\omega$. 
	Therefore the vanishing of the cohomology class 
	$[\tilde{\omega}] \in H^{n+1}(C_\mathfrak{g}^\bullet)$
	is a necessary and sufficient condition for the existence of a \comoment for the infinitesimal action of $\mathfrak{g}$ on $M$.
\end{proof}
The upshot is that \momaps pertaining to a certain multisymplectic action $\vAct:\mathfrak{g}\to \mathfrak{X}_{\msy}(M,\omega)$ are in a one-to-one correspondence with the primitives of the cochain $\tilde{\omega}$, hence they form an affine space modelled on $B^{n+1}(C_\mathfrak{g})$.

\begin{remark}\label{rk:kuenneth}
By the K\"{u}nneth theorem, the cohomology $H^\bullet(C_\mathfrak{g}^\bullet)$ is isomorphic to $H^{\geq 1}(\mathfrak g)\otimes H^\bullet(M)$, we will give a geometric interpretation to this fact in the next section.
\end{remark}	

\begin{remark}\label{rk:cp_obsruction}
In Chapter \ref{Chap:MauroPaper} we will consider a slightly different obstruction to the existence of a \comoment for $v:\mathfrak{g}\to \mathfrak{X}(M,\omega)$.
Namely, we are going to consider the following cocycle in the Chevalley-Eilenberg complex of $\mathfrak{g}$
\begin{displaymath}
	\morphism{c^{\mathfrak{g}}_{p}=(\iota^{n+1}_\mathfrak{g}\omega)\big\vert_p}
	{\Lambda^{n+1} \mathfrak{g}}
	{\mathbb{R}}
	{x_{1} \wedge \dots \wedge x_{n+1}}
	{(\iota( v_{1} \wedge \dots \wedge v_{n+1})\omega)\big\vert_p}
\end{displaymath}
where $p\in M$ is a fixed point in $M$.
Lemma \ref{lemma:multicartan} guarantees that $\delta_\text{CE} c^{\mathfrak{g}}_{p} = 0$.
Note that, when $M$ is connected, the cohomology class $[c_p^{\mathfrak{g}}] \in H^{n+1}(\mathfrak{g})$ is independent of the point $p$, see \cite[Proposition 9.1]{Callies2016}.
\\
The vanishing of $[\tilde{\omega}]$ implies in particular that
$(\iota^{n+1}_\mathfrak{g}\omega) \in C^{n+1}_\mathfrak{g}$ must be a boundary, hence the vanishing of $[c_p^{\mathfrak{g}}]$ is a necessary condition for the existence of a \comoment.
\\
Moreover, it follows from Remark \ref{rk:kuenneth} that when $H^{i}_{\mathrm{dR}}(M) =0$ for $1 \leq i \leq n-1$ the vanishing of $[c_p^{\mathfrak{g}}]$ is also a sufficient condition for the existence of a \comoment.
\end{remark}

\begin{remark}[Relations between $C_\mathfrak{g}$ and equivariant cohomology]
At this stage, this cohomological construction might appear to be a simple formal rearrangement of the moment map definition.
Its real power is manifested when one can draw a correspondence between the cohomology groups of $C$ and those of the other geometric structures involved, namely the manifold $M$, the Lie group $G$ and the action $\vartheta$.
	\\
 When talking about $G$-spaces, \ie triples $(M,G,\vartheta)$, the cohomology theory of interest is called \emph{equivariant cohomology}.
 It can be shown, see \cite[\S 6]{Callies2016}, that the complex $C_\mathfrak{g}$ can be embedded into the Bott-Shulman-Stasheff model $\Omega(G\ltimes M)$ pertaining to the action. 
 Hence the statement of theorem \ref{Thm:CoboundaryinCGRANDE} can be framed as a coboundary condition in equivariant cohomology since $\Omega(G\ltimes M)$ is a model of equivariant cohomology (see proposition 6.4 in \cite{Callies2016}).
 In the same article, it has also been proved how to frame the same cohomological obstruction in another complex called \emph{Cartan model}, which is a model for equivariant cohomology only in the case of compact group actions.
\end{remark}

 In the following subsections we will focus on compact groups trying to give a geometric interpretation of such obstruction which will not depend on the model used to represent the equivariant cohomology of the action.

\subsection{A geometric interpretation of the obstruction class}\label{subsec:geomint}
In symplectic geometry  the existence of a comoment map implies that $\omega$ can be extended to a cocycle in equivariant de Rham cohomology  (\cite{MR721448}). Following \cite{Callies2016}, we illustrate this fact by giving a geometric interpretation to the obstruction class $[\tilde \omega]$ defined above and explain its analogue in the multisymplectic setting.

When the Lie algebra action $v$ comes from a Lie group action, we can interpret the complex $C^\bullet_\mathfrak g$ and the cocycle $\tilde \omega$ in terms of the de Rham cohomology of invariant forms.
\begin{definition}
Let $\vartheta:G\times M\to M$ be a Lie group action. We denote by $\Omega^\bullet(M,\vartheta)$ the subcomplex of $\vartheta$-invariant differential forms. The cohomology of this complex is called \emph{invariant de Rham cohomology} and denoted by $H^\bullet(M,\vartheta)$.
\end{definition}
\begin{remark}
It is more common in the literature to denote these invariant spaces by $\Omega^\bullet(M)^G$ and $H^\bullet(M)^G$. We use the above notation to emphasize the specific Lie group action involved.
\end{remark}
\begin{remark}\label{comactintegration}
The invariant cohomology is not the same as the equivariant cohomology, which we will define later. 
For example, whenever $G$ is compact, the natural map $H^\bullet(M,\vartheta)\to H^\bullet(M)$ induced by the inclusion of the subcomplex is an isomorphism, as in this case any form can be made invariant by averaging. Pullbacks along equivariant maps lead to homomorphisms of the invariant cohomology groups. 
For details we refer to \cite[\textsection 4.3]{MR0336651}.
\end{remark}
\note{Recall that  $C_\mathfrak{g}$ is built from $\wedge^{\geq 1}\mathfrak{g}^*\otimes \Omega^\bullet(M)$ \ie neglecting $\wedge^0\mathfrak g^*$
}
\begin{lemma}[Lemma 6.3 in \cite{Callies2016}]\label{old_infmom}
Let $\vartheta:G\times M\to M$ be a right Lie group action, denote by 
\begin{displaymath}
	\morphism{(r\times \id)}
	{G\times (G\times M)}
	{(G\times M)}
	{(h,(g,m))}
	{(gh,m)}
\end{displaymath}
 the right multiplication action on the second factor.
 \\
 The complex ~$\Omega^\bullet(G\times M,{r\times \id})$ is naturally isomorphic to
 ~$C_\mathfrak g^\bullet\oplus\Omega^\bullet(M)$. 
\end{lemma}

\begin{proof}  
	We have a natural morphism in the category of cochain complexes: 
	\begin{equation}\label{eq:naturalmap}
	 \varphi:~
	 \tot(\Lambda^\bullet\mathfrak g^* \otimes \Omega^\bullet(M))
	 \to \tot(\Omega^\bullet(G,r) \otimes \Omega^\bullet(M))
	 \to \Omega^\bullet(G\times M, {r\times \id}).
	\end{equation}
	The first arrow comes from the isomorphism $\Lambda^k \mathfrak{g}^\ast \to \Omega^k(G,r)$, which associates to any element of $\Lambda^k\mathfrak g^*=\Lambda^kT^*_eG$ its right-invariant differential form extension on $G$.
	This is in particular an isomorphism of cochain complex obtained dualizing the Chevalley-Eilenberg chain complex pertaining to right-invariant vector fields on $G$ that is $CE(\mathfrak{g})\cong CE(\mathfrak{X}(G,r))$. 
	Note that from lemma \ref{lemma:multicartan}, the  cochain differential $\delta_{CE}$ on $(CE(\mathfrak{X}(G,r)))^\ast \cong \Omega(G,r)$ corresponds to the usual de Rham differential.
	\\
	The second arrow comes from the exterior wedge product \ie from the map
	\begin{displaymath}
		\begin{tikzcd}[column sep= small,row sep=0ex]
				\Omega^q(G,r) \otimes \Omega^p(M) 	\arrow[r]
				& \Omega^{q+p}(G\times M, {r\times \id})
				\\
				\alpha\otimes\beta 	\ar[r, mapsto]
				&	\pi_1^*\alpha \wedge \pi_2^* \beta
		\end{tikzcd}		
	\end{displaymath}
	where $\pi_i$ are the projections on the $i$-th factor of $G\times M$.
	Regarding the complexes involved as graded vector spaces, the previous map can be extended to a degree 0, bilinear map
	\begin{equation}\label{eq:kunnetqiso}
		\begin{tikzcd}[column sep= small,row sep=0ex]
				\Omega^\bullet(G,r) \otimes \Omega^\bullet(M) 	\arrow[r]
				& \Omega^{\bullet}(G\times M, {r\times \id})
				\\
				\alpha\otimes\beta 	\ar[r, mapsto]
				&	(-1)^{|\alpha|}\pi_1^*\alpha \wedge \pi_2^* \beta
		\end{tikzcd}		
	\end{equation}
	where the extra signs comes from the Koszul convention we employed in definition \ref{def:comomentbicomplex} when defining the differential in the total complex. 
	\\
	The map $\varphi$ defined in \eqref{eq:naturalmap} admits an inverse given by restricting a form $\alpha\in \Omega^\bullet(G\times M, r\times \id)$ to the submanifold ${\{e\}\times M}$.
	Explicitly one has $\varphi^{-1} := (\iota_e)^\ast $, where $\iota_e : \{e\}\times M \hookrightarrow G \times M$ is the standard inclusion, since
	\begin{displaymath}
		\begin{aligned}
		\Gamma(\iota_e^* \Lambda^k T(G\times M))
		&\cong \Gamma(\Lambda^k((T_e G)_M \otimes TM )) \cong
		\\
		& \cong \Gamma \left(\bigoplus_{n=0}^k \Lambda^n T_e G \otimes \Lambda^{k-n} TM \right) \cong
		\\
		&\cong \bigoplus_{n=0}^k \Lambda^n\mathfrak{g}\otimes \Omega^{k-n}(M)
		~.
		\end{aligned}
	\end{displaymath}
	(In the above equality, $(T_e G)_M$ stands for the trivial vector bundle $(T_e G)\times M \to M$.)
	\\
	The statement follows from the observation that $C_\mathfrak g^\bullet\oplus( \Lambda^0\mathfrak g^*\otimes \Omega^\bullet(M))\cong C_\mathfrak g^\bullet\oplus \Omega^\bullet(M)$ is the total complex of $\Lambda^\bullet\mathfrak g^*\otimes \Omega^\bullet(M)$ and that the second arrow defined above is precisely the function inducing the K\"unneth isomorphism \cite{Bott-Tu82}.
\end{proof}

\begin{proposition}[Proposition 6.4 in \cite{Callies2016}]\label{infmom}
	Assume that $G$ preserves a pre-$n$-\-plec\-tic form $\omega$. 
Call $v$ the infinitesimal action induced by $G$.
	\\
	Then the cocycle $\tilde \omega \in C_\mathfrak g^{n+1}\hookrightarrow\Omega^\bullet(G\times M, {r\times \id})$ with respect to the infinitesimal action $v$ is given by 
	\begin{displaymath}
		\tilde{\omega}= \varphi^{-1}\left(\vartheta^*\omega-\pi^*\omega\right)
		~,
	\end{displaymath}
	where $\pi:G\times M\to M$ is the projection onto the second factor and $\varphi$ is the natural isomorphism defined by equation \eqref{eq:naturalmap}.
\end{proposition}
\begin{proof}
	First must be checked that $\vartheta^* \omega - \pi^* \omega$ is a well-defined $(r\times \id)$-invariant form.
	Being an action, the map $\vartheta\colon G \times M \to M$ is equivariant (with respect to $r\times \id$ in the domain and $\vartheta$ in the target). 
	Hence, the cochain-map $\vartheta^*:\Omega^\bullet(M)\to \Omega^\bullet(G\times M)$ 
	restricts to a well-defined map on the invariant subcomplexes 
	$\vartheta^*:\Omega^\bullet(M,{\vartheta})\to \Omega^\bullet(G\times M,{r \times \id}) $,
	and in particular we have a well-defined map in cohomology.
	\\
	Proceeding by inspection, we consider, without loss of generality,  the point $(e,p)\in G\times M$.
	From the general theory of product manifolds one has that
	\begin{displaymath}
		T_{(e,p)}G\times M \cong T_e G \oplus T_p M \cong \mathfrak{g}\oplus T_p M
		~.
	\end{displaymath}
	Let $X_1,...,X_{n+1-i}\in T_pM$ and $\xi_1,...,\xi_{i}\in \mathfrak g$.
	For all $ 1<i\leq n+1 $ we get
	\begin{align*}
		\vartheta^*\omega(\xi_1,...,\xi_{i},X_1,...,X_{n+1-i}) & = \omega(\vartheta_*\xi_1,...,\vartheta_*\xi_{i},\vartheta_*X_1,...,\vartheta_*X_{n+1-i}) =\\
		&=\omega(v(\xi_1),...,v(\xi_{i}),X_1,...,X_{n+1-i})=\\
		&= \left[(\iota^i_\mathfrak g\omega)(\xi_1,...,\xi_{i})\right](X_1,...,X_{n+1-i})
	\end{align*} 
	and for $i=0$ we get
	\begin{displaymath}
		\vartheta^*\omega(X_1,\dots,X_{n+1} )=
		\omega(X_1,\dots,X_{n+1})
		= \pi^* \omega (X_1,\dots, X_{n+1})
		~.
	\end{displaymath}		
Observe then that by the sign convention in equation \eqref{eq:kunnetqiso} follows that
\begin{displaymath}
	\varphi (\tilde{\omega})
	=
	\sum_{k=1}^{n+1}(-)^{k-1} \varphi (\iota^k_{\mathfrak{g}}\omega )
	=
	-\sum_{k=1}^{n+1}\iota^k_{\mathfrak{g}}\omega  
\end{displaymath}
with the last $\iota^k_{\mathfrak{g}}\omega$ is seen as a differential form on $G\times M$.
Being $\xi_1,\dots,\xi_i,X_1,\dots,X_{n+1-i}$ arbitrary elements, it turns out that
\begin{displaymath}
	\vartheta^* \omega = - \varphi(\tilde{\omega}) + \pi^*(\omega)
\end{displaymath}
and the invertibility of $\varphi$ (see lemma \ref{old_infmom}) finishes the proof.
One should note, that this is true although $\pi$ is not an equivariant map with respect to $( r\times \id, \vartheta)$.
\end{proof}

\begin{remark}
The previous results are subsumed by the following sequence in the category of cochain complexes:
\begin{displaymath}
	\begin{tikzcd}
		\Omega(M,\vartheta) \ar[r,"\vartheta^\ast-\pi^\ast"'] &
		\Omega(G\times M,r\times id) \ar[r,equal,"\varphi^{-1}"',"\sim"] &
		\Lambda^{\geq 0}\mathfrak{g}^\ast\otimes\Omega(M) \ar[r,equal,"\sim"] &[-3em]
		C_{\mathfrak{g}}\oplus\Omega(M)
	\end{tikzcd}
\end{displaymath}
for any smooth action $\vartheta:G\action M$. 
Given $\omega\in \Omega(M,\vartheta)$ multisymplectic, are implied the following cohomological conditions to the existence of a \momap:
\begin{displaymath}
	\left(
		\underset{\text{admits \Hcmm}}{\vartheta:G\action (M,\omega)}
	\right)
	~ \Leftrightarrow ~
	\Big(
		[\tilde{\omega}]=0 \in H(C_{\mathfrak{g}})
	\Big)
	~ \Leftrightarrow ~
	\Big(
		[(\vartheta^\ast-\pi^\ast)\omega]=0 \in H(G\times M, r\times id)	
	\Big)
\end{displaymath}
Notice, in particular, that the right-hand side provides an obstruction in terms of the $G$-invariant cohomology, which is closer to the problem's geometric data and less "ad hoc" than $C_{\mathfrak{g}}$.
\\
All of this does not require the group to be compact. 
In case of compact groups, the obstruction can be read directly in de Rham cohomology, see corollary  \ref{cor:core}.
\end{remark}

\begin{remark}
We think that the above proposition is central to the understanding of multisymplectic \comoments, as it enables an elementary and unified treatment of many phenomena in multisymplectic geometry.
	\begin{itemize}
		\item For symplectic manifolds, this result gives a nice interpretation for the result of \cite{MR721448} that moment maps are in correspondence to equivariant extensions of $\omega$ and also explains why this correspondence fails in the general multisymplectic setting (\cf Example \ref{exnongen} and \cite[\textsection 4]{weinstein1977lectures}). 

		\item A \comoment exists, whenever the multisymplectic form $\omega$ can be lifted to a class in the equivariant cohomology $H_G^{n+1}(M)$.
		(See Theorem \ref{cpteqcom}). 
		
		\item Let $G_i$ act on the multisymplectic manifolds $(M_i,\omega_i)$ for $i\in\{1,2\}$. 
		If there exists a \comoment for $\vartheta_i:G_i\action (M_i,\omega_i)$, then there exists a \comoment for the product $\vartheta=(\vartheta_1\times \vartheta_2): (G_1\times G_2) \action (M_1\times M_2)$ with respect to the multisymplectic structure $\omega=(\tau_1^*\omega_1\wedge \tau_2^*\omega_2)\in \Omega(M_1\times M_2)$. 
		(Here $\tau_i:M_1\times M_2 \to M_i$ denotes the standard projection of the cartesian product). 
		This theorem from \cite{Shahbazi2016}, proven by explicit exhibition of the sought \comoment, can be derived from proposition \ref{infmom} observing that
		\begin{align*}
					 [ (\vartheta^\ast - & \pi^\ast) \omega]
					 =
					 \\ 
			=&~		
			 [ (\vartheta_1^\ast\times\vartheta_2^\ast - \pi_1^\ast \times \pi_2^\ast) (\tau_1^*\omega_1\wedge \tau_2^*\omega_2)] 
			 =
			\\
			=&~
			[\tau_1^*\vartheta_1^\ast \omega_1 \wedge \tau_2^*\vartheta_2^\ast \omega_2
			-
			\tau_1^*\pi_1^\ast \omega_1 \wedge \tau_2^*\pi_2^\ast \omega_2
			]
			~=
			\\
			=&~
			[ (\tau_1^*(\vartheta_1^*-\pi_1^*)\omega_1)\wedge (\tau_2^*\vartheta_2^* \omega_2) +
			  (\tau_1^\ast\pi_1^\ast\omega_1)\wedge (\tau_2^\ast(\vartheta_2^\ast-\pi_2^\ast)\omega_2)]		
			~=
			\\
			=&~
			\kappa\left(
				\cancel{[(\vartheta_1^\ast -\pi_1^\ast) \omega_1]},[\vartheta_2^\ast \omega_2] 
			\right) +
			\kappa\left(
				[\pi_1^\ast \omega_1], \cancel{[(\vartheta_2^\ast -\pi_2^\ast )\omega_2]}
			\right)
			~=
			\\
			=&~ 0
			~,
		\end{align*}
		where $\kappa$ denotes the K\"unneth isomorphism and the last cancellations are due to the existence of a \momap with respect to the components $\vartheta_i$.
		A similar result also holds endowing $(M_1\times M_2)$ with the multisymplectic form $\pi_1^\ast \omega + \pi_2^\ast \omega$.
		\note{Maggiori dettagli nel file LeoQ.tex}
	\end{itemize}
\end{remark}

\begin{remark}\label{rk:invpot}
In particular, Proposition \ref{infmom} immediately implies that a $\vartheta$-invariant potential of $\omega$ induces a \comoment, 
as an invariant potential $\alpha$ of $\omega$ would be mapped to a potential $(\vartheta^*\alpha-\pi^*\alpha) \in \Omega^\bullet( G\times M , r\times \id)$ of $\tilde{\omega}=\vartheta^*\omega-\pi^*\omega$.
Note that $\omega$ being exact is not a sufficient condition, 
because a primitive $\vartheta^*\alpha-\pi^*\alpha$ need not to be an element in the invariant cochain complex $\Omega^\bullet(G\times M, r\times \id)$ in general. \\ 
\end{remark}

When such an invariant potential exists, it is fairly easy to give an explicit expression for the components of a \comoment, as illustrated by the following Lemma:
\begin{lemma}[Section 8.1 in \cite{Callies2016}]\label{lem:extexact}
	Let $M$ be a manifold with a $G$-action, 
	let $\alpha\in \Omega^n(M,G)$ be a $G$-invariant $n$-form and consider the pre-$n$-plectic form $\omega=d\alpha$ on $M$.\\
	The action $G\action \left(M,d\alpha\right)$ admits a $G$-equivariant \momap $(f):\mathfrak{g} \to L_{\infty}(M,\omega)$, given by $(k=1,\dots,n)$:
\begin{displaymath}
	\morphism{f_{k}}
	{\Lambda^k\mathfrak{g}}
	{\Omega^{n-k}(M)}
	{q}
	{(-1)^{k-1}\varsigma(k)\iota(v_q)(\alpha)~.}
\end{displaymath}

\end{lemma}
\begin{proof}
	A direct proof of this statement can be given by showing that equation \eqref{eq:fk_hcmm} is satisfied.
	Upon employing Lemma \ref{lemma:multicartan}, we have:
	\begin{align*}
		\textrm{d} f_m (p) &= (-1)^{m-1} \varsigma(m) \textrm{d} \iota_{v_p} \alpha =
		\\
		&=
		-\varsigma(m) (-1)^{m} \textrm{d} \iota(v_1\wedge\dots\wedge v_m) \alpha =\\
		&= -\varsigma(m) \left(\iota_{v_p} \textrm{d}\alpha + \iota_{\partial v_p} \alpha +
		\sum_{k=1}^{m} (-1)^k  \iota( x_1\wedge\, \hat{x_k}\, \wedge\dots\wedge x_m) \cancel{\mathcal{L}_{x_k} \alpha}
		 \right) =\\
		&= -\varsigma(m) \iota_{v_p} \omega +  (-1)^{m-1} \varsigma(m-1) \iota_{\partial v_p} \alpha  =
		-\varsigma(m) \iota_{v_p} \omega - f_{m-1}(\partial v_p) ,
	\end{align*}	
	thus $G$-equivariance follows from lemma \ref{lemma:multicartan}.
\end{proof}
\begin{corollary}[$SO(n)$-action on $\mathbb{R}^n$, Example 8.4 in \cite{Callies2016}]\label{cor:sorn}
	The canonical action 
	$$SO(n) \action \left( \mathbb{R}^{n}, dx^{1\dots n}\right),$$
	where~$x = (x^i)$ are the standard coordinates on~$\mathbb{R}^n$ 
and~$dx^{1\dots n} = d x^1\wedge\dots \wedge dx^{n}$ is the standard volume form of $\mathbb{R}^n$, admits a \comoment given by $(k=1,\dots,n)$:
	\begin{displaymath}
		\morphism{f_{k}}
		{\Lambda^k\mathfrak{g}}
		{\Omega^{n-1-k}(M)}
		{q}
		{(-1)^{k-1} \frac{\varsigma(k)}{n} \iota(E \wedge v_q)\left(dx^{1\dots n}\right)}
	\end{displaymath}
 	where $E =\sum_i x^i\partial_i$ is the Euler vector field.
\end{corollary}
\begin{proof}
    The proof follows from Lemma \ref{lem:extexact} noting that the standard volume form admits the following $G=SO(n)$-invariant primitive
    $$\alpha = \dfrac{\iota_E\left(dx^{1\dots n}\right)}{n}~.$$
\end{proof}

\subsection{Cohomological obstructions for compact group actions}
Let us specialize our description of the obstructions to the existence of a \momap pertaining to a compact Lie group action. 
In this case,  we do not have to care about the invariance of forms and Proposition \ref{infmom} reads as follows:
\begin{corollary}\label{cor:core}
Let $\vartheta:G\times M\to M$ be a compact Lie group acting on a pre-$n$-plectic manifold, preserving the pre-multisymplectic form $\omega$. 
A \comoment exists if and only if $[\vartheta^*\omega-\pi^*\omega]=0\in H^{n+1}_{dR}(G\times M)$. 
\end{corollary}
\begin{proof}
From Proposition \ref{infmom} we get the following sequence of maps between complexes together with the induced maps in cohomology:
\begin{center}
\begin{tikzcd}
 \Omega^\bullet(M,\vartheta) \ar[d,"\vartheta^\ast-\pi^\ast"] &\quad
 &[-2em] H_\text{dR}(M) \ar[d,"\vartheta^\ast-\pi^\ast"]  
 &[-2em] \lbrack \omega \rbrack \ar[d,mapsto]
 \\ 
 \Omega^\bullet(G\times M, r\times \id) \ar["\cong",leftrightarrow]{d} &\quad
 & H_\text{dR}(G\times M) \ar[leftrightarrow,"\cong"]{d}[swap]{\text{\tiny (K\"unneth)}} 
 & \lbrack \vartheta^\ast \omega - \pi^\ast \omega \rbrack \ar[ddd,mapsto]
 \\ 
 \Omega^\bullet(G,r) \otimes \Omega^\bullet(M) \ar["\cong",leftrightarrow]{d}[swap]{} &\quad
 & H_\text{dR}(G) \otimes  H_\text{dR}(M) \ar[d,"\cong",leftrightarrow]
 \\ 
 \Lambda^\bullet \mathfrak{g}^* \otimes \Omega^\bullet(M)\ar["\cong",leftrightarrow]{d} &\quad
 & H_\text{CE}(\mathfrak{g}) \otimes  H_\text{dR}(M) 
 \ar["\cong",leftrightarrow]{d} & 
 \\
 C_\mathfrak{g}^\bullet \oplus ( \mathbb{R}\otimes \Omega^\bullet(M))&\quad 
 & H(C_\mathfrak{g}^\bullet)\oplus H_\text{dR}(M)
 & \lbrack \tilde{\omega}\rbrack
\end{tikzcd}
\end{center}
The statement follows by resorting to Remark \ref{comactintegration}, \ie noting that 
$H_\text{dR}(G) \cong H(G,r)$ and $H_\text{dR}(G\times M) \cong H(G\times M, r\times \id)$ via the averaging trick on compact connected Lie groups (see reminder \ref{rem:Averaging}).
\end{proof}

\begin{reminder}[Averaging trick]\label{rem:Averaging}
Let $\vartheta: G\action M$ be a compact and connected Lie group action.
By \emph{averaging trick} we mean the procedure associating a $\vartheta$-invariant differential form to any differential form on $M$.
\\
Recall at first that Lie groups are orientable (since the tangent bundle $TG\cong G\times \mathfrak{g}$ is trivial) and in particular it is always possible to find left or right invariant volume forms thereon.
Such volume forms are uniquely defined by a constant.
When $G$ is a compact group, such constant can be fixed via normalization. Namely, there exists a unique $dg \in (\Omega^{\text{top}}(G),R)$ right-invariant volume form, called \emph{Haar volume}, such that $\int_G \d g = 1$ (see for instance \cite[\S 3.13]{Duistermaat2000}).
\\
Consider now the aforementioned compact group action $\vartheta:G\action M$.
For any $0\leq k \leq \text{dim}(M)$, we define the \emph{averaging} as the mapping
\begin{displaymath}
	\morphism{A}
	{\Omega^k(M)}
	{(\Omega^k(M),\vartheta)}
	{\omega}
	{\displaystyle \int_G (\vartheta^\ast_g \omega) dg}~,
\end{displaymath}
	where $dg$ denotes the Haar measure, $\vartheta_g:x\mapsto x.g$ is the action of $g$ and $(\vartheta_g^\ast \omega)$ is the pullback of $\omega$ along the action seen as a vector-valued smooth function on $G$:
	\begin{displaymath}
		\vartheta^\ast_\blank \omega : G \to \Omega^k(M)~.
	\end{displaymath}
A simple consequence of this construction is that any $\vartheta$-invariant exact form admits a $\vartheta$-invariant primitive.
Namely, considering $\omega= \d \beta$, one has
$$
	\omega = A(\omega) = \int_G(\vartheta_g^\ast \d \beta) dg = \d A(\beta)~.
$$

\end{reminder}

\medskip
\noindent
We will now investigate the connection between \comoments and equivariant cohomology.
The latter is the most general cohomology theory pertaining to group actions on topological space (see \cite{Tu2011} for a gentle introduction).
\begin{definition}[{Equivariant Cohomology (\cite{Borel1960}, see also \cite[\S 3]{Hsiang1975}\cite[\S 1]{Guillemin1999})}] \label{def:equivariantcoho}
Consider a smooth Lie group action $\vartheta:G\action M$ on a manifold $M$. 
Let $EG$ be a contractible space on which $G$ acts freely by $\vartheta^{EG}$. 
Then we define the \emph{equivariant cohomology of $M$} as $H^\bullet_G(M):=H^\bullet((M\times EG)/G)$, where $G$ acts on $M\times EG$ diagonally.
\end{definition}

\begin{remark}
	As $EG$ might not be a manifold, we have to interpret $H^\bullet_G(\cdot)$ as a suitable cohomology theory (\eg singular cohomology with real coefficients) in the above definition.
	\\
	When the group $G$ is compact, any action $\vartheta:G\times M\to M$ is in particular proper. 
If the action $\vartheta$ is also free then the quotient manifold theorem holds and $H_G^\bullet(M)=H^\bullet_{dR}(M/G)$. 
\end{remark}

For a not necessarily free action $\vartheta$, we have the following diagram in the category of topological spaces
\begin{center}
\begin{tikzcd}[column sep=large]
	G\times (M\times EG)  \ar[r,shift left=1.5ex,"\vartheta\times \vartheta^{EG}"]\ar[r,shift right=1.5ex,,"\pi"]
	& M\times EG \ar[r,"q"]& (M\times EG)/G,
\end{tikzcd}
\end{center}
where $q$ is the projection to the orbits, which induces a sequence in cohomology:
\begin{center}
\begin{tikzcd}[column sep=large]
	H^\bullet(G\times M) & H^\bullet(M) \ar[l,"\vartheta^*-\pi^*"]
	& \ar[l,"q^*"]H^\bullet_G(M)
\end{tikzcd}
\end{center}
The latter diagram is obtained using the definition of equivariant cohomology $H(\frac{M\times EG}{G})\cong H_G(M)$ on the rightmost term and the contractibility of $EG$, \ie
	\begin{displaymath}
		\begin{split}
		H(M\times EG) &\underset{\text{de Rham thm.}}{\cong} 
		H_{dR}(M\times EG) \cong \\
		&\underset{\text{K\"unneth thm.}}{\cong} 		
		H_{dR}(M)\otimes H_{dr}(EG) \cong \\ 
		&\underset{\text{$EG$ conctractible}}{\cong} 	
		H_{dR}(M)~.
		\end{split}
	\end{displaymath}	 
Being $q$ a mapping into orbits one has $q\circ \vartheta=q\circ\pi$ (where $\theta$ here is shorthand for the action $\vartheta\times \vartheta^{EG}$ on $M\times EG$.
Therefore $(\vartheta^*-\pi^*)\circ q^*=0$.
Using Remark \ref{comactintegration} and Proposition \ref{infmom} we get the following statement:
\begin{theorem}[\cite{Callies2016}]\label{cpteqcom}
Let $G\times M\to M$ be a compact Lie group preserving a pre-$n$-plectic form $\omega$. If $[\omega]\in H^\bullet(M)$ lies in the image of $q^*:H^\bullet_G(M)\to H^\bullet(M)$, then $\vartheta$ admits a \comoment.
\end{theorem}
\begin{remark}
	The advantage of our approach to the theorem is that it is much simpler and more intrinsic, in particular we do not need to choose a model for equivariant cohomology.
In Section 6 of \cite{Callies2016} one can find similar results framed in the Bott-Shulman-Stasheff and in the Cartan model.
 Furthermore, in Section 7.5 of \cite{Callies2016}, a generalization of this statement to non-compact groups is discussed.
 \end{remark}

Observe that the vector space $\text{Im}(q^\ast) \subset H_{dR}(M)$ consists of closed forms (it is better to say cohomology classes) coming from a class in equivariant cohomology; in other words, elements in $\text{Im}(q^\ast)$ can be extended to an equivariant cohomology class.
The condition for $[\omega]$ to come from an equivariant cohomology class is sufficient for the existence of a \momap.
\\
Unfortunately, unlike the symplectic case (\cf \cite{MR721448}), the converse statement does not hold in general. 
Even if a (pre-)multisymplectic action of $G$ on $(M,\omega)$ admits a \comoment, $[\omega]$ does not need to come from an equivariant cocycle. 
We will illustrate this fact by an example.

\begin{example} \label{exnongen}
Consider the action $\vartheta: S^1\times S^3\to S^3$ given by the Hopf fibration $S^1\hookrightarrow S^3 \twoheadrightarrow S^2$.
Let $\omega$ be the standard volume on $S^3$. By Remark $\ref{rk:kuenneth}$ the obstructions to a \comoment sit in $H^1(\mathfrak u(1))\otimes H^2(S^3) $, $ H^2(\mathfrak u(1))\otimes H^1(S^3)$ and  $H^3(\mathfrak u(1))\otimes H^0(S^3) $, which are trivial. Hence, a \comoment exists. 
\\
As the action is free (and the quotient is $S^2$), we have $H^3_{S^1}(S^3)= H^3(S^3/S^1)=H^3(S^2)=0$. 
But $[\omega]\neq 0$ in $H^3(S^3)$, so the class $[\omega]$ cannot come from an equivariant cocycle.
\end{example}

\begin{remark}
We note that this example has a different character than the ones provided in Section 7.5 of \cite{Callies2016}. 
They exhibit cases where individual \comoments do not come from equivariant cocycles, whereas in our case no equivariant cocycle can be found for any of the possible \comoments.
\end{remark}

\section{Non-transitive multisymplectic group actions on spheres}\label{secintrans}
The goal of this section is to prove the existence of \comoments for non-transitive actions and construct an explicit \comoment for the $SO(n)$-action on $S^n$.
\begin{remark}[Canonical multisymplecticity of the sphere]
	Recall that the $n$-dimensional sphere $S^n$ is a connected, simply connected manifold; therefore is orientable.
	Recall also that $S^n$ can be canonically embedded in the Euclidean space $\R^{n+1}$ as the standard unit sphere centered on the origin.
	Thus, it possesses a canonical Riemannian structure induced by the standard metric on $\R^{n+1}$.
	One can induce a volume form on $S^n$ by pulling-back contraction of the Euclidean volume with a nowhere vanishing normal vector field.
	At this scope, it is customary to choose the outward-facing vector field $E = x^i \partial_i$ called \emph{Euler vector field}.
	The unique Riemannian volume form singled out of this procedure will be the standard multisymplectic form considered on spheres. (See example \ref{Ex:VolumesAreMultiSymp}.)
\end{remark}

Recall first that the \emph{orbit map of $p\in M$} with respect to the smooth action $\vartheta:G\action M$ is the smooth function
\begin{displaymath}
	\morphism{\vartheta_p}
	{G}
	{M}
	{g}
	{g.p = \vartheta(g,p)}
\end{displaymath}
determining an immersed submanifold of $M$ simply called \emph{the orbit of $p$}.
\\
The existence of a \momap for actions on a sphere can be ascertained by studying the restriction of the multisymplectic form on orbits:
\begin{lemma}\label{lem:core}
	Let $\vartheta: G\times S^n\to S^n$ be a compact Lie group acting multisymplectically on $S^n$ equipped with the standard volume form $\omega\in \Omega^n(S^n)$. 
	Let $p\in S^n$ be any point.
	\\
	A \momap exists if and only if $\vartheta_p^*[\omega]\in H^{n}(G)$ vanishes.
\end{lemma}
\begin{proof}
	By Corollary \ref{cor:core}, we only have to check that
	$$ 
	[\vartheta^*\omega-\pi^*\omega]=0\in H^{n}(G\times S^n)
	\quad\Leftrightarrow\quad 
	\vartheta_p^*[\omega]= 0\in H^{n}(G)~.$$
	The direct implication follows by considering the map 
	\begin{displaymath}
	 \morphism{i}
	 {G}
	 {G\times S^n}
	 {g}
	 {(g,p)}
	\end{displaymath}
	and its induced linear map in cohomology $i^* : H^\bullet(G \times S^n) \to H^\bullet(G)$ which acts on $[\vartheta^*\omega - \pi^*\omega] \in H^n(G\times S^n)$ as
\begin{displaymath}
	i^* [\vartheta^* \omega - \pi^* \omega ] = 
	[(\vartheta \circ i)^\ast \omega - \cancel{(\pi \circ i)^\ast \omega}] = \vartheta^*_p[\omega] ,
\end{displaymath}
	because $\vartheta \circ i = \vartheta_p$ and $\pi \circ i $ is the constant map valued in $p\in S^n$.
	\\
	For the converse implication, note at first that the cohomology of the sphere implies
	$$H^{n}(G\times S^n)=(H^n(G)\otimes H^0(S^n))\oplus (H^0(G)\otimes H^n(S^n)) .	$$ 
	 Therefore, recalling Proposition \ref{infmom}, the obstruction $[\vartheta^*\omega-\pi^*\omega]$ lies entirely in $H^n(G)\otimes H^0(S^n)$ as $[\tilde{\omega}]$ has null component in $ H^0(G)\otimes H^n(S^n)$.
	 Since the restriction of $i^*$ to $H^n(G)\otimes H^0(S^n)$ is an isomorphism (the $0$-th cohomology group of a connected manifold is isomorphic to $\mathbb{R}$), one can conclude that $[\vartheta^*\omega-\pi^*\omega]$ vanishes if and only if 
	$$i^*[\vartheta^*\omega-\pi^*\omega]=\vartheta_p^*[\omega]=0\in H^n(G) .$$

\end{proof}
\begin{remark}
	Note that the direct implication in the proof of Lemma \ref{lem:core} does not depend from being the base manifold a sphere.	
	In other words, a necessary condition for the existence of a \comoment is that the pullback of $\omega$ with respect to any orbit map vanishes.
\end{remark}
\begin{remark}
Observe that a result similar to Lemma \ref{lem:core} can be stated for any compact multisymplectic action $G \action (	M,\omega)$, with $\omega$ in degree $n+1$, such that the following cohomological condition holds
$$
\bigoplus_{k=1}^n H^k(G)\otimes H^{n-k}(M) = 0 .
$$
In particular, the same result is true for the action of any compact, connected and semisimple Lie group $G$ (\ie such that $H^1(\mathfrak{g})=H^2(\mathfrak{g})=0$) acting on a $2$-plectic manifold.
See \cite[Proposition 7.1]{Callies2016} for an existence result for \comoments related to this kind of actions in presence of a fixed point.
\end{remark}
	An action $G\action M$ is said \emph{transitive} if the orbit map $\vartheta_p$ is surjective for any $p\in M$.
	When this condition is not met we have the following result:
\begin{proposition}\label{prop:intransitive}
	Let $G$ be a compact Lie group acting non-transitively on $S^n$ and preserving the standard volume form. 
	Then the action $G\action M$ admits a \comoment.
\end{proposition}
\begin{proof} 
	If $G$ acts non-transitively then there exists an orbit $O\subset S^n$ of dimension strictly less than $n$. Let $p\in O$. 
	Then we have $\vartheta_p^*[\omega]=\vartheta_p^*[\omega|_O]$, but $\omega|_O\in \Omega^n(O)$ is zero due to dimension reasons. Hence, the action admits a \comoment, due to Lemma \ref{lem:core}.
\end{proof}

\subsection{The action of $SO(n)$ on $S^n$}\label{subsecson}
The goal of this section is to give an explicit construction for the \comoment of the action $SO(n) \action S^n$ by resorting to Corollary \ref{cor:inducedmachinery}.

\vspace{1ex}
In what follows, we will consider the standard $SO(n+1)$-invariant embedding $ j: S^n \to \mathbb{R}^{n+1}$ of $S^n$ as the sphere with unit radius
and consider the linear action of the group $SO(n)$ on $\mathbb{R}^{n+1}$ as the subgroup of special orthogonal linear transformations fixing the axis $x^0$.
\\
In~$\mathbb{R}^{n+1}$ we consider the standard coordinates~$x = (x^0,\dots,x^n)$ and the corresponding volume form~$\Vol = d x^0\wedge\dots \wedge dx^{n}$. 
Furthermore, we will make use of the following notation for the cylindrical radius (with axis $x^0$) and the central radius respectively:
\begin{displaymath}
r = \sqrt{(x^1)^2 + \dots + (x^n)^2} \qquad R = \sqrt{(x^0)^2 + r^2} 
~.
\end{displaymath}

\vspace{1ex}
Recall that the volume form on the unit sphere embedded in $\mathbb{R}^{n+1}$ is given by 
$$ \omega = j^\ast \iota_E ~\Vol$$
where $E$ is the Euler vector field.
$E$ can be seen as the fundamental vector field of the action
\begin{align*}
	\vartheta: \mathbb{R} \times \mathbb{R}^{n+1}
	&\to \mathbb{R}^{n+1}\\
	(\lambda,x)&\mapsto e^{\lambda} x
\end{align*} 
of $\mathbb{R}$ on the Euclidean space via dilations,
that is the linear action generated by the identity matrix $\id_{n+1} \in \mathfrak{gl}(n,\mathbb{R}^{n+1})$, \ie $ E= v_{\id_{n+1}} = \sum_i x^i\, \partial_i$.

Let us call $H = SO(n)\times\mathbb{R}$ the subgroup of $GL(n,\mathbb{R}^{n+1})$ generated by
\begin{displaymath}
	\mathfrak{h} = 
	\left\lbrace
		\begin{bmatrix} 
			1 & 0 \\ 
			0 & a 
		\end{bmatrix}
		\; \vert \: a \in \mathfrak{so}(n)
	\right\rbrace 
	\oplus \langle \id_{n+1} \rangle \simeq \mathfrak{so}(n)\oplus \mathbb{R}
	~.
\end{displaymath}
The group $H$ acts linearly on $\mathbb{R}^{n+1}$ through the standard infinitesimal action
\begin{displaymath}
	\morphism{v}
	{\mathfrak{h}}
	{\mathfrak{X}(\mathbb{R}^{n+1})}
	{A}
	{\displaystyle\sum_{i,j} [A]_{i j} ~ x^j \partial_i}
\end{displaymath}
where $[A]_{i j}$ denotes the $i j$ entry of the matrix $A$.

\begin{lemma}\label{lem:rescaledvolume}
The differential form
\begin{displaymath}
 \eta = \rescaling \,\Vol \in \Omega^{n+1}(\mathbb{R}^{n+1} \setminus \{0\})
 \qquad \text{with} \qquad \rescaling = \dfrac{1}{R^{n+1}}
\end{displaymath}	
is multisymplectic on $N = (\mathbb{R}^{n+1} \setminus \{0\})$, invariant under the action $H \action N$ and restricts to the standard volume form on the unit sphere.
\end{lemma}
\begin{proof}
	Multisymplecticity follows from the closure and non-degeneracy of $\Vol$ together with the property that $\rescaling$ never vanishes.
	\\ 
	The form is clearly $SO(n)$-invariant because $\rescaling$ depends only on $R$. The $H$-invariance follows from the invariance along the Euler vector field:
$$
	\mathcal{L}_E (\rescaling~\Vol) = (\mathcal{L}_E \rescaling + (n+1) \rescaling ) \Vol
$$
and
$$
	\mathcal{L}_E \rescaling = \sum_i x^i \dfrac{\partial}{\partial x^i} R^{-(n+1)}	
	= -\dfrac{n+1}{2}\sum_i 2 (x^i)^2 R^{-(n+1)-1} =-(n+1) \rescaling	
$$
Finally, $\eta$ restricts to the volume form on $S^n$ because $j^\ast \rescaling = 1$.
\end{proof}

The function $\rescaling \in C^\infty(\mathbb{R}^{n+1}\setminus\{0\})$ is precisely the scaling factor that makes the Euclidean volume invariant with respect to the extended group $SO(n+1)\times\mathbb{R}$.
The problem to find explicitly a \comoment:
\begin{displaymath}
s \colon \mathfrak{h} \rightarrow L_{\infty} \left(\mathbb{R}^{n+1}\setminus\{0\},\eta = \rescaling \, \Vol\right)
\end{displaymath}
can be solved by exhibiting an $H$-invariant primitive of  the rescaled volume $\eta$ and then resorting to Lemma \ref{lem:extexact}.

\begin{lemma}\label{lem:uglyprimitive}
	The differential (n+1)-form $\eta$ admits an $H$-invariant potential $n$-form: 
	\begin{displaymath}
		\beta = (\hat{\varphi}~x^0)~dx^1\wedge\dots\wedge dx^n
	\end{displaymath}	
where $\hat{\varphi}\in C^\infty (\mathbb{R}^{n+1}\setminus{0})$ depends only on the cylindrical coordinates $(x_0,r)$ and it is given by
	\begin{equation}\label{eq:uglyprimitive}
	\hat{\varphi}(x^0,r)  = 	
	\begin{cases}
		\frac{1}{\left((x^0)^2 + r^2\right)^{\frac{n+1}{2}}}
		\left(x^0 (n+1) - r \arctan\left(\dfrac{x^0}{r}\right)\right)
		&\quad r \neq 0
		\\
		(n+1)\frac{1}{\vert x^0 \vert^n} 
		&\quad r=0,\, x^0 > 0
		\\  
		-(n+1)\frac{1}{\vert x^0 \vert^n} 
		&\quad r=0,\, x^0 < 0
	\end{cases}
	\end{equation}
\end{lemma}
\begin{proof}
	Let us start from the following ansatz 
		$$\beta =\iota_{(x^0\partial_0)} \varphi(x^0,r) \eta $$
	for a potential $n$-form of the scaled volume $\eta$ as defined in Lemma \ref{lem:rescaledvolume}.
	At this point, $\varphi$ is an arbitrary smooth function depending only on the cylindrical parameters $(x^0,r)$.
	\\
	Being $x^0 \partial_0$ the fundamental vector field of 
	\begin{displaymath}
		\zeta =
			\begin{bmatrix} 
				1 & 0 \\ 
				0 & 0_n 
			\end{bmatrix}
		\in \mathfrak{gl}(n+1,\mathbb{R})
		,
\end{displaymath}
one gets
\begin{displaymath}
\mathcal{L}_{v_\xi} \beta = \left(
\iota(\cancel{v_{[\xi,\zeta]}}) + \iota({v_\zeta}) \mathcal{L}_{v_\xi}
\right)
\varphi \, \eta = 0
\qquad \forall \xi \in \mathfrak{so}(n),
\end{displaymath}
because the $(n+1)$-form $\varphi \, \eta $ depends only on $(x^0,r)$,
\ie $\beta$ is $SO(n)$ invariant. On the other hand, one has:
\begin{displaymath}
	\begin{aligned}
		\mathcal{L}_E \beta &= \left(
		\cancel{\iota_{[\id,\zeta]}} + \iota_{x_0 \partial_0} \mathcal{L}_{E}
		\right)
		\varphi \, \rescaling \, \Vol = 
		\\
		&= \iota_{x_0 \partial_0} \left[
		\left(\mathcal{L}_E \varphi \right)\, \rescaling \, \Vol 
		+ \varphi \cancel{\mathcal{L}_E\,\rescaling\Vol}
		\right]
	\end{aligned}
\end{displaymath}
and
\begin{displaymath}
	\begin{aligned}
		\textrm{d} \beta &=
		\left(
		\dfrac{\partial \varphi}{\partial x^0}\,\rescaling\,x^0 +
		\varphi\, x^0 \, \dfrac{\partial \rescaling}{\partial x^0} + 
		\varphi \rescaling \right) \Vol 
		\\
		&= \left( \partial_0 \varphi \, x^0 - (n+1)\dfrac{(x^0)^2}{(r^2 + (x^0)^2)} 
		+ \varphi \right)\, \rescaling\,\Vol 
		~.
	\end{aligned}
\end{displaymath}
Therefore, in order for $\beta$ to be $G-$invariant primitive of $\omega$, the following conditions on $\varphi$ have to be met:
\begin{displaymath}
	\begin{cases}
		\mathcal{L}_E \varphi = r\,\partial_r\,\varphi + x^0\,\partial_0\,\varphi = 0 
		\\
		x^0\,\partial_0 \varphi   - (n+1)\dfrac{(x^0)^2}{(r^2 + (x^0)^2)} + \varphi = 1
	\end{cases}
\end{displaymath} 
The general solution of this system reads:
\begin{displaymath}
	\varphi(x^0,r) = - \dfrac{r}{x^0}\arctan\left(\dfrac{x^0}{r}\right)(n+1) + n + 2
\end{displaymath}
which is a smooth function defined on 
$\left\{ x \in \mathbb{R}^{n+1} ~\vert x^0 \neq 0, r \neq 0\right\}$.
Recalling that 
$$\lim_{y\to 0}\dfrac{\arctan(y)}{y}=1 \qquad \lim_{y\to \infty}\dfrac{\arctan(y)}{y}=0 ~,$$
one can see that the limits to all the critical points of $\varphi$, except $(x^0,r)=(0,0)$, are finite.
Hence, considering the unique smooth extension $\hat{\varphi}\in C^\infty(\mathbb{R}^{n+1}\setminus\{0\})$ of $\varphi$, given explicitly by equation (\ref{eq:uglyprimitive}), we conclude that 
	\begin{displaymath}
		\beta = (x^0~\hat{\varphi})~dx^1\wedge\dots\wedge dx^n
	\end{displaymath}	
is the sought $H$-invariant primitive.
\end{proof}

\begin{proposition}\label{Prop:SonSn}
	A \comoment for the action $SO(n) \action \left( S^{n}, \omega\right)$, for $n \geq 2$, is given by
	\begin{displaymath}
		\morphism{f_i}
		{\Lambda^i\mathfrak{so}(n)}
		{\Omega^{n-1-i}(S^n)}
		{q}
		{-j^\ast\iota(v_q)(\iota_E \beta)}
		~.
	\end{displaymath}
	where $\beta$ is the primitive given in Lemma \ref{lem:uglyprimitive}.
\end{proposition}
\begin{proof}
	The statement follows directly from Corollary \ref{cor:inducedmachinery} upon considering

	\begin{align*}
		(N,\eta) = (\mathbb{R}^{n+1},\rescaling\Vol ) \qquad\quad 
		&(M,\omega) = (S^n, j^\ast \iota_E \Vol) 
		\\
		p= \id_{(n+1)} \in Z_1(\mathfrak{h}) \qquad\quad
		& H = SO(n)\times \mathbb{R} = H_E \supset SO(n)		
	\end{align*}
	and noting that an explicit \comoment for the $H$-action is given by Lemma \ref{lem:extexact} via employment of the primitive constructed in Lemma \ref{lem:uglyprimitive}.
\end{proof}

\note{Check if this result gives the known result for $n=3$. (Angular momentum)}

\begin{remark}
	Proposition \ref{Prop:SonSn} extends to spheres of arbitrary dimension a similar result given in \cite[Paragraph 8.3.2]{Callies2016} up to dimension $5$.
\end{remark}

\section{Transitive multisymplectic group actions on spheres}\label{sectra}
Recall at first an important property of transitive compact group actions:
\begin{proposition}[Isotropies of transitive compact group actions \emph{(see e.g \cite{encyclopedia:Isotropy})}]
	Be $G\action M$ a transitive compact group smooth action and let $H\subset G$ be an isotropy subgroup.
	Then $H$ is a closed subgroup, all isotropy subgroups are isomorphic to $H$, \ie
	\begin{displaymath}
		G_x:= \left\lbrace g \in G ~\vert~ g.x = x \right\rbrace\cong H
		\qquad \forall x \in M
	\end{displaymath}		
	and there exists a canonical diffeomorphism $G/H \cong M$.
\end{proposition}
\begin{remark}\label{Rem:ActionasBundle}
	Synthetically, we will denote any transitive action of a compact group $G$ on the manifold $M$ with isotropy subgroup $H$ simply by the corresponding canonical isomorphism $G/H = M$.
	
	Recall also that there exists a corresponding $H$-principal bundle over $M$ given by fixing a certain point $p\in M$ and considering the corresponding orbit map
	\begin{displaymath}
		\begin{tikzcd}[column sep= small,row sep=0ex]
			H \cong G_p \ar[r,hook,"i"]
			& G \ar[r,two heads, "\vartheta_p"]
			& M
		\end{tikzcd}
		~.
	\end{displaymath}
\end{remark}

Focussing to our case, all effective transitive compact connected group actions on spheres have been classified (\cf \cite{MR2371700} for an overview of the results):
\begin{proposition}[Classification of transitive compact group action on spheres, \cite{MR0008817, MR0034768,MR0029915} ]\label{Prop:ClassificActionSpheres}
	The only compact groups $G$ acting transitively and effectively on $S^n$ with isotropy group $H$, hence $G/H=S^n$, are given by the following list:
\begin{itemize}
\item $SO(n)/SO(n-1)=S^{n-1}$
\item $SU(n)/SU(n-1)=U(n)/U(n-1)=S^{2n-1}$
\item $Sp(n)Sp(1)/Sp(n-1)Sp(1)=Sp(n)U(1)/Sp(n-1)U(1)=Sp(n)/Sp(n-1)=S^{4n-1}$
\item $G_2/SU(3)=S^6$
\item $Spin(7)/G_2=S^7$
\item $Spin(9)/Spin(7)=S^{15}$.
\end{itemize}
\end{proposition}
Observe that effectiveness is not a particularly stringent requirement since the action of any group $G$ on $M$ is effective when restricted to the quotient group $G/N$ given by the closed subgroup $\displaystyle N :=\bigcap_{x\in M} G_x$.

The goal of this section is to prove the following theorem:
\begin{theorem}\label{thm:surprise} 
Let $G$ be a compact Lie group acting multisymplectically, transitively and effectively on $S^n$ equipped with the standard volume form, then the action admits a \comoment if and only if $n$ is even. 
\end{theorem}
\begin{proof}
According to the classification given by proposition \ref{Prop:ClassificActionSpheres}, it will suffice to prove the following statements:
\begin{enumerate}
\item The action of $SU(n)$ on $S^{2n-1}$ does not admit a \comoment. As $SU(n)\subset U(n)\subset SO(2n)$, from this we will automatically get the statements, that $U(n)$ and $SO(2n)$ do not admit a \comoment when acting on $S^{2n-1}$. Moreover, as $SU(4)\cong Spin(6)\subset Spin(7)$, this implies that $Spin(7)$ does not admit a \comoment, when acting on $S^7$.
\item The action of $Sp(n)$ on $S^{4n-1}$ does not admit a \comoment. Hence, neither $Sp(n)U(1)$ nor $Sp(n)Sp(1)$ do. 
\item Spin(9) does not admit a \comoment, when acting on $S^{15}$.

\item $SO(2n+1)$ has a \comoment when acting on $S^{2n}$. As $G_2\subset SO(7)$, this implies that $G_2$ admits a \comoment when acting on $S^6$.
\end{enumerate}
Here we have employed the standard inclusions and isomorphisms coming from the well-known classification of compact connected Lie groups (see for example \cite{fulton2013representation,knapp2013lie}).

The theorem can, therefore, be proved resorting to Lemma \ref{lem:core}, after proving the following  Proposition \ref{prop:stiefel} and Proposition \ref{prop:annoying}. 
\end{proof}

\begin{proposition}\label{prop:stiefel}
Let $\omega_n$ be the volume form of $S^n$, $N$ the north pole and consider the orbit map $\vartheta_N$ of $N$ for a certain group acting on the sphere.
\begin{itemize}
\item Let $\vartheta_N:SU(n)\to S^{2n-1}$. Then $\vartheta_N^*[\omega_{2n-1}]\neq 0$.
\item Let $\vartheta_N:Sp(n)\to S^{4n-1}$. Then $\vartheta_N^*[\omega_{4n-1}]\neq 0$.
\item Let $\vartheta_N:Spin(9)\to S^{15}$. Then $\vartheta_N^*[\omega_{15}]\neq 0$.
\end{itemize}
\end{proposition}
\begin{proof}
 Starting from the first case, consider the principal bundle
 	\begin{displaymath}
		\begin{tikzcd}[row sep=0ex]
			SU(n-1) \ar[r,hook,"i"]
			& SU(n) \ar[r,two heads, "\vartheta_N"]
			& S^{2n-1}
		\end{tikzcd} 	
 	\end{displaymath}
 corresponding to the action $SU(n)\action S^{2n-1}$ with respect to the orbit of the north pole $N$ as introduced in remark \ref{Rem:ActionasBundle}.
 This bundle has a $(2n$-$2)$-connected base manifold, hence the pair $(SU(n),SU(n$-$1)$ is $(2n$-$2)$-connected and the pull-back maps $i^\ast: H^k(SU(n))\to H^k(SU(n$-$1))$ are isomorphism for $k\leq (2n$-$2)$ (see \cite{MR1867354} for extra details).
 The cohomology of $SU(n)$ is completely described by its generators, namely one has
	\begin{displaymath}
 		H^\bullet(SU(n),\mathbb{R})\simeq \Lambda_{\mathbb{R}}[u_3,u_5,\dots,u_{2n-1}]
	\end{displaymath}
	where the right-hand side denotes the exterior algebra generated by elements $u_k$ in odd degree $k$ (\cf \cite[Corollary 4D.3]{MR1867354}).
 In particular one can see that all generators $H^\bullet(SU(n-1))$ comes from the restriction of the generators of $H^\bullet(SU(n))$ , \ie are images under $i^\ast$.
 \note{Hint MZ: Since you mention this theorem several times, make sure you can explain it at the defense}
 This means that the Leray-Hirsch-theorem \cite[Theorem 4D.1]{MR1867354} can be applied, hence the following isomorphism holds:
 \begin{displaymath}
  	H^\bullet(S^{2n-1})\otimes H^\bullet ( SU(n-1)) \xrightarrow{\sim} H^\bullet(SU(n)) ~.	
 \end{displaymath}
 Inserting the explicit values for the real cohomology of the higher dimensional sphere, the left-hand side in degree $(2n-1)$ results in
 \begin{displaymath}
 	\mathclap{
  \left(H^{(2n-1)}(S^{2n-1})\otimes H^0(SU(n-1)\right)
 	\oplus 
 	\left(H^0(S^{2n-1})\otimes H^{2n-1}(SU(n-1))\right)
 	\cong 
 	H^{2n-1}(S^{2n-1})}
\end{displaymath}  
 and the Leray-Hirsch isomorphism restrict to $\vartheta_N^\ast$
 \begin{displaymath}
 	H^{2n-1}(S^{2n-1}) \xrightarrow{\vartheta_N^\ast}
 	H^{2n-1}(SU(n))
 \end{displaymath}
 Being the volume $[\omega_{2n-1}]$ the (degree $(2n-1)$) generator of the cohomology of $S^{2n-1}$, that means that $\vartheta_N^*[\omega_{2n-1}]$ is a generator of the cohomology of $SU(n)$ and thus non-zero in $H^{2n-1}(SU(n))$.
 
 The same reasoning holds almost verbatim to the second case. The expression for the generators of $H^\bullet(Sp(n))$ can be read again in \cite[Corollary 4D.3]{MR1867354}.
 
 Regarding the third case, one has only to observe that, $Spin(n)$ is the universal covering of $SO(n)$ and, in particular, the two groups are locally isomorphic, hence the real-valued cohomologies of $Spin(n)$ and $SO(n)$ are isomorphic (see \cite[\textsection 11.1]{Borel1955}).
 The cohomology groups can be explicitly given through the set of generators
 \begin{equation}\label{Eq:SOn-cohomology}
 	H^\bullet(Spin(k))
 	\simeq H^\bullet(SO(k))
 	\simeq
 	\begin{cases}
 		\Lambda_{\mathbb{R}}[a_3,a_5,\dots,a_{4n-1}]
 		&~ k = 2n+1
 		\\
 		\Lambda_{\mathbb{R}}[a_3,a_5,\dots,a_{4n-1},a'_{2n+1}]
 		&~ k = 2n + 2 		
 	\end{cases}
 \end{equation}
 where, as before, subscripts denote degrees (see, for instance, \cite[Theorem 3D.4]{MR1867354}).
 As before,the pair $(Spin(9),Spin(7))$ is 14-connected and the maps $H^i(Spin(9))\to H^i(Spin(7))$ are isomorphisms for $i\leq 13$. 
 Hence Leray-Hirsch-theorem holds and $\vartheta_\omega^\ast$ maps $[\omega_{15}]$ to the generator $a_{15}\in H^{15}(Spin(9))$ which is necessarily non-zero.
\end{proof}

In the case of $SO(2n+1)$, the Leray-Hirsch theorem can not be applied since is clear from equation \eqref{Eq:SOn-cohomology} that $SO(2n)$ has a class in degree $2n-1$ which does not come from any class in $SO(2n+1)$. 
In fact we have the following:
\begin{proposition}\label{prop:annoying}
Let $\omega_{2n}$ be the volume form of $S^{2n}$ and $N$ the north pole. Let $\vartheta_N:SO(2n+1)\to S^{2n}$. Then $\vartheta_N^*[\omega_{2n}]= 0$.
\end{proposition}
\begin{proof}
Let $i:SO(2n)\hookrightarrow SO(2n+1)$ be the inclusion. 
According to equation \eqref{Eq:SOn-cohomology}, cohomologies of $SO(2n+1)$ and $SO(2n)$ are isomorphic up to degree $2n-2$ since their generators are isomorphic through  $i^*$ up to the same degree.
Nevertheless, also $i^*:H^{2n}(SO(2n+1))\to H^{2n}(SO(2n))$ is an isomorphism because the first different generator is $a'_{2n-1}\in (SO(2n))$ and it cannot produce any element in degree $2n$ since the lowest degree generator has degree equal to $3$. 
\\
Observe than that the class $[i^*\vartheta_N^*\omega_{2n}]$ is the obstruction to 
 the existence of a \momap pertaining to the $SO(2n)$-action on $S^{2n}$. 
We know from Proposition \ref{lem:core}, that this action admits a \comoment, \ie $[i^*\vartheta_N^*\omega_{2n}]=0\in H^{2n}(SO(2n))$. 
But as $i^*$ is an isomorphism, this implies that $[\vartheta_N^*\omega_{2n}]=0\in H^{2n}(SO(2n+1))$.
\end{proof}

\note{
\textbf{Marco} How to get rid of effectiveness?\\
$\forall s: G \twoheadrightarrow	 SO(n+1)$ surjective Lie group morphism, the effective and transitive action $\theta:SO(n+1)\action S^{n}$ induce a transitive but generally non-effective action $\tilde{\theta}:G\action S^{n}$ given by $\tilde{\theta}= \theta \circ (\id\times s) : S^n \times G \to S^n$
\\
Effective means that $\hat{\theta}_g = \Id_M \Leftrightarrow g = e$ (\ie $\hat{\theta}:G\to \text{Diff}(M)$ is not injective). 

Defining
$N = \text{ker}(\theta) = \{g \in G\; \vert\: \hat{\theta}_g = \Id_M \}$
(subgroup of G) one can see that every non-effective action can be made effective, \ie
$$ \dfrac{G}{\text{ker}(\theta)} \twoheadrightarrow SO(n) \action S^{n}$$
is transitive and effective. 

In general, consider $G \twoheadrightarrow \frac{G}{N}$
(Check that $\frac{G}{N}$ is an honest Lie group. $N$ is always a closed normal subgroup and if $G$ is closed it should be ok)
you get $j:\mathfrak{g}\twoheadrightarrow \frac{\mathfrak{g}}{\mathfrak{n}}$, in general not injective.
If $\frac{\mathfrak{g}}{\mathfrak{n}}$ admits a \comoment you get a \comoment for $\mathfrak{g}$ via pre-composition with $j$

If you can give an example of \comoment for $\mathfrak{g}$ which does not descend from a \comoment of $\frac{\mathfrak{g}}{\mathfrak{n}}$ you should be able to complete the argument.
}

\subsection{Explicit \comoment for $SO(n+1)$ on $S^{n}$}\label{subsectra}
Giving an explicit expression for a \comoment of $SO(n+1) \action S^{n}$ requires to find iteratively, for $k \in (1,\dots,n-1)$ and  for all $p \in \Lambda^k \mathfrak{so}(n)$, a primitive, denoted as $f_k(p)$, of the closed differential $(n-k)$-form
\begin{equation}\label{eq:comomentasprimitive}
	\mu_k (p) = 	-f_{k-1} (\partial p) - \varsigma(k) \iota(v_p) \omega 
\end{equation}
with $f_0 = 0$ and satisfying the following constraint
\begin{equation} \label{eq:comomentconstraints}
	f_n(\partial p ) = - \varsigma(k+1) \iota_{v_p} \omega ~.
\end{equation}

As $H^0(S^n)$ and $H^{n}(S^n)$ are the only non-trivial cohomology groups, it is always possible to find primitives of $\mu_k (p)$.
The only thing that could fail, and actually fails when $n$ is odd, is the fulfilment of condition \eqref{eq:comomentconstraints}.
In the latter case, it is however possible to consider an extension of $\mathfrak g$ to a suitable Lie-$n$ algebra and consider \emph{Lie-$n$ \momaps} instead of our notion of \comoments. (See \cite{Callies2016} or \cite{Mammadova2019} for the explicit construction in the $2$-plectic case. 
The latter applies to the action $SO(4)\action S^3$ in the situation we are considering here.)

\vspace{1em}
Instead of dealing with the analytical problem of finding explicit potentials for the form $\mu_k (p)$, let us translate the problem in a more algebraic fashion focusing on the particular structure of the Chevalley-Eilenberg complex of $\mathfrak{so}(n)$.
\\
In general, it is fairly easy to express the action of a \comoment on boundaries:
\begin{lemma}[\Momaps on boundaries]
	Let $v:\mathfrak{g}\to (M,\omega)$ be a multisymplectic infinitesimal action.\\
	Let $F^k: B_k(\mathfrak{g}) \to \Lambda^{k+1}\mathfrak{g}$ such that $\partial \circ F^k = \id_{B^k}$, 
	\ie $F^k$ gives a choice of representative of a primitive for every $k$-boundary of $\mathfrak{g}$.\\
	Then, the function $ 	f_k \colon B_k(\mathfrak{g}) \to\Omega^{n-k}(M)$ defined as
	\begin{displaymath}
		f_k(p) = -\varsigma(k+1) \iota(v_{F^k(p)}) \omega
	\end{displaymath}
	satisfies equation (\ref{eq:fk_hcmm}) defining the $k$-th component for a \comoment of the infinitesimal action, for every boundary $p$.
\end{lemma}
\begin{proof}
	It is a straightforward application of Lemma \ref{lemma:multicartan} together with the multisymplecticity of the infinitesimal action:
	\begin{align*}
		d f_k (p) & = 
		-\varsigma(k+1) d \iota(v_{F^k(p)}) \omega =
		\\
		&=
		-(-1)^{k+1}\varsigma(k+1) \iota(v_{\partial F^k(p)}) \omega =\\
		& = \cancel{- f_{k-1}(\partial p)} - \varsigma(k) \iota_{v_p} \omega.
	\end{align*}
\end{proof}
\begin{remark}
Note that the constraint given by equation (\ref{eq:comomentconstraints}) is precisely the requirement that action on boundaries of the highest component $f_n$ of the \comoment $(f)$ is independent from the choice of $F^n$.
\end{remark}
In other words, finding the action of the component $f_k$ of the \comoment $(f)$ on boundaries is tantamount to finding a function $F^k: B_k(\mathfrak{g}) \to \Lambda^{k+1}\mathfrak{g}$ mapping a boundary $p$ to a specific primitive $q$.
\\
Note that this is nothing more than replacing the problem of finding a potential of an exact differential form to the one of finding a primitive of a boundary in the  CE-complex.

It follows from the previous lemma that the $k$-th component of the \comoment is completely determined by its action on boundaries when the $k$-th Chevalley-Eilenberg homology group vanishes:
\note{Hint Leo: Maybe connection to leyli (concerning hermann sfuff)}
\begin{corollary}
	Consider a Lie algebra with $H_k(\mathfrak{g})=0$ and fix a choice of representatives $F^k$ and $F^{k-1}$ as before.
	Then, the function $ 	f^k \colon \Lambda^k \mathfrak{g} \to\Omega^{n-k}(M)$ defined as
	\begin{displaymath}
	f_k(q) = \varsigma(k+1)\iota(v_{F^k(F^{k-1}(\partial q) -q)})\omega
	\end{displaymath}
	satisfies equation (\ref{eq:fk_hcmm}) for every chain $q \in \Lambda^k\mathfrak{g}$.
\end{corollary}
\begin{proof}
	Consider a function $f_{k-1}$ defined by means of $F^{k-1}$ according to the previous lemma.
	For every cycle $q \in \Lambda^k\mathfrak{g}$ we get
	\begin{displaymath}
		\begin{split}
		- f_{k-1}(\partial q) -\varsigma(k)\iota_{v_q}\omega =&~
		\varsigma(k)[\iota(v_{F^{k-1}(\partial q)}) - \iota(v_{q})] \omega=
		\\
		=&~
		\varsigma(k)\iota( v_r )\omega		
		\end{split}
	\end{displaymath}
	where $r =(F^{k-1}(\partial q) - q)$ is closed, hence exact.
	Again from Lemma \ref{lemma:multicartan} follows that the right-hand side it is equal to
	\begin{displaymath}
		\begin{split}
		\varsigma(k)\iota(v_{\partial F^k(r)} \omega)	=&~
		\varsigma(k)(-1)^{k+1} d \iota(v_{F^k(r)}) \omega =
		\\
		=&~
		d f_k(q)
		~.
		\end{split}
	\end{displaymath}
\end{proof}

An explicit construction of a \comoment is generally more delicate in presence of cycles that are not boundaries.
Nevertheless, we know from equation \eqref{Eq:SOn-cohomology} that the first two homology group of $\mathfrak{so}(n)$ are trivial, therefore it is easy to give the first two components of the \comoment.
From the linearity of $f_k$, it is clear that we only need to give its action on the standard basis of the finite-dimensional vector space $\mathfrak{so}(n)$.
\begin{remark}[Standard basis of $\mathfrak{so}(n)$]
Recall that $\mathfrak{so}(n)$ is the Lie sub-algebra of $\mathfrak{gl}(n,\mathbb{R})$ consisting of all skew-symmetric square matrices. A basis can be constructed as follows:
\begin{equation}\label{eq:standard-basis}
	\mathcal{B}\coloneqq \big\lbrace 	A_{a b} = (-1)^{1+a+b} \left( E_{a b} - E_{b a}\right)
	\quad \vert \quad 1\leq a<b\leq n \big\rbrace
\end{equation}
where $E_{a b}$ is the matrix with all entries equal to zero and entry $(a,b)$ equal to one.
\\
The fundamental vector field of $A_{a b}$ associated to the linear action of $SO(n)$ on $\mathbb{R}^n$ reads as follows:
\begin{displaymath}
	v_{A_{a b}}= \sum_{i,j}[A_{a b}]_{i j}x^j \partial_i  = (-1)^{1+a+b}\left(x^a \partial_b - x^b \partial_a\right)
\end{displaymath}
\noindent Using such a basis, the structure constants read as follows:

	\begin{displaymath}
	\begin{aligned}
		[A_{a b}, A_{c d}] =& (-1)^{(b+c+1)}\delta_{b c} A_{a d} +
		(-1)^{(a+d+1)}\delta_{a d} A_{b c} +\\
		&	(-1)^{(d+b+1)}\delta_{d b} A_{a c} +
		(-1)^{(a+c+1)}\delta_{c a} A_{d b}	
	\end{aligned}
	\end{displaymath}
		in particular:
	\begin{equation}\label{eq:reductionformula}
		[A_{k a}, A_{k b}] = A_{a b} \qquad \forall k \neq a,b ~.
	\end{equation}
\end{remark}

\paragraph{$f_1$ for any $SO(n)$.}

Since all 1-chains in the CE complex are automatically cycles, $H^1(\mathfrak{so}(n))=0$ implies that all elements of $\mathfrak{so}(n)$ are boundaries.
\\
Formula \eqref{eq:reductionformula} suggests a natural choice for the function $F^1$ when acting on elements of the standard basis:
\begin{equation}\label{eq:primitiveMapF}
	F^1(A_{a b}) = -\sum_{k=1}^n \dfrac{1}{n-2}  A_{k a}\wedge A_{k b} ~.
\end{equation}
Therefore the first component of the \comoment is given by
\begin{displaymath}
	\begin{split}
	f_1 (A_{a b}) =&~ - \iota(v_{F^1(A_{a b})}) \omega 
	=
	\\
	=&~
	\dfrac{1}{n-2}\sum_{k=1}^n \iota(v_{A_{k b}})\iota(v_{A_{k a}})\omega~.
	\end{split}
\end{displaymath}
\begin{example}
In the three-dimensional case, denoting the three generators of $\mathfrak{so}(3)$ as $l_x,l_y,l_z$, namely

\begin{displaymath}
\mathclap{
	l_x = A_{1\, 2} = \begin{bmatrix} 0 & 1 & 0 \\ -1 & 0 & 0 \\ 0 & 0 & 0 \end{bmatrix} \quad
	l_y = A_{1\, 3} = \begin{bmatrix} 0 & 0 & -1 \\ 0 & 0 & 0 \\ 1 & 0 & 0 \end{bmatrix} \quad
	l_z = A_{2\, 3} = \begin{bmatrix} 0 & 0 & 0 \\ 0 & 0 & 1 \\ 0 & -1 & 0 \end{bmatrix}
	~,
}
\end{displaymath}
 one gets
\begin{displaymath}
	F^1 (l_i) = - \dfrac{1}{2}\sum_{j,k=1}^3\epsilon_{i j k} l_j \wedge l_k,
\end{displaymath}
where $\epsilon_{i j k}$ is the Levi-Civita symbol, and
\begin{displaymath}
	\begin{split}
	f_1(l_z) =&~ \iota(v_{l_y \wedge l_z}) \omega =
	\\
	=&~
	j^\ast \iota(E \wedge v_{l_y}\wedge v_{l_z}) dx^{123}.
	\end{split}
\end{displaymath}
\end{example}

\paragraph{$f_2$ for any $SO(n)$.}
In this case there are two subsets of generators of $\Lambda^2 \mathfrak{so}(n)$ to consider:
\begin{displaymath}
	\begin{cases}
    p = A_{a b} \wedge A_{c d} \xmapsto{\,\partial\,} 0 & \quad\text{for } a,b,c,d \text{ different} \\
    q = A_{j a} \wedge A_{j b} \xmapsto{\,\partial\,} -A_{a b} & \quad\text{for } j,a,b \text{ different} 
  \end{cases}
\end{displaymath}
The elements of the first set are boundaries and a primitive can be given as follows
\begin{displaymath}
F^2 (A_{a b} \wedge A_{c d}) = \dfrac{n-2}{4}\left( F^1(A_{a b}) \wedge A_{c d} - A_{a b}\wedge F^1(A_{c d}) \right)
\end{displaymath}
where $F^1$ is given again by equation \eqref{eq:primitiveMapF}.
In the second case, we need to find a primitive of
\begin{displaymath}
	\begin{aligned}
		 F^1(\partial q) -q &= -F^1(A_{a b}) - A_{j a} \wedge A_{j b} 
		 =
		 \\
		 &= \dfrac{1}{n-2}\sum_{k=1}^n( A_{k a}\wedge A_{k b} - A_{j a} \wedge A_{j b})
		 =\\
		 &=\dfrac{1}{n-2}\sum_{k=1}^n \partial (A_{k a}\wedge A_{j b}\wedge A_{k j}) 
		 =\\
		 &=\partial \left(\dfrac{1}{n-2}\sum_{k=1}^n (A_{k a}\wedge A_{j b}\wedge A_{k j}) \right)
	\end{aligned}
\end{displaymath}
The last equality suggests the following choice
\begin{displaymath}
 F^2(F^1(\partial q) -q) = \left(\dfrac{1}{n-2}\sum_{k=1}^n (A_{k a}\wedge A_{j b}\wedge A_{k j}) \right).
\end{displaymath}
Finally, one gets
\begin{displaymath}
\begin{aligned}
	f_2(A_{a b} \wedge A_{c d}) &= \dfrac{n-2}{4}\left(
	\iota(v_{F^1(A_{a b}) \wedge A_{c d}}) - \iota(v_{A_{a b}\wedge F^1(A_{c d})})\right)\omega\\
	f_2(A_{j a} \wedge A_{j b}) &= \dfrac{-1}{n-2}
	\sum_{k=1}^n	
	\left(\iota(v_{A_{k a}\wedge A_{j b}\wedge A_{k j}})\right)\omega.
\end{aligned}
\end{displaymath} 

\note{
	\textbf{$f_k$ for any $SO(n)$ and $k\geq 3$}:\\
	We know from Theorem \ref{thm:son-cohomology} that $H^3(\mathfrak{so}(n))$ never vanishes... how can we proceed? Tony drafted a code in Python. The repository is still private
}

\subsection{Explicit \comoment for $G_2$ on $(S^6,\phi)$}
We finish by providing a nice example of \comoments for non-volume forms on spheres.\\

Recall that	$G_2$ is a subgroup of $SO(7)$ acting transitively and multisymplectically on $S^6$ with the standard volume. Therefore, according to Theorem \ref{thm:surprise}, the action $G_2\action (S^6,\omega)$ admits a \comoment.
	
This group can be explicitly defined as the subgroup of $GL(7,\mathbb{R})$ preserving the multisymplectic $3$-form
\begin{equation}
	\phi
	=d x^{123}+ d x^{145}+ d x^{167}+ d x^{246}- d x^{257}- d x^{356}-d x^{347},
\end{equation}
where~$x = (x^i)$ are the standard coordinates on~$\mathbb{R}^7$ 
and~$d x^{ijk} = d x^i \wedge d x^j \wedge d x^k$. 
(See \cite{MR1939543} for further remarks on $G_2$-homogeneous multisymplectic forms and \cite{MR2253159} for details on the $G_2$-manifold $S^6$).

Considering the multisymplectic structure $j^\ast\phi$ on $S^6$, where $j$ is the inclusion of the sphere in $\mathbb{R}^7$, instead of the standard volume, it is possible to give an explicit \comoment for the action of $G_2$:
\begin{lemma}
	The action $G_2 \circlearrowright (S^6, j^\ast\phi)$ admits an equivariant \comoment given by:
	\begin{displaymath}
		\morphism{f_k}
		{\Lambda^k\mathfrak{g}_2}
		{\Omega^{2-k}(M)}
		{q}
		{\dfrac{(-1)^{k-1}}{3}~j^\ast\left(\iota_{v_q}~\iota_E~ \phi\right)}
		~,
	\end{displaymath}
for $(k=1,2)$, where $E=x^i \partial_i \in \mathfrak{X}(\R^7)$ is the Euler vector field.
\end{lemma}	
\begin{proof}
	It follows from Lemma \ref{lem:extexact}, noting that $(\frac{1}{3} \iota_E \phi)$ is a $G_2$ invariant primitive of $\phi$ in $\mathbb{R}^3$.
\end{proof}

\note{
Hint Marco: if there exists a multisymplectic 3-form $\beta$ such that $(j^*\phi) \wedge \beta = j^\ast\iota_E Vol=\omega$ and if you have an explicit \comoment for $G_2 \action (S^6, \beta)$ you can cook up a \comoment for $G_2 \action (S^6, \omega)$ resorting to the main theorem of \cite{Shahbazi2016}.\\
}

\ifstandalone
	\bibliographystyle{../../hep} 
	\bibliography{../../mypapers,../../websites,../../biblio-tidy}
\fi

\cleardoublepage


%% file: chapters/higherrogersembedding/higherrogersembedding.tex
\chapter{Gauge transformations and \momaps}\label{Chap:MarcoPaper}
A natural theme that arises when dealing with both symplectic and multisymplectic structures is to investigate what relationship exists between gauge-related multisymplectic manifolds, \ie manifolds endowed with multisymplectic forms lying in the same cohomology class (see section \ref{Section:GaugeTransformations} in chapter \ref{Chap:MultiSymplecticGeometry}).

To date, no canonical correspondence between the $L_\infty$-algebras of observables of two gauge related multisymplectic manifolds is known\footnote{Clearly they are not isomorphic as graded vector spaces. In particular, they differ in their degree $0$ component $\Omega^{n-1}_{\ham}(M,\omega)$.}.
In this chapter, we will exhibit a compatibility relation between those observables that are momenta of corresponding homotopy moment maps (the higher analogues of a moment map in the multisymplectic setting).
Although this construction is essentially algebraic in nature, it admits also a geometric interpretation when declined to the particular case of pre-quantizable symplectic forms. 
For a symplectic manifold $(M,\omega)$, a choice of prequantization circle bundle with connection gives a Lie algebra embedding
of the Poisson algebra of functions on $M$ into the invariant vector fields on the prequantization  bundle. 
The latter form the sections of a Lie algebroid over $M$, called Atiyah algebroid, which is isomorphic to a central extension $(TM\oplus \RR)_{\omega}$ of the tangent bundle $TM$. Hence we obtain a Lie algebra embedding 
\begin{equation}\label{eq:intropreqmap}
  C^{\infty}(M)\to \Gamma((TM\oplus \RR)_{\omega}),
\end{equation} which we call \emph{prequantization map}.

Assume $(M,\omega)$ is endowed with a moment map for the action of some Lie group $G$, whose corresponding comoment map we denote $J^*\colon\g\to C^{\infty}(M)$.
Any choice of $\alpha\in \Omega^1(M)^G$ provides another $G$-invariant symplectic form $\omega+d\alpha$ (assuming this is non-degenerate), and $\alpha$ can be used to twist some of the above data, obtaining in particular a moment map   $J_{\alpha}$  for $\omega+d\alpha$. A geometric argument shows that the following diagram   commutes:
	\begin{equation}
	\label{intro:diag:main}
		\begin{tikzcd}[column sep=huge]
			&
			  C^{\infty}(M)_{\omega}   \ar[r,""] 
			&
			 \Gamma(TM\oplus \RR)_{\omega}  \ar[dd,"\tau_\alpha"]
			\\[-1em]
			\mathfrak{g}\ar[ru,"J^*"] \ar[dr,"J^*_{\alpha}"']
			\\[-1em]
			&
			 C^{\infty}(M)_{\omega+d\alpha} \ar[r,""] 
			&
			\Gamma(TM\oplus \RR)_{\omega+d\alpha}
		\end{tikzcd}
	\end{equation}
Here $\tau_{\alpha}$ is the  gauge transformation of the Atiyah algebroids induced by $\alpha$, which, in particular, is a Lie algebroid isomorphism.
We interpret this commutativity by saying that  the twisting {of the moment map by $\alpha$} is compatible with the twisting of the {Atiyah algebroid}. 
\\
The prequantization bundle over $(M,\omega)$  exists only when $[\omega]$ satisfies an integrality condition.
Notice that the above Lie algebra embedding \eqref{eq:intropreqmap} and the commutative diagram \eqref{intro:diag:main} make no reference to the prequantization bundle. 
Indeed, one can check that they make sense and hold for any arbitrary symplectic manifold $(M,\omega)$. 

\bigskip
In this chapter, based on joint work with Marco Zambon \cite{Miti2020}, 
we show that the existence of the above Lie algebra embedding and commutative diagram extends to the setting of higher geometry, \ie replacing $\omega$ by a  multisymplectic $(n+1)$-form (no integrality condition is required). 
In that case the Poisson algebra of functions $C^{\infty}(M)$ is replaced by a $L_{\infty}$-algebra \cite{LadaMarkl}\cite{LadaStasheff}, the Atiyah Lie algebroid   by the higher Courant algebroid $TM \oplus \Lambda^{n-1} T^\ast M$, and $\alpha$ by an invariant $n$-form $B$.
\\
Our previous discussion around diagram \eqref{intro:diag:main} provides some evidences that this construction may be related to the higher analogue of geometric quantization for integral multisymplectic forms.

Our main results are the following.
\begin{itemize}
\item In theorem   \ref{thm:iso} {and corollary \ref{cor:Psi}} we establish the existence of the embedding for $n\le 4$. The method is based on the description of $L_{\infty}$-algebras as suitable coderivations, and we expect the proof to extend to the case of arbitrary $n$.
\item Building on this, in theorem \ref{thm:comm}, we establish that the higher version of the above diagram \eqref{intro:diag:main} commutes, for $n\le 4$.
  \begin{equation}
  \label{intro:eq:pentagonDiagram}
	\begin{tikzcd}
		&
		L_{\infty}(M,\omega) \ar[r]
		&
		L_{\infty}(TM \oplus \Lambda^{n-1} T^\ast M,\omega) \ar[dd,"\tau_B"]
		\\[-1em]
		\mathfrak{g}\ar[ru]
		 \ar[dr]
		\\[-1em]
		&
		L_{\infty}(M,\widetilde{\omega}) \ar[r]
		&
		L_{\infty}(TM \oplus \Lambda^{n-1} T^\ast M,\widetilde{\omega})
	\end{tikzcd}
\end{equation}
\end{itemize}

We relate our first result above with the literature.
 In the special case $n=2$, the Atiyah algebroid is an instance of Courant algebroid, and the embedding was established by Rogers \cite[Theorem 7.1]{Rogers2013}. 
 For arbitrary $n$ the existence of the embedding is stated by S\"amann-Ritter in their preprint \cite[Theorem 4.10.]{Ritter2015a}. 
 They provide a proof in which the embedding is constructed recursively, but not all steps are worked out explicitly. 
They do not give a closed formula for the embedding. 
For a different approach in the case of integral multisymplectic forms, involving a choice of open cover on the manifold $M$, see   Fiorenza-Rogers-Schreiber 
 \cite[\S 5]{Fiorenza2014a}.
\\ 
 {Our second result above, to the best of our knowledge, has not been addressed in the literature yet}. 
\\
{We expect to be able to extend both results above to arbitrary values of $n$. The proof of theorem \ref{thm:comm} suggests that  theorem   \ref{thm:iso} can be extended by choosing the coefficients there to depend suitably on the Bernoulli numbers, and in that case theorem \ref{thm:comm} would hold for all $n$ as a consequence. 
}

\medskip
\noindent The outline of the chapter is as follows.
\\
First, in section \ref{Sec:GeoMotiv}, we describe the possible motivation for studying the commutation of a diagram like \eqref{intro:eq:pentagonDiagram} coming from geometric quantization.
\\
Then, we move to prepare the ground for delivering our explicit construction for the embedding of the $L_\infty$-algebra associated to a multisymplectic manifold into the $L_\infty$-algebra associated to the corresponding Vinogradov algebroid.
Our construction is based on the observation that both structures - when restricted to a suitable subcomplex - can be seen as being generated by the same set of few multibrackets. 
Accordingly, in section \ref{Sec:RogersL1prop} we review again the construction of the Rogers $L_\infty$-algebra of observable, discussed in section \ref{Section:RogersObservables}, with the language of \RN products.
\\
In section \ref{Sec:Vinoids} we present the notion of Vinogradov algebroid and of the corresponding $L_\infty$-algebra working out certain recurrence properties of multibrackets that are similar to what we found for the Rogers' $L_\infty$-algebra of observables.
\\
In section \ref{Section:ExtendedRogersEmbedding} we reap the fruits of this preliminary work providing an explicit expression of the sought embedding of the Rogers' $L_\infty$-algebra into the Vinogradov's one.
\\
Finally, in section \ref{Sec:DiagramGaugeTransf} we discuss the compatibility of the previous construction under gauge transformations in presence of symmetries admitting \momaps.

\section{Geometric motivation: the symplectic case}\label{Sec:GeoMotiv}
In this section, we explain in detail the considerations on symplectic geometry outlined in the first part of the introduction. 
\begin{remark}
On the  symplectic manifold  $(M,\omega)$ we adopt the conventions that $\iota_{X_f}\omega=-df$ and $\{f,g\}=\omega(X_f,X_g)$ (hence $f\mapsto X_f$ is a Lie algebra morphism). To shorten the notation, we denote by  $C^{\infty}(M)_{\omega}$ the Lie algebra $(C^{\infty}(M),\{\cdot,\cdot\})$.
\end{remark} 

We first introduce  in \S \ref{subsec:prequantizationReminder} the notion of "prequantization" needed to discuss the geometric interpretation of the embedding \eqref{eq:intropreqmap} and the diagram \eqref{intro:diag:main}.
The latter will be explained in \S \ref{subsec:first} and \S \ref{subsec:second}.
Finally, in \S \ref{subsec:higher} we provide some motivation for the higher case.

\subsection{Reminders on Atiyah algebroids and Prequantization}\label{subsec:prequantizationReminder}
	In this subsection, we succinctly review the language required to deliver a geometric interpretation of the -purely algebraic- problem of ascertaining the commutativity of diagram \eqref{intro:diag:main}.
	We will need the notion of Atiyah algebroid and prequantization.
	The former is a certain Lie algebroid uniquely associated to any principal bundle and the latter is a construction involving $S^1$-principal bundles on certain symplectic manifolds.

	\subsubsection{Lie algebroids}
	Informally, Lie algebroids are infinite dimensional Lie algebras controlled by a geometric data, namely, elements are sections of given vector bundle.
	\begin{definition}[Lie algebroid]
		We call \emph{Lie algebroid} a triple $(E,[\cdot,\cdot]_E,\rho)$ consisting of
		\begin{itemize}
			\item a vector bundle $\pi:E \twoheadrightarrow M$;
			\item a Lie algebra structure $[\cdot,\cdot]_E$ on the space of section $\Gamma(E)$;
			\item a vector bundle morphism $\rho:E\to TM$ (over the identity $\id_M$), called \emph{anchor};
			\end{itemize}
		such that:
		\begin{enumerate}
			\item[(a)] $\rho$ induces a Lie algebra morphism at the level of sections\footnote{Note that this condition follows directly from condition (b).}
				 $$\rho: (\Gamma(E),[\cdot,\cdot]_E)\to (\mathfrak{X}(M),[\cdot,\cdot])~;$$
			\item[(b)] $[\cdot,\cdot]_E$ is compatible with the anchor in the sense of the \emph{Liebniz rule}\footnote{Algebraically, it tells that $[X,\cdot ]$ is a derivation on the $C^\infty(M)$-module of sections $\Gamma(E)$.}:
				\begin{displaymath}
					[X, f~Y ]_E = (\rho(X) f)~ Y + f~ [X,Y]_E
					\qquad \forall f \in C^\infty(M);~X,Y\in \Gamma(E)~.
				\end{displaymath}
				($\rho(X)$ in the above equation has to be interpreted as the unique derivation on the associative commutative algebra $C^\infty(M)$, with respect to the point-wise product, associated to the smooth vector field $\rho(X)\in\mathfrak{X}(M)$.)				
		\end{enumerate}
	\end{definition}
	\begin{notation}
		If the anchor is a surjective map $ \rho:E \twoheadrightarrow TM$, the Lie algebroid is said to be  \emph{transitive}.	
	\end{notation}	
	
	\begin{example}[Lie Algebras]
		Any Lie algebra $\mathfrak{g}$ can be seen as a Lie algebroid over a $0$-dimensional manifold (\ie a point) $\{\ast\}$.
		In this sense, a Lie algebroid is a "many points version" (see "horizontal categorification" \cite{nlab:horizontal_categorification}) of a Lie algebra.
	\end{example}
	\begin{example}[Tangent bundle]
		Given a manifold $M$, the corresponding tangent bundle is Lie algebroid with anchor given by the identity bundle map.
	\end{example}

	\begin{example}[Standard Lie algebroid]\label{ex:StandardLieAlgbroid}
		Given any smooth manifold $M$, the vector bundle $E= \RR_M\oplus TM$ together with the standard projection $\rho: \RR_M\oplus TM \twoheadrightarrow TM$ and the binary bracket
		\begin{displaymath}
			\morphism{[\cdot,\cdot]}
			{\Gamma(TM \oplus\RR_M)\otimes \Gamma(TM \oplus\RR_M)}
			{\Gamma(TM \oplus\RR_M)}
			{\pair{X_1}{f_1}\otimes\pair{X_2}{f_2}}
			{\pair{[X_1,X_2]}
			{\mathcal{L}_{X_1} f_2 - \mathcal{L}_{X_2} f_1}}
		\end{displaymath}
		constitute a Lie algebroid called \emph{standard Lie algebroid}.
	\end{example}
	\begin{example}[$\omega$-Twisted (standard) Lie algebroid]\label{ex:TwistedLieAlgbroid}
		Consider a smooth manifold $M$.	
		Let be $\omega\in\Omega^2(M)$	 a closed $2$-form , \ie the manifold $(M,\omega)$ is a pre-$1$-plectic manifold.
		Then the vector bundle $E= \RR_M\oplus TM$ together with the standard projection $\rho: \RR_M\oplus TM \twoheadrightarrow TM$ and the binary bracket
		\begin{displaymath}
			\morphism{[\cdot,\cdot]_{\omega}}
			{\Gamma(TM \oplus\RR_M)\otimes \Gamma(TM \oplus\RR_M)}
			{\Gamma(TM \oplus\RR_M)}
			{\pair{X_1}{f_1}\otimes\pair{X_2}{f_2}}
			{\pair{[X_1,X_2]}
			{\mathcal{L}_{X_1} f_2 - \mathcal{L}_{X_2} f_1- \omega(X_1,X_2)} }
		\end{displaymath}
		constitute a Lie algebroid called \emph{$\omega$-twisted (standard) Lie algebroid}.
	\end{example}

\subsubsection{Principal connections and Atiyah algebroids}
		Given a finite dimensional Lie group $G$, consider a principal bundle $G\hookrightarrow P \twoheadrightarrow M~.$	
		We denote by $\hat{\xi}\in \mathfrak{X}(P)$ the fundamental vector field of $\xi \in \mathfrak{g}$ with respect to the action $R:G\action P$ of the group on $P$ from the right, namely
		\begin{displaymath}
			\hat{\xi}\eval_p := \dfrac{\d}{\d t} R_{\exp{t \xi}}(p) \eval_{t=0}
			~;
		\end{displaymath}
		(\cf remark \ref{rem:RightActionMess}).
		Le us recall the following two definitions:
	\begin{definition}[Connection of a Principal bundle]\label{def:connections}
		Given a $G$-principal bundle $G\hookrightarrow P \twoheadrightarrow M$	, we call a \emph{connection on $P$} any $\mathfrak{g}$-valued differential $1$-form $\theta\in \Omega^1(P,\mathfrak{g})$ satisfying the two following property
		\begin{itemize}
			\item $\theta$ reproduces the fundamental vector fields, \ie $\theta(\hat{\xi})= \xi$;
			\item $\theta$ is $G$-equivariant: $R_g^\ast \theta = \text{Ad}_{g^{-1}}\theta$.
		\end{itemize}
	\end{definition}

	\begin{definition}[Curvature $2$-form]
		Given a connection $\theta\in \Omega^1(P,\mathfrak{g})$, we call \emph{curvature of $\theta$} the $\mathfrak{g}$-valued differential 2-form $\d_{\cA}\theta \in \Omega^2(P,\mathfrak{g})$ such that
		\begin{displaymath}
			\d_{\cA}\theta (v_1,v_2) 
			= 
			\d\theta (v_1,v_2) + [\theta(v_1),\theta(v_2)]
			\qquad \forall v_i \in \mathfrak{X}(P)~.
		\end{displaymath}
		The connection $\theta$ is said \emph{flat} if $\d_{\cA}\theta = 0$.
	\end{definition}
	
	Principal connections are special cases of Ehresmann connections:
	\begin{reminder}[Ehresmann connections]\label{rem:ehresmann}
	Let $\pi:E\to M$ be a smooth fibre bundle. 
	Consider the tangent bundle $\tau:TE\to E$.
	One can introduce the unique \emph{vertical bundle} $V:=\ker(\d\pi:TE\to TM)$ where $\d \pi$ is the differential of the smooth map $\pi$ (the tangent map $T\pi$).
	$V$ is a subbundle of $TE$ whose fibres $V_e=T_e (E_{\pi(e)})$ consist of vectors on $TE$ which are tangent to the fibres of $E$.
	\\
	An \emph{Ehresmann connection} of $E$ is any smooth subbundle $H$ of $E$ such that 
	$$ TE = H \oplus V~.$$
	\\
	Consider now a principal bundle $\pi:P\to M$ with connection encoded by a $1$-form $\theta$ as in definition \ref{def:connections}. The subbundle 
	$H_{\theta}:=\ker(\theta)$ defines a invariant Ehresmann connection on $P$ corresponding by $\theta$. Hence $TP= \ker(\theta)\oplus \ker(\d \pi)$.
\end{reminder}

	\begin{reminder}[Horizontal lifts]\label{rem:horlifts}
		Let be $\pi:P\to M$ a vector bundle and consider an Ehresmann connection encoded by an horizontal subbundle $H$ of $P$.
		Vector fields on $P$ decompose uniquely in a vertical and a horizontal part.
	Given a vector field $X\in \mathfrak{X}(M)$, we call the \emph{horizontal lift} of $X$ the unique (horizontal) vector field $X^H\in \mathfrak{X}(P)$ such that the following diagram commutes in the category of smooth manifolds:
	\begin{displaymath}
		\begin{tikzcd}[column sep = huge]
			P\ar[r,bend left=30,"X^H"] \ar[d,two heads, "\pi"']&
			H \ar[l,two heads] \ar[d,dashed] \ar[r,hook] &
			TP \ar[dl,"\d\pi"]
			\\
			M \ar[r,bend right=30,"X"'] & TM \ar[l,two heads] 
		\end{tikzcd}
		~.
	\end{displaymath}	
	
	\end{reminder}

		Consider a $G$-principal bundle $G\hookrightarrow P \twoheadrightarrow M$, denoted simply as $\pi:P\to M$. 
		The action $R:G\action P$ is free and proper, therefore, there is a well-defined quotient.
		The same property holds for the lift of the action $R$ to $G\action TP$.
		In other words, we have the following commutative diagram in the category of smooth manifold
		\begin{displaymath}
			\begin{tikzcd}
				& TP \ar[r,"\d\pi"]\ar[dd,"\tau_P"] & TM \ar[dd,"\tau_M"] 
				\\[-1.5em]
				G \ar[ru,hook] \ar[rd,hook] & &
				\\[-1.5em]
				& P \ar[r,two heads,"\pi"] & M 
				~,
			\end{tikzcd}
		\end{displaymath}
		where the vertical arrows $\tau_P$ and $\tau_M$ denote the canonical fibrations of the tangent bundles over their corresponding base manifolds.
		Observe that, while the lower horizontal line is a $G$-principal bundles, the upper one never inherits the structure of $G$-principal bundle since the Lie group does not act transitively on the fibers of $\d \pi$.
		\\
		This justify the following definition:
		\begin{definition}[Atiyah algebroid]\label{def:ati}
		Given a $G$-principal bundle $\pi:P\to M$, 
		we call \emph{Atiyah algebroid} the Lie algebroid obtained by taking the quotient, with respect of $G$, of the tangent bundle $\tau_P: TP \to P$.
		Namely, it is given by the vector bundle
		\begin{displaymath}
			\begin{tikzcd}
				A_P \cong \dfrac{TP}{G} \ar[r,two heads,"\tau_P"] & \dfrac{P}{G}\cong M
			\end{tikzcd}
		~,
		\end{displaymath}
		whose sections correspond to $G$-invariant vector fields over $P$, \ie $\Gamma(A_P)\cong \mathfrak{X}(P)^{G}$, together with the restriction of the standard Lie bracket on $\mathfrak{X}(P)$ to $G$-invariant vector fields and with anchor given by $\d\pi: TP \to TM$.
		\end{definition}
		The upshot is that Atiyah algebroids are certain Lie algebroids naturally associated to principal bundles.
		
		The following lemma gives another characterization of the Atiyah algebroid that will be used in the following.
		\begin{lemma}[Atiyah exact sequence {(\cite[Thm.1]{Atiyah1957})}]\label{lem:atiexactseq}
			Given a principal bundle $\pi:P\to M$, its corresponding Atiyah algebroid $A_P$ fits in a short exact sequence in the category of Lie algebroids:
			\begin{displaymath}
				\begin{tikzcd}
					0 \ar[r] & P \times_G \mathfrak{g} \ar[r] & A_P \ar[r,"\d \pi"] & TM \ar[r] & 0
				\end{tikzcd}
			\end{displaymath}
			where $P \times_G \mathfrak{g}$ is the \emph{adjoint bundle of $P$}, \ie the vector bundle $(P\times \mathfrak{g})/\sim$ where $(p, ad_{g} \xi)\sim (R_g(p),\xi)$ (see \cite[\S 17.6.]{Kolar1993}).
		\end{lemma}

\subsubsection{Prequantization}\label{sec:Prequantum}
In this subsection we concisely review some basic notions related to geometric quantization, tailored to our needs (see \cite{Kostant70,Souriau66} for the original articles or \cite{Woodhouse97,Bry,Carosso2018} for a more recent review), in a finite dimensional environment.

Let $(M,\omega)$ be a connected symplectic manifold.
\begin{definition}[Prequantum bundle]\label{Def:PrequantumBundle}
	We call a \emph{prequantum bundle} of the symplectic manifold $(M,\omega)$ the pair $(P,\omega)$ consisting of a $S^1$-principal bundle 
	\begin{displaymath}
		S^1 \hookrightarrow P \twoheadrightarrow M
		~,
	\end{displaymath}
	together with a connection $\theta\in \Omega^1(P,\mathfrak{g})\cong\Omega^1(P)$ such that 
	\begin{displaymath}
		\pi^\ast \omega = \d \theta
	\end{displaymath}
	where $\pi:P\to M$ denotes the fibre projection encoding $P$.
	\\
	A prequantum bundle $(P,\theta)$ is also called a \emph{prequantization of $(M,\omega)$}. 
	When a given symplectic manifold admits a prequantum bundle it is said to be \emph{prequantizable}.
\end{definition}
Not all symplectic manifolds admit a prequantum bundle. 
The following celebrated theorem provides a cohomological condition to the existence of a prequantization.

\begin{theorem}[Weyl-Kostant integrality condition \emph{(see \cite{Kostant70} or \cite[Thm. 8.3.1]{Woodhouse97})}]\label{thm:integralitycondition}
	Consider a symplectic structure $\omega$ on the connected manifold $M$.
	The symplectic manifold $(M,\omega)$ is "prequantizable" (in the sense of definition \ref{Def:PrequantumBundle}) if and only if ${\frac{1}{2\pi}}[\omega]$ is an integral class, \ie it lies in the image of the mapping
	\begin{displaymath}
		\begin{tikzcd}
			H^2_{\text{sing}}(M,\mathbb{Z}) \ar[r] & H^2_{\text{sing}}(M,\R) \ar[r,equal,"\sim"] & H^2_{dR}(M)
		~.
		\end{tikzcd}
	\end{displaymath}
\end{theorem}
	\begin{remark}[Connections on circle bundles]
		Specializing definition \ref{def:connections} to the case of circle bundles, like in definition \ref{Def:PrequantumBundle}, has the following implications:
		\begin{itemize}
			\item since $\mathfrak{g}\cong \RR$, the connection is an honest $1$-form $\theta \in \Omega^1(P)$;
			\item since $\mathfrak{g}$ is generated by $1\in \RR$, one has that $\theta(\vAct_1)=1\in \RR$ where $\vAct_1$ denotes the fundamental vector field of the generator $1$;
			\item recalling that $S^1\cong U(1)$ and $\mathfrak{g}^\ast\cong \RR$, it results that $Ad_g^\ast = \id$ for any $g\in S^1$.
			Therefore $R_g^\ast \theta = \theta$, \ie $\theta$ is $S^1$-invariant.
		\end{itemize}
	\end{remark}
\begin{notation}
	We denote by $E\in \vX(P)$ the fundamental vector field, pertaining to the action of the $1$-dimensional group $S^1$ on $P$, corresponding to generator $1$.
\end{notation}

Once one fixes a "preferred" differential form, it is natural to select the class of smooth maps that preserve such a structure.
\begin{definition}[Infinitesimal quantomorphisms]
	Consider a prequantum bundle $(P,\theta)$ of $(M,\omega)$.
	 We call an \emph{infinitesimal quantomorphism} any vector field on $P$ preserving the connection $\theta$.
	We denote by $$Q(P,\theta):=\{Y\in \vX(P)~\vert~\cL_Y\theta=0\}$$ 
	the Lie subalgebra of $\mathfrak{X}(P)$ consisting of infinitesimal quantomorphisms.
\end{definition}
\begin{lemma}[{\cite[\S 2.2]{VaughanJ}}]\label{lem:preqInclusionSinvariant}
	Consider a prequantum bundle $(P,\theta)$. 
	The infinitesimal quantomorphisms are automatically $S^1$-invariant, \ie
	\begin{displaymath}
		Q(P,\theta) \subset \mathfrak{X}(P)^{S^1}~.
	\end{displaymath}
\end{lemma}
\begin{proof}
	Consider $E\in\mathfrak{X}(P)$, the generator of the action of $S^1$ on $P$.
	Notice that $E$ is determined by $\theta$ being the unique vector fields such that
$\iota_E\theta=1$ and $\iota_E \d\theta=0$. 
	Hence if $X\in\mathfrak{X}(P)$ preserves $\theta$, then it will preserve $E$ too.
	More precisely, let be $X\in\mathfrak{X}(P)$ such that $\mathcal{L}_X \theta =0$.
	One has that ${[X,E]}$ is an horizontal vector field since
	\begin{displaymath}
	 	\iota_{[X,E]}\theta = \mathcal{L}_X \iota_E \theta - \iota_E \mathcal{L}_X \theta = 0
	 	~,
	\end{displaymath}
	hence $[X,E]$ is projectable and completely determined by its projection onto $M$.
	Similarly, one has that $\iota_{[X,E]} \d \theta  =0$.
	Noticing that
	\begin{displaymath}
		\iota_{[X,E]} \d \theta = \iota_{[X,E]} \pi^\ast \omega = \iota_{\pi_\ast [X,E]} \omega
		~,
	\end{displaymath}
	the non-degeneracy of $\omega$ implies $[X,E]= \mathcal{L}_E X =0$ for any $X\in Q(P,\theta)$.
\end{proof}
Assume that  ${\frac{1}{2\pi}}[\omega]$ is an integral class and fix a prequantization circle bundle $P\to M$. 
As above, we denote by $E\in \vX(P)$ the unique infinitesimal generator pertaining to the action of the $1$-dimensional group $S^1$ and by $H_{\theta}:=\ker(\theta)$ the invariant Ehresmann connection on $P$ corresponding by $\theta$.
It is known that a prequantization provides a Lie algebra isomorphism between the observables Poisson algebras and the infinitesimal quantomorphisms.

\begin{lemma}[Kostant {\cite{Kostant70}} \emph{(see also \cite[Thm. 2.8]{VaughanJ})}]\label{lem:preqMap}
	Consider a prequantizable presymplectic manifold $(M,\omega)$, let be $(P,\theta)$ a prequantum bundle.
	One has a Lie algebra isomorphism
	\begin{equation}\label{eq:preq}
		\isomorphism{\Preq_{\theta}}
		{C^{\infty}(M)_{\omega}}
		{Q(P,\theta)}
		{f}
		{X_f^{H_\theta}+(\pi^*f)\cdot E},
	\end{equation}
where $X^{H_\theta}$ denotes the horizontal lift  of a vector field $X$ on $M$ using the Ehresmann connection $H_\theta$ (see reminder \ref{rem:ehresmann} and \ref{rem:horlifts}) and $\pi^\ast f \in C^{\infty}(P)$ is the pullback of $f$ along $\pi:P\to M$. 
\end{lemma}
The Lie algebra isomorphism $\Preq_\theta$ is the first ingredient to realize the map \eqref{eq:intropreqmap} anticipated in the introduction of this chapter.

\begin{remark}
 In the definition of $\Preq_\theta$ is implied the commutation of the following diagram in the category of Lie algebras:
 \begin{equation}\label{eq:preqiso}
 	\begin{tikzcd}[column sep = huge]
	 C^{\infty}(M)_{\omega} \ar[r,"\Preq_{\theta}"',"\sim"] \ar[d,"X_{\bullet}"'] &
	 Q(P,\theta)  \ar[dl,shift left=.5em,"\pi_*"]  
 		\\
 		\vX(M) & 
	\end{tikzcd}
	~.
\end{equation}
	Notice that the vertical map has a one-dimensional kernel given by the constant functions on $M$, \ie $\ker(X_\bullet)=\RR$. 
	Hence the same holds for the map $\pi_*$ in the diagram giving the pushforward of projectable vector fields via $\pi$, i.e
	\begin{displaymath}
		\Preq_\theta(\ker(X_\bullet)) = \ker(\pi_\ast\vert_{Q(P,\theta)})~,
	\end{displaymath} 
	where the right-hand side consists of vertical vector fields preserving $\theta$, \ie elements of $\ker(\pi_\ast\vert_{Q(P,\theta)})$ are constant multiples of $E$.
\end{remark} 

\begin{remark}[About the \emph{Quantum} name]
	We briefly explain why the term "quantum" appears in this context.
	\\
	Associating a prequantum bundle to a prequantizable symplectic manifold is the first step of a 3-step procedure called \emph{"geometric quantization (scheme)"} essentially due to Kostant, Kirillov and Souriau (KKS).
	\\
	Roughly, a "quantization scheme" is procedure to associate to any symplectic manifold (prequantizable in some sense), taken together with the corresponding Poisson algebra of observables $(C^{\infty}(M),\{\cdot,\cdot\})$, a Hilbert $\mathbb{C}$-vector space $\mathscr{H}$, taken together with the algebra of self-adjoint operators on $\mathscr{H}$.
	\\
	If one understands states of a classical mechanical system as points on a (finite dimensional) manifold $M$, the corresponding "quantum" states will be  vectors of an (infinite dimensional) complex vector spaces $\mathscr{H}$ with unitary norm.
	The keypoint is the linearity of $\mathscr{H}$. This property provides a framework making possible to encompass the phenomenon of "superposition of states".
	In particular, given $\varphi\in \mathscr{H}$, the vectors $e^{i\lambda}\varphi$ represents the same "physical state" for any $\lambda\in \RR$.
	\\
	Once this point is clear, the reason why we have named "prequantization" the act of associating a $S^1$-bundle $P$ over $M$ should begin to emerge.
	Intuitively, a $S^1$ bundle over $M$ is the attachment of a $1$-dimensional circle to any classical state $p\in M$. 
	On the other hand, it is well-known that $S^1$ is diffeomorphic to the unitary circle in the complex plane $\mathbb{C}$, \ie the Lie group $U(1)$. 
	Hence points in $P$ can be seen as classical states taken together with a certain \emph{phase factor} $e^{i \lambda}\in\mathbb{C}$.
	\\
	Nevertheless, being $P$ a generic smooth manifold, is not yet the sought linear space. 
	According to the (KKS) procedure, the \emph{prequantum Hilbert space} is represented by a certain subclass of complex-valued smooth function over $P$.
	Furthermore, being possible to regard vector fields on $P$ as derivations $C^\infty(P)\to C^\infty(P)$, the images of the morphism given in equation \eqref{eq:preq} can be seen as prequantum versions of the classical observables of $C^\infty(M)$, hence the name "prequantization map"\footnote{We mention that the prequantization procedure could be equivalently expressed in terms of $1$-dimensional Hermitian bundles and integrable sections, see for instance \cite[\S 1.1]{Weinstein2005a}.}.
\note{	Furthermore, the morphism given in equation \eqref{eq:preq} can be suitably extended to give an operator $\Gamma(P)\to \Gamma(P)$ rather then a vector field $\mathfrak{X}(P)$ ( which in turn can be seen as a derivation $C^\infty(P)\to C^\infty(P)$). The former are prequantum version of the classical observables of $C^\infty(M)$, hence the name "prequantum map".
}
 	We do not insist here on further details and we refer the interested reader to the fundamental manuals of geometric quantization \cite{Bry} and \cite{Woodhouse97}.	
	We only stress that what we loosely described here are quantization schemes for ordinary, \ie point-like, mechanical systems. 
	The mathematical foundation of quantization procedures for $\infty$-dimensional mechanical systems, \ie field theories, is still largely incomplete.   
\end{remark}

\subsection{Embedding of the observables in the Atiyah algebroid}\label{subsec:first}\label{subsec:At}
%
In this subsection we derive the map \eqref{eq:intropreqmap}. {The material reviewed here can be found also in \cite[\S 2]{Rogers2013}.}

Consider again a principal circle bundle $\pi\colon P\to M$.
According to definition \ref{def:ati}, the corresponding Atiyah Lie algebroid $A_P$ is the transitive Lie algebroid over $M$ with space of sections given by $\vX(P)^{S^1}$, the invariant vector fields on $P$, and anchor given by $\pi_*$.

Lemma \ref{lem:atiexactseq} says that $A_P$ fits in a short exact sequence of Lie algebroids
	\begin{displaymath}
		\begin{tikzcd}
			0 \ar[r] & \RR \ar[r] & A_P \ar[r] & TM \ar[r]& 0
		\end{tikzcd}
	\end{displaymath}
	where 
	{$\RR$ denotes the trivial rank-1 bundle\footnote{This is usually denoted as $\RR_M$, we are employing a short-hand notation here.} over $M$} (a bundle of Abelian Lie algebras, {with the constant section $1$ mapping to $E\in \Gamma(A_P)$}) and {the second map  is the anchor}.
A principal connection $\theta$ on $P$ provides a linear splitting of the above short exact sequence.
	\begin{lemma}[$TM \oplus  \RR$ is isomorphic to $A_P$]\label{lem:sigmatheta}
		Let be $\pi: P \to M$ a $S^1$-principal bundle. Consider a connection $1$-form $\theta$.
		There is an isomorphism of vector bundles
	\begin{displaymath}
		\isomorphism{\sigma_\theta}
		{TM \oplus  \RR}
		{A_P}
		{\pair{v}{c}}
		{v^{H_\theta}+ c\cdot E}
		~,
	\end{displaymath}
		where the superscript $H_\theta$ denotes the horizontal lift on vectors and $E$ is the fundamental vector field of the generator $1$ in the Lie algebra of $S^1$.	
	\end{lemma}
	\begin{remark}[Recovering the standard $\omega$-twisted Lie algebroid]
		Considering the inverse of the isomorphism $\sigma_\theta$ introduced in lemma \ref{lem:sigmatheta}, one can pull back the Lie algebroid structure from $A_P$ to $TM\oplus \RR$.
	As a result, we obtain the Lie algebroid  
$(TM\oplus \RR)_{{\omega}}$, with anchor map given by the first projection onto $TM$ and Lie bracket  on sections given by:
	\begin{displaymath}
		\begin{aligned}
			\left[\pair{X}{f},\pair{Y}{g} \right]_{\omega}
			:=&~
			\sigma_\theta^{-1} \left( \left[\sigma_\theta\pair{X}{f},\sigma_\theta\pair{Y}{g}\right] \right) =
			\\
			=&~\pair{[X,Y]}{X(g)-Y(f) + {\iota_X\iota_Y\omega} }
			~.
		\end{aligned}
	\end{displaymath}
	To prove the last equality, observe that
	\begin{displaymath}
		\mathclap{
		\begin{aligned}
				&[X^{H_\theta} +(\pi^\ast f) E, Y^{H_\theta} + (\pi^\ast g) E ]
				=
				\\
				&=~
				[X^{H_\theta},Y^{H_\theta}] + [	(\pi^\ast f) E, Y^{H_\theta}] 
				+ [X^{H_\theta}, (\pi^\ast g) E]
				+ \cancel{[(\pi^\ast f) E,(\pi^\ast g) E] 		}
				=
				\\
				&=~ [X^{H_\theta},Y^{H_\theta}] + 
				\cancel{ (\pi^\ast f) [E,Y^{H_\theta}]} +
				\cancel{(\pi^\ast g) [X^{H_\theta},E]}
				+
				\left(	
				-\cL_{Y^{H_\theta}}(\pi^\ast f)	
				+\cL_{X^{H_\theta}}(\pi^\ast g)
				\right) E
				=
				\\
				&=~ [X^{H_\theta},Y^{H_\theta}] + 
				\pi^\ast\left( X(g)-Y(f)	\right) E
		\end{aligned}
		}
	\end{displaymath}	
	where the first cancellation occurs because $E$ is vertical and $\pi^\ast f$ and $\pi^\ast g$ are constant along the fibres, and the second one follows from $\pi_\ast ([X^{H_\theta},E])=\theta([X^{H_\theta},E])=0$.
	The claim is proved by noticing that $$[X^{H_\theta},Y^{H_\theta}]= [X,Y]^{H_\theta} + \pi^\ast(\iota_X\iota_Y\omega) E$$ since
	\begin{displaymath}
		\mathclap{
		\begin{aligned}
		\theta([X^{H_\theta},Y^{H_\theta}] 
		=&~
		\cL_{X^{H_\theta}}\cancel{\iota_{Y^{H_\theta}}\theta}
		-\iota_{Y^{H_\theta}}\iota_{X^{H_\theta}} \d \theta 
		+ \iota_{Y^{H_\theta}}\d \cancel{\iota_{X^{H_\theta}}\theta} =
		\\
		=&~
		-\iota_{Y^{H_\theta}}\iota_{X^{H_\theta}} \pi^\ast \omega =
		\\
		=&~
		-\pi^\ast (-\iota_{Y}\iota_{X} \omega) 		
		\end{aligned}
		}
	\end{displaymath}
	and
	\begin{displaymath}
		\pi_\ast [X^{H_\theta},Y^{H_\theta}] = [X,Y]
		~.
	\end{displaymath}
	\end{remark}	
	The upshot is that $\sigma_\theta$ is a Lie algebroid isomorphism
$(TM\oplus \RR)_{{\omega}} \cong A_P$ between the standard $\omega$-twisted Lie algebroid (see example \ref{ex:TwistedLieAlgbroid}) and the Atiyah algebroid of the $S^1$-principal bundle. 	

Finally, noticing that $Q(P,\theta)\subset \vX(P)^{S^1} =\Gamma(A_P)$, we conclude that the sought embedding (see equation \ref{eq:intropreqmap}) can be obtained by composing\footnote{We will tend, with a slight abuse of notation, to regard the map $\Preq_{\theta}$ introduced in equation \eqref{eq:preq} as a morphism $C^{\infty}(M)_{\omega}\to \vX(P)^{S^1}$.
} 
 $\Preq_{\theta}$ and   $\sigma_\theta^{-1}$.
\begin{proposition}[The Poisson algebra embeds into the sections of $(TM\oplus\RR)_\omega$]
	Let be $(M,\omega)$ a prequantizable symplectic manifold.
	Consider the $\omega$-twisted standard Lie algebroid $(TM\oplus\RR)_\omega$.
	The Lie algebra morphism $\Psi$ obtained by the composition of the maps introduced in lemmas \ref{lem:preqMap}, \ref{lem:preqInclusionSinvariant} and \ref{lem:sigmatheta},
	\ie the map obtained from the following diagram
\begin{displaymath}
	\begin{tikzcd}[row sep = small, column sep = small]
		C^\infty(M)_\omega \ar[dr,equal,sloped,"\sim","\Preq_{\theta}"'] \ar[rrr,dashed,"\Psi"] 
		& &[2em] & \Gamma(TM\oplus \RR)_{{\omega}} 
		\\
		& Q(P,\theta) \ar[r,hook] 
		& \Gamma(A_P) \ar[ur,equal,sloped,"\sim","\sigma_\theta^{-1}"']
	\end{tikzcd}
	~,
\end{displaymath}
	is the Lie algebra embedding \eqref{eq:intropreqmap} appearing in the introduction.
	\\
	Namely one has:
		\begin{equation}\label{eq:chris}
			\morphism{\Psi}
			{C^{\infty}(M)_{\omega}}
			{\Gamma(TM\oplus \RR)_{{\omega}}}
			{f}
			{\pair{X_f}{f}}
			~,
		\end{equation}
		where $X_f$ denotes the Hamiltonian vector field pertaining to $f$.
\end{proposition}
	For the rest of this section, we will simply denote $\Psi$ as
	$$\sigma_{\theta}^{-1}\circ \Preq_{\theta}\colon  C^{\infty}(M)_{\omega} \to
\Gamma(TM\oplus \RR)_{\omega}
~.$$
 \begin{remark}[Independence from the choice of prequantization]\label{rem:IndiPreQuantum}
  We point out that the above expression \eqref{eq:chris} is a Lie algebra embedding  
even when $\omega$ does not satisfy the integrality condition.
	Furthermore, it does not depend on the connection $\theta$ implied by the prequantization procedure.
\end{remark}
%

\subsection{Commutativity after twisting}\label{subsec:second}
{In this subsection we show, by geometric arguments, the commutativity of the diagram \eqref{intro:diag:main} from the introduction.}

Assume we have an action of a Lie group $G$ on $M$, and denote by $v:\g\to \vX(M)$ the corresponding infinitesimal action (a Lie algebra morphism).
Assume the existence of an equivariant moment map  $$J\colon M\to \g^*.$$ This means that $J$ satisfies $\iota_{v_x}\omega=-d(J^*(x))$ (\ie $v_x=X_{J^*(x)}$) and that $J^*\colon \g\to  C^{\infty}(M)_{\omega}$ 
is a Lie algebra morphism\footnote{This is equivalent to infinitesimal equivariance, \ie 
$\cL_{v_y}J^*(x)=J^*([y,x])$.} 
(see reminder \ref{Rem:SymplecticMomaps} in chapter \ref{Chap:MultiSymplecticGeometry}).

Therefore, diagram \eqref{eq:preqiso} is extended to
 \begin{equation}\label{diag:gQ}
 	\begin{tikzcd}[column sep = huge, row sep = large]
 		& C^{\infty}(M)_{\omega} \ar[r,"{\sim}","\Preq_\theta"'] \ar[d,"X_{\bullet}"]  
 		& Q(P,\theta)  \ar[dl,"\pi_*"]  
 		\\
 		\g \ar[r,"{v_{\bullet}}"'] \ar[ru,"{J^*}"]
 		& \vX(M) 
 		&	
 	\end{tikzcd}
\end{equation} 
In particular, we obtain a Lie algebra morphism 
\begin{equation}\label{eq:L0map}
	\morphism{L_0:=\Preq_{\theta}\circ J^*}
	{\g}
	{Q(P,\theta)}
	{x}
	{v_x^{H_\theta}+\left((\pi^*J^*)(x)\right)\cdot E}~,
\end{equation}
lifting the infinitesimal action in the sense of the commutativity of the following diagram in the category of Lie algebras  
\begin{equation*}
	\begin{tikzcd}[column sep = large]
		& Q(P,\theta)  \ar[d,"{\pi_*}"] 
		\\
		\g \ar[r,"v_{\bullet}"'] \ar[ru,"L_0"]
		& \vX(M)
	\end{tikzcd}~.
\end{equation*}

 \subsubsection{Twisting by an invariant one-form}
{Notice that the difference between any two connection 1-forms on the circle bundle $\pi\colon P\to M$ is basic, \ie is the pullback of a 1-form on $M$.}
\\
Now we take  $\alpha\in \Omega^1(M)^G$ and use it to twist some of the above data, keeping the $G$-action fixed:
$\omega+d\alpha$ is an invariant symplectic form on $M$ (assuming it is non-degenerate), with moment map $J_{\alpha}$ determined by\footnote{Indeed it can be checked easily that $\iota_{v_x}(\omega+d\alpha)=-d(J^*(x)+\iota_{v_x}\alpha)$ using the $G$-invariance of $\alpha$ (expressed as $\cL_{v_x}\alpha=0$
for all $x\in \g$).} 
\begin{equation}\label{eq:jalpha}
	\morphism{J_{\alpha}^*}
	{\g}
	{C^{\infty}(M)_{\omega+d\alpha}}
	{x}
	{\mapsto J^*(x)+\iota_{v_x}\alpha}
	~.   
\end{equation}
 Furthermore, a prequantization of the symplectic manifold $(M,\omega+d\alpha)$ is given by the same circle bundle $P$ but with connection $\theta+\pi^*\alpha$.

We can repeat the  procedure outlined above (see in particular equation \eqref{eq:L0map}), 
obtaining a  Lie algebra morphism $L_\alpha\colon \g \to Q(P,\theta+\pi^*\alpha) $ lifting the infinitesimal action.
Since $\alpha$ is $G$-invariant,  any lift to $P$ of a generator $v_x$ preserves $\pi^*\alpha$, hence we can view both $L_0$ and $L_{\alpha}$ as maps $$ \g \to Q(P,\theta)\cap Q(P,\theta+\pi^*\alpha)$$ 
which are \emph{Lie algebra morphisms
lifting the infinitesimal action}.
There   are ``few'' such Lie algebra morphisms. (They are in bijection with moment maps for $(M,\omega)$, by diagram \eqref{diag:gQ}; if $H^1(\g)=0$ then
the moment map is unique
\cite[Theorem 26.5]{CannasdaSilva2001}.) Hence the following is not a surprise.
\begin{proposition}\label{prop:L}
The   Lie algebra morphisms $L_0$ and $L_\alpha$ coincide.
\end{proposition}
\begin{proof}
Fix $x\in \g$  and write $f_x:=J^*(x)$. We have to show that $L_0(x)=L_{\alpha}(x)$, \ie
	$$v_x^{H_\theta}+\left(\pi^*f_x\right)\cdot E=v_x^{H_{\theta+\pi^*\alpha}}+\pi^*(f_x+\iota_{v_x}\alpha)\cdot E~.$$
We do so decomposing $TP$ as $H_\theta\oplus \RR E$. Since both the left-hand side and the right-hand side $\pi$-project to the same vector field (namely $v_x$), we have to check  that   we obtain the same function applying $\theta$ to both vector fields. This is indeed the case, since
applying $\theta$ to the vector field on the right  we obtain
$$\pi^*(f_x+\iota_{v_x}\alpha)-(\pi^*\alpha)(v_x^{H_{\theta+\pi^*\alpha}})=\pi^*f_x.$$
\end{proof}

We can also repeat the construction of \S \ref{subsec:At}
using the connection $\theta+\pi^*\alpha$, {yielding a Lie algebroid isomorphism $$\sigma_{\theta+\pi^*\alpha}\colon (TM\oplus \RR)_{\omega+d\alpha} \cong A_P.$$}
The composition $(\sigma_{\theta+\pi^*\alpha})^{-1}\circ \sigma_{\theta}$ reads
\begin{equation}\label{eq:tau}
	\isomorphism{\tau_{\alpha}}
	{(TM\oplus \RR)_{\omega}}
	{(TM\oplus \RR)_{\omega+d\alpha}}
	{\pair{v}{c}}
	{\pair{v}{c+ \iota_{v}\alpha}}
\end{equation}
and is often referred to as  \emph{gauge transformation}.

\subsubsection{Commutativity}
We end up with the following commutative diagram
\begin{displaymath}
	\begin{tikzcd}
		& C^{\infty}(M)_{\omega} \ar[dr,"{\Preq_{\theta}}"]
		&&  \Gamma(TM\oplus \RR)_{\omega} \ar[dd,"{\tau_{\alpha}}"']
		\\
		\g \ar[rd,"{J_{\alpha}^*}"'] \ar[ru,"{J^*}"]
		&& \vX(P)^{S^1} \ar[ur,"{\sigma_{\theta}^{-1}}"] \ar[dr,"{\sigma^{-1}_{\theta+\pi^*\alpha}}"'] 
		&
		\\
		& C^{\infty}(M)_{\omega+d\alpha} \ar[ur,"{\Preq_{\theta+\pi^*\alpha}}"']   
		&& \Gamma(TM\oplus \RR)_{\omega+\d\alpha}	
	\end{tikzcd}
\end{displaymath}
where the left square commutes by  proposition \ref{prop:L} and the right one by the very definition of $\tau_\alpha$ (see equation \eqref{eq:tau}).

As we emphasized in remark \ref{rem:IndiPreQuantum}, the composition $\sigma_{\theta}^{-1}\circ \Preq_{\theta}\colon  C^{\infty}(M)_{\omega} \to
\Gamma(TM\oplus \RR)_{\omega}$ does not depend on $\theta$. 
Hence, after removing $\vX(P)^{S^1}$ from the above diagram, we obtain a commutative diagram that makes no reference to the prequantization bundle $P$:
	\begin{equation}\label{diag:main}
		\begin{tikzcd}[column sep=huge]
			&
			  C^{\infty}(M)_{\omega}   \ar[r,""] 
			&
			 \Gamma(TM\oplus \RR)_{\omega}  \ar[dd,"\tau_\alpha"]
			\\[-1em]
			\mathfrak{g}\ar[ru,"J^*"] \ar[dr,"J^*_{\alpha}"']
			\\[-1em]
			&
			 C^{\infty}(M)_{\omega+d\alpha} \ar[r,""] 
			&
			\Gamma(TM\oplus \RR)_{\omega+d\alpha}
		\end{tikzcd}
	\end{equation}
 
 \begin{remark} For a given $\alpha$, in general, there is no linear map  $C^{\infty}(M)_{\omega}\to C^{\infty}(M)_{\omega+d\alpha}$ making the left part of diagram \eqref{diag:main} commute. 
 Indeed such a map exists if, and only if, for all $f\in C^{\infty}(M)$ we have
 $$X_f^{\omega}=X_{f+\iota_{X^{\omega}_f}\alpha}^{\omega+d\alpha}$$
 where $X_g^{\nu}$ denotes the Hamiltonian vector field pertaining to the function $g$ with respect to the symplectic form $\nu$.
 The latter is equivalent to say that $\cL_{X_f^{\omega}}\alpha=0$. 
 This explains why it is necessary to consider moment maps for a $\alpha$-preserving action.
\end{remark}
 
 \begin{remark}\label{rem:symcomm}
Diagram \eqref{diag:main} commutes for any symplectic form $\omega$, even for those that do not satisfy the integrality condition and therefore do not admit a prequantization bundle.
{This is immediate using the explicit expressions for the maps involved in eq. \eqref{eq:chris}, 
\eqref{eq:jalpha} and \eqref{eq:tau}.}
 The discussion of this subsection -- in particular proposition \ref{prop:L} -- provides a geometric argument for the commutativity of diagram \eqref{diag:main} in the integral case.
\end{remark}

\subsection{Motivation for the higher case}\label{subsec:higher}
In the rest of this chapter we will consider a multisymplectic form $\omega$, and 
 show that the higher analogue of diagram \eqref{diag:main}  commutes too. 
 This provides some evidence that, in the integral case, one can expect a global geometric picture (higher prequantization) as the one we outlined for the symplectic case in the previous section. 

	We first address \S\ref{subsec:first}. 
	In the integral case the analogue of the prequantization map  $\Preq_{\theta}$   of eq. \eqref{eq:preq} was already established  for all $n$ in Fiorenza-Rogers-Schreiber \cite[Thm. 4.6]{Fiorenza2014a}; there however the higher prequantization bundle is described by means of an open cover on the manifold $M$. 
	For $n=2$,  the analogue of the prequantization map was established on a higher prequantization bundle that admits a global description (without choosing a cover) in Krepski-Vaughan \cite[\S 5.1]{krepski2020multiplicative} -- but not as an $L_{\infty}$-algebra morphism -- and later in Sevestre-Wurzbacher
\cite[Thm. 3.5]{SevestreWurzbacherPreq}. 
	See \cite[Remark 5.2]{krepski2020multiplicative} 
	and \cite[Rem. 3.6]{SevestreWurzbacherPreq}  
for a comparison. 
For the fact that ``higher Atiyah algebroids'' (also known as Vinogradov algebroids) can be obtained from $S^1$-gerbes and higher analogues, see \cite[\S 2.3]{Gualtieri2004}.


\section{Algebraic Properties of Rogers' $L_\infty$-algebra}\label{Sec:RogersL1prop}
In this section, we review again the construction of the $L_\infty$-algebra of observables introduced by Rogers, focusing on its presentation in terms of symmetric multibrackets ($L_\infty[1]$-algebra, \cf  definition \ref{Def:LInfinityShifted}).
Namely, we establish certain relations between the multibrackets of different degrees. We do so by means of concise computations using the \RN operation $\cs$   introduced in eq. \eqref{eq:compsymm}. 

In sections \ref{Sec:Vinoids} and \ref{Section:ExtendedRogersEmbedding}, these relations will be necessary to relate $L_\infty(M,\omega)$ with the $L_\infty$-algebra associated to the corresponding "higher Courant algebroid" (more precisely they will be used to make explicit the expression \eqref{eq:pi'} appearing in proposition \ref{prop:main}).

\subsection{Rogers'$L_{\infty}[1]$-algebra}\label{sec:RogersL1}
Let $(M,\omega)$ be an $n$-plectic manifold, denote by $L_{\infty}(M,\omega):=(L,\{l_{k} \})$ the associated $L_\infty$-algebra prescribed by the Rogers' construction\footnote{Notice the small change of notation: in chapter \ref{Chap:Linfinity} we denoted the $k$-ary multibrackets as $[\cdots]_k$.}.
In order to understand the relationship of $L_\infty(M,\omega)$ with the Vinogradov's $ L_\infty $ structure (see section \ref{subsec:Vinogradov}), it will be more convenient to consider the graded vector space $\cA$ given by the following components
\begin{equation}\label{eq:Aspace}
	\cA^i:=
	\begin{cases}
		 \left.\left\lbrace
		\pair{X}{\alpha}\in \mathfrak{X}(M)\oplus \Omega^{n-1}(M)
		~ \right\vert ~
		\iota_X \omega = -\d \alpha\right\rbrace
		 & ~\text{if } i=0
		 \\
		 ~\Omega^{n+i-1}(M) &  ~\text{if } 1-n \leq i\leq -1
		 \\
		 ~0 & ~\text{otherwise}
		 ~.
	\end{cases}
\end{equation}
\begin{remark}
	Observe that the graded vector space $\cA$ coincides with the graded vector space
\begin{displaymath}
	Ham_\infty(M,\omega):= Ham^{n-1}(M,\omega)\oplus \trunc_{0}(\Omega(M)[n-1])
\end{displaymath} mani
	introduced in remark \ref{Rem:DegenerateCase}.
	Accordingly, $\cA^0 = Ham^{n-1}(M,\omega)$ consists of all the Hamiltonian pairs pertaining to $\omega$, \ie Hamiltonian forms together with their corresponding Hamiltonian vector field.
\end{remark}

When $\omega$ is $n$-plectic, $\cA\cong L$ as graded vector spaces. Let us denote $\pi_k$ the pullback of $\ell_k$ from $L$ to $\cA$ along the isomorphism acting as the projection $\cA^0\twoheadrightarrow \Omega^{n-1}_{\ham}(M,\omega)$ in degree $0$ and as the identity in other degrees.
Being $\omega$ non-degenerate, the latter projection is in particular bijective. We will denote as $\vartriangle$ its inverse, given, in degree $0$, by mapping Hamiltonian forms into the corresponding Hamiltonian pairs.
Namely, the map acts like the identity in  degrees lesser than $0$ and by
\begin{equation}\label{eq:trianglemap}
	\morphism{(\blank)^\vartriangle}
	{L^0=\Omega^{n-1}_{\ham}(M,\omega)}
	{\cA^0 \subset \mathfrak{X}(M)\oplus\Omega^{n-1}(M)}
	{\alpha}
	{\pair{\vHam_\alpha}{\alpha}}
	~,
\end{equation}
in degree $0$. Accordingly we will also denote $(\Omega^{n-1}_{\ham}(M,\omega))^{\Delta}=\cA^0$.

For the sake of clarity, we reiterate the explicit expression for the higher observables multibrackets in this slightly different setting.
Denoting by $\varv= f\oplus e$ a generic element in $\cA$, where $f\in \bigoplus_{k=0}^{n-2}\Omega^k(M)$ and $e=\pair{X}{\alpha}\in \cA^0$, the unary multibracket reads as the following linear map
\footnote{To be completely consistent with the framework introduced in section \ref{Sec:ConventionsMultiLinearAlgebras}, we should replace $\cA$ with $\mathcal{A^\oplus}$. Here we are choosing to lighten that notation. Everything should be clear from the context. 
}
\begin{displaymath}
	\morphism{\pi_1}
	{\cA}
	{\cA}
	{f\oplus\pair{X}{\alpha}}
	{\d f}		
	 ~;					
\end{displaymath}
and all the other non-trivial $k$-multibrackets $(2\leq k \leq n$+$1)$ result in :
\begin{displaymath}
	\morphism{\pi_k}
	{\cA^{\otimes k}}
	{\cA}
	{e_1\otimes\dots \otimes e_k}
	{\varsigma(k) \iota(X_1\wedge\dots\wedge X_k) \omega}
	~.
\end{displaymath}

Recall that, given a graded vector space $\cA$, an $L_{\infty}[1]$-algebra structure on $\cA[1]$ is equivalent to a $L_{\infty}$-algebra structure on $\cA$, via the d\'ecalage isomorphism (see equation \eqref{eq:deca}), reading in the present case as follows
\begin{equation}\label{deca}
	\isomorphism{\dec}
	{(\wedge^n \cA)[n]}
	{\odot^n(\cA[1])}
	{(u_1\wedge \cdots \wedge u_n)_{[n]}}
	{u_{1 [1]}\cdots u_{n [1]}\cdot (-1)^{(n-1)|u_{1}|+\dots+2|u_{n-2}|+|u_{n-1}|}}
\end{equation}
 where $u_1,\dots,u_n\in \cA$ are homogeneous vectors nd $|u_i|$ denote the degrees of $u_i\in \cA$. 
 \\
 We denote  by $(\cA[1],\{\bpi_k\})$ the $L_{\infty}[1]$-algebra corresponding to Rogers' Lie infinity algebra on $\cA$.

\subsection{{Properties} of Rogers' $L_{\infty}[1]$-algebra}\label{sec:L1RogersProp}
	We consider now the $L_\infty[1]$-algebra $(\cA[1],\{\bpi_k\})$.
	The crucial property is that, whatever the degree of the multisymplectic form $\omega$ considered on $M$, all Rogers' multibrackets can be expressed as combination of the three lowest arity multibrackets together with the following auxiliary operators:
	\begin{definition}[Symmetric and skew-symmetric Pairings $\pairing_\pm$]\label{def:pairing-}
		Let be $M$ a smooth manifold and $n\in \ZZ$ a fixed integer.
		Consider the graded vector space 
		\begin{equation}\label{eq:VtildeSpace}
			\widetilde{\cV}:= \mathfrak{X}(M)\oplus\left( \Omega(M)[n]\right)~.
		\end{equation}
		We call \emph{symmetric (resp. skew-symmetric) pairing} the, degree $-1$, binary brackets
		$\pairing_+\in M^{\sym}(\widetilde{\cV})$ (resp. $\pairing_-\in M^{\skew}(\widetilde{\cV})$)  defined as
		\begin{displaymath}
			\morphism{\pairing_\pm}
			{\widetilde{\cV}\otimes\widetilde{\cV}}
			{\widetilde{\cV}}
			{(x_1+f_1, x_2 +f_2)}
			{\dfrac{1}{2}(\iota_{x_1}f_2 \pm \iota_{x_2} f_1)}
		\end{displaymath}
		for any $x_i \in \mathfrak{X}(M)$ and $f_i\in \Omega(M)$.
	\end{definition}
	
	\begin{remark}\label{rem:pair-}
		In the following, we will be mostly interested in the restriction of $\pairing_-$ to the graded subspaces $\cA$ (or to $\cV$, that is a truncation of $\widetilde{\cV}$, see \S \ref{subsec:Vinogradov} below).
		This determines a graded skew-symmetric bilinear map $\cA\otimes \cA\to \cA$ of degree $-1$.		
		Namely, for any $f_i\in \bigoplus_{k=0}^{n-2}\Omega^k(M)$ and $e_i=\pair{X_i}{\alpha_i}\in \cA^0$, the pairing reads as follows
\begin{equation}\label{Eq:PairingExtensions}
	\left\langle f_1 \oplus \pair{X_1}{\alpha_1}, f_2 \oplus \pair{X_2}{\alpha_2} \right\rangle_\pm = 
	\frac{1}{2}{\left(
	\iota_{X_1}( \alpha_2 + f_2) -\iota_{X_2}( \alpha_1 + f_1)\right)}\oplus
	\pair{0}{0}
	~.
\end{equation} 
In turn, the latter defines by d\'ecalage a degree zero map  $S^{ 2}(\cA[1])\to \cA[1]$ which we  denote by {$\pairing$}.
Observe that, for any $\pair{X_i}{\alpha_i}_{[1]}\in (\cA[1])^{-1} = \cA^0$, one has
	\begin{displaymath}
		\begin{aligned}
		\left\langle \pair{X_1}{\alpha_1}_{[1]}, \pair{X_2}{\alpha_2}_{[1]}\right\rangle
		&=
		\left(\iota_{X_1}\alpha_2 + (-)^{|X_{2[1]}||\alpha_{1[1]}|}\iota_{X_2}\alpha_1\right) =
		\\
		&=
		\left(\iota_{X_1}\alpha_2 - \iota_{X_2}\alpha_1\right) =
		\\
		&= \left(\left\langle \pair{X_1}{\alpha_1}, \pair{X_2}{\alpha_2}\right\rangle_-\right)_{[1]}
		~.
		\end{aligned}
	\end{displaymath}
This map vanishes if both entries of $\cA[1]$ lie in degrees $\le -2$. By extending trivially we obtain a map $\pairing \colon S^{\ge 1}(\cA[1])\to \cA[1]$. 
\end{remark}

According to the next lemma, the whole Rogers' $L_\infty[1]$ structure on $\cA$ is completely generated by the following four multibrackets $\{\pairing, \bpi_1, \bpi_2, \bpi_3\}$:
\begin{lemma}[Higher Rogers multibrackets recursive formula]\label{lem:rogersRecurFormula}
 \begin{displaymath}
	[\pairing,\bpi_{k-1}] _{\cs}
	= 
	\frac{k}{2} \bpi_k \qquad \forall k \geq 4
 \end{displaymath}
 	where $\ca$ and $\cs$ denote respectively the \RN product of graded skew-symmetric and graded symmetric multibrackets introduced in section \ref{sec:RNProdMB}	(see appendix \ref{App:RNAlgebras} for all the details).
\end{lemma}
{We observe that $\bpi_3$ can be expressed in terms of $\bpi_2$ too, see the proof of Lemma \ref{Lemma:BoringAssociator}.}

\begin{proof}
	Inspecting elements $e_i=f_i + \pair{X_i}{\alpha_i}\in \cA$ one gets
	\allowdisplaybreaks
	\begin{align*}
			\pi_k(e_1,\dots,e_k) 
			&= \varsigma(k)\omega(X_1,\dots,X_k) 
			= 
			\\
			&=
			\varsigma(k)\sum_{\sigma\in \ush{k-1,1}} \frac{1}{k}\iota_{X_{\sigma_k}}\omega(X_{\sigma_1},\dots,X_{\sigma_{k-1}})
			=
			\\
			&=
			\left(\frac{\varsigma(k)\varsigma(k-1)}{k}\right)
			\sum_{\sigma\in \ush{k-1,1}}
			\iota_{X_{\sigma_k}}~\pi_{k-1}(e_{\sigma_1},\dots,e_{\sigma_{k-1}})
			=
			\\
			&= -\frac{2}{k} (-)^k 
			\sum_{\sigma\in \ush{k-1,1}}
			\Big\langle \pi_{k-1}(e_{\sigma_1},\dots,e_{\sigma_{k-1}}), e_{\sigma_k} \Big\rangle_-
			=
			\\
			&= -\frac{2}{k} (-)^k (-)^{|\pi_{k-1}|}
			\pairing_- \ca \pi_{k-1} (e_1,\dots e_k)
			=
			\\
			&=
			\frac{2}{k}			\pairing_- \ca \pi_{k-1} (e_1,\dots e_k)		
	\end{align*}
	\allowdisplaybreaks[0]
	The claim follows after d\'ecalage.
\end{proof}
The upshot is that the  the graded subalgebra of $M^{sym}(\cA[1])$ generated by the Rogers multibrackets is isomorphic to the subalgebra generated by $\bpi_1,\bpi_2,\bpi_3$ together with the pairing $\pairing$ \footnote{The same property could also be expressed in term of the $\ca$ product in the spaces $M(V)$ of multilinear maps without any symmetry. In that terms, one can see that are sufficient the three generators $\pi_1,\pi_2$ and $(\pairing_+ + \pairing_-)$ to generate the entire $L_\infty$-structure.}.
\\
All possible multibrackets generated by $\{\bpi_1,\bpi_2,\bpi_3,\pairing\}$
can be reconstructed from the $L_\infty$-algebra axioms and by working out the iterated {commutators} of powers 
$\pairing^{\cs l}$ with $\bpi_k$.

{
We now proceed to compute some of this commutators in some relevant cases.
}
\begin{remark}\label{rem:two}
 a) $\bpi_k$ with $k\geq 2$ is non-zero only when evaluated on degree $0$ elements, hence
 $$ [\pairing, \bpi_k ]_{\cs} = \pairing \cs \bpi_k.$$
 
 b) We carry out many of the proofs in terms of the multilinear maps $\pi$ on $\cA$, rather than using the corresponding graded-symmetric maps $\bpi$ on $\cA[1]$. This is possible thanks to the graded algebra isomorphism given by the D\'ecalage \eqref{eq:Dec}, and it is convenient because it allows us to apply the identities of Cartan calculus easily.
 
 c) {In this section, we will sometimes make use of the symmetric pairing 
$\pairing_+ : \cV \otimes \cV \to \cV$,
which is defined  analogously to $\pairing_-$ in eq. \eqref{Eq:PairingExtensions} but replacing the minus sign there with a plus sign.}
\end{remark}

\begin{proposition}[Commutators of arity $3$]\label{Prop:TernaryCommutator}
 \begin{displaymath}
 	[\pairing,\bpi_2]_{\cs} = [\pairing,[\pairing, \bpi_1]_{\cs}]_{\cs}
 \end{displaymath}
\end{proposition}
	{We point out that the
		  computations in the proof below are similar to -- but more concise than --  those found in \cite[Lemmas 7.2, 7.3 7.4]{Rogers2013}.}
\begin{proof}
	First note by inspecting on elements $e_i=f_i + \pair{X_i}{\alpha_i}\in \cA$ that:
	\begin{displaymath}
		2 \big(\pairing_- \ca \pi_2\big)~(e_1,e_2,e_3)
		=
		\iota_{[X_1,X_2]} e_3 - \omega(X_1,X_2,X_3) + \cyc			
	\end{displaymath}
	where $\cyc$ denotes sum on all cyclic permutations.
	Using Cartan's magic formula twice:
	
	\allowdisplaybreaks
	\begin{align}
		&\iota_{[X_1,X_2]} e_3 +\cyc
		=
		\notag
		\\
		&=
		\mathcal{L}_{X_1} \iota_{X_2} e_3 - \iota_{X_2} \mathcal{L}_{X_1} e_3 +\cyc=
		\notag
		\\
		&=
		(\iota_{X_1} \dd + \dd \iota_{X_1}) \iota_{X_2} e_3 
		- \iota_{X_2} (\dd \iota_{X_1} + \iota_{X_1} \dd ) e_3 +\cyc=
		\label{eq:bruttocoso}
		\\
		&=
		\dd \iota_{X_1} \iota_{X_2} e_3 
		+ \iota_{X_1} \dd \iota_{X_2} e_3 
		- \iota_{X_2} \dd \iota_{X_1} e_3 
		+ \iota_{X_2}\iota_{X_1}\iota_{X_3} \omega 
		- \iota_{X_2} \iota_{X_1} \mu_1 e_3 +\cyc =
		\notag
		\\
		&=
		(\mu_1 \iota_{X_1} \iota_{X_2} e_3) 
		+ (\iota_{X_1} \mu_1 \iota_{X_2} e_3 - \iota_{X_2} \mu_1 \iota_{X_1} e_3) 
		- (\iota_{X_2} \iota_{X_1} \mu_1 e_3) 
		+ \omega(X_1,X_2,X_3)	+\cyc  	\notag
	\end{align}
	\allowdisplaybreaks[0]
	where, in the penultimate equation, is employed that:
	\begin{displaymath}
		\dd~ e_3 = \dd (\alpha_3 + f_3) = -\iota_{X_3} \omega + \mu_1 e_3
		~.
	\end{displaymath}
	{The first three terms on the r.h.s. of} equation \eqref{eq:bruttocoso} can be recast as follows:
	\allowdisplaybreaks
	\begin{align*}
		\mu_1 \iota_{X_1} \iota_{X_2} e_3 + \cyc
		&=
		\mu_1 \iota_{X_3} \iota_{X_1} e_2 + \cyc 
		=
		\\
		&= 
		\mu_1 \iota_{X_3} (\langle e_1 , e_2\rangle_+ + \langle e_1 , e_2\rangle_- ) + \cyc
		=
		\\
		&=
		- 2 \mu_1 \langle (\langle e_1 , e_2\rangle_+ + \langle e_1 , e_2\rangle_- ), e_3 \rangle_- + \cyc
		=
		\\
		&=
		 2 \mu_1 \pairing_- \ca (\cancel{\pairing_+} + \pairing_- ) (e_1,e_2,e_3) =
		\\
		&=
		2 \Big(\mu_1 \pairing_- \ca \pairing_-\Big)~(e_1,e_2,e_3)		
		=
		\\
		\iota_{X_1} \mu_1 \iota_{X_2} e_3 - \iota_{X_2} \mu_1 \iota_{X_1} e_3 + \cyc 
		&=
		\iota_{X_3} \mu_1 \iota_{X_1} e_2 
		- \iota_{X_3} \mu_1 \iota_{X_2} e_1 + \cyc =
		\\
		&= 
		2 \iota_{X_3} \mu_1 \langle e_1,e_2 \rangle_- + \cyc 
		=
		\\
		&= 
		-4 \langle \mu_1 \langle e_1,e_2 \rangle_-, e_3 \rangle_- + \cyc 
		=
		\\
		&=
		-4 \Big(\pairing_- \ca \mu_1 \ca \pairing_-\Big) (e_1,e_2,e_3)
		\\
		-\iota_{X_2} \iota_{X_1} \mu_1 e_3 + \cyc
		&=
		-\iota_{X_3} \iota_{X_2} \mu_1 e_1 + \cyc 
		=
		\\
		&= 
		- \dfrac{1}{2}\iota_{X_3}
		\big(\iota_{X_2} \mu_1 e_1 - \iota_{X_1} \mu_1 e_2) + \cyc
		=\\
		&=
		+ \iota_{X_3}
		\big(\langle \mu_1 e_1, e_2 \rangle_- - \langle \mu_1 e_2,e_1 \rangle_-\big) + \cyc		
		=\\
		&=
		- 2 \langle
		\big(\langle \mu_1 e_1, e_2 \rangle_- - \langle \mu_1 e_2,e_1 \rangle_-\big), e_3 \rangle_- + \cyc			
		=\\
		&= 2 \Big(\pairing_- \ca (\pairing_- \ca \mu_1)\Big)~(e_1,e_2,e_3)
	\end{align*}
	\allowdisplaybreaks[0]
	Hence, after d\'ecalage, one gets:
	\begin{equation}\label{Eq:pairing-pi2}
	\begin{aligned}
		[\pairing,\bpi_2]_{\cs} 
		=&~\pairing \cs \bpi_2 
		=
		\\
		=&~
		\bpi_1 \cs \pairing\cs \pairing
		-2 \pairing \cs \bpi_1 \cs \pairing 
		+ \pairing \cs \pairing \cs \bpi_1
		=\\
		=&~ [\pairing,[\pairing, {\bpi_1}]_{\cs}]_{\cs} 
	\end{aligned}
\end{equation}
\end{proof}

\begin{proposition}[Commutators of arity $4$]\label{Prop:QuaternaryCommutator}
	\begin{align}
			[\pairing,\bpi_3] 
			=&~ 
			2 \bpi_4 \label{eq:commPairingMu3}
			\\[1em]
			[\pairing^{\cs 2}, \bpi_2 ]
			=&~
			[\pairing^{\cs 2},[\pairing,\bpi_1]]
			= 
			\\[-.5em]
			=&~
			[\pairing,[\pairing^{\cs 2},\bpi_1]] \nonumber 
			\\[1em]
			[\pairing,[\pairing,\bpi_2]] 
			=&~
			[\pairing,[\pairing,[\pairing,\bpi_1]]] = 
			\\[-.5em]
			=&~[\pairing^{\cs 2}, \bpi_2 ] -	2 \alpha (\pairing,\pairing,\bpi_2) 
			= \nonumber
			\\[-.5em]
			=&~  3 \bpi_4 \nonumber
	\end{align}
	where $\alpha$ denotes the associator (see definition \ref{def:gradedAssociator}).
\end{proposition}
\begin{proof}
 Rather straightforward algebraic computations together with the following Lemma \ref{Lemma:BoringAssociator}.
\end{proof}
\begin{lemma}[{A} recurrent associator]\label{Lemma:BoringAssociator}
	\begin{displaymath}
		\alpha(\cs; \pairing,\pairing, \bpi_2) = \frac{1}{2}\Big( [\pairing^{\cs 2},\bpi_2] -3 \bpi_4\Big)
	\end{displaymath}
\end{lemma}
\begin{proof}
	Observe that
	\begin{equation}\label{eq:obscure}
		2 \pairing_- \ca \pi_2 = K + 3 \pi_3
	\end{equation}
	where the auxiliary operator $K$ is given by the following equation
	\begin{equation*}
		\begin{aligned}
			K(e_1,e_2,e_3)=&~
		(\pairing_+ + \pairing_-)\ca \pi_2 (e_1,e_2,e_3) =
		\\
		=&~
		\iota_{[X_1,X_2]}e_3 + \iota_{[X_2,X_3]}e_1 + \iota_{[X_3,X_1]}e_2
		~.
		\end{aligned}
	\end{equation*}	
	According to equation \eqref{Eq:explicitassociators},  the following  holds:
	\begin{equation}\label{eq:Kalpha}
		\begin{aligned}
			\alpha(\ca;\pairing_-,&\pairing_-,\pi_2) (e_1,e_2,e_3,e_4)
			=
			\\
			=& (-)^2\frac{1}{2}\iota_{[X_3,X_4]}\langle e_1,e_2 \rangle_- + \unsh{(2,2)}
			=\\
			=&
			\frac{1}{2}\biggr(			
			+ \iota_{[X_3,X_4]}\langle e_1,e_2\rangle_- 
			+ \iota_{[X_1,X_2]}\langle e_3,e_4\rangle_-
			- \iota_{[X_2,X_4]}\langle e_1,e_3\rangle_-
			+
			\\
			&\phantom{\frac{1}{2}\biggr(}- \iota_{[X_1,X_3]}\langle e_2,e_4\rangle_-
			+ \iota_{[X_2,X_3]}\langle e_1,e_4\rangle_-
			+ \iota_{[X_1,X_4]}\langle e_2,e_3\rangle_-
			\biggr)~,
		\end{aligned}
	\end{equation}	
	{where $\unsh{(2,2)}$ denotes sum over all $(2,2)$-unshuffles.}
	\\
	{On the other hand, one finds that}
	\begin{align*}
		\big(\pairing_-&\ca K\big) (e_1,e_2,e_3,e_4)
		=
		\\
		=&
		(-)^{|K|} \langle K(e_1,e_2,e_3),e_4 \rangle_- + \unsh{(3,1)} =
		\\
		=&
		-\langle K(e_1,e_2,e_3),e_4\rangle_-
		+ \langle K(e_1,e_2,e_4),e_3\rangle_- +
		\\
		&
		- \langle K(e_1,e_3,e_4),e_2\rangle_-
		+ \langle K(e_2,e_3,e_4),e_1\rangle_-
		=\\
		=&\frac{1}{2}\Big\{
			+\iota_{X_4}(\iota_{[X_1,X_2]}e_3+\iota_{[X_2,X_3]}e_1+\iota_{[X_3,X_1]}e_2)+
		\\
		&\phantom{\frac{1}{2}\big\{}
			-\iota_{X_3}(\iota_{[X_1,X_2]}e_4+\iota_{[X_2,X_4]}e_1+\iota_{[X_4,X_1]}e_2)+
		\\
		&\phantom{\frac{1}{2}\big\{ }
			+\iota_{X_2}(\iota_{[X_1,X_3]}e_4+\iota_{[X_3,X_4]}e_1+\iota_{[X_4,X_1]}e_3)+
		\\
		&\phantom{\frac{1}{2}\big\{}
			-\iota_{X_1}(\iota_{[X_2,X_3]}e_4+\iota_{[X_3,X_4]}e_2+\iota_{[X_4,X_2]}e_3)
		\Big\}
		=
		\\
		=&+ \iota_{[X_1,X_2]}\langle e_3, e_4\rangle_- 
			- \iota_{[X_1,X_3]}\langle e_2, e_4\rangle_-
			+ \iota_{[X_1,X_4]}\langle e_2, e_3\rangle_-
		+\\&
			+ \iota_{[X_2,X_3]}\langle e_1, e_4\rangle_-
			- \iota_{[X_2,X_4]}\langle e_1, e_3\rangle_-
			+ \iota_{[X_3,X_4]}\langle e_1, e_2\rangle_-
		=\\=&
			2~ \alpha(\ca;\pairing_-,\pairing_-,\pi_2) (e_1,e_2,e_3,e_4)
			~,
	\end{align*}
{using eq. \eqref{eq:Kalpha} in the last equality. In other words,}
	\begin{displaymath}
		\pairing_- \ca K = 2~ \alpha(\ca;\pairing_-,\pairing_-,\pi_2)
		~.
	\end{displaymath}
	Plugging this last result into the equation obtained by {composing} equation \eqref{eq:obscure} with $\pairing_-$ from the left one gets:
	\begin{displaymath}
		2~\pairing_- \ca (\pairing_-\ca \pi_2 ) =
		2~\alpha(\ca;\pairing_-,\pairing_-,\pi_2) +
		3~ \pairing_- \ca \pi_3
		~.
	\end{displaymath}			
	Applying on the l.h.s. the definition of the associator, using equation \eqref{eq:commPairingMu3},  and remark \ref{rem:two}, one gets
	\begin{displaymath}
		[ \pairing_-^{\ca 2}, \pi_2]  = 2~\alpha(\ca; \pairing_-, \pairing_-, \pi_2 ) + 3\,\pi_{4}
		~.
	\end{displaymath}		
	The claim follows after d\'ecalage.
\end{proof}

The following technical lemma will be used in remark \ref{Rem:DiagramSituation} to express the gauge transformation of a \momap in terms of the pairing.
The key idea is to take advantage of the operator $\pairing_-$ on the space $\cA$ by expressing the operation of inserting several vector fields in a given differential form as a ``power'' of the pairing:
	\begin{lemma}[Insertions as pairing]\label{lemma:InsertionsAsPairing}
		Consider the graded vector space $\tilde{\cV}:=\mathfrak{X}(M)\oplus\Omega(M)[k]$. 
		Denote by $\rho$ the standard projection $\rho:\tilde{\cV} \twoheadrightarrow\mathfrak{X}(M)$.
		Given a $k$-form  $B$, \ie an element in $\ker(\rho)\subset \tilde{\cV}$, and given  vector fields
		 $x_i$,   
		the following equation holds {for all $m$}:
		\begin{displaymath}
			\left(
			\pairing_-^{\ca m}
			\right)(B, x_1,\dots, x_m)
			=
			\left(-\varsigma(m) \cdot \frac{m!}{2^m}
			\right)~
			\iota_{x_m}\dots \iota_{x_1} B
			~.
		\end{displaymath}
Here the left-hand side denotes the evaluation of operator $\pairing_-^{\ca m}$, see definition \ref{def:pairing-},
on the element $ B\otimes x_1\otimes \dots \otimes x_m \in \tilde{\cV}^{\otimes m+1}$.
	\end{lemma}
	\begin{proof}
		By induction. From equation \eqref{Eq:RNProducts-explicit}, given two vector fields $x_1,x_2$, which are degree $0$ elements in $\widehat{\cV}$ and a differential form $B$ one has
		\allowdisplaybreaks
		\begin{align*}
				\pairing_-\ca \pairing_- ~& (B,x_1,x_2)=
				\\
				&=
				(-)^{|\pairing_-|}
				\mkern-20mu
				\sum_{\sigma \in \ush{1,1}}
				\mkern-20mu
				\chi(\sigma)
				\langle\langle B,x_{\sigma_1}\rangle_-, x_{\sigma_2} \rangle_-
				=
				\\
				&=
				\left(- \frac{2!}{(-2)^2}\right)~\iota_{x_2}\iota_{x_1} B
				=
				\\
				&=
				\left(-\varsigma(2) \dfrac{2!}{2^2}\right)~\iota_{x_2}\iota_{x_1}
				~,
		\end{align*}
		\allowdisplaybreaks[0]
		where $\chi(\sigma)=(-)^\sigma \epsilon(\sigma)$ denotes the odd Koszul sign.
		Assuming now
		\begin{displaymath}
			\pairing_-^{\ca(m-1)}\ca \pairing_-~(B,x_1,\dots,x_{m})
			=
			\left(-\varsigma(m) \frac{m!}{2^m}\right) \iota_{x_m}\dots \iota_{x_1} B
			~,
		\end{displaymath}
		it follows that
		
		\begin{displaymath}
			\mathclap{
			\begin{aligned}
				\pairing_-^{\ca(m)}&\ca \pairing_- 
				~(B,x_1,\dots,x_{m+1})
				=
				\\
				=&
				\pairing_- \ca 
				\left(
					\pairing_-^{\ca (m-1)}\ca \pairing_- 
				\right) 
				~ (B,x_1,\dots,x_{m+1})
				= 
				\\
				=&
				(-)^{|\pairing_-^{\ca (m)}|}
				\mkern-20mu\sum_{\sigma \in \ush{m,1}} \mkern-20mu
				\chi(\sigma) 
				\left\langle 
					\left( 
						\pairing_-^{\ca(m-1)}\ca \pairing_-~(B,x_1,\dots,x_{m})
					\right)
					, x_{\sigma_{m+1}}				
				\right \rangle_-
				=
				\\
				=&
				(-)^m
				\mkern-20mu\sum_{\sigma \in \ush{m,1}} \mkern-20mu
				\chi(\sigma)
				\left(-\varsigma(m) \frac{m!}{2^m}\right)\left(-\frac{1}{2}\right) 
				\iota_{x_{\sigma_{m+1}}}\iota_{x_{\sigma_m}}\dots \iota_{x_{\sigma_1}} B
				=
				\\
				=&
				\left(
					(-)^m \varsigma(m) \frac{(m+1)!}{2^{m+1}}
				\right)
				\iota_{x_{m+1}}\dots \iota_{x_1} B
				~,
			\end{aligned}
			}
		\end{displaymath}
		hence the claim.
		Formula \eqref{Eq:explicitassociators} ensures associativity in the first equality.
	\end{proof}
	\begin{remark}\label{Rem:piccoloerroredisegnoneldraft}
		Notice that the above statement can be re-expressed by singling out the contraction with the differential from $B$.
		For any given differential form $B$, seen as a degree $|B|-k$ element in $\widetilde{\cV}$, one has:
		\begin{displaymath}
			\pairing_-^{\ca m} \ca \langle B, \cdot \rangle_- = (-)^{m(|B|-k)} \pairing_-^{\ca (m+1)} B \otimes \Unit_{m+1}
			~.
		\end{displaymath}
		This follows from inspecting on vector fields $x_i$, seen as degree $0$ elements in $\widetilde{\cV}$:
		\begin{displaymath}
			\begin{aligned}
				\pairing_-^{\ca (m)} \ca \langle B, \cdot \rangle_- &(x_1,\dots,x_{m+1}) 
				=
				\\
				&=
				(-)^{m(|\langle B, \cdot \rangle_-|)}
				\sum_{\sigma \in S_{m+1}}\chi(\sigma) \langle\dots\langle B, x_{\sigma_1}\rangle_-,\dots,x_{\sigma_{m+1}}\rangle_-
				=
				\\
				&=
				(-)^{m(|B|-k+1)}(-)^{m|\pairing_-|}
				\pairing_-^{\ca(m+1)} (B,x_1,\dots, x_{m+1})				
				~.			
			\end{aligned}
		\end{displaymath}
	
	\end{remark}
	
	As a corollary\footnote{This corollary is not essential in this chapter, but we expect that it should be useful to extend theorem \ref{thm:iso} to arbitrary values of $n$, just as Lemma \ref{lem:mun} below.},
	one can express all higher multibrackets $\bpi_k$ in terms of $\bpi_3$:	
		\begin{corollary}\label{cor:higherpi}
		\begin{displaymath}
			\bpi_n = \left( \frac{2^{n-3}}{n!} 3!\right) \pairing^{\cs(n-3)} \cs \bpi_3
		\end{displaymath}
	\end{corollary}
	\begin{proof}
		This can be proved iterating lemma \ref{lem:rogersRecurFormula}, or employing lemma \ref{lemma:InsertionsAsPairing} in the following way.
	{
	Consider Hamiltonian pairs $e_i= \pair{x_i}{\alpha_i}\in \cA$, one has:}
		\begin{displaymath}
			\begin{aligned}
				\pi_{n}& (e_1,\dots e_n ) 
				= 
				\varsigma(n) \iota_{x_n}\dots \iota_{x_1} \omega 
				=
				\\
				&=
				\varsigma(n)\varsigma(3) \frac{1}{\# \ush{3,n-3}}
				\mkern-50mu\sum_{\quad\qquad\sigma \in \ush{3,n-3}}\mkern-50mu
				\chi(\sigma)~
				\iota_{x_{\sigma_n}}\dots \iota_{x_{\sigma_4}} ~ \pi_3(x_{\sigma_1},x_{\sigma_2}x_{\sigma_3})
				=
				\\
				&=
				- \varsigma(n)\frac{3!(n-3)!}{n!}
				\left(-\varsigma(n-3)\frac{2^{n-3}}{(n-3)!}\right)\cdot
				\\
				&\qquad \cdot
				\mkern-50mu\sum_{\quad\qquad\sigma \in \ush{3,n-3}}\mkern-50mu
				\chi(\sigma)~
				\pairing_-^{\ca (n-3)} \cdot \pi_3 \otimes \mathbb{1}_{n-3} ~(x_{\sigma_1}\dots x_{\sigma_n})
				=
				\\
				&=
				\left(
				\varsigma(n)\varsigma(n-3)(-)^{n-3}
					\frac{2^{n-3}}{n!}
					3!				
				\right)
				\left(\pairing_-^{\ca (n-3)} \ca \pi_3 \right){(e_1,\dots e_n)}
				=
				\\
				&=
				\left(
					\frac{2^{n-3}}{n!}
					3!				
				\right)
				\left(\pairing_-^{\ca (n-3)} \ca \pi_3\right){(e_1,\dots e_n)}
				~.
			\end{aligned}
		\end{displaymath}
		The statement follows after d\'ecalage.
	\end{proof}

\section{Vinogradov algebroid and $L_\infty$-algebra}\label{Sec:Vinoids}
To any smooth manifold $M$, and index $n\in \NN$, one can associate an auxiliary geometric structure called \emph{Vinogradov algebroid}.
In this section, we will recall this definition and some basic properties.
In particular we will be interested in describing how one can produce a $L_\infty$-algebra out of this geometric data and we will discuss the algebraic properties enjoyed by the corresponding set of multibrackets with the language of the \RN product.

A (standard) Vinogradov algebroid is a slightly generalized version of a Courant algebroid \cite{Courant1990}.
\begin{definition}[(Standard ) $n$-Vinogradov algebroid]\label{def:standardVinalgoid}
	Given a smooth manifold $M$, fixed a natural number $n\geq 1$, the \emph{(standard) Vinogradov algebroid} consists of the data $(E^n ,  \rho, \pairing_{-}, [\cdot,\cdot]_C)$ where:
	\begin{itemize}
		\item $E^n$ denotes the vector bundle over $M$ given by
			\begin{displaymath}
				E^n := TM \oplus \Lambda^{n-1} T^\ast M ~,
			\end{displaymath}	
				we will denote elements of $E^n$ as $e = \pair{X}{\alpha}$;
		\item $\rho \colon E^n \to  TM$ is a vector bundle morphism given by the first projection, also called \emph{the anchor};
		\item $\pairing_{\pm} \colon E^n \otimes E^n \to \Lambda^{n-2} T^\ast M $ are binary bundle maps given by 
		\begin{displaymath}
		\left(\pair{X_1}{\alpha_1},\pair{X_2}{\alpha_2}\right) \mapsto 
	 \frac{1}{2}\left( \iota_{X_1}\alpha_2 \pm \iota_{X_2}\alpha_1 \right)~;
		\end{displaymath}
		\item $[\cdot,\cdot]_C$ is a skew-symmetric bracket on the vector space of sections  $\Gamma(E^n)$, called \emph{higher Courant bracket}, given by 
		\begin{displaymath}
		\left(\pair{X_1}{\alpha_1},\pair{X_2}{\alpha_2}\right) \mapsto 	
		\pair{[X_1,X_2]}{\mathcal{L}_{X_1}\alpha_2 - \mathcal{L}_{X_2}\alpha_1 - \dd 
		\left\langle \pair{X_1}{\alpha_1},\pair{X_2}{\alpha_2}\right\rangle_-}	
		~.
		\end{displaymath}
	\end{itemize}
\end{definition}
\begin{notation}
	In the following, we will drop the $n$-prefix; everything should be clear from the context.
	We will also drop the "standard" term most of the times. 
	This adjective is due to the existence of an "abstract" notion of Vinogradov algebroid, see remark \ref{rem:abstractVino} below, that will not be needed in the following.
\end{notation}

Consider now a pre-$n$-plectic form $\omega\in \Omega^{n+1}(M)$. The degree of the form select a specific integer and thus a specific standard $n$-Vinogradov algebroid.
Furthermore, the higher Courant bracket can be ``twisted'' by the closed form $\omega\in \Omega^{n+1}(M)$.
\begin{definition}[(Standard) $\omega$-twisted Vinogradov algebroid]\label{def:Vinalgoid}
	Let be $M$ a smooth manifold, and $\omega\in\Omega^{n+1}(M)$ a closed differential form. 
	The  \emph{Vinogradov algebroid twisted by $\omega$} consists of the  data\\
 $(E^n,  \rho, \pairing_{-}, [\cdot,\cdot]_\omega)$
	where $[\cdot,\cdot]_\omega: \Gamma(E^n)\otimes \Gamma(E^n) \to \Gamma(E^n)$ is defined by 
	\begin{displaymath}
		[e_1,e_2]_\omega =  [e_1,e_2]_C  + \pair{0}{\iota_{X_1}\iota_{X_2} \omega}
		~.
	\end{displaymath}
\end{definition}
\begin{example}[Lie and Courant algebroids]
	When $n=1$, the operators $\pairing_\pm$ are trivial. 
	Hence one recovers, as a particular case, the definition of standard and twisted Lie algebroid (\cf  examples \ref{ex:StandardLieAlgbroid} and \ref{ex:TwistedLieAlgbroid}).
	\\
	If $n=2$, definition \ref{def:standardVinalgoid} gives the definition of the standard Courant algebroid with Courant bracket.
\end{example}

\begin{remark}[About the naming]
	The choice to name "Vinogradov algebroid" this, somewhat natural, higher generalization of the Courant algebroid is borrowed from Ritter and Saemann \cite{Ritter2015a}.
	They pointed out that the key structure, 
	given in particular by the bracket on sections, 
	was first studied by Vinogradov in \cite{zbMATH04172022}.
	In other sources, see for example \cite{Zambon2012} or \cite{Bi2011a}, the same object is simply called \emph{higher Courant algebroid}.
\end{remark}

\begin{remark}["Abstract" Vinogradov algebroids]\label{rem:abstractVino}
	When introducing the Vinogradov algebroid we are mostly employing a constructive, "hands-on", approach here.
	Our definition appears as a "standard" example, \cf  the notion of "standard" vs. "abstract" Courant algebroid, in \cite[\S 5, Def. 5.6]{Grutzmann2015}.
	More conceptually, this notion can be framed in the language of graded geometry as a prototypical $NQ$-manifold, see \cite[\S 3.3.]{Deser2018b} and  \cite{Ritter2015a}.
	In turn, the latter concept can be interpreted in terms of \emph{horizontal categorification}\footnote{This justifies the "-oids" suffix in the name.} since $NQ$-manifolds are related to \emph{symplectic Lie-$n$ algebroids}.
\end{remark}
\begin{remark}[Vinogradov algebroid morphisms]\label{rem:VinoidsMorphism}
	The precise notion of morphism between Vinogradov algebroids can be given by mimicking the corresponding definition for the Courant ($n=2$) algebroid case (see \cite{Courant1990} for the original definition or also \cite[\S 1.3]{Li-Bland2009} and \cite[\S 2.2]{Bursztyn2008}).
	\\
	We will not need the full fledged apparatus here.
	In this chapter, precisely in \ref{Sec:VinoGauge}, we will only deal with morphisms between Vinogradov algebroid on the same smooth manifold $M$ and twisted by differential forms $\omega$ and $\widetilde{\omega}$ of the same degree.
	In the latter situation, a Vinogradov morphism 
		\begin{displaymath}
			\Psi:~
			(E^n, \rho, \pairing_\pm, [\cdot,\cdot]_\omega) \to 
			(E^n, \rho, \pairing_\pm, [\cdot,\cdot]_{\widetilde{\omega}})		
		\end{displaymath}
	will be simply given by a vector bundle automorphism $\Psi:E^n\to E^n$ that is compatible with the anchor, \ie the following diagram commutes in the category of smooth manifolds
	\begin{displaymath}
		\begin{tikzcd}[column sep = small, row sep = small]
				E^n \ar[rr,"\Psi"] \ar[ddr,bend right = 45,"\rho"'] \ar[dr,two heads]& &
				E^n \ar[ddl,bend left = 45,"\rho"] \ar[dl,two heads]
				\\
				& M 
				\\
				& TM \ar[u,two heads]
		\end{tikzcd}
		~,
	\end{displaymath}
	and such to preserve the symmetric pairing and the higher Courant bracket in the following sense:
	\begin{displaymath}
		\begin{aligned}
				\Psi \circ \pairing_+ &= \pairing_+ \circ (\Psi\otimes \Psi) ~;
				\\
				\Psi \circ [\cdot,\cdot]_\omega &= [\cdot,\cdot]_{\widetilde{\omega}} \circ (\Psi \otimes \Psi) 				~.
		\end{aligned}
\end{displaymath}		
\end{remark}
For the sake of completeness, we mention here some basic properties enjoyed by the twisted higher Courant bracket.
\begin{proposition}[\emph{\cite[Thm 2.2]{{Bi2011a}}}]\label{prop:VinoAlgoidsProperties}
	Let be $(M,\omega)$ be a pre-$n$-plectic manifold. 
	Consider the associated Vinogradov algebroid 
	$$E^n=(TM\oplus\wedge^{n-2}T^\ast M, \rho, \pairing_\pm,[\cdot,\cdot]_\omega) ~.$$
	The higher Courant bracket $[\cdot,\cdot]_\omega$ satisfies the following properties
	\begin{enumerate}
		\item $[\cdot,\cdot]_\omega$ is skew-symmetric;
		\item $[\cdot,\cdot]_\omega$ satisfies the Jacobi equation up to an exact term. Namely
		\begin{equation}\label{eq:CourantFailsJacobi}
			[[e_1,e_2]_\omega,e_3]_\omega +\cyc = \d T_\omega(e_1,e_2,e_3)
			\qquad \forall e_i \in \Gamma(E^n)
		\end{equation}
		where
		\begin{equation}\label{eq:CourantTernaryOp}
			\morphism{T_\omega}
			{(\Gamma(E^n))^{\otimes 3}}
			{\Omega^{n-2}(M)}
			{(e_1,e_2,e_3)}
			{\frac{1}{3}\langle [e_1,e_2]_\omega, e_3 \rangle_+ + \cyc} ~.
		\end{equation}
		\item Regard $\Gamma(E^n)$ as a $C^\infty(M)$-module. $[\cdot,\cdot]_\omega$ is not $C^\infty(M)$-linear. In the standard case, one has
		\begin{displaymath}
			[e_1,f e_2]_C = f [e_1,e_2]_C + (\mathcal{L}_{X_1}f) e_2 - \dd f \wedge \langle e_1,e_2\rangle_+  \qquad \forall e_i \in \Gamma(E^n)~.
		\end{displaymath}
		\item $[\cdot,\cdot]_\omega$ is compatible with the anchor:
			\begin{displaymath}
				\rho\left([e_1,e_2]_\omega \right) = [\rho(e_1),\rho(e_2)]_\omega
				~.
			\end{displaymath}
	\end{enumerate}
\end{proposition}
\begin{proof}
	The proof relies on the same computations routinely performed for the standard Courant algebroid case. 
	See for instance \cite[Prop. 4.7]{Ritter2015a}.
\end{proof}

\begin{remark}[Comparison with the literature]
	We stress that there are several different conventions that one could adopt when defining Vinogradov algebroids and, in turn, Courant algebroids.
	We briefly compare our choices to the references in bibliography.
	\begin{itemize}
		\item The pairing operators $\pairing_\pm$ are often defined without the $1/2$ prefactor (see for example \cite[eq. (1)]{Zambon2012} and \cite[ex. 1]{Rogers2013}).
		Here, we are adopting the same convention of \cite[eq. (2.2.1),(2.2.2)]{Courant1990} and \cite[eq. (1)]{Bi2011a}.
		\item The choice of the sign used to "twist" the standard Courant bracket also differs across the literature.
		Here, like in \cite[eqn. 5.6]{Rogers2013} and \cite[eqn. 4.12]{Ritter2015a}, we are choosing to twist the Courant bracket $[e_1,e_2]$, for any given $e_i\in \Gamma(E^n)$ with $\rho(e_i)=x_i$, by adding the term $\iota_{x_1}\iota_{x_2}\omega$.
		With the notation\footnote{Recall that
			\begin{displaymath}
				\morphism{\iota_\rho^k\omega}
				{(\Gamma(E^n))^{\otimes k}}
				{\Omega^{n+1-k}(M)}
				{(e_1,\dots,e_k)}
				{\varsigma(k) \iota_{x_k}\dots \iota_{x_1} = (-)^{k+1}\iota_{x_1}\dots\iota_{x_k}\omega ~,}
\end{displaymath}	
		hence, in particular, one has $\iota_{x_1}\iota_{x_2}\omega = - \iota^2_\rho \omega$.}
		introduced in remark \ref{Rem:SignedMultiContraction}, this convention can compactly stated as:
		\begin{equation}
			[\cdot,\cdot]_\omega = [\cdot,\cdot]_C - \iota^2_\rho \omega ~.
		\end{equation}
		Notice that in \cite[\S 3.7]{Gualtieri2004}, \cite[eq (2.3)]{Gu2} and \cite[\S 2]{Zambon2012} it is preferred the opposite sign.
		\item The previous choices also affect the definition of the ternary operator $T_\omega$ introduced in proposition \ref{prop:VinoAlgoidsProperties}.
		In the above notation, this can be compactly written as
		\begin{displaymath}
			T_\omega = T_0 + \dfrac{1}{2}\iota^3_\rho \omega ~,
		\end{displaymath}
		where $T_0$ denotes the ternary operator pertaining to to the standard higher Courant bracket.
 		Our expression can be easily compared with \cite[Def 4.1]{Rogers2013}, \cite[eq (4.12)]{Ritter2015a} and \cite[Thm. 4.2]{Bi2011a}.
	\end{itemize}
\end{remark}
The following remark supports the observation that multisymplectic manifolds are deeply related to Vinogradov algebroids.
\begin{remark}[Higher Dirac structures \cite{Zambon2012}]
	Exactly as symplectic structures are special cases of Dirac structures, given by the graph subbundle inside the standard Courant algebroid, also multisymplectic structures can be regarded as special case of \emph{higher Dirac structures}.
	The latter has been identified in \cite[Def. 3.1, Prop. 3.2 and 3.7]{Zambon2012} as a 
	\emph{involutive, Lagrangian subbundles} of $E^n$. (See also \cite[\S 4]{Bi2011a}.)
\end{remark}

Proposition \ref{prop:VinoAlgoidsProperties}, especially items $1.$ and $2.$, 
leads us naturally to replicate the reasoning described in section \ref{Section:RogersObservables}, \ie the construction of the multisymplectic observables $L_\infty$-algebra, to the case of Vinogradov algebroids.

\subsection{Vinogradov's $L_{\infty}$-algebra}\label{subsec:Vinogradov}
Consider a $\omega$-twisted Vinogradov algebroid $E^n$.
The upshot of proposition \ref{prop:VinoAlgoidsProperties} is that, similarly to what has been observed for the vector space $\Omega^{n-1}_{\ham}(M,\omega)$ in section \ref{Section:RogersObservables},
the higher Courant bracket on the space of section $\Gamma(E^n)$
fails to be a Lie algebra structure.
Namely, the higher Courant bracket $[\cdot,\cdot]_\omega$ is skew-symmetric but satisfies the Jacobi equation only modulo the exterior derivative of the ternary operator $T_\omega$ (see equation \eqref{eq:CourantFailsJacobi}).
This naturally prompts to look for a suitable completion of $\Gamma(E^n)$ to give a full-fledged $L_\infty$-algebra.

Roytenberg and Weinstein showed how this $L_\infty$-algebra can be constructed in the case of ordinary Courant algebroids \cite[Thm 4.3]{Roytenberg1998}. 
This result has been extended by Zambon in \cite[Prop. 8.1 and 8.4]{Zambon2012} to  Vinogradov algebroids.
Summing up, for any $\omega$-twisted Vinogradov algebroid there is an associated $L_n$-algebra, given by the following definition:
 \begin{definition}[$L_\infty$-algebra of a twisted Vinogradov algebroid]\label{def:vinolinfty}
Given a twisted Vinogradov algebroid 
$$ (E^n,\omega)=(TM\oplus \Lambda^{n-1}T^\ast M , \rho, {\pairing_{-}}, [\cdot,\cdot]_\omega) ~,$$
we denote the the associated $L_n$-algebra structure as 
$$L_{\infty}(E^n,\omega) = ({\cV}, \{\mu_k\})$$
Its underlying graded vector space is given by
\begin{equation}\label{eq:VSpace}
	{\cV^i}:=
	\begin{cases}
		\mathfrak{X}(M)\oplus \Omega^{n-1}(M)
		 & \quad ~\text{if } i=0
		 \\
		 ~\Omega^{n+i-1}(M) 
		 &  \quad~\text{if } 1-n \leq i\leq -1
		 \\			
		~0 & \quad ~\text{otherwise}
		 ~.
	\end{cases}
\end{equation}	
The actions of non-vanishing multi-brackets (up to permutations of the entries) on arbitrary vectors $\varv_i=f_i\oplus e_i \in {\cV}$,
with $e_i = \pair{X_i}{\alpha_i} \in \mathfrak{X}(M)\oplus \Omega^{n-1}(M) = \Gamma(E^n)$ and 
$f_i \in \bigoplus_{k=0}^{n-2}\Omega^k(M)$,
are given as follows:
\begin{itemize}
	\item unary bracket:
		$$\mu_1 \left(f\right) =   \dd f ~;$$
	\item binary bracket:
		\begin{displaymath}
			\begin{aligned}
				\mu_2 \left(e_1,e_2\right) 	=& [e_1,e_2]_\omega 
				= \pair{[X_1,X_2]}{\dd \langle e_1, e_2\rangle_- 
				+ (\iota_{X_1}\dd\alpha_2 - \iota_{X_2}\dd\alpha_1 + \iota_{X_1}\iota_{X_2}\omega)}
				~;
				\\
				\mu_2 \left(e_1,f_2\right) =& -\mu_2(f_2,e_1) = 
				\frac{1}{2} \mathcal{L}_{X_1} f_2 = \langle e_1, \dd f_2 \rangle_-
				~;
			\end{aligned}
		\end{displaymath}
	\item ternary bracket:
		\begin{displaymath}
			\begin{aligned}
				\mu_3 (e_1, e_2, e_3) =& -T_\omega(e_1,e_2,e_3) = 
				-\frac{1}{3} \langle[e_1,e_2]_\omega,e_3 \rangle_+ + \cyc \\
				\mu_3 (f_1, e_2, e_3) =& -\frac{1}{6}
				\left(
				\frac{1}{2}(\iota_{X_1}\mathcal{L}_{X_2} 
				- \iota_{X_2}\mathcal{L}_{X_1}) + \iota_{[X_1,X_2]}
				\right)f
			\end{aligned}
		\end{displaymath}	
	\item $k$-ary bracket for $k \ge 3$ an \emph{odd} integer:
	 
		\begin{equation}\label{eq:VinoMultibrakAllaZambon_1}
			\begin{aligned}
				\mu_k(\varv_0,\cdots,\varv_{k-1})
				=&
				\left(\sum_{i=0}^{k-1} {(-)^{i-1}\mu_k(f_i+\alpha_i,X_0,\dots,\widehat{X_i},\dots,X_{k-1})}\right)
				+\\
				&+(-)^{\frac{k+1}{2}} \cdot k \cdot B_{k-1} \cdot 
				\iota_{X_{k-1}}	\dots \iota_{X_{0}} \omega			
				~;
			\end{aligned}
		\end{equation}
		where 
		\begin{equation}\label{eq:VinoMultibrakAllaZambon_2}
			\begin{aligned}
			\mu_k  &(f_0+\alpha_0, X_1,\dots,X_{n-1}) =
			\\
			&=~c_k
			\mkern-30mu\sum_{\quad1\le i<j\le k-1}\mkern-20mu
			(-1)^{i+j+1}
			\iota_{X_{k-1}}\dots   
  			\widehat{\iota_{X_{j}}}\dots \widehat{\iota_{X_{i}}}\dots
				\iota_{X_{1}} ~ [f_0+\alpha_0,X_i,X_j]_3~.
			\end{aligned}
		\end{equation}			
In the above formula,		$[\cdot,\cdot,\cdot]_3 = -T_0$ denotes the ternary bracket 
associated to the untwisted ($\omega=0$) Vinogradov algebroid, and $c_k$ is a numerical constant
		\begin{equation}\label{eq:UglyCoefficient}
			c_k= (-)^{\frac{k+1}{2}}\frac{12~B_{k-1}}{(k-1)(k-2)}~,
		\end{equation}
	containing the \emph{Bernoulli numbers}\footnote{{Hence $B_0=1, B_1=-\frac{1}{2}, B_2=\frac{1}{6}$, and $B_k=0$ for odd $k\neq 1$.}} $B_{n}$.
				\end{itemize}
\end{definition}

\begin{remark}\label{Remark:semplifica-conti}
	Notice that Vinogradov's brackets $\mu_k$ with $k\geq 2$ vanishes unless $k-1$ entries are elements of degree zero.
For $k\ge 2$, Rogers' brackets $\pi_k$ vanish unless all entries are elements in degree zero (\cf  remark \ref{rem:two}).
\end{remark}
\begin{remark}[On the origin of definition \ref{def:vinolinfty}]
	The definition $L_\infty(M,\omega)$ has been worked out in  \cite[Prop. 8.1 and 8.4]{Zambon2012} relying on a result by Getzler \cite{Getzler1991} (see also theorem \ref{Thm:Getzler} in chapter \ref{Chap:Linfinity}) and on some observations in graded geometry.
	Briefly, they noticed that $\cV$ can be obtained as a certain truncation of the graded vector space of smooth functions on the graded manifold $T^\ast[n]T[1]M$.
	Namely, denoting by $\mathscr{C}= C^{\infty}(T^\ast[n]T[1]M)$ the space of smooth functions on the above-mentioned graded manifold, one can show that
	\begin{displaymath}
			\cV[1] = \trunc_{n} \mathscr{C}[n]~.
	\end{displaymath}
	As such, $\cV[1]$ inherits from $\mathscr{C}$ a binary bracket $\lbrace\cdot,\cdot\rbrace$, given by the canonical Poisson bracket on the cotangent bundle (\ie the graded analogue of example \ref{Ex:Multicotangent} in the $1$-plectic case, see \cite{CATTANEO2011}),
	and a unary bracket given by the de Rham differential.
	\\
	Summing up $(\mathscr{C},\d,\lbrace\cdot,\cdot\rbrace)$ constitute a DGLA (see definition \ref{def:DGLA}).
	Thence, one could readily apply the machinery given by theorem \ref{Thm:Getzler}  endowing $\cV=\mathscr{C}[n][-1]$ with the above $L_\infty$-structure.
\end{remark}

\begin{remark}[Vinogradov algebroids as $L_\infty$-algebroids]
		Recall that, in layman terms, a Lie algebroid can be thought of as the geometric data encoding a certain "well-behaving" infinite dimensional Lie algebra.
		In other words, it can be seen as a $\infty$-dimensional Lie algebra whose elements are sections of a vector bundles, taken together with a Lie algebra morphism $\rho$ into the Lie algebra $\mathfrak{X}(M)$ giving a sort of "representation" in terms of vector fields. 
		\\
		Definition \ref{def:vinolinfty} tells us that the same point of view could be also loosely applied to a Vinogradov algebroid. 
		In this spirit, one could say that definition \ref{def:Vinalgoid} is the geometric data giving a certain $\infty$-dimensional $L_n$-algebra whose elements are differential forms of a smooth manifold $M$ up to the degree $n-1$.
		\\
		This observation drops a hint about the interpretation of Vinogradov algebroids as "$L_n$-algebroids" mentioned in remark \ref{rem:abstractVino}.
\end{remark}

\subsection{Hamiltonian subalgebroid}
Let be $(M,\omega)$ be a $n$-plectic manifold, consider the corresponding $\omega$-twisted Vinogradov algebroid $E^n$.
Consider also the two corresponding $L_\infty$-algebras $L_\infty(M,\omega)$ and $L_\infty(E^n,\omega)$.
We denoted respectively by $L$ and $\cV$ the underlying graded vector spaces (see equations \eqref{eq:VSpace} above and \eqref{eq:Lspace} in chapter \ref{Chap:MultiSymplecticGeometry}).
There is an obvious sequence of inclusions at the level of graded vector spaces
\begin{equation}\label{eq:GVSinclusions}
	\begin{tikzcd}
		L \ar[r,equal,"\sim"] & \cA \ar[r,hook,"h"] & \cV \ar[r,hook] & \widetilde{\cV}
	\end{tikzcd}
\end{equation}
where $\cA$ is the graded vector space introduced in equation \eqref{eq:Aspace} and $\widetilde{\cV}$ has been introduced in equation \eqref{eq:VtildeSpace}.
The first isomorphism is the map $\vartriangle$ defined in equation \eqref{eq:trianglemap}, the second one is the standard inclusion of $\cA$ as a subspace of $\cV$ and the last one is a truncation of the identity map.

With the following lemma, we want to point out  that there are two $L_n$-structures naturally induced by the $n$-plectic form $\omega$ on the cochain complex $(\cA,\pi_1)$. Namely, the structure $\{\pi_k\}$ given by Rogers, and the multibrackets $\{\mu_k\}$ obtained by restricting the Vinogradov's $L_\infty$-algebra.
\begin{lemma}
	\noindent
	(i) The inclusions expressed by equation \eqref{eq:GVSinclusions} extend at the level of cochain complexes to
	\begin{displaymath}
	\begin{tikzcd}
		(L, [\cdot]_1) \ar[r,equal,"\sim"] & (\cA,\pi_1) \ar[r,hook] & (\cV,\mu_1)
		~.
	\end{tikzcd}		
	\end{displaymath}
	\noindent
	(ii)
	The $L_\infty$-algebra structure $\{\mu_k\}$ on $\cV$ restricts to a $L_\infty$-subalgebra structure on $\cA$, namely
	\begin{displaymath}
		(h)=(h,0,\dots)~:~ (\cA,\{\mu_k \circ h^{\otimes k}\}) \to (\cV,\{\mu_k\})
	\end{displaymath}
	is a strict $L_\infty$-morphism.
\end{lemma}
\begin{proof}
	\noindent
	(i) 
	By their very definition (\cf equations \eqref{eq:Lspace}, \eqref{eq:Aspace}, \eqref{eq:VSpace}), the considered cochain complexes coincide in degrees lesser than $0$.
	Therefore, one has only to prove that the following diagram commutes in the category of ordinary vector spaces
	\begin{displaymath}
		\begin{tikzcd}[column sep= large]
 			& &\quad \mathfrak{X}(M)\oplus \Omega^{n-1}(M) 
 			\\
			\cdots \ar[r,"d"] 
 			& \Omega^{n-2}(M) \ar[ru,sloped,"\mu_1"]\ar[r,"\pi_1"]
 			\ar[rd,sloped,"\dd"']
 			& \cA^{n-1}\ar[u,hook,"h "'] 
 			\\
 			& & \Omega^{n-1}_{\text{Ham}}(M,\omega)\ar[u,"\Delta "']
		\end{tikzcd}
		~.
	\end{displaymath}	
	The commutativity of the two rightmost triangles follows immediately remembering that exact $(n-1)$-forms are Hamiltonian with trivial Hamiltonian vector field (see remark \ref{Rem:ClosedformsTrivialHamiltonian}), hence
	\begin{displaymath}
		\pi_1\ (\alpha) = \mu_1 (\alpha)  = h(\d\alpha) = \pair{0}{\d \alpha}
		\qquad \forall \alpha \in \Omega^{n-2}(M)~.
	\end{displaymath}
	
	\noindent
	(ii) 
 	One has to prove that
		 \begin{displaymath}
 			\text{Im} (\mu_k \vert_{\cA}) = 
 			\text{Im} (\mu_k \cdot h^{\otimes k}) \subset
 			\text{Im} ( h)
 			~.
		 \end{displaymath}
	Being $\cA^i\equiv \cV^i$ for all $i\leq -1$, the previous condition holds automatically for any $\mu_k$ with $k\geq 3$ since  $|\mu_k| = 2-k$.
	The argument of (i) implied that that $\mu_1=\pi_1$.
	It remains only to check the case of $\mu_2$ restricted to the degree $0$ sector. 
	Consider $e_i=\pair{x_i}{\alpha_i}\in \cA^0\subset \cV^0$, where $x_i$ is the Hamiltonian vector field pertaining to $\alpha_i$, one gets:
	\begin{displaymath}
		\begin{aligned}
			[e_1,e_2]_\omega 
			=&~
			\pair{[x_1,x_2]}{\mathcal{L}_{x_1}\alpha_2 - \mathcal{L}_{x_2}\alpha_1 - \d \langle e_1,e_2\rangle_- + \iota_{x_1}\iota_{x_2}\omega}
			=
			\\
			=&~
			\pair{[x_1,x_2]}{\d (\iota_{x_1}\alpha_2 - \iota_{x_2}\alpha_1) - \d \langle e_1,e_2\rangle_- + \iota_{x_1}\d \alpha_2 - \iota_{x_2}\d\alpha_1 + \iota_{x_1}\iota_{x_2}\omega}=
			\\
			=&~
			\pair{[x_1,x_2]}{\d \langle e_1,e_2\rangle_- + \iota_{x_2}\iota_{x_1}\omega}
		\end{aligned}
	\end{displaymath}	
	employing the Cartan formula on the second equality and the Hamilton-DeDonder-Weyl equation (see the definition of the "Hamiltonian condition" in \ref{Def:Hamiltonianform}) in the last one.
	From the above equation follows that Hamiltonian pairs in
 $(\Omega^{n-1}_{\text{Ham}}(M,\omega))^{\Delta}$
are closed under the bracket $[\cdot,\cdot]_\omega$ since, according to lemma \ref{Lem:BinBrackofHamFormsisHamiltonian}, $[x_1,x_2]$ is the Hamiltonian vector field pertaining to  $\iota_{x_2}\iota_{x_1}\omega$.	
\end{proof}
\begin{notation}
	From now on, we will simply denote as $\mu_k$ the restriction of the $k$-ary brackets of definition \ref{def:Vinalgoid}, to $\mathcal{A}$.
	Let us stress that, according to previous definitions, the grading convention reads as follows:
\begin{align*}
	|\mathfrak{X}(M)\otimes \Omega^{n-1}(M)| 
	&= |\cA^{n-1}| = |\Omega^{n-1}_{\text{Ham}}| = 0
	; 
	\\
	|\Omega^k(M)| &= k - (n-1)
	~.
\end{align*}		
Summing up,
\begin{displaymath}
	L_\infty(M,\omega) \cong (\mathcal{A},\pi) ~;\quad
	L_\infty(E^n,\omega) \supset (\mathcal{A},\mu)~.
\end{displaymath}

\end{notation}
%
\subsection{{Properties} of Vinogradov's $L_\infty[1]$-algebra}\label{sec:L1VinoProp}
Is a result stated by Ritter and Saemann \cite{Ritter2015a} that the two aforementioned $L_\infty$-algebras  $(\cA,\pi)$ and $(\cA,\mu)$ are isomorphic.
In section \ref{Section:ExtendedRogersEmbedding} we will manage to provide an explicit construction of such isomorphism up to the $4$-plectic case.
In order to do so, we need to understand the relations between the multibrackets $\{\pi_k\}$ and $\{\mu_k\}$ and we are going to express them using the \RN product introduced in eq. \eqref{eq:compsymm}. 
%

Denote  by $(\cA[1],\{\bmu_k\})$ the $L_{\infty}[1]$-algebra corresponding to  Lie infinity algebra on $\cA$ obtained restricting Vinogradov's Lie-$n$ algebra $ L_{\infty}(E^n,\omega)$. 
We establish some relationship satisfied by the multibrackets $\bmu_k$. 
\begin{proposition}[{Binary bracket}]\label{Prop:mu2}
	\begin{displaymath}
		{\bmu_2} = \bpi_2 - [\pairing,\bpi_1]_{\cs}	
	\end{displaymath}
\end{proposition}
\begin{proof}
 First, we introduce the following {auxiliary} binary bracket on $\cV$:
\begin{displaymath}
\begin{aligned}
	\eta_2(e_1,e_2) 
	=& 
	\pair{[X_1,X_2]}{\iota_{X_1}\dd \alpha_2 - \iota_{X_2}\dd \alpha_1 + \iota_{X_1}\iota_{X_2}\omega}
	\\
	\eta_2(e,f) 
	=& 
	\eta_2(f,e)= 0
	~.
\end{aligned}
\end{displaymath}
 Recalling the natural extension of the pairing operator (see remark \ref{rem:pair-}), the binary bracket in definition \ref{def:vinolinfty} can then be written as
$$ 	\mu_2 = \eta_2 + \mu_1 \ca \pairing_- - \pairing_- \ca \mu_1 ~.$$
 The latter equality can be checked by inspection on arbitrary elements on $\cV$:
 
 \begin{displaymath}
 \mathclap{
 \begin{aligned}
	\mu_2&(f_1\oplus e_1  ,f_2\oplus e_2) = 
	\mu_2(e_1,e_2) + \mu_2(f_1,e_2)-\mu_2(e_1,f_2) + \cancel{\mu_2(f_1,f_2)} =\\
	&= \frac{1}{2} d \langle	 e_1,e_2\rangle_- + \eta_2(e_1,e_2)-\frac{1}{2}\mathcal{L}_{X_2}f_1 + \frac{1}{2}\mathcal{L}_{X_1}f_2 = 
	\\
	&= \frac{1}{2}\big[
		\dd (\iota_{X_1}\alpha_{2} - \iota_{X_2}\alpha_1	)	
		- (\dd \iota_{X_2}f_1 + \iota_{X_2} \dd f_1)
		+ (\dd \iota_{X_1}f_2 + \iota_{X_1} \dd f_2)	\big]
		+  \eta_2(e_1,e_2) =\\
	&= \frac{1}{2}\big[
		\dd \iota_{X_1}(\alpha_{2} + f_2)
		-\dd \iota_{X_2}(\alpha_{1} + f_1)
		+ \iota_{X_1} \dd f_2 -\iota_{X_2}\dd f_1		
		\big]
		+ \eta_2(f_1\oplus e_1 , f_2\oplus e_2) =\\
	&=\big[
		\mu_1 \ca \pairing_- - \pairing_- \ca \mu_1 + \eta_2
	\big](f_1\oplus e_1 , f_2\oplus e_2)
 \end{aligned}
 }
 \end{displaymath}
Restricting to $\cA{\subset\cV}$, 
one finds that
%
 $\eta_2  = \pi_2$, since {$\eta_2(e_1,e_2) 	= \pair{[X_1,X_2]}{\iota_{X_2}\iota_{X_1}\omega} 
	= \pi_2 (e_1,e_2)$.}
 Hence {on $\cA$ we have}
\begin{equation}\label{Eq:mu2}
	\mu_2 = \pi_2 + \mu_1 \ca \pairing_- - \pairing_- \ca \mu_1 
\end{equation} 	
and thus on $\cA[1]$:
\begin{equation}\label{Eq:mu2-shifted}
 {\bmu_2 = \bpi_2 + \bmu_1 \cs\pairing - \pairing \cs \bmu_1 ~.}
\end{equation} 
\end{proof}	

\begin{corollary}\label{cor:pairmu2} 
	The Vinogradov binary bracket commutes with the pairing:
 \begin{displaymath}
 	[\pairing,\bmu_2]_{\cs} = 0
 \end{displaymath}
\end{corollary}
\begin{proof}
	Computing the commutator of equation \eqref{Eq:mu2-shifted} with the pairing yields
 \begin{displaymath}
 	\begin{aligned}
 	[\pairing,\bmu_2]_{\cs} 
 	=&
 	[\pairing, \bpi_2 - [\pairing,\bpi_1]_{\cs}]_{\cs}
 	=
 	\\
 	=&
 	[\pairing,\bpi_2]_{\cs} - [\pairing,[\pairing,\bpi_1]_{\cs}]_{\cs}
 	= 0
 	~,
 	\end{aligned}
 \end{displaymath}
 where the last equality is given by equation \eqref{eq:commPairingMu3}.
\end{proof}
	
\begin{proposition}[{Ternary bracket}]\label{Prop:mu3}
	\begin{displaymath}
		{\bmu_3} =
		\bpi_3
		-\frac{1}{2} [\pairing, \bpi_2]_{\cs}
		- \frac{1}{6} [\pairing^{\cs 2},\bpi_1]_{\cs} 
	\end{displaymath}
\end{proposition}
\begin{proof}
 Employing the definition of the Richardson-Nijenhuis product, one can express the ternary bracket in definition \ref{def:vinolinfty} as
\begin{equation}\label{Eq:mu3}
	\mu_3 = -\dfrac{1}{3} \pairing_+\ca \mu_2 ~.
\end{equation}
 The explicit definition of $\ca$, see equation \eqref{Eq:RNProducts-explicit}, ensures that multiplying on the left by a binary bracket, not necessarily skew-symmetric, is a well-defined operation valued in graded skew-symmetric multilinear maps. 
 More precisely, equation \eqref{Eq:mu3} can be deduced by inspection on homogeneous elements. 
 When evaluated on degree $0$ elements, $\mu_3$ reads as:
	$$
		\mu_3(e_1,e_2,e_3) = -T_\omega(e_1,e_2,e_3) = -\frac{1}{3} \langle[e_1,e_2]_\omega,e_3 \rangle_+ + \cyc
	~,
	$$
	{while in other degrees, for any $f$ such that $|f|\neq 0$, it reads}
	\begin{align*}
		\mu_3(f,e_1,e_2) &=
		-\frac{1}{6}\left[
			\iota_{X_1}(\frac{\mathcal{L}_{X_2}}{2}f)
		- \iota_{X_2}(\frac{\mathcal{L}_{X_1}}{2}f)
		+ \iota_{[X_1,X_2]}f
		\right] =
		\\
		&=
		-\frac{1}{6}\left[
			\iota_{X_1}\mu_2(e_2,f)	- \iota_{X_2}\mu_2(e_1,f)
		+ \iota_{[X_1,X_2]}f
		\right] =
		\\
		&=
		-\frac{1}{3}\left[
			\pairing_+ \circ ( \mu_2 \otimes \mathbb{1})
			\right]
			\big(
			(f,e_1,e_2)-(f,e_2,e_1)+(e_1,e_2,f)
			\big) =
		\\
		&=
		-\frac{1}{3}\left[
			\pairing_+ \circ ( \mu_2 \otimes \mathbb{1})
			\right]
			(f,e_1,e_2)
			+ \cyc ~.
	\end{align*} 
	Equation \eqref{Eq:mu3} can be further expressed as:
	\begin{align*}
		\mu_3
		\equal{Eq: \eqref{Eq:mu3}}&
		-\dfrac{1}{3} \pairing_+ \ca \mu_2  
		=
		\\
		\equal{Eq: \eqref{Eq:mu2}}&
		-\dfrac{1}{3} \pairing_+ \ca \Big( \pi_2 + \pi_1\ca \pairing_- 
		- \pairing_- \ca \pi_1 \Big)  
		=
		\\
		\equal{\phantom{Eq: \eqref{Eq:mu2}}}&
		-\frac{1}{3} \pairing_+ \ca \pi_2 
		- \frac{1}{3}\Big( \pairing_+\ca\pi_1\ca \pairing_- 
		- \pairing_+\ca\pairing_- \ca \pi_1 \Big) 
		=		
		\\
		\equal{Eq: \eqref{Eq:pi3}}&
		\pi_3
		-\frac{1}{3} \pairing_- \ca \pi_2 + \frac{1}{3}
		\Big( 
			\pairing_-\ca\pi_1\ca \pairing_- 
			- \pairing_-\ca\pairing_- \ca \pi_1 
		\Big)
	\end{align*} 
   where in the last equation we used that 
   $\pairing_-\ca \eta = -\pairing_+ \ca \eta$ for any multilinear map $\eta$ in degree non-zero, and that
   \begin{equation}\label{Eq:pi3}
	\pairing_+ \ca \pi_2 = \pairing_- \ca \pi_2 - 3\, \pi_3
	~.
	\end{equation}
	Equation \eqref{Eq:pi3} {can be checked by probing with} elements $e_i=f_i + \pair{X_i}{\alpha_i}$, \ie
		\begin{displaymath}
			\begin{aligned}
			\Big[\big(\pairing_+ - \pairing_-\big) \ca \pi_2\Big] (e_1,e_2,e_3)
			=&~
			\iota_{X_3} \pi_2(e_1,e_2) + \cyc 
			~=
			\\
			=&~
			\omega(X_1,X_2,X_3) + \cyc 
			~=
			\\
			=&~
			3 \omega(X_1,X_2,X_3)
			~=
			\\
			=&~
			- 3 \pi_3 (e_1,e_2,e_3)
			~.
			\end{aligned}
		\end{displaymath}
	The claim {of the proposition} follows after applying the d\'ecalage:
	\begin{displaymath}
	 \begin{aligned}
		{\bmu_3}
		\equal{\phantom{Eq: \eqref{Eq:mu2}}}&
		\bpi_3
		-\frac{1}{3} \pairing \cs \bpi_2 + \frac{1}{3} \pairing\cs\bpi_1\cs \pairing
		-\frac{1}{3} \pairing\cs\pairing \cs \bpi_1 
		=
		\\
		\equal{Eq: \eqref{Eq:pairing-pi2}}&
		\bpi_3
		-\frac{1}{3} [\pairing, \bpi_2] 
		+ \frac{1}{6}\Big( 
		-[\pairing,\bpi_2] + \bpi_1\cs\pairing^{\cs 2}
		+\pairing^{\cs 2}\cs \bpi_1
		\Big)~+
		\\
		&
		-\frac{1}{3} \pairing^{\cs 2} \cs \bpi_1
		=
		\\
		\equal{\phantom{Eq: \eqref{Eq:mu2}}}&		
		\bpi_3
		-\frac{1}{2} [\pairing, \bpi_2] 
		+ \frac{1}{6}\Big( 
		 \bpi_1\cs\pairing^{\cs 2}
		-\pairing^{\cs 2}\cs \bpi_1
		\Big)
		=
		\\		
		\equal{\phantom{Eq: \eqref{Eq:mu2}}}&		
		\bpi_3
		-\frac{1}{2} [\pairing, \bpi_2] 
		- \frac{1}{6} [\pairing^{\cs 2},\bpi_1] 
		~.
	 \end{aligned}
	\end{displaymath}
\end{proof}
Observe that there is a common pattern for constructing $\pi_k$ and $\mu_k$ (when $k\geq 4$).
Both of them are realized via multiple insertions $\iota_{x_k},\iota_{x_{k-1}},\dots$ inside a fixed differential form $B\in \Omega(M)$.
In the first case, 
$B$ coincides with the $n$-plectic form $\omega$, in the second $B=\mu_3(\xi_1,\xi_2,\xi_3)$ where $x_i$ are generic elements of $\mathcal{A}$.
	In lemma \ref{lemma:InsertionsAsPairing} we stated how similar constructions can be phrased in terms of the pairing operator $\pairing$.
	We finish this section by showing the analogue of corollary \ref{cor:higherpi} in the case of the Vinogradov $L_\infty$-algebra\footnote{This will not be used in this chapter but we believe it will be useful in extending theorem \ref{thm:iso} to arbitrary values of $n$.}.	
\begin{lemma}\label{lem:mun}
	\begin{displaymath}
		\bmu_{{n}} = 			
		3~
		\left(
			\frac{2^{n-1}}{(n-1)!}B_{n-1}	
		\right)
			~\cdot~ \pairing^{\cs ({n}-3)} \cs \bmu_3
	\end{displaymath}
\end{lemma}	
\begin{proof}
	According to equation 	\eqref{eq:VinoMultibrakAllaZambon_1}, the explicit value of $\mu_n(\varv_1,\dots,\varv_n)$  is a sum of two terms which can be rewritten employing the anchor operator $\rho$, such that $\rho(\varv_i) = X_i$ for any  $\varv_i=f_i\oplus e_i \in \cA$. {We consider the two terms separately.}
	
	i) The first summand reads as
	\begin{displaymath}
		\begin{aligned}
			\sum_{i=1}^n (-)^{i-1}\mu_n & \left( \varv_i,\rho(\varv_1),\dots,\widehat{\rho(\varv_i)},\dots,\rho(n)\right)
			=
			\\
			=&~
			\left(
				\mu_n \circ 
				\big(\mathbb{1}\otimes \rho^{\otimes(n-1)}\big) \circ 
				P_{1,n-1}
			\right)~(\varv_1,\dots,\varv_n)				
		\end{aligned}
	\end{displaymath}	
	noticing  that $\sigma=(i,1,\dots,\hat{i},\dots,n)\in \ush{1,n-1}$ and $|\sigma|=(-)^{i+1}$.
	Explicitly, see equation \eqref{eq:VinoMultibrakAllaZambon_2}, the term   {on the r.h.s. between the large brackets} can be given by contraction with several vector fields:
	\begin{displaymath}
		\mathclap{
		\begin{aligned}
			\mu_n &\circ\big(\mathbb{1}\otimes \rho^{\otimes(n-1)}\big)~
			(\varv_0,\varv_1,\dots,\varv_{n-1})=
			\\
			=&~			
			\mu_n( \varv_0,X_1\dots,X_{n-1} )
			\quad=
			\\
			=&~
			c_n			
			\mkern-20mu\sum_{1\leq i < j \leq n-1}\mkern-20mu 
			(-)^{i+j+1}
			\iota_{\rho(n-1)}\dots\widehat{\iota_{\rho(\varv_j)}}\dots
			\widehat{\iota_{\rho(\varv_i)}}\iota_{\rho(\varv_1)}
			~[\xi_0,\rho(\varv_i),\rho(\varv_j)]_3 =
			\\
			=&~
			c_n
			\left(
				-\varsigma(n-3)\frac{2^{n-3}}{(n-3)!}
			\right)	\cdot	
			\mkern-50mu \sum_{\qquad~1\leq i < j \leq n-1} \mkern-50mu 
			(-)^{i+j+1}
			\pairing_-^{\ca (n-3)}
			\circ
			\left([\cdot,\cdot,\cdot]_3\otimes\mathbb{1}_{n-3}\right)
			\circ 
			\left(\mathbb{1}\otimes \rho^{\otimes n-1}\right)\\
			&~~~
			(\varv_0,\varv_i,\varv_j,\varv_1,\dots,\hat{\varv_i},\dots,\hat{\varv_j}\dots,\varv_{n-1}) =
			\\
			=&~
			3~d_n
			~
			\pairing_-^{\ca (n-3)}
			\circ
			([\cdot,\cdot,\cdot]_3\otimes\mathbb{1}_{n-3})
			\circ 
			(\mathbb{1}\otimes \rho^{\otimes n-1})
			\circ 
			(\mathbb{1}\otimes P_{2,n-3})
			~(\varv_0,\varv_1,
			\dots,\varv_{n-1}) ~.
		\end{aligned}
		}
	\end{displaymath}	
Here $c_n$ denotes the coefficient defined in equation \eqref{eq:UglyCoefficient} and, in the last equality, we noticed that $(-)^{i+j+1}=|\sigma|$ with $\sigma=(i,j,1,\dots,\hat{i},\dots,\hat{j},\dots,n-1)\in \ush{2,n-3}$.
	
	Further, $d_n$ is given by
	\begin{displaymath}
		\begin{aligned}
		d_n =&~		\frac{c_n}{3} 		
		\left(
				-\varsigma(n-3)\frac{2^{n-3}}{(n-3)!}
		\right)
		\\	
		=&~ 
		\left(-\varsigma(n-3)(-)^{\frac{n+1}{2}}\right)		
		\left(\frac{4~B_{n-1}}{(n-1)(n-2)} \right)
		\left(\frac{2^{n-3}}{(n-3)!} \right)
		=
		\\
		=&~
		\frac{2^{n-1}}{(n-1)!}B_{n-1}
		~.
		\end{aligned}
	\end{displaymath}
	Therefore
	\begin{align*}
			[\cdot,\dots,\cdot]_n &\circ \big(\mathbb{1}\otimes \rho^{\otimes(n-1)}\big)\circ P_{1,n-1}
			=
			\\
			=&
			3~d_n~		
			\pairing_-^{\ca(n-3)}
			\circ \left( [\cdot,\cdot,\cdot]_3\otimes\mathbb{1}_{n-3} \right)
			\circ \left(\mathbb{1}\otimes \rho^{\otimes n-1} \right)
			\circ \left(\mathbb{1}\otimes P_{2,n-3}\right)
			\circ P_{1,n-1}=
			\\
			=&
			3~d_n~		
			\pairing_-^{\ca(n-3)}
			\circ \left( [\cdot,\cdot,\cdot]_3\otimes\mathbb{1}_{n-3} \right)
			\circ \left(\mathbb{1}\otimes \rho^{\otimes n-1} \right)
			\circ P_{1,2,n-3}=
			\\
			=&
			3~d_n~	
			\pairing_-^{\ca(n-3)}
			\circ 
			\left(			
			\left(
				[\cdot,\cdot,\cdot]_3
				\circ \mathbb{1}\otimes \rho^{\otimes 2} 
				\circ P_{1,2}
			\right)\otimes \rho^{\otimes (n-3)}\right)
			\circ P_{3,n-3}=
			\\
			=&
			3~d_n~	
			\pairing_-^{\ca(n-3)}
			\circ \left([\cdot,\cdot,\cdot]_3
			\otimes \rho^{\otimes (n-3)}\right)
			\circ P_{3,n-3}=
			\\
			=&
			3~d_n~		
			\left(
			(-)^{n-3}		
			\pairing_-^{\ca(n-3)}
			\circ \left(
				[\cdot,\cdot,\cdot]_3
				\otimes \mathbb{1}_{n-3}\right)
				\circ P_{3,n-3}
			\right)
			=
			\\			
			=&
			3~d_n~	
			\left(
			\pairing_-^{\ca(n-3)}
			\ca
			[\cdot,\cdot,\cdot]_3
			\right)
			~.
		\end{align*}
	The last three equalities follow respectively from the observations that
	 $\mu_3\circ \big(\mathbb{1}\otimes \rho^{\otimes 2}\big) \circ P_{1,2} = \mu_3$ (see the definition of the ternary bracket), 
	that $\pairing_-^{n-3}\circ \big(\alpha \otimes \rho^{\otimes (n-3)}\big) = \pairing_-^{n-3}\circ (\alpha \otimes \mathbb{1}_{n-2})$ for any element $\alpha\in \cV$ such that $\rho(\alpha) = 0$, 
	and that $(-)^{n-3}=1$ for any $n\geq 3$ odd. 

ii) The second summand, {see  equation 	\eqref{eq:VinoMultibrakAllaZambon_1}}, {can be re-expressed via lemma \ref{lemma:InsertionsAsPairing} to give:}
\begin{displaymath}
	\begin{aligned}
		\Big(
			&(-)^{\frac{n+1}{2}}\cdot n \cdot B_{n-1} 
		\Big)
		~\iota_{\rho(\varv_n)}\dots\iota_{\rho(\varv_1)} \omega
		=
		\\
		=&
		-\varsigma(n)
		\left(
			-\frac{2^{n-3}}{n!}3!
		\right)
		\left(
			(-)^{\frac{n+1}{2}}\cdot n \cdot B_{n-1} 
		\right)
		\pairing_-^{\ca (n-3)} \ca \pi_3
		~ (\varv_1,\dots,\varv_n)
		=
		\\
		=&
		\left(\varsigma(n) (-)^{\frac{n+1}{2}}\right)
		\cdot
		\frac{3}{2}
		\cdot
		\left(
			\frac{2^{n-1}}{(n-1)!} B_{n-1}
		\right)
		~
		\pairing_-^{\ca (n-3)} \ca \pi_3
		~ (\varv_1,\dots,\varv_n)
		=
		\\
		=&
		-\frac{3}{2} d_n
		~
		\pairing_-^{\ca (n-3)} \ca \pi_3
		~ (\varv_1,\dots,\varv_n)
	\end{aligned}
	\end{displaymath}
	Summing up, one can conclude that
	\begin{displaymath}
		\begin{aligned}
			\mu_n =&
			3 ~d_n ~
			\pairing_-^{\ca (n-3)} \ca \left([\cdot,\cdot,\cdot]_3 -\frac{1}{2}\pi_3 \right)
			=
			\\
			=&
			3 ~d_n ~
			\pairing_-^{\ca (n-3)} \ca \mu_3
			~
		\end{aligned}	
	\end{displaymath}		
	{using proposition \ref{Prop:mu3}.}
	The claim follows after d\'ecalage.
\end{proof}
The upshot of the previous lemmas is that the two subalgebras of $M^{\sym}(\cA[1])$ respectively generated by $\{\bpi_k\}$ and $\{\bmu_k\}$
are actually generated by the same subset of only $4$ generators  $\{\bpi_1,\bpi_2,\bpi_3,\pairing \}$.
This will be crucial when trying to ascertain the existence of a $L_\infty$-morphism between $(\cA,\{\bpi_k\})\to (\cA,\{\bmu_k\})$.

\section{Extending Rogers' embedding}\label{Section:ExtendedRogersEmbedding}
Let $(M,\omega)$ be a $n$-plectic manifold. 
In this section we provide an explicit $L_{\infty}$-embedding from the $L_{\infty}$-algebra of observables on $(M,\omega)$ into the $L_{\infty}$-algebra associated to the Vinogradov algebroid (higher Courant algebroid) $E^n=TM \oplus \Lambda^{n-1} T^\ast M$  with bracket twisted by $\omega$; see theorem \ref{thm:iso} and corollary \ref{cor:Psi}.
\\
This construction can be seen as the higher analogue of mapping \eqref{eq:chris}. The case with $n=2$ has been originally carried out by Rogers in \cite{Rogers2013}.

\subsection{An $L_{\infty}$-isomorphism}
Recall that we can transfer the multibrackets of Rogers' $L_{\infty}$-algebra  $L_{\infty}(M,\omega)$ to $\cA$, by using
the isomorphism $\Omega^{n-1}_{\text{Ham}}(M,\omega) \cong (\Omega^{n-1}_{\text{Ham}}(M,\omega))^{\Delta}\cong \cA^0$
in degree zero and the identity in negative degrees.
This way we obtain an $L_{\infty}$-algebra structure on $\cA$, which we denoted by $(\cA,\{\pi_k\})$.
The latter differs from the $L_\infty$-algebra $(\cA,\{\mu_k\})$ we associated after Def. \ref{def:vinolinfty} to the $\omega$-twisted Vinogradov algebroid, but the underlying chain complex (\ie the unary brackets) are the same. 

We  show that these two $L_{\infty}$-algebra structures are $L_{\infty}$-isomorphic, by a morphism whose first component is $\Id_{\cA}$.
We will prove the following theorem in \S \ref{subsec:isoproof}.

\begin{theorem}\label{thm:iso}
For $n\le 4$ there is an $L_{\infty}$-isomorphism  ${\Phi}\colon (\cA,\{\pi_k\})\to(\cA,\{\mu_k\})$. Its non-vanishing components are   given by  
{
 \begin{align*}
{\Phi}_1&=\Id_{\cA}\\ 
{\Phi}_2&=-\pairing_-\\
{\Phi}_3&=\frac{1}{3}\pairing_-\ca \pairing_- 
 \\
{\Phi}_4&=0
\end{align*}
where $\ca$ is the skew-symmetric \RN product defined in equation \eqref{Eq:RNProducts-explicit}.
}
\end{theorem} 

\begin{remark}\label{rem:Phientries}
	By construction  $\Phi$ enjoys the   property  that
		$ \Phi_k (\varv_1 \dots \varv_k) =0$ unless $|\varv_i|=0$ for  at least $k-1$ entries.
\end{remark}

Recall the sequence of graded vector space morphisms 
$$L \cong \cA \hookrightarrow \cV~, $$
given by the identity on $\Omega^{\le n-2}(M)$ in negative  degrees, and by
$\Omega^{n-1}_{\text{Ham}}(M,\omega)\cong (\Omega^{n-1}_{\text{Ham}}(M,\omega))^{\Delta}\subset \Gamma(E^n)$ in degree zero.
As an immediate corollary of theorem \ref{thm:iso} we obtain the following embedding, which for $n=2$ is due to Rogers \cite[Theorem 7.1]{Rogers2013}.
\begin{corollary}\label{cor:Psi}
For $n\le 4$ there is an $L_{\infty}$-embedding ${\Psi}: L_{\infty}(M,\omega)\hookrightarrow L_\infty(E^n,\omega)$ of  Rogers' $L_{\infty}$-algebra into
into Vinogradov's.
	The embedding consists of only three non-trivial components given, for any $\varv_i=f_i\oplus \alpha_i \in L_{\infty}(M,\omega)$ with $f_i \in \bigoplus_{k=0}^{n-2}\Omega^k(M)$ and $\alpha \in \Omega^{n-1}_{ham}(M,\omega)$,
by the following equations:
	\begin{displaymath}
		\begin{aligned}
			{\Psi}_1(\varv) 
			=& 
			f\oplus \pair{v_\alpha}{\alpha}
			\\
			{\Psi}_2(\varv_1,\varv_2)
			=&
			-\dfrac{1}{2}\left(
				\iota_{v_{\alpha_1}}(f_2\oplus \alpha_2) - \iota_{v_{\alpha_2}}(f_1\oplus \alpha_1)
			\right)
			\\
			{\Psi}_3(\varv_1,\varv_2,\varv_3)
			=&
			\dfrac{1}{6}\iota_{v_{\alpha_1}}\iota_{v_{\alpha_2}}(f_3\oplus \alpha_3) + \cyc
		\end{aligned}
		~.
	\end{displaymath}
\end{corollary}

\begin{remark}[Comparisons with the Ritter-Saemann Theorem]
	We emphasize that a result similar to corollary \ref{cor:Psi} was already stated, slightly less explicitly, by Ritter and Saemann in \cite[Thm. 4.10]{Ritter2015a}.
	Basically, their idea was to deform $\mu_k$ into $\pi_k$ via a sequence of approximating $L_\infty$-morphisms.
	\\
	Let us briefly paraphrase their results in our notation. Given an $n$-plectic manifold $(M,\omega)$, consider the graded vector space $\cA$ as defined in equation \eqref{eq:Aspace} above.
	They claimed the existence of $(n$-$1)$ $L_\infty$-morphisms
	\begin{displaymath}
		\Psi^{(m)}= (id,0,\dots,\varphi^{(m)}_{m+1},0,\dots): (\mathcal{A},\mu_k^{(m)}) \to (\mathcal{A}, \mu_k^{(m+1)}) ~,
	\end{displaymath}
	where $\mu_k^{(\ell)}$ is a $L_\infty$-structure on $\mathcal{A}$ that agrees to Roger's up to order $\ell$, i.e
 $$\mu_k^{(\ell)} = \pi_k \qquad \forall k\leq \ell ~,$$
and in particular $\mu_k^{(1)}$ coincides with the restriction of Vinogradov's $L_\infty$-algebra on $\mathcal{A}$, 
	such that $\Psi = \Psi^{(n-1)}\circ\dots\circ\Psi^{(1)}$ is a $L_\infty$-isomorphism between $(\mathcal{A},\mu_k^{(1)})$ and $(\mathcal{A},\pi_k)$.
	This construction can be alternatively subsumed by the following diagram:
	\begin{displaymath}
		\begin{tikzcd}[column sep = large]
			L_{\infty}(M,\omega) \ar[d,equal,sloped,"\sim "]& & L_{\infty}(E^n,\omega)
			\\
			(\mathcal{A},\pi_k) & & (\mathcal{A},\mu_k^{(1)}) \ar[u,hook,"h"']
			\ar[d,"{\Psi^{(1)}=(id,\varphi_2^{(1)},0,\dots)}"]
			\\	
			& & (\mathcal{A},\mu_k^{(2)})
			\ar[d,"{\Psi^{(2)}=(id,0,\varphi_3^{(2)},0,\dots)}"]
			\\
			& & \vdots 
			\ar[d,"{\Psi^{(n-1)}=(id,0,\dots,\varphi_n^{(n-1)})}"]
			\\
			& & (\mathcal{A},\mu_k^{(n)}) \ar[d,"\Psi^{(n)}=id"]
			\\
			& & (\mathcal{A},\mu_{k}^{(n+1)})	\ar[lluuuu,equal]
		\end{tikzcd}
	\end{displaymath}
	Notice that in \cite{Ritter2015a} also the single components $\Psi^{(k)}$, for any $k$, are claimed to be obtained from subsequent approximations. Namely $\mu_{k+1}^{(k)}$ is deformed into $\mu_{k+1}^{(k+1)}=\pi_{k+1}$ degree by degree, with respect to the grading of $\cA$.
	\\
	The interest in finding an explicit expression for each of the above components $\varphi_n$ has been our original motivation for performing the computations involved in the results of sections \ref{sec:L1RogersProp} and \ref{sec:L1VinoProp}.
\end{remark}

\subsection{The proof of Theorem \ref{thm:iso}}\label{subsec:isoproof} 
{For the proof, it is convenient to work with $L_{\infty}[1]$-algebras, by applying the d\'ecalage isomorphism \eqref{eq:deca}.} 
\\
In the following, for any given homogeneous linear map $m\in\underline{\Hom}(S^{\ge 1}(V),V)$, we denote by $C_m$ the corresponding coderivation of $S^{\ge 1}(V)$ given by the lift, \ie $C_m=\widetilde{L}_{\sym}(m)$.

Denote  by $(\cA[1],\{\bpi_k\})$ the $L_{\infty}[1]$-algebra corresponding to $(\cA,\{\pi_k\})$.  
Notice that applying the d\'ecalage isomorphism  to obtain $\bpi$ does not introduce any extra signs, since the higher multibrackets in Rogers' $L_{\infty}$-algebra vanish unless all entries have degree $0$.
\\
Similarly, denote  by $(\cA[1],\{\bmu_k\})$ the $L_{\infty}[1]$-algebra corresponding to $(\cA,\{\mu_k\})$,
and write $\bmu \colon  S^{\ge 1}(\cA[1])\to \cA[1]$ for the map with components $\bmu_k$ ($k\ge 1)$.

{We want to construct an $L_{\infty}[1]$-isomorphism from $(\cA[1],\{\bpi_k\})$ to $(\cA[1],\{\bmu_k\})$.}
The idea is to apply the key remark \ref{rem:exp} in section \ref{SubSection:studycase}, \ie to construct the sought isomorphism as the exponential of a degree $0$ coderivation.
Denote by $Q_{\bpi}$ the codifferential on $S^{\ge 1}(\cA[1])$ corresponding to $(\cA[1],\bpi)$. 
For any   degree $0$ linear map $p\colon S^{\ge 1}(\cA[1])\to \cA[1]$ {such that $e^{C_p}$ converges,} we know that
\begin{itemize}
\item $ e^{C_p}\circ Q_{\bpi}\circ e^{-C_p}$ is a new codifferential, which corresponds to a new 
 $L_{\infty}[1]$-algebra structure $\bpi'$ on $\cA[1]$;
 \item $e^{C_p}$  corresponds to
an $L_{\infty}[1]$-isomorphism $f$ from $(\cA[1],\bpi)$ to $(\cA[1],\bpi')$.
\end{itemize}
The explicit formulae for $\bpi'$ and $f$ were given in Remark \ref{rem:exp}.
We will show that $p$ can be chosen in such a way that $\bpi'=\bmu$, at least when $n\le 4$.

In the following we will  employ the straightforward extension of the pairing operator  $\pairing_-$ of \S \ref{subsec:Vinogradov} from $\mathfrak{X}(M)\oplus \Omega^{n-1}$(M) to the whole {graded vector space $\cV$}.
(In the terms of remark \ref{rem:pair-}, one can equivalently say that we are considering the restriction of $\pairing_-$ from $\widetilde{\cV}$ to $\cV$.)
\\
The following proposition implies immediately the first part of  theorem \ref{thm:iso}. 
The proof relies on our previous computations about $\bpi$ and $\bmu$ in terms of the \RN products (see  sections \ref{sec:L1RogersProp} and \ref{sec:L1VinoProp}).
\begin{proposition}\label{prop:main}
Let $n\le 4$.
  Let $p\colon S^{\ge 1}(\cA[1])\to \cA[1]$ be the degree $0$ linear map given by\footnote{Thus $p=-\langle \;,\; \rangle-\frac{1}{6}\langle \;,\; \rangle\cs \langle \;,\; \rangle$.}
$$p=c_1\langle \;,\; \rangle+c_2\langle \;,\; \rangle\cs \langle \;,\; \rangle+c_3\langle \;,\; \rangle\cs \langle \;,\; \rangle\cs \langle \;,\; \rangle$$
for $c_1=-1, c_2=-\frac{1}{6}, c_3=0$,
where $\cs$ denotes the operation introduced in eq. \eqref{eq:compsymm}. 
Then 
\begin{equation}\label{eq:pi'}
\bpi':=
 \bpi+[p,\bpi]_{\cs}+\frac{1}{2!}[p, [p,\bpi]_{\cs}]_{\cs}+\dots
\end{equation}
agrees with $\bmu$.
\end{proposition}

\begin{remark}\label{rem:degrees}
{Although proposition \ref{prop:main} is phrased only for $n\le 4$, we make some remarks that hold for arbitrary $n$ and for any linear map $p=\sum_{i=1}^{n-1}c_i \langle \;,\; \rangle^{\cs i}$ with arbitrary string of coefficients $c_1,\dots,c_{n-1}$.}

a) For every $k\ge 1$, the component $\bpi'_k=\bpi'|_{S^k(\cA[1])}$ is a finite sum: more precisely, only the first $k$ summands on the r.h.s. of eq. \eqref{eq:pi'} contribute to it.
The reason is that $p_1=0$ by definition, and the only non-vanishing component of the maps $p_i\cs\bpi_j$ and $\bpi_j\cs p_i$ is the component defined on $S^{i+j-1}(\cA[1])$.

b) The elements of $\cA[1]$ of maximal degree are those of degree $-1$, and are those lying in $(\Omega^{n-1}_{\text{Ham}}(M,\omega))^{\Delta}[1]$. 
We now make some considerations about the vanishing of $\bpi'_k$ 
when applied to {homogeneous} elements of $\cA[1]$ of lower degree (by which we mean:
of degree that is non-maximal, \ie $\le -2$). First notice the following, {where $p_k=p|_{S^k(\cA[1])}$:}
\begin{itemize}
\item $p_1=0~.$
\item for $k\ge 2$: $p_k$ applied to  {homogeneous} elements of $\cA[1]$ might be non-vanishing only when \emph{all but possibly one} elements are in degree $-1$. The result is a lower degree element.
\item for $m\ge 2$: $\bpi_m$ applied to  {homogeneous} elements of $\cA[1]$ might be non-vanishing only when \emph{all} elements are in degree $-1$. The result is a lower degree element, except in the case of $\bpi_2$.
\end{itemize}

As a consequence, for all $k\ge 2 $ we have:
\begin{itemize}
\item[i)]  $[p_k,\bpi_2]_{\cs}$ applied to  {homogeneous} elements of $\cA[1]$ might be non-vanishing only when \emph{all but possibly one} elements are in degree $-1$.\\
The same holds for iterated brackets $[p_{k_1},\dots,[p_{k_l},\bpi_2]_{\cs}]_{\cs}$.
\item[ii)]  for $m\ge 3$, $[p_k,\bpi_m]_{\cs}$ applied to  {homogeneous} elements of $\cA[1]$ might be non-vanishing only when \emph{all} elements are in degree $-1$.\\
The same holds for iterated brackets $[p_{k_1},\dots,[p_{k_l},\bpi_m]_{\cs}]_{\cs}$.
\end{itemize}

We consider separately the case of iterated brackets involving $\bpi_1$.
Notice that $\bpi_1(\xi)$ is a degree $-1$ element only when  $\xi\in \cA[1]$ has degree $-2$. 
In that case $\bpi_1(\xi)= \pair{0}{\dd \xi}\in(\Omega^{n-1}_{\text{Ham}}(M,\omega))^{\Delta}[1]$, \ie the vector field component vanishes. 
Further for $k\ge 2$, refining a statement above, results that $p_k$ applied to elements of $\cA[1]$ may be non-vanishing only if \emph{all but possibly one} entries are are in degree $-1$ with \emph{non-vanishing vector field component}. 
Consequently we have:
\begin{itemize}
\item[iii)]  $[p_k,\bpi_1]_{\cs}$ applied to  {homogeneous} elements of $\cA[1]$ might be non-vanishing only when \emph{all but possibly one} elements are in degree $-1$. \\
The same holds for iterated brackets $[p_{k_1},\dots,[p_{k_l},\bpi_1]_{\cs}]_{\cs}$. 
\end{itemize} 
{The conclusion we  draw from i), ii), iii) is that the  $L_{\infty}[1]$-algebra structure $\bpi'$ on $\cA[1]$, defined as in eq. \eqref{eq:pi'}, has the following property:  the evaluation of multibrackets with arity $k\ge 2$ on  homogeneous elements
   might be non-vanishing only when \emph{all but possibly one} elements are of top  degree (\ie degree $-1$). Notice that the $L_{\infty}[1]$-algebra associated to the $\omega$-twisted Vinogradov algebroid has the same property, as recalled in Remark \ref{Remark:semplifica-conti}.
   }
\end{remark}

\begin{proof}[Proof of  proposition \ref{prop:main}]
We write out explicitly the first $\bpi_k$'s:
\begin{align*}
\allowdisplaybreaks
 \bpi_1'&=\bpi_1
 \\
 \bpi_2'&=\bpi_2+[p_2,\bpi_1]_{\cs}
 \\
 \bpi_3'&=\bpi_3+[p_3,\bpi_1]_{\cs}+[p_2,\bpi_2]_{\cs}+\frac{1}{2} [p_2, [p_2,\bpi_1]_{\cs}]_{\cs}
 \allowdisplaybreaks[0]
\\
\bpi_4'&=\bpi_4+[p_4,\bpi_1]_{\cs}+[p_3,\bpi_2]_{\cs}+[p_2,\bpi_3]_{\cs}\\
&+
\frac{1}{2} [p_3, [p_2,\bpi_1]_{\cs}]_{\cs}+\frac{1}{2} [p_2, [p_3,\bpi_1]_{\cs}]_{\cs}+
\frac{1}{2} [p_2, [p_2,\bpi_2]_{\cs}]_{\cs}+
\\
&+
\frac{1}{6} [p_1,[p_1,[p_1,\bpi_1]_{\cs}]_{\cs}]_{\cs}. 
\end{align*}
We now compute explicitly the right-hand sides.
According to proposition \ref{Prop:mu2}, one has that
	\begin{displaymath}
		\bpi'_2 = \bpi_2 + c_1 [\pairing,\bpi_1]_{\cs}	
	\end{displaymath}
{equals	$\bmu_2$} if and only if $c_1 = -1$.
Similarly, according to proposition \ref{Prop:TernaryCommutator}, one has
	\begin{displaymath}
	\begin{aligned}
 \bpi'_3 =&~ \bpi_3 + c_2 [\pairing^{\cs 2}, \bpi_1] + c_1 [ \pairing, \bpi_2] + \frac{c_1^2}{2}[\pairing,[\pairing,\bpi_1]]
		\\
		=&~\bpi_3 + c_2 [\pairing^{\cs 2}, \bpi_1] + c_1\left(\frac{2+c_1}{2}\right) [ \pairing, \bpi_2],
	\end{aligned} 
	\end{displaymath}
{and this equals 		$\bmu_3$} 
	if and only if $c_1=-1$ and $c_2 = -1/6$, {by proposition  \ref{Prop:mu3}}.\\
	Every even multibracket in the Vinogradov's $L_\infty$-algebra with arity greater than 2 is trivial, in particular $\bmu_4 = 0$. On the other hand, proposition \ref{Prop:QuaternaryCommutator} implies that 
	\begin{align*}
			\bpi'_4 
			=& \bpi_4 +
			\\
			&+ c_1 \pairing \cs \bpi_3 + c_2 \pairing^{\cs 2} \cs \bpi_2 + c_3 [ \pairing^{\cs3}, \bpi_1] +
			\\
			&
			+ \frac{c_1^2}{2} [ \pairing, \pairing\cs \bpi_2] 
			+ \frac{c_1 c_2}{2} \left(
				[\pairing^{\cs 2},[\pairing,\bpi_1]] + [\pairing,[\pairing^{\cs 2},\bpi_1]]
			\right)+
			\\
			&+ \frac{c_1^3}{6} [\pairing,[ \pairing , [ \pairing , \bpi_1 ] ] ]
			=
			\\
			=& 
			\left(1+ 2c_1 + \frac{3}{2}c_1^2 + \frac{{c_1^3}}{2}\right)\bpi_4 +
			\\
			&+ (c_2 + c_1 c_2)[\pairing^{\cs 2},\bpi_2]
			\\
			&+ c_3 [\pairing^{\cs 3}, \bpi_1].
		\end{align*}
Substituting the values of $c_1$ and $c_2$ just found implies that
	\begin{displaymath}
		\bpi_4' = c_3 [\pairing^{\cs 3}, \bpi_1]
	\end{displaymath}
	hence, $\bpi_4' = 0$
	when $c_3=0$.
\end{proof}

 \begin{proof}[Proof of theorem \ref{thm:iso}]
{We apply  remark \ref{rem:exp} to the $L_{\infty}[1]$-algebra $(\cA[1],\bpi)$
as indicated at the beginning of this subsection, choosing $p\colon S^{\ge 1}(\cA[1])\to \cA[1]$ as in proposition \ref{prop:main}. Notice that $p$ satisfies the condition explained in Remark \ref{rem:codermor}, hence $e^{C_p}$ is convergent.}
Thanks to proposition \ref{prop:main}, we just have to make explicit the $L_{\infty}$-isomorphism $f\colon (\cA[1],\{\bpi_k\})\to (\cA[1],\{\bmu_k\})$ provided by   remark \ref{rem:exp}.
\\
Due to the restriction $n\le 4$, we know that $\cA[1]$ is concentrated in degrees $-4,\dots,-1$. 
Since  $\pairing \colon S^{ 2}(\cA[1])\to \cA[1]$ has degree zero, we have $\pairing^{\cs 4}=0$. 
Computing $p^{\cs 2}=\pairing^{\cs 2} +\frac{1}{3}\pairing^{\cs 3}$ and 
$p^{\cs 3}=-\pairing^{\cs 3}$
we obtain
	\begin{align*}
		f &=pr_{\cA[1]}+p+\frac{1}{2!}p^{\cs 2}+\frac{1}{3!}p^{\cs 3} ~=\\
	&=pr_{\cA[1]}- \pairing +\frac{1}{3}\pairing^{\cs 2}. 
	\end{align*}
Applying the d\'ecalage isomorphism \eqref{eq:Dec}, and a fortiori of remark \ref{rem:pair-}, we obtain an  $L_{\infty}$-isomorphism  ${\Phi}\colon (\cA,\{\pi_k\})\to(\cA,\{\mu_k\})$.

\end{proof}

\begin{remark}\label{rem:arbitraryn}
We expect to be able to extend proposition \ref{prop:main} to arbitrary values of $n$.
Thanks to Remark \ref{rem:degrees} b), to do so it suffices to apply $\bpi'$ on strings of elements of $\cA[1]$ which are of maximal degree except for possibly one. 

{For the resulting $L_{\infty}$-isomorphism  ${\Phi}\colon (\cA,\{\pi_k\})\to(\cA,\{\mu_k\})$,  extending theorem \ref{thm:iso}, we expect each component  $\Phi_k$ to be given as in eq. \eqref{eq:Phik} in the next section. This expectation is suggested by the proof of theorem \ref{thm:comm}.}
\end{remark}	 

\section{Gauge transformations}\label{Sec:DiagramGaugeTransf}
Given two $n$-plectic forms $\omega$ and $\widetilde{\omega}$ on the same manifold $M$, we say that the two are \emph{gauge related} if there exists a $n$-form $B$ such that 
$$
	\widetilde{\omega} = \omega  {+}\dd B ~,
$$
\ie if they lie in the same de Rham cohomology class.
Recall that in section \ref{Section:GaugeTransformations} we already described the relationship between the \momaps pertaining to a group action that is multisymplectic with respect to two gauge related $n$-plectic form. 
{The aim of this section is to show that the $L_{\infty}$-embedding constructed in corollary \ref{cor:Psi} is compatible with gauge transformations, see theorem \ref{thm:comm}.}

\subsection{Vinogradov algebroids and gauge transformations}\label{Sec:VinoGauge}
 {Let $\omega$ and $\widetilde{\omega} = \omega + \dd B$ be gauge-related closed $(n+1)$-forms on $M$.}
The two corresponding twisted Vinogradov algebroids $(E^n, \omega)$ and $(E^n,\widetilde{\omega})$ are isomorphic. Indeed, the vector bundle isomorphism
\begin{displaymath}
	\begin{tikzcd}[column sep= small,row sep=0ex,
				/tikz/column 1/.append style={anchor=base east}]
		\tau_B \colon &
		E^n = TM \oplus \Lambda^{n-1} T^\ast M \ar[r]& E^n ~, \\
	 &\pair{X}{\alpha} \ar[r, mapsto] & 
	 \pair{X}{\alpha + \iota_X B} =
 \pair{X}{\alpha} + \pair{0}{\iota_X B}
	\end{tikzcd}		
	 ~,		
\end{displaymath}
{preserves the anchor $\rho$, the pairing  $\pairing_{-}$, and maps the bracket $[\cdot,\cdot]_{\omega}$ to $[\cdot,\cdot]_{\widetilde{\omega}}$. (Cf.  remark \ref{rem:VinoidsMorphism}).}
\\
Hence this bundle isomorphism induces a strict $L_\infty$-isomorphism at the level of the corresponding Vinogradov's $L_\infty$-algebras (\cf  \cite[Prop. 8.5]{Zambon2012}) given by
\begin{displaymath}
	(\tau_B)_m = 
	\begin{cases}
		\id_{\cV} + \iota_{\rho(\cdot)} B	
		& m=1
		\\
		0 & m\geq 2
	\end{cases}
	~.
\end{displaymath}

Notice that there is a natural diagram in the category of $L_\infty$-algebras
\begin{equation}\label{diag:3sides}
	\begin{tikzcd}
		L_{\infty}(M,\omega) \ar[r,"\Psi"] 
		&
		L_{\infty}(E^n,\omega) \ar[d,"\tau_B"]
		\\
		L_{\infty}(M,\widetilde{\omega}) \ar[r,"\Psi"] 
		&
		L_{\infty}(E^n,\widetilde{\omega})
	\end{tikzcd}
\end{equation}
where the horizontal arrows are the extensions of the Rogers's embedding 
{constructed in corollary \ref{cor:Psi}  (for $n\le 4$), which, in particular, are given by formulae that do not depend on the multisymplectic form $\omega$ or $\widetilde{\omega}$.}
\\ 
When considering two gauge-related multisymplectic manifolds, it is not possible  to define a canonical $L_\infty$-morphisms between the two corresponding
observables $L_\infty$-algebras\footnote{The two objects are already substantially different at the level the underlying vector space. An Hamiltonian form with respect to $\omega$ is not in general Hamiltonian with respect to $\tilde{\omega}$.}.
In particular there is no canonical way to close this diagram on the left to give a commutative square.
In what follows, we will look for a suitable pullback $\mathfrak{g}$ in the category of $L_\infty$-algebras:
\begin{displaymath}
	\begin{tikzcd}
		\mathfrak{g}\ar[r] \ar[d] \arrow[dr, phantom, "\scalebox{1.5}{$\lrcorner$}" , very near start, color=black]
		& L_\infty(M,\omega) \ar[d,"\Psi"]
		\\
		L_\infty(M,\tilde{\omega}) \ar[r,"\Psi"] 
		& L_\infty(E^n,\omega)\cong L_\infty(E^n,\tilde{\omega})		
	\end{tikzcd}
\end{displaymath}

	%

\subsection{Commutativity}
{Consider an infinitesimal action $v:\mathfrak{g}\to \mathfrak{X}(M)$ 
	preserving the $n$-plectic form $\omega$. 
	Suppose that $B\in \Omega^n(M)$ is also preserved, and   that $\widetilde{\omega}:=\omega + dB$ is non-degenerate.}
Further,  assume there exists a  homotopy comoment map $(f):\mathfrak{g}\to L_\infty(M,\omega)$ for $\omega$. 
Then we obtain a homotopy comoment map $(\widetilde{f})$ for $\widetilde{\omega}$, by Lemma \ref{lem:momaps} of section \ref{Section:GaugeTransformations}.

 {We will show that the diagram \eqref{diag:3sides} can be completed} 
 to a commutative pentagon, at least when $n\le 4$.
 \begin{theorem}\label{thm:comm}
{Assume that $n\le 4$.}
The following diagram of $L_{\infty}$-algebra morphisms commutes,
where $\Psi$ is the morphism introduced in corollary \ref{cor:Psi}. 
  \begin{equation}\label{eq:pentagonDiagram}
	\begin{tikzcd}
		&
		L_{\infty}(M,\omega) \ar[r,"\Psi"] 
		&
		L_{\infty}(E^n,\omega) \ar[dd,"\tau_B"]
		\\[-1em]
		\mathfrak{g}\ar[ru,"f"] \ar[dr,"\widetilde{f}"']
		\\[-1em]
		&
		L_{\infty}(M,\widetilde{\omega}) \ar[r,"\Psi"] 
		&
		L_{\infty}(E^n,\widetilde{\omega})
	\end{tikzcd}
\end{equation}
\end{theorem}
{We interpret this commutativity by saying that the twisting of the homotopy comoment moment map by $B$ is compatible with the twisting of the Vinogradov algebroid.} 



{For the proof of theorem \ref{thm:comm} we}
{make use of the strict isomorphism $L_{\infty}(M,\omega)\cong (\cA,\{\pi_k\})$ (see equations \eqref{eq:Aspace} and \eqref{eq:trianglemap}), 
given by  the identity in negative degrees, and  $\alpha\mapsto \pair{X_{\alpha}}{\alpha}$ in degree zero.
\\
The commutativity of diagram \eqref{eq:pentagonDiagram} is then equivalent to the commutativity of this diagram, where $\Phi$ is the $L_{\infty}$-morphism constructed in theorem \ref{thm:iso}.
\\
%
	\begin{equation}\label{diag:HamiltonianpentagonDiagram}
		\begin{tikzcd}[column sep=huge]
			&
			(\cA,\{\pi_k\}) \ar[r,"\Phi"] 
			&
			{L_{\infty}(E^n,\omega)}  \ar[dd,"\tau_B"]
			\\[-1em]
			\mathfrak{g}\ar[ru,"f"] \ar[dr,"\widetilde{f}"']
			\\[-1em]
			&
			(\widetilde{\cA},\{\widetilde{\pi_k}\})) \ar[r,"\Phi"] 
			&
			{L_{\infty}(E^n,\widetilde{\omega})}  
		\end{tikzcd}
	\end{equation}
	\\

\begin{remark}\label{Rem:DiagramSituation}
All the arrows involved in diagram \eqref{diag:HamiltonianpentagonDiagram} can be expressed in terms of the pairing $\pairing_-$ {and the operation $\ca$  defined in eq. \eqref{Eq:RNProducts-explicit}}:
	\begin{align}\label{eq:Phim}
		\Phi_m =&
		\begin{cases}
			\id_{\cA} & m=1
			\\
			\varphi_m ~ \pairing_-^{\ca (m-1)} & m\geq 2
		\end{cases}	
		\\[.5em]
		(\tau_B)_m =&
		\begin{cases}
			\id_{\cA}  -  2 ~\langle B, \cdot \rangle_- & m=1
			\\
			0 & m\geq 2
		\end{cases}
		\\[.5em]	
		(f)_m =&
		\begin{cases}
			f_1: \xi \mapsto \pair{v_\xi}{f_1(\xi)} \in \cA_0 & m=1
			\\
			f_m & m\geq 2
		\end{cases}
		\\[.5em]
		(\widetilde{f})_m =&
		\begin{cases}
			f_1   -  2~ \langle B, \cdot \rangle_- \circ f_1
			& m=1
			\\
			f_m  -  d_m ~\left(\pairing_-^{\ca (m-1)} \ca \langle B,\cdot \rangle_-\right) \circ f_1^{\otimes m} & m\geq 2.
		\end{cases}		
	\end{align}
{The coefficients $\varphi_m$ for $\Phi_m$ are provided by theorem \ref{thm:iso}, which holds for   $n\le 4$: they are given by
$\varphi_2=-1$, $\varphi_3=\frac{1}{3}$, and  $\varphi_m=0$ for $m\ge 4$.}
The coefficients $d_m$ are given by
$$d_m= \left( \frac{2^m}{m!}\right),$$ as follows {from lemma \ref{lem:momaps}}, lemma \ref{lemma:InsertionsAsPairing} and remark \ref{Rem:piccoloerroredisegnoneldraft} noting that
	\begin{equation}\label{eq:bm}
		\begin{aligned}
		\mathsf{b}_m(x_1,\dots x_m) &= 
		\varsigma(m+1)
		\iota_{v_{x_n}}\dots \iota_{v_{x_1}} B
		=
		\\
		&=
		- \left(\varsigma(m+1)\varsigma(m) \frac{2^m}{m!}\right)
		~\left(\pairing_-^{\ca (m)} \right) 
		\circ \left(B\otimes f_1^{\otimes m}\right)
		=
		\\
		&=
		- \dfrac{2^m}{m!} \left(\pairing_-^{\ca(m-1)}\ca \langle B, \cdot \rangle_- \right) 
		\circ f_1^{\otimes m}
		~.
	\end{aligned}
	\end{equation}
\end{remark}

\subsection{Proof of Theorem \ref{thm:comm}}
In this subsection we provide the proof of theorem \ref{thm:comm}.
The condition $n\le 4$ will be used only at the very end.

	To ascertain the strict commutativity of diagram \eqref{diag:HamiltonianpentagonDiagram} one has to make sure that 
	\begin{displaymath}
		(\tau_B \circ \Phi \circ f)_m -(\Phi\circ \widetilde{f})_m = 0 \qquad \forall m\geq 1
		~.
	\end{displaymath}
	The case $m=1$ is straightforward:
	\begin{displaymath}
		\begin{aligned}
		(\tau_B \circ \Phi \circ f)_1 -(\Phi\cdot \widetilde{f})_1 
		=&~
		(\tau_B)_1 \circ \Phi_1 \circ f_1 -\Phi_1 \circ \widetilde{f}_1
		=
		\\
		=&~
		(\tau_B)_1  \circ f_1 - f_1 + 2 \langle B, \cdot\rangle_- \circ f_1
		=
		\\
		=&~ 0
		~.
		\end{aligned}
	\end{displaymath}
	{Alternatively, one can adapt the argument given for the symplectic case $n=1$ in remark \ref{rem:symcomm}.}

{Now let $m\ge 2$.}	Since $\tau_B$ is a strict morphism and $(\tau_B)_1$ acts as the identity on any element of $\cA$ in degree different than $0$, 
	the higher cases requires to check that
	\begin{displaymath}
		{(\Phi \circ f)_m- (\Phi \circ \widetilde{f})_m = 0}.
	\end{displaymath}
Recall now the explicit expression for the composition of of two $L_{\infty}$-morphisms in the skew-symmetric framework (obtained by working out definition \ref{Def:CompositionFormula} together with remark \ref{rem:operatorSln}) :
 \begin{reminder}
	Given two $L_\infty$-morphisms $f\colon V \to V'$ and $g\colon V' \to V''$	,
	the components of their  {composition}   $g\circ f$ are
	\begin{equation}\label{eq:LinfinityMorphismsComposition}
		(g \circ f)_m = \sum_{\ell=1}^m g_\ell \circ \mathcal{S}_{\ell,m} (f)
	\end{equation}
	where the operator $\mathcal{S}_{\ell,m} (f)$ is the component $\otimes^m V \to \otimes^\ell V'$ of the lift of $f$ to a coalgebra morphism  $T(V)\to T(V')$ . 
	Explicitly,
	\begin{displaymath}
		\mathcal{S}_{\ell,m} (f) =
		\left( 
			\sum_{\substack{k_{1}+\cdots+k_{\ell}=m\\1\leq k_{1}\leq\cdots\leq k_{\ell}}}
			(-)^{\sum_{i=1}^{\ell-1}(|f_{k_i}|)(\ell-i)}
			(f_{k_1}\otimes\cdots\otimes f_{k_\ell})\circ P_{k_1,\ldots,k_\ell}^<		
		\right)~,
		\end{displaymath}
		with $P_{k_1,\ldots,k_\ell}^<$ denoting the sum over all the (odd) action of permutations $\sigma$ in $(k_1,\cdots,k_\ell)$-unshuffles on $V^{\otimes (k_1+\dots+k_\ell)}$ satisfying the extra condition
	\begin{displaymath}
		\sigma(k_1+\dots+k_{j-1}+1)<\sigma(k_1+\dots+k_{j}+1) 
		\quad \text{if}~k_{j-1}=k_j ~.
	\end{displaymath}	
\end{reminder}
{Thanks to remark \ref{rem:Phientries},} equation \eqref{eq:LinfinityMorphismsComposition}, defining the  the composition of   $L_{\infty}$-morphisms, takes the simpler form
		\begin{displaymath}
			(\Phi \circ f )_m =
			\Phi_m \circ f_1^{\otimes m} + \left[
				\sum_{\ell=1}^{m-1} \Phi_\ell \circ \left( f_1^{\otimes (\ell -1)} \otimes f_{m-\ell+1} \right) \circ P_{\ell -1, m-\ell +1}
			\right]
			~,
		\end{displaymath}
 and similarly for 		$(\Phi \circ \widetilde{f} )_m $.
	Observe now that one can compare the images of two gauge related \momaps $(f)$ and $(\widetilde{f})$ by viewing them as elements in $\cV$ in virtue of the diagram (yet to proven to be commutative):
		\begin{displaymath}
		\begin{tikzcd}[row sep=small]
			\mathfrak{g} \ar[r,"f"]\ar[rd,"\widetilde{f}"'] &
			L_{\infty}(M,\omega) \ar[r,"\cong"]&
			(\cA,\pi) \ar[r,hook] &
			L_\infty(E^n,\omega)
			\\
			&
			L_{\infty}(M,\widetilde{\omega}) \ar[r,"\cong"]&
			(\widetilde{\cA},\widetilde{\pi}) \ar[ru,hook] &
		\end{tikzcd}
	\end{displaymath}
	For instance, fixing $\xi \in \mathfrak{g}$, one could conclude that
	\begin{displaymath}
		\widetilde{f}(\xi)-f(\xi) = \pair{0}{\iota_{v_\xi}B}=\mathsf{b}_1(\xi)
	\end{displaymath}
	and read it as an element lying in the kernel of the anchor $\rho$, \ie the map that projects onto the vector fields component of any given Hamiltonian pair.
	In particular one could make sense of the following equation:
	\begin{displaymath}
		\Phi_m \circ \widetilde{f}_1^{\otimes m} -  		
  		\Phi_m \circ f_1^{\otimes m}
  		=
	 		\Phi_m \circ\left( \widetilde{f}_1^{\otimes m} - f_1^{\otimes m}\right)
	 		~.
	\end{displaymath}
	According to remark \ref{rem:Phientries},  all terms in the r.h.s. of the above equation which involve more than one occurrence of $\mathsf{b}_1$ vanish by degree reasons\footnote{$\Phi_m$ is constructed out of the pairing. The cancellation of these terms happens since they involve $m-1$ contractions with less than $m-1$ non-trivial vector fields.},
	thence 
	\begin{displaymath}
		\begin{aligned}
		\Phi_m \circ \widetilde{f}_1^{\otimes m} - \Phi_m \circ f_1^{\otimes m}
			=&~  
			\sum_{\ell=0}^{m-1}\Phi_m \circ \left(f_1^{\otimes(m-1-\ell)}\otimes \mathsf{b}_1 \otimes f_1^{\otimes \ell}\right) =
			\\
			=&~
			\Phi_m \circ \left( f_1^{\otimes (m-1)}\otimes \mathsf{b}_1 \right) \circ P_{m-1, 1}
			~.
			\end{aligned}
	\end{displaymath}
	Summing up, we conclude that
	\begin{equation}\label{eq:square}
		(\Phi \circ f)_m
		-(\Phi\circ \widetilde{f})_m = 
		- \left[
				\sum_{\ell=1}^{m} 
				\Phi_\ell \circ \left( f_1^{\otimes(\ell -1)} \otimes \mathsf{b}_{m-\ell+1} \right) \circ P_{\ell -1, m-\ell +1} 			\right]
		~.
	\end{equation}

{The following lemma allows to compute the summands on the right-hand side of equation \eqref{eq:square}, using equation \eqref{eq:Phim} to write 
$\Phi_{\ell}=\varphi_\ell \circ \pairing_-^{\ca (\ell-1)}.$
}	
	\begin{lemma}\label{lem:rhsbm}
{For all $\ell\ge 1$ we have}
	 \begin{displaymath}
	 	\pairing_-^{\ca \ell -1} \circ \left( f_1^{\otimes(\ell -1)} \otimes \mathsf{b}_{m-\ell+1} \right) \circ P_{\ell -1, m-\ell +1} 
	 	= \binom{m}{\ell-1} \left[\frac{(\ell-1)!}{2^{\ell-1}}\right] \cdot \mathsf{b}_m
	 \end{displaymath}
	 where $\binom{m}{\ell-1}$ is the Newton binomial.
	\end{lemma}
	\begin{proof}
Recall that $\mathsf{b}_m = - d_m ~\left(\pairing_-^{\ca (m-1)} \ca \langle B,\cdot \rangle_-\right) \circ f_1^{\otimes m}$ {(see eq. \eqref{eq:bm})}.
{We use this in the first and last equalities below,} to write
 the left-hand side above as
 
	\begin{displaymath}
		\mathclap{
		\begin{aligned}
			&(l.h.s.)~=
			\\ 
			&=
			 - d_{m-\ell+1} \cdot \pairing_-^{\ca (\ell -1)} \circ
			\left( 
				\mathbb{1}_{\ell-1} \otimes 
				\left(
					\pairing_-^{\ca (m -\ell)}\ca \langle B,\cdot \rangle_-
				\right)
			\right) \circ
			P_{\ell-1, m-\ell +1} \circ {f_1}^{\otimes m}	
			=	
			\\
			&=
			 - (-)^{(\ell-1)(m-\ell)} ~			
			d_{m-\ell+1} \cdot 
			\\
			&\phantom{=-}\cdot\left[
			\pairing_-^{\ca (\ell -1)} \circ			
			\left( 
				\left(
					\pairing_-^{\ca (m -\ell)}\ca \langle B,\cdot \rangle_-
				\right) \otimes			
				\mathbb{1}_{\ell-1} 
			\right) \circ
			P_{m-\ell +1, \ell-1 }\right] \circ {f_1}^{\otimes m}
			=
			\\
			&=
			 - (-)^{(\ell-1)(m)} ~			
			(-)^{(\ell-1)(|B|-k-1-m+\ell)} ~		
			d_{m-\ell+1} \cdot 
			\left( 
				\pairing_-^{\ca (m -1)}\ca \langle B,\cdot \rangle_-
			\right) \circ
			{f_1}^{\otimes m}
			=
			\\
			&=
			\frac{d_{m-\ell+1}}{d_m} \cdot 
			\mathsf{b}_m
			~.
	 	\end{aligned}
	 	}
	 \end{displaymath}
	 \note{cos'è quel $k$ che appare? Correggere!}
	 The sign term in the second equality comes from noting that, for any graded $b$-multilinear map $\nu_b$ on the graded vector space $V$, one has
	\begin{displaymath}
		\mathbb{1}_a \otimes \nu_b =
		(-)^{a(b+1)}~ \mathsf{C}_{(a+1)} \circ
		\left(\nu_b\otimes\mathbb{1}_a\right) \circ 
		\mathsf{C}^{-1}_{(a+b)},
	\end{displaymath}
	where $\mathsf{C}_{(i)}^j$ denotes the odd action of the cyclic permutation on $V^{\otimes i}$ repeated $j$ times.
	The sign term in the third equality comes from the sign convention in the definition of $\ca$.
	\\ 
	Finally, the claim now follows from an explicit computation of the coefficients
		\begin{displaymath}
			\dfrac{d_{m-\ell+1}}{d_m} 
			=
			\dfrac{ 2^{m-\ell+1}}{(m-\ell+1)!}~
			\dfrac{m!}{2^{m}} 
			=
			\dfrac{1}{2^{\ell-1}} \dfrac{m!}{(m-\ell+1)!}
			=
			\binom{m}{\ell-1} \dfrac{(\ell-1)!}{2^{\ell-1}} 
			~.
		\end{displaymath}
	\end{proof}
{Thanks to Lemma \ref{lem:rhsbm}, 	
 we can write eq. \eqref{eq:square} as follows, for any $m\geq 2$:}
	\begin{equation}\label{eq:nicesum}
		(\Phi \circ f)_m 
		-(\Phi\circ \widetilde{f})_m = 
		- \left[
				\sum_{\ell=1}^{m}
				\binom{m}{\ell-1}
				~
				\dfrac{(\ell-1)!}{2^{\ell-1}}\varphi_\ell
			\right]
			\circ 	\mathsf{b}_m
	\end{equation}	
		{where $\varphi_1=1$.}
		
	Observe now that the Bernoulli numbers $B_k$, given $B_0=1$,
	are completely defined by the the following summation formula for all $m\ge 2$  (see for instance \cite{Agoh}\cite{Weisstein}):
	\begin{equation}
  \sum_{j=0}^{m-1}
				\binom{m}{j}B_{j}=0.
\end{equation}	
	Using this recurrence relation, it is clear that the r.h.s. of equation \eqref{eq:nicesum} vanishes {provided the coefficients $\varphi_{\ell}$ in eq. \eqref{eq:Phim} are such that}
{$\dfrac{(\ell-1)!}{2^{\ell-1}}\varphi_{\ell}=B_{\ell-1}$  for $\ell=1,\dots,m$}.

We  conclude that diagram \eqref{diag:HamiltonianpentagonDiagram} commutes  {provided}
the following equation   holds true for any ${k}\geq 2$:
	\begin{equation}\label{eq:Phik}
		\Phi_k = \left(\dfrac{2^{k-1}}{(k-1)!}B_{k-1} \right) \circ \pairing_-^{\ca (k-1)}	
	~,
	\end{equation}
{The following table displays the values of the coefficient for low values of $k$:} 	
	\begin{center}
				\begin{tabular}{|l|llllllllll|}
					\hline
					$k$ & 1 & 2 & 3 & 4 & 5 & 6 & 7 & 8 & 9 & 10 \\
					\hline
					$\varphi_k = \frac{2^{k-1}}{(k-1)!}B_{k-1}$ & 1 & -1 & 1/3 & 0 & -1/45 & 0 & 2/945 & 0 & -1/4725 & 0\\
					\hline
				\end{tabular}
		\end{center}	
	In {theorem \ref{thm:iso}}} we proved that equation \eqref{eq:Phik} holds  true {for $k=2,3,4$ when $n\le4$}. 	\\
	{This concludes the proof of theorem \ref{thm:comm}.}

\ifstandalone
	\bibliographystyle{../../hep} 
	\bibliography{../../mypapers,../../websites,../../biblio-tidy}
\fi

\cleardoublepage


%% file: chapters/hydromomaps/hydromomaps.tex
\chapter{Hydrodynamical \Hcmm and linked vortices}\label{Chap:MauroPaper}
%
In this chapter, we discuss some applications of multisymplectic techniques in a hydrodynamical context. 
The possibility of applying symplectic techniques therein ultimately comes from Arnol'd's pioneering work culminating in the geometrization of fluid mechanics (\cite{Arnold66, Abraham1978,Arn-Khe,MW83}).
In particular, in this connection we may mention the paper \cite{Rasetti-Regge75}, with its symplectic reinterpretation \cite{Pe-Spe89,Pe-Spe92,Pe-Spe00}, and the general portrait depicted in \cite{Bry}. 
Here we wish to apply some recently emerged concepts in multisymplectic geometry (mostly building on \cite{Callies2016,Ryvkin2016,Ryvkin2018}) and 
construct an explicit {\it \momap} (\cite{Callies2016}) in a hydrodynamical setting, leading to a multisymplectic interpretation of the so-called {\it higher order linking numbers}, viewed \`a la Massey (\cite{Pe-Spe02,Spe06,Hebda-Tsau12}). 
The construction is generalized to cover connected compact oriented Riemannian manifolds (with a specified volume form) having vanishing intermediate de Rham groups. 
Moreover, a {\it covariant phase space} interpretation of the multisymplectic setting is outlined.
\par
We make clear from the outset that our constructions, together with the covariant phase space portrait, will not adhere to the standard multisymplectic approach to continuum mechanics set forth \eg in \cite{Gimmsy1, MPSW} but they will be based instead on the peculiar structure of an ideal fluid, whose configuration space is the ``Lie group" of  diffeomorphisms preserving a volume form, which will be directly taken as a multisymplectic form (\cite{CatIbort}).
\par
The content of this chapter is a joint work with Mauro Spera appeared in \cite{Miti2018} and \cite{Miti2019a}.

\medskip
The layout of the chapter is the following. 
First, in Section \ref{Sec:MSHydroFluids}, we give an example of \momap in fluid mechanics - in the sense of Callies-Fr\'egier-Rogers-Zambon (\cite{Callies2016}) 
- transgressing to Brylinski's symplectic structure on loop spaces and descending, in turn, to the manifold of smooth oriented knots, see \cite{Bry,BeSpe06} and below for precise definitions.  
We briefly discuss the (non-)equivariance of the above construction with respect to the group of volume-preserving diffeomorphisms of 3-space and we outline a generalization thereof in a Riemannian framework, signalling  potential topological obstructions. 
Moreover, covariant phase space aspects will be analysed.
In Section \ref{Sec:Ham1FormLinks} we prepare the ground for the forthcoming
applications by depicting a hydrodynamical multisymplectic portrait of basic knot theoretic objects, used,
in Section \ref{Sec:MasseyMess}, to reinterpret the Massey higher order linking numbers in multisymplectic terms: the $1$-forms appearing in the hierarchical Massey construction (viewed, in turn, differential geometrically \`a la Chen) provide an example of {\it first integrals in involution} in a multisymplectic framework. 
Appropriate background material is provided within the various sections in order to ease readability. 

\begin{remark}[Hydrodynamical brackets]
In this chapter, we are slightly departing from the \emph{hydrodynamical bracket convention} employed in \cite{Miti2018}.
	In that convention, the infinitesimal action pertaining to a left action is an anti-homomorphism, and the Lie bracket on the space of vector fields $\mathfrak{X}(M)=\Gamma(M)$ is given by minus the standard one (see \cite[pag. 6]{Arn-Khe}).
	In the convention employed throughout this text, the Lie bracket of vector fields is the standard one; in particular, $\mathcal{L}_X Y = [X, Y]$ for any $X,Y\in \mathfrak{X}(M)$.
\end{remark}

\section{Multisymplectic geometry and hydrodynamics of perfect fluids}\label{Sec:MSHydroFluids}
By \emph{hydrodynamics} we mean the discipline study the motion of liquids;
it can be seen as a branch of fluid mechanics and, in turn, of continuum mechanics.
In what pertains to us, we will only focus on the specific ideal model embodied by the perfect, inviscid, incompressible continuous medium (fluid) in space, \eg water in ideal conditions.
	%
	\begin{figure}[h!]
	 	 \begin{center}
	 	   \includegraphics[width=0.38\textwidth]{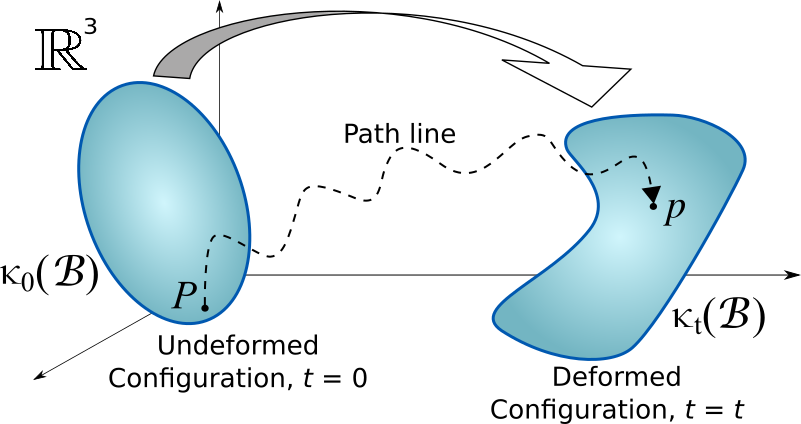}
		  \end{center}
		  \caption{Spatial displacement of a continuous bulk in the Euclidean space.\\
		  \small(\href{https://commons.wikimedia.org/wiki/File:Displacement_of_a_continuum.svg}{Wikimedia Commons})}
		  \label{Fig:FluidDisplacement}		
	\end{figure}

		In a rather general stance, fluid mechanics can be framed in the context of Riemannian geometry regarding the displacement of a bulk of material at a given time as a submanifold $B$, possibly with a regular boundary $\partial B$, embedded into the Riemannian "physical" space $(M,g)$ (see Figure \ref{Fig:FluidDisplacement}).
		Given a generic fluid system, its kinematics can be encoded by a "mass density" function $\rho$ and a "velocity field" $\mathbf{u}$, to be thought of as distributional sections in the sense of de Rham (see Remark \ref{Rem:ReminderDeRhamCurrents}).
		In the case that the fluid fills the whole ambient space and turbulence phenomena are negligible, this quantities can be seen as honest smooth objects $\rho \in C^\infty(M)$ and $\mathbf{u}\in \mathfrak{X}(M) $.
		\\
		The dynamics can therefore be formulated through evolutionary equations involving this pair of quantities.
		As far as we are concerned, we will only be interested in the following example:
		\begin{example}[Ideal fluids dynamics]\label{Ex:EulerEquations}
			The dynamics of an ideal fluid occupying an open region $\Omega \subset M$, with a regular boundary $\partial \Omega$, is ruled by the following system of equations:
		\begin{displaymath}\tag{EE}
			\label{Eq:EulerEquation}
			\begin{cases}
				\frac{\partial \mathbf{u}}{\partial t} + \nabla_\mathbf{u} \mathbf{u} & = -\nabla P  
				\\[.5em]
				\div (\mathbf{u}) 
				&= 0 \qquad\qquad \textrm{in } \Omega 
				\\[.5em]
				\mathbf{u} \cdot \hat{n} 
				&= 0 \qquad\qquad \textrm{on } \partial \Omega\\
			\end{cases}
			~,
		\end{displaymath}
		called \emph{Euler equation}, where $\nabla$ denotes the Levi-Civita connection, $\mathbf{u}$ is the velocity field of the fluid and $P$ is a function encoding the force exerted on an infinitesimal element of the fluid by its surroundings (pressure).
		\end{example}

		An alternative description of continuous systems, more akin to the framework of \emph{geometric mechanics}, can be achieved as follows.
		Note that, if one fixes a reference displacement $\kappa_0: B\hookrightarrow M$ of a given continuous system,
		any other spatial configuration can be encoded by a diffeomorphism $F= \kappa_t \cdot \kappa_0^{-1} : M \to M$ (see again Figure \ref{Fig:FluidDisplacement}).
		In other words, the configuration space of said continuous system, \ie the set of all its spatial displacements (not to be confused with the "physical states" of the system), is given by a suitable "Lie" subgroup $Q\subset \Diff(M)$. 
		As it is often done, we shall gloss over analytic subtleties coming from the infinite dimensionality of the group $\Diff(M)$ (see \eg \cite{Arnold66,Arn-Khe,Eb-Mar,Kri-Mich} for more information).
			We just recall here that $\Diff(M)$ is a {\it regular Lie group} in the sense of Kriegl-Michor \cite[43.1]{Kri-Mich} and that its associated exponential map is not even locally surjective (a quite general phenomenon). Further details can be found in the survey \cite{Roger2005}. 
		\\
		Two notable examples are given by the \emph{rigid body} and the \emph{incompressible fluid} filling the whole three-dimensional Euclidean space.
		In these cases, the configuration spaces are given respectively by the space of Euclidean isometries $\text{Iso}(\R^3)\cong \R^3 \rtimes SO(3)$ and the group of volume-preserving diffeomorphisms $\sDiff(\R^3)$.
		\\
		Formally speaking, physical states of a fluid mechanical system should be given by a point in the tangent bundle of $Q$ (or in the cotangent bundle when adopting the Hamiltonian framework). 
		Namely, fixed a spatial displacement $F\in Q$, a vector $v \in T_F Q$ can be thought of as vector field over $M$ integrating to an infinitesimal flow sitting in $Q$.
		\\
		In this framework, the Euler equation given in Example \ref{Ex:EulerEquations}
		can be deduced as an extremal of a certain Lagrangian (see for instance \cite{Arn-Khe,Abraham1978,Kri-Mich}) and its solutions can be recast as geodesics pertaining to a right-invariant metric on the group $\Diff(M)$ (see \cite{Arn-Khe} or \cite[Appendix A]{Khesin2020a}).

		\medskip
    In a preprint appeared in 1998 \cite{Gimmsy1}, Gotay, Isenberg, Marsden, Montgomery,  Sniatycki and Yasskin introduced a recipe to associate to any classical (relativistic but not quantum) first order field theory a multisymplectic manifold called \emph{multiphase space}. 
    This construction mimics the well-known prescription of an ordinary phase space starting from an associated configuration space, which is usually done in the case of point-like particles, in the context of continuous mechanics.
		(See section \ref{Sec:CPSMess} for further details.)
		\\
		In particular, this would also apply to any ideal fluid system. In the following, however, we will not adhere to this construction.
Instead, we will focus on the other "multisymplectic" feature peculiar of ideal fluid systems; more specifically, our key point will be that their configuration space is a Lie group consisting of multisymplectic (\ie preserving a multisymplectic form) diffeomorphisms.

\subsection{The hydrodynamical Poisson bracket}\label{Sec:IdroPoisson}
		In the present subsection we briefly review, for motivation and further applications, the symplectic geometrical portrait underlying the theory of hydrodynamics in its simplest instance.
		Namely we will focus on the free motion of an ideal fluid filling the standard "physical" space.
		More concretely, let us assume the following conditions:
	\begin{itemize}
		\item The ambient space is given by $M=\R^3$ taken with the standard volume form $\nu = \d x\wedge \d y \wedge \d z$.
		\item The configuration space is given by the "Lie" group $G = \sDiff({\mathbb R}^3)$ of volume-preserving diffeomorphisms of ${\mathbb R}^3$.
		\item The dynamics is ruled by equation \eqref{Eq:EulerEquation} with $P=0$.
	\end{itemize}	 
	
	We denote by ${\mathfrak g}$ the ({\it infinite dimensional}) Lie subalgebra of ${\mathfrak X}({\mathbb R}^3)$ consisting of divergence-free vector fields on ${\mathbb R}^3$. These are the vector fields integrating to volume-preserving diffeomorphisms. 
	In this sense $\mathfrak{g}$ is the ``Lie algebra" of the ``Lie group"  $G$.
	It is customary to tacitly assume that our fields {\it rapidly vanish at infinity}. 
	This is justified by the physical assumption that the motion of the system is localized in a finite region of the ambient space.
	Such condition ensures that convergence problems are avoided and boundary terms are absent.
	Briefly, we also assume the following:
	\begin{equation}\label{Eq:IdrodynamicalLieAlgebra}
		\begin{aligned}
		\mathfrak{g} :=&~ \rm{sdiff}_0\,(\R^3)~=
		\\
		=&~ \left\lbrace  X \in \mathfrak{X}(\R^3) ~\left\vert~ \rm{div} X = 0, \textrm{\emph{ rapidly vanishing at }}\infty \right.\right\rbrace
		\end{aligned}~.
	\end{equation}
	
	It is crucial to observe that equation \eqref{Eq:EulerEquation} is exclusively stated in terms of the velocity field of the system, hence the physical states of the system are completely encoded by $\mathfrak{g}$
	(this can be expected a priori since the incompressibility condition would imply that the density function must by a constant.)
	\begin{remark}[Vorticity]\label{Rem:Vorticity}
			{\it Euler evolution} can be read, among others, in the so-called {\it vorticity form}:
			\begin{equation}
			\label{Eq:EulerVorticityForm}
					\begin{cases}
	  				\frac{\partial \mathbf{w}}{\partial t} &= [\mathbf{w}, \mathbf{u}] 
	  				\\
						\mathbf{w} &= \curl \left( \mathbf{u} \right)  \\
					\end{cases}
  		\end{equation}
		We call $\mathbf{w}=\curl \mathbf{u}$ the \emph{vorticity field} pertaining to the velocity field $\mathbf{u}\in \mathfrak{g}$ (hence divergence-free).
		\\
		Recall also the standard result about  vector analysis in $\R^3$ that the Lie bracket of two divergence vector field can be expressed via the curl of their cross product
	\begin{equation}\label{eq:BrackAsCurl}
		[\mathbf{v},\mathbf{w}] = -\curl (\mathbf{v}\times \mathbf{w})
		~.
	\end{equation}
	\end{remark}

	Following \eg \cite{Arn-Khe}, we shall consider the so-called  {\it regular dual} ${\mathfrak g}^*$ of ${\mathfrak g}$
consisting of all $1$-forms modulo exact $1$-forms:
\begin{equation}\label{eq:gastPoisson}
	{\mathfrak g}^* := \Omega^1({\mathbb R}^3)/ d \Omega^0({\mathbb R}^3)
\end{equation}
together with the standard pairing ($\omega \in {\mathfrak g}^*$, $\xi \in {\mathfrak g} $)
$$
(\omega, \xi) = \int \langle \omega (x), \xi(x) \rangle \, d^3x .
$$
	Observe that this is different from the the full topological dual that, in principle, would contain also suitable genuine distributional elements as well (\ie {\it currents}, in the sense of de Rham, see \cite{dR} and remark \ref{Rem:ReminderDeRhamCurrents} below).

	\begin{theorem}[\cite{Arn-Khe,Kuznetsov-Mikhailov80,MW83,Pe-Spe89,Pe-Spe92,Pe-Spe00,Spera16}]\label{Thm:HydroBracket}
		The (regular) dual ${\mathfrak g}^*$ is naturally interpreted as a Poisson manifold with respect to the hydrodynamical Poisson bracket  (Arnol'd --Marsden-Weinstein Lie-Poisson structure)
		\begin{equation}\label{Eq:HydroBracket}
			\{ F, G \} ([{\mathbf v}]) 
			= 
			\int_{ \R^3} \left\langle {\mathbf v}, \,  \left[\frac{\delta F}{\delta {\mathbf v}}, \frac{\delta G}{\delta {\mathbf v}}\right]  \right\rangle \,\d^3x
		\end{equation}
		with ${\mathbf v} \in {\mathfrak g}$ (velocity field), ${\mathbf w} := {\rm curl}\,{\mathbf v}$,
		its {\it vorticity},  with $[{\mathbf v}]$ denoting the ``gauge" class of ${\mathbf v}$: $[{\mathbf v}] = 
		\{ {\mathbf v}+ \nabla f \}$.	
	\end{theorem}

	The {\it Euler evolution}, given for instance by equation \eqref{Eq:EulerVorticityForm}, is naturally {\it volume-preserving} and it also preserves the {\it symplectic leaves} of ${\mathfrak g}^*$ given by the $G$-{\it coadjoint orbits} ${\mathcal O}_{[{\mathbf v}]} \equiv {\mathcal O}_{{\mathbf w}}$.
	The symplectic structure on ${\mathcal O}_{\mathbf w}$, given by the celebrated  Kirillov-Kostant-Souriau (KKS) construction
(\cite{Kirillov01,Kostant70,Souriau70}),  is explicitly given by
\begin{align*}
	{\Omega_{KKS}}{([{\mathbf v}])} (ad^*_{{\mathbf b}}([{\mathbf v}]),ad^*_{{\mathbf c}})([{\mathbf v}]) 
	=&~ \int_{{\mathbb R}^3} \langle {\mathbf v}, [{\mathbf b},{\mathbf c}]\rangle \, d^3x ~=
	\\
	=&~
	\int_{{\mathbb R}^3} \langle {\mathbf w}, {\mathbf b} \times {\mathbf c}\rangle\, d^3x
\end{align*}
with the coadjoint action reading, explicitly, {\it up to a gradient} (not influencing calculations)
$$
ad^*_{{\mathbf b}}  ({\mathbf v}) = - {\mathbf w} \times {\mathbf b} 
\quad(\equiv ad^*_{{\mathbf b}}  ([{\mathbf v}])
~.
$$

The {\it Hamiltonian algebra} $\Lambda$ pertaining to ${\mathcal O}_{{\mathbf w}}$ consists of the so-called {\it Rasetti-Regge currents} originally introduced in \cite{Rasetti-Regge75} and further developed in
\cite{Pe-Spe89,Pe-Spe92,Pe-Spe00,Spera16,Bry}:
	\begin{definition}[Rasetti-Regge currents algebra]\label{def:RR-current}
		We call \emph{Rasetti-Regge current} pertaining to a given $\mathbf{b} \in \mathfrak{g}$ the linear function $\lambda_{\mathbf{b}}$ in the topological dual of $\mathfrak{g}$ given by
		\begin{displaymath}
			\morphism{\lambda_{\mathbf{b}}}
			{\mathfrak{g}}
			{\R}
			{\mathbf{v}}
			{\displaystyle\int_{\R^3}\left\langle \mathbf{b} , 
			\mathbf{v}\right\rangle d^3x 
			\,= \int_{\R^3}\left\langle {\mathbf B},
			\curl({\mathbf v}) \right\rangle d^3x} 
		\end{displaymath}	
		where $\mathbf{B}\in \mathfrak{X}(\R^3)$ is an arbitrary vector field such that
		$\curl({\mathbf B}) = {\mathbf b}$.
		\\
		We call \emph{Rasetti-Regge currents algebra}  the vector space $\Lambda = \left\lbrace \lambda_{\mathbf b} \right\rbrace_{{\mathbf b}\in\mathfrak{g}}$ endowed with the skew-symmetric bilinear structure $\lbrace\cdot,\cdot\rbrace$ given by
		\begin{displaymath}
			\{\lambda_{\mathbf b}, \lambda_{\mathbf c} \} = \lambda_{[{\mathbf b},{\mathbf c}] }
			\qquad \forall ~{\mathbf b}, {\mathbf c} \in {\mathfrak g}
			~.
		\end{displaymath}
	\end{definition}
	\begin{theorem}[\cite{Pe-Spe92}]
		\begin{enumerate}[label=(\roman*)]
			\item $\Lambda$ is a Lie algebra and the map 
				\begin{displaymath}
					\morphism{\lambda}
					{\mathfrak{g}}
					{\Lambda}
					{\mathbf{b}}
					{\lambda_{\mathbf{b}}}
				\end{displaymath}							
			gives a {\it $G$-equivariant co-momentum map}. 
			(Observe in particular that
$\frac{\delta \lambda_{{\mathbf b}}}{\delta {\mathbf v}} = {\mathbf b}$).
			\item The Euler equation \eqref{Eq:EulerEquation} can be recast in term of the following Hamilton equation
			\begin{displaymath}
				\partial_t ~ \lambda_{\mathbf{b}} = - \langle H, \lambda_{\mathbf{v}}\rangle
			\end{displaymath}
			with Hamiltonian given by $H(\mathbf{v})=\frac{1}{2}\langle \mathbf{v},\mathbf{v}\rangle$. 
		\end{enumerate}
		\note{
			Nelle note si accenna a "Kutnetsiv-Mikhailov" Poisson brackets coincidente con KKS Poisson bracket. Devo spiegare qualche fonte?

			sarebbe l'articolo \cite{Kuznetsov-Mikhailov80}, citato esplicitamente da Mauro in \cite{Pe-Spe92}.
		}
	\end{theorem}
	\begin{remark}[Helicity]\label{Rem:Helicity}
		We call \emph{helicity}, pertaining to an ideal fluid with velocity ${\mathbf v}$ and vorticity ${\mathbf w} = {\rm curl}\,{\mathbf v}$ in $\R^3$,
		the quantity
		\begin{displaymath}
			{\mathcal H} = \int_{\R^3} \langle {\mathbf v},  {\mathbf w}\rangle \d^3x 
			= \int_{\R^3} v \wedge dv \d^3x
		\end{displaymath}	
		where the last expression is the differential form counterpart.
		\\
		Helicity is preserved along the Euler evolution of the fluid \cite{Moffatt-Ricca92}:
		\begin{align*}
			\partial_t \mathcal{H} 
			=&~ \partial_t \int \langle \mathbf{v},\mathbf{w}\rangle
			=
			\\
			=&~  \int \langle \partial_t\mathbf{v},\mathbf{w}\rangle + \int \langle \mathbf{v},\partial_t\mathbf{w}\rangle
			=
			\\
			=&~ 
			- \lbrace H, \lambda_{\mathbf{w}}\rbrace + \int \langle\mathbf{v},[\mathbf{w},\mathbf{v}]\rangle = 0
			~.
		\end{align*}
	\end{remark}

\subsubsection*{Singular vortices}
	The preceding portrait carries through to the {\it singular} vorticity case, in particular when the vorticity field is $\delta$-like and concentrated on a two-dimensional patch or -possibly knotted- filament.
	We like to stress that the latter hydrodynamical configurations have been extensively studied and are recently proved to exists both analytically \cite{Enciso2015} and experimentally \cite{Kleckner2013}.
	It is crucial to observe that the support of the vorticity $\mathbf{w}$ of an ideal fluid is preserved along the motion in the sense that $\text{supp}(\mathbf{w}(0))$ and $\text{supp}(\mathbf{w}(t))$ are diffeomorphic for any time $t$.

	\begin{figure}[h!]
	  \begin{center}
			\href{https://www.nature.com/articles/nphys2560}{\includegraphics[width=0.32\textwidth]{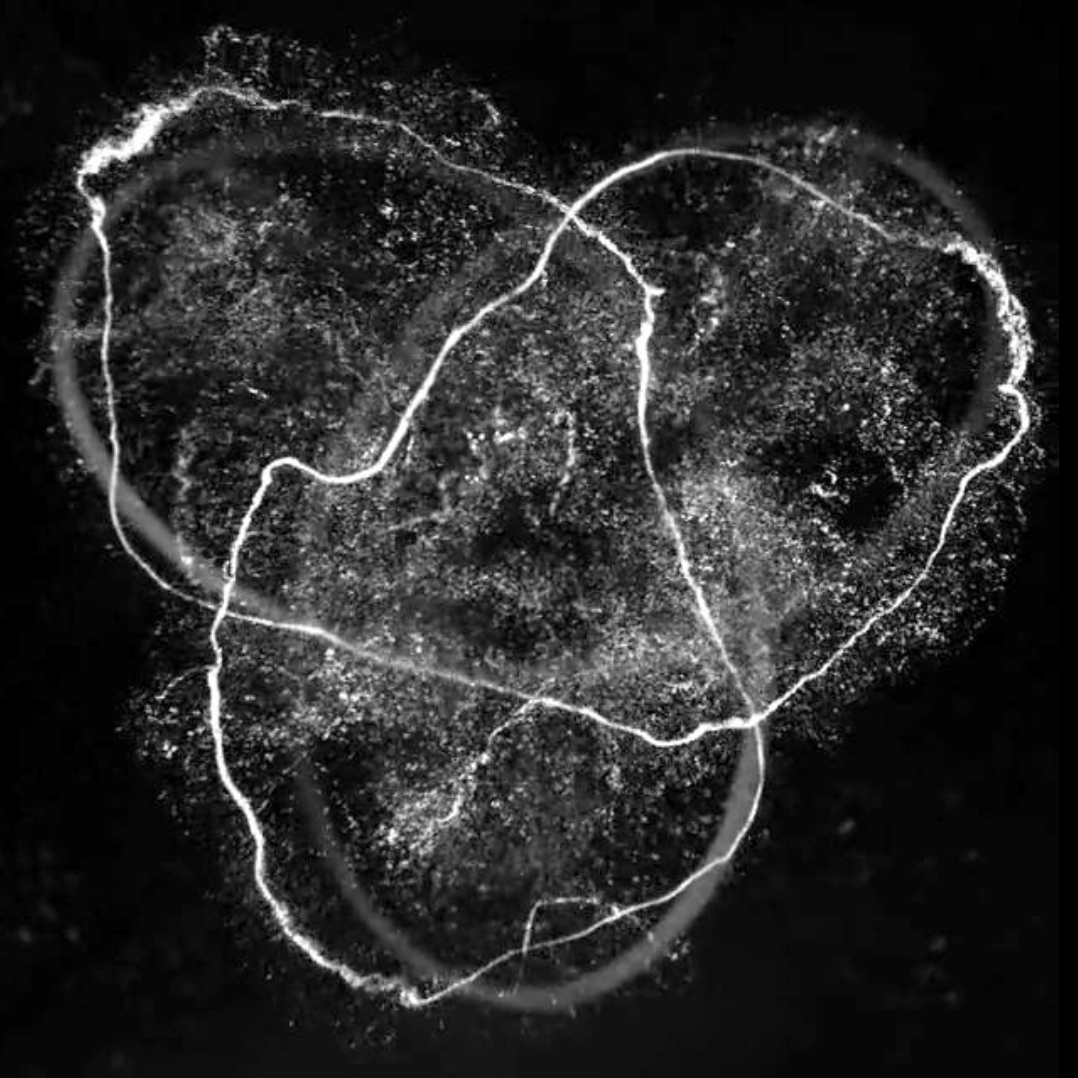}}
	  \end{center}
		\caption{Experimental realization of a knotted vortex in water. (Klenecker \& Irvine \cite{Kleckner2013})}
	  \label{Fig:KnottedVortexExp}
	\end{figure}

	When dealing with the embeddings of loops in the Euclidean space we will make use of the following nomenclature due to Brylinski:
	\begin{notation}[Loop spaces and Brylinski manifolds]\label{Rem:BryLoopSpaces}
		Let $M$ be a smooth, paracompact, oriented manifold of dimension $n$.
		(Throughout all this chapter we will mostly assume $M=\R^3$.)
		\begin{itemize}
			\item 
			We call \emph{(free) loop space} the set of smooth maps from the circle into the manifold, it will be denoted by $LM=C^{\infty}(S^1,M)$.
		This set can be made into a (infinite dimensional) smooth Fr\'echet manifold modelled on the topological vector space $C^\infty(S^1,\R^n)$ (which is Lind\"elof and paracompact)\cite[\S{3.1}]{Bry}.
			One can make sense of the tangent space to a $\gamma \in LM$ as
			\begin{displaymath}
				T_{\gamma}LM = C^\infty(S^{1},\gamma^{\ast}TM)= \Gamma^\infty(\gamma^\ast TM)	
				~.
			\end{displaymath}
			Observe that the Lie\footnote{
				Given a manifold $M$ with a volume form $\Omega$, $\sDiff(M,\Omega)$ is a bona fide Lie group at least in the case that $M$ is compact \cite[prop. 2.1]{Molitor2007}. 
				In particular, $\sDiff(S^1)$ is a rather well-behaved Lie group, modelled on a Fréchet space.
				When $M$ is a non-compact manifold, things are more subtle.
			} group $\sDiff(S^1)$ of oriented diffeomorphisms of $S^1$ acts smoothly and on the right on $LM$ via precomposition. 
			This action, determines all possible re-parametrization of a given loop.			
			
			\item 
			We call \emph{space of mildly singular knots} the subset $\hat{X}\subset LM$ consisting of immersions $\gamma$ which have the property to induce an embedding of $S^1\setminus A$,
			for $A$ a finite subset of $S^1$, and such that the branches of $\gamma$ at any distinct points $x_1$ and $x_2$ of $A$  have at most finite-order tangency.
			
			\item
			We call \emph{Knot space} of $M$ the subset $X\subset \hat{X}$ consisting of bona fide embeddings.
		(When $M=\R^3$, the image of $\gamma \in X$ is called \emph{knot} and $\gamma:S^1\to \R^3$ is called a \emph{parametrization}.)
			\item 
			We call \emph{space of oriented singular knots} and \emph{space of oriented knots} in $M$ the quotients $\hat{Y} = \hat{X}/\sDiff(S^1)$ and $Y = X/\sDiff(S^1)$ respectively.
			Note that $\hat{Y}$ is taken with a Fr\'echet structure such that the projection $\hat{X}\to \hat{Y}$ is a $\sDiff(S^1)$-principal bundle. Thence, $Y$ is an open subset of $\hat{Y}$ such that the projection $X\to Y$ is a locally trivial $\sDiff(S^1)$-principal bundle. 
		\end{itemize}				
	\end{notation}
	Considering the case of vorticity concentrated on a knot in (smooth embedded loops modulo orientation-preserving reparametrizations in $\R^3$), we ultimately retrieve the symplectic structure ${\Omega_Y}$ on the the {\it Brylinski manifold $Y$}
	(see \cite{Bry,BeSpe06} for more details):
\begin{equation}
\Omega_Y (\cdot , \cdot) ({\gamma}) := \int_{\gamma} \nu (\dot{\gamma}, \cdot, \cdot) =
\int_{\gamma} \langle \dot{\gamma}, \cdot \times \cdot \rangle
~.
\end{equation}
Indeed, given a volume form $\nu$ on a 3-dimensional $M$, one gets, by {\it transgression} (see Example \ref{Ex:TransgressionOnLoops}), a $2$-form $\Omega$ on $LM$ via the formula
	\begin{equation}
	\Omega = \int_{S^1} ev^* (\nu ) ~,
	\end{equation}
where $ev:  LM \times S^1\rightarrow M$ given by $ev (\gamma , t ) := \gamma(t)$ is the evaluation map of a loop $\gamma \in LM$ at a point $t \in S^1 \equiv [0,1] / \,\,{}_{\widetilde{}}\,$ (endpoint identification) and $\int_{S^1}$ denotes integration along the $S^1$-fibres (see \cite[\S 3.5]{Bry}). 
More explicitly, given  tangent vectors $u$ and  $v$ at $\gamma$, the symplectic form reads  (\cf \cite{Bry}, formula 6 - 8, p. 238):
\begin{equation}
\Omega_{\gamma} (u , v)  = \int_{0}^1  \nu ({\dot \gamma}(t), u(t), v(t))\, dt
\end{equation}
(where we set ${\dot \gamma} = {d\gamma \over d t}$).
The above construction carries through to $Y$. In this case, the coadjoint orbits are labelled by the equivalence types of knots (via ambient isotopies),
by virtue of a result of Brylinski\cite{Bry}.

	\begin{remark}
			Rasetti-Regge currents has been originally introduced in the context of the theory of singular vortices\cite{Rasetti-Regge75}.
			Considering a velocity field configuration $\mathbf{v}\in \mathfrak{g}$ with vorticity $\mathbf{w}= \delta_\gamma$ 
			concentrated on a closed loop $\gamma: S^1 \to \mathbb{R}^3$,
			one could denote the Rasetti-Regge current on the vortex filament as
			\begin{displaymath}
				\lambda_{\mathbf b} (\gamma)
				:=
				\lambda_{\mathbf{b}} (\mathbf{v})	
				=	\int_{\R^3} \mathbf{v} \cdot \mathbf{b} 
				= \int_{\R^3} \mathbf{v} \cdot \curl(\mathbf{B}) 
				= -\int_{\R^3} \mathbf{w} \cdot \mathbf{B}
				= - \oint_\gamma \mathbf{B}		
			\end{displaymath}				
			where $\mathbf{B}$ is a vector field such that $\mathbf{b}=\curl{B}$ 
			(vector potential of $b$, see \ref{Rem:ConcreteSolutionBiotSavart}).

			\begin{figure}[h!]
				\begin{center}
					\begin{tikzpicture}[use Hobby shortcut]
						\begin{knot}[
					 	 	consider self intersections=true,
					 		ignore endpoint intersections=false,
					 		flip crossing=5,
						]
							\strand[blue] ([closed]0,0) .. (0.5,1) .. (-0.5,2) .. (-0.65,2.5) .. (0,3) .. (0.5,2) .. (-0.5,1) ..  (0,0);
							\strand (-0.65,2.5) circle[x radius=0.3, y radius=0.15];
						\end{knot}
						\draw[blue,-stealth](-0.65,2.5)--(-0.65,3) node[anchor=north east]{$\omega$};
						\draw[black,-stealth](-0.65,2.35)--(-0.3,2.35) node[anchor= west]{$v$};
						\draw[black,-stealth](-0.95,2.5)--(-0.95,2.3) node[anchor= west]{};
						\draw[black,-stealth](-0.35,2.5)--(-0.35,2.7) node[anchor= west]{};
						\node[blue,anchor= north] at (0,0) {$\gamma$};
					\end{tikzpicture}		
				\end{center}
			  \caption{Singular vorticity concentrated on a closed curve $\gamma$.}
			\end{figure}
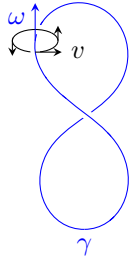
	\end{remark}

\subsection{A hydrodynamical \momap}\label{Sec:IdroMoMap}
In this Subsection we elaborate on the previous discussion by introducing an explicit {\it \momap} (\Hcmm) emerging in fluid dynamics.
We depart from the standard setting since we are going to consider actions of an infinite dimensional group $G = \sDiff(\R^3)$.
\\
More precisely, we consider an infinitesimal action pertaining to a subalgebra $\mathfrak{g}$ of the Lie algebra of $G$ given by equation \eqref{Eq:IdrodynamicalLieAlgebra}.
We ought to notice that the study of the deformation quantization problem for this Lie algebra, carried out by Claude Roger in \cite{Roger2012}, brought to the construction of a $L_\infty$-algebra anticipating the one proposed by Chris Rogers in the general multisymplectic case.

In example \ref{Ex:VolumesAreMultiSymp}, it has been noted how any oriented $n$-dimensional manifold $M$ can be interpreted as multisymplectic interpreting a volume form determining its orientation as a $(n$-$1)$-plectic form.
In the case we are considering, the volume form takes the canonical expression $\nu :=  dx \wedge dy \wedge dz $.
%
Denoting by $\alpha$ the flat operator pertaining to the multisymplectic form $\nu$ (see equation \eqref{Eq:OmegaFlat}), one has, in coordinates for any $\xi = (\xi^i)$, that
$$
\alpha(\xi) = \iota_{\xi} \nu = \xi^1 dy \wedge dz + \xi^2 dz \wedge dx + \xi^3 dx \wedge dy
$$
which is manifestly bijective.

It is also important to recall that the standard volume on $\R^3$ is in particular coming from a (standard) Riemannian structure.
Hence, we have at our disposal the Hodge $*$ operator which allows us to compare differential forms of different degree.
\begin{reminder}[Hodge calculus]\label{Rem:HodgeCalculus}
	For the sake of completeness, let us briefly recall the basic operators involved in the Hodge calculus (the reader can refer to \cite{Warner,Kri-Mich} for a complete account).
	\\
	Recall that, given a $n$-dimensional vector space $V$ endowed with an inner product $\pairing$,
	there exists an induced inner product on the skew-symmetric tensor space $\Lambda( V)$, denoted with the same symbol,
	given on homogeneous $k$-vectors by
	\begin{displaymath}
		\morphism{\pairing}
		{\Lambda^p V \otimes \Lambda^q V}
		{\R}
		{u_1\wedge\dots\wedge u_p,~ v_1\wedge\dots\wedge v_q}
		{\text{det}\left[ \langle u_i,v_j\rangle\right] \delta_{p, q}}
	\end{displaymath}		
	and extended by linearity.
	Considering on the same vector space the orientation, specified by the volume form $\nu$, induced by the inner product, one can construct an invertible linear operator $\ast: \Lambda^k (V)\to \Lambda^{n-k}(V)$ defined completely, for any $0\leq k \leq n$, by the following property
	\begin{displaymath}
		\alpha \wedge (\ast \beta ) = \langle \alpha, \beta \rangle \nu
		\qquad \forall \alpha \in \Lambda^k (V)
		~.
	\end{displaymath}
In particular $**\eval_{\Lambda^k(V)} = (-)^{k(n-k)} \id$.
	\\
	Given a Riemannian manifold $(M,g)$, oriented with the corresponding Riemannian volume, such construction can be replicated locally on the tangent and cotangent space for any points $p\in M$.
This yields a well-defined isomorphism called \emph{Hodge dual}, or "star",
	\begin{displaymath}
		\morphism{\ast}
		{\Omega^k(M)}
		{\Omega^{n-k}(M)}
		{\sigma}
		{\ast \sigma}
	\end{displaymath}
	uniquely defined by the following property:
	\begin{displaymath}
	 \omega \wedge \ast \sigma = \langle \omega, \sigma \rangle \nu 
	 \qquad \forall \omega \in \Omega^k(M)
	 ~.
	\end{displaymath}			 
Here $\langle \omega, \sigma \rangle$ denotes the smooth function on $M$ given point-wise on $p\in M$ by the inner product $\langle \omega_p, \sigma_p \rangle$. 
When $M$ is compact, the integration of the previous expression defines an inner product on $\Omega^k(M)$.
	\\
	Through the Hodge dual it is possible to introduce the following differential operators:
	\begin{displaymath}
		\begin{aligned}
			\delta :=~ (-)^{n(k+1)+1}\ast\circ\d\circ\ast \quad & : \Omega^k(M) \to \Omega^{k-1}(M)
			\\
			\Delta :=~ \delta\circ\d + \d\circ\delta \quad & : \Omega^k(M) \to \Omega^{k}(M)
			\\
			\grad :=~ \sharp \circ \d \quad & : C^\infty(M) \to \mathfrak{X}(M)
			\\
			\div :=~ -\delta \circ \flat = \ast\circ \d\circ \ast \circ \flat \quad & : \mathfrak{X}(M) \to C^\infty(M)
			\\
			\curl :=~ \sharp \circ \ast \circ \d  \circ \flat	\quad & : \mathfrak{X}(M) \to \mathfrak{X}^{n-2}(M)
			~,		
		\end{aligned}
	\end{displaymath}
	where $\flat$ and $\sharp$ denote the Riemannian "musical isomorphisms" associated to the metric.
	When $M=\R^3$ these operators assume the more familiar expression involving the gradient operator $\nabla$.
\end{reminder}
Upon  introducing the Hodge $*$ relative on $(\R^3,\nu)$, one can easily recast the operator $\alpha$ as the invertible linear map
\begin{equation}\label{eq:alphacontractionmap}
	\morphism{\alpha = (\ast \circ \flat)}
	{\mathfrak{X}(\R^3)}
	{\Omega^1(\R^3)}
	{\xi}
	{\ast(\xi^{\flat})}
	~.
\end{equation}
Then we have, for $\xi \in \mathfrak {g}$ (via Cartan's formula)
$$
0 = {\mathcal L}_{\xi} \nu = d  \iota_{\xi} \nu  +  \iota_{\xi} d \nu  = d \iota_{\xi} \nu = \div(\xi)~\nu
$$
and thus we have an isomorphism $\mathfrak {g} \cong Z^2(\R^3)$ (closed $2$-forms on ${\mathbb R}^3$). 
This will be important in the sequel. 
The above also expresses the fact that $\nu$ is a {\it strictly conserved} $3$-form (see definition \ref{Def:conservedQuantities}).
\\
We are now ready to give the promised example of {\rm \momap} pertaining to the action of $\sDiff_0$ on $\R^3$.
\begin{theorem}[Explicit \Hcmm for $\sDiff_0 \circlearrowright (\R^3,\nu)$]\label{Thm:HydrodynamicalComoment}
	Consider the infinitesimal action of $v:\mathfrak{g}\to \mathfrak{X}(\mathbb{R}^3)$ concretely given by the inclusion of divergence free fields in the set of all vector fields.
	The previous action admits a \momap $(f)$ with components $f_j: \Lambda^j {\mathfrak g} \to  \Omega^{2-j} (\R^3)$ given by
			\begin{displaymath}
				\begin{aligned}
					f_1 =&~ \flat\circ {\rm curl}^{-1}
					\\
					f_2 =&~ \Delta^{-1} \circ \delta~ \circ\mu_2		
				\end{aligned}
			\end{displaymath}
		where $\mu_2(p):= f_{1} (\partial p) +  \iota(v_p) \omega$ is the term introduced in remark \ref{Rem:TermMuByMauro}.
		The inverse of the vector calculus operators involved, $\curl$  and Laplacian $\Delta$ (see reminder \ref{Rem:HodgeCalculus}),  are to be interpreted as their corresponding Green operators. Hence they are not unique.
In remark \ref{Rem:ConcreteSolutionBiotSavart} we provide a more explicit expression for the components $f_k$ in the standard coordinates of $\R^3$ .
\end{theorem}
\begin{proof}
	$(1)$
In this case, the observables are given by a graded vector space concentrated in degrees $-1$ and $0$,
	 $$L= \Omega^1_{\textrm{ham}}(\mathbb{R}^3)\oplus\left(\Omega^0(\mathbb{R}^3)[1]\right)~,$$
	hence, a \momap, consists of a pair of functions:
	\begin{align*}
		f_1 &\colon \mathfrak{g} \rightarrow \Omega^1_{\textrm{ham}}(\mathbb{R}^3) \\
		f_2 &\colon \mathfrak{g}\wedge\mathfrak{g} \rightarrow C^\infty(\mathbb{R}^3)
	\end{align*}
	sitting in the following (non-commutative) diagram:
\begin{displaymath}
	\begin{tikzcd}[column sep=large,
	   execute at end picture={
	     \node[label=below:{\tiny CE complex},dashsquare=(L1)(L2)]{};
	     \node[label=below:{\tiny Observables},dashsquare=(R1)(R2)]{};
			}]
		&	& \Omega^3(M) \ar[ddd,leftrightarrow,green!60!black,bend left=60,"\ast"]\\
		& \mathfrak{X}(M) \ar[r,red,"\alpha"] \ar[dr,blue,"\flat",shift left=1ex]\ar[dr,blue,leftarrow,"\sharp"',shift right=1ex]& \Omega^2(M) \ar[u,"d"']  \ar[d,leftrightarrow,green!60!black,bend left=60,"\ast"] \\
		|[alias=L1]| \mathfrak{g} \ar[ru,hookrightarrow,"v"] \ar[rr,purple,"f_1"]
		& 
		& |[alias=R1]| \Omega^1_{(ham)}(M) \ar[u,"d"'] \\
		|[alias=L2]| \mathfrak{g} \wedge \mathfrak{g}  \ar[u,"\partial"]\ar[rr,purple,"f_2"]
		& \color{purple}\underbrace{\qquad\qquad}_{\Hcmm} 
		& |[alias=R2]| \Omega^0(M) \ar[u,"d"']
	\end{tikzcd}
\end{displaymath}	
	and satisfying the following system of three equations:
	\begin{subequations}
	   \begin{equation}\label{eq_1Idro}
					\d f_1(\xi) = -\iota_\xi \nu = -\alpha(\xi) 
		 \end{equation}
		 \begin{equation}\label{eq_2Idro}
				\d f_2(\xi_1 \wedge \xi_2) = f_1\left([\xi_1,\xi_2]\right) - \iota_{\xi_2}\iota_{\xi_1} \nu 
								 := \mu_2(\xi_1,\xi_2)
		 \end{equation}
		 \begin{equation}\label{eq_3Idro}
		 				f_2\left(\partial \xi_1 \wedge \xi_2 \wedge \xi_3 \right) = \iota_{\xi_3}\iota_{\xi_2}\iota_{\xi_1} \nu
		 \end{equation}
	\end{subequations}
%
	(Compare with Lemma \ref{Lem:ExplicitHCCM} fixing, in the case $k=2$, $\xi_i \in {\mathfrak g}$ ($i=1,2$),  $p = \xi_1 \wedge \xi_2$, and $\partial p = -[\xi_1, \xi_2]$.)
\\
	The first two equations imply that the components of the \momap are  primitives of the closed forms on the right-hand side.
	Closedness of the $1$-form $\mu_2(\xi_1,\xi_2) := f_1([\xi_1,\xi_2]) - \iota_{\xi_1 \wedge \xi_2} \nu$ can be checked using lemma \ref{lemma:multicartan}, hence
$$
df_1 ([\xi_1, \xi_2]) = d (\iota_{\xi_1 \wedge \xi_2} \nu) = -\iota_{[\xi_1,\xi_2]} \nu
$$
where $\iota_{\xi_1 \wedge \xi_2} \nu = \nu(\xi_1,\xi_2, \cdot)$.
According to Poincar\'e Lemma, such primitives can always be found but are in principle determined only up to a constant.

Fixing $\mathbf{b} \in \mathfrak{g}$, as noted below remark \ref{Rem:HodgeCalculus}, equation \eqref{eq_1Idro} can be easily recast as
\begin{displaymath}
	\d f_1 (\mathbf{b}) = - \ast \circ \flat \circ \mathbf{b} 
	~.
\end{displaymath}
Inverting the Hodge and flat operators one gets
\begin{displaymath}
	- \left(~\sharp \circ \ast \circ \d \circ f_1\right)(\mathbf{b}) = \mathbf{b}
	~.
\end{displaymath}
Hence, introducing the vector field $\mathbf{A}= -~\sharp \circ f_1(\mathbf{b})$,
the first component of the \momap can be given, modulo a constant $c_1(\xi)$, by solving the following equation on smooth vector fields
\begin{equation}\label{Eq:MagnetoEq}
	\curl \mathbf{A} = \mathbf{b}
\end{equation}
which is the well-known equation of magnetostatic, see remark \ref{Rem:ConcreteSolutionBiotSavart} for the explicit solution.

Applying the $\delta$ operator, equation \eqref{eq_2Idro} yields that
\begin{equation}\label{Eq:PoissonEq}
	\Delta (f_2(\xi_1\wedge\xi_2)) = \delta~ \mu_2(\xi_1,\xi_2)
\end{equation}
hence $f_2$ can be expressed, modulo a constant $c_2(\xi_1,\xi_2)$, as a solution of a Poisson equation in $\R^3$ with source $\delta~ \mu_2(\xi_1,\xi_2)$.

 In order to prove that we have a bona fide \momap, we must have, in particular, for $q = \xi_1 \wedge \xi_2 \wedge \xi_3$, the explicit formula
\begin{equation}\label{eq:condition3_hccm}
	f_2(\partial q) = \nu(\xi_1, \xi_2, \xi_3) 
\end{equation}
which is a priori true up to a constant $c_3(\xi_1, \xi_2, \xi_3)$ by virtue of (\ref{eq_1Idro}) and \cite[Lemma 9.2]{Callies2016}. 
However, the constant is in fact zero since 
 $\nu(\xi_1, \xi_2, \xi_3)$ vanishes at infinity and the same is true for $f_2(\partial q)$ upon solving the related Poisson equation
 $$
 \Delta f_2(\partial q) = \Delta \nu(\xi_1, \xi_2, \xi_3)
 $$
 (obtained via a straightforward computation; notice that we use the Riemannian Laplacian, which is {\it minus} the standard one).
 \\
An alternative derivation uses $x$-independence of the class $[c_x]$ defined in remark \ref{rk:cp_obsruction}.
Upon taking $S^3 = {\mathbb R}^3 \cup \{\infty \}$ , we have $c_{\infty } = 0$, hence
$c_x  = \delta_{CE} (b)$, with 
$$
b =  -\int_{\gamma_{\infty}} \iota (v_1\wedge v_2) \nu
$$
($\gamma_{\infty}$ being a path connecting $ x $ to ${\infty}$, compare with the proof of \cite[Prop. 9.1]{Callies2016}.
In this case, the expression is meaningful in view of the assumed decay at infinity of our objects).
This is equivalent to the previous equation (\ref{eq:condition3_hccm}).
\end{proof}

We may also naturally ask the question of whether the above \momap is equivariant.
\begin{proposition}
	The map $(f)$ is {\rm not $G$-equivariant} in the sense of definition \ref{Def:EquivariantMomap}.
\end{proposition}
\begin{proof}
	To ascertain the infinitesimal {\it G-equivariance} of the above map $(f)$, one should check the validity of the formula
$$
{\mathcal L}_{\xi} f_1({\mathbf b}) = f_1 ([\xi, {\mathbf b}])
\qquad \forall \xi, {\mathbf b} \in {\mathfrak g} ~.
$$
Considering in particular $\xi = {\mathbf b}$, the above equation easily yields that
\begin{displaymath}
				f_1\left( \cancel{[\xi, \mathbf{b}]} \right)
				=
				{\mathcal L}_{\xi} f_1({\mathbf{b}})  = -{\mathcal L}_{\xi} \mathbf{A}^\flat =
					- d \iota_\xi \mathbf{A}^\flat = - d \langle \mathbf{A}, \mathbf{b} \rangle_g
				\qquad \forall \xi \in \mathfrak{g}
\end{displaymath}
where $\mathbf{A}$ is defined as within the previous proof.
Hence $\langle \mathbf{A}, \mathbf{b} \rangle_g$ ought to be a constant vanishing at infinity, \ie must be equal to zero.

	\begin{figure}[h!]
  \begin{center}
    \includegraphics[width=0.32\textwidth]{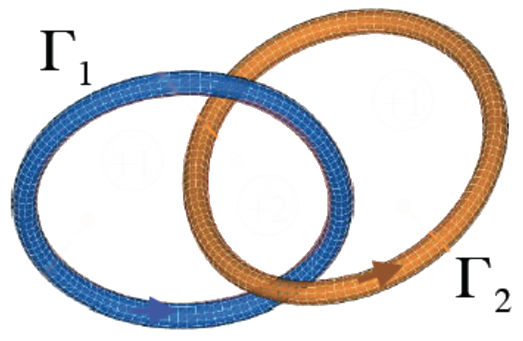}
  \end{center}
  \caption{Linked flux tubes.}
  \label{Fig:VortexLink}
	\end{figure}
%
Consider now a solenoidal field $\mathbf{b} = \curl(\mathbf{A})$ supported on a domain consisting of two disjoint, unknotted but linked, closed  \emph{flux tubes}
($\text{supp}(\mathbf{b}) = \Gamma_1 \cup \Gamma_2$, see figure \ref{Fig:VortexLink}).
One gets
		\begin{displaymath}
			\int \langle \mathbf{A}, \mathbf{b} \rangle_g = 2 n ~\Phi_1 \Phi_2 
			\qquad \text{with} \quad 
			\Phi_i = \int_{S_i} \mathbf{n} \cdot \mathbf{b}\, d\sigma
		\end{displaymath}
		where $n$ (\emph{Gauss linking number}, see theorem \ref{Thm:Moffatt-Ricca} below), is an integer different than zero, hence we get a contradiction. 
	%
\end{proof}
	\begin{remark}
Interpreting $\mathbf{A}$ in the proof of theorem as the velocity field of the fluid and $\mathbf{v}$ as the associated vorticity, the configuration given in figure \ref{Fig:VortexLink}) exemplifies the case of a configuration with non-zero {\it helicity}.
See \cite{Moffatt-Ricca92,BeSpe06,Spe06} and below for further elucidation of this train of concepts. 
Notice that the argument does not depend on the choice of the vector field $\mathbf{A}$ pertaining to $\mathbf{b}$.
The lack of $G$-equivariance is not surprising, since our construction involves Riemannian geometric features.
	\end{remark}

\note{In the symplectic case, $G$ connected implies $G$-equivariance. This does not hold true in the multisymplectic case. (ask Leyli!)}

\begin{remark}	
	Observe that the condition for the fields to vanish at infinity plays a crucial role when proving that the solution fulfils equation \eqref{eq_3Idro}.
	Relaxing this condition arises an obstruction to the existence of a \momap, however the component $f_1$ and $f_2$ constructed constitute a \emph{weak \momap} (see \cite[Def. 3.11]{Herman2017}\cite[Def. 3.8]{Herman2018} or \cite[Def. 1.16]{Mammadova2020}).
	\note{La definizione di weak moment map non è stata aggiunta da nessuna parte!
		Pensavo di metterla ma desisto per mancanza di tempo	
	}
\end{remark}

\begin{remark}\label{Rem:ConcreteSolutionBiotSavart}
	To give a completely explicit expression of the components $f_1$ and $f_2$ in theorem \ref{Thm:HydrodynamicalComoment}, it is required to solve two differential problems. Notice that the solution would not be unique, thus one can say that theorem \ref{Thm:HydrodynamicalComoment} yields a family of \momaps pertaining to the same action of the Lie algebra $\mathfrak{g}_0$ on $\R^3$.
	\\
	Given a vector field $\mathbf{b} \in \mathfrak{g}$ and defining $\mathbf{A}^\flat := - f_1(\mathbf{b})$ for a certain vector field $\mathbf{A}\in \mathfrak{X}(\R^3)$, the first component requires to solve the differential equation \eqref{Eq:MagnetoEq}.
	Such equation is the well-known equation of magnetostatic which admits a solution given by the so-called  Biot-Savart law
	\begin{displaymath}
					\tag{Biot-Savart law}
					\mathbf{A}(r) = \int\frac{\mathbf{b}(r')\times(\vec{r}-\vec{r}')}{|\vec{r}-\vec{r}'|^3}\d r'
	\end{displaymath}							
	defined up to a gradient \emph{(gauge freedom)}. It is customary to chose $\div(\mathbf{b})=0$ (Coulomb gauge).
	When $\mathbf{b}$ denotes a magnetic induction field, $\mathbf{A}$ is interpreted as the magnetic vector potential.
	In the context of hydrodynamics, $\mathbf{A}$ can be interpreted as a \emph{velocity field} and $\mathbf{b}$ as the corresponding vorticity.
	\\
	On the other hand, the computation of the second component requires to solve the Poisson equation \eqref{Eq:PoissonEq} with source given by $\rho:=-\delta \mu_2(\xi_1,\xi_2)$.
	A general solution is given by
	\begin{displaymath}
		\varphi(r) = \int\frac{-\rho(r')}{|\vec{r}-\vec{r}'|}\d r'
		~.
	\end{displaymath}
	Observe furthermore that $\rho$ can be explicitly given by
	\begin{displaymath}
		\begin{aligned}
		\rho :=~ \delta~ \mu_2(\xi_1,\xi_2) \equal{}&~ 
		\delta (f_1([\xi_1,\xi_2]) - \nu(\xi_1,\xi_2,\cdot)=
		\\
		\equal{}&~
		\delta \cdot \flat \cdot (\curl^{-1}([\xi_1,\xi_2]) - \xi_1 \times \xi_2)=
		\\
		\equal{}&~ -\div( \curl^{-1}([\xi_1,\xi_2]) - \xi_1 \times \xi_2)	
		=
		\\
		\equal{\eqref{eq:BrackAsCurl}}&~ 2~ \div(\xi_1 \times \xi_2) =
		\\
		\equal{}&~  2~\langle \xi_2, \curl \xi_1\rangle_g - \langle \xi_1, \curl \xi_2 \rangle_g
		\end{aligned}
	\end{displaymath}
	where $\times$ denotes the vector (cross) product on $\R^3$ and the last equality comes from another well-known result in vector analysis.
	\\
	For further informations see \cite[Appendix A1]{Miti2019a}.
\end{remark}


\subsubsection*{Hydrodynamical reinterpretation}
	Up to now, the main connection between the construction given in theorem \ref{Thm:HydrodynamicalComoment} and hydrodynamics consists in the fact that the considered Lie algebra $\mathfrak{g}$ can be interpreted as the phase space of a freely moving ideal fluid (see Subsection \ref{Sec:IdroPoisson}).
	
	Another important reason why the above construction is relevant in hydrodynamics is that the above \momap can be related, after transgression along the evaluation map ${\rm ev}: L{\mathbb R}^3 \times {\mathbb R} \ni (\gamma, t) \mapsto \gamma(t) \in {\mathbb R}^3$, to the hydrodynamical co-momentum map of Arnol'd and Marsden-Weinstein, defined on the Brylinski (infinite-dimensional) manifold $Y$ of oriented knots\cite{MW83}.
\begin{theorem}
	The above \Hcmm for $G\circlearrowright(\mathbb{R}^3,\nu)$ induces an ordinary (\ie "symplectic") co-moment map for $G\circlearrowright (L{\mathbb R}^3,\nu^{\ell})$, where $\nu^\ell$ denotes the trangression of the volume form along loops (see example \ref{Ex:TransgressionOnLoops}).
	Namely $(f) \colon \mathfrak{g} \to L_{\infty}(\mathbb{R}^3,\nu)$,
	given in theorem \ref{Thm:HydrodynamicalComoment}, 
	transgresses to
			\begin{displaymath}
				\begin{tikzcd}[column sep= small,row sep=0ex]
					\lambda \colon& \mathfrak{g}	\arrow[r]& C^\infty(L{\mathbb R}^3) \\
					& {\mathbf b}	\arrow[r, mapsto]
					& \displaystyle \biggr( \gamma \mapsto \lambda_{\mathbf b}(\gamma) = - \oint_{\gamma} f_1({\mathbf b})  \biggr)	
				\end{tikzcd}	
			\end{displaymath}
			that is a  moment map for the induced action $G$ on the pre-symplectic loop space $(L{\mathbb R}^3,\nu^{\ell})$
			(Smooth space in the sense of Brylinski, see remark \ref{Rem:BryLoopSpaces}).
\end{theorem}
\begin{proof}
Observe that the relevant piece of the \momap is  $f_1$ which, under transgression, becomes
$$
\mu_{\mathbf b} = - \int_{\gamma} B = -\lambda_{\mathbf b}
$$
\ie, up to sign, the {\it Rasetti-Regge current}  $\lambda_{\mathbf b}$ pertaining to ${\mathbf b} \in {\mathfrak g}$ (see definition \ref{def:RR-current}).
In particular, $\mu_{\mathbf b}$ is independent of the choice of $B$.
This is in accordance with the general result in \cite{Ryvkin2016} asserting that, roughly speaking, \momaps transgress to \momaps on loop (and even mapping) spaces.
{\it Actually, the ansatz for the $f_1$ term was precisely motivated by this phenomenon}.
\end{proof}

Note that it is natural to define a {\it "Poisson" bracket} on momenta, \ie on the image of $f_1$, via the expression:
\begin{equation}\label{Eq:MSPB}
\{ f_1({\mathbf b}), f_1({\mathbf c}) \} (\cdot):= \iota_{\mathbf c} \iota_{\mathbf b} \nu (\cdot)= \nu({\mathbf b}, {\mathbf c}, \cdot)
\end{equation}
satisfying the formula
\begin{equation}\label{eq:bracketsformula}
\{ f_1({\mathbf b}), f_1({\mathbf c}) \}  - f_1([{\mathbf b},{\mathbf c}]) = -df_2 ({\mathbf b} \wedge {\mathbf c}).
\end{equation}
that one gets after rewriting equation \eqref{eq_2Idro}. 
This construction will be employed below and in section \ref{Sec:MasseyMess}.

\subsection{A generalization to Riemannian manifolds}
 We ought to notice that a hydrodynamically flavoured \momap can be similarly construed also for an $(n+1)$-dimensional, connected, compact, orientable Riemannian manifold $(M,g)$, upon taking a Riemannian volume form $\nu$ as a multisymplectic form and again the group $G$ of volume-preserving diffeomorphisms as symmetry group. 

\begin{theorem}\label{Thm:RiemannGeneralization}
Let $(M,g)$ be a connected compact oriented Riemannian manifold of dimension $n+1$, $n\geq1$, with multisymplectic
form $\nu$ given by its Riemannian volume form, and such that the {\it de Rham}
cohomology groups $H_{dR}^{k}(M)$ vanish for $k=1,2,\dots n-1$ (one has necessarily $H_{dR}^{0} (M) = H_{dR}^{n+1} (M) = {\mathbb R}$). 
Let ${\mathfrak g}_0$ be the Lie subalgebra of ${\mathfrak g}$ consisting of divergence-free vector fields vanishing at a point $x_0 \in M$.
	\\
	Then there exists an associated family\footnote{They are a family in the sense that the "inverse" of $\Delta$ in equation \eqref{eq:RiemannCase} is a Green operator, hence it does not yield an unique solution}	
		 of \momaps for the infinitesimal action of $\mathfrak{g}_0$ on $(M,\nu)$ 
	given by the following compact formulae:
	\begin{equation}\label{eq:RiemannCase}
	f_1(\xi) := -\Delta^{-1} \delta (\iota_{\xi} \nu); \quad\quad f_k = \Delta^{-1} \delta \mu_k, \qquad k=2\dots n
~,
	\end{equation}
	where $\Delta$ denotes the Riemannian Laplace operator (see reminder \ref{Rem:HodgeCalculus}) and $\Delta^{-1}$ denotes the corresponding Green operator.
\end{theorem}
\begin{proof}
	We already noticed in the proof of theorem \ref{Thm:HydrodynamicalComoment}, how the defining equations of a \momap (see lemma \ref{Lem:ExplicitHCCM}) triggers a recursive construction starting from $f_1$, up to topological obstructions.
	Namely, we have a sequence of closed forms, which must be actually exact, together with the constraint $ f_n(\partial q) = (-1)^{\frac{(n+1)(n+2)}{2}}\nu(\xi_1,\dots\xi_{n+1})$, with $q= \xi_1 \wedge \dots \xi_{n+1}$.
	As before, the Riemannian structure allows us to connect elements of the Rogers $L_\infty$-algebra in different degrees:
			\begin{displaymath}
				\begin{tikzcd}[column sep=huge]
					&	
					& \Omega^{n+1}(M) \ar[dddddd,leftrightarrow,green,bend left=60,"\ast"]
					\\
					& \mathfrak{X}(M) \ar[r,red,"\alpha"] \ar[ddddr,blue,"\flat",shift left=1ex]\ar[ddddr,blue,leftarrow,"\sharp"',shift right=1ex]
					& \Omega^{n}(M) \ar[u,"d"']  \ar[dddd,leftrightarrow,green,bend left=60,"\ast"] 
					\\
					\mathfrak{g} \ar[ru,hookrightarrow,"v"] \ar[rr,purple,"f_1",crossing over]
					& 
					& \Omega^{n-1}_{(ham)}(M) \ar[u,"d"']  \ar[dd,leftrightarrow,green,bend left=60,"\ast"]
					\\
					\vdots \ar[u,"\partial"]
					& 
					& \vdots \ar[u,"d"'] 
					\\
					\bigwedge^{n-3}\mathfrak{g} \ar[u,"\partial"]\ar[rr,purple,"f_{n-3}",crossing over]
					& 
					& \Omega^2(M) \ar[u,"d"'] 
					\\
					\bigwedge^{n-2}\mathfrak{g} \ar[u,"\partial"]\ar[rr,purple,"f_{n-2}",crossing over]
					& 
					& \Omega^1(M) \ar[u,"d"'] 
					\\
					\bigwedge^{n-1}\mathfrak{g} \ar[u,"\partial"]\ar[rr,purple,"f_{n-1}",crossing over]
					& 
					& \Omega^0(M) \ar[u,"d"'] 
				\end{tikzcd}
			\end{displaymath}

In the present case, a natural candidate for the $(n$-$1)$-form $f_1$ can be readily manufactured via Hodge theory (see remark \ref{Rem:HodgeCalculus}):
\begin{equation}\label{eq:f1_hcmm_riemann}
f_1(\xi) := -\Delta^{-1} \delta (\iota_{\xi} \nu)
\end{equation}
(the direct generalization of the preceding case)
after imposing $\delta f_1({\xi}) = 0$ (the analogue of the Coulomb gauge condition), provided one can safely invert the
Hodge Laplacian $\Delta = d\delta + \delta d$.
According to the Hodge theorem (see \cite[\S 6]{Warner}) $H_{dR}^k(M) \cong \ker(\Delta\eval_{\Omega^k(M)})$, hence the latter holds true for any $1 \leq k \leq n-1$.	
\\
The topological assumption of vanishing of all the middle cohomology groups ensure that the entire procedure goes through unimpeded due to the formula
$$
df_k (\xi_1 \wedge\dots \wedge \xi_k) = \mu_k (\xi_1 \wedge\dots \wedge \xi_k), \qquad \qquad k=2,3,\dots n
~,
$$
where $\mu_k$ is the auxiliary form defined in remark \ref{Rem:TermMuByMauro} (see equation \eqref{eq:MauroMukForm}).
Hence, one can compactly state that:
\begin{equation}\label{eq:fk_hcmm_riemann}
	f_k = \Delta^{-1} \delta \mu_k, \qquad k=2\dots n,
\end{equation}
since
\begin{displaymath}
	\d ( \Delta f_k - \delta \mu_k) = \Delta ( \d f_k - \mu_k) = 0
	~.
\end{displaymath}
	One has finally to check that
	$$
	f_n (\partial (\xi_1 \wedge\dots \wedge \xi_{n+1})) =  -\varsigma(n+1) \iota (\xi_1 \wedge\dots \wedge \xi_{n+1})\nu
	$$
	but this is guaranteed by the vanishing of all fields at $x_0$, hence $c_{x_0} = 0$ and the obstruction class $[c_x]$ vanishes
	 (\cf remark \ref{rk:cp_obsruction}).
	%
	Note that one can also alter the above solution, given by equations \eqref{eq:f1_hcmm_riemann}  and \eqref{eq:fk_hcmm_riemann}, by addition of an exact form. 
	This freedom explain why the claimed solution has to be thought of as a family of \momaps.
\end{proof}

\begin{remark}
We notice that the above result holds, in particular, for  {\it homology spheres} such as, for instance, the celebrated Poincar\'e dodecahedral space\cite{Dror1973}. 
We point out that the case in which the intermediate homology groups are at most torsion (hence not detectable by de Rham techniques) is also encompassed: this is \eg the case of {\it lens spaces}. Notice that $G$-equivariance cannot be expected a priori. 
Also notice that one could restrict to the natural symmetry group provided by the isometries of $(M,g)$. See \eg \cite{zbMATH06448534} for a general discussion of topological constraints to existence and uniqueness of \momaps.

\end{remark}

\begin{remark}[Relation with weak homotopy moment maps]
	Recall that any \momap induces a {\it weak homotopy moment map}
(see \eg \cite{Herman2018,Mammadova2020}) by restriction on cycles in the Chevalley-Eilenberg chain complex (known also as \emph{Lie kernels}).
	Namely, if in equation \eqref{eq:fk_hcmm} we set $\partial p = 0$, we get
\begin{equation}
d f_k (p) + \varsigma(k) \iota(v_p) \omega = 0
\end{equation}
for $k=1,\dots,n+1$, which is the very property defining a weak homotopy moment map. 
For the time being we just 
%
	notice that, in theorems \ref{Thm:HydrodynamicalComoment} and \ref{Thm:RiemannGeneralization} above, the  vanishing of all fields at $x_0$ (that is infinity in the first case) is crucial for the existence of a \momap. 
	Upon relaxing such condition {\it we get a weak homotopy moment map which, in general, is not a \momap} since, for $k=n+1$, the value of the ensuing constant is not specified. 
\end{remark}

\subsection{Covariant phase space aspects}\label{Sec:CPSMess}
We are now going to propose a multisymplectic interpretation of the hydrodynamical bracket (see theorem \ref{Thm:HydroBracket}), which ties neatly with the topics discussed in previous sections, via {\it covariant phase space} ideas but without literally following the standard recipe, as we shall see shortly.

Let us briefly recall first the standard recipe for associating a multisymplectic manifold to any classical field system (see \cite{KS,Gimmsy1, Forger2005,Zuckerman87,Crnkovic,Ryvkin2018}).
Assume that one can encode all the configurations of a given field-theoretic (continuous) mechanical system as sections of a \emph{"configuration"} bundle $E\to M$ on a \emph{"parameter"} space $M$ (usually a globally hyperbolic spacetime) and assume that the dynamics of the system is governed by a $I$-order differential operator on $\Gamma(E)$ coming from a certain Lagrangian density $\mathcal{L}$.
Then, there is an associated multisymplectic manifold, called \emph{(Lagrangian) multiphase space}, given by the first jet space $J^1 E$ of the configuration bundle together with a $m$-plectic form $\omega_{\mathcal{L}}$, where $m$ is given by the dimension of the parameter space $M$, depending on the Lagrangian density.( See \cite[Eq. (3B.2)]{Gimmsy1}.)
\\
For instance, in the case of a classical (\ie non-relativistic) continuous medium in Galilean spacetime, one would have a tower of trivial bundles
	\begin{displaymath}
		\begin{tikzcd}[column sep = small,row sep =small]
			M \ar[r,phantom,"\cong"]
			& T \times \Sigma \ar[r,phantom,"\cong"]
			& \R^4 \ar[r,phantom,"\sim"] 
			&(t,x^i)
			\\
			E \ar[u] \ar[r,phantom,"\cong"]
			& M \times \Sigma \ar[r,phantom,"\cong"]
			& \R^{4+3} \ar[r,phantom,"\sim"] 
			&(t,x^i;y^j)
			\\
			J^1 E \ar[u] \ar[r,phantom,"\cong"]
			& E \times \Sigma' \ar[r,phantom,"\cong"]
			& \R^{7+4\cdot 3} \ar[r,phantom,"\sim"] 
			& (t,x^i,y^j; v^j_t, v^j_{, i})
		\end{tikzcd}
	\end{displaymath}
	where $T$ denotes the "absolute" time line, $\Sigma$ denotes the physical space and $\Sigma'$ is the space of generalized velocities. Hence $(J^1 E, \omega_{\mathcal{L}})$ is a $19$-dimensional $4$-plectic manifold.\footnote{In the particular case of an ideal fluid, the Lagrangian density would be given, employing Einstein summation convention, by $\mathcal{L}= \frac{1}{2} g_{i,j} v^i_t ~ v^j_t$.}

	Recall also the notion of \emph{covariant phase space} for a classical field theory, given as above, defined as the (usually $\infty$-dimensional) submanifold $\text{Sol}\subset\Gamma(E)$ of all configurations consisting of only the solutions of the motion equation. 
	It is well-known how to induce a presymplectic structure on $\text{Sol}$ from the multisymplectic structure $\omega_{\mathcal{L}}$ on $J^1 E$ (see for example \cite{Forger2005,Helein2011c}).
	Namely one has that for any $\sigma \in \text{Sol}$ and $\frac{\delta\phi}{\delta x^i}\in T_{\sigma}\text{Sol}$, (vertical) variation of the solution $\sigma$, the $2$-form given by
	\begin{displaymath}
		\Omega_\sigma \left(\frac{\delta\phi}{\delta x^i}, \frac{\delta\phi'}{\delta x^i}\right)
		= \int_{\Gamma} (j^1 \sigma)^\ast 
		\left( 
		\iota_{\frac{\delta\phi}{\delta x^i}}\iota_{\frac{\delta\phi'}{\delta x^i}}
		\omega_{\mathcal{L}}\right)
\end{displaymath}		
	is closed.

	Now we want to show how the Poisson structure given by equation \eqref{Eq:HydroBracket} can be obtained as a covariant phase space bracket obtained from a suitably reduced version of the multisymplectic manifold associated the ideal fluid.
	Observe that any divergence-free vector field  can be viewed as an initial condition ${\mathbf v} (x,0)$ for the (volume-preserving) Euler evolution (at least for small times, but as we previously said, we do not insist on refined analytical nuances) ${\mathbf v} (x,t)$, yielding a section of $E$. 

We set 
	\begin{displaymath}
		{\mathcal J}^1 {\mathbf v}  := {\mathbf w} \quad (:= {\rm curl}\, {\mathbf v} )
	\end{displaymath}
which has to be understood as the natural "covariant" jetification of the section ${\mathbf v}$.
This object is "covariant" in the sense that, in contrast with the standard jetification $j^1$ of a section (called \emph{first jet prolongation} in \cite[eq. (2A.4)]{Gimmsy1}), does not yield an object depending on the choice of coordinates in $E$.
Also, since  the section of the first jet bundle $J^1 E\to E$ can be interpreted as Ehresmann connections\cite[Rem. 1]{Gimmsy1}, this allows to look at ${\mathbf v} $ as the vector space counterpart of a connection $1$-form.
\\	
Using the 3-volume form $\nu$, orienting fibres (notice that, when viewed on $E$, it is only pre-$2$-plectic, namely closed but degenerate), we can rewrite the hydrodynamical bracket, 
mimicking \cite{Forger2005}, 
as 
\begin{equation}\label{Eq:FormerStar}
	\begin{aligned}
	\{ F, G \} ([{\mathbf v}]) 
	=&~ 
	\int_{\Sigma = {\mathbb R}^3} \left\langle {\mathbf w}, \,  \frac{\delta F}{\delta {\mathbf v}} \times \frac{\delta G}{\delta {\mathbf v}}\right\rangle \,d^3x
	\\
	=&~
	\int_{\Sigma = {\mathbb R}^3} \nu ({\mathcal J}^1{\mathbf v}, \, \frac{\delta F}{\delta {\mathbf v}}, \frac{\delta G}{\delta {\mathbf v}}) \,d^3x
	\end{aligned}
\end{equation}
since the variations $\frac{\delta F}{\delta {\mathbf v}}$ and $\frac{\delta G}{\delta {\mathbf v}}$ are {\it vertical} and divergence-free: $\delta F/ \delta {\mathbf v} = {\rm curl} \, (\delta F/ \delta {\mathbf w})$.
 Taking again ${\mathbf b} = {\rm curl}\, {\mathbf B}$ et cetera,  and setting finally $F = \lambda_{\bullet}$ (see \eg in particular \cite{Pe-Spe92,Spera16}) we see that the expression \eqref{Eq:FormerStar} can be manipulated to yield the expressive layout (with slight abuse of language)
\begin{displaymath}
	\begin{aligned}
		\{ F, G \} ([{\mathbf v}]) 
		=&~ 
		\int_{\Sigma} ({{\mathcal J}^1}^*\nu) ({\mathbf v},{\mathbf B},{\mathbf C}) \,d^3x =
		\\
		=&~ \int_{\Sigma}\nu ( {\mathcal J}^1{\mathbf v}, {\mathcal J}^1{\mathbf B}, {\mathcal J}^1{\mathbf C}) \,d^3x =
		\\
		=&~ \int_{\Sigma}\nu ( {\mathbf w}, {\mathbf b}, {\mathbf c}) \,d^3x
	\end{aligned}
\end{displaymath}
(in full adherence with the discussion carried out in Section \ref{Sec:IdroPoisson}).
The same portrait can be depicted, {\it mutatis mutandis}, for the singular case.
Ultimately, we reached the following conclusion: 
\begin{theorem}
(i) The Poisson manifold ${\mathfrak g}^*$, introduced in equation \eqref{eq:gastPoisson} and in theorem \ref{Thm:HydroBracket}, can be naturally  interpreted as a (generalized) covariant phase space pertaining to the volume-preserving Euler evolution. 
The latter indeed preserves the symplectic leaves of
${\mathfrak g}^*$ given by the $G$-coadjoint orbits ${\mathcal O}_{[{\mathbf v}]}$.\par
(ii) The above construction reproduces the symplectic structure of the Brylinski manifold ${Y}$, see remark \ref{Rem:BryLoopSpaces}, upon taking singular vorticities concentrated
on a smooth oriented knot. 
The covariant phase space picture is fully retrieved upon passing to
a 2-dimensional space-time $S^1 \times {\mathbb R}\rightsquigarrow (\lambda, t) $, with $\lambda \in S^1 \equiv \Sigma$ being a knot parameter (and staying of course with the same $\nu$). \par
\end{theorem}

\begin{remark}
\begin{enumerate}
	\item We stress the fact that we did not literally follow the standard "multisymplectic to covariant" recipe developed in \cite{Forger2005}. In fact the multisymplectic manifold we consider is not the one prescribed by \cite{Gimmsy1} since we directly took the standard volume form 
$\nu$ on ${\mathbb R}^3$ as a $2$-plectic structure (or pre-$2$-plectic when pulled back to $E$), \cf \cite{CatIbort}.
This neatly matches  Brylinski's theory  and fits with the stance long advocated, among others, by Rasetti and Regge and Goldin (see \eg \cite{Rasetti-Regge75,Goldin71,Goldin87,Goldin12}, and \cite{Spera16} as well) pinpointing the special and ubiquitous role played by the group of orientation-preserving diffeomorphism $G = \sDiff({\mathbb R}^3)$.
Another motivation for considering $\nu$ is its pivotal role in the formulation of conservation theorems (see \cite{Ryvkin2016}). We shall pursue this aspect in what follows.

	\item
	In line with the preceding remark, notice that the above portrait can, in principle, be generalized to any volume form (on an orientable manifold), with its attached group $G$. 
	The covariant phase space picture should basically persist in the sense that one might construct, in greater generality, an $n$-plectic structure out of an $(n+1)$-plectic one via an expression akin to equation \eqref{Eq:FormerStar}. 
	The (non) $G$-equivariance issue should be relevant in this context.
\end{enumerate}
\end{remark}

\section{A Hamiltonian $1$-form for links}\label{Sec:Ham1FormLinks}
This section wants to show how it is possible to build a bridge between knot theory and multisymplectic geometry, exploiting the close connection between them and hydrodynamics.
Vortex theory, together with the ubiquitous role of the group of volume-preserving diffeomorphisms, will be the cornerstone of this relationship.
	The basic and quite natural idea is {\it to associate to a knot (or link) a perfect fluid whose vorticity in concentrated thereon}
(\cf the preceding discussion on the Brylinski manifold).
\begin{displaymath}
	\begin{tikzpicture}[thick,node distance = 8em, auto]
    \node [->,rectangle, draw, fill=blue!20,
    text width=7em, text centered, rounded corners, minimum height=4em] (A)
    {Multisymplectic Geometry};
    \node[inner sep=1em,minimum size=4em, text width=4em,right of=A] (k) {}; 
    \node [rectangle,above of=k, draw,
    text width=7em,text centered, rounded corners, minimum height=4em] (B)
    {Hydrodynamics};
    \node [rectangle,right of=k, draw, fill=blue!20, 
    text width=7em, text centered, rounded corners, minimum height=4em] (C)
    {Knot theory};
    
		\path[every node/.style={font=\sffamily\small}]
    	(A) edge[bend left=45,-Latex] node [left,text width=5em,align=center] {Geometric\\ Fluid\\ Mechanics} (B)
    	(B) edge[bend left=45,-Latex] node [right,text width=5em,align=center] {Vortex\\ dynamics} (C) ;
		{\path[every node/.style={font=\sffamily\small}]
		    (A) edge[bend right=45,dashed,-Latex,gray] node [below,text width=10em] {Geometric Mechanics \\ of classical strings} (C);}
	\end{tikzpicture}
\end{displaymath}

	Let us recall that Knot theory is a rich and ramified theory, born as a pure mathematical curiosity and naturally framed in the field of low-dimensional topology, but yielding innumerable implications in the most disparate areas of both pure and applied mathematics.
	We do not venture to review all the multifaceted aspects of this broad topic, and, for what follows, it is enough to recall just a few basic ideas regarding knots and links.
	(As general references on knot theory we quote, among others, \cite{Rolfsen,Lickorish1997},
together with \cite{Bott-Tu82} for the algebraic-topological tools employed below.)
	\begin{reminder}[Geometry of $n$-Links]\label{Rem:ReminderOnKnots}
		As already anticipated in remark \ref{Rem:BryLoopSpaces}, a (smooth) \emph{knot} can be formalized as a compact submanifold of codimension $2$ embedded in $\R^3$. 
		A slight generalization of this concept is given by the notion of \emph{links}:
		\begin{itemize}
			\item We call \emph{$n$-link} any  smooth embedding
		$\gamma: {\coprod_{i=1}^n S^1 } \to \mathbb{R}^3 $		
		of $n$ disjoint copies of the circle into the Euclidean space $\R^3$.
		
			\item 
				Links are studied modulo \emph{ambient isotopies}. 
				Specifically, given two embeddings $f,g: {\coprod_{i=1}^n S^1 }\to \R^3$, an ambient isotopy between $f$ and $g$ is a smooth isotopy $h : \R^3 \times [0,1] \to \mathbb{R}^3 $ of the "ambient" space,
				\ie $\hat{h}(t) : \R^3 \to \mathbb{R}^3$ is a diffeomorphism for any $t\in [0,1]$,
				such that $\hat{h_1}\circ f = g $.
				In layman terms, an ambient isotopy describes a continuous deformation of the parametrized link $\hat{h}(0)$ into $\hat{h}(1)$. 
				Hence, to study $n$-links modulo ambient isotopies, means to study the equivalence classes of those links that can be transformed continuously one into the other (in particular without cuts or other surgeries).
				\\
				It is then implicit that the Brylinski manifold $Y$ (see remark \ref{Rem:BryLoopSpaces}), and its obvious generalization to $n$-links, must be composed of several connected components, one for each equivalence class of links. 
				Clearly all these classes are invariant with respect to volume-preserving diffeomorphisms.
			
			\item The ultimate goal of knot theory is to completely classify all non-equivalent (by ambient isotopies) classes of $n$-links. 
			Note that it is not simply an instance of studying the topology of loops but the ambient space takes a significant role.
			\\
			To date, such a classification is still missing. 
			Nonetheless, numerous "partial" criteria to distinguish non-equivalent knots are known. 
			These are mostly achieved introducing the so-called \emph{knot polynomials} (see section \ref{sec:prelim} below).
		\end{itemize}
	\end{reminder}
	
	Let us now specialize the considerations in Subsection \ref{Sec:IdroMoMap} to the  case of links.
	Namely, we show how one can naturally associate to any link in $\R^3$ an \emph{Hamiltonian $1$-form} that in particular turns out to be conserved under the Euler evolution, \ie under the infinitesimal action of the Lie algebra $\mathfrak{g}=\sDiff(\R^3)$ on $\R^3$ (action defined in \ref{Sec:IdroPoisson}).	
	Such construction can be achieved introducing the concept of \emph{Poincar\'e} dual of a given closed, oriented submanifold.
	To make sense of the latter, we briefly recall the notion of \emph{de Rham currents}:
\begin{reminder}[de Rham currents (\cite{dR})]\label{Rem:ReminderDeRhamCurrents}
	It is possible to mimic the basic definition of "distribution" or "generalized function", originally introduced on $\R^n$, on any smooth manifold $M$ by mean of differential forms:
	\begin{itemize}
		\item We call \emph{(de Rham) $k$-current} any (sequentially) continuous\footnote{Continuity is meant in the "sequential" sense.
			More explicitly, the functional $T$ is continuous if, for any given sequence $\omega_{k}$ of smooth forms, all supported in the same compact set,  such that all derivatives of all their coefficients tend uniformly to $0$ when 
		$k$ tends to infinity, $T(\omega_{k})$ tends to 0.
		\\
		Recall the definition of \emph{support} of a differential form as the closure of the non-vanishing locus of $\omega$. 
		Namely
		$
			\text{supp}(\omega) := \overline{\lbrace p \in M \; \vert \: \omega_p \neq 0 \rbrace}
		$.
		}
		linear functional $\eta$ from the space $\Omega^k_c(M)$ of compactly supported $k$-forms to $\R$.
		We denote by $\mathcal{D}_k(M)$ the vector space of all $k$-current $M$.
		Clearly one has $\mathcal{D}_k(M)\cong\left(\Omega^k_c (M) \right)^\ast$ understanding the space on the right-hand side as the \emph{topological dual} of $\Omega^k(M)$.
		\item
			We call \emph{annihilation set} of $\eta\in \mathcal{D}$, an open subset $A \subset M$ such that 	$\langle \eta, \varphi \rangle = 0$ for any $\varphi\in\Omega_c(M)$ compactly supported in $A$.
			(here $\langle \eta, \varphi \rangle$ denotes evaluation of the given current $\eta$ on the "test" form $\varphi\in \Omega_c^k(M)$.)
			\\
			We call \emph{support} of a given current $\eta\in \mathcal{D}$ the complement of the union of all open annihilation sets of $T$.
		\item 
			There exists an analogue of a \emph{regular distributions} given by the following mapping:
			\begin{displaymath}
				\morphism{D}
				{\Omega^k(M)}
				{\mathcal{D}^{n-k}(M)}
				{\alpha}
				{\displaystyle\left(	D_\alpha : \varphi \mapsto \int_M \alpha\wedge \varphi
				\right)}~.
			\end{displaymath}
		\item The de Rham differential can be extended to currents by duality
			\begin{displaymath}
				\morphism{\partial}
				{\mathcal{D}_k(M)}
				{\mathcal{D}_{k-1}(M)}
				{\eta}
				{\displaystyle\left(	\partial \eta : \varphi \mapsto (-)^k\int_M \eta\wedge \d \varphi
				\right)}~.
			\end{displaymath}
			\ie $\langle \partial \eta, \blank \rangle = (-)^{k} \langle \eta, \text{d}\blank \rangle$.
		Hence, de Rham currents build up a chain complex which is dual (modulo a sign) to the de Rham cochain complex.
		In a similar way, one can make sense of the wedge of de Rham currents.
	\end{itemize}
\end{reminder}

To any given a compact, oriented embedded $k$-dimensional submanifold one can associates a unique de Rham current
\begin{definition}[Poincar\'e dual]
	Given a smooth embedding $\left( i : 	\Sigma \hookrightarrow M \right) \in \text{Emb}_c(k)$ of a compact, oriented, $k$-dimensional submanifold $\Sigma$ of the $n$-dimensional manifold $M$, we call \emph{Poincar\'e dual} of $\Sigma$ the de Rham current $D_\Sigma$ uniquely defined by the following equation
	\begin{displaymath}
		\langle D_\Sigma, \omega\rangle = \int_\Sigma i^\ast (\omega) \qquad 
		\forall \omega \in \Omega^k(M)
		~.
	\end{displaymath}
\end{definition}
	These can be understood as \emph{generalized} differential $(n-k)$-forms concentrated on $\Sigma$. As such, the Poincar\'e dual $D_\Sigma$ is the analogue of a Dirac delta function localized on the submanifold $\Sigma$, in that sense are "non regular" or "singular".
	One can then be interested in a \emph{regular approximation} of such singular behaviour:
	\begin{definition}[Smooth (regular) Poincar\'e dual]\label{Def:RegularPoinDual}
			We call \emph{Smooth Poincar\'e} dual of an embedding $\left( i : 	\Sigma \hookrightarrow M \right) \in \text{Emb}_c(k)$ as above, any (smooth) differential form $\eta_\Sigma \in \Omega^k_c$ supported on a tubular neighbourhood $T$ of $\Sigma$ s.t.
		 \begin{displaymath}
				\langle D_{\eta_\Sigma},\omega\rangle \equiv \int_M \omega \wedge \eta_\Sigma =
				\int_T i^\ast \omega 
				~.
		 \end{displaymath}
	\end{definition}
	Clearly, smooth Poincar\'e duals 	are not uniquely defined. 
	The distribution associated to the differential form $\eta_\Sigma$ represents a possible regular approximation of $D_\Sigma$ in the sense that $\langle D_{\eta_\Sigma},\omega\rangle \sim 
				\langle D_{\Sigma},\omega\rangle $.
		\\
		Poincar\'e duals satisfy the following properties:

		\begin{lemma}\label{Lem:PropertiesPoinDuals}
			Poincar\'e duals satisfy the following properties for any given compact oriented submanifold $\Sigma_i$:
			\begin{align}
				& \partial D_\Sigma = (-)^{k-1} D_{\partial \Sigma} \label{Eq:PoinProp1}\\
				& D_{\Sigma_1} \wedge D_{\Sigma_2} = D_{\Sigma_1 \cap \Sigma_2} \label{Eq:PoinProp2}
			\end{align}
		\end{lemma}
		\begin{proof}
			The proof of \eqref{Eq:PoinProp1} is straightforward:
			
		\begin{displaymath}
			\mathclap{
			(-)^{k-1} \langle \partial D_\Sigma , \omega \rangle = 
			\langle D_\Sigma, d\omega\rangle =
			\int_\Sigma i^\ast d \omega = 
			 \int_\Sigma d i^\ast \omega = 
			 \int_{\partial \Sigma} i^\ast \omega =
			\langle D_{\partial \Sigma}, \omega \rangle
			~.
			}
		\end{displaymath}
		Claim \eqref{Eq:PoinProp2} is better understood by means of smooth Poincar\'e duals:
		\begin{displaymath}
			\text{supp}(\eta_1 \wedge \eta_2) \subset
			\text{supp}(\eta_1) \cap \text{supp}(\eta_2) \subset
			T_{\Sigma_1 \cap \Sigma_2} \quad\Rightarrow\quad
			\eta_1 \wedge \eta_2 = \eta_{\Sigma_1 \cap \Sigma_2}	
			~.
		\end{displaymath}				
		For instance, take as $\Sigma_1$ the $z$-line in $\mathbb{R}^3$, $\eta_1 = \delta_{\{x=y=0\}}	dx \wedge dy$, and as $\Sigma_2$ the $xy$-plane, $\eta_2= \delta_{\{z=0\}} dz$.
		You get $\eta_1 \wedge \eta_2 = \delta_{\{x=y=z=0\}} dx \wedge dy \wedge dz = \eta_{\Sigma_1 \cap \Sigma_2}$.
		\end{proof}
	\begin{remark}[	Poincar\'e duals as Thom classes]
		Generally speaking, in a more algebraic topological flavour, one can make sense of the 
		Poincar\'e dual of a $k$-dimensional closed oriented submanifold $\Sigma$ of an n-dimensional manifold $M$ as a cohomology class $[\eta_\Sigma] \in H^{n-k}(M)$ characterized by the property:
		$$
		\int_M \omega \wedge \eta_\Sigma = \int_\Sigma i^* \omega
		$$
		for any {\it closed, compactly supported $k$-form} $\omega$ on $M$ 
		($i: \Sigma \hookrightarrow M$ being the inclusion map).  	 
		\\
		One can therefore see a (regular) Poincar\'e dual as 
		$k$-form localized in a  cross-section of a  suitable tubular neighbourhood around $\Sigma$ \ie a \emph{Thom class} (see \cite{Bott-Tu82}).
		In \cite{Miti2019a} the authors indifferently view Poincar\'e duals as genuine forms or currents in the sense of de Rham.		
	\end{remark}

\begin{figure}[h!]
		\centering
		\includegraphics[width=\textwidth]{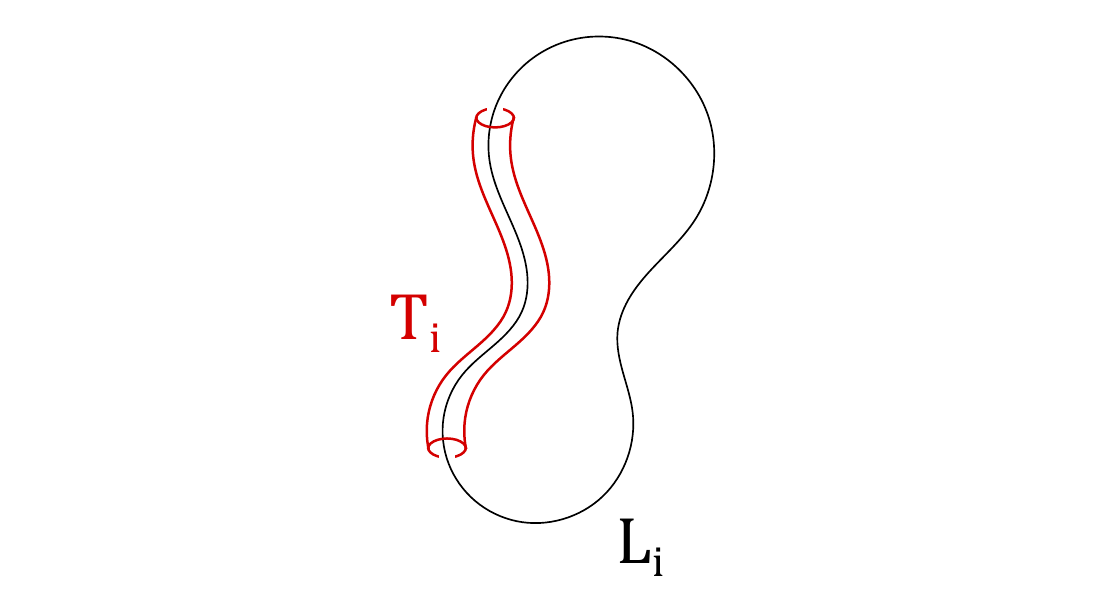}
		\caption{Tubular neighbourhood for the $i$-th component of $L$.}
		\label{fig: tubes}
	\end{figure}

\begin{figure}[h!]
		\centering
		\includegraphics[width=\textwidth]{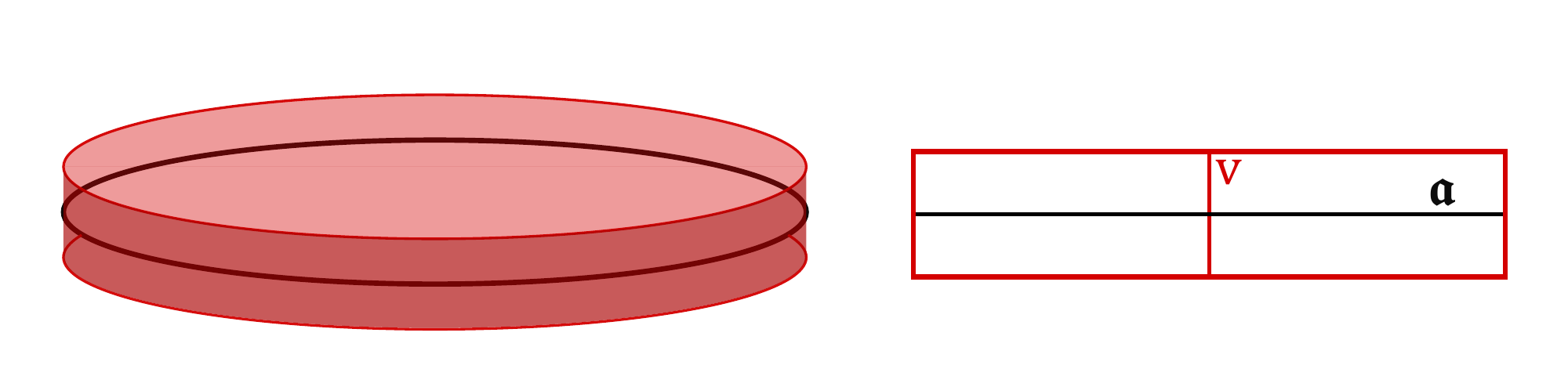}
		\caption{Poincar\'e duals}
		\label{fig: thom}
	\end{figure}

	Building on \cite{Pe-Spe02,BeSpe06,Spe06}, let $ L = \coprod_{i=1}^n L_i$ be an oriented link in ${\mathbb R}^3$ with components $L_i$, $i=1,\dots,n$ required to be  {\it trivial} knots.
	Let us choose a suitable tubular neighbourhood $T_i$ around any component $L_i$ (see Figure \ref{fig: tubes}).
	We introduce the following differential forms:
 \begin{definition}[Vorticity $2$-form]\label{Def:Vorticity2Form}
	 	We call \emph{vorticity $2$-form} pertaining to $L$ the smooth differential form
	 	\begin{displaymath}
			\omega_{L} := \sum_{i=1}^n \omega_{L_i} \in \Omega^2(\R^3)
	 	\end{displaymath}
		with $\omega_{ L_i}$ denoting a smooth Poincar\'e dual associated to $L_i$ with support given by $T_i$.
 \end{definition}
	\begin{definition}[Velocity $1$-form]\label{Def:Velocity1Form}
		We call \emph{velocity $1$-form} pertaining to $L$ the smooth differential form
		\begin{displaymath}
			v_{ L} = \sum_{i=1}^n v_{L_i} \in \Omega^1(\R^3)
		\end{displaymath}
		with $v_{L_i}:= \omega_{{\mathfrak a}_i}$ denoting a smooth Poincar\'e dual associated to a disc ${\mathfrak a}_i$ 	bounded by 	$L_i$ (Seifert surface in the knot theory jargon, see Figure \ref{fig: thom}). 
	\end{definition}		
 	\begin{remark}[Non-uniqueness of  $v_L$ and $\omega_L$]\label{Rem:KnotFormsNotUnique}
		Notice that $\omega_L$ and $v_L$ 
 are not uniquely determined but rather heavily depend on the choice of the concrete differential form used to represent the regular Poincar\'e dual of each component, and corresponding Seifert surface, of the links.
	 (Singular de Rham currents would be a correct language for restoring uniqueness.)
	 	Although this feature may appear unpleasant at first glance, it actually fits neatly the philosophy of studying knots and links regardless their parametrization and modulo ambient isotopies.
	 	In fact, the property of being "knotted" of a certain link configuration is indifferently detected whether it is considered as a singular object or as something "smeared" in a region close to a possible spatial displacement of the link.
		An exemplification of this phenomenon will be discussed in the next subsection. 
	\end{remark}

	\begin{theorem}\label{thm:VorticityFormExact}
		The vorticity $2$-form $\omega_L$ associated to the link $L$ is an exact form.
		The velocity $1$-form $\omega_L$ is a primitive  of $\omega_L$ and is an Hamiltonian form for the $2$-plectic manifold $(\R^3,\nu)$ in the sense of definition 	\ref{Def:Hamiltonianform}.
	\end{theorem}
\begin{proof}
	According to equation \eqref{Eq:PoinProp1} in Lemma \ref{Lem:PropertiesPoinDuals}, the vorticity form $\omega_L$ is closed since the boundary of $L$ is clearly empty. In particular, $\omega_L$ is exact by virtue of the Poincar\'e lemma.
	\\
	By the same lemma, and recalling that $\partial\mathfrak{a}_i = L_i$, one gets for each component of the link, that
	\begin{displaymath}
		\d v_{L_i} = \d \eta_{\mathfrak{a}_i} = \eta_{\partial \mathfrak{a}_i} = \eta_{L_i} = \omega_{L_i}
	\end{displaymath}
	where $\eta_{\Sigma}$ denotes a regular Poincar\'e duals of $\Sigma$.
	\\
	Finally, notice that, for each component $L_i$, the Hamiltonian vector field
	of $v_{L_i}$ is given by 	
	 $\xi_{L_i}=-\alpha^{-1}(\omega_{L_i})$(via the map $\alpha$ of Section \ref{Sec:IdroMoMap}). 
	 Explicitly, one has (setting $\xi_L = \sum_{i=1}^n \xi_{L_i}$)
\begin{equation}\label{eq:poincaredual}
dv_L + \iota_{\xi_L} \nu = 0 ~.
\end{equation}
\end{proof}

	\begin{remark}[Cohomological interpretation of $v_L$]\label{Rem:CohoIntLinkForm}
		Recall that there is a canonical cohomology theory associated to a $n$-link $L$
		given by the cohomology of $S^3 \setminus L$.
		(Note that there is no harm in considering the compactification $S^3$ instead of the standard Euclidean space since the image of any link configuration would be in any case sits in a compact subset of $\R^3$.)
		For the sake of clarity, the (de Rham) cohomology and relative homology groups of $S^3 \setminus L$ with real coefficients, reads as follows
		\begin{displaymath}
			\begin{aligned}
				H^0 (S^3 \setminus L) &\cong H_3(S^3, L) \cong {\mathbb R}
				\\
				H^1 (S^3 \setminus L) &\cong H_2(S^3, L) \cong {\mathbb R}^{n} 
				\\
				H^2 (S^3 \setminus L) &\cong H_1(S^3, L) \cong {\mathbb R}^{n-1}
				\\
				H^3 (S^3 \setminus L) &\cong H_0(S^3, L) \cong 0
			\end{aligned}
		\end{displaymath}
		which can be computed iterating the Mayer-Vietoris on $\R^3$ deprived of a line.\par
		The de Rham classes of the forms $v_{L_i}$ generate the cohomology group $H^1(S^3 \setminus L, {\mathbb R})$ (or, better, that of $S^3 \setminus T$, with $ T = \cup_{i=1}^n T_i$). 
		Their homological counterparts are given by the (classes of) the discs ${\mathfrak a}_i$. 
		One can also interpret the other groups:
in particular, elements in $H_1(S^3, L) $ can be represented by classes $[\gamma_{ij}]$ of (smooth) paths $\gamma_{ij}$ connecting two components $L_i$ and $L_j$, subject to the relation 
$[\gamma_{ij}] + [\gamma_{jk}] = [\gamma_{ik}]$. (See \cite{Spe06})
	\end{remark}

\begin{remark}
	Inspection of the very geometry of Poincar\'e duality shows that the velocity $1$-forms $v_i$ correspond (upon approximation of the associated Euler equation) to the so-called LIA (Linear Induction Approximation) or  {\it binormal evolution}
of the ``vortex ring" $L_i$ 
(``orthogonal" to the discs ${\mathfrak a}_i$ - an easy depiction, \cf figure \ref{fig: thom}), see \cite{Khe} for more information.
	Formula (\ref{eq:poincaredual}) will be the prototype for the calculations in section \ref{Sec:MasseyMess}.
\end{remark}

\subsection{Relation to Gauss linking number}\label{subsec:ReltoGauss}
	We want now to show how the previous constructions can determine quantities relevant to the classification of a link despite their heavy reliance on arbitrary choices.
	
	Consider a link $L$ embedded in $(\R^3,\nu)$. Let be $v_L$ and $\omega_L$ the velocity and vorticity forms defined in \ref{Def:Velocity1Form} and \ref{Def:Vorticity2Form} respectively.
	One can therefore introduce the {\it Chern-Simons (helicity) } $3$-form:
	\begin{equation}\label{Eq:CSHeli}
		CS({L}) :=  v_{L} \wedge  \omega_{ L}
	~.
	\end{equation}
	Integration of $CS({L})$ over ${\mathbb R}^3$ (or $S^3$ in case of compactification of the ambient space) yields a number ${\mathcal H}(L) $ called  the {\it helicity} of $L$.
	The naming is clearly not incidental. Interpreting a $L$ as an ideal fluid configuration with singular vorticity concentrated in a filament given by $L$, one recovers the \emph{helicity} of the fluid as introduced in remark \ref{Rem:Helicity}.
	
	It is a celebrated result that the quantity $\mathcal{H}(L)$ is indeed an integer, depending on the knot and not on its possible parametrization, measuring the mutual knotting of two generic flow lines (we point to \cite{Moffatt-Ricca92, Arn-Khe,Pe-Spe89,Pe-Spe92,Spe06} for a more extensive discussion) 
	\begin{theorem}[Moffatt-Ricca \cite{Moffatt-Ricca92}]\label{Thm:Moffatt-Ricca}
		Given a link $L$ in $\R^3$, one has
		\begin{displaymath}
			\int_{S^3} CS({L}) =: {\mathcal H}(L)  = \sum_{i,j=1}^n \ell(i,j) ~\in \mathbb{N},
		\end{displaymath}	
		where $\ell(i,j)\in\mathbb{Z}$ are quantities that can be algorithmically computed on indented diagrams\footnote{
A useful way to visualise and manipulate knots is to project the knot onto a plane.
In order to be able to recreate the original knot, one must distinguish between the over-strand and the under-strand at every crossing. 
			This is usually done by creating a break in the strand going underneath, called "indentation". 
The resulting diagram is an immersed plane curve with the additional data of which strand is over and under each crossing.
			}
		 (see \eg \cite{Rolfsen,Moffatt-Ricca92,Ricca2011,Spe06}). 		
		Namely:
		\begin{itemize}
			\item when $i\neq j$, $\ell(i,j) =\ell(j,i)= \ell(L_i,L_j)$ is the {\it Gauss linking number}\footnote{Informally, it represents the number of time that  the curve $L_i$ winds around the curve $L_j$.
			Operatively, it can be computed by summing all the signed crossings of the indented diagram (see \cite{Rolfsen,Lickorish1997}).
			} 
			pertaining to the components $L_i$ and $L_j$; 
			\item when $i=j$, $\ell(j,j)$ is equal to the Gauss linking number $\ell(L_j, L_j^{\prime})$ with $L_j^{\prime}$ being a section of the normal bundle of $L_j$ (think $L_j^\prime$ as an another loop winding around $L_j$, arbitrarily close without intersecting it).
		\end{itemize}			
	\end{theorem}
	\begin{proof}
	(\textit{Sketch})\\	
		Choosing a parametrization $\mathbf{r}_i$ (in standard coordinates) for each component $L_i$, one gets the Gauss' linking integral formula\cite{Ricca2011} 
		\begin{displaymath}
			{\mathcal H}(L)  = 
			\sum_{i,j=1}^n
			\,\frac{1}{4\pi}
			\oint_{\gamma_i}\oint_{\gamma_j}
			\frac{\mathbf{r}_i - \mathbf{r}_j}{|\mathbf{r}_i - \mathbf{r}_j|^3}
			\cdot (d\mathbf{r}_i \times d\mathbf{r}_j)
			~.
		\end{displaymath}
		\begin{figure}[h!]
			\vspace{-1em}
		 	 \begin{center}
		 	   \includegraphics[width=0.3\textwidth]{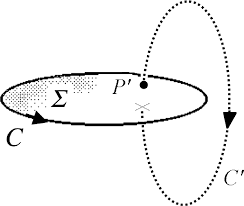}
			  \end{center}
			  \caption{Hopf link, component $C'$ crossing the Seifert surface of component $C$. \cite{Ricca2011}.}
			  \label{Fig:HopfRicca}		
		\end{figure}
		Solving this integral could seem a daunting task at first; however, one could convince oneself that it yields an integer by using the notion of Poincar\'e duals in a simple example.
		\\
		Let $L=C\coprod C'$ be a Hopf Link as in figure \ref{Fig:HopfRicca}, consider the de Rham current 
		\begin{displaymath}
			\mathbf{H}(L) = (D_C\wedge D_\Sigma + D_{C'}\wedge D_{\Sigma'})
			+ (D_C \wedge D_{\Sigma'} + D_{C'}\wedge D_{\Sigma'})
		\end{displaymath}
		where $D_C$ denotes the (singular) Poincar\'e dual of $C$ and $D_{\Sigma}$ denotes the (singular) Poincar\'e dual of a certain Seifert surface $\Sigma$ of $C$. 
		Notice that this can be seen as a singular version of the Chern-Simons $3$-form given in equation \eqref{Eq:CSHeli}.
		According to Lemma \ref{Lem:PropertiesPoinDuals}, we have
			\begin{displaymath}
				\mathbf{H}(L) = D_{P'} + D_{P}
			\end{displaymath}
		where $D_P$ and $D_P'$ are singular currents localized in the intersection point of $C$ with $\Sigma'$ and viceversa.
		Hence, the integration of $\int \mathbf{H}(L)$ simply counts the times that a component cross another Seifert Surface \ie the Gauss linking number of the link.
	\end{proof}
	It is important to recall that ${\mathcal H}(L)$ is invariant under ambient isotopies but non-ambient isotopic links do not necessarily yield different linking numbers. Hence, it does not solve the ultimate goal sketched in remark \ref{Rem:ReminderOnKnots}.
%
\section{A multisymplectic interpretation of Massey products }\label{Sec:MasseyMess}
	In this section we resort to the techniques developed in the sections above and propose a reformulation of the so-called
higher order linking numbers in multisymplectic terms.
	
	\begin{figure}[h!]
		\vspace{-1em}
	 	\begin{center}
	 		\includegraphics[width=0.43\textwidth]{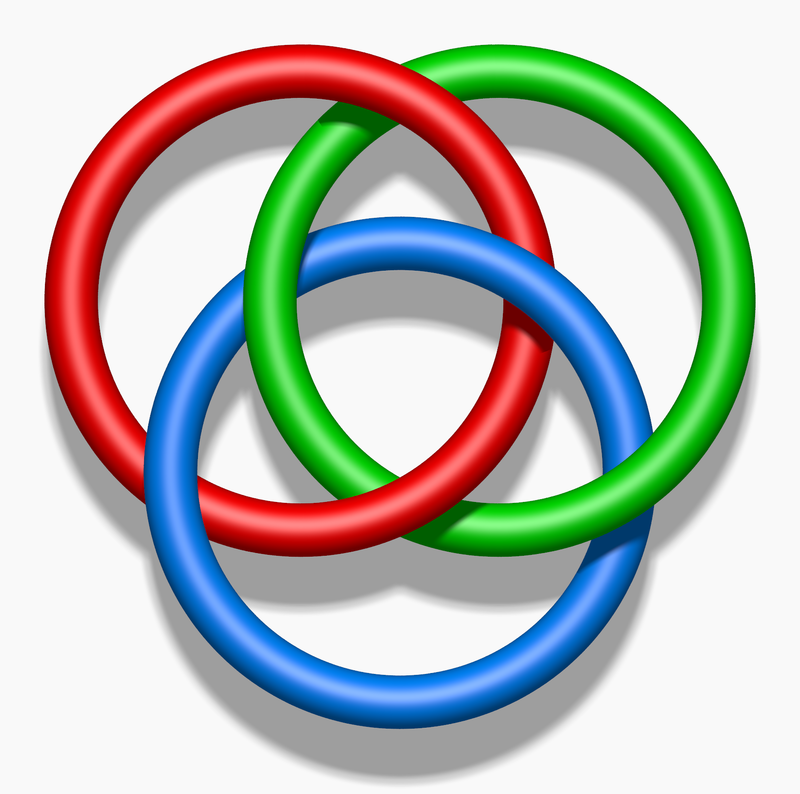}
		\end{center}
		\vspace{-1em}		
		\caption{
			\emph{Borromean rings} is a prototypical example of \emph{Brunnian link}. Removing a component of the link it yields a pair of unknots.
			\small(\href{https://commons.wikimedia.org/wiki/File:Borromean_Rings_Illusion.png}{Wikimedia Commons})
			\smallskip	
			}
		\label{Fig:BorromeanRings}		
	\end{figure}

	Ordinary and higher order linking numbers provide, among others,  a quite useful tool for the investigation of Brunnian phenomena in knot theory: recall that a link is {\it almost trivial} or {\it Brunnian} if upon removing any component therefrom one gets a trivial link (see figure \ref{Fig:BorromeanRings}).
They can be defined recursively in terms of Massey products, or equivalently, Milnor invariants,
by the celebrated Turaev-Porter theorem (see \cite{Fenn,Pe-Spe02,Spe06,Hebda-Tsau12}). We are going to review, briefly and quite concretely, the basic steps of the Massey procedure, read differential geometrically as in \cite{Pe-Spe02, Spe06,Hebda-Tsau12}, presenting at the same time our novel multisymplectic interpretation thereof. 

	As a first thing, let us recall how one can associate to any pair of knots in a link a certain cohomology class:
\begin{remark}[Cohomological reinterpretation of the ordinary linking number]\label{Rem:CohoReintOrdinaryLinking}
	Consider a pair of linked knots $L_1$ and $L_2$, possibly components of a bigger link.
	The cohomological reinterpretation of the ordinary linking number $\ell(1,2)$ of the components $L_1$ and $L_2$, starts from consideration of the $2$-form
	\begin{equation}
		\Omega_{1 2} := - v_{1}\wedge v_{2}
	\end{equation}
	with $v_i$ the velocity $1$-forms of the component $L_i$.
	Observe that $\d \Omega_{1 2} = -\omega_1\wedge v_2 + v_1 \wedge \omega_2$ is equal to the CS-form minus a "self-linking" term, \ie $\int \d \Omega_{1 2} = \ell(1,2)$ (\cf Theorem \ref{Thm:Moffatt-Ricca}).
	
	Although this form is not uniquely defined, it determines a unique (integral) de Rham class in the cohomology of the link (see Remark \ref{Rem:CohoIntLinkForm}) completely independent from the various choices
	\begin{displaymath}
		\langle L_1, L_2 	\rangle := 
		\left[\Omega_{1 2}\Big\rvert_{S^3\setminus L} \right] \in H^2(S^3\setminus L)
	\end{displaymath}
	The closedness of $\Omega_{1 2}\eval_{S^3\setminus L}$ follows from observing that $\text{supp}(\d \Omega_{1 2})$ is contained in a tubular neighbourhood centered along the link that can be taken arbitrarily small, hence, one can assume that $\text{supp}(\d \Omega_{1 2})\subset L$ hence $\d \Omega_{1 2}$ vanishes out of $L$. (See \cite[\S 2.3]{Pe-Spe02} or \cite[\S 6.2]{Spe06} for the complete argument.)
	
	The linking number $\ell(1,2)$ is non-zero precisely when $\langle L_1, L_2 \rangle $, which in $H_1(S^3,L)$ equals
$\ell(1,2) [\gamma_{12}]$, is non-trivial.
	If the latter class vanishes (\ie $\Omega_{1 2}\eval_{S^3\setminus L}$ is exact), we have
	\begin{equation}\label{eq:massey1}
		dv_{12} + v_{1}\wedge v_{2}  = dv_{12} + \Omega_{12} = 0.
	\end{equation}
	for some $1$-form $v_{12}$.
\end{remark} 
	Let $L$ be an oriented link with three or more components $L_j$ such that all the ordinary mutual linking numbers of the components under consideration vanish.
	Out of the primitives obtained in equation \eqref{eq:massey1}, one can manufacture another closed form:
	\begin{definition}[Third order linking number (class)]
		Let $v_{i j}$ be the primitive obtained through equation \eqref{eq:massey1} from the vanishing of the cohomology class $\langle L_i, L_j 	\rangle$ defined in remark \ref{Rem:CohoReintOrdinaryLinking}.
		We call {\it third order linking number} (as a class) for the three components $L_1, L_2,L_3$ the cohomology class:
	$$
		\langle L_1, L_2, L_3 \rangle :=   [\Omega_{123}\rvert_{S^3\setminus L}  ] \in H^2(S^3 \setminus L).
	$$
	of the (closed) \emph{Massey} $2$-form 
	\begin{displaymath}
		\Omega_{123} =  v_{1}\wedge v_{23}  + v_{12}\wedge v_{3}	
		~.
	\end{displaymath}
	\end{definition}
	It is then easy to devise a general pattern, 
	if the latter class vanishes, one can find a $1$-form $v_{123}$ such that
	\begin{equation}\label{eq:massey2}
		dv_{123} + v_{1}\wedge v_{23}  + v_{12}\wedge v_{3}  =  dv_{123} + \Omega_{123}  = 0.
	\end{equation}
	The iteration of this procedure, eventually obstructed by the non-vanishing of certain higher linking number classes, yields an hierarchy of pairs
	\begin{displaymath}
		v_I \in \Omega^1(S^3\setminus L) \qquad \Omega_{I} \in Z^2(S^3\setminus L)
	\end{displaymath}	
	with $I$ a general multi-index constructed out of the set $\{1,\dots,n\}$ of the $n$-link components.
	\begin{remark}
		The previous construction can be organised - via Chen's calculus of iterated path integrals \cite{Chen,Chen4}- in terms of sequences of {\it nilpotent connections} ${\mathbf v}^{(k)}$, $k=1,2...$ on a trivial vector bundle over $S^3 \setminus L$ and their attached  	{\it curvature forms} ${\mathbf w}^{(k)}$ (ultimately, the $\Omega_I$, \cite{Pe-Spe02,Spe06,Tavares,Hain}), everything stemming from the {\it Cartan structure equation}
	$$
		d {\mathbf v}^{(k)} + {\mathbf v}^{(k)} \wedge {\mathbf v}^{(k)} = {\mathbf w}^{(k)}
	$$
	together with the ensuing Bianchi identity
	$$
		d {\mathbf w}^{(k)} +  {\mathbf v}^{(k)} \wedge {\mathbf w}^{(k)} - {\mathbf w}^{(k)} \wedge {\mathbf v}^{(k)} = 0
	$$
	(the latter implying closure of the forms $\Omega_I$).
	In order to give a flavour of the general argument, start from the nilpotent connection ${\mathbf v}^{(1)}$ with its corresponding curvature ${\mathbf w}^{(1)}$:
	\begin{equation*}
		{\mathbf v}^{(1)} =
		\begin{pmatrix}
			0 & v_1 &  0 & 0 \\
			0 & 0    & v_2   & 0    \\
			0 & 0     &  0    &   v_3 \\
			0 & 0     &  0    &   0 \\
		\end{pmatrix}, \quad
		{\mathbf w}^{(1)} =
		\begin{pmatrix}
			0 & 0 &  \Omega_{12} = v_1\wedge v_2 & 0 \\
			0 & 0    & 0   & \Omega_{23} = v_2\wedge v_3    \\
			0 & 0     &  0    &   0 \\
			0 & 0     &  0    &   0 \\
		\end{pmatrix}
	\end{equation*}
	Then proceed similarly with
	$$
		{\mathbf v}^{(2)} =
		\begin{pmatrix}
			0 & v_1 &  v_{12} & 0 \\
			0 & 0    & v_2   & v_{23}    \\
			0 & 0     &  0    &   v_3 \\
			0 & 0     &  0    &   0 \\
		\end{pmatrix},\quad  {\mathbf w}^{(2)} =
		\begin{pmatrix}
			0 & 0 &  0 & \Omega_{123} = v_{1} \wedge v_{23} + v_{12} \wedge v_{3} \\
			0 & 0    & 0   & 0    \\
			0 & 0     &  0    &   0 \\
			0 & 0     &  0    &   0 \\
		\end{pmatrix}
	$$
	(we made use of $dv_{12} + \Omega_{12} = dv_{23} + \Omega_{23} = 0$), and so on.
	\end{remark}
	\begin{remark}
		Recall that all forms $\Omega_{I}$ can be neatly interpreted, via Poincar\'e duality, as auxiliary (trivial) knots $L_I$, and $v_I$ as discs bounded by  $L_I$, in adherence to the considerations in Section \ref{Sec:Ham1FormLinks}, see \cite{Pe-Spe02,Spe06} for more details and worked out examples,
		including the {\it Whitehead link} (involving fourth order linking numbers - with repeated indices) and the {\it Borromean rings} (exhibiting a third order linking number). 
		Just notice here that, for instance, formula (\ref{eq:massey1}) becomes, intersection theoretically
		$$
			\partial {\mathfrak a}_{12} + {\mathfrak a}_1 \cap {\mathfrak a}_2 = 0,
		$$
		see Figure \ref{fig: chen}.
		\begin{figure}[h!]
			\centering
			\includegraphics[width=\textwidth]{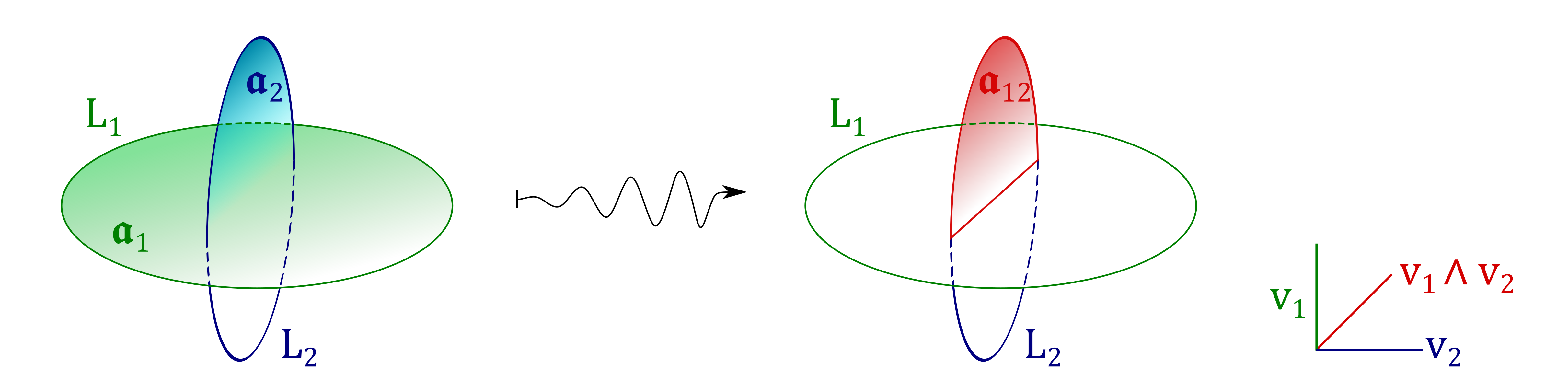}
			\caption{Starting the Chen procedure}
			\label{fig: chen} 
		\end{figure}
	\end{remark}
	
	We are now ready to interpret these quantities in the framework of multisymplectic higher observables:
	
	\begin{proposition}\label{Prop:MasseyMess}
		\begin{enumerate}[label=(\roman*)]
			\item All Massey $2$-forms are globally conserved in the sense of definition \ref{Def:conservedQuantities};
			\item Given an exact Massey $2$-form $\Omega_I$, the corresponding primitive $v_I$ will be Hamiltonian with respect to the volume $\nu\in \Omega^3(S^3\setminus L)$.
			The latter is obtained by extending first the standard Euclidean volume $\nu$ from $\R^3$ to its compactification $\R^3\cup \{\infty\} \cong S^3$ and then restricting it again from $S^3$ to $S^3\setminus L$.
			The corresponding Hamiltonian vector field is given by
			\begin{equation}
				\xi_I = \alpha^{-1}(\Omega_I)
				~,
			\end{equation}
			where $\alpha$ is the "flat" map introduced in equation \eqref{eq:alphacontractionmap}.
		\end{enumerate}
	\end{proposition}
\begin{proof}
	Formula (\ref{eq:massey2}) can be rewritten as
	$$
		dv_{123} + \iota_{\xi_{123}} \nu = 0 
	$$
	where ${\xi \equiv \xi_{123}}  = \alpha^{-1}(\Omega_{123})$.
	The above (``vorticity") vector field $\xi_{123}$ can be thought of as being concentrated on the knot corresponding to  $\xi_{123}$, or, alternatively, in a thin tube around it, when considering a bona fide Poincar\'e dual, \cf 
(\ref{eq:poincaredual}).

	This tells us that $v_{123}$ is a {\it Hamiltonian $1$-form} with respect to the volume $\nu$ and the formula
	$$
		{\mathcal L}_{\xi} \Omega_{123} =  
		\d  \iota_{\xi} \Omega_{123}  +  \iota_{\xi} \d \Omega_{123}  = 
		\d \iota_{\xi} \Omega_{123} 
	$$
 expresses the fact that $\Omega_{123}$ is a globally conserved $2$-form and the same holds for $\Omega_{12}$.
 
 This discussion can be carried out verbatim for a general  multi-index $I$:
	$$
		d v_I + \iota_{\xi_I} \nu = 0
	$$
(an extension of (\ref{eq:poincaredual}))
  and, in general, $\Omega_I$ is globally conserved.
\end{proof}
	We stress that, by construction, the momenta associated to the divergence-free field $\xi_I$ through the \momap constructed in section \ref{Sec:IdroMoMap} corresponds to $v_I$.
	
	The following is the main result of this section.
\begin{theorem}\label{Thm:MasseyMess}
	The $1$-forms $v_I$ are {\rm first integrals in involution} with respect to the flow generated by the Hamiltonian vector field $\xi_{ L}$ pertaining to the velocity $1$-form introduce in definition \ref{Def:Velocity1Form} , namely
	$$
	{\mathcal L}_{\xi_{ L}} v_I = 0
	$$
	(\ie the $v_I$'s are  {\rm strictly conserved})
	and the Poisson brackets given in equation \eqref{Eq:MSPB} above yields
	$$
	\{v_I, v_J \} = 0
	$$ 
	(for multi-indices $I$ and $J$).
\end{theorem}
\begin{proof}
Using Cartan's formula, we get
$$
{\mathcal L}_{\xi_{ L}} v_I = d \iota_{\xi_{ L}} v_I + \iota_{\xi_{ L}} d v_I = 
d \iota_{\xi_{ L}} v_I  -  \iota_{\xi_{ L}} \iota_{\xi_{ I}} \nu,
$$
but the second summand vanishes in view of the general expression
$$
\{v_\xi, v_\eta \}(\cdot) = \nu (\xi, \eta, \cdot)
$$
and of the peculiar structure of the vector fields involved (they either partially coincide or have disjoint supports). By the same argument,
one gets $\iota_{\xi_{ L}} v_I = 0$, in view of the Poincar\'e dual interpretation of $v_I$ (\cf Section \ref{Sec:Ham1FormLinks}), together with the second assertion;
a crucial point to notice  is that the auxiliary links obtained via Chen's procedure may be suitably split from their ascendants, this leading to
$$
\iota_{\xi_{ L}} v_I = 0,
$$
the consequent {\it strict} conservation of the $v_I$'s being then immediate.\par
\smallskip
Notice that, in particular, from
$$
\iota_{\xi_{ L}} v_L = 0
$$
(Poincar\'e dual interpretation again) we also get
$$
{\mathcal L}_{\xi_{ L}} v_L = 0
$$
(this is {\it not} to be expected a priori in multisymplectic geometry, \cf  \cite{Ryvkin2016}).
\end{proof}

We ought to remark that, upon altering the $v_I$'s by an exact form, we may lose strict conservation, but in any
case global conservation is assured (the Poisson bracket is an exact form, by equation \eqref{eq:bracketsformula} in Section \ref{Sec:MSHydroFluids} and in view of commutativity of the
vector fields $\xi_I$ and $\xi_J$).\par
\smallskip
Ultimately, we can draw the conclusion that  {\it the Massey invariant route to ascertain the Brunnian character of a link can be mechanically understood as a recursive test of a kind of {\rm knot theoretic integrability}: the Massey linking numbers provide obstructions to the latter}. \par
\smallskip
Thus, somewhat curiously, higher order linking phenomena receive an interpretation in terms of  multisymplectic geometry, which is a sort of higher order symplectic geometry. Also, integrability comes in with a twofold meaning: first, higher order linking numbers emerge from the construction of a sequence of  flat, \ie integrable nilpotent connections; second, this very process yields first integrals in involution in a mechanical sense. \par
%

\section{A symplectic approach to the HOMFLYPT  polynomial}\par
In this section, extracted from \cite{Miti2019a}, 
we present a novel {\it interpretation} of the HOMFLYPT 
polynomial (\cite{Freyd-etal,PT}) as a WKB-wave function via geometric quantization of the Brylinski manifold of singular knots (and links), drawing inspiration from the  {\it ad hoc} helicity-based hydrodynamical procedures devised in \cite{Liu-Ricca12,Liu-Ricca15}.

This section is divided in two parts.
Subsection \ref{sec:prelim}, gathers together basic symplectic geometric and knot theoretic tools which are necessary for understanding the problem.
Subsection \ref{Sec:HomWkb} contains the main theorem \ref{Thm:MainHomflyThm} in which the HOMFLY polynomial is read in the geometric quantization framework.
Before proving theorem \ref{Thm:MainHomflyThm}, we recall, and extend to links, some results contained in \cite{BeSpe06}.
These results are the foundation on which our -geometric quantization flavoured- approach is built.

\subsection{Preliminaries}\label{sec:prelim}

\subsubsection{The HOMFLYPT polynomial}\label{sub:HomflyPrelim}
In this subsection we recall a few basic knot theoretic notions, referring, among others, to
\cite{Rolfsen, Kauffman,Moffatt-Ricca92,BeSpe06,Spe06} for further information.\

First recall that, given a {\it framed}\footnote{Roughly speaking, the framing is tubular neighbourhood center around the knots of the knot.
	It is chosen to be "small" as desired so as to never intersect with itself in the vicinity of adjacent branches.
	See \cite{nlab:framed_link} for a concise description of this concept.
}
 oriented link, its {\it helicity} ${\mathcal H}(L)$  is given as
\begin{equation}\label{eq:helicity}
	{\mathcal H}(L) = \sum_{i,j =1}^n \ell(i,j)
\end{equation}
with $\ell(i,j) = \ell(j,i)$ being the {\it Gauss linking number} of components $L_i$ and $L_j$ if $i\neq j$ and where $\ell(j,j)$ is the {\it framing} of $L_j$, equal to $\ell(L_j, L_j^{\prime})$ with $L_j^{\prime}$ being a section of the normal bundle of $L_j$, see subsection \ref{subsec:ReltoGauss} above.
\\
A regular projection of a link $L$ onto a plane produces a natural framing called the {\it blackboard framing}\footnote{
Think of it as an indented projection of framing. In other words, it is the extension of the knot diagram from a curve to an infinitesimally thin ribbon.}, and ${\mathcal H}(L) = w(L)$, the {\it writhe} of $L$, given by
	\begin{equation}\label{eq:Writhe}
		w(L) = \sum_{\rm all \,\,crossings} \pm 1
	\end{equation}

\begin{figure}[h!]
	\begin{center}
		\includegraphics[width=0.43\textwidth]{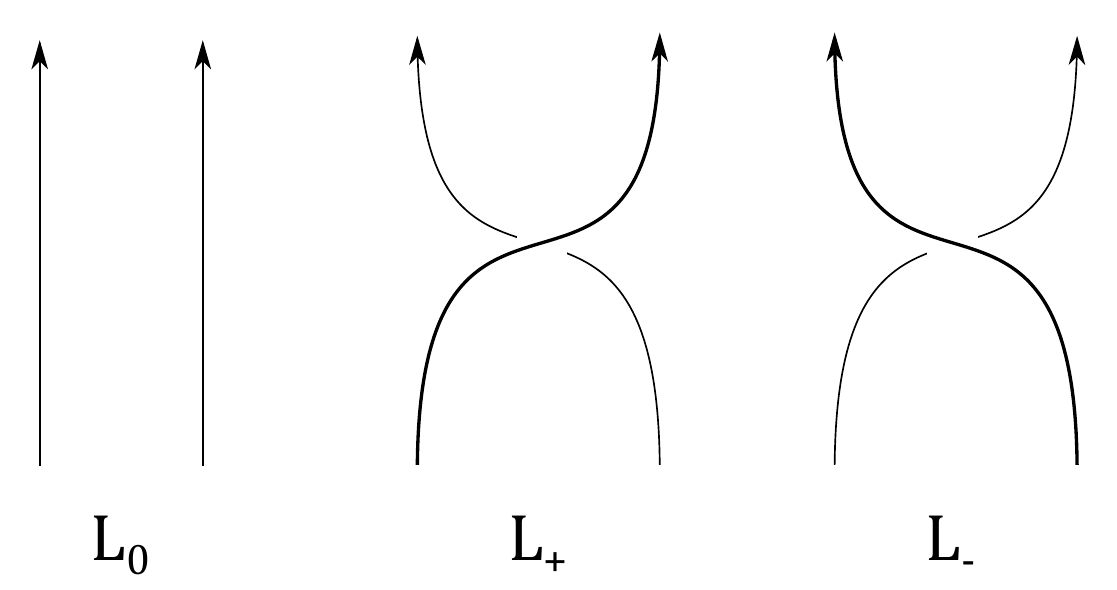}
	\end{center}
	\vspace{-3mm}
	\caption{Crossings}
	\label{fig: crossings}
\end{figure}
Let be $L$ an arbitrary link and consider a certain crossing pertaining to its indented diagram (\ie regularly projected onto a plane, $z=0$, say).
We denote, as usual, by $L_+$, $L_-$ and $L_0$ three oriented links  differing at a single crossing ($(\pm 1)$-crossing, no crossing, respectively), see Figure \ref{fig: crossings}.
Denoting by $m_\pm$ and $m_0$ their respective number of components, one has $m_+ = m_-  = m_0 \pm 1$.
\begin{remark}\label{Rem:L0pm}
	Let us stress the fact that the three links $L_{\pm}$ and $L_0$ may well be mutually inequivalent (thus they belong to different connected components of the manifold $Y$ defined below): 
	for instance, if $L_+$ is a trefoil knot, $L_-$ is a trivial knot and $L_0$ is a Hopf link. 	
	This situation is not general, considering
$L_+ = E_+$ and $L_- = E_-$ (``figure of eight" knots, \ie trivial knots with $\pm 1$ writhe, see figure \ref{fig:surgery}), both equivalent to the unknot $\bigcirc$, the corresponding $L_0$ is a trivial 2-component link.
\end{remark}

We already mentioned in reminder \ref{Rem:ReminderOnKnots} the notion of \emph{knot polynomials} as a tool to provide knot invariants. 
With the notation introduced in remark \ref{Rem:BryLoopSpaces}, one can see a knot polynomial as a function
\begin{displaymath}
	\hat{Y} \to R[X]
\end{displaymath}
from the space of (possibly) mildly singular knots to the $R$-algebra of polynomials on a certain ring $R$ and set of indeterminates $X$.
An example is given by the so-called \emph{HOMFLYPT polynomial} (see \cite{Freyd-etal, PT, Kauffman}).
\begin{definition}[HOMFLYPT polynomial]
	For any given link $L$, the corresponding \emph{HOMFLYPT polynomial}  is a $2$-variables (knot) polynomial $P = P(\alpha,z)$ defined by the {\it skein relation}
	\begin{displaymath}
		\alpha \,P ({L_+}) - \alpha^{-1}\, P({L_-}) = z \, P({L_0})~,
	\end{displaymath}
	for any crossing of a given diagram of $L$ (see above remark \ref{Rem:L0pm}),
	and by the normalization relation on the unknot $\bigcirc$
	 \begin{displaymath}
		 	P(\bigcirc) = 1~.
	 \end{displaymath}
\end{definition}
\begin{remark}
	The {\it HOMFLYPT polynomial} is an {\it ambient isotopy}.
	Diagrammatically, this amounts to invariance under the {\it Reidemeister moves}
 $R_j$, $j=0,1,2,3$.
The actual calculation may be performed, for instance, via the skein-template algorithm (see \eg \cite{Kauffman}, p.57).
	We recall also that $P$ is not a universal invariant, meaning that one can find inequivalent knots determining the same HOMFLYPT polynomial (see for example \cite{RAMADEVI1994}).
	\note{Domanda per Mauro:\\ ci sono fonti preferibili a questa?}
\end{remark}
The HOMFLYPT polynomial can be obtained from the so-called $H$-{\it polynomial}:
\begin{definition}[$H$-polynomial]\label{Def:Hpoly}
	For any given link $L$, the corresponding \emph{$H$-polynomial}  is a $2$-variables (knot) polynomial $H = H(\alpha,z)$ defined by the {\it skein relations}
	\begin{align*}
				H(L_+) - H(L_-) =&~ z \, H(L_0) ~,
				\\ 
				H(L \sharp {\,8_\pm}) =&~ \alpha^{\pm 1} H(L), 
	\end{align*}
	for any crossing of a given diagram of $L$.
	With $L \sharp {\,8_\pm}$ we mean the topological glueing of a "figure of eight"-shaped ``curl"  (see also section \ref{Sec:HomWkb} below for further interpretation of the notation).
\end{definition}
\begin{proposition}[``Kauffman principle" (see {\cite[p.55]{Kauffman}})]\label{prop:KaufPrin}
	It is possible to recover the HOMFLYPT polynomial out of the $H$-polynomial  via the standard writhe correction :
	\begin{equation}
		P(L) = \alpha^{-w(L)}\,H(L) = \alpha^{-{\mathcal H}(L)}\,H(L) ~.
	\end{equation}
	In particular, the $H$-polynomial is not an ambient isotopy invariant but only a regular isotopy invariant. Namely, one drops invariance under $R_1$, \ie addition of a ``curl".
\end{proposition}

\subsubsection{Lagrangian submanifolds revisited}
Recall that a submanifold $\Lambda$ of a symplectic manifold $(M, \omega)$ is {\it Lagrangian} when the symplectic form 
$\omega$ vanishes thereon and it is of maximal dimension with respect to this property. 

Consider now the cotangent space $T^*Q$ of a given smooth manifold $Q$ taken with its canonical symplectic form (see example \ref{Ex:Multicotangent}). 
A Lagrangian submanifold $\Lambda \subset T^*Q$ in general position can be described in the following way:
\begin{theorem}[Maslov-H\"ormander Morse family theorem \emph{(see \eg \cite{Mas,Hor,Gui-Ste,McD-Sal})}]\label{thm:morseFam}
	Consider a Lagrangian submanifold $\Lambda$ of the canonically symplectic manifold $T^\ast Q$.
	There exists (locally) a smooth function $\phi = \phi(q, a), (q, a) \in Q \times {\mathbb R}^k$ (${\mathbb R}^k$ being a space of auxiliary parameters)
and a submanifold
	\begin{equation}
		C_{\phi} = \{ (q, a) \in   Q \times {\mathbb R}^k \, \mid \, d_{a} \phi = 0  \}
\end{equation}
	with $d(d_{a})$ of maximal rank thereon (here $d = d_q + d_{a}$) such that the map
\begin{equation}
		\lambdamorphism{C_{\phi}}
		{T^*Q  }
		{(q, a)}
		{(q,d_q \phi )}
\end{equation}
is an immersion with image $\Lambda$. 
\end{theorem}
If the Hessian $H_{a}$ (with respect to the auxiliary variables $a$) is non-degenerate, one can solve $a = a(q)$ and introduce the following objects:
\begin{itemize}
	\item we call \emph{phase function} the smooth function $F\in C^{\infty}(Q)$, also denoted $F(q):= \phi(q, a(q))$, such that $(q, dF(q) ) \in \Lambda$;
	\item given a phase function, we call \emph{momentum} at $q\in Q$ the covector $dF(q) =: p(q)$.
\end{itemize}
The non-degeneracy of the Hessian fails at the singular points of the obvious projection $\Lambda \rightarrow Q$.
However, the singular locus $Z$ (the {\it Maslov cycle}) turns out to be orientable and of codimension $1$ in $\Lambda$ with $\partial Z$ of codimension $\geq 3$ (see \eg \cite[\S II.7]{Gui-Ste}).


\subsubsection{WKB wave functions}
Recall that, given a {\it prequantizable} symplectic manifold $(M,\omega)$,  the Weil-Kostant theorem (see theorem \ref{thm:integralitycondition} above) implies the existence of a complex line bundle ${\mathcal L} \to M$ (prequantum bundle), equipped with a Hermitian metric and compatible connection $\nabla$ with curvature $\Omega_{\nabla} = -2\pi i \omega$. 

Since the symplectic 2-form $\omega$ vanishes on any Lagrangian submanifold $\Lambda \subset M$, any (local) symplectic potential $\vartheta$ 
($d\vartheta = \omega$) becomes a closed form thereon, giving a (local) connection form  pertaining to the restriction of the prequantum connection $\nabla$. 
The former is a {\it flat} connection that will be denoted by the same symbol $\nabla$.
When $\nabla$ has trivial holonomy, one can defined the so-called \emph{WKB-wave functions}:
\begin{definition}[WKB-wave function]
	We call a \emph{WKB-wave function} any global covariantly constant section of the restriction to $\Lambda$ of the prequantum line bundle ${\mathcal L} \to M$.
	In other terms, the latter is a section $s\in \Gamma(\mathcal{L},\Lambda)$ such that $\nabla s = 0$.
\end{definition}
A WKB-wave function $s$ takes the local form (neglecting the so-called half-form correction, see \eg \cite{Woodhouse97})
 \begin{equation}
 s(m) := hol_\gamma(\nabla) \cdot  s(m_0) = e^{i \int_{\gamma} \theta} \cdot s(m_0) 
 \end{equation}
with $\gamma$ denoting any path connecting a chosen point $m_0$ in a (connected) symplectic manifold $M$ with a generic point $m \in M$, ${\rm hol}_{\gamma}(\nabla)$ being the holonomy of the (restriction to $\Lambda$ of the) prequantum connection $\nabla$ along $\gamma$. 
The right-hand side tacitly assumes a trivialization of $L \to M$  around $m_0$, and $m$ in a corresponding local chart.


\subsection{The HOMFLYPT polynomial as a WKB wave function}\label{Sec:HomWkb}
Let be $(M,\nu)$ a $3$-dimensional smooth manifold oriented by the volume form $\nu$.
Consider the free loop space $LM$ and the Brylinski's manifold $\widehat{Y}_{M}$ of (mildly) singular {\it knots} embedded in $M$ introduced in reminder \ref{Rem:BryLoopSpaces}.

Recall that, by transgression, one gets a 2-form $\Omega$ on $LM$ via the formula
\begin{equation}
\Omega = \int_{S^1} ev^* (\nu )
\end{equation}
where $ev: S^1 \times LM \rightarrow M$, given by $ev (\lambda , \gamma ) := \gamma(\lambda)$, is posthe evaluation map (of a loop $\gamma \in LM$
at a point $\lambda \in S^1$). 
More explicitly, given  tangent vectors $u$ and  $v$ at $\gamma$, it reads
\begin{equation}
\Omega_{\gamma} (u , v)  = \int_{0}^1  \nu ({\dot \gamma}(\lambda), u(\lambda), v(\lambda))
\end{equation}
(where we set ${\dot \gamma} = {d\gamma \over d\lambda}$).

The 2-form $\Omega$ is basic with respect to the ${\rm Diff}^+(S^1)$-principal bundle $\widehat{X}_{M} \rightarrow \widehat{Y}_{M}$,
namely $i_{\xi} \Omega = i_{\xi} d\Omega = 0$, with $\xi$ any vertical vector field (\ie generating an orientation preserving
reparametrization of the loop).
Hence, $\Omega$  descends to a closed, non-degenerate 2-form on $\widehat{Y}_{M}$, \ie a (weak) {\it symplectic form}.

The keypoint is given by the following proposition:
\begin{proposition}[Pre-quantum bundle for the Brylinksi's manifold $\widehat{Y}_{M}$]
	Let be $(M,\nu)$ a $3$-dimensional manifold with volume form $\nu$.
	If $[\nu]$ is an integral $3$-form then $(\widehat{Y}_{M},\Omega)$ is prequantizable.
\end{proposition}
\begin{proof}
	Recall that, in general, the transgression gives rise to a (degree shifting) morphism of complexes 
	$\Omega^{\bullet} (M) \rightarrow \Omega^{\bullet  - 1} (LM)$, mapping closed (resp. exact) forms to closed (resp. exact) ones
	in view of the general formula (direct calculation, or see \cite{Chen}):
	\begin{equation}
	d \int \omega = - \int d \omega
	\end{equation}
	where, of course, the l.h.s. differential pertains to $LM$ and the r.h.s. one pertains to $M$. 
\\
Consequently, integral cohomology classes on $M$ are mapped to integral cohomology classes on $LM$.
Therefore, if $[\nu ]$ is integral, then $[\Omega]$ is integral as well, this ensuring, via the Weil-Kostant theorem,
the existence of a prequantum bundle.
\\ A subtle though explicit construction can be given via the integral class  $[\nu ] \in H^3(M, {\mathbb Z})$, defining a {\it gerbe}, see \cite{Bry,Spe11}.
\end{proof}

We  can naturally extend the above discussion to oriented links and accordingly define the symplectic structure on the {\it generalized Brylinski space of oriented mildly singular links} $\widehat{Y}$ (same notation)
via the
same formula above, by replacing a knot $K$ by a link $L$:
\begin{equation}\label{eq:LinkSymplecticForm}
	\Omega_L (\cdot , \cdot )= (\int_L \nu) \, \,(\cdot,\cdot):= \sum_{i=1}^n \int_{L_i} \nu (\dot{\gamma}_i, \cdot,\cdot)
\end{equation}
\begin{remark}
	We also point out, for completeness, that the weak symplectic manifold $(\widehat{Y}_{M}, \Omega)$  can be endowed with several other natural structures:
	\begin{itemize}
		\item $\widehat{Y}_{M}$ can be naturally
equipped with a (formally) integrable compatible almost complex structure making it a K\"ahler manifold in an appropriate sense, see \eg \cite{Bry,Pe-Spe89,Arn-Khe}.
		\item Each connected component of $\widehat{Y}_{M}$ is (up to technical subtleties, see \cite{Bry} and \cite{Pe-Spe89} as well) a {\it coadjoint orbit} of the group of unimodular diffeomorphisms of $M$, \ie those preserving a volume form.
	\end{itemize}		
\end{remark}

\medskip
We shall  deal with the case $M = {\mathbb R}^3$; the ensuing manifold 
$\widehat{Y} := \widehat{Y}_{{\mathbb R}^3}$ is called the manifold of {\it oriented singular links} in ${\mathbb R}^3$, whereas $Y := {Y}_{{\mathbb R}^3}$
is called the manifold of {\it oriented  links} in ${\mathbb R}^3$.

\begin{remark}[$\widehat{Y}$ is prequantizable]\label{rem:OmegaExact}
	Observe that the standard volume form $\nu$ in ${\mathbb R}^3$ can be portrayed as
	\begin{equation}
		\nu = dx \wedge dy \wedge dz = d (z\, dx \wedge dy) \equiv d\hat{\theta}
	\end{equation}
	in terms of the (multisymplectic) potential $\hat{\theta}$; the latter transgresses to a (symplectic) potential $\theta$ for $\Omega$, which vanishes identically when restricted on the plane $z=0$. 
	\\
	In other terms, $\Omega$ is exact.
	Thence the assumptions of the Weil-Kostant theorem are fulfilled and we have a trivial prequantization bundle ${\mathcal L} \to \widehat{Y}$ (Brylinski's line bundle).
\end{remark}

\subsubsection{The Chern-Simons Lagrangian submanifold}
Consider the weak symplectic manifold  $T^*{\widehat{Y}}$ that is the cotangent space associated to $\widehat{Y}$ together with its canonical symplectic structure.

One can naturally introduce an appropriate Morse family\footnote{To be formally interpreted in the sense of H\"ormander (see \eg \cite{Mas},\cite{Hor}).}, \cf  with theorem \ref{thm:morseFam}, treating the space of $U(1)$-connections ${\cal A}$ as a set of auxiliary parameters\footnote{
It may be identified with ${\cal D}_{{\mathbb R}} ({\mathbb R}^3)\otimes {\mathbb R}^3$, the space of compactly supported (real) vector fields
on ${\mathbb R}^3$, and standardly topologized accordingly. 
We regard it as an infinite dimensional manifold modelled on itself
(\cf \eg \cite{Kri-Mich}, p.439).
}.
This family is given by the (link analogue) of the {\it Abelian Chern-Simons action with source}, \ie is  the smooth function\footnote{Concretely, one first checks G\^ateaux differentiability of $\Phi =\Phi (K, \cdot )$,
then observes that   $\Phi$ is indeed Fr\'echet differentiable (for instance by checking continuity of the G\^ateaux derivative,
see \eg \cite[p.128]{Kri-Mich}). }
 $\Phi \in C^{\infty}(Y\times \cA)$ given by
	\begin{equation}
		\Phi (L,A) := {k\over{8\pi}}\int_{{\mathbb R}^3}A \wedge dA + \int_L A 
	\end{equation}
Accordingly to theorem \ref{thm:morseFam}, $\Phi$ defines locally a Lagrangian submanifold $\widetilde{\Lambda}$ of the cotangent space $T^*{\widehat{Y}}$  via the position
	\begin{equation}
		d_{\cal A}\Phi \mid_{(K,A)} = 0
	\end{equation}
where $d_{\cal A}$ denotes the differential of ${\cal A}$.
Furthermore, the Lagrangian submanifold $\widetilde{\Lambda} \hookrightarrow T^*{\widehat{Y}}$ admits a local phase $\phi\in C^\infty(\widehat{Y})$
given by
\begin{equation}\label{eq:phasephi}
\phi (L) = -\frac{2\pi}{k}  {\mathcal H}(L) \equiv 2\pi \lambda \,{\mathcal H}(L) 
\end{equation}
where ${\mathcal H}(L)$ is the helicity of the link (see equation \eqref{eq:helicity})  and $\lambda := -1/k$ is non-zero number
\footnote{The generic value taken by $\lambda$ (in particular, it can be taken equal to a root of unity) avoids trivialities.}
(see \cite[\S 3]{BeSpe06} and references therein for details).
The previous discussion in subsumed by the following statement

\begin{theorem}[Helicity as a phase function \emph{(\cite[Thm 3.1]{BeSpe06})}]
	The writhe (helicity) of a link can be interpreted as a phase function pertaining to $\widetilde{\Lambda}$, looked upon
as a Lagrangian submanifold of $T^*\widehat{Y}$.
\end{theorem}

\subsubsection{Link invariants via geometric quantization}
Consider now the submanifold $\Lambda \subset {\widehat Y}$ consisting of the links on a plane with transversal intersections.
It has been proved in \cite{BeSpe06} that $\Lambda$ is a Lagrangian submanifold with respect the symplectic given in equation \eqref{eq:LinkSymplecticForm}.
Henceforth we tacitly identify a link in $Y$ with its projection onto $\Lambda$ (whilst retaining crossing information via a $\pm$-marking). 
\\
Since the symplectic potential of Brylinski's form can be taken equal to zero (see remark \ref{rem:OmegaExact}) the phase is (locally) constant.
This is in accordance with the fact that the phase is essentially given by the helicity and the latter is a topological invariant.  
The CS-Lagrangian can be used, as in \cite{BeSpe06}, as a Morse family also for $\Lambda$.
The Lagrangian submanifold $\Lambda$ is thence locally given by the graph
\begin{equation}
(L, d \, {\mathcal H}(L)) = (L, 0)
\end{equation}
($d {\mathcal H}(L) = 0$ is the so-called {\it eikonal equation}, see \cite{BeSpe06,Spe06}).

	In our context the assumptions of the Weil-Kostant theorem are fulfilled (see remark \ref{rem:OmegaExact}).
	Hence we have a trivial prequantization bundle ${\mathcal L} \to \widehat{Y}$,
restricting to $\Lambda$ (and denoted by the same symbol).
A covariantly constant section  on ${\mathcal L} \to \Lambda$ is then just a {\it locally constant function} on $\Lambda \subset \widehat{Y}$  since, locally $\theta = d{\mathcal H} = 0$ and we neglect the so-called ``half-form" correction (see \eg \cite{BeSpe06} and \cite{Woodhouse97}). 
\\
Then, observe that the phase function in equation \eqref{eq:phasephi}  gives rise precisely to the exponent of the -regular isotopy- \emph{Witten invariant}
\begin{equation}\label{eq:wittenInv}
 \psi = \psi (L) := e^{2\pi i \lambda{\mathcal H}(L)} 
\end{equation}
(obtained in \cite{Witten89} via a path integral computation, see also \cite{Gua,Kauffman,BeSpe06}).
Furthermore $\psi$ arises as a WKB wave function for the prequantum line bundle ${\mathcal L} \to \Lambda$ in a purely topological theory.
\smallskip
	The discussion developed in the preceding subsections, extending to links the 	results obtained for knots in \cite{BeSpe06}, 
	can be summarized via the following theorem (suitably merging and extending Theorems 3.1, 4.1, 5.1 of the above paper):\par
\begin{theorem}[Witten invariant as a WKB-phase function]
Consider the Brylinski' manifold $\widehat{Y}$ of mildly singular links.
The regular isotopy Witten invariant $\psi$ can be interpreted as WKB wave function in Brylinski's framework.
		Namely, it corresponds to a covariantly constant section of the restriction of the prequantum line bundle to the submanifold $\Lambda$ of singular knots on a plane, viewed as a Lagrangian submanifold of $\widehat{Y}$.
\end{theorem}

In this spirit, we propose a geometric quantization interpretation of the HOMFLYPT polynomial building on the Besana-Spera symplectic approach to framing via Brylinski's manifold of mildly singular links.
%
%
%
 \begin{theorem}[HOMFLYPT as a WKB-phase function]\label{Thm:MainHomflyThm}
The HOMFLYPT polynomial  can be recovered from the geometric quantization procedure applied to the Brylinski manifold $\widehat{Y}$ and to its Lagrangian subspace $\Lambda$.
\\
	Namely, it coincides (after normalization) with a suitable covariantly constant section $\Psi = \Psi(\alpha,z)$.
	The coefficient $\alpha$ is a  phase factor related to the helicity of a standard ``eight-figure" and $z$ comes from accounting for the variation of the number of components of a link.  
\end{theorem}
\begin{proof}
\begin{figure}[h!]
	\begin{center}
		\includegraphics[width=0.43\textwidth]{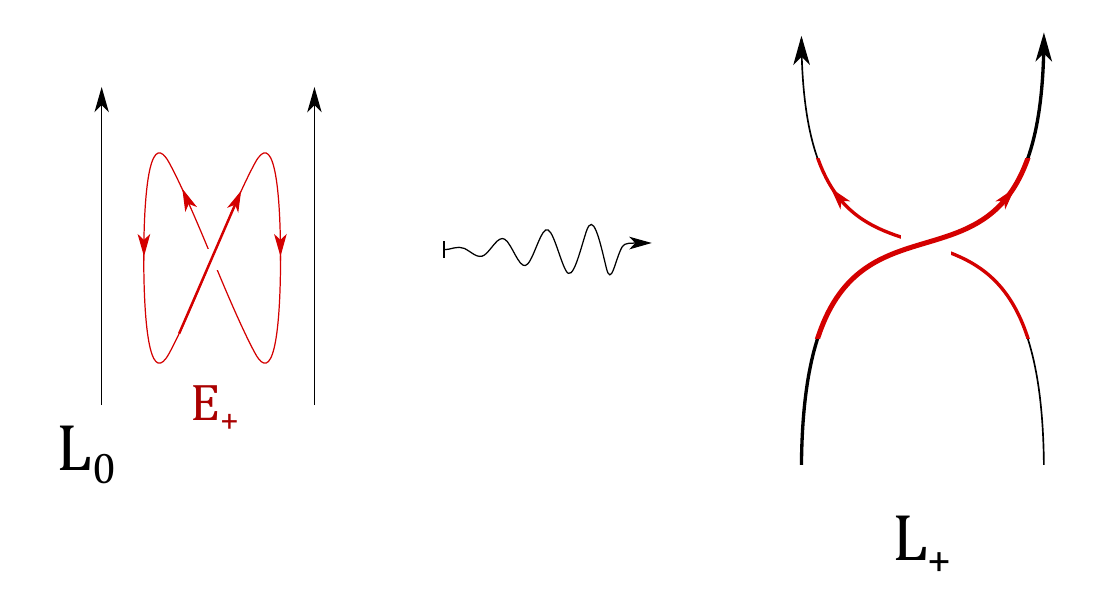}
	\end{center}
	\vspace{-3mm}
	\caption{Surgery via $E_+$}
	\label{fig:surgery}  
\end{figure}
Let be $L$ a link. Consider three links $L_0$, $L_\pm$ pertaining to a given crossing in the indented diagram of $L$ (see subsection \ref{sub:HomflyPrelim}).
\\
Then, inspired by the Liu-Ricca approach (\cite{Liu-Ricca12,Liu-Ricca15}), let us consider the above ``figures of eight" $E_{\pm}$. Having the figure of eight a single crossing, one has that ${\mathcal H}(E_{\pm}) = \pm 1$.
Starting, for instance, from $L_0$, one can ``add" $E_+$ to the two coherently oriented parallel strands of $L_0$  in such a way that $E_+$ comes with the opposite orientation: a partial cancellation occurs and the net result is 
$L_+$.  
Conversely, proceeding backwards we can, by adding appropriately an $E_-$, produce $L_0$ from $L_+$ and so on. 
Therefore, addition of $E_{\pm}$ allows one to pass from one local configuration to the other\footnote{ We explicitly notice that the two strands may indifferently belong to the same or to two different components of the link involved.}, see Figure \ref{fig:surgery}. 
\\
Now set: 
\begin{equation}
	\begin{aligned}
		\alpha :=&~ e^{2\pi i \lambda {\mathcal H}(E_+)} = e^{2\pi i \lambda}
		~,
		\\
		\alpha^{-1} =& e^{-2\pi i \lambda} = e^{2\pi i \lambda {\mathcal H}(E_-)}  ~,
	\end{aligned}
\end{equation}
 so that, trivially (compare with equation\eqref{eq:Writhe})m
  \begin{equation}\label{eq:vecchia4.2}
  	\begin{aligned}
 \psi(L_{\pm}) =&~ \alpha^{\pm 1} \psi(L_{0})
 		~,
 		\\
 		 \psi(L_{\pm}) =&~ \alpha^{\pm 2} \psi(L_{\mp})
 		~,
 		\\
 		\psi (L_+)  - \psi (L_-) =&~ (\alpha  - \alpha^{-1}) \psi (L_0)
 		~.
 		\end{aligned}
  \end{equation}
Thus we see that $\alpha^{\pm 1}$ arises as the local contribution to the WKB wave function $\psi$ upon addition (surgery) of an eight figure (or ``curl") - which can be applied to a single branch as well (first Reidemeister move: this explains our notation $ L \sharp {\,8_\pm}$ for the link obtained from $L$ by adding a $\pm$-curl on one of its strands)
and $\alpha^{\pm 2}$ as the corresponding contribution upon crossing the Maslov cycle $Z$.

We now wish to modify the above formula \eqref{eq:vecchia4.2} so as to produce a genuine {\rm ambient isotopy} link invariant, while keeping the above interpretation, that is, we shall prove the following 
\\
First of all, the prequantization bundle ${\mathcal L} \to \Lambda$ can be trivialized
via the trivial link invariant $\Psi_0 := {\mathbb 1}$. 
We wish to alter $\Psi_0$ beyond the connected component of the unknot $\bigcirc$ in order to get a non-trivial invariant.
This will be accomplished via a minimal alteration of the r.h.s. of  equation \eqref{eq:vecchia4.2}. 
Given links $L_{\pm}$ and $L_0$, with respect to one of their crossings (we already observed that they may all be mutually inequivalent, \ie they may lie in different connected components of $Y$), we can nevertheless identify their respective fibres ${\mathcal L}_{\pm}$ and ${\mathcal L}_0 
$ of the line bundle ${\mathcal L}$ via  the given trivialization. 
Let us then slightly modify the above formula \eqref{eq:vecchia4.2} and first look for a regular isotopy invariant wave function $\widetilde{\Psi}$ such that
\begin{equation}
	\begin{aligned}
	\widetilde{\Psi}(L_+) - \widetilde{\Psi}(L_-) =&~ 
	z \,\widetilde{\Psi}(L_0)
	~,
	\\
	\widetilde{\Psi}(L \sharp {\,8_\pm}) =&~ 
	\alpha^{\pm 1} \widetilde{\Psi}(L)
	\end{aligned}
\end{equation}
with $z$ now {\it independent} of $\alpha$. 
But this is precisely the skein relation for the $H$-polynomial (see definition \ref{Def:Hpoly}).
Hence  $\widetilde{\Psi}$ exists and can be promoted to the looked for HOMFLYPT polynomial wave function $\Psi$ via Kauffman's principle (proposition \ref{prop:KaufPrin}).
	Namely
	\begin{equation}
		\Psi(L) := \alpha^{-w(L)}\,\widetilde{\Psi}(L) = \alpha^{-{\mathcal H}(L)}\,\widetilde{\Psi}(L),
	\end{equation}
fulfilling
 \begin{equation}
 	\begin{aligned}
		\alpha \Psi ({L_+}) - \alpha^{-1} \Psi ({L_-}) =&~ z \,\Psi ({L_0}) ~,
		\\
		\Psi(\bigcirc) =&~ 1
		~.
	\end{aligned}
\end{equation}
	The above procedure, applied to $\psi$, yields the trivial invariant $\Psi_0 = \mathbb{1}$. 
\end{proof}

\begin{remark}
	\begin{enumerate}
			\item The skein relation for the $H$-polynomial, and thence for the HOMFLYPT polynomial, has been {\it used} in order to guarantee the representation of the latter as a covariantly constant section of the above line bundle.  
 			As such, our interpretation is consistent and {\it crossing independent}.	
 			\item The position $z= \alpha^{-\frac{1}{2}} - \alpha^{\frac{1}{2}}$ reproduces the {\it Jones} polynomial.
For $\alpha = 1$ and $z \neq 0$ we recover the {\it Conway} polynomial. The case $z = 0$
yields the trivial invariant $\Psi = \Psi_0 = {\mathbb 1}$.
			\item The skein relation (5.3) can be equivalently  written in the form
\begin{equation}
 \Psi(L_{-}) =  \alpha^{2} \,\Psi(L_{+}) - z \alpha \,\Psi(L_{0})
\end{equation}
which tells us that $\Psi(L_{-})$ can be obtained by suitably adding $\Psi(L_{+})$, corrected  by a {\it Maslov type} transition (switching term: local surgery via   $\alpha^{2}$ - one has the same number of link components) and $\Psi(L_{0})$,   corrected by a splicing term (``component transition")   $\alpha$ (and multiplied by an extra coefficient $-z$).
The latter contribution was absent in   \cite{BeSpe06} since that paper dealt with   {\it knots} only. 
Notice that the apparent notational clash (one would naively expect a switch $\alpha \leftrightarrow\alpha^{-1}$ in (5.3) and (5.4)) is simply due to the Kauffman principle.
			\item Passage from $L_{\pm}$ to $L_0$ (and conversely) in $Y$ - abutting, as already remarked, at a change in the number of the link components - involves coalescence
of {\it two} opposite crossings into one and corresponding tangent alignment
(a sort of ``dipole", related to the second Reidemeister move). This is a sort of ``higher order" contribution beyond the Maslov one.
			\item In this way we essentially recover the hydrodynamical portrait of Liu and Ricca \cite{Liu-Ricca12,Liu-Ricca15}, essentially stating that `` $P= t^{\mathcal H}$ " via a different (and more conceptual) interpretation. In particular, the meanings of the two parameters used in HOMFLYPT are not quite the same. 
The local surgery operation involves helicity, as in Liu-Ricca, but we portray the latter as yielding a local phase function, governing a component  transition or Maslov, upon squaring it, as in \cite{BeSpe06}.
			\item The Chern-Simons (CS) 3-form $A_L \wedge dA_L $
can be interpreted, in adherence to  \cite{Spe11}, as a {\it connection 3-form} for a {\it 2-gerbe}, having {\it zero curvature} on ${\mathbb R}^3 \setminus L$. The provisional wave function $\psi = \exp 2\pi i \lambda {\mathcal H}(L)$ then essentially becomes the ``parallel transport" of this connection ``along"  ${\mathbb R}^3$, and it is already an important topological invariant.
			\item The considerations made in the previous remark may provide the starting point of a {\it multisymplectic} reformulation, involving $({\mathbb R}^3, \nu)$ instead of its transgressed symplectic manifold  $({L}{\mathbb R}^3, \Omega)$, possibly casting further light on the Jones-Witten theory.
	\end{enumerate}
\end{remark}

\ifstandalone
	\bibliographystyle{../../hep} 
	\bibliography{../../mypapers,../../websites,../../biblio-tidy}
\fi

\cleardoublepage


%% file: chapters/gradedmultilinear/gradedmultilinear.tex
\note{
	Hint MZ non c'è tempo:\\
Before appendix A, maybe you can add a few lines recalling that the parts of app A and B that are necessary to read the body of the thesis, are condensed in Ch 1. You said it earlier, but about 200 pages earlier....Also emphatize what is "new" in the appendix
}

\chapter{Graded multilinear algebra}\label{App:GradedMultilinearAlgebra}
	This chapter introduces some basics definitions of linear algebra on graded vector spaces and establishes the notation adopted in the thesis's body.
	\\
	Although this material can now be considered standard, it is useful to include it for the sake of fixing unambiguously the several conventions related to graded objects.
	\\
	Most of the content is drawn from the following sources  \cite{Manetti-website,Bandiera2016,Schatz2009,Reinhold2019,Doubek2007,Delgado2015,Fiorenza2006}.
	
To clarify how certain conventions and sign rules emerge, we thought it worthwhile to arrive at the definition of graded vector space starting from the more abstract notion of \emph{category of graded objects}. 
	However, the reader may decide to jump directly to the section \ref{Section:GradedVectorSpaces} where the category of graded vector spaces, the main setting of this thesis, is introduced.
		
\section{Categorical prelude: categories of graded objects}\label{Section:GradedObjects}
	In this section will be assumed some basic notions in category theory, 
	see \cite{Riehl2016} or \cite{MacLane1978} for an introductory exposition. 
	Note in particular that all the category considered are \emph{locally small}.

	\subsection{Basic definitions}
	Consider an arbitrary set $S$ and a category $\cat$. We call $S$-graded object any family of objects of $\cat$ indexed by $S$.
	\begin{definition}[Category of $S$-Graded $\cat$-objects]\label{Def:GradedObjects}\index{$S$-Graded $\cat$-objects}
	We define the \emph{Category of $S$-Graded $\cat$-objects} as the functor category
		\begin{displaymath}
			\cat^{S}:= 
			\left\lbrace
					\hat{S}\to \cat \text{ functors}
			\right\rbrace
		\end{displaymath}
		where $\hat{S}$ is the set $S$ regarded as a discrete category (objects are the elements on the set $S$ and the only arrows are identities).
	\end{definition}
	\begin{remark}
	Definition \ref{Def:GradedObjects} implies that a morphism $\phi$ between two graded object $A$ and $B$ in $\cat^S$ is a natural transformation
	\begin{displaymath}
		\begin{tikzcd}
			\left(\phi: A \xrightarrow[]{~\cdot~} B  \right)
			& \equiv &
			 S \arrow[r, bend left=50, "A"{name=U, above}]
			\arrow[r, bend right=50, "B"{name=D, below}]
			& \cat
			\arrow[shorten <=1pt,shorten >=1pt,Rightarrow, from=U, to=D, "\phi"]
		\end{tikzcd}	
	\end{displaymath}
	\end{remark}
	\begin{notation}
		A $S\times S$ graded $\cat$-object is also called \emph{S-bi-graded}.
		Such terminology could be iterated straightforwardly.
	\end{notation}
	\begin{remark}
	Concretely, any object $A \in \cat^{S}$ is equivalent to a family $\{A_s\}_{s\in S}$ of objects in $\cat$ labelled by the elements of $S$.
	Similarly,
	any morphism $\phi\in \Hom_{\cat^S}(A,B)$ can be seen as a collection $\phi = \lbrace \phi_s : A_s \to B_s \rbrace_{s\in S}$ of morphisms in $\cat$ indexed by $S$.
	\\
	Informally, both can be seen as functions from $S$ to, respectively, objects and morphisms of $\cat$. 
	This definition is more correctly phrased in terms of functors since $ob(\cat)$ is not in general a set.
	\end{remark}
	\begin{remark}
		Observe that $\hat{S}=\hat{S}^{op}$, therefore graded objects are in particular presheaves.
	\end{remark}
	\begin{notation}
		When one specializes the category $\cat$ to be the category of sets (resp. $R$-modules, $\mathbb{K}$-vector spaces or $\mathbb{K}$-algebras; being $R$ a ring and $\mathbb{K}$ a field), 
		$S$-graded objects are called $S$-graded sets (resp. $S$-graded $R$-modules, $S$-graded $\mathbb{K}$-vector spaces or $S$-graded $\mathbb{K}$-algebras).
	\end{notation}
	Since $\cat^S$ is a category of functors, it inherits several properties from the categorical structure of $\cat$:
	\begin{proposition}\label{prop:limitsofFunctorCat}
		Consider two categories $\cat$ and $\cat[D]$ and denote by $[\cat,\cat[D]]$ the category of functors from $\cat$ to $\cat[D]$.
		\begin{enumerate}
			\item If $\cat[D]$ has limits or colimits of a certain shape, then so does $[\cat,\cat[D]]$.
			\item If $\cat$ is small and $\cat[D]$ is Cartesian closed and complete, then $[\cat,\cat[D]]$  is Cartesian closed. 
			\item If $\cat[D]$ is pre-additive so is $[\cat,\cat[D]]$, in particular if $\cat[D]$ is Abelian so is $[\cat,\cat[D]]$.
		\end{enumerate}
	\end{proposition}	
	\begin{proof}
		\begin{enumerate}
			\item Limits in $[\cat,\cat[D]]$ can be computed "point-wise" from the limits in $\cat[D]$. 
			Note that if $\cat[D]$ is not complete, then can exist other limits in $[\cat,\cat[D]]$ that are not defined point-wise.
			(See \href{https://ncatlab.org/nlab/show/functor+category#properties}{\cite{nlab:functor_category}} for details).
			\item The sketch of the proof can be found in \href{https://ncatlab.org/nlab/show/cartesian+closed+category#exponentials_of_cartesian_closed_categories}{\cite{nlab:cartesian_closed_category}}	.
			\item The $Ab$-enrichment property on $[\cat,\cat[D]]$, as well as the $ker$ and the $coker$ of a morphism, can be defined again point-wise from the pre-additive structure on $\cat[D]$. 
			(The specific proof for $\mathbb{Z}$-graded objects in an Abelian category can be found in \href{https://stacks.math.columbia.edu/tag/09MF}{\cite[Lemma 12.16.2.]{stacks-project}}.)
		\end{enumerate}
	\end{proof}
	
	In practice, there are three main choices for the grading set $S$ that appear most frequently in the literature: 
	the natural numbers $\mathbb{N}$, the integers $\mathbb{Z}$ and the group $\mathbb{Z}/2$ \footnote{A $\mathbb{Z}/2$-graded object it is also called a \emph{super-object}.}.
	Note that these three objects are not simply sets but they carry a canonical algebraic structures: $\mathbb{N}$ is a monoid, $\mathbb{Z}$ and $\mathbb{Z}/2$ are Abelian groups.
	In the following, we will see how the algebraic structure of $S$ reverberates in the properties of $\cat^S$.

	\subsection{Degree shifts}\label{sec:degreeshifts}
	Let be $\cat$ an arbitrary category and be $S$ a grading set.
	Observe that, given a $S$-graded object $A\in \cat^S$, exchanging the components associated to two fixed labels define a different graded object.
	More generally, an $S$-indexed family of objects can be re-indexed giving a different $S$-graded object.
	This operation can be encoded as a functor:
	\begin{definition}[Degree shift functor]\label{Def:ShiftFunctor}\index{Shift functor}
		Given a (set-theoretic) function $f:S\to S$, denote by $\hat{f}:\hat{S}\to\hat{S}$ its promotion to an endofunctor on the discrete category $\hat{S}$. 
		We define the \emph{degree shift pertaining to $f$} as the endofunctor $[f]:\cat^S \to \cat^S$ defined on a graded object $A$ as the precomposition of the associated functor with  $\hat{f}$, \ie
		\begin{displaymath}
			[f](A) := A[f] = A \circ \hat{f}
			~,
		\end{displaymath}
		diagrammatically denoted as
		\begin{displaymath}
		\begin{tikzcd}
			\hat S \ar[r,"\hat f"]
			\ar[dr,dashed,"A\lbrack f \rbrack"'] 
			& \hat S \ar[d,"A"] 
			\\
			 & \cat			
		\end{tikzcd}
		~.
		\end{displaymath}				
		The action on a morphism $\phi:A\to B$ is given by the (horizontal) precomposition of the natural transformation $\phi$ with the functor $\hat{f}$.
	\end{definition}
	Concretely, for any two graded objects $A$ and $B$, for any graded morphism $\phi:A\to B$, and for any $i\in S$, one has that $(A[f])_i = A_{f(i)}$ and	 $\phi[f]_i =\phi_{f(i)}$, \ie the following diagram commutes in the category $\cat$
	\begin{displaymath}
	\begin{tikzcd}
		(A{[f]})_i \ar[r,equal] \ar[d, "\phi{[f]}_i"'] & A_{f(i)} \ar[d, "\phi_{f(i)}"] 
		\\
		(B{[f]})_i \ar[r,equal] & B_{f(i)}
	\end{tikzcd}
	~.
	\end{displaymath}		
	
	\begin{remark}[Grading given by a monoid]\label{rem:GradingGivenByMonoid}
		A little bit more abstractly, one could say that there exists a monoid action\footnote{One must pay extra care when considering the collection of natural transformations since it is not a set in general.
However, when $\cat$ is small, $\cat^S$ is locally small \cite{Freyd1995}.}
		\begin{displaymath}
			(\Hom_{Set}(S,S),\circ) \xrightarrow{[~\cdot~]} (\Hom_{Cat}(\cat^S,\cat^S),\circ))
			~,
		\end{displaymath}
		\ie there is a sort of "action" of the monoid of endomorphism of $S$ on the category of $S$-graded $\cat$-object.
%
	\\
	In particular, when $S$ is a monoid, with product denoted by $\cdot$, it is natural to consider the subset of (left) actions of $S$ on itself 
	\begin{displaymath}
		\left\{
			\left[
				\morphism{L_s}{S}{S}{q}{s \cdot q}
			\right]
		\right\}_{s \in S} \subset
		(\Hom_{Set}(S,S),\circ)
	\end{displaymath}
	which is in a one-to-one correspondence to $S$.
	In other words, there is an action of $S$ on $\cat^S$, that is a monoid homomorphism
	\begin{displaymath}
		\morphism{[\cdot]}
		{S}
		{\Hom_{Cat}(\cat^S,\cat^S)}
		{s}
		{
			\left(
				\Almorphism{[s]}
				{\cat^S}
				{\cat^S}
				{(k\mapsto A_k)}
				{(k\mapsto A_{s\cdot k})}
			\right)
		}
	\end{displaymath}
	Being a monoid morphism guarantees that $[s'][s]=[s'\cdot s]$.
	\\
	When the monoid is commutative, there is no need to distinguish between left or right action; hence we simply talk about $[s]$ as the \emph{shift by $s \in S$}.
	\\
	When $S$ is also a group, $[L_s]$ functors are invertible with $[L_s]^{-1} = [L_{s^{-1}}]$.
	\end{remark}
	\begin{remark}
		Note that the definition of shift can be seen as a particular instance of the action of $\End(\cat[S])=\Hom_{cat}(\cat[S],\cat[S])$ on $[\cat[S],\cat[D]]$ given by "precomposition" or "pullback". In other terms,
		\begin{displaymath}
			[\cdot]: \End(\cat[S])
			=\Hom_{cat}(\cat[S],\cat[S]) \to \Hom_{cat}(\cat[S],\cat[D])
			~.
		\end{displaymath}
	\end{remark}
	\begin{remark}
		Let be $A \in \cat^S$ a $S$-graded object.
		When $S$ is a group, it is trivial to observe that, taken all together,  the shifts define a $S$-graded $\cat^S$ object 
		\begin{displaymath}
			\left(k \mapsto A[k]\right) ~.
		\end{displaymath}
		Therefore, one could alternative regard the shift operation as a functor $[\bullet]: \cat^S \to (\cat^S)^S$.
	\end{remark}
	\begin{notation}[Graded-morphisms vs homogeneous maps]
		Given any $k\in S$ and $A,B$ graded objects, it is customary to call any map in $\Hom_{\cat^S}(A,B[k])$ \emph{homogeneous map in degree $k$}\index{Homogeneous maps}. 
		In particular \emph{morphisms in $\cat^S$} are \emph{degree $0$ homogeneous maps}. 
	\end{notation}

	\subsection{Graded objects in categories with (co)limits}
		Recall that, according to proposition \ref{prop:limitsofFunctorCat}, any limit\footnote{In the categorical sense, \ie universal cone over a certain diagram. See \cite[Chap. V]{MacLane1978}.} $Lim$ in $\cat$ can be easily extended to a limit in $\cat^S$ simply defining $(Lim(D))_i= Lim (D_i)$ for any diagram $D$ in the category of $S$-graded $\cat$-objects
		
		Any category $\cat$ with terminal object can be easily regarded has a full subcategory of $\cat^S$:
		\begin{lemma}\label{Lemma:CembedsinCS}
			Given $k\in S$, a category $\cat$ and a fixed object $T\in ob(\cat)$,
			the functor
			\begin{displaymath}
				G_k : \cat \to \cat^S ~,
		\end{displaymath}			 
		acting on objects as
		\begin{displaymath}
			(G_k(A))_i = 
			\begin{cases}
				A & i=k \\
				T & i\neq k			
			\end{cases}
		\end{displaymath}
		and on morphisms as
		\begin{displaymath}
			\biggr(G_k\big(A\xrightarrow \phi B\big)\biggr)_i = 
			\begin{cases}
				(A\xrightarrow \phi B\big) & i=k \\
				(T \xrightarrow {id_T} T ) & i\neq k			
			\end{cases}~,
		\end{displaymath}
		is faithful. If $T$ is terminal object in $\cat$ the functor $G_k$ is also full.
		\\
		Given two indices $k,k'$ in $S$, one has $G_{k'} = [f] \circ G_k$ where $f$ is the unique invertible function on $S$ such that $f(k)= k', f(k')= k$ and $f(i)=i$ when $i\neq k,k'$ (\ie the function swapping index $k$ with $k'$). 
		In a diagram, the following triangle of functors commutes
		\begin{displaymath}
		\begin{tikzcd}
			\cat \ar[r,"G_k"]
			\ar[dr,"G_{k'}"'] 
			& \cat^S \ar[d,"{[f]}"] 
			\\
			 & \cat^S			
		\end{tikzcd}
		~.
		\end{displaymath}	

		\end{lemma}
		\begin{proof}
			The definition is well-posed, functoriality is enforced by the composition rule of morphisms of $\cat$.
			\\
			Faithfulness follows easily from the definition; namely  
			$G_k(f) = G_k(g)$ if and only if $f = g$ for any $f,g \in \Hom_{\cat}(A,B)$.
			\\
			Fullness is a consequence of the universal property of the terminal object since the only possible morphism between $(G_k(A))_i = (G_{k'}(A))_i$ when $i\neq k,k'$ is the terminal arrow $T \to T$.			
			\\
			The last property follows from the definition of shift, for any $A \in \cat$ one has
			\begin{displaymath}
				\big( G_k(A)[f] \big)_i 
				= \big( G_k(A)\big)_{f(i)} 
				= \begin{cases}
					0 & f(i) \neq k\\
					A & f(i) = k
				\end{cases}
				= \begin{cases}
					0 & i \neq k'\\
					A & i = k'
				\end{cases}
				= \big( G_{k'}(A) \big)_i 
			\end{displaymath}

		\end{proof}
		\begin{remark}\label{Rem:TrivialExtension}
			Observe that, given $A\in\cat^S$ and  $T\in\cat$ terminal object,
			any morphisms $f\in \Hom_{\cat}(A_k,B)$ can be uniquely extended to a morphism $f \in \Hom_{\cat^S}(A,G_k(B))$ in the trivial way
			\begin{displaymath}
			f_i = \begin{cases}
				A_i \xrightarrow ! T & i\neq k \\
				A_k \xrightarrow f B & i = k
			\end{cases}
			~.
			\end{displaymath}			 
			Uniqueness is a simple consequence of the universal property for $T$.		
		\end{remark}

		It follows from Lemma \ref{Lemma:CembedsinCS} that any object in a category $\cat$ with terminal object can be seen as a "special" object in $\cat^S$. 
		When $\cat$ possesses all coproducts $\coprod$, and the terminal object $T$ is initial (\ie zero object), also the "converse" statement holds:
		\begin{lemma}\label{lemma:totalfunctor}
			Given a category with all coproducts $\coprod$, consider the functor
			\begin{displaymath}
				\morphism{F}{\cat^S}{\cat}{A}{\coprod_{i\in S} A_i}
				~,
			\end{displaymath}
			acting on morphisms as
			\begin{displaymath}
				\biggr(F\big(A\xrightarrow \phi B\big)\biggr) = \coprod_{i\in I} \varphi_i 
				~, 
			\end{displaymath}
			where $\coprod_{i\in S} \varphi_i$ is the unique arrow provided by the universal property of coproducts. 
			Namely, given two arrows $(f:A_1\to B_1)$ and $(g:A_2\to B_2)$ one could define their sum as the unique arrow $f\sqcup g $ provided by the universal property of the coproduct $\sqcup$:
			\begin{displaymath}
				\begin{tikzcd}[column sep=small, row sep =tiny,ampersand replacement=\&]
					A_1 \ar[dr]\ar[dd,"f"']\& \& A_2\ar[dl]\ar[dd,"g"] \\
					\& A_1 \sqcup A_2 \ar[dd,"f\sqcup g ","!"'] \& \\
					B_1 \ar[dr] \& \& B_2 \ar[dl] \\
					\& B_1 \sqcup B_2 \& \\				
				\end{tikzcd}
			~.
			\end{displaymath}	
			When $T$ defining $G_k$ is an initial object, $F$ is a left inverse of $G_k$, \ie there exists a canonical natural isomorphism 
			$$F \circ G_k \Rightarrow id_{\cat}~.$$
		\end{lemma}
		\begin{proof}
			$F$ is well-defined since the universal property of $\coprod$ guarantees functoriality.
			The natural isomorphism is given on every component by the canonical isomorphism implieds by the universal property of coproducts and initial objects
			$$
				A \coprod I \simeq A \simeq I \coprod A ~.
			$$
		\end{proof}
		\begin{remark}
			Note that the two functors $F$ and $G_k$ are not adjoint.
		\end{remark}		
	\begin{example}\label{Example:GradedSets}
		When $\cat$ is the category of sets with functions and $\coprod$ is the disjoint union of sets, a \emph{S-graded set} $A$ can be regarded as a pair:
		\begin{displaymath}
			\left(
				\tilde{A}= \coprod_{s\in S} A_s 
				~,~
				\morphism{\vert\cdot\vert}{\tilde{A}}{S}{A_s\ni a}{s}
			\right)
		\end{displaymath}
		where 
		$\tilde{A}=\{(a,s)\vert s\in S, a\in A_s\}$ is called the \emph{(total) set of homogeneous elements} and $\vert\cdot\vert$ is the grading function.
	\end{example}

\begin{remark}\label{Remark:GradedObjectasCoproduct}
	The reasoning of example \ref{Example:GradedSets} can be extended to any \emph{concrete category}, \ie categories of structured sets.
	Recall that a concrete category can be formalized as category $\cat$ together with a faithful functor $J:\cat \to Set$, called \emph{forgetful functor}, which associates to any object $A \in \cat$ the underlying sets of points $\underline{A} \in Set$. 
	The forgetful functor implies a functor from the category of $S$-graded $\cat$-objects to the category of $S$-graded sets simply given by post-composition with $J$,
	\begin{displaymath}
		\morphism{J'}{\cat^S}{Set^{S}}{A}{J\circ A}
		~.
	\end{displaymath}
	Therefore, one can define a "grading" as a function on the set of homogeneous elements $\coprod_{s\in S}\underline{A_s}$.
	\\
	When $\cat$ is a concrete category with coproducts $\coprod$ one could mimic the definition in example \ref{Example:GradedSets} identifying any graded object $A$ as a pair
	\begin{displaymath}
		\left(
			\coprod_{s\in S} A_s 
			~,~
			\morphism{\vert\cdot\vert}{	\coprod_{s\in S} \underline{A_s}}{S}{\underline{A_s}\ni a}{s}
		\right)
		~.
	\end{displaymath} 
	Note that one can not define the grading as a function on $\underline{\coprod_{s\in S} A_s }$, or as morphism in $\cat^S$, since, in general, the subset of homogeneous elements $\coprod_{s\in S} \underline{A_s}$ is only a subset of $\underline{\coprod_{s\in S} A_s}$. 
	When the two coincide the category is said to \emph{admit concrete coproducts} \cite[1.3.3]{Bourles2017}. 
	%
\end{remark}
		\begin{remark}[About the total functor]\label{rem:totaloplus}
			When $\cat$ is pre-additive, \ie finite products and coproducts coincide (see for example the category $\Vect$ of vector spaces in section \ref{Section:GradedVectorSpaces}), it is customary to denote coproducts as $\oplus$ and the functor $F$ of lemma \ref{lemma:totalfunctor} as $(-)^\oplus$.
			Given a graded object $A$, $A^\oplus$ is sometimes called the \emph{total object} or \emph{total space}.
		\end{remark}

	\subsection{Monoidal structure}\label{Section:CategoricalGradedMonoidalStructure}
		In simple terms, a monoidal category is a category equipped with a notion of "tensor product" (see \cite{MacLane1978} or also \cite{nlab:monoidal_category}).
		Most of the subtleties that can be encountered when dealing with graded spaces derive from consistently using the monoidal structure.
		The keypoint is that, when certain conditions are met, a category of graded objects $\cat^S$ inherit the monoidal structure from $\cat$.

		Observe first that that the Cartesian structure $\times$ on Cat and Set implies that there exists a canonical functor
		\begin{displaymath}
			\morphism{\times}
			{\cat^S \times \cat^S}
			{(\cat\times\cat)^{S\times S}}
			{((k\mapsto A_k),(s\mapsto B_s))}
			{((k,s)\mapsto A_k\times B_s)}			
			~.
		\end{displaymath}
		When $\cat$ is equipped with the monoidal structure $\otimes$,
		the postcomposition of the previous functor with $\otimes$ yields a functor from pairs of $S$-graded objects to $S$-bigraded objects:		
		\begin{displaymath}
			\morphism{(\otimes) \circ (\times)}
			{\cat^S \times \cat^S}
			{(\cat)^{S\times S}}
			{((k\mapsto A_k),(s\mapsto B_s))}
			{((k,s)\mapsto A_k\otimes B_s)}			
			~.
		\end{displaymath}
		In the case that the grading set $S$ is equipped with a monoid product $\cdot$, and given a monoidal structure $\oplus$ on $\cat$ (possibly different from $\otimes$) one can also introduce the canonical functor
		\begin{equation}\label{eq:totalbigradedcomplex}
			\morphism{\tot}
			{\cat^{S\times S}}
			{\cat^S}
			{((k,s)\mapsto A_{(k,s)})}
			{(g\mapsto \bigoplus_{s\cdot k =g}A_{(k,s)})}
			~.		
		\end{equation}
		In general, when two suitable compatible monoidal structures occur, 
		it is possible to introduce a canonical monoidal structure on the category of graded objects.
		Namely, one can consider the consecutive composition of the three previous functors:
		\begin{proposition}\label{Prop:GradedMonoidalStructure}
			Given a monoidal structure $(\cdot,1)$ on $S$ and be $\cat$ a "bimonoidal"
			\footnote{A \emph{Rig-category} \cite{nlab:rig_category} is 
			  a category with two monoidal structures $(\cat,\oplus,\mathbb{0})$ and 
			 $(\cat,\otimes,\mathbb{1})$, where the first one is symmetric, 
			 together with left and right distributivity natural isomorphisms
			\begin{align*}
				d_\ell: A\otimes(B\oplus C) \to (A \otimes B) \oplus (A \otimes C) \\
				d_r: (A\oplus B) \otimes C \to (A \otimes C) \oplus (B \otimes C)		
			\end{align*}
			and absorption/annihilation isomorphisms
			\begin{align*}
				a_\ell: A\otimes \mathbb{0} \to \mathbb{0} \\
				a_r : \mathbb{0} \otimes A \to \mathbb{0}		
			\end{align*}
			satisfying a set of coherence laws (see for example \cite{Laplaza1972}).	
		}
		category $(\cat,\oplus,\mathbb{0},\otimes,\mathbb{1})$,
		the endofunctor $ \tilde{\otimes} = (\tot)\circ(\otimes)\circ(\times): \cat^S\times \cat^S \to \cat^S$, defined on objects as
		\begin{displaymath}
			(A \,\tilde{\otimes}\, B)_s = \bigoplus_{g\cdot" h =s} A_g \otimes B_h
			~,
		\end{displaymath}
		yields a monoidal structure on $\cat^S$.
		\end{proposition}
		\begin{proof}
			Defining the unit object in $\cat^S$ as
			\begin{displaymath}
			  \begin{tikzcd}[
			    column sep=2em,
			    row sep=-1ex,
			  ]
			  \tilde{\mathbb{1}}\colon &[-3em]
			  S \arrow[r] & \cat \\
			  & 1 \arrow[r,mapsto] & \mathbb{1} \\
			  & s\neq 1 \arrow[r,mapsto] & \mathbb{0}
				\end{tikzcd}
			\end{displaymath}
			the definition of the associator, left unitor and right unitor of $\tilde{\otimes}$ follows straightforwardly from the their counterparts related to $\otimes$ on each degree.			
		\end{proof}
		\begin{remark}\label{Rem:CanonicalInducedBraiding}
		 In the case that $\cat$ is equipped with a braiding with respect to $\otimes$ (\ie $\cat$ is a "braided monoidal" category), there is a natural induced braiding pertaining to $\tilde{\otimes}$,
		 \begin{displaymath}
		 	(B_{V,W})_n = \bigoplus_{k+\ell = n} B_{V_k,W_\ell}
		 	: (V\otimes W)_n \to (W \otimes V)_n
		 	~,
		 \end{displaymath}
		 where $B_{V_k,W_\ell}$ is the braiding provided on $\cat$ and $(V,W)$ is an arbitrary pair of $S$-graded $\cat$-objects.
		\end{remark}
		\begin{remark}\label{Rem:AnticipazionediKoszul}
			When the grading set is a ring $(S,+,\cdot)$, there are two possible monoidal structures to consider when constructing $\tilde{\otimes}$.
			In the case that $S=\mathbb{Z}$, it is customary to choose the sum operation as the monoid structure defining $\tilde{\otimes}$ while the product is used to twist the induced braiding with an extra sign.
			(Compare with the so-called "Koszul convention" in section \ref{sec:bloodyKoszulConvention}).
		\end{remark}

\section{The category of graded vector spaces}\label{Section:GradedVectorSpaces}
	In what concern this thesis, we are only interested in the case where $S=\mathbb{Z}$ and $\cat=\Vect$ is the category of vector spaces over the field of real numbers $\mathbb{R}$.
	Most of what follows could be extended to vector spaces over a generic field $\mathbb{K}$ of characteristic $0$ and $R$-modules with a relatively small effort.
	\begin{definition}[Graded vector spaces]\label{Def:GradedVectorSpaces}
	 We call \emph{graded vector space} any object in the $\mathbb{Z}$-graded vector spaces category $\Vect^{\mathbb{Z}}$. 
	 In other words, a graded vector space is a collection $\{V^k\}_{k\in\mathbb{Z}}$ of vector spaces over $\mathbb{R}$ indexed by $\mathbb{Z}$.
	 A morphism of graded vector spaces $V \to W$ is a family of linear maps $\{f^k : V^k \to W^k \}_{k\in \mathbb{Z}}$.
	\end{definition}
	\begin{notation}
		Notice that, when dealing with graded vector spaces, we slightly changed the notation introduced in section \ref{Section:GradedObjects}.
		Namely, we denote the degree $k$ sector of a graded vector space as $V^k$, \ie we put a superscript $k$ instead of a subscript.
		This choice is in accordance with the "cohomological convention" (see section \ref{sec:HomologicalAlgebrasConventions}) employed throughout the text.	
	\end{notation}

	\begin{remark}[Categorical properties of $\Vect$]\label{rem:Vectpdvs}
		Recall that $\Vect$ is an extremely rich category. 
		we will be interested in the following properties:
		\begin{itemize}
			\item $\Vect$ is a complete and cocomplete category. 
			In particular the trivial 0-dimensional vector space $\mathbb{0}$ is a \emph{zero object}, \ie it is both initial and terminal:
				\begin{displaymath}
					\begin{tikzcd}[row sep =tiny]
						V \ar[r,two heads,"!"] & \mathbb{0} \ar[r,hook,"!"]& W\\[-.7em]
						v \ar[r,mapsto] & 0 \ar[r,mapsto] & 0_W
					\end{tikzcd}									
					~.	
				\end{displaymath}
				$\Vect$ admits products, called \emph{direct products}, for any set of indices $J$, given by:
				\begin{displaymath}
					\mathclap{
					\prod_{j\in J} W_j = 
					\left(
						\underbrace{\dots\times W_j \times \dots}_{\text{Set-theoretic Cartesian product}},
						\underbrace{
							\begin{aligned}
(\omega_j)_{j\in J} + (\omega_j')_{j'\in J} &= (\omega_j + \omega_{j'})_{j\in J}
							\\
							\lambda(\omega_j)_{j\in J} &= (\lambda \omega_j)_{j \in J}									
							\end{aligned}
						}_{\text{index-wise \\linear structure}}					
					\right)
					}
				\end{displaymath}
				with surjective "projectors":
				\begin{displaymath}
					\morphism{{\pi_k}}
					{\displaystyle\prod_{j\in J} W_j }
					{W_k}
					{(\omega_j)_{j\in J}}
					{\omega_k}
					~.
				\end{displaymath}		
				$\Vect$ admits infinite coproducts, called \emph{direct sums}, given by:
				\begin{displaymath}
					\bigoplus_{j\in J} W_j = 
					\left\lbrace\left.
						(\omega_j)_{j\in J} \in \prod_{j \in J} W_j \right\vert 
						\quad\omega_j \neq 0_{W_j} \quad\parbox{9em}{only for a finite\\ numbers of indices $j$} 
					\right\rbrace
				\end{displaymath}
				with injective "injectors":
				\begin{displaymath}
					\morphism{{\iota_k}}{W_k }{\bigoplus_{j\in J} W_j}
					{\omega}{\left(\omega_j= \begin{cases}
						0_j & j \neq k \\
						\omega & j = k				
						\end{cases}
						\right)_{j \in J}}
					~.
				\end{displaymath}
				By definition $\displaystyle \bigoplus_{j \in J} W_j \subseteq \prod_{j\in J} W_j $ is a linear subspace, the equality holds when $\#J<\infty$.
			
				\item $\Vect$ is a symmetric monoidal category. 
				The tensor product functor is given by the usual tensor product of vector spaces:
			\begin{displaymath}
				V_1 \bigotimes_{\mathbb{K}} V_2 
					= \frac{\Free( V_1 \oplus V_2)}{\sim}	
			\end{displaymath}
			where $\Free( V_1 \oplus V_2)$ denotes the free vector space generated by the elements of set $V_1\times V_2$ and $\sim$ is the vector subspace encoding distributivity 
			\begin{displaymath}
			\begin{aligned}
				(v+v',w) \sim (v,w) + (v',w) \\
				(v,w+w') \sim (v,w) + (v,w')
			\end{aligned}
			\end{displaymath}
			and scalar multiplication
			\begin{displaymath}
				(\lambda v , w) \sim \lambda (v,w) \sim (v,\lambda w)~.
			\end{displaymath}
			Such object is a vector space by the definition of quotient of a vector space with respect to a  linear subspace.
			The unit object is given by the 1-dimensional vector space $\mathbb{R}$.
			The corresponding associator, left unitor, and right unitor isomorphisms, are trivial and usually treated as simple identifications\footnote{Technically, $\Vect$ is not "strict" monoidal (see \cite{Muger2010}), \eg ${V_1 \otimes (V_2 \otimes V_3)}$ and ${(V_1 \otimes V_2 ) \otimes V_3}$ are isomorphic but -in principle- different.
			However, most of the time they will be considered "equal" implying the underlying isomorphism $\alpha$. The same idea is also applied to $\lambda$ and $\rho$.}:
			\begin{displaymath}
				\isomorphism{\alpha}
				{V_1 \otimes (V_2 \otimes V_3)}
				{(V_1 \otimes V_2 ) \otimes V_3}
				{x_1\otimes(x_2 \otimes x_3)}
				{(x_1\otimes x_2)\otimes x_3}
			\end{displaymath}
			\begin{displaymath}
				\isomorphism{\lambda}
				{V \otimes \mathbb{R}}
				{V}
				{v\otimes \lambda}
				{\lambda v}
			\end{displaymath}
			\begin{displaymath}
				\isomorphism{\rho}
				{\mathbb{R} \otimes V}
				{V}
				{\lambda \otimes v}
				{\lambda v}
				~.
			\end{displaymath}
			There is also a trivial symmetric braiding:
			\begin{displaymath}
				\morphism{B_{V_1,V_2}}
				{V_1 \otimes V_2}
				{V_2 \otimes V_1}
				{v_1\otimes v_2}
				{v_2 \otimes v_1}
				~.
			\end{displaymath}

			\item $\Vect$ is a \href{https://ncatlab.org/nlab/show/distributive+monoidal+category}{distributive monoidal category} \cite{nlab:distributive_monoidal_category}.
			That means that the monoidal structure $\otimes$ distributes over the (cartesian) monoidal structure given by direct sums, thus we have two canonical isomorphisms
			\begin{displaymath}
				\isomorphism{d_\ell}
				{X \otimes (Y \oplus Z)}
				{(X\otimes Y)\oplus(X\otimes Z)}
				{x\otimes(y+z)}
				{x\otimes y + x \otimes z}
			\end{displaymath}
			\begin{displaymath}
				\isomorphism{d_r}
				{(X\oplus Y) \otimes Z}
				{(X\otimes Z)\oplus(Y\otimes Z)}
				{(x+y)\otimes z}
				{x\otimes z + y \otimes z}
				~.
			\end{displaymath}			
			This is a particular case of \emph{Rig-category} \cite{nlab:rig_category}.
			
			\item $\Vect$ is a $\Vect$-enriched category (also known as $\mathbb{R}$-linear category or \emph{algebroid} \cite{nlab:algebroid}).
			That means that there is a canonical $\mathbb{R}$-structure on the space of linear maps $\Hom_{\Vect}(V,W)$ or, in other words, there exists an "internal hom functor"
			\begin{displaymath}
				[-,-] \,:~ {\Vect^{op}\times \Vect} \to {\Vect}
			\end{displaymath}
			defined on a pair of objects $V,W\in\Vect$ as
			\begin{displaymath}
				[V,W] =
					\left(
						\underbrace{\Hom_{\Vect}(V,W)}_{\text{Set of linear maps}},
						\underbrace{
							\begin{aligned}
							(A+B)(v) &= A(v)+B(v)
							\\
							(\lambda A)(v) &= \lambda A(v)									
							\end{aligned}
						}_{\substack{\text{vector-wise linear structure}\\ \forall v \in V \\ \forall A, B \in \Hom_{\Vect}(V,W)}}					
					\right)
					~.
			\end{displaymath}
			 			 
			 Clearly the usual $\Hom$ functor factors through the internal one and the forgetful functor neglecting the linear structure.
			It follows easily from the definition  that the composition map of linear maps is bilinear with respect to the linear structure on the hom-space,
			therefore there is natural morphism
			\begin{displaymath}
				\morphism{-\circ-}
				{{[V,W]\otimes [W,V']}}
				{{[V,V']}}
				{(A:V\to W)\otimes (B:W\to V')}
				{(B\circ A : V \to V')}
				~.
			\end{displaymath}						 
			
			\item $\Vect$ is \href{https://ncatlab.org/nlab/show/closed+monoidal+category}{monoidal closed} \cite{nlab:closed_monoidal_category}.
			The functor $[-,-]$ defined above is also \emph{internal} in the sense of monoidal closed categories;
			\ie, for any $V \in \Vect$, $[V,-]$ is right adjoint to $-\otimes V$. 
			The latter can be proved by exhibiting an explicit "currying" \cite{nlab:currying} natural isomorphism in the category of sets:
			\begin{equation}\label{eq:Curry}
				\mathclap{
				\isomorphism{"Curry"}
				{\Hom_{\Vect}(V\otimes W, L)}
				{\Hom_{\Vect}(V,[W,L])}
				{(h:V\otimes W \to L)}
				{
					\left(
						\Almorphism{\hat{h}}
						{V}
						{[W,L]}
						{v}
						{(h(v,-):W\to L)}
					\right)				
				}}
			\end{equation}
			Note that this isomorphism realizes precisely the \emph{universal property of tensor products of vector spaces}.
			It is not difficult to see that $\Hom_{\Vect}(V,[W,L])$ consists of functions from $V\times W$ to $L$ which are separately linear in both of the two entries \ie coincides with the set of $L$-valued bilinear maps $\text{Multi}(V,W;L)$	.
			
			\item $\Vect$ is an Abelian category. 
			Being enriched in vector spaces, $\Vect$ is in particular an ab-enriched category, also called pre-additive.
			The further condition that one can always define a kernel and a cokernel for any linear maps, ultimately makes $\Vect$ Abelian.
		\end{itemize}
	\end{remark}

	\begin{remark}[Specializing the general theory of graded objects]\label{Rem:GVectCategoricalProps}
	Most of the structures explained in section \ref{Section:GradedObjects} specialize to the category $\Vect^{\mathbb{Z}}$; 
	in particular:
	\begin{itemize}
		\item The category $\Vect^{\mathbb{Z}}$ is complete and cocomplete. Limits are defined component-wise. 
		For instance there is a zero object
		\begin{displaymath}
			\mathbb{0} = (k \mapsto	\mathbb{0}\in \Vect)
			~,
		\end{displaymath}
		product and coproduct are given by
		\begin{displaymath}
			\prod_{j\in J }V_j = \left(
				k \mapsto \prod_{j \in J}(V_j)^k 
			\right)
		\end{displaymath}
		and
		\begin{displaymath}
			\bigoplus_{j\in J }V_j = \left(
				k \mapsto \bigoplus_{j \in J}(V_j)^k 
			\right)~.
		\end{displaymath}
		
		\item The category $\Vect^{\mathbb{Z}}$ is Abelian. Kernels and cokernels are defined component-wise:
		\begin{align*}
		 	\ker(f:V\to W) 
		 	=&~ 
		 	\big(k \mapsto \ker( f^k: V^k \to W^k) \big)
			\\
		 	\coker(f:V\to W) 
		 	=&~ 
		 	\big(k \mapsto \coker( f^k: V^k \to W^k) \big)
		 	~.
		\end{align*}
		
		\item The category $\Vect^{\mathbb{Z}}$ is $\Vect^{\mathbb{Z}}$-enriched.
		For any two $V,W \in \GVect$, the internal hom-functor is defined component-wise as the graded vector space\footnote{We ought to notice that the convention we are employing here is slightly non-standard.
	A popular choice (see \cite{nlab:internal_hom}) is to consider the graded vector space homogeneous graded maps, defined below in remark \ref{rem:neglectinginternalgrading}, as the internal hom-space for the category of graded vector spaces.}
		\begin{displaymath}
			[V,W] = \big( k \mapsto [V^k,W^k] \big) ~.
		\end{displaymath}
		In particular it is also monoidal closed \cite[ex 1.1.]{Delgado2015}.
		Considering the direct sums of all the components one gets an enrichment over $\Vect$.

		\item The grading set $\ZZ$ is a field and in particular a ring.
		 We thus may consider the \emph{shift by $j$} functor for any $j\in\mathbb{Z}$:
		\begin{displaymath}
			V[j] = \big( k \mapsto V^{k+j} \big)
			~.
		\end{displaymath}
	\end{itemize}	
	\end{remark}

	\begin{remark}\label{Remark:GvsAsTotalSpace}
		Similarly to what hinted in \ref{Remark:GradedObjectasCoproduct}, it is customary to understand a graded vector space as its \emph{associated total vector space} $V^\oplus \cong \bigoplus_{k\in\mathbb{Z}} V^k$ (see also remark \ref{rem:totaloplus})
		keeping implied the choice of a particular $\mathbb{Z}$-labelled decomposition, or \emph{grading}, and referring to elements completely contained in $V^p$ as "\emph{homogeneous} of degree $p$".
		\\
		In particular, the standard isomorphism $A\oplus 0 \simeq A$ implies that the category defined in Definition \ref{Def:GradedVectorSpaces} includes any $S$-graded vector space where $S$ is a countable set. 
		For instance, if $S=\{0,1\}$, any $S$-graded vector space $A$ can be regarded as a $\mathbb{Z}$-graded vector space simply by imposing that $A_i=0$ for all $i\notin S$. 
		\\
		Furthermore, $\GVect$ can be viewed a the (non-full) subcategory of $\Vect$ consisting of decomposable vector spaces $V^\oplus = \oplus_{k\in \mathbb{Z}} V^k$ with morphisms given by degree-preserving linear maps.
		According to this, one can easily define the direct sum and tensor product of two graded vector spaces out of the corresponding operators $\oplus$ and $\otimes$ on $\Vect$.
		Namely it suffices to specify a grading on the following "decomposable" ordinary vector spaces
		\begin{displaymath}
			\begin{aligned}
			(V\oplus W)^\oplus &= V^\oplus \oplus W^\oplus ~,\\
			(V\otimes W)^\oplus &= V^\oplus \otimes W^\oplus ~.
			\end{aligned}
		\end{displaymath}
		For instance, one can enforce that $(V\oplus W)_k = \oplus_{i+j=k} V^i\otimes W^k$ for any $k\in\mathbb{Z}$.
	\end{remark}
	\begin{notation}[Concentrated graded vector spaces]
		A graded vector space $V$ is said to be \emph{concentrated in degrees $S\subset \mathbb{Z}$} if $V^k=0$ for all $k \not\in S$.		
		Note that we can regard ordinary (in the sense of "ungraded") vector spaces as $\mathbb{Z}$-graded vector spaces concentrated in degree $0$.
	\end{notation}
	\begin{definition}[Bi-graded vector spaces]
	 	We call \emph{bi-graded vector space} any object in the $(\mathbb{Z}\times\mathbb{Z})$-graded vector space category $\Vect^{(\mathbb{Z}\times\mathbb{Z})}$. 
	 	In other words, is is a collection $\{V_k\}_{k\in\mathbb{Z}\times\ZZ}$ of vector spaces over $\mathbb{R}$ indexed by pairs of integers.		
		(The definition extends trivially	 to any multi-grading on generic graded objects).
	\end{definition}
	\begin{remark}\label{Rem:TotGradedVecSpace}
		To any bi-graded vector space one can associated a graded one through the \emph{total space} construction,
		\ie, for any given $\overline{V}=( s,t \mapsto V_{s, t})  \in \Vect^{\mathbb{Z}\times\mathbb{Z}}$, one can introduce
		\begin{displaymath}
			\tot(\overline{V}) = (q \mapsto \bigoplus_{s+t=q}V_{s,t})
			~.
		\end{displaymath}
	\end{remark}

	\begin{definition}[Homogeneous maps]
		It is customary to call any graded-morphism $f \in \Hom_{\Vect^\mathbb{Z}}(V,W[k])$ a \emph{homogeneous map from $V$ to $W$ in degree $k$}.
		In particular any graded morphism $V \to W$ is a \emph{degree $0$ homogeneous map}.
	\end{definition}
	\begin{remark}[Neglecting the internal grading]\label{rem:neglectinginternalgrading}
		The set of all homogeneous maps from $V\to W$ constitutes a $\mathbb{Z}$-graded object in the category $\Vect^\mathbb{Z}$.	In other words, there is a bi-graded vector space
		\begin{displaymath}
			\left(
				(k,s) \mapsto \Hom_{\Vect}(V^s,W[k]^s)
			\right)
			~.
		\end{displaymath}
		It is customary to neglect the "internal grading" given by the index $s$.
		Namely, one introduces the \emph{$\ZZ$-graded} vector space of homogeneous maps as
		\begin{displaymath}
			\underline{\Hom}_{\Vect^{\ZZ}}(V,W):=
			\left(
				k \mapsto \left(\Hom_{\Vect^{\ZZ}}(V,W[k])\right)^{\oplus}
			\right)
			~,
		\end{displaymath}
		where the superscript $\oplus$ means the direct sums of all the components of the graded vector space, \ie $V^\oplus=\oplus_{k\in\ZZ}V^k$, as defined in remark \ref{rem:totaloplus}.
		Accordingly, we will often refer to elements in $\underline{\Hom}_{\Vect^{\ZZ}}^k(V,W)$ simply as \emph{linear maps} in degree $k$.
		This is completely consistent with the fact that $f\in \underline{\Hom}^k_{\Vect^{\ZZ}}(V,W)$ is, by its very definition, a linear map between ordinary vector spaces $f:V^{\oplus}\to W^\oplus$ such that $\im(f\eval_{V^i})\subset W^{i+k}$.
	\end{remark}

	\begin{notation}[Shifted elements]\label{not:shiftedelements}
		Given $v \in V^{|v|}$, we denote by $v_{[k]}$ the element $v$ seen as an element in $V[k]^{|v|-k}$. 
		In other words,  $|v_{[k]}|= |v|-k$.
		An homogeneous map $f$ in degree $|f|=k$ from $V$ to $W$ is a graded morphism (\ie degree $0$ linear map) given by:
		\begin{displaymath}
			\morphism{f}
			{V}
			{W[k]}
			{v}
			{(f(v))_{[|f|]}} ~.
		\end{displaymath}	
		The latter will be often denoted as a linear map $f:V\to W$ implying the passage to $(f)^\oplus$ explained in remark \ref{rem:neglectinginternalgrading}.
		\\
		The action of the shift functor $[\ell]$ on an homogeneous map $f$ (graded morphism $f:V\to W[|f|]$) is given by
		\begin{equation}\label{eq:shiftHomoMaps}
			\morphism{f[\ell]}
			{V[\ell]}
			{W[|f|][\ell]}
			{v_{[\ell]}}
			{(f(v))_{[|f|][\ell]}}
			~.
		\end{equation}	
	\end{notation}

	\begin{definition}[Composition of homogeneous maps]\label{Def:compositionofhomogeneousmaps}
		Given  two homogeneous maps $f \in \Hom_{\GVect}(V,W[k])$ and $g \in \Hom_{\GVect}(W,X[\ell])$ in degree $k$ and $\ell$ respectively, we define their composition as the homogeneous map
		$$ g\circ f := g[k] \circ f \in \Hom_{\GVect}(V,X[k+\ell]) ~.$$
	\end{definition}
	\begin{remark}
		According to the notation introduced in remark \ref{rem:neglectinginternalgrading},
	one has
		\begin{displaymath}
			\underline{\Hom}^i_{\GVect}(V,W[l])\cong \underline{\Hom}^{i+l}_{\GVect}(V,W)
			~.
		\end{displaymath}
		Notice also that linearity ensures the following compatibility rule on homogeneous maps
		\begin{equation}\label{Eq:HomOfDirectSum}
			\underline{\Hom}_{\GVect}(A\oplus B, C) \cong \underline{\Hom}_{\GVect}(A,C)\oplus \underline{\Hom}_{\GVect}(B,C)
			~.
		\end{equation}	
	
	\end{remark}

	\begin{notation}[Graded maps vs. homogeneous maps]\label{Not:GradedvsHomogeneous-maps}
		From now on, we will only work in the category of graded vector spaces, therefore we will omit the subscript $\GVect$ when denoting the hom-sets.
		\\
		Namely, we denote as $\Hom(V,W)$ the graded vector space (due to the internal hom-functor) of \emph{graded morphisms} between $V$ and $W$.
		We also denote as $ \underline{\Hom}^k(V,W) = \Hom (V,W[k])$ the graded vector space of \emph{homogeneous maps in degree $k$}.
		Therefore we call $\underline{\Hom}(V,W)=\oplus_{k\in\mathbb{Z}}\Hom (V,W[k])$ the graded vector space of \emph{homogeneous maps in any degree}.
		\\
		Recall that $\Hom(V,W) = \underline{\Hom}^0(V,W)$, we stress that in our convention \emph{arrows between graded vector spaces} are always homogeneous maps in degree $0$.
	\end{notation}

\subsection{Monoidal structure and Koszul convention}\label{sec:bloodyKoszulConvention}
	Some slightly more subtle conventions appear when dealing with the induced monoidal structure on $\GVect$.
		
	Since $\Vect$ is a bi-monoidal\footnote{Namely with respect to the monoidal structures given by the tensor product $\otimes$ and by the direct sum $\oplus$.} category and $\mathbb{Z}$ is, in particular, a monoid, 
	proposition \ref{Prop:GradedMonoidalStructure} assures the existence of an induced monoidal structure on graded vector spaces. 
	The action on objects is given by
	\begin{displaymath}
		V\otimes W = \left(
			k \mapsto \bigoplus_{i+j=k} V^i\otimes W^j
		\right)
		~,
	\end{displaymath}		
	and the action of morphism $f:X\to Y$ and $g:W\to Z$ is given by
	\begin{displaymath}
		f \otimes g
		\left(
		k \mapsto \bigoplus_{i+j=k}\left(f^i\otimes g^j : X^i\otimes W^j \to Y^i\otimes Z^j
		\right)
		\right)~.
	\end{displaymath}
	The associator and unitor isomorphisms come from the associativity and unity isomorphisms in $\Vect$.
	Note that $\Vect$ is also braided, therefore, according again to proposition \ref{Prop:GradedMonoidalStructure}, there is a canonical braiding induced on the category $\Vect^\mathbb{Z}$ of graded vector spaces.
	 
	As it has been anticipated in remark \ref{Rem:AnticipazionediKoszul}, it is customary to consider on $\Vect^\mathbb{Z}$ a "twisted" version of the canonical Braiding called "Koszul braiding".
	
	\begin{definition}[Koszul braiding]\label{Def:KoszulBraiding}
		We call \emph{Koszul braiding} the braiding natural transformation defined on homogeneous elements by the isomorphism
		\begin{displaymath}
			\morphism{B_{V,W}}
			{V\otimes W}
			{W \otimes V}
			{x\otimes y}
			{(-)^{|x||y|}y\otimes x}		
		\end{displaymath}
		for any $V,W \in \Vect^\mathbb{Z}$.
	\end{definition}
	\begin{remark}[]\label{not:shortenBraiding}
		The braiding is clearly symmetric since 
		\begin{displaymath}
			B_{V,W} \circ B_{W,V} = \id_{W\otimes V} ~,\qquad	
			B_{W,V} \circ B_{V,W} = \id_{V\otimes W}	 
			~.
		\end{displaymath}
	\end{remark}
	
	\begin{notation}
		We will often omit the subscript $V,W$ in $B_{V,W}$ when there is no ambiguity about the domain $V\otimes W$ of the braiding operator.
	\end{notation}	
	
	\begin{remark}[Koszul convention]\index{Koszul convention}\label{rem:KoszulRuleofThumb}
		Informally, the choice of the Koszul braiding implies that the exchange of two homogeneous elements keeps track of their degree of the two exchanged elements.
		Namely, whenever two elements in degree $m$ and $n$ are swapped, a sign $(-)^{m n }$ is introduced. 
	\end{remark}
	%
	
	This convention has several tricky consequences mostly coming from the tensor product of shifted spaces.
	\begin{remark}[Suspension]\index{Suspension isomorphism}\label{Rem:suspension}
		A first important observation is that the three graded vector spaces, $V[1]$, $\mathbb{R}[1]\otimes V$, and $V \otimes \mathbb{R}[1]$, coincides components-wise
	\begin{displaymath}
		(V[1])_k \equiv (\mathbb{R}[1]\otimes V)_k \equiv (V \otimes \mathbb{R}[1])_k = V_{k+1}
		\qquad \forall k \in \mathbb{Z} ~.
	\end{displaymath}
	Nevertheless, we cannot assume that the above three spaces coincide as graded vector spaces because $\mathbb{R}[1]\otimes V$ and $V \otimes \mathbb{R}[1]$ are isomorphic through the braiding but, in principle, different.
	Therefore there is the freedom to choose which of the two spaces can be identified with $V[1]$.
	In the wording of \cite{Fiorenza2006}:
	\begin{quotation}
		\underline{ we adopt the convention that \emph{"degrees are shifted on the left"}}.
	\end{quotation}
	Namely we understand the following natural identification $V[1]\cong \mathbb{R}[1]\otimes V$ realized by the  \emph{suspension isomorphism}\footnote{Would be better to refer to it as \emph{suspension natural transformation}, $\susp:\mathbb{R}[1]\otimes \blank \Rightarrow [1]$.},
	\begin{displaymath}
		\isomorphism{\susp}
			{V[1]}
			{\mathbb{R}[1]\otimes V}
			{v_{[1]}}
			{1_{[1]}\otimes v}
			~.
	\end{displaymath}
	According to this convention, there is also a corresponding "suspension on the right" by composing the suspension with the Braiding:
	\begin{displaymath}
		\begin{tikzcd}[row sep=small]
			&
			\mathbb{R}[1]\otimes V \ar[dd,"B"]
			\\
			V[1] \ar[ru,"\susp"] \ar[dashed,dr] &
			\\
			&
			V \otimes \mathbb{R}[1]  
		\end{tikzcd}
	\end{displaymath}	
	that is :
	\begin{displaymath}
		\isomorphism{B \circ \susp}
			{V[1]}
			{V\otimes \mathbb{R}[1]}
			{v_{[1]}}
			{(-)^{|v|} v\otimes 1_{[1]}}	
			~.
	\end{displaymath}
	Basically, we are imposing that the shift functor is equivalent to left tensor product with $\mathbb{R}[1]$.
	Iterating the suspension, one obtains the following isomorphism
		\begin{displaymath}
		\isomorphism{\susp}
			{V[k]}
			{\mathbb{R}[k]\otimes V}
			{v_{[k]}}
			{1_{[k]}\otimes v}		
			~,
	\end{displaymath}	
	denoted again by $\susp$ with a slight abuse of notation.
	In particular we assume that
\begin{displaymath}
	\isomorphism{\text{susp}}
	{\RR[k]\otimes \RR[\ell]}
	{\RR[k+\ell]}
	{1_{[k]}\otimes 1_{[\ell]}}
	{(1_{[\ell]})_{[k]}=1_{[\ell+k]}}
	~.
\end{displaymath}
The latter implies the following canonical identification\footnote{In the sense that no extra sign arises.} $\RR[k]\otimes \RR[\ell]\equiv \RR[\ell+k] \equiv \RR[\ell]\otimes\RR[k]$ and, according to remark \ref{rem:GradingGivenByMonoid}, the following shifted spaces are also identified
		$$
			V[k][\ell]\equiv V[k+\ell] \equiv V [\ell][k] ~.
		$$
		We stress that the latter convention differs from some references in the bibliography, see \eg \cite{Delgado2018b}, where the shift is precisely defined by tensor multiplication on the left and  $V[k][\ell]$ and $V[\ell][k]$ are not identified but considered isomorphic (implying in particular an extra sign).
	 	\end{remark}
	More in general, the previous choice of a suspension implies the following isomorphism defined on any pair of graded vector spaces:
	\begin{definition}[D\'ecalage]\label{Def:DecIso}
		\begin{displaymath}
			\isomorphism{\dec}
			{V[k]\otimes W[\ell]}
			{(V\otimes W)[k+\ell]}
			{v_{[k]}\otimes w_{[\ell]}}
			{(-)^{\ell\cdot |v|}(v\otimes w)_{[k+\ell]}}
		\end{displaymath}
	\end{definition}
	\begin{remark}[Construction of $\dec$]\label{rem:decConstruction}
		The d\'ecalage operator introduced in definition \ref{Def:DecIso} can be explicitly 
	 constructed by a suitable composition of the "building block" introduced so far (suspension, braiding, associator and unitors).
		For instance, $dec$ could be given in four "steps" as follows:
		\begin{displaymath}
			\begin{tikzcd}
				V[k]\otimes W[\ell] \ar[r,"\susp\otimes\susp"] &[2em]
				\RR[k]\otimes V \otimes \RR[\ell]\otimes W \ar[r] &[-2em] 
				\cdots
				\\[-2em]
				v_{[k]}\otimes w_{[\ell]} \ar[r,mapsto]&
				1_{[k]}\otimes v \otimes 1_{[\ell]}\otimes w \ar[r,mapsto] &
				\cdots
				\\
				\cdots \ar[r,"\Unit\otimes B \otimes \Unit"] &
				\RR[k]\otimes \RR[\ell] \otimes V \otimes W \ar[r]&
				\cdots
				\\[-2em]
				\cdots \ar[r,mapsto] &
				(-)^{\ell\,|v|}~1_{[k]}\otimes 1_{[\ell]}\otimes v \otimes w \ar[r,mapsto] &
				\cdots 				
				\\
				\cdots \ar[r,"\susp^{-1} \otimes \Unit"] &
				\RR[k+\ell] \otimes (V \otimes W) \ar[r]&
				\cdots
				\\[-2em]
				\cdots \ar[r,mapsto] &
				(-)^{\ell\,|v|}~1_{[k+\ell]}\otimes (v \otimes w) \ar[r,mapsto] &
				\cdots 		
				\\
				\cdots \ar[r,"\susp^{-1}"] &
				(V \otimes W)[k+\ell] 
				\\[-2em]
				\cdots \ar[r,mapsto] &
				(-)^{\ell\,|v|}~(v \otimes w)_{[k+\ell]} 		
			\end{tikzcd}			
		\end{displaymath}			 
	(In the above construction, we understood the isomorphisms $\alpha,\lambda$ and $\rho$ given by the monoidal structure as identifications, see remark \ref{rem:Vectpdvs}). 
	\end{remark}
	Iterating the previous definition, one can introduce the \emph{d\'ecalage} isomorphism, defined on any $n$-tuple of graded vector spaces $(V_1,\dots, V_n)$ as:
	\begin{displaymath}
		\isomorphism{dec}
		{V_1[1]\otimes\dots\otimes V_n[1]}
		{(V_1\otimes\dots \otimes V_n)[n]}
		{v_{1[1]}\otimes\dots v_{n[1]}}
		{(-)^{\sum_{i=1}^{n}(n-i)|v_i|}(v_1\otimes\dots\otimes v_n)_{[n]}}
		~.
	\end{displaymath}

	\begin{remark}[Comparison with other conventions]
		Observe that with the building blocks in our hands one could make other slightly different choices on how to build a canonical isomorphism $V[k]\otimes W[\ell]\cong (V\otimes W)[k+\ell]$.
		For instance, one could also introduce the operator
		\begin{equation}\label{eq:decstrano}
			\morphism{\overline{\dec}}
			{V[k]\otimes W[\ell]}
			{(V\otimes W)[k+\ell]}
			{v_{[k]}\otimes w_{[\ell]}}
			{(-)^{\ell\cdot |v_{[k]}|}(v\otimes w)_{[k+\ell]}}
		\end{equation}
		which can be obtained substituting the second arrow in the diagram of remark \ref{rem:decConstruction}, given by $(\Unit_{\RR[k]}\otimes B_{V,\RR[\ell]} \otimes \Unit_W)$, with
		$(B_{\RR[k]\otimes V, \RR[\ell]}\otimes \Unit_W)$.
		Such a choice is rather common in the literature (see remark \ref{rem:comparesignmesswithliterature}). 
		In our convention, the operator $\overline{\dec}$ will only play a role in the definition of the tensor product of graded homogeneous maps (see remark \ref{Rem:TensorHomogeneousMaps} below).	
	\end{remark}

	\begin{remark}\label{rem:decVsBraiding}
		It is clear, from the characterization of the shift functor in terms of the suspension, that $[k]$ is not a "monoidal functor".
		Namely, one has that
		\begin{displaymath}
			(V\otimes W)[k] = \mathbb{R}[k]\otimes V \otimes W 
			~\cancel{\cong}~
			\mathbb{R}[k]\otimes V \otimes \mathbb{R}[k]\otimes W = V[k]\otimes W[k]
			~.
		\end{displaymath}	
		However, the d\'ecalage isomorphism defines a sort of "compatibility rule" between $[k]$ and $\otimes$. 
		In other terms, it determines an iso-natural transformation
		\begin{displaymath}
			\otimes \circ ( [k]\times [l] ) \Rightarrow [k+l]\circ \otimes
			~.
		\end{displaymath}
		Accordingly, the shift functor is not compatible with the braiding structure.
		This implies that, without loss of generality, the following diagram do not plain commutes:
		\begin{displaymath}
			\begin{tikzcd}
				{V[k]\otimes W[l]} \ar[r,"\dec"]\ar[d,"{B_{V[k],W[l]}}"'] 
				\ar[rd,phantom,"\cancel{\circlearrowleft}"]&
				{(V\otimes W)[k+l]} \ar[d,"{[k+l]B_{V,W}}"]
				\\
				{W[l]\otimes V[k]} \ar[r,"\dec"']&
				{W\otimes V [k+l]}			
			\end{tikzcd}		
		\end{displaymath}			
		but it is required an extra sign on the right vertical arrow in order to achieve the commutation:
		\begin{equation}\label{Eq:DecalageInterwiningPermutations}
			\begin{tikzcd}
				{V[k]\otimes W[l]} \ar[r,"\dec"]\ar[d,"{B_{V[k],W[l]}}"']&
				{(V\otimes W)[k+l]} \ar[d,"{(-)^{k l}[k+l]B_{V,W}}"]
				\\
				{W[l]\otimes V[k]} \ar[r,"\dec"']&
				{W\otimes V [k+l]}			
			\end{tikzcd}
			~.
		\end{equation}	
		The appearance of this extra sign determines the crucial property of $\dec$ to intertwine odd and even permutations of homogeneous vectors (see remark \ref{Rem:abstractpermutationofGVS}).
	\end{remark}

	According to this choice, there is also a corresponding non-trivial formula for the tensor product of graded maps.
	\begin{remark}[Tensor product of homogeneous maps]\label{Rem:TensorHomogeneousMaps}
		Recall at first that an homogeneous map $f:V\to W$ in degree $|f|=k$ is simply a graded morphism $f:V\to W[k]$.
		\\
		Consider now any two homogeneous maps $f \in \underline{\Hom}^{|f|}(V,W), f' \in \underline{\Hom}^{|f'|}(V',W')$ seen as graded morphisms in the above sense.
		Applying the monoidal product functor $\otimes$ defined on $\GVect$ one obtains the following graded morphism
		\begin{equation}\label{eq:CompletelyAccurateTensorOfMaps}
		 \morphism{f\otimes f'}
		 {V~\otimes~ V'}
		 {W[k]\otimes W'[k']}
		 {v \otimes v'}
		 {{f(v)}_{[|f|]}\otimes {f'(v')}_{[|f'|]}}
		\end{equation}	
		which -technically speaking- is not an homogeneous map\footnote{
			According to the previous conventions:
			$$
				f\otimes f' \in \Hom(V\otimes V', W[k]\otimes W[k'])
			\cong \underline{\Hom}^k(V\otimes V',W\otimes (W'[k']))
			~.
			$$
		} from $V\otimes V'$ to $W\otimes W'$.		
		\\
		To get an an honest homogeneous map in $\underline{\Hom}(V\otimes V',W\otimes W')$ one should postcompose the above map $f\otimes f'$ with a suitable isomorphism $W[|f|]\otimes W'[|f'|]\cong (W\otimes W')[|f|+|f'|]$.
		\\
	Accordingly, 
	we define, with a provisionally decorated notation, the tensor product of the above homogeneous map $f$ and $f'$ to be the homogeneous map 
	$f~\widetilde{\otimes}~f' \in \underline{\Hom}^{|f|+|f'|}(V\otimes V', W \otimes W')$ acting on homogeneous elements as
		\begin{equation}\label{eq:KoszulConvTensorProducts-appendix}
			\morphism{f~\widetilde{\otimes}~ f'}
		 {V\otimes V'}
		 {(W\otimes W')[k+k']}
		 {v \otimes v'}
		 {(-)^{|v||f'|}(f(v)\otimes f'(v'))_{[|f|+|f'|]}}
		 ~.
		\end{equation}
	In equation \eqref{eq:KoszulConvTensorProducts-appendix} is implied to choose a postcomposition of the above $\otimes$ with with $\overline{\dec}$ (see equation \eqref{eq:decstrano}) since $|{f(v)}_{[|f|]}|=|v|$.
	\\
	This choice is perfectly consistent with the \emph{Koszul convention} (see remark \ref{rem:KoszulRuleofThumb}), which, as a rule of thumb, prescribes that a sign $(-)^{|\alpha||\beta|}$ appears whenever two graded elements $\alpha$ and $\beta$.
	In case of \eqref{eq:KoszulConvTensorProducts-appendix}, the sign come from the permutation $(f,f',x,x')\mapsto(f,x,f',x')$.
		\\
		Equation \eqref{eq:KoszulConvTensorProducts-appendix} implies also a sign rule for the composition of tensor products of homogeneous maps:
		\begin{equation}\label{Eq:TensorHomogeneousMaps-appendix}
			(f'~\widetilde{\otimes}~ g') \circ (f ~\widetilde{\otimes}~ g) 
			= (-)^{|g'||f|}(f'\circ f)~\widetilde{\otimes}~ (g'\circ g)
			~.
		\end{equation}
	\end{remark}
	
	\begin{notation}[Dropping the decorated notation for $\widetilde{\otimes}$]\label{not:abuseotimes}
		It is an almost universally accepted practice to drop the decorated notation $\widetilde{\otimes}$ and denote with $\otimes$, with a slight abuse of notation, the tensor product of homogeneous map
		\begin{displaymath}
			\otimes:~
			\underline{Hom}^k(V,W)\times \underline{Hom}^{k'}(V',W')
			\to			
			\underline{Hom}^{k+k'}(V\otimes V',W\otimes W')
		\end{displaymath}
		whose image is given by equation \eqref{eq:KoszulConvTensorProducts-appendix}
	\end{notation}

	\begin{remark}[Considering homogeneous maps as morphisms in $\GVect$]\label{rem:homomapsasmorphism}
		Some sources in the literature propose to consider as morphisms in the category of graded vector spaces all homogeneous maps, therefore not only maps in degree $0$ as we are considering here.
		In that case, equation \ref{Eq:TensorHomogeneousMaps-appendix} must be imposed in order to define the action of functor $\otimes$ on arrows properly.
		\\
		We stress again that this is not the convention that we are employing here. In our diagrams, arrows must always be interpreted as degree $0$ maps.
	\end{remark}

	\begin{remark}[Shift of homogeneous maps]\label{rem:shiftHomogMapsTrickySign}
	There is a possible source of confusion coming from notation \ref{not:abuseotimes} that we shall clarify.
	\\
	Recall that we mentioned in remark \ref{Rem:suspension} that the convention "shift is suspension from the left" can be read as a natural transformation.
	The naturality condition implies that the below (on the right) square diagram ought to commute, in the category of graded vector spaces, for any arrow (below, on the left)
		\begin{displaymath}
			\begin{tikzcd}
				V \ar[d,"f"]&[2em] & 
				V[\ell] \ar[r,"\text{susp}"] \ar[d,"{f[\ell]}"]& 
				\mathbb{K}[\ell]\otimes V \ar[d,"{\id_{\mathbb{K}[\ell]}\otimes f}"]
				\\
				W[|f|] & &
				W[|f|][\ell] \ar[r,"\text{susp}"] &
				\mathbb{K}[\ell]\otimes W[|f|]
			\end{tikzcd}		
			~.
		\end{displaymath}
	Observe that the naturality is with respect to the monoidal structure of $\Vect^\ZZ$, hence, even if the arrow $f$ is an homogeneous map in degree $|f|$, the rightmost arrow $({\id_{\mathbb{K}[\ell]}\otimes f})$ has to be interpreted in the sense of equation \eqref{eq:CompletelyAccurateTensorOfMaps} \underline{not} in the sense of \eqref{eq:KoszulConvTensorProducts-appendix}.
	Chasing a generic homogeneous element $v\in V$ incorrectly, \ie interpreting $\otimes$ in the sense of notation \ref{not:abuseotimes}, would implies that
		\begin{displaymath}
			(f[\ell])(v_{[\ell]}) = (\id_{\mathbb{K}[\ell]} \otimes f) ~ (1_{[\ell]}\otimes v
			= (-)^{\ell k } 1_{[\ell]}\otimes f(v) = (-)^{\ell k } (f(v))_{[\ell]}
			~,
		\end{displaymath}
		which contradicts what we stated in equation \eqref{eq:shiftHomoMaps}.
		\\
		Notice that different conventions, see for example remark \ref{rem:homomapsasmorphism}, would lead to different interpretation of the symbol $\otimes$.
		\\
		We emphasize that in the body of the thesis the tensor product of homogeneous maps will be always interpreted in the sense of notation \ref{not:abuseotimes}
	\end{remark}		

	\begin{remark}[Suspension maps]\label{Rem:suspensionmaps}
		Observe that any homogeneous element $v$ of degree $|v|$ in the graded vector space $V$ is also an element in $V[i]$ with shifted homogeneous degree.
		Denoting $v_{[i]}$ the corresponding element in the shifted vector space one has $|v_{[i]}| = |v| -i$ because if $v \in V^k$ therefore $v_{[i]}\in (V[i])^{k-i}$.
		Accordingly, the identity morphism $\id_V=(k\mapsto \id_{V^k})$ can be regarded as an invertible degree $k$ homogeneous map between $V$ and $V[-k]$.
		\\
		Sometimes can be useful to introduce the \emph{suspension map} $\uparrow \in \underline{\Hom}^1(V,V[-1])$ and the \emph{desuspension map} $\downarrow  \in \underline{\Hom}^{-1}(V,V[1])$ (to not be confused with the suspension natural transformation defined in \ref{Rem:suspension}) defined as follows:
		\begin{displaymath}
			\uparrow(v) = 1_{[-1]}\otimes v = v_{[-1]}
			\qquad
			\downarrow(v)= 1_{[1]}\otimes v = v_{[1]}
			~.
		\end{displaymath}
		\\
		The suspension convention implies that $(\downarrow\otimes \downarrow)\circ (\uparrow \otimes \uparrow)= -\mathbb{1}$ or, more generally,
		\begin{displaymath}
			\downarrow^{\otimes n}\circ \uparrow^{\otimes n} 
			= 
			\uparrow^{\otimes n} \circ \downarrow^{\otimes n}
			=
			(-)^{\frac{n(n-1)}{2}}\mathbb{1}
			~.
		\end{displaymath}		 
		Finally, observe that if one chooses homogeneous maps as morphisms in the graded vector spaces category, as hinted in remark \ref{rem:homomapsasmorphism}, a graded vector space $V$ will be isomorphic to all its shifts $V[n]$ through the suspension and desuspension maps. 
	\end{remark}
	\begin{remark}[About the convention ${[k][\ell]=[\ell][k]}$]\label{rem:compsingShiftsvsSuspensions}
		Our convention about shift endofunctors to identify $[k][\ell]$, $[\ell][k]$ and $[\ell+k]$ for any $k,\ell \in \mathbb{Z}$ comes quite naturally when regarding a graded vector space as a functor $\mathbb{Z}\to \Vect$.
		In fact, since $\mathbb{Z}$ is in particular a group, the action $\mathbb{Z}\action\mathbb{Z}$ given by $+$ naturally translates to a family of endofunctor $\mathbb{Z}\to \mathbb{Z}$ where $\mathbb{Z}$ is interpreted as discrete category. 
		The shift functors on $\Vect^{\mathbb{Z}}$ are therefore obtained by precomposition with the latter functors.
		\\		
		More precisely, since $[k][\ell],[\ell][k]$ and $[\ell+k]$ are different functors, it is implied the following natural transformation
		\begin{displaymath}
			\begin{tikzcd}[ column sep=2em,row sep=-1ex]
				V[k][\ell] \ar[r,"\sim"] &
				V[k+\ell]=V[\ell+k] \ar[r,"\sim"] &
				V[\ell][k]
				\\
				(v_{[k]})_{[\ell]} \ar[r,mapsto] &
				v_{[k+\ell]} \ar[r,mapsto] &
				(v_{[\ell]})_{[k]}
			\end{tikzcd}
		\end{displaymath}				
		for any homogeneous vector $v\in V$. 
		\\
		This choice is perfectly compatible with the "(suspension) shift on the left" convention:
		\begin{displaymath}
				\begin{tikzcd}
				(v_{[k]})_{[\ell]} \ar[ddd,mapsto]&[-7em] 
				\\[-2em]
				&
				V[k][\ell] \ar[r,equal]\ar[d,"\text{susp}"] &
				V[k+\ell]=V[\ell+k] \ar[r,equal]\ar[d,"\text{susp}"] &
				V[\ell][k]\ar[d,"\text{susp}"]
				\\
				&
				\mathbb{K}[\ell]\otimes\mathbb{K}[k]\otimes V \ar[r,"\text{susp}^{-1}"] 		&
				\mathbb{K}[\ell+k]\otimes V \ar[r,"\text{susp}"] &
				\mathbb{K}[k]\otimes\mathbb{K}[\ell]\otimes V
				\\[-2em]
				1_{[\ell]}\otimes 1_{[k]}\otimes v \ar[rr,mapsto]&&
				1_{[\ell+k]}\otimes v \ar[r,mapsto]&
				1_{[k]}\otimes 1_{[\ell]} \otimes v					
				\end{tikzcd}
		\end{displaymath}
		We notice that is also common to make different choices. For instance in \cite[pag. 4]{Delgado2018} one can find the opposite \ie $V[k] \equiv \mathbb{R}[k]\otimes V$ and
		\begin{displaymath}
			\lambdaisomorphism{V{[k][\ell]}}
			{V{[\ell][k]}}
			{(v_{[k]})_{[\ell]}}
			{\mathit{\nrt{[(-)^{k \ell}]}}(v_{[\ell]})_{[k]}}
			~.
		\end{displaymath}
		The latter is in particular compatible with a different convention regarding the action of "suspended maps" (see \eg \cite{Fiorenza2006}). Namely, they define the shift functor $[\ell]$ on an homogeneous map (graded morphism $f:V\to W[|f|]$) as given by
		\begin{displaymath}
			\morphism{f[\ell]}
			{V[\ell]}
			{W[|f|][\ell]}
			{v_{[\ell]}}
			{\mathit{\nrt{[(-)^{|f|k}]}} (f(v))_{[|f|][\ell]}}
			~.
		\end{displaymath}
	In the previous two expressions we emphasized (with square brackets) the extra signs that are not present in our conventions.		
	\end{remark}

	\begin{remark}[Koszul convention - Take away message]\label{Rem:KoszulTAM}
		Summing up, the so-called "Koszul convention" rigorously consists of the following three choices:
		\begin{enumerate}
			\item Endowing the category $\GVect$ with the Koszul braiding defined in definition \ref{Def:KoszulBraiding}, that is different from the canonical one induced from the grading in $\Vect$ according to \ref{Rem:CanonicalInducedBraiding}.
			\item Adopting the convention that degrees are "shifted on the left", that is imposing that the shift functor $[k]$ coincides (is naturally identified) with $\mathbb{R}[k]\otimes-$. (See remark \ref{Rem:suspension}).
			\item Understanding the tensor product of homogeneous maps as obtained by postcomposition of the monoidal structure $\otimes$ with the operator $\overline{\dec}$ given in equation \eqref{eq:decstrano}. (See remark \ref{Rem:TensorHomogeneousMaps})
		\end{enumerate}
		Informally, these are subsumed by the following rule of thumb (already anticipated in \ref{rem:KoszulRuleofThumb}):
		\begin{quote}
		whenever two objects of degrees $p$ and $q$ are permuted, one has also to multiply by the sign prefactor $(-)^{p q}$.
		\end{quote}
		\end{remark}

	\begin{remark}[Permutations of tensor products]\label{Rem:abstractpermutationofGVS}
		Observe that, out of the symmetric braiding $B$, one can construct two natural transformations (an odd and an even one) pertaining to a permutation $\sigma\in S_n$, for any $n\geq 1$.
		Recall at first that any permutation can be expressed ---not uniquely--- as a composition of transpositions $\tau_i$, \ie, by permutations exchanging the $i$-th index with the $(i+1)$-th only.
		The sought transformations can be given on the generating set $\{\tau_i\}_{0< i < n}$ of $S_n$.
		The "even" transformation is given by
	\begin{displaymath}
		\mathclap{
		\morphism{B_{s_i}}
		{W_1\otimes\dots\otimes W_n}
		{W_1\otimes\dots\otimes W_{i-1}\otimes W_{i+1}\otimes W_i\otimes W_{i+2}\otimes\dots W_n}
		{w_1\otimes\dots\otimes w_n}
		{(-)^{|w_i||w_{i+1}|}w_1\otimes\dots\otimes 
		w_{i-1}\otimes w_{i+1}\otimes w_i\otimes w_{i+2}\otimes\dots w_n}	
		}	
	\end{displaymath}
	and the "odd" transformation is given by:
	\begin{displaymath}
		\mathclap{
		\morphism{P_{s_i}}
		{W_1\otimes\dots\otimes W_n}
		{W_1\otimes\dots\otimes W_{i-1}\otimes W_{i+1}\otimes W_i\otimes W_{i+2}\otimes\dots W_n}
		{w_1\otimes\dots\otimes w_n}
		{-(-)^{|w_i||w_{i+1}|}w_1\otimes\dots\otimes 
		w_{i-1}\otimes w_{i+1}\otimes w_i\otimes w_{i+2}\otimes\dots w_n}	
		~.
		}
	\end{displaymath}			
		More succinctly, this means that
	\begin{displaymath}
		\begin{aligned}
		B_{s_i} =& \mathbb{1}_{W_1}\otimes\dots\otimes\mathbb{1}_{W_{i-1}}\otimes B_{W_i,W_{i+1}}\otimes \mathbb{1}_{W_{i+2}}\otimes \dots \otimes \mathbb{1}_{W_{n-i-1}}
		\\
		P_{s_i} =& -	B_{s_i}
		\end{aligned}
		~.
	\end{displaymath}		
	Accordingly, to any given permutation $\sigma\in S_n$, decomposed in transpositions as $\sigma=s_{i_k}\dots s_{i_1}$, one can associate
	\begin{displaymath}
		\morphism{B_\sigma = B_{s_{i_k}}\dots B_{s_{i_1}}}
		{W_1\otimes\dots\otimes W_n}
		{W_{\sigma_1}\otimes\dots\otimes W_{\sigma_n}}
		{w_1\otimes\dots\otimes w_n}
		{\epsilon(\sigma,w) w_{\sigma_1}\otimes\dots\otimes w_{\sigma_n}}		
	\end{displaymath}	 
	and
	\begin{displaymath}
		P_\sigma= (-)^\sigma B_\sigma ~.
	\end{displaymath}
	We will often call this last two operators "even (resp. odd) permutators". They will be more extensively studied in appendix \ref{App:UnshuffleAtors}.
	\\	
	The coefficient $\epsilon(\sigma,w)$  arising from a generic permutation is called \emph{(symmetric) Koszul sign}.
	The corresponding sign $\chi(\sigma,w)=(-)^\sigma\epsilon(\sigma,w)$ pertaining to $P$ is called \emph{skew-symmetric Koszul sign}.
	See remark \ref{rem:aboutKoszulSign} for further details.	
	\end{remark}

	\begin{remark}[On the Koszul sign]\label{rem:aboutKoszulSign}
		The Koszul sign $\epsilon(\sigma,w)$  is the sign one gets by transposing elements in $w=w_1\otimes\dots\otimes w_n$ complying with the Koszul convention. 
		Namely it is the sign of the subpermutation $\bar{\sigma}$ obtained by restricting $\sigma$ to the subset of odd-degree elements only (see also remark \ref{rem:practicalComputKoszulSign}).
		For instance, the Koszul sign pertaining to the transposition of two elements $\tau_{1,2}\in S_2$, is given by
\begin{displaymath}
	\epsilon(\tau_{1,2};x_1,x_2)= (-)^{|x_1||x_2|}
\end{displaymath}
that is the sign given by the (Koszul) braiding. In other words $B_{\tau_{1,2}}\equiv B_{V,V}$.
One can then extend multiplicatively this definition to an arbitrary permutation using a decomposition into transpositions to get the general sign.
	\end{remark}

	\begin{notation}[Shortening the notation of the Koszul sign]
		Although the odd and even Koszul signs depend on the specific sequence of homogeneous vectors that are permuted, we will often omit the dependence on the list of graded vectors to shorten the notation.
		Everything should be clear from the context. 
		Namely we will often abuse the notation $\epsilon(\sigma;x_1,\dots,x_n)$ by writing $\epsilon(\sigma)$, and writing
 $\chi(\sigma):= \epsilon(\sigma)(-)^{\sigma}$ for the sign involved in the definition of $P_\sigma$.
	\end{notation}

\subsection{Tensor spaces}
	Considering the tensor product of a graded vector space $V$ with itself $n$-times, one gives rise to the so-called \emph{$n$-th tensors space}:
	\begin{displaymath}
		\otimes^n V = V^{\otimes n} = \underbrace{V\otimes\dots\otimes V}_{\text{n times}}
		~.
	\end{displaymath}		

	\subsubsection{Action of symmetric groups on Tensor spaces}\label{Section:ActionsonTensorSpaces}
		When specialized to $n$ copies of the graded vector space $V$, the iso-natural transformations $B_\sigma$ and $P_\sigma$, 
		defined in remark 	\ref{Rem:abstractpermutationofGVS} for a generic permutation $\sigma \in S_n$, can be regarded as canonical representations of the symmetric group $S_n$ on $V^{\otimes n}$ (for any $V\in \GVect$) induced by the Koszul Braiding.
		Namely, the graded vector space $\otimes^n V$, for any $n>1$, carries two natural actions of the group of permutations $S_n$: a canonical \emph{even representation}
		\begin{displaymath}
			\morphism{B}
			{S_n \times V^{\otimes n}}
			{V^{\otimes n}}
			{(\sigma,v_1\otimes\dots\otimes v_n)}
			{\epsilon(\sigma,v) v_{\sigma_1}\otimes\dots \otimes v_{\sigma_n}}
			~,
		\end{displaymath}
		and an \emph{odd representation}
		\begin{displaymath}
			\morphism{P}
			{S_n \times V^{\otimes n}}
			{V^{\otimes n}}
			{(\sigma,v_1\otimes\dots\otimes v_n)}
			{\chi(\sigma,v) v_{\sigma_1}\otimes\dots \otimes v_{\sigma_n}}
			~,
		\end{displaymath}
		where $\epsilon(\sigma,v)$ is the (even) Koszul sign and $\chi(\sigma,v)=(-)^\sigma\epsilon(\sigma,v)$ is the odd Koszul sign 
		(see remark \ref{rem:aboutKoszulSign}).
		Conventionally, we consider this actions as left-actions.
		\\
		Recall that this actions are natural in $V$ therefore they can be equally encoded as group morphisms $ S_n \to \text{Iso.Nat.}(\otimes^2,\otimes^2)$
	from the group of permutations to the group of natural isomorphisms from the functor $\otimes^n$ to itself.

		\begin{remark}[Practical computation of the Koszul sign]\label{rem:practicalComputKoszulSign}
			Consider an element $x_1\otimes\dots\otimes x_n$ in $\otimes^n V$. 
			Call $\zeta \in S_n$ the unique permutation that accumulates all odd degree elements on the left, 
			that is
			$$
				x_{\zeta_1}\otimes \dots \otimes x_{\zeta_n} = y_1\otimes\dots\otimes y_k\otimes z_1 \otimes \dots \otimes z_{n-k}
			$$ 
			with $|y_i|$ odd and $|z_i|$ even.
			\\
			Operatively, the permutation $\zeta$ can be constructed by scanning $x_1 \otimes \dots \otimes x_n$ element by element from the left and, whenever the considered $x_i$ is even, cycling it with all subsequent elements. 
			Since all of these permutations involve only swapping with an element in even degree, there is no Koszul sign to take into account.
			\\
			Given any permutation $\sigma \in S_n$, 
			there are two unique permutations $\bar{\sigma}\in S_k$ and $\bar{\bar{\sigma}} \in S_{n-k}$ defined by the commutativity of the following diagram:
			\begin{displaymath}
				\begin{tikzcd}
					x_1\otimes \dots \otimes x_n \ar[r,"B_\zeta"]
					\ar[d,"B_\sigma"]
					&
					y_1\otimes\dots\otimes y_k\otimes z_1 \otimes \dots \otimes z_{n-k}
					\ar[d,"B_{\bar{\sigma}}\otimes B_{\bar{\bar{\sigma}}}"]
					\\
					\epsilon(\sigma)~x_{\sigma_1}\otimes \dots \otimes x_{\sigma_n} \ar[r,"B_\zeta"]
					&
					(-)^{\bar{\sigma}}~y_{\bar{\sigma}_1}\otimes\dots\otimes y_{\bar{\sigma}_k}
					\otimes z_{\bar{\bar{\sigma}}_1} \otimes \dots \otimes z_{\bar{\bar{\sigma}}_{n-k}}
				\end{tikzcd}
				~.
			\end{displaymath}
			These are respectively the subpermutation of $\sigma$ involving odd degree elements and the one involving even degree elements only.
			One can conclude that $\epsilon(\sigma;x_1\otimes\dots\otimes x_n) = (-)^{\bar{\sigma}}$~.		
		\end{remark}		

		\begin{proposition}\label{Proposition:PermutingHomogeneousMaps}
			The permutation action on $(\underline{\Hom}(V,W))^{\otimes n}$ is compatible with the conjugate of the twist action on $\underline{\Hom}(V^{\otimes n},W^{\otimes m})$, \ie
			\begin{displaymath}
				 B_\sigma(f_1\otimes\dots\otimes f_n) =
				 B_\sigma \circ (f_1\otimes\dots \otimes f_n) \circ B_\sigma^{-1}
			\end{displaymath}
			where $B$ on the left-hand side denotes the action of the permutation group on $(\underline{\Hom}(V,W))^{\otimes n}$ and on the right-hand side denotes the action on $V^{\otimes n}$ and $W^{\otimes n}$.
		\end{proposition}
		\begin{proof}
			Without loss of generality, consider two homogeneous maps $f,g$ and two arbitrary graded homogeneous vectors $x,y$, one has
			\begin{align*}
				B \circ (f\otimes g) \circ ( x \otimes y) 
				=&~ 
				(-)^{|x||g|} B(f(x)\otimes  g(y)) =
				\\
				=&~ 
				(-)^{|f||g|+|f||y| + |x||y|} g(y)\otimes f(x) =
				\\
				=&~
				(-)^{|f||g|+ |x||y|} (g\otimes f) \circ (y\otimes x) =
				\\
				=&~
				(-)^{|f||g|} (g\otimes f) \circ B \circ (x \otimes y)
				~.
			\end{align*}
			Hence
			\begin{displaymath}
				B \circ (f \otimes g) = B(f\otimes g) \circ B
				~,
			\end{displaymath}
			where $B(f\otimes g)$ has to be carefully interpreted as the braiding on $\underline{\Hom}(V,W)^{\otimes 2}$ and it is not a postcompostion.
		\end{proof}
		
		\begin{lemma}[Composition of the Koszul signs]
			Given $\sigma,\tau \in S_n$, we have:
		\begin{displaymath}
			B_{\tau \sigma} (v_1\otimes\dots\otimes v_n) = \epsilon(\sigma;v_1,\dots,v_n) B_{\tau}(v_{\sigma_1}\otimes\dots\otimes v_{\sigma_n})
		\end{displaymath}
		and 
		\begin{displaymath}
			\epsilon(\tau \sigma;v_1,\dots,v_n)= \epsilon(\sigma;v_1,\dots,v_n)\epsilon(\tau;v_{\sigma_1},\dots,v_{\sigma_n})
			~.
		\end{displaymath}		
		\end{lemma}
		
		Specializing the d\'ecalage isomorphism of definition \ref{Def:DecIso} to several tensor products of a graded vector space with itself, it results the following expression
		\begin{displaymath}
			\morphism{\dec}
			{\otimes^n (V[1])}
			{(\otimes^n V)[n]}
			{(v_{1~[1]}\otimes \dots \otimes v_{n~[1]})}
			{(-)^{\sum_{k=1}^n |v_k|(n-k)}(v_1\otimes\dots\otimes v_n)_{[n]}}
			~.
		\end{displaymath}

		\begin{proposition}[D\'ecalage intertwines even and odd permutation actions]\label{Prop:DecalageOfPermutation}
			\begin{displaymath}
				\begin{tikzcd}
				(V[1])^{\otimes n} \ar[r,"dec"] \ar[d,"B_{\sigma}"'] 
				&
				(V^{\otimes n})[n] \ar[d,"(P_{\sigma}){[n]}"] 
				\\
				(V[1])^{\otimes n} \ar[r,"dec"] 
				&
				(V^{\otimes n})[n]
				\end{tikzcd}
			\end{displaymath}	
		\end{proposition}
		\begin{proof}
		 Specialize the equation given by diagram \eqref{Eq:DecalageInterwiningPermutations} to two copies of the same vector space and then iterate $n$ times.
		\end{proof}
		\begin{remark}
			Employing the language of the suspension maps rather then the suspension functors (see remark \ref{Rem:suspensionmaps}), the previous lemma can be subsumed by the following equation:
		\begin{displaymath}
			B_\sigma\circ \downarrow^{\otimes n} = \downarrow^{\otimes n} \circ P_\sigma~.
		\end{displaymath}		
		\end{remark}

		\subsubsection{(Skew)-symmetric tensor spaces}
		Due to the presence of the aforementioned canonical actions of $S_n$ on $V^{\otimes n}$, it is natural to consider the corresponding subspaces of coinvariants elements:
		\\
		\begin{paracol}{2}
			\begin{definition}[Symmetric tensor space]
				We call $\odot^n V$ the space of coinvariants with respect to the (even) canonical action of $S_n$ on $\otimes^n V$. 
				In other words, it is the unique subspace $\odot^n V \subset \otimes^n V$ such that:
				\begin{displaymath}
					\begin{tikzcd}[ampersand replacement=\&]
						\otimes^n V \ar[r,"B_\sigma"] \& \otimes^n V
						\\
						\odot^n V \ar[u,hook]\ar[ur,hook]
					\end{tikzcd}
				\end{displaymath}
				commutes for all $\sigma \in S_n$.
			\end{definition}		
		\switchcolumn
			\begin{definition}[Skew-symmetric tensor space]
				We call $\Lambda^n V$ the space of coinvariants with respect to the the twisted canonical action of $S_n$ on $\otimes^n V$.
				In other words, it is  the unique subspace $\Lambda^n V \subset \otimes^n V$ such that:
				\begin{displaymath}
					\begin{tikzcd}[ampersand replacement=\&]
						\otimes^n V \ar[r,"P_\sigma"] \& \otimes^n V
						\\
						\Lambda^n V \ar[u,hook]\ar[ur,hook]
					\end{tikzcd}
				\end{displaymath}
				commutes for all $\sigma \in S_n$.
			\end{definition}		
		\end{paracol}

		Exploiting the Abelian structure of $\GVect$, one can introduce the projectors on the (skew)-symmetric subspaces:
		\\
		\begin{paracol}{2}
			\begin{definition}[Symmetrizator]\label{Def:Symmetrizator}
				\begin{displaymath}
					\mathcal{S}= \left(\sum_{\sigma \in S_n} \frac{1}{n!} B_\sigma \right)
					: \otimes^n V \to \otimes^n V
					~.
				\end{displaymath}
			\end{definition}

		\switchcolumn
			\begin{definition}[Skew-symmetrizator]\label{Def:SkewSymmetrizator}
				\begin{displaymath}
					\mathcal{A}= \left(\sum_{\sigma \in S_n} \frac{1}{n!} P_\sigma \right)
					: \otimes^n V \to \otimes^n V
					~.
				\end{displaymath}
			\end{definition}		
		\end{paracol}			
		\begin{lemma}\label{lemma:splitsequencebrutta}
			There are several equivalent characterizations of $V^{\odot n}$ and $V^{\wedge n}$
			\begin{enumerate}
				\item 	$\odot^n V =\im(\mathcal{S}) = \ker(\mathcal{A})$
				$\qquad$
				$\Lambda^n V = \im(\mathcal{A}) = \ker(\mathcal{S})$~.
				\item $\Lambda^n V = \frac{V^{\otimes n}}{\ker(\mathcal{A})} = \frac{V^{\otimes n}}{\im(\mathcal{S})}$
				$\qquad$
				$\odot^n V = \frac{V^{\otimes n}}{\ker(\mathcal{S})} = \frac{V^{\otimes n}}{\im(\mathcal{A})}$.
				\item $V^{\otimes n} = V^{\wedge n}\oplus V^{\odot n}$.
			\end{enumerate}
		\end{lemma}
		\begin{proof}
			Consider the following sequence in the category of graded vector spaces
			\begin{equation}\label{eq:symskewsplitting-appendix}
				\begin{tikzcd}
					0 \ar[r, shift left =.5ex] 
					&[-5ex]
					\Lambda^n V \ar[r, shift left =.75ex,hookrightarrow,"N_a"] 
					&
					\bigotimes^n V 
					\ar[r,two heads, shift left =.75ex,two heads,"\pi_s"] 
					\ar[l,two heads, shift left =.75ex,two heads,"\pi_a"]
					&
					\bigodot^n V 
					\ar[r]			
					\ar[l, shift left =.75ex,hookrightarrow,"N_s"]
					&[-5ex] 0
				\end{tikzcd}~,
			\end{equation}				
			where the mapping are explicitly defined on homogeneous elements as follows
	\begin{equation}\label{eq:SymSkewOperatorsDef-appendix}
		\mathclap{	
		\begin{aligned}[c]
			\pi_s:&~ x_1\otimes\dots\otimes x_n \mapsto x_1\odot\dots\odot x_n ~,
			\\
			\pi_a:&~ x_1\otimes\dots\otimes x_n \mapsto x_1\wedge\dots\wedge x_n ~;
			\\
			{N_s}:&~ x_1\odot\dots\odot x_n \mapsto 
			\left(\sum_{\sigma\in S_n}x_{\sigma_1}\otimes\dots \otimes x_{\sigma_n}\right) ~,
			\\
			{N_a}:&~ x_1\wedge\dots\wedge x_n \mapsto 
			\left(\sum_{\sigma\in S_n}(-)^\sigma x_{\sigma_1}\otimes\dots \otimes x_{\sigma_n}\right)	~;
			\\
			\symAtor_{(n)}:&= \frac{1}{n!}N_s\circ \pi_s \equiv \left(\sum_{\sigma \in S_n} \frac{1}{n!} B_\sigma \right)~,
			\\
			\skewAtor_{(n)}:&= \frac{1}{n!}N_a\circ \pi_a \equiv \left(\sum_{\sigma \in S_n} \frac{1}{n!} P_\sigma \right) ~;
		\end{aligned}
			}
	\end{equation}
			and trivially extended by linearity.			
			The diagram in equation \eqref{eq:symskewsplitting-appendix} is a short exact sequence (both reading it from left to right than the converse).
			In fact, for any $x_1\wedge\dots\wedge x_n  \in \Lambda^n V$, one gets that
			\begin{align*}
				\pi_s \circ N_a (x_1\wedge\dots\wedge x_n) 
				=&~
				\pi_s\left(
					\sum_{\sigma \in S_n}(-)^\sigma x_{\sigma_1}\otimes\dots x_{\sigma_n}\right)
				=
				\\
				=&~
				\left(n!\sum_{\sigma \in S_n}(-)^\sigma\right) x_1\odot\dots \odot x_n
				=
				\\
				=&~
				0~.
			\end{align*}
			In particular one has $ \mathcal{S}\circ \mathcal{A} = \mathcal{A}\circ \mathcal{S} = 0$ and this proves the first claim.
			\\
			The sequence is both left and right splitting, \ie $\frac{1}{n!}\pi_a \circ N_a = \id_{\Lambda^n V}$ and $\frac{1}{n!}\pi_s \circ N_s = \id_{V^{\odot n}}$.
			The other two claims follows from the splitting lemma (see \cite{Weibel} or \cite{nlab:split_exact_sequence}).
		\end{proof}

		\begin{example}
			When $n=2$, one has
			\begin{align*}
				V\odot V &= \left\lbrace 
					v \in \otimes^2 V \; \left| \; B_{\sigma} v = v \quad \forall \sigma \in S_2
				\right.\right\rbrace =
				\\
				&=
				\frac{V\otimes V}{\langle v_1\otimes v_2 - (-)^{|v_1||v_2|}v_2\otimes v_1 \rangle_{v_i\in V}}~,
			\end{align*}	
			in the symmetric case, and	
			\begin{align*}
				V\wedge V  &= \left\lbrace 
					v \in \otimes^2 V \; \left| \; P_{\sigma} v = v \quad \forall \sigma \in S_2
				\right.\right\rbrace =
				\\
				&=
				\frac{V\otimes V}{\langle v_1\otimes v_2 + (-)^{|v_1||v_2|}v_2\otimes v_1 \rangle_{v_i\in V}}
				~,
			\end{align*}					
			in the skew-symmetric case.
		\end{example}		
	
	\begin{lemma}\label{Lemma:DecalageRestrictToSymTens}
	The D\'ecalage isomorphism restrict nicely on the symmetric and skew-symmetric subspaces
	\begin{displaymath}
		\begin{tikzcd}[column sep = huge]
			(V[1])^{\wedge n}
			\ar[r,"dec\big|_{V^{\wedge n}}"]\ar[d,hook]&
			(V^{\odot n})[n] 
			\ar[d,hook]
			\\
			(V[1])^{\otimes n}
			\ar[r,"dec"] &
			(V^{\otimes n}) [n] 
			\\
			(V[1])^{\odot n} \ar[u,hook]
			\ar[r,"dec\big|_{V^{\odot n}}"'] \ar[u,hook]&
			(V^{\wedge n})[n] 
		\end{tikzcd}
	\end{displaymath}
	\end{lemma}
	\begin{proof}
		Consider, without loss of generality, the case $n=2$.
		The commutation of the upper and lower squares follows from the following equation:
		\begin{displaymath}
			\begin{aligned}
				dec( v_{1~[1]}\otimes v_{2~[1]} &\pm (-)^{|v_{1~[1]}||v_{2~[1]}|} v_{2~[1]} \otimes v_{1~[1]} ) ~=
				\\
				=&
				(-)^{|v_1|} v_{1}\otimes v_{2} 
				\pm (-)^{|v_{2~[1]}||v_{1~[1]}|+|v_2|} v_{2}\otimes v_{1} =
				\\
				=& (-)^{|v_1|} (v_{1}\otimes v_{2} 
				\pm (-)^{|v_{2~[1]}||v_{1~[1]}|+|v_2|+|v_1|} v_{2}\otimes v_{1}) =
				\\
				=& (-)^{|v_1|} (v_{1}\otimes v_{2} 
				\pm (-)^{|v_{2}||v_{1}|+1} v_{2}\otimes v_{1}) =
				\\
				=& (-)^{|v_1|} (v_{1}\otimes v_{2} 
				\mp (-)^{|v_{2}||v_{1}|} v_{2}\otimes v_{1})
				~.
			\end{aligned}
		\end{displaymath}
	\end{proof}

\section{Graded algebras and coalgebras}\label{section:GradedAlgebras}
	In the course of the thesis, several different algebraic structures over graded vector spaces are considered.
	In this section, we review the basic definitions of (co)associative (co)algebra in the context of graded vector spaces and recall the notion of tensor (co)algebra associated with any given graded vector space.
		
	\subsection{(co)Associative (co)Algebras}\label{sec:abstractAlgebras}
In this subsection, we review the basic notions of algebras and coalgebras. Existing a notion of duality between these two concepts, we will present them in parallel columns.
What follows can be considered standard material; hence most of the proofs will be omitted.
Further details can be found, for instance, in \cite{Loday} and \cite{Manetti-website-coalgebras}.
		\\
		Let be $\mathbb{K}$ a field in characteristic $0$. Consider the category of graded vectors space over the field $\mathbb{K}$. 
		\globalcounter*
		
	\sidebyside{
		\begin{definition}[(Graded) Algebra]
			We call a \emph{Graded Algebra} a pair $(A,\mu)$ composed of a graded vector space $A$ and a graded morphism 
			\begin{displaymath}
				\morphism{\mu}
				{A\otimes A}
				{A}
				{a\otimes b}
				{a\cdot b}
			\end{displaymath}
			(\ie an homogeneous bilinear map in degree $0$).
		\end{definition}
	}{
		\begin{definition}[(Graded) coAlgebra]
			We call a \emph{Graded coAlgebra} a pair $(C,\Delta)$ composed of a graded vector space $C$ and a graded morphism (expressed in the Sweedler notation)
			\begin{displaymath}
				\morphism{\Delta}
				{C}
				{C\otimes C}
				{x}
				{\sum x_{(1)}\otimes x_{(2)}}
				~.
			\end{displaymath}
		\end{definition}
	}
	
	\sidebyside{
		\begin{definition}[Algebra morphism]
			Given two graded algebras $A=(A,\mu)$ and $A'=(A',\mu')$ a graded morphism $\phi:A\to A'$ is a \emph{algebra morphism} if the following diagram commutes:
			\begin{displaymath}
				\begin{tikzcd}[ampersand replacement=\&]
					A\otimes A \ar{r}{\mu} \ar{d}[swap]{\phi\otimes\phi} \& A \ar{d}{\phi}
					\\
					A'\otimes A' \ar{r}[swap]{\mu} \& A'				
				\end{tikzcd}
				~.
			\end{displaymath}
		\end{definition}
	}{
		\begin{definition}[coAlgebra morphism]
			Given two graded coalgebras $C=(C,\Delta)$ and $C'=(C',\Delta')$ a graded morphism $\psi:C\to C'$ is a \emph{coalgebra morphism} if the following diagram commutes:
			\begin{displaymath}
				\begin{tikzcd}[ampersand replacement=\&]
					C \ar{r}{\Delta} \ar{d}[swap]{\psi} \& C\otimes C \ar{d}{\psi\otimes\psi}
					\\
					C' \ar{r}[swap]{\Delta} \& C'\otimes C'				
				\end{tikzcd}
				~.
			\end{displaymath}
		\end{definition}	
	}
	\begin{remark}\label{Rem:CoalgebrasNotLinear}
		Observe that the space of (co)-algebra morphisms is not linear, \eg the sum of two coalgebra morphisms is not a morphism
		\begin{displaymath}
			\Delta (f + g ) = (f\otimes f + g \otimes g ) \Delta
			~. 
		\end{displaymath}
		However, one can define a deformation of a coalgebra morphism $f$ as the linear map $h$ such that $(f+h)$ is a coalgebra morphism.
		Accordingly, a deformation $h$ must satisfy the following equation:
		\begin{displaymath}
			\Delta h = \left( h\otimes h + h \otimes f + f \otimes h \right) \Delta ~.
		\end{displaymath}
	\end{remark}

	\sidebyside{
		\begin{definition}[Associative algebra]
			A given algebra $(A,\mu)$ is said \emph{associative} if the following diagram commutes
			\begin{displaymath}
				\begin{tikzcd}[ampersand replacement=\&]
					A \otimes A \otimes A \ar{r}{\mu\otimes\id_A}\ar{d}[left]{\id_A\otimes\mu}
					\& A \otimes A \ar{d}{\mu}
					\\
					A \otimes \ar{r}{\mu} \& A 
				\end{tikzcd}
				~.
			\end{displaymath}
		\end{definition}

	}{
		\begin{definition}[coAssociative coalgebra]
			A given coalgebra $(C,\Delta)$ is said \emph{coassociative} if the following diagram commutes
			\begin{displaymath}
				\begin{tikzcd}[ampersand replacement=\&]
					C \ar{r}{\Delta} \ar{d}[left]{\Delta}
					\& C \otimes C \ar{d}{\Delta\otimes\id_C}
					\\
					C \otimes C \ar{r}{\id_C \otimes \Delta}
					\&
					C\otimes C\otimes C
				\end{tikzcd}
				~.
			\end{displaymath}
		\end{definition}	
	}

	\sidebyside{
		\begin{notation}[Iterated multiplication]
		 On an associative algebra one will often consider the \emph{iterated product} 
		 $\mu^n: A^{\otimes(n+1)}\to A$ defined as follows:
		 \begin{displaymath}
		 	\mu^n := \mu(\mu^{n-1}\otimes \id)
		 \end{displaymath}
		 with $\mu^0= \id$ and $\mu^1=\mu$.
		\end{notation}
	
	}{
		\begin{notation}[Iterated comultiplication]
		 On a coassociative coalgebra one will often consider the \emph{iterated coproduct} 
		 $\Delta^n: C \to C^{\otimes(n+1)}$ defined as follows:
		 \begin{displaymath}
		 	\Delta^n := (\Delta\otimes\id\otimes\dots\otimes\id)\circ\Delta^{n-1}
		 \end{displaymath}
		 with $\Delta^0= \id$ and $\Delta^1=\Delta$.
		\end{notation}	
		\begin{lemma}\label{Lemma:IteratedCoalgebra}
			Given a coalgebra $(C,\Delta)$,  one has
			\begin{displaymath}
				(\Delta^p \otimes \Delta^q)\Delta = \Delta^{p+q+1}
				~.
			\end{displaymath}
			Furthermore, for any given coalgebra morphism $f:(C,\Delta_C)\to (D,\Delta_D)$, one has
			\begin{displaymath}
				f^{\otimes(n+1)} \Delta^n_C = \Delta^n_D \circ f
				~.
			\end{displaymath}
		\end{lemma}	
	}

	\sidebyside{
		\begin{definition}[(Anti)commutative algebra]
			A given algebra $(A,\mu)$ is said \emph{(anti)commutative} if the following diagram commutes
			\begin{displaymath}
				\begin{tikzcd}[ampersand replacement=\&]
					A \otimes A \ar{rr}{B_{A,A}}\ar{rr}[swap]{(-B_{A,A})} \ar{dr}[swap]{\mu}
					\&\& A\otimes A \ar{dl}{\mu}
					\\
					\& A \&
				\end{tikzcd}
				~,
			\end{displaymath}
			where $B$ and $P=-B$ denote the even and odd braiding operators (see notation \ref{not:shortenBraiding}).
		\end{definition}
	
	}{
		\begin{definition}[(Anti)cocommutative coalgebra]
			A given coalgebra $(C,\Delta)$ is said \emph{(anti)cocommutative} if the following diagram commutes
			\begin{displaymath}
				\begin{tikzcd}[ampersand replacement=\&]
					\&C \ar{dl}[swap]{\Delta} \ar{dr}{\Delta} \&
					\\
					C \otimes C \ar{rr}{B_{C,C}}\ar{rr}[swap]{(-B_{C,C})} \& \& C\otimes C 
				\end{tikzcd}
				~,
			\end{displaymath}
			where $B$ and $P=-B$ denote the even and odd braiding operators (see notation \ref{not:shortenBraiding}).			
		\end{definition}
	}

	\sidebyside{
		\begin{definition}[Unital algebra]
			A given algebra $(A,\mu)$ is said \emph{unital} if there exists an element $e\in A$, 
			called \emph{unit}, with $u:\mathbb{K}\to A$ the corresponding generated subspace, 
			such that the following diagram commutes
			\begin{displaymath}
			\hspace{-.1\textwidth}
				\begin{tikzcd}[ampersand replacement=\&]
					\mathbb{K}\otimes A \ar{r}{e\otimes\id_A}\ar{dr}[sloped,swap]{\sim} 
					\& A\otimes A \ar{d}{\mu}
					\& A \otimes \mathbb{K} \ar{l}[swap]{\id_A\otimes e}  \ar{dl}[sloped,swap]{\sim}]
					\\
					\& A \&
				\end{tikzcd}
				~.
			\end{displaymath}
		\end{definition}

	}{
		\begin{definition}[coUnital coalgebra]
			A given coalgebra $(C,\Delta)$ is said \emph{counital} if there exists a graded morphism $\epsilon:C\to\mathbb{K}$, called \emph{counit} (or \emph{augmentation map}),
			such that the following diagram commutes
			\begin{displaymath}
				\begin{tikzcd}[ampersand replacement=\&]
					\& C \ar{d}{\Delta}  \ar{dl}[sloped]{\sim} \ar{dr}[sloped]{\sim} \&
					\\
					\mathbb{K}\otimes C
					\& C \otimes C \ar{l}{\epsilon\otimes\id_C"} \ar{r}[swap]{\id_C\otimes\epsilon}
					\& C \otimes \mathbb{K}
				\end{tikzcd}
				~.
			\end{displaymath}		
		\end{definition}
	}

	\sidebyside{
		\begin{remark}
			Given two unital algebras $(A,\mu,e)$ and $(A',\mu',e')$, a \emph{unital morphism} is an algebra morphism $f:A\to A'$ such that $e' = f(e)$.
			Unital algebras form a non-full subcategory of the category of graded algebras with algebra morphisms.
		\end{remark}	

	}{
		\begin{remark}
			Given two counital coalgebras $(C,\Delta,\epsilon)$ and $(C',\Delta',\epsilon')$, a \emph{counital morphism} is a coalgebra morphism $f:C\to C'$ such that $\epsilon'\circ f = \epsilon$.
			Counital coalgebras form a non-full subcategory of the category of graded coalgebras with coalgebra morphisms.
		\end{remark}		
	}

	\begin{remark}
		Observe that the ground field $\mathbb{K}$ itself, in addition to constituting an unital associative algebra with its product, it is also a counital coassociative coalgebra for $\Delta(1_{\mathbb{k}})=1_{\mathbb{k}} \otimes 1_{\mathbb{k}}$.
	\end{remark}
	\begin{notation}
		From now on, we will often abbreviate associative unital graded algebras as "algebras" and coassociative counital graded coalgebras as "coalgebras".
	\end{notation}
	
	\sidebyside{
		\begin{definition}[Augmented Algebra]
		 A unital associative graded algebra is called \emph{augmented} when there is a given morphism of algebras $\epsilon:A\to \mathbb{K}$, called \emph{augmentation map}.
		 \footnote{In particular $\epsilon(1_A)=1_{\mathbb{K}}$.}
		\end{definition}

	}{
		\begin{definition}[coAugmented coAlgebra]
		 A counital coassociative graded coalgebra is called \emph{coaugmented} when there is a given morphism of coalgebras $e:\mathbb{K}\to C$ called \emph{coaugmentation map}.
		 \footnote{The coproduct on $\mathbb{K}$ is given  by $1_{\mathbb{K}}\mapsto 1_{\mathbb{K}}\otimes 1_{\mathbb{K}} $, in particular $\epsilon \circ e = \id_{\mathbb{K}}$.}
		 \note{A morphism of coaugmented coalgebras}
		\end{definition}
	}
	
	\sidebyside{
		\begin{theorem}
			If $A$ is augmented, then $A$ is canonically isomorphic to 
			$\overline{A}\oplus\mathbb{K}$
			where $\overline{A}=\ker(\epsilon)$ is called the \emph{augmentation ideal}.
		 Furthermore, the category of nonunital augmented associative algebras is equivalent to the category of unital augmented associative algebras.
		\end{theorem}
	
	}{
		\begin{theorem}
			If $C$ is coaugmented, then $C$ is canonically isomorphic to 
			$\overline{C}\oplus\mathbb{K}$
			where $\overline{C}=\ker(\epsilon)$.
			The reduced coproduct $\overline{\Delta}:\overline{C}\to\overline{C}\otimes\overline{C}$ is given by
			$$
				\overline{\Delta}(x):= \Delta(x) - x\otimes 1 - 1 \otimes x
				~.
			$$
		 Furthermore, the category of coaugmented coalgebras is equivalent to the category of coassociative coalgebras.
		\end{theorem}
		(proof, see \cite{Reinhold2019}.)
	}	
	
	\sidebyside{
		\begin{definition}[Nilpotent Algebra]
			A coassociative coalgebra $(A,\mu)$ is said \emph{nilpotent} if there exists some positive integer $n\in\mathbb{N}$ such that
			$\mu^n=0$.
		\end{definition}
	
	}{
		\begin{definition}[coNilpotent coalgebra]
		A coassociative coalgebra $(C,\Delta)$ is called \emph{conilpotent} if for any $x\in \overline{C}$ there exists some positive integer $n\in\mathbb{N}$ such that $\overline{\Delta}^n(x) =0 $.
		\end{definition}
		\note{Forse questa definitione è locally nilpotent}
	}
	
	\vspace{.5em}
	\sidebyside{
		Consider a graded vector space $W$ and a subcategory $\cat[A]$ of the category associative algebras.
		\begin{definition}[Free algebra generated by $W$ in the subcategory ${\cat[A]}$]\label{Def:FreeAlgebra}
			The \emph{Free associative algebra} generated by $W$ in the subcategory ${\cat[A]}$ is an associative algebra $\Free(W)\in\cat[A]$ equipped with a graded morphism $i:W\to \Free(W)$ which satisfies the following universal property:
			\\
			for any graded vector spaces morphism $f:W\to A$ where $A$ is an associative algebra in $\cat[A]$, there exists an unique  morphism $\tilde{f}:\Free(W)\to A$ of $\cat[A]$ such that the following diagram commutes:
			\begin{displaymath}
				\begin{tikzcd}[ampersand replacement=\&]
					W \ar{r}{i}\ar{dr}[swap]{\forall f} \& \Free(W)\ar{d}{\exists!\tilde{f}} \\
					\& A
				\end{tikzcd}
				~.
			\end{displaymath}
		\end{definition}
	
	}{
		Consider a graded vector space $W$ and a subcategory $\cat[C]$ of the category of  coassociative coalgebras.
		\begin{definition}[coFree coalgebra generated by $W$]\label{Def:FreeCoAlgebra}
			The \emph{coFree coassociative coalgebra} generated by $W$
			in the subcategory ${\cat[C]}$ is a coalgebra $\Free^c(W)\in\cat[C]$ equipped with a graded morphism $p:\Free^c(W)\to W$
			which satisfies the following universal property:
			\\
			for any graded vector spaces morphism $\varphi:C\to W$ where $C$ is a coassociative coalgebra in $\cat[C]$, there exists an unique morphism $\tilde{\varphi}:C\to\Free^c(W)$ of $\cat[C]$ such that the following diagram commutes:
			\begin{displaymath}
				\begin{tikzcd}[ampersand replacement=\&]
					C \ar{dr}{\forall\varphi} \ar{d}[swap]{\exists! \tilde{\varphi}} \& 
					\\
					\Free^c(W) \ar{r}{p} \&
					W
				\end{tikzcd}
				~.
			\end{displaymath}
		\end{definition}
	}
	
	\begin{proposition}
		The (co)free (co)associative (co)algebra on $W$ is uniquely determined up to isomorphisms.
	\end{proposition}

	\begin{remark}
		Putting the previous definitions side to side emphasize that the notion of coalgebra is obtained by formally dualizing, \ie reversing arrows in the diagrams, the notion of algebra.
	However, these two notions are not equivalent. 
	Despite the fact that a coalgebra gives an algebra by linear dualization,
	the dual of an algebra is not in general a coalgebra.
	The latter is guaranteed when considering finite dimensional algebras.
		\\	
		Heuristically, this asymmetry comes from the fact that, given any (possibly infinite-dimensional) graded vector space $V$, and considering the corresponding linear dual $V^\ast := \Hom(V,\mathbb{K})$,
		there exists a canonical map 
		\begin{displaymath}
			\morphism{\omega}
			{V^\ast\otimes V^\ast}
			{(V\otimes V)^\ast}
			{(f\otimes g)}
			{
				\left(
					\Almorphism{\omega(f\otimes g)}
					{V\otimes V}
					{\mathbb{K}}
					{x\otimes y}
					{f(x)g(y)}
				\right)
			}
		\end{displaymath}			
		while does not exists any natural map $(V\otimes V)^\ast \to V^\ast \otimes V^\ast$ in general.
		\\
		When $V$ is finite dimensional, $\omega$ is an isomorphism.
		If $(C,\Delta)$ is a coalgebra, then $(C^\ast,\Delta^\ast \circ\omega)$ is an algebra (no need of finiteness hypothesis).
		If $(A,\mu)$ is an algebra which is finite dimensional, then $(A^\ast,\omega^{-1}\circ \mu^\ast)$ is a coalgebra.
		See \cite{Loday,Manetti-website} for further details.
	\end{remark}
	
	\subsection{Tensor algebras}
		Consider a graded vector space $V$, by direct sum of all the possible tensor spaces one gets the so-called \emph{tensor algebra}.
	\begin{definition}[Tensor algebra]\label{Def:TensorAlgebra}
		Given a graded vector space $V\in \GVect$, we call \emph{tensor algebra of $V$} the graded algebra $(T(V),\otimes)$ where 
		\begin{displaymath}
			T(V) = \oplus_{k\geq 0} V^{\otimes k}
			~,
		\end{displaymath}
		conventionally it is assumed that $V^{\otimes 0} = \mathbb{R}$,
		and $\otimes$ denotes, with a slightly (yet standard) abuse of notation, the standard projection 
		\begin{displaymath}
			\morphism{\otimes}
			{V\times V}
			{V\otimes V = \frac{\Free(V\oplus V)}{\sim}}
			{(v_1,v_2)}
			{v_1\otimes v_2}
		\end{displaymath}
		which extends naturally to $T(V)$ yielding a bilinear map (concatenation product) and therefore defining an algebra structure on $T(V)$.
		\\
		We call \emph{reduced tensor algebra}, the subalgebra
		\begin{displaymath}
			\overline{T(V)} = \oplus_{k > 0} V^{\otimes k} \hookrightarrow T(V) ~.
		\end{displaymath}
	\end{definition}
	\begin{remark}
		Recall that we do not consider $V^{\otimes n}$ as the $\ZZ^n$-graded vector space
		\begin{displaymath}
			\left(
			(i_1,\dots,i_n) \mapsto	V_{i_1}\otimes\dots\otimes V_{i_n}
			\right)		
		\end{displaymath}
		but it is implicit in the definition of $\otimes$ to pass to the total space.
		Similarly, we do not regard $T(V)$ as the bigraded vector space
		\begin{displaymath}
			\left(
				(s,n)\mapsto (V^{\otimes n})_s = \bigoplus_{i_1+\dots+i_n=s} V_{i_1}\otimes\dots\otimes V_{i_n}
			\right)
		\end{displaymath}
		in other terms, the $n$-th component  of the Tensor algebra reads a follows
		\begin{displaymath}
			T^n(V) = \bigoplus_{q=0}^n (V^{\otimes q})_n \neq V^{\otimes n}
			~.
		\end{displaymath}
	\end{remark}	
	\begin{proposition}
		The tensor algebra of a graded vector space $V$ satisfies the following properties:
		\begin{enumerate}
			\item  $T(V)$ is an associative unital graded algebra, $\overline{T}(V)$ is an associative non-unital sub algebra of $T(V)$;
			\item $T(V)$ is augmented by $\epsilon (v_1,\dots,v_n)=0$ for any $ n\geq 1$ and $\epsilon(1)=1$, therefore $T(V)=\mathbb{R}\oplus \overline{T}(V)$.
			\item $T(V)$ is free in the category of unital associative algebra and $\overline{T(V)}$ is free in the category of non-unital associative algebra.
			That means that holds the following universal property:
			\begin{itemize}
				\item for any given graded unital associative algebra $A$, for each linear degree preserving map $f:V\to A$, there exists a unique homomorphism of unital graded algebras $\varphi: T(V)\to A$ that agrees on $V$ to $f$.
				I.e.\  the following diagram commutes
				\begin{displaymath}
					\begin{tikzcd}
						T(V) \ar[d] \ar[dr,"\exists ! F"] &
						\\
						V \ar[r,"\forall f"] & A
					\end{tikzcd}
					~,
				\end{displaymath}			
				where the vertical arrow denotes the standard projection.
			\end{itemize}
		\end{enumerate}
	\end{proposition}
	Similar constructions can be repeated on the subspaces of symmetric and skew-symmetric tensors.
	We present them side by side to emphasize the parallelism of the construction. However, unlike section \ref{sec:abstractAlgebras}, the two columns are not implying a duality in the categorical side. 
	\sidebyside{
		\begin{definition}[Symmetric tensor algebra]
			Given a graded vector space $V\in \GVect$, we call \emph{symmetric tensor algebra of $V$} the graded algebra $(S(V),\odot)$ where 
			\begin{displaymath}
				S(V) = \oplus_{k\geq 0} V^{\odot k}
				~,
			\end{displaymath}
			conventionally it is assumed that $V^{\odot 0} = \mathbb{R}$,
			and $\odot$ is the symmetrization of $\otimes$
		 	\begin{displaymath}
		 		\hspace{-.1\textwidth}
		 		\begin{tikzcd}[row sep=-1ex, column sep= small, ampersand replacement=\&]
		 			V\times V \ar{r}{\otimes} \ar[bend left = 30]{rr}{\odot} \&
		 			V\otimes V \ar{r}{\mathcal{S}} \&
		 			V \wedge V
		 			\\
		 			(v,w) \ar[mapsto]{r} \& 
		 			v\otimes w \ar[mapsto]{r} \&
		 			\displaystyle\sum_{\sigma \in S_2} B_\sigma(v\otimes w)
		 		\end{tikzcd}
		 	\end{displaymath}	
			which extends naturally to $S(V)$ yielding a bilinear map (concatenation product) and therefore defines an algebra structure on $S(V)$.
		\end{definition}
	}{
		\begin{definition}[Exterior algebra]
			Given a graded vector space $V\in \GVect$, we call \emph{exterior algebra of $V$} the graded algebra $(\Lambda(V),\wedge)$ where 
			\begin{displaymath}
				\Lambda(V) = \oplus_{k\geq 0} V^{\wedge k}
				~,
			\end{displaymath}
			conventionally it is assumed that $V^{\wedge 0} = \mathbb{R}$,
			and $\wedge$ is the skew-symmetrization of $\otimes$
		 	\begin{displaymath}
		 		\begin{tikzcd}[row sep=-1ex, column sep= small,ampersand replacement=\&]
		 			V\times V \ar{r}{\otimes} \ar[bend left = 30]{rr}{\wedge} \&
		 			V\otimes V \ar{r}{\mathcal{A}} \&
		 			V \wedge V
		 			\\
		 			(v,w) \ar[mapsto]{r} \& 
		 			v\otimes w \ar[mapsto]{r} \&
		 			\displaystyle\sum_{\sigma \in S_2} P_\sigma(v\otimes w)
		 		\end{tikzcd}
		 	\end{displaymath}	
			which extends naturally to $\Lambda(V)$ yielding a bilinear map (concatenation product) and therefore defines an algebra structure on $\Lambda(V)$.
		\end{definition}
	}

	\sidebyside{
		\begin{definition}
			We call \emph{reduced symmetric tensor algebra}, the subalgebra
			\begin{displaymath}
				\overline{S(V)} = \oplus_{k > 0} V^{\odot k} \hookrightarrow S(V) ~.
			\end{displaymath}		
		\end{definition}		
	}{
		\begin{definition}
			We call \emph{reduced exterior algebra}, the subalgebra
			\begin{displaymath}
				\overline{\Lambda(V)} = \oplus_{k > 0} V^{\wedge k} \hookrightarrow \Lambda(V) ~.
			\end{displaymath}		
		\end{definition}				
	}

	\sidebyside{	
 		\begin{proposition}[Symmetric tensor algebra]
			$S(V)= \frac{T(V)}{I_S}$ where $I_S$ is the two-sided homogeneous ideal generated by elements of the form $$v_1\otimes v_2 -(-)^{|v_1||v_2|}v_2\otimes v_1~.$$
		\end{proposition}			
	}{
 		\begin{proposition}[Skew-symmetric tensor algebra]
			$\Lambda(V)= \frac{T(V)}{I_A}$ where $I_A$ is the two-sided homogeneous ideal generated by elements of the form $$v_1\otimes v_2 +(-)^{|v_1||v_2|}v_2\otimes v_1~.$$		
		\end{proposition}
	}

	\sidebyside{
		\begin{proposition}
			$S(W)$ is the free graded commutative associative algebra over $W$, i.e it satisfies the analogue of the universal property but in the graded unital associative commutative algebras subcategory.
			That means that holds the following universal property:
			\begin{itemize}
				\item 			for any given graded unital associative commutative algebra $A$, for each linear degree preserving map $f:V\to A$, there exists a unique morphism of graded unital associative commutative algebras $\varphi: T(V)\to A$ that agrees on $V$ to $f$.
				I.e.\  the following diagram commutes
				\begin{displaymath}
					\begin{tikzcd}[ampersand replacement=\&]
						S(V) \ar{d} \ar{dr}{\exists ! F} \&
						\\
						V \ar{r}{\forall f} \& A
					\end{tikzcd}
					~,
				\end{displaymath}
				where the vertical arrow denotes the standard projection.
			\end{itemize}
		\end{proposition}
	
	}{
		\begin{proposition}
			$\Lambda(W)$ is the free graded anti-commutative associative algebra over $W$, i.e it satisfies the analogue of the universal property but in the
			graded unital associative anti-commutative algebras subcategory.
			That means that holds the following universal property:
			\begin{itemize}
				\item 			for any given graded unital associative anti-commutative algebra $A$, for each linear degree preserving map $f:V\to A$, there exists an unique homomorphism of graded unital associative anti-commutative algebras $\varphi: T(V)\to A$ that agrees on $V$ to $f$.
				I.e.\  the following diagram commutes			
				\begin{displaymath}
					\begin{tikzcd}[ampersand replacement=\&]
						\Lambda(V) \ar{d} \ar{dr}{\exists ! F} \&
						\\
						V \ar{r}{\forall f} \& A
					\end{tikzcd}
					~,
				\end{displaymath}
				where the vertical arrow denotes the standard projection.
			\end{itemize}
		\end{proposition}
	}
	\begin{notation}[Inclusion $N$]\label{Remark:ManettiNotation}
		For later reference, let us note that, from the very definition of the symmetrizator operator (definition \ref{Def:Symmetrizator} ), the direct sum of all symmetrizators
		\begin{displaymath}
				\mathcal{S}= \sum_{n\geq 0} \mathcal{S}_{(n)} ~ : T(V) \to T(V)
		\end{displaymath}
		admits the following trivial factorization
		\begin{displaymath}
			\begin{tikzcd}[column sep = huge]
				T(V) \ar[rr,"\mathcal{S}"] \ar[d,two heads,"\pi"] 
				&& T(V) \ar[dll,two heads,"\pi"]
				\\
				S(V) \ar[rr,"\sum_{n\geq 0} \frac{1}{n!} \id_{V^{\odot n}}"']
				&& S(V) \ar[u,hook,"N"]
			\end{tikzcd}
		\end{displaymath}
		where $\pi$ and $N$ are given in equation \eqref{eq:SymSkewOperatorsDef-appendix}.
		We omit the similar construction that can be made in the skew-symmetric case.
	\end{notation}
	\begin{remark}[Tensor functors]\label{Remark:TensorFunctors}
		One can easily assemble an endofunctor $T$ which maps a graded vector space $V$ to the corresponding tensor algebra assigning the action on any $f \in \Hom_{\text{Alg}}(V,W)$ to be 
		\begin{displaymath}
			\morphism{	T(f) := \sum_{n\geq o} f^{\otimes n} ~}
			{T(V)}{T(W)}
			{x_1\otimes\dots\otimes x_n}
			{f(x_1)\otimes\dots\otimes f(x_n)}
			.
		\end{displaymath}
		Similarly it is possible to introduce a \emph{symmetric tensor endofunctor} $S$ given on morphisms by
		\begin{displaymath}
			\morphism{	S(f)}
			{S(V)}{S(W)}
			{x_1\odot\dots\odot x_n}
			{f(x_1)\odot\dots\odot f(x_n)}
			.
		\end{displaymath}
		Note that the following factorization holds
		\begin{displaymath}
			\begin{tikzcd}
				S(V) \ar[r,hook,"N"] \ar[d,"S(f)"]& T(V) \ar[d,"T(f)"]
				\\
				S(W) \ar[r,hook,"N"] & T(W)
			\end{tikzcd}
		\end{displaymath}
		since $	(f\otimes f) \circ B = B \circ (f\otimes f)$ on any degree $0$ map, hence $S(f) = \pi \circ T(f) \circ N$.
		\\
		Furthermore, associativity of the tensor product of vectors ensures that $T$ and $S$ can be safely regarded as functors valued in the category of graded associative algebras.
	\end{remark}

	\subsection{Tensor coalgebras}\label{Section:TensorCoalgebras}
		In addition to the natural notion of associative algebra, the tensor space $T(V)$ (definition \ref{Def:TensorAlgebra}) can be also equipped with a canonical coalgebra structure.
		\begin{definition}[Deconcatenation coproduct]\label{def:Decon}
			Given a graded vector space $V\in \GVect$, we call \emph{deconcatenation coproduct} the graded morphism
			\begin{displaymath}
				\Delta = \sum_{n \geq 0} \sum_{a=0}^{n} \mathfrak{d}_{a,n-a} ~:~ T(V) \longrightarrow T(V)\otimes T(V)
			\end{displaymath}
			given by the following operators 
			\begin{displaymath}
				\morphism{\mathfrak{d}_{a,n-a}}
				{V^{\otimes n}}
				{V^{\otimes a}\otimes V^{\otimes n-a}}
				{x_1\otimes \dots \otimes x_n}
				{(x_1\otimes\dots\otimes x_a)\otimes (x_{a+1}\otimes \dots \otimes x_n)}
				.
			\end{displaymath}
			Similarly, we call \emph{reduced deconcatenation coproduct}, the graded morphism
			\begin{displaymath}
				\overline{\Delta} = \sum_{n > 0} \sum_{a=1}^{n-1} \mathfrak{d}_{a,n-a} ~:~ T(V) \longrightarrow T(V)\otimes T(V)
			\end{displaymath}
			~.
		\end{definition}
		\begin{notation}[Tensor coalgebra]
			We call \emph{tensor coalgebra of $V$} and \emph{reduced tensor coalgebra} the pairs $T^c(V)= (T(V),\Delta)$ and $\overline{T^c}(V)= (\overline{T(V)},\overline{\Delta})$  respectively.
			The latter 	acts on decomposable elements as follows
	\begin{displaymath}
		\morphism{\overline{\Delta}}
		{\overline{T(V)}}{\overline{T(V)}\otimes \overline{T(V)}}
		{x_1\otimes\dots \otimes x_n}
		{\displaystyle\sum_{i=1}^{n-1}(x_1\otimes\dots\otimes x_i)\otimes (x_{i+1}\otimes\dots\otimes x_n)}~.
	\end{displaymath}	 
		\end{notation}
		\begin{lemma}\label{Lemma:TensorCoalgebraProperties}
			The tensor coalgebra of a graded vector space $V$ satisfies the following properties:
			\begin{enumerate}
				\item  $T^c(V)$ is a coassociative counital graded coalgebra.
				 $\overline{T^c(V)}$ is a non-counital coassociative sub-coalgebra of $T(V)$, \ie
				 $$\overline{T^c(V)} = \oplus_{k > 0} V^{\otimes k} \hookrightarrow T^c(V) ~;$$
				\item $T^c(V)$ is coaugmented by the standard inclusion $u: \mathbb{R} \hookrightarrow T(V)$, therefore $T^c(V)=\mathbb{R}\oplus \overline{T^c(V)}$.
				\item $T^c(V)$ and $\overline{T^C(V)}$ are conilpotent, the iterated coproducts reads as follows:
					\begin{displaymath}
						\mathclap{
						\begin{aligned}
							\Delta^s(v_1\otimes\dots\otimes v_n) =&
							\sum_{0\leq i_1 < i_2 < \dots < i_{s} \leq n} 
							(v_1\otimes\dots\otimes v_{i_1})\otimes\dots\otimes (v_{i_s+1}\otimes\dots\otimes v_n)
							~;
							\\
							\overline{\Delta}^s(v_1\otimes\dots\otimes v_n) =&
							\sum_{1\leq i_1 < i_2 < \dots < i_{s} \leq n-1} 
							(v_1\otimes\dots\otimes v_{i_1})\otimes\dots\otimes (v_{i_s+1}\otimes\dots\otimes v_n)
							~.
						\end{aligned}
						}
					\end{displaymath}
			\end{enumerate}
		\end{lemma}		
		Tensor coalgebras enjoy the property to be cofree in a suitable subcategory of graded coalgebras:
		\begin{proposition}[Universal property of graded coalgebras]\label{Prop:UniversalPropertyGradedCoalgebras}
					${T^c(V)}$ is cofree in the category of counital coaugmented conilpotent coassociative coalgebras.
					\begin{enumerate}
						\item The universal property reads as follows:
						\\
						given a  conilpotent coaugmented coassociative coalgebra $(C,\Gamma)$ and a graded vector space $V$, 
						for any $f: C\to V$ graded morphism, there exists a unique morphism of graded coaugmented coalgebras $F:C\to T(V)^c$ 
						such that the following diagram commutes
						\begin{displaymath}
							\begin{tikzcd}
								C \ar[r,"\exists! F"] \ar[dr,"\forall f"']
								& T(V) \ar[d,"p"]
								\\
								& V \\
							\end{tikzcd}
							~,
						\end{displaymath}
						where $p:T(V)\to V$ denotes the standard projection.
						\\
						Similarly,
						$\overline{T^c(V)}$ is cofree in the category of conilpotent coassociative coalgebras.
						\item In the reduced case, the unique map $F$ associated to $f$ is explcitly given by the following commutative diagram in the category of graded coalgebras:
						\begin{displaymath}
							\begin{tikzcd}
								\overline{T}^c(V) \ar[rrd,sloped,"\sum_{n>0} f^{\otimes n} := \overline{T(f)}"]
								\\
								C \ar[u,"\sum_{n\geq 0}\Gamma^n"] \ar[rr,"F"] & & \overline{T^c(V)	}						
							\end{tikzcd}
							~.
						\end{displaymath}
							
					\end{enumerate}
	
					Namely, $F$ is given as follows
					\begin{displaymath}
						F = \sum_{n=1}^\infty (f^{\otimes n})\circ \Gamma^{n-1} 
						~: 
						C \to \overline{T}(C) \to \overline{T}(V)
						~.
					\end{displaymath}		
		\end{proposition}
		\begin{proof}
			The construction follows by noting that $\sum_{n\geq 0} \Gamma^n~: C \to \overline{T(C)}$ is a canonical coalgebra morphism and noting that lemma \ref{Lemma:IteratedCoalgebra} implies that $f^{\otimes n} \circ \Gamma^{n-1} = \Gamma^{n-1} \circ f$.
		\end{proof}
	
		Let us now focus on the case where the considered coalgebra is a tensor coalgebra itself, \ie consider $C=\overline{T}(U)$ for some graded vector space $U$.
		\begin{definition}[Corestriction]
			Given any homogeneous map $F\in \underline{\Hom}(\overline{T}(V),\overline{T}(W))$, we call \emph{$n$-th corestriction} the homogeneous map $f_n \in \underline{\Hom}^k(\overline{T}(V),W)$ given by
			\begin{displaymath}
				f_n := \pr \circ F \big\vert_{V^{\otimes n}}
			\end{displaymath}
			where $\pr: \overline{T}(W) \to W$ denotes the canonical projection.
			\\
			We can encode the process of taking all possible corestrictions as the single graded morphism
			\begin{equation}
				\morphism{P}
				{\underline{\Hom}^k (\overline{T}(V),\overline{T}(W))}
				{\underline{\Hom}^k (\overline{T}(V),W)\cong \bigoplus_{n>0}\underline{\Hom}^k(V^{\otimes n},W)}
				{F}
				{(f_1,f_2,\dots)}
				~.
			\end{equation}
		\end{definition}

	\begin{remark}[Unique lift to a coalgebra morphism]\label{Rem:LiftToMorphism}
		According to proposition \ref{Prop:UniversalPropertyGradedCoalgebras}, given a graded morphism $f:\overline{T}(U)\to V$ the corresponding unique coalgebra morphism $F: \overline{T}(U)\to\overline{T}(V)$ is given by
		\begin{displaymath}
			F(v_1\otimes\dots\otimes v_n) = 
			\sum_{s=1}^n \mkern-30mu
			\sum_{\mkern45mu 1\leq i_1<\dots< i_s =n}\mkern-40mu
			f(v_1\otimes\dots\otimes v_{i_1})\otimes\dots\otimes f(v_{i_{s-1}+1}\otimes\dots\otimes v_{i_s})
		\end{displaymath}
		and is called \emph{unique lift (to a graded morphism) of $f$}.
		The uniqueness condition of such a lift allows to encode it as a "lift" mapping
		\begin{equation}\label{Eq:LiftMorphismMapping}
			\morphism{L}
			{{\Hom}(\overline{T}{(V)},W)}
			{\Hom_{\text{coAlg}}(\overline{T^c}{(V)},\overline{T^c}{(W)})}
			{f}
			{\displaystyle\sum_{n>0}\sum_{s=1}^n 	\left(\sum_{j_1+\dots+j_s = n} f_{j_1}\otimes \dots \otimes f_{j_s}\right)~.}
		\end{equation}
		Note that the function is manifestly non-linear compatibly with remark \ref{Rem:CoalgebrasNotLinear}.				
	\end{remark}
	%
	%
	It is worth to note that the universal property of tensor coalgebras can be reviewed as an isomorphism of graded sets (they are not graded vector spaces):
	\begin{theorem}[Universal property as an isomorphism between hom-spaces]\label{Theorem:HomCoAlgISO}
		The following graded sets are in a one-to-one correspondence
		\begin{displaymath}
			{\Hom}_{\text{coAlg}}(\overline{T^c(V)},\overline{T^c(W)}) 
			\cong {\Hom}(\overline{T(V)},W) 
			\cong \bigoplus_{n>0} \Hom(V^{\otimes n},W)
			~.
		\end{displaymath}
		In particular, the isomorphism is induced from the corestriction $P$ according to the following diagram:
		\begin{displaymath}
			\begin{tikzcd}
				\Hom(\overline{T(V)},\overline{T(W)}) \ar[r,"P"] &
				\Hom(\overline{T(V)},W) \ar[d,equal]
				\\
				\Hom_{\text{coAlg}}(\overline{T^c(V)},\overline{T^c(W)}) \ar[u,hook]&
				\Hom(\overline{T(V)},W) \ar[l,"L"]
			\end{tikzcd}
			~.
		\end{displaymath}
	\end{theorem}
	\begin{proof}
		The fact that $P \circ L = \id_{\Hom(\overline{T}(V),W)}$ follows by noting that the term between brackets in equation \eqref{Eq:LiftMorphismMapping} is an operator $V^{\otimes n}\to V^{\otimes s}$, therefore $P\circ F = \sum_{n>0} f_n$.
		On the other hand, lemma \ref{Lemma:IteratedCoalgebra} implies
		\begin{align*}
			L (P (F)) 
			=&~
			\sum_{n>0} f^{\otimes n}\circ \Delta^{n-1} = 
			\\
			=&~
			\sum_{n>0} p^{\otimes n} \circ F^{\otimes n} \circ \Delta^{n-1} =	
			\\
			=&~
			\sum_{n>0} p^{\otimes n} \circ \Delta^{n-1} \circ F =	
			\\
			=&~ F ~,		
		\end{align*}
		since $p^{\otimes n} \circ \Delta^{n-1}$ coincides with the projector $T(V) \to V^{\otimes n}$. Hence $L\circ P = \id_{\Hom(\overline{T}(V),\overline{T}(W))}$.	
	\end{proof}
	The upshot is that any morphism of graded coalgebras is entirely determined by its corestriction.
	In particular, isomorphisms can be characterized as follows:
	\begin{lemma}\label{Lemma:CoalgebraIsosCondition}
		A coalgebra morphism $F: \overline{T}(V) \to \overline{T}(W)$ is invertible if and only if its first corestriction $f_1 = \pr\circ F \eval_V$ is a graded vector spaces isomorphism.
	\end{lemma}
	\begin{proof}
		Note first that, by its very definition, the unique lift of the canonical projection $\pr: T(V) \to V$ is the identity $\id\vert_{T(V)}$.
		Hence, given a coalgebra isomorphism $F: \overline{T(V)} \to \overline{T(W)}$ with inverse denoted by $G$, one has
		\begin{displaymath}
			G\circ F = L (g \circ F) = L (\pr) = \id_{\overline{T(W)}}
			~.
		\end{displaymath}
		According to remark \ref{Theorem:HomCoAlgISO}, $L$ is injective. Therefore:
		\begin{displaymath}
			\id_V = \pr \eval_V = g\circ F \eval_V = g_1 \circ f_1
		\end{displaymath}
		hence $f_1$ is an isomorphism with inverse $g_1$.
		\\
		On the converse, one can check that, given a coalgebra morphism $F$ with invertible first corestriction $f_1$, the lift $L(g)$ with respect to $g_n = f_1^{-1} \circ f_n \circ (f_1^{-1})^{\otimes n}$, yields an inverse of $F$.		
	\end{proof}
	We will see in appendix \ref{App:RNAlgebras} how one can define a similar universal property also in the case of homogeneous maps in degree different than zero.
	
	\subsubsection{Symmetric tensor coalgebras}
		We already seen how a given tensor space $T(V)$ can be decomposed  in a symmetric and skew-symmetric part, \ie $T(V) = S(V) \oplus \Lambda(V)$.
		Let us now focus on the graded symmetric tensor space (we omit the analogue constructions that can be straightforwardly carried out in the skew-symmetric case taking care of extra signs).
		\\
		$S(V)$ admits a canonical coassociative coalgebra  structure:
	\begin{definition}[Unshuffle coproduct]\label{def:unshuffleDecon}
			Given a graded vector space $V\in \GVect$, we call \emph{unshuffle coproduct} the graded morphism
			\begin{displaymath}
				\Xi = 	\sum_{n \geq 0} \sum_{a=0}^{n} \dfrac{1}{n!} (\pi\otimes \pi) \circ \mathfrak{k}_{a,n-a}
				\circ N ~:~ S(V) \longrightarrow S(V)\otimes S(V)
			\end{displaymath}
			given by the following operators 
			\begin{displaymath}
					\mathfrak{k}_{a,n-a}:= \mathfrak{d}_{a,n-a} \circ B_{a,n-a}~:~ T(V) \longrightarrow T(V) \otimes T(V)
			\end{displaymath}
			where  $\mathfrak{a}_{i,j}$ is the operator of definition \ref{def:Decon}, $\pi$ and $N$ have been introduced in remark \ref{Remark:ManettiNotation} and $B_{a,n-a}$ denotes the sum of all possible unshuffle permutation (unshuffleators, see appendix \ref{App:UnshuffleAtors}).
			Similarly, we call \emph{reduced unshuffle coproduct}, the graded morphism
			\begin{displaymath}
				\overline{\Xi} = 	\sum_{n > 0} \sum_{a=1}^{n-1} \dfrac{1}{n!} (\pi\otimes \pi) \circ \mathfrak{k}_{a,n-a}
				\circ N ~:~ \overline{S(V)} \longrightarrow \overline{S(V)} \otimes \overline{S(V)}
				~.
			\end{displaymath}
		\end{definition}
		\begin{notation}[Symmetric tensor coalgebra]
			We call \emph{symmetric tensor coalgebra of $V$} and \emph{reduced symmetric tensor coalgebra} the pairs $S^c(V)= (S(V),\Xi)$ and $\overline{S^c}(V)= (\overline{S(V)},\overline{\Xi})$  respectively.
		\end{notation}	
		\begin{notation}[Explicit action of the unshuffle coproduct]
			The action of $\Xi$ on homogeneous elements is given explicitly by
			\begin{displaymath}
				\Xi(v_1\odot\dots\odot v_n) = \sum_{i=0}^{n}
				\sum_{\sigma \in S_{i,n-i}} \epsilon(\sigma) 
				(v_{\sigma_1}\odot\dots\odot v_{\sigma_i})
				\otimes
				(v_{\sigma_{(i+1)}}\odot\dots\odot v_{\sigma_n})
			\end{displaymath}
			where $S_{i,n-i}\subset S_{n}$ denotes the set of $(i,n-i)$ unshuffles (see appendix \ref{App:UnshuffleAtors}).
			
		\end{notation}

	\begin{proposition}
		$S^c(V)=(S(V),\Xi)$ forms a cocommutative coalgebra, it forms a sub-coalgebra of  the tensor coalgebra $T^c(V)$ with respect to the injection $N$ defined in remark \ref{Remark:ManettiNotation}, \ie
		\begin{displaymath}
			\begin{tikzcd}
				N:~S^c(V) 	\ar[r,hook] & T^c(V)~.
			\end{tikzcd}
		\end{displaymath}
		A similar result holds in the reduced case for $(\overline{S(V)},\overline{\Xi})$.
	\end{proposition}
	\begin{proof}
		Cocommutativity follows from
		\begin{displaymath}
			\begin{aligned}
			B \circ \Xi =&~ 
			\sum_{n \geq 0} \sum_{a=0}^{n} \dfrac{1}{n!} (\pi\otimes \pi) \circ \mathcal{C}_{(n)}^a \circ \mathfrak{d}_{a,n-a}\circ B_{a,n-a}	\circ N =
			\\
			 =&~ 
					\sum_{n \geq 0} \sum_{a=0}^{n} \dfrac{1}{n!} (\pi\otimes \pi) \circ \mathfrak{d}_{n-a,a}\circ B_{n-a,a}	\circ \cancel{\mathcal{C}_{(n)}^a} \circ N = 
					\\
				=&~
				\Xi
				~,
			\end{aligned}
		\end{displaymath}
		where $\mathcal{C}_{(n)}^a$ denotes the cyclic permutation of $n$ elements $a$ times (see appendix \ref{App:UnshuffleAtors}).
		The latter can be easily ascertained by inspection on homogeneous elements.
		\\
		The property of being coassociative follows by proving that $N$ is an injective map that preserve coproducts, i.e $\Delta\circ N = (N\otimes N) \circ \Xi$. This follows from the next equation
		\begin{displaymath}
				\mathfrak{d}_{a,n-a} \circ N =
				(N \otimes N) ~ \circ
				\left[\pi\otimes \pi \left(\dfrac{1}{n!} \circ\mathfrak{k}_{n,n-a}\right) \circ N \right]		
		\end{displaymath}
		which can be again ascertained by inspection on homogeneous elements
		\begin{displaymath}
			\begin{aligned}
				\mathfrak{d}_{a,n-a} \circ N &(x_1\odot\dots\odot x_n) \equal{} \\
				\equal{Rem \ref{Remark:ManettiNotation}}&
					n! ~\mathfrak{d}_{a,n-a} \circ \mathcal{S}_{(n)}\circ \mathcal{S}_{(n)} (x_1\otimes\dots \otimes x_n)
				=
				\\
				\equal{Rem: \ref{Rem:UnshufflesAsCoset}}&
				a!~(n-a)!~ ( \mathcal{S}_{(a)}\otimes\mathcal{S}_{(n-a)})\circ \mathfrak{d}_{a,n-a} \circ B_{a,n-a} \circ \mathcal{S}_{(n)}
				(x_1\otimes\dots \otimes x_n)
				=
				\\
				\equal{Rem \ref{Remark:ManettiNotation}}&
				(N \otimes N) \circ \left[
				\dfrac{1}{n!} ~ (\pi\otimes \pi) \circ  \mathfrak{d}_{a,n-a} \circ B_{a,n-a} \circ N
				\right] (x_1\odot\dots\odot x_n)
				~.
			\end{aligned}
		\end{displaymath}
		In other words, the restriction of $\Delta$ to $S(V)$ has image in $S(V)\otimes S(V)$ and it is invariant under the even action of the symmetric group. In particular one has $B_\sigma \circ \Delta = \Delta$ for all $\sigma \in S_2$.
	\end{proof}

	Therefore, the conilpotency and coaugmentation property of $T^c(V)$ stated in lemma \ref{Lemma:TensorCoalgebraProperties} are automatically reflected on $S^c(V)$ and $T^c(V)$.
		Similarly, they satisfy an analogous universal property:	
		\begin{proposition}[Universal property of cocommutative graded coalgebras]\label{Prop:UniversalPropertyCocommutativeGradedCoalgebras}
					${S^c(V)}$ is cofree in the category of counital coaugmented conilpotent coassociative coalgebras.
					\begin{enumerate}
						\item The universal property reads as follows:
						\\
						given $(C,\Gamma)$ a cocommutative conilpotent coaugmented coassociative coalgebra and given a graded vector space $V$, for each graded morphism $f: C\to V$ there exists a unique morphism of graded coaugmented coalgebras $F:C\to S(V)^c$ 
						such that the following diagram commutes
						\begin{displaymath}
							\begin{tikzcd}
								C \ar[r,"\exists! F"] \ar[dr,"\forall f"]
								& S(V) \ar[d,"p"]
								\\
								& V \\
							\end{tikzcd}
							~,
						\end{displaymath}
						where $p:S(V)\to V$ denotes the standard projection.
						Similarly,
						$\overline{S^c(V)}$ is cofree in the category of cocommutative conilpotent coassociative coalgebra.
						\item In the reduced case, the unique map $F$ associated to $f$ is explicitly given by the following commutative diagram in the category of graded coalgebras:
						\begin{displaymath}
							\begin{tikzcd}[column sep = huge]
								\overline{S(V)} \ar[hook,d,"N"] \ar[drr,"S(f)"]&&
								\\
								\overline{T(V)}\ar[drr,"T(f)"] && \overline{S(V)}\ar[d,hook,"N"]
								\\
								C \ar[rr,"L(f)"'] \ar[uu,bend left=60,"\sum_{n\geq 0}\frac{1}{n!}~\pi\circ\Gamma^n"] \ar[u,"\sum_{n\geq 0}\Gamma^n"']\ar[drr,"f"']
								&& \overline{T(V)} \ar[d]
								\\
								&& V
							\end{tikzcd}
							~.
						\end{displaymath}	
						Namely, $F$ is given by
						\begin{displaymath}
							F= \sum_{n=1}^\infty  \frac{\pi}{n!}\circ(f^{\otimes n})\circ \Gamma^{n-1} : C \to \overline{S(C)} \to \overline{S(V)}
							~.
						\end{displaymath}											
					\end{enumerate}
		\end{proposition}
		\begin{proof}
			The leftmost factorization follows by noting that the canonical coalgebra morphism $\sum_{n\geq 0} \Gamma^n$ has image in permutation invariant elements when the coalgebra $C$ is cocommutative, hence
			\begin{displaymath}
				\mathcal{S} \circ \left(\sum_{n\geq 0} \Gamma^n \right)
				= 
				N \circ \frac{\pi}{n!} \circ  \left(\sum_{n\geq 0}\Gamma^n\right) 
				.
			\end{displaymath}
			The uppermost square has been explained in remark \ref{Remark:TensorFunctors} and the other two triangles are precisely given by proposition \ref{Prop:UniversalPropertyGradedCoalgebras}.
		\end{proof}

		Continuing the parallelism with what explained in the previous subsection, we now focus in the case in which the coalgebra considered is the reduced symmetric coalgebra itself, \ie $C=\overline{S(V)}$
	\begin{remark}[Unique symmetric lift]\label{Rem:SymmetricLiftToMorphism}
		According to proposition \ref{Prop:UniversalPropertyCocommutativeGradedCoalgebras}, given a graded morphism $f:\overline{S(V)}\to W$ the corresponding unique coalgebra morphism $F: \overline{S^c(V)}\to\overline{S^c(W)}$ is given by

		\begin{displaymath}
			\mathclap{
			\begin{aligned}
			F(x_1&\odot\dots\odot x_n) =
			\\
			=&~ 
			\sum_{s=1}^n \mkern-30mu
			\sum_{\qquad\substack{i_1+\dots+i_s = n \\ 0<i_1\leq \dots \leq i_s}}\mkern-60mu
			\sum_{\qquad\qquad\sigma \in \ush{i_1,\dots,i_s}^{<}}\mkern-60mu
			f_{i_1}(x_{\sigma_1}\odot\dots\odot x_{\sigma_i})\odot \dots \odot
			f_{i_s}(x_{n-i_s -1}\odot \dots \odot x_n)			
			\end{aligned}}
		\end{displaymath}
		and is called \emph{unique lift (to a graded morphism) of $f$}.
		The uniqueness condition of such lift allows to encode it as a \emph{symmetric lift mapping}:
		
		\begin{equation}\label{Eq:LiftMorphismMapping-Sym}
		  	\hspace{-.1\textwidth}
			\morphism{L_{\sym}}
			{{\Hom}(\overline{S(V)},W)}
			{\Hom_{\text{coAlg}}(\overline{S^c(V)},\overline{S^c(W)})}
			{f}
			{\displaystyle
			\sum_{n>0}\sum_{s=1}^n 	
			\pi \circ \left[\mkern-30mu\sum_{\qquad\substack{i_1+\dots+i_s = n \\ 0<i_1\leq \dots \leq i_s}}\mkern-40mu (f_{i_1}\otimes \dots \otimes f_{i_s}) \circ B^{<}_{i_1,\dots,i_s} \right]\circ N~}
		\end{equation}	
		where $B^<_{i_1,\dots,i_n}$ denotes the sum on all the ordered $(i_1,\dots,i_n)$-unshuffles (ordered unshuffleator, see appendix \ref{App:UnshuffleAtors}) and the term between square brackets is an operator $V^{\otimes n}\to V^{\otimes s}$.
	\end{remark}
%
%
%
	The universal property of symmetric tensor coalgebras can be regarded as an isomorphism of graded vector spaces. The following is the analogue of theorem \ref{Theorem:HomCoAlgISO} in the case of symmetric tensor coalgebras.
	\begin{theorem}[Universal property as an isomorphism between hom-spaces (symmetric case)]\label{Theorem:HomCoAlgISOSymmetric}
		The following graded sets are isomorphic:
		\begin{displaymath}
			{\Hom}_{\text{coAlg}}(\overline{S^c(V)},\overline{S^c(W)}) 
			\cong {\Hom}(\overline{S(V)},W) 
			\cong \bigoplus_{n>0} \Hom(V^{\odot n},W)
			~.
		\end{displaymath}
		In particular, the isomorphism is induced from the corestriction according to the following diagram:
		\begin{displaymath}
			\begin{tikzcd}[column sep = huge, row sep = huge]
				\Hom(\overline{S(V)},\overline{S(W)}) \ar[r,"P"] &
				\Hom(\overline{S(V)},W) \ar[d,equal]
				\\
				\Hom_{\text{coAlg}}(\overline{S^c(V)},\overline{S^c(W)}) \ar[u,hook]&
				\Hom(\overline{S(V)},W) \ar[l,"L_{\sym}"] \ar[dl,"L"]
				\\
				\Hom_{\text{coAlg}}(\overline{S^c(V)},\overline{T^c(W)}) \ar[u,"\tiny\substack{\sum\\{n>0}} \frac{1}{n!} \pi \eval_{V^{\otimes n}}"]
			\end{tikzcd}
			~.
		\end{displaymath}
	\end{theorem}
	\begin{proof}
		The commutation of the lower triangle is precisely given by proposition \ref{Prop:UniversalPropertyCocommutativeGradedCoalgebras}. Commutation of the upper square follows from the condition that $\pr \circ \pi = \pr : T(V) \to V$, denoting with a slightly abuse of notation the canonical projection $\pr$ from both $T(V)$ and $S(V)$ to $V$.
	\end{proof}

	\begin{remark}
		Note that functors $T$ and $S$ introduced in remark \ref{Remark:TensorFunctors} can be safely regarded as functors valued in the category of coalgebras since
		\begin{align*}
			f^{\otimes n} \circ \Delta =&~
			\sum_{a=0}^n f^{\otimes n} \circ \mathfrak{d}_{a,n-a} =
			\\
			=&~ 
			\sum_{a=0}^n (f^{\otimes a} \otimes f^{\otimes n -a})\circ\mathfrak{d}_{a,n-a} = 
			\\
			=&~
			\sum_{a=0}^n\mathfrak{d}_{a,n-a}\circ f^{\otimes n} = 	
			\\
			=&~
			\Delta\circ f^{\otimes n}
			~.
		\end{align*}
	\end{remark}

	\begin{notation}\label{not:droppingthesuperc}
		From now on we will drop the superscript $c$ when dealing with tensors (co)algebras.
		It will be clear from the context whether we are focusing on the algebra structure or the coalgebra structure.
	\end{notation}

\ifstandalone
	\bibliographystyle{../../hep} 
	\bibliography{../../mypapers,../../websites,../../biblio-tidy}
\fi

\cleardoublepage


%% file: chapters/multibracketscoderivations/multibracketscoderivations.tex
\chapter{Algebraic structure of multibrackets and coderivations}\label{App:RNAlgebras}
This appendix discusses some algebraic properties possessed by the space of homogeneous multilinear maps between graded vector spaces.
	More precisely, we will deal with defining and providing some basic properties of the Gerstenhaber \cite{Gerstenhaber1963a}\cite{Gerstenhaber1964} and \RN \cite{Nijenhuis1967} products.
These algebraic operations enjoy the crucial property to determine a Lie algebra structure despite being not associative (they are pre-Lie, see appendix \ref{App:PreLie}).
	This framework will be useful in chapter \ref{Chap:Linfinity} for introducing the definition of \emph{$(L_\infty)$-structures}.
	
	This material hereby presented can be considered standard. It has been recorded here to provide a compact reference in the body of the thesis. Several proofs are included for the sake of completeness.
	We point out that some statements specific to the context of graded multilinear algebra (see \eg propositions \ref{Prop:SymmetricGerstenhaberAssociators} and \ref{Prop:DecalageAsCoalgebrasISO}) are not easy to found in the literature.
 	In section \ref{section:AppBConclusions} we synthetically recap all the results in a single commutative diagram.
	Our presentation will be mostly tailored to our subsequent needs; its inspiration can be tracked in several sources 
	\cite{Lecomte1992,Doubek2007,Manetti-website,Delgado2015,Bandiera2016}.

\section{Algebra of multilinear operators (Gerstenhaber algebra)}
In this section, we deal with some algebraic structures canonically associated with the graded vector space of homogeneous multilinear maps; namely, the Gerstenhaber product \cite{Gerstenhaber1964} and the \RN product \cite{Nijenhuis1967}.
	\\
	Let $V$ be a graded vector space. We start by introducing the following notation:
	\begin{notation}\label{Notation:MultilinearMapsSpaces}
		We denote the vector space of $a$-multilinear homogeneous maps in degree $k$ valued in $W$ as
		\begin{displaymath}
			M_{a,k}(V,W)
			:= \left\lbrace
			m:\underbrace{V\times\dots\times V}_{a-\text{times}}\to W 
			\;\left\vert\; 
			\stackanchor[5pt]{\text{$m$  is linear in each argument}}
			{$|m(x_1,\dots,x_a)| = k + \sum_{i=1}^a |x_i|$} 
		\right.\right\rbrace
		~.
		\end{displaymath}		
		The linear structure of $M_{a,k}(V,W)$ is inherited "point-wise" from the vector space structure of $W$ (\cf  with the internal hom-functor in remark \ref{rem:Vectpdvs}).
		 When $W=V$ we will lighten the notation omitting the second entry.		
	\end{notation}

	\begin{remark}[Universal property of multilinear maps]\label{Rem:UnivPropMultilinMaps}
		Recall that the tensor product  $\otimes$ of graded vector spaces is completely characterized by the following universal property:
		\begin{itemize}
			\item given three graded vector space $V,W,Z$, for any homogeneous bilinear map $\varphi : V \times W \to Z$ in degree $k$, there exists an unique homogeneous linear map $\hat{\varphi}: V \otimes W \to Z$, with same degree, such that the following diagram commutes
			\begin{displaymath}
			\begin{tikzcd}
				V \times W \ar[r] \ar[dr,"\varphi"']
				& V \otimes W \ar[d,"\hat{\varphi}"]
				\\
				& Z
			\end{tikzcd}		
			~.			
			\end{displaymath}
		\end{itemize}
		Accordingly, for any graded vector spaces $V,W$, and for any $n\geq 0$ and $k \in \ZZ$, one has that 
	\begin{displaymath}
		M_{n,k}(V,W)\cong \underline{\Hom}^k (V^{\otimes n},W)
		~.	
	\end{displaymath}
	We will usually treat this isomorphism as an identification.
	Hence we will be free to understand elements of $M_{n,k}(V,W)$, \ie $n$-ary functions $V\times\dots \times V \to W[k]$ with the extra property of being separately linear in each entry, as graded homogeneous map from $V^{\otimes n}$ to $W$.
	Homogeneous linear maps from $V^{\otimes n}$ to $W$ will be said \emph{of arity $n$}, and we will often denote the image of a multilinear map $\mu_n$ on $x_1\otimes\dots\otimes x_n$ as $\mu_n(x_1,\dots,x_n)$, separating the input elements by commas and omitting the symbol $\otimes$.
		In particular we have $M_{a,k}(V)\subseteq V[k]\otimes(V^{\otimes a})^\ast$, the equality holds when $dim(V)<\infty$.
	\end{remark}

	\begin{remark}[The (bi-)graded vector space of multilinear maps]\label{rem:BigradedSpaceofMultilinearMaps}
	Considering all the possible arities and degrees collectively, and neglecting the "internal" grading of the graded vector space ${\Hom}(V^{\otimes n},W[k])$ (see remark \ref{rem:neglectinginternalgrading}), one obtains the $(\NN_0\times \ZZ)$-graded (bi-graded) vector space
	\begin{displaymath}
		M(V,W)_{\bullet\bullet} :=
		\left(
			(a,d) \mapsto	\underline{\Hom}^d(V^{\otimes a},W)
		\right)
	\end{displaymath}
	which is respectively graded by the arity and the degree of the homogeneous maps.	
	\\
	In what follows, it will be crucial to identify a correct, in the sense of most suitable to our needs, $\ZZ$-graded (single grading) vector space built out of the above bi-graded vector space.
	The most common choice is to understand the degree of multilinear maps as their degree in the sense of homogeneous maps.
	Namely, we introduce the following graded vector space:
 	\begin{displaymath}
 		\begin{aligned}
 			{M(V,W)} :=  M_{\oplus,\bullet} =
 			\left(
 				k \mapsto \bigoplus_{n\geq 1}M_{n,k}(V,W)
 			\right)
			~. 			
 		\end{aligned}
 	\end{displaymath}
 	We will refer to it as "the" \emph{graded vector spaces of (homogeneous) multilinear maps}.
	According to remark \ref{Rem:UnivPropMultilinMaps} and definition \ref{Def:TensorAlgebra}, one has the following isomorphism of graded vector spaces
	\begin{displaymath}
		{M(V,W)} \cong \underline{\Hom}(\overline{T(V)},W)
	\end{displaymath}
	that we will interpret as an identification without further noticing.
	\\
	In remark \ref{rem:DecasGradedMorph} we will introduce another graded vector space that one can construct out of $M_{\bullet \bullet}(V,W)$ by means of the $\tot$ functor defined in equation \eqref{eq:totalbigradedcomplex}.
	\end{remark}		
	
	Without loss of generality, we now focus to the case where $V=W$. Everything can be easily extended to multilinear maps with any codomain\footnote{When defining the composition, one has to recall that it will be an operation only defined on pairs of multilinear maps with matching domains and codomains.}.
	The graded vector space of multilinear maps can be endowed with a non-associative algebra structure (see \eg \cite{Manetti-website-coalgebras}):
	\begin{definition}[$i$-th Gerstenhaber product]\label{Def:ithGerstenhaberProduct}
		We call \emph{$i$-th Gerstenhaber product} the bilinear map
		\begin{displaymath}
			\morphism{-\gerst_i-}
			{M_{a,k}(V)\times M_{a',k'}(V)}
			{M_{(a+a'-1),(k+k')}(V)}
			{(f,g)}
			{f\gerst_i g = f \circ (\mathbb{1}_{i-1}\otimes g \otimes \mathbb{1}_{a-i})}
		\end{displaymath}
		%
		where $\circ$ denotes the usual composition of graded linear maps and $1\leq i \leq a$, and $\Unit_k$ denotes the identity automorphism on $V^{\otimes k}$.
		\\
		The latter definition can be extended by linearity to the whole $M(V)$ keeping implied that $f\gerst_i g = 0$ when $i$ is greater than the arity of $f$.
	\end{definition}
	\begin{definition}[(Full) Gerstenhaber product]\label{Def:FullGerstenhaberProduct}
		We call \emph{(full) Gerstenhaber product} the bilinear map
		\begin{displaymath}
			\morphism{-\gerst-}
			{M_{a,k}(V)\times M_{a',k'}(V)}
			{M_{(a+a'-1),(k+k')}(V)}
			{(f,g)}
			{\displaystyle f\gerst g = \sum_{i=1}^{a} f \gerst_i g}
			~.
		\end{displaymath}
		As before, this can be extended to the space $M(V)$ of all multilinear maps by linearity.
		%
	\end{definition}
	%
	%
	
	\begin{notation}[Unfolding the Koszul convention]\label{not:GerstProdUnfolded}
		Slightly more explicitly, given $f\in \Hom(V^{\otimes a},V[|f|])$ and $g\in \Hom(V^{\otimes b},V[|g|])$ one has 
		\begin{displaymath}
			(f\gerst_i g)\in	\Hom(V^{\otimes(a+b-1)},V[|f|+|g|])
			~.
		\end{displaymath}
		According to the Koszul convention (see remark \ref{Rem:KoszulTAM}), one gets that
		\begin{displaymath}
			\mathclap{
				\begin{aligned}
					&(f  \gerst_i g) (x_1\dots x_{a+b-1}) =
					\\
					&=~ (-)^{|g|(|x_1|+\dots +|x_{i-1}|)} f(x_1,\dots,x_{i-1},g(x_i,\dots,x_{i+b-1}),x_{i+b},\dots,x_{a+b-1})			
				\end{aligned}
			}
		\end{displaymath}
		and
		
		\begin{displaymath}
			\mathclap{
				\begin{aligned}
				&(f\gerst g) (x_1 \otimes\dots\otimes x_{a+b-1})=
				\\
				&~=
				\mkern-5mu
			\sum_{k=0}^{a-1}(-)^{|g|(|x_1|+\dots+|x_k|)} 
			f(x_1,\dots, x_k , g(x_{k+1} , \dots , x_{k+b}), x_{k+n+1}, \dots , x_{a+b-1})
				\end{aligned}
			}
		\end{displaymath}
		for any given $(a+b-1)$-tuple of homogeneous elements $(x_1\dots x_{a+b-1})$ of $V$.
	\end{notation}	
	The space $M(V)$, taken together with $\gerst$, forms a non-associative monoid with unit given by the identity automorphism $\Unit \in M_{1,0}(L)$.
	More precisely, the Gerstenhaber product yields a pre-Lie structure 
	(see appendix \ref{App:PreLie} for the definition):
	\begin{lemma}
		The graded vector space $M(V)$ together with the Gerstenhaber product $\gerst$ forms a graded right pre-Lie algebra.
	\end{lemma}
	\begin{proof}
		Recall that, according to our definition, the grading of $M(V)$ is given by the degree of the homogeneous map. 
		We thus have to show that the associator $\alpha(f,g,h)= f\gerst(g\gerst h) - (f\gerst g)\gerst h$ is graded symmetric in the two rightmost entries, \ie
		\begin{displaymath}
			\alpha(f,g,h) = (-)^{|h||g|} \alpha(f,h,g)
			~.
		\end{displaymath}
		We show how this computation works in a fairly simple case, a complete and conceptual proof follows from remark \ref{Remark:GerstPreLie}.
		Let us consider for simplicity $f\in \underline{\Hom}(V^{\otimes 2},V)$ and $g,h \in\underline{\Hom}(V,V)$, one has
		\begin{displaymath}
			f\gerst (g\gerst h) =
			\sum_{i=1}^2 f \gerst_i (g\circ h) 
			= f \circ (gh\otimes \mathbb{1} + \mathbb{1}\otimes gh)
		\end{displaymath}
		and
		\begin{displaymath}
			\begin{aligned}
			(f\gerst g) \gerst h =&
			\sum_{j=1}^2 (\sum_{i=1}^2 f \gerst_i g ) \gerst_j h 
			=
			\\
			=&
			f \circ (g\otimes \mathbb{1}+ \mathbb{1}\otimes g) \circ (h\otimes \mathbb{1} + \mathbb{1}\otimes h)
			=
			\\
			=&
			f \circ ( g\otimes h + (-)^{|g||h|} h\otimes g + gh \otimes \mathbb{1} + \mathbb{1}\otimes gh)
			\end{aligned}
		\end{displaymath}
		where the sign coefficient comes from the Koszul convention (see equation \eqref{Eq:TensorHomogeneousMaps}).
		The corresponding associator reads
		\begin{align*}
			\alpha(f,g,h) =&~ f \circ (g\otimes h + (-)^{|h||g|} h \otimes g ) =
			\\ 
			=&~			
			f \circ (g\odot h) ~,
		\end{align*}
		where $\odot$ denotes the symmetric tensor product of linear maps.
		Hence, the latter equation is manifestly graded symmetric in the second and third entries.
	\end{proof}
	Therefore, there is a well-defined Lie bracket:
	\begin{definition}[Gerstenhaber bracket]\label{def:GerstenhaberBracket}
		We call \emph{Gerstenhaber bracket} the Lie bracket $[\cdot,\cdot]$ on the graded vector space of multilinear maps $M(V)$ given by the graded pre-Lie product $\gerst$, \ie
		\begin{displaymath}
		 [f,g] := f\gerst g - (-)^{|f||g|} g \gerst f
		 \qquad \forall f,g \in M_{\bullet}(V)
		 ~.
		\end{displaymath}
	\end{definition}		
	It goes without saying that the Gerstenhaber bracket satisfy the Jacobi equation.		
	\begin{remark}\label{Rem:signsProblemwithGerstenhaberProducts}
		Note that definition \ref{def:GerstenhaberBracket} differs from the original one given by Gerstenhaber \cite{Gerstenhaber1964} where the grading of $M(V)$ is given by the arity minus one:
		\begin{displaymath}
			[f,g]_G := f \gerst g - (-)^{(a-1)(b -1)} g \gerst f
			~,
		\end{displaymath}
		for any $f$ and $g$ as in remark \ref{not:GerstProdUnfolded}.
	\end{remark}

		\subsection{D\'ecalage of multilinear maps}
		It is possible to establish a relationship between multilinear operators on the graded vector space $V$ and those on its shift composing these operators with suitable d\'ecalage isomorphisms (see definition \ref{Def:DecIso}).
		Recall that we introduced the following d\'ecalage isomorphism:
		\begin{displaymath}
			\dec: (W[1])^{\otimes n} \xrightarrow{~\sim~}  (W^{\otimes n}) [n] 
			~,
		\end{displaymath}
		given on arbitrary elements $\omega_1,\dots,\omega_n$ in $W$ by
		\begin{equation}\label{eq:decStandard-appendix}
			\begin{aligned}
			\dec \left( \omega_{1\,[1]}\otimes\dots\otimes \omega_{n\,[1]}	\right) =&~ (-)^{\sum_{k=1}^n(n-k)(|\omega_{k\,[1]}|+1)}  (\omega_1\otimes\dots\otimes \omega_n)_{[n]}
			\\[1em]
			\dec^{-1} \left( (\omega_1\otimes\dots\otimes \omega_n)_{[n]}\right) =&~ (-)^{\sum_{k=1}^n(n-k)|\omega_k|} \omega_{1\,[1]}\otimes\dots\otimes \omega_{n\,[1]}
			~.
			\end{aligned}
		\end{equation}
		The d\'ecalage isomorphism  induces by precomposition a graded linear isomorphism between the spaces of multilinear maps:
	\begin{definition}[D\'ecalage of multilinear maps]\label{Def:MultiBrackDecalage}
		Given two graded vector spaces $V$ and $W$, we call \emph{d\'ecalage} of multilinear maps from $V$ to $W^{\otimes m}$ the invertible homogeneous map:
		\begin{displaymath}
			\Dec :  \lineHom^{k}( V^{\otimes n}, W^{\otimes m}) \to \lineHom^{k +n-m}(V[1]^{\otimes n},W[1]^{\otimes m})
		\end{displaymath}
		given by
		\begin{displaymath}
			\begin{aligned}
			\Dec(\mu) = \dec^{-1}[|\mu|+n-m] \circ \mu[n] \circ \dec 
			&
			\qquad \forall \mu\in \lineHom^{|\mu|}( V^{\otimes n}, W^{\otimes m})
			\\
			\Dec^{-1}(\varphi) = \left(\dec[|\varphi|] \circ \varphi \circ \dec^{-1} \right)[-n] 
			&
			\qquad
			\forall \varphi \in \lineHom^{|\varphi|}( V[1]^{\otimes n}, W[1]^{\otimes m})
			~.
			\end{aligned}
		\end{displaymath}
	\end{definition}
	
	\begin{remark}[D\'ecalage of multilinear maps as a diagram]\label{Rem:DecasDiagram}
			Observe that the d\'ecalage of a given multilinear map is given by the following commutative diagram in the category of graded vector spaces:
		\begin{center}
			\begin{tikzcd}
				(V[1])^{\otimes n} \ar[r,dashed,"\Dec(\mu)"] \ar[d,"\dec"']
				&
				(W[1])^{\otimes m}[|\mu| + n -m ] \ar[d,"\dec{[|\mu|+n-m]}"]
				\\
				V^{\otimes n}{[n]}  \ar[r,"\mu{[n]}"'] 
				& W^{\otimes m} [|\mu| + n] = W^{\otimes m} [m][|\mu| + n-m] 
			\end{tikzcd}
		\end{center}
		 which in particular makes manifest the invertibility of the map $\Dec$. 
		 Namely, the inverse can be read out of the following diagram:
		\begin{center}
			\begin{tikzcd}[column sep = huge]
				(V[1])^{\otimes n} \ar[r,"\varphi"] \ar[d,"\dec"']
				&
				(W[1])^{\otimes m}[|\varphi| ] \ar[d,"\dec{[|\varphi|]}"]
				\\
				V^{\otimes n}{[n]}  \ar[r,dashed,"\Dec^{-1}(\varphi){[-n]}"'] 
				&[2em] W^{\otimes m} [m][|\varphi|] 
			\end{tikzcd}
			~.
		\end{center}
	\end{remark}

	\begin{remark}[D\'ecalage of $M_{\bullet,\bullet}(V)$]\label{Rem:MultiBracketDecalage}
		In the following we will be mainly concerned with the space $M(V)$ given at the beginning of this section.
		In the case $m=1$, the d\'ecalage of multilinear maps $\Dec$ is given by graded morphisms
		\begin{equation}\label{eq:DecComp-appendix}
			\begin{aligned}
				\Dec &: M_{n,k}(V) \xrightarrow{\quad\sim\quad} M_{n,k+n-1}(V[1])
				\\
				\Dec^{-1} &: M_{n,k}(V[1]) \xrightarrow{\quad\sim\quad} M_{n, k +1 -n}(V)
			\end{aligned}~.
		\end{equation}
		The situation is encoded by the following commutative diagram for any $\mu_n\in M_{n,|\mu_n|}(V)$:
		\begin{equation}\label{eq:DecalageDiagram-appendix}
			\begin{tikzcd}
				V^{\otimes n}[n] \ar[r,"\mu_n{[n]}"] & W[|\mu_n|][n] \ar[r,equal] & W[1][|\mu_n|+n-1]
				\\
				(V[1])^{\otimes n} \ar[u,"\dec"] \ar[urr,dashed,"\Dec(\mu_n)"']
			\end{tikzcd}
		\end{equation}
		meaning that $\Dec(\mu_n)= \mu_n[n]\circ \dec$. 
		The d\'ecalage operator of multilinear maps is explicitly given as follows:
			\begin{equation}\label{eq:ExplicitDECA-appendix}
				\morphism{\Dec(\mu_n)}
				{\mkern-50mu V[1]^{\otimes n}}
				{W[1]}
				{v_{1~[1]}\otimes\dots\otimes v_{n~[1]}}
				{(-)^{ \sum\limits_{j=1}^{n}(n-j)|v_j|}\mu_n(v_1\otimes\dots\otimes v_n)_{[n]}}
				~.
			\end{equation}
	\end{remark}

	\begin{remark}[D\'ecalage $\Dec$ as a graded morphism]\label{rem:DecasGradedMorph}
		Observe that the operators given in equation \eqref{eq:DecComp-appendix}, due to their property of intertwining the grading and the arity of any given multilinear maps, 
		cannot be read as the components of a bi-graded morphism on $M_{\bullet \bullet}(V)$.
		\\
		However, it is possible to read the operator $\Dec$ as a genuine graded isomorphism (not bi-graded) by appropriately contracting indices $n$ and $k$ to give a certain $\ZZ$-grading.
		Namely, by its very definition, the linear operator $\Dec$ descends to a graded morphism  
 	\begin{displaymath}
 		\Dec: \overline{\overline{M(V,W)}} \to {M(V[1],W[1])}
 	\end{displaymath}
 	between the $\ZZ$-graded vector spaces constructed out of $M_{\bullet\bullet}$ according to remark \ref{rem:BigradedSpaceofMultilinearMaps} and the graded vector space obtained from the total space of $M_{\bullet \bullet}$.
 	More explicitly, one has:
 	\begin{align*}
 			{M(V,W)} &:= \underline{\Hom}(\overline{T(V)},W) = M_{\oplus,\bullet} =
 			\left(
 				k \mapsto \bigoplus_{n\geq 1}M_{n,k}(V,W)
 			\right)
 			~,
 			\\
 			\overline{\overline{M(V,W)}} &:= 
 			\left(\tot(M_{\bullet \bullet}(V,W))\right)[1]
 			= 
 			\left(
 				k \mapsto \bigoplus_{n+m=k+1}M_{n,k}(V,W)
 			\right)~. 			
 	\end{align*}
	In other words, ${M(V,W)}$ and $\overline{\overline{M(V,W)}}$ denote the $\ZZ$-graded vector space of $n$-multilinear map, for any $n\geq 1$, from $V$ to $W$ taken with two different gradings.
	Namely, given an homogeneous multilinear map $\mu \in M^{n,k}(V,W)$, we denote by 
	\begin{displaymath}
		|\mu|=k ~,\qquad  ||\mu||=k+n-1 ~,
\end{displaymath}	
	the degree of $\mu$ when regarded respectively as an element of ${M(V,W)}$ and $\overline{\overline{M(V,W)}}$. 	
 	Notice that $|\mu|$ coincides, as expected, with the degree of $\mu_k$ as an homogeneous map.
 	\end{remark}	
	
	\begin{remark}[$\Dec$ in the literature]\label{Rem:DecalageinLiterature}
		We ought to notice that our definition of the d\'ecalage of multilinear maps differs by a sign prefactor from the definition often found in literature (e.g \cite[Eq. 3]{LadaStasheff}, \cite[\S 1]{Fiorenza2006},\cite[Rem 1.7]{Bandiera2016},\cite[Prop. 1.5]{Delgado2018b}). 
		The different sign comes from our convention about shift endofunctors to identify $[k][\ell]$, $[\ell][k]$ and $[\ell+k]$ for any $k,\ell \in \mathbb{Z}$ (see remark \ref{rem:compsingShiftsvsSuspensions}).
		Namely, in our convention, the isomorphism $V[k][\ell]\cong V[\ell][k]$ is treated as a natural identification and accordingly denoted with an "$=$" in diagram \eqref{eq:DecalageDiagram-appendix}.
		\\
		However different choices can be made. 
		For instance, if one replaces that identification with natural isomorphism given by the braiding
		\begin{displaymath}
				\begin{tikzcd}[ column sep=4em,row sep=-1ex]
					\mathbb{K}[\ell]\otimes\mathbb{K}[k]\otimes V \ar[r,"B_{\mathbb{K}[\ell],\mathbb{K}[k]}"] 		&
				\mathbb{K}[k]\otimes\mathbb{K}[\ell]\otimes V
				\\
				1_{[\ell]}\otimes 1_{[k]}\otimes v \ar[r,mapsto]&
				(-)^{\ell k}~1_{[k]}\otimes 1_{[\ell]} \otimes v					
			\end{tikzcd}		
			~,
		\end{displaymath}
		the resulting d\'ecalage would get an extra overall sign prefactor.
	\end{remark}

	A crucial property of the d\'ecalage operator $\Dec$ (see definition \ref{Def:DecIso}) is to preserve the Gerstenhaber product modulo a sign. The proof is based on the following lemma:

\begin{lemma}\label{Lemma:CalcoloDecalageGerstenhaber}
	Let $\mu_b\in \underline{Hom}^{|\mu_b|}(V^{\otimes b}, V)$ be a multilinear map  of arity $b$.
	Consider the graded morphism 
	\begin{displaymath}
		(\Unit_a \otimes \mu_b \otimes \Unit_c) \,:~ V^{\otimes (a+b+c)}\to V^{\otimes (a+1+c)}[|\mu_b|]~.
	\end{displaymath}
	Then the d\'ecalage of the latter, computed according to definition \ref{Def:MultiBrackDecalage}, satisfy the following equation
	$$ \Dec(\Unit_a\otimes \mu_b \otimes \Unit_c) = (-)^{S}\Unit_a\otimes \Dec(\mu_b)\otimes \Unit_c~,$$	
	where $\Unit_i$ denotes the identity on $V^{\otimes i}$ on the left-hand side and the identity on $(V[1])^{\otimes i}$ on the right-hand side, and
	$$
	S= |\mu_b|(a+c) + a(b-1)
	~.
	$$
\end{lemma}
\begin{proof}
	According to definition \ref{Def:MultiBrackDecalage}, the left-hand side of the above equation results
	$$
	 \dec_{(a+1+c)}^{-1} \circ \left( \Unit_a\otimes\mu_b\otimes\Unit_c\right) \circ \dec_{{(a+b+c)}}
	$$
	and the right-hand side reads
	$$
		\Unit_a \otimes \left( \dec_{(1)}^{-1}\otimes \mu_b \otimes \dec_{(b)} \right) \otimes \Unit_c ~.
	$$
	In the previous equations we decorated $\dec$ with a subscript to stress the tensor power corresponding to the domain.
	When inspecting on generic elements $x_{i[1]}\in V[1]$, one can check that both sides are proportional to
	$$
	 \left(
	  x_{1[1]},\dots,x_{a[1]},
	  \left(
	  	\mu_a(x_{a+1},\dots,x_{a+b})
	  \right)_{[1]},
	  x_{a+b+1[1]},\dots,x_{a+b+c [1]},	  
	 \right)
	$$
	and only differ by a sign prefactor.
	The left-hand side sign prefactor is $(-)^{s_1 + s_2 +s_3}$ where
	$$
		s_1 = \sum_{k=1}^{a+b+c-k} |x_k|
	$$
	comes from the application of the isomorphism $\dec_{(a+b+c-1)}$ on the given elements,
	$$
		s_2= |\mu_b| \sum_{k=1}^a |x_k|
	$$
	comes from the Koszul convention for evaluating tensor product of maps (see equation \eqref{eq:KoszulConvTensorProducts-appendix}), and
	\begin{align*}
		s_3 
		&=~
		\sum_{k=1}^{a+1+c-k} |y_k| =
		\\
		&=~
		\sum_{k=1}^{a}(a+1+c-k)|x_k| +
		c\left(|\mu_b| + \sum_{k=1}^{b}|x_{a+k}|\right) +
		\sum_{k=1}^{c}(c-k)|x_{a+b+k}|
		\\
		&=~
		\sum_{k=1}^{a}(a+1+c-k)|x_k| +
		c|\mu_b| +
		c  \left( \sum_{k=a+1}^{a+b}|x_{k}|\right) +
		\sum_{k=a+b+1}^{a+b+c}(a+b+c-k)|x_{k}|
	\end{align*}
	comes from evaluating $\dec_{(a+1+c)}^{-1}$ on $(y_1,\dots,y_{a+1+c})$ where
	\begin{displaymath}
		y_k=
		\begin{cases}
			x_k & k\leq a\\
			\mu_b(x_{a+1},\dots,x_{a+b}) & k = a+1\\
			x_{b+k-1} & k\geq a+2
		\end{cases}
		~.
	\end{displaymath}
	The sign prefactor on the right-hand side is given by $(-)^{s_4+s_5}$ where
	$$
		s_4 = | \Dec(\mu_b)| \left(\sum_k=1^a\right) =
		(|\mu_b| + b -1) \left(a + \sum_{k=1}^a |x_k|\right)
	$$
	is the Koszul sign coming from evaluating $\left(\Unit_a\otimes \mu_b \otimes \Unit_c\right)$ on elements, and
	$$
		s_5 = \sum_{k=1}^b (b-k) |x_{a+k}| =
		\sum_{k=a+1}^{a+b} (b+a-k)|x_k|
	$$
	comes from the application of $\dec_{(b)}$ on $(x_{a+1 [1]},\dots, x_{a+b [1]})$.
	Then sought sign prefactor $S$ appearing in the statement results from computing $\sum_{i=1}^5 s_i \mod{2}$.
\end{proof}

\note{Era Sbagliato per un segno nella versione stampata!}

	In particular, when $a=0$ one has
		\begin{displaymath}
			\Dec(\mu_m\otimes\Unit_n) = (-)^{|\mu_m|\cdot n} \Dec(\mu_m)\otimes \Unit_n~.
	\end{displaymath}

	The following corollary is used in subsection \ref{subsec:LinftyMorphi}.
	\begin{corollary}\label{Cor:DecalageofTensorsProducts}
		Consider a $n$-tuple of multilinear maps $f_{k_i} \in\underline{Hom}^{|f_{k_i}|}(V^{\otimes k_i}, V)$, with $1\leq i \leq n$, of arity $k_i$.
		The d\'ecalage of their tensor product is obtained as follows:
			\begin{displaymath}
				\Dec(f_{k_1}\otimes \dots \otimes f_{k_\ell}) =
				(-)^{S}\Dec(f_{{k_1}})\otimes \dots \otimes \Dec(f_{{k_\ell}})
			\end{displaymath}	
			where
			\begin{displaymath}
				S = \sum_i^n
				\left(
				 |f_{k_i}|(n-i) + 
				 \left(\sum_{j=1}^{i-n} k_j\right)
				 \left( |f_{k_i}| + k_i -1 \right)
				\right)
			\end{displaymath}
	\end{corollary}
	\begin{proof}
		Note that
		\begin{align*}
			\Dec\big(f_{k_1} & \otimes \dots \otimes f_{k_n}\big) =
			\\ 
			&=~
			\Dec\left(f_{k_1}\otimes \Unit_{n-1}\right) 
			\circ \dots 	\circ 
			\Dec\left( \Unit_{k_1+\dots+k_{i-1}}\otimes f_{k_i}\otimes \Unit_{n-i}\right)
			\circ \dots \circ
			\Dec\left( \Unit_{k_1+\dots k_{n-1}} \otimes f_{k_n}\right)
			=
			\\
			&=~
			(-)^S
			\Dec(f_{k_1})\otimes \dots \otimes \Dec(f_{k_n})
		\end{align*}
		 where 
		 $$
		 S = \sum_{i=1}^n \Big(
		 	|f_{k_i}|\big(k_1+\dots + k_{i-1} +n -i\big)
		 	+
		 	\big(k_1+\dots +k_{i-1}\big)(k_i-1)
		 \Big)
		 $$ 
		 is the sign coming from multiple applications of lemma \ref{Lemma:CalcoloDecalageGerstenhaber}.

	\end{proof}

	\begin{proposition}[D\'ecalage of Gerstenhaber product]\label{Prop:DecalageGerstenhaberProducts}
		Given any $\mu_n \in M_{n,|\mu_n|}(V)$ and $\ell_n \in M_{n, |\ell_n|}(V[1])$ the following equations holds:
		\begin{displaymath}
			\begin{aligned}
			\Dec(\mu_n \gerst_k \mu_m ) =&
			 (-)^{|\mu_m|(n-k)}
			 \Dec(\mu_n) \gerst_k \Dec(\mu_m)
			\\
			\Dec^{-1}(\ell_n \gerst_k \ell_m)
			=&
			(-)^{|\Dec^{-1}(\ell_m)|(n-k)}
			 \Dec^{-1}(\ell_n) \gerst_k \Dec^{-1}(\ell_m)
			~.
			\end{aligned}		
		\end{displaymath}
	\end{proposition}
	\begin{proof}
		In the first equation we are in the following situation
		\begin{displaymath}
			\begin{tikzcd}[row sep = large]
				(V^{\otimes(n+m-1)})[n+m-1] 
				\ar[d,"(\Unit_{k-1}\otimes \mu_m \otimes \Unit_{n-k})_{[n+m-1]}"']
				&
				(V[1])^{\otimes(n+m-1)} \ar[d,shift right = 1em,"\Dec(\Unit_{k-1}\otimes \mu_m \otimes \Unit_{n-k})"]
				\ar[l,"\dec"]
				\\
				(V^{\otimes n})[n][|\mu_m|][m-1] 
				\ar[d,"{\mu_n [n][|\mu_m|][m-1]}"']
				&
				(V[1])^{\otimes n}[|\mu_m|][m-1]
				\ar[ddl,"\Dec(\mu_n){[|\mu_m| + m -1]}",end anchor=east]
				\ar[l,"\dec{[|\mu_m|][m-1]}"]
				\\
				V[|\mu_n|][n][|\mu_m|][m-1]
				\ar[d,phantom,sloped,"\cong"]
				\\[-2ex]
				(V[1])[|\mu_n| + |\mu_m| + (n+m -1) -1]
			\end{tikzcd}
		\end{displaymath}
		The composition of the two leftmost vertical arrows gives  $(\mu_n \gerst_k \mu_m){[n][m-1]}$ while the composition of the two rightmost gives its d\'ecalage.
		In particular, recalling definition \ref{Def:ithGerstenhaberProduct}, one gets that 
		\begin{displaymath}
			\mathclap{
			\begin{aligned}
			\Dec( \mu_n \gerst_k& \mu_m) 
			= \Dec(\mu_n)[|\mu_m + m -1] \circ \Dec (\Unit_{k-1}\otimes\mu_m \otimes \Unit_{n-k})
			=
			\\
			=& (-)^{|\mu_m|(n-k)} \Dec(\mu_n)[|\mu_m + m -1] \circ (\Unit_{k-1}\otimes \Dec (\mu_m) \otimes \Unit_{n -k})
			=
			\\
			=&
			(-)^{|\mu_m|(n-k)}
			 \Dec(\mu_n) \gerst_k \Dec (\mu_m)
			\end{aligned}
			}
		\end{displaymath}
		where in the last two lines has been employed Lemma \ref{Lemma:CalcoloDecalageGerstenhaber} and definition \ref{Def:ithGerstenhaberProduct}.
		\\
Regarding the second claim, one could draw a diagram similar to the previous one. Alternatively, one can notice that the equality can be derived directly from the invertibility of $\Dec$ together with the following equation:
		\begin{displaymath}
			\Dec \left( \Dec^{-1}(\ell_n) \gerst_k \Dec^{-1}(\ell_m) \right)
			=
						(-)^{|\Dec^{-1}(\ell_m)|(n-k)}
			  (\ell_n \gerst_k \ell_m)
			~.
		\end{displaymath}
	\end{proof}
	\begin{remark}\label{rem:gerstStrano}
		The previous proposition could be read more compactly as follows.
		Introduce the operator
		\begin{displaymath}
			\morphism{-\overline{\overline{\gerst_i}}-}
			{M_{a,k}(V)\times M_{a',k'}(V)}
			{M_{(a+a'-1),(k+k')}(V)}
			{(\mu_a,\mu_{a'})}
			{(-)^{k'(a-i)}\mu_a\gerst_i \mu_{a'}}
		\end{displaymath}	
		that corresponds to $\gerst_i$ with an extra signs introduced in accordance to proposition \ref{Prop:DecalageGerstenhaberProducts}.
		Consider then the operator $\overline{\overline{\gerst}}=\sum_i \overline{\overline{\gerst_i}}$.
		Therefore, proposition \ref{Prop:DecalageGerstenhaberProducts} implies that
		\begin{displaymath}
			\begin{tikzcd}
				Dec:~ (\overline{\overline{M(V)}},\overline{\overline{\gerst}}) \ar[r,"\sim"]& (M(V),\gerst)
			\end{tikzcd}~,
		\end{displaymath}
		where $\overline{\overline{M(V)}}$ is given by remark \ref{rem:DecasGradedMorph} and $M(V)$ is given by remark \ref{rem:BigradedSpaceofMultilinearMaps},
		is an isomorphism in the category of pre-Lie algebras.
	
	\end{remark}

	\subsection{Algebra of multibrackets (Nijenhuis$-$Richardson algebra)}\label{Section:MultibracketsAlgebra}
	We now focus on graded (skew)symmetric multilinear operators. 		
	To distinguish these particular multilinear maps from those without any particular symmetry, we will often refer to them as "multibrackets".
	\begin{notation}
		We denote the spaces of symmetric and skew-symmetric  $n$-multilinear homogeneous maps on $V$ in degree $k$ as $M_{n,k}^{\sym}(V)$ and $M_{n,k}^{\skew}(V)$ respectively. 
		Namely, the two spaces are defined as
		\begin{align*}
				M_{n,k}^{\sym}(V)
				=&
				\left\{
					\ell_n \in M_{n,k} ~|~ \ell_n \circ B_{\sigma} = \ell_n \quad \forall \sigma \in S_n
				\right\}
				~,
				\\
				M_{n,k}^{\skew}(V)
				=&
				\left\{
					\mu_n \in M_{n,k} ~|~ \mu_n \circ B_{\sigma} = (-)^{\sigma}\mu_n \quad \forall \sigma \in S_n
				\right\}
				~.
		\end{align*}		 
	\end{notation}
	\begin{remark}[Universal property of (skew)symmetric multilinear maps]\label{Rem:UnivPropSkewSymMaps}
		In the case of (skew)symmetric maps a universal property, similar to \ref{Rem:UnivPropMultilinMaps}, holds:
		\begin{displaymath}
			\begin{aligned}
				M_{n,k}^{\sym}(V)
				\cong&
				\lineHom^k(V^{\odot n},V)
				\\
				M_{n,k}^{\skew}(V)
				\cong&
				\lineHom^k(V^{\wedge n},V)
				~.
			\end{aligned}
		\end{displaymath}
		It also follows from the general properties of $\Hom$ functors, namely the property to map colimits in limits, that
		\begin{displaymath}
			\begin{aligned}
			M_{n,k} =&~ \lineHom^k(V^{\otimes n}, V) =
			\\
			=&~ \lineHom^k(V^{\odot n}\oplus V^{\wedge n}) \cong
			\lineHom^k(V^{\odot n})\oplus 			\lineHom^k (V^{\wedge n}) =
			\\
			=&~ M_{n,k}^{\sym}\oplus M_{n,k}^{\skew}
			~.
			\end{aligned}
		\end{displaymath}

	\end{remark}
	\begin{remark}\label{Remark:MSymEmbedding}
		Applying the controvariant endofunctor
		\begin{displaymath}
			\underline{Hom}(-,V) = \oplus_{k\in \mathbb{N}} (-,V[k]) ~:~ \GVect \to \GVect
		\end{displaymath}
		to the diagram introduced in remark \ref{Remark:ManettiNotation} we get the following commutative diagram (where $-$ denotes a blank entry, \eg $-\circ f$ means precomposition)
		\begin{displaymath}
			\begin{tikzcd}[column sep = huge]
				M(V)  
				&& M(V) \ar[ll,"\mathcal{-\circ S}"'] \ar[d,two heads,"-\circ N"]
				\\
				M^{\sym}(V) \ar[u,hook,"-\circ\pi"]  \ar[urr,hook,"-\circ\pi"]
				&& M^{\sym}(V) 				\ar[ll,"\sum_{n\geq 0} \frac{1}{n!} \circ - \eval_{V^{\odot n}}"]
			\end{tikzcd}
		\end{displaymath}
		which formalizes how $M^{\sym(V)}$ can be embedded in $M(V)$. As before, we omit the analogous result in the skew-symmetric case.	
	\end{remark}
	
	By its very definition, the subspace of (skew)symmetric multibrackets does not form a subalgebra of $M(V)$ (it does not close with respect to $\gerst$). This can be cured by considering a suitable (skew)symmetrization of $\gerst$.

	\subsubsection{The \RN product on $M(V)$}\label{subsec:RNMV}
	A preliminary step to understand how one can equip the spaces of symmetric and skew-symmetric multilinear maps with a meaningful algebra product (or "composition") is to investigate how the Gerstenhaber product $\gerst$ behaves under symmetrization.
	Observe first the following:
	\begin{lemma}\label{Lemma:GerstenhaberSymmetrization}
		Given any two $\ell_n,\ell_m \in M(V)$ one has
		\begin{displaymath}
			\left(\ell_n \gerst \ell_m\right) \circ \mathcal{S} =
			\left(\frac{n }{\#\ush{m,n-1}} \right) ~ 
			(\ell_n  \circ \mathcal{S}) \circ ((\ell_m \circ \mathcal{S}) \otimes \Unit_{n-1}) \circ B_{m,n-1}
		~,
		\end{displaymath}
		where $\mathcal{S}$ denotes the symmetrizator (see definition \ref{Def:Symmetrizator}) and $B_{m,n-1}$ is the even permutator with respect to the subgroup of unshuffles (see definition \ref{Def:Unshuffleator}).
		Similarly, one has
		\begin{displaymath}
			\left(\mu_n \gerst \mu_m\right) \circ \mathcal{A} = 
			\left(\frac{n}{\#\ush{m,n-1}} \right)
			~ \left(\mu_n \circ \mathcal{A}\right)
			 \circ ((\mu_m\circ\mathcal{A}) \otimes \Unit_{n-1}) \circ P_{m,n-1}
		~,
		\end{displaymath}
		where $\mathcal{A}$ denotes the skew-symmetrizator (see definition \ref{Def:SkewSymmetrizator}) and 
		$P_{m,n-1}$ is the odd permutator with respect to the subgroup of unshuffles (see definition \ref{Def:Unshuffles}).
	\end{lemma}
	\begin{proof}
		Consider any two multilinear maps $\ell_n,\ell_m \in M^{\sym}(V)$. 
		Via a simple inspection, one can see that 
		\begin{displaymath}
		\begin{aligned}
			\ell_n \gerst_i \ell_m 
			=&~
			\ell_n \circ (\Unit_{i-1}\otimes \ell_m \otimes \Unit_{n-i}) =
			\\
			=&~
			\ell_n \circ (\cycPermutator_{(i)}\otimes \Unit_{n-i}) \circ (\ell_m\otimes \Unit_{n-1}) \circ (\cycPermutator_{(i+m)}^i\otimes \Unit_{n-1})
			=
			\\
			=&~
			\left(
			(\ell_n \circ (\cycPermutator_{(i)}\otimes \Unit_{n-i})) \gerst_1 \ell_m \right)  
			\circ (\cycPermutator_{(i+m)}^i\otimes \Unit_{n-1})
		\end{aligned}
		\end{displaymath}
		where $\cycPermutator_{(k)}^i$ denotes the cyclic permutation of $k$ elements $i$-times (see remark \ref{Notation:Cyclic permutator} in the following appendix).
		Therefore, one has
		
		\begin{displaymath}
		\mathclap{
		\begin{aligned}
			&\left(\ell_n \gerst \ell_m\right) \circ \mathcal{S}
			~=
			\\
			&=
			\sum_{i=1}^n \ell_n \gerst_i \ell_m \circ \mathcal{S}
			\\
			&=
			\sum_{i=1}^n \left(\ell_n \circ (\cycPermutator_{(i)}\otimes \Unit_{n-i})) \gerst_1 \ell_m \right)  \circ \mathcal{S}
			=
			\\
			&=
				\left(\frac{m! (n-1)!}{(n+m-1)!}\right)
				\left( \ell_n \circ 
					\left(\sum_{i=1}^n
						(\cycPermutator_{(i)}\otimes \Unit_{n-i})\circ ( \Unit \otimes \mathcal{S}_{(n-1)})
					 \right) 
				 \right)
				 \gerst_1 
				 (\ell_m \circ \mathcal{S}_{(m)})
				 \circ B_{m,n-1}
			=
			\\		
			&=
				\frac{1}{\#\ush{m,n-1}} \frac{n!}{(n-1)!}
				(\ell_n \circ \mathcal{S}_{(n)})
								 \gerst_1 
				 (\ell_m \circ \mathcal{S}_{(m)})
				 \circ B_{m,n-1}
			=
			\\
			&=
				\frac{n}{\#\ush{m,n-1}}
				(\ell_n \circ \mathcal{S}_{(n)})
								 \gerst_1 
				 (\ell_m \circ \mathcal{S}_{(m)})
				 \circ B_{m,n-1}
			\end{aligned}
			}
		\end{displaymath}
		where we used 
		the decomposition of the symmetrizator operator through unshuffles (see remark \ref{Rem:UnshufflesAsCoset}) in the third equality and
		the computation of the cardinality of the subgroup of unshuffles (see  remark \ref{Rem:UnshufflesCardinality}) and the decomposition of the symmetrizator operator through unshuffles (see remark \ref{Rem:UnshufflesAsCoset}) in the fourth one.
		\\
		The same computation holds in the case of skew-symmetric multibrackets by substituting $B_\sigma$ with $P_\sigma$ and $\mathcal{S}$ with $\mathcal{A}$.
	\end{proof}

	\note{MZ\\
	OK. Only now I understand that "Gerstenhaber project" is for maps without any symmetry, while "NR" is the name for the product on (skew)symm maps. Is this made clear in the introduction?}
	Lemma \ref{Lemma:GerstenhaberSymmetrization} suggests the following two definitions:
	\note{metti sidebyside?}
		\begin{definition}[(Symmetric) \RN product]\label{Def:SymGerstProd}
			We call \emph{(symmetric) \RN product} (or \emph{symmetric Gerstenhaber product}) the bilinear operator
			\begin{displaymath}
				- \symgerst - : M_{n, k}\otimes M_{m, d}(V) \to M_{n+m -1, k +d}(V)
			\end{displaymath}
			defined on any $\mu_n,\mu_m \in M(V)$ as
			\begin{displaymath}
				\mu_n \symgerst \mu_m = (\mu_n \gerst_1 \mu_m) \circ B_{m, n-1}
				~,
			\end{displaymath}
		\end{definition}
		where $B_{n,m-1}$ denotes the operator that sum on all possible even actions of the $(n,m-1)$-unshuffles (see appendix \ref{App:UnshuffleAtors}) and $\gerst_1$ is the first Gerstenhaber product (see definition \ref{Def:ithGerstenhaberProduct}).
	%
		\begin{definition}[(Skew-symmetric) \RN product]\label{Def:SkewGerstProd}
			We call \emph{(skew-symmetric) \RN product} the bilinear operator
			\begin{displaymath}
				- \skewgerst - : M_{n, k}\otimes M_{m, d}(V) \to M_{n+m -1, k +d}(V)
			\end{displaymath}
			defined on any $\mu_n,\mu_m \in M(V)$ as
			\begin{displaymath}
				\mu_n \skewgerst \mu_m =(-)^{|\mu_m|(n-1)}(\mu_n \gerst_1 \mu_m) \circ P_{m, n-1}
				~,
			\end{displaymath}
			where $P_{n,m-1}$ denotes the operator that sum on all possible odd actions of the $(n,m-1)$-unshuffles (see appendix \ref{App:UnshuffleAtors}) and $\gerst_1$ is the first Gerstenhaber product (see definition \ref{Def:ithGerstenhaberProduct}).
		\end{definition}	
		The role of the extra sign in the previous definition is justified a posteriori by the following lemma:
	\begin{proposition}\label{Prop:DecalageAsCoalgebrasISO}
		The graded map $\Dec$ defined in remark \ref{Notation:FixingDec} is a graded algebra isomorphism:
		\begin{displaymath}
			\Dec: (\overline{\overline{M(V)}}, \skewgerst) \xrightarrow{\quad \sim \quad} (M(V[1]),\symgerst)
			~.
		\end{displaymath}
	\end{proposition}
	\begin{proof}
		Inspecting on any one two given $\mu_n, \mu_m \in M(V)$ one gets
		\begin{displaymath}
			\begin{aligned}
			\Dec( \mu_n \skewgerst \mu_m ) 
			\equal{Def: \ref{Def:SkewGerstProd}}&
			(-)^{|\mu_m|(n-1)}\Dec(\mu_n \gerst_1 \mu_m \circ P_{m,n-1})
			=
			\\
			\equal{Def: \ref{Def:MultiBrackDecalage}}&
			(-)^{|\mu_m|(n-1)}\Dec(\mu_n \gerst_1 \mu_m) \circ  \Dec(P_{m,n-1})
			=
			\\
			\equal{Prop: \ref{Prop:DecalageGerstenhaberProducts}}&
			\Dec(\mu_n) \gerst_1 \Dec(\mu_m) \circ  \Dec(P_{m,n-1})
			=
			\\
			\equal{Prop: \ref{Prop:DecalageOfPermutation}}&
			\Dec(\mu_n) \gerst_1 \Dec(\mu_m) \circ  B_{m,n-1}
			=
			\\
			\equal{Def: \ref{Def:SymGerstProd}}&
			\Dec(\mu_n) \symgerst \Dec(\mu_m)
			~.
			\end{aligned}
		\end{displaymath}
	\end{proof}

	\begin{notation}
 	More explicitly, evaluating on homogeneous element $x_i\in V$, the products read as follows:
 	
	\begin{equation}\label{Eq:RNProducts-explicit-appendix}
		\mathclap{
		\begin{aligned}
		 \mu_n \cs \mu_m &(x_1,\dots,x_{m+k-1}) =
		 \\
		 =&~
		 \mkern-30mu	
		 \sum_{\sigma \in unsh(m,n-1)}
		 \mkern-30mu		 
		  \epsilon(\sigma) 
		 \mu_n\Big(\mu_m(x_{\sigma_1},\dots,x_{\sigma_m}),x_{\sigma_{m+1}}\dots,x_{\sigma_{m+k-1}}	\Big)
		 \\[2em]
		 \mu_n \ca \mu_m &(x_1,\dots,x_{m+k-1}) =
		 \\
		 =&~ 
		 (-)^{|\mu_m|(n-1)}\mkern-30mu
		 \sum_{\sigma \in unsh(m,n-1)}\mkern-30mu	
		  (-)^\sigma \epsilon(\sigma) 
		 \mu_n\Big(\mu_m(x_{\sigma_1},\dots,x_{\sigma_m}),x_{\sigma_{m+1}}\dots,x_{\sigma_{m+k-1}}\Big)
		\end{aligned}
		}
	\end{equation}
 where $\epsilon(\sigma) $ is the Koszul sign.
	\end{notation}
	
	The bigraded vector space $M(V)$, taken together with $\symgerst$ or $\skewgerst$, forms a non-associative graded algebra. The following proposition gives an explicit expression for the associator:
	\begin{proposition}\label{Prop:SymmetricGerstenhaberAssociators}
		Given any three multilinear operators $\mu_\ell,\mu_m,\mu_n \in M(V)$ the corresponding associators, with respect to $\symgerst$ and $\skewgerst$ respectively, results
		\begin{align*}
			\alpha(\symgerst;\mu_a,\mu_b,\mu_c) 
			=&~
			\mu_a \circ \left((\mu_b\otimes\mu_c \circ B_{b,c})\otimes \Unit_{a-2}\right) \circ B_{b+c,a-2}
			\\
			\alpha(\skewgerst; \mu_a,\mu_b,\mu_c) 
			=&~
			(-)^{\mathcal{s}}
			\mu_a \circ \left((\mu_b\otimes\mu_c \circ P_{b,c})\otimes \Unit_{a-2}\right) \circ P_{b+c,a-2}		
			~.
		\end{align*}
		where the latter sign prefactor is given by:
		\begin{displaymath}
			{\mathcal{s}}={{|\mu_c|(b+a)} +{|\mu_b|(a-1)} +b(c+1)}
			~.
		\end{displaymath}
	\end{proposition}
	\begin{proof}
		A complete proof of these statements can be obtained employing the combinatorial properties of the unshuffles. We refer to appendix \ref{Section:AppendixProofPreLie} for a complete argument.
		We only observe here that the second claim follows from the first since, from the equality $\Dec(\mu_a\ca\mu_b) = \Dec(\mu_a)\cs \Dec(\mu_b)$, one has
		$$
			\alpha\Big( 
				\cs; \Dec(\mu_a),\Dec(\mu_b),\Dec(\mu_c)
			\Big)
			=
			\Dec\Big(\alpha\big(
				\ca; \mu_a, \mu_b, \mu_c
			\big)\Big)~.
		$$ 
		Namely, the latter equality reads as
		\begin{align*}
			\Dec(\mu_a) \circ &
			\left(
				\Dec(\mu_b)\otimes\Dec(\mu_c)\otimes\Unit_{a-2}
			\right)
			\circ
			B_{b,c,a-2}	
			=
			\\
			=&~
			(-)^{\mathcal{s}}
			\Dec(\mu_a) \circ
			\Dec\left(
				\mu_a \otimes \mu_b \otimes \Unit_{a-2}
			\right)\circ
			B_{b,c,a-2}
			=
			\\
			=&~
			(-)^{\mathcal{s}+S}
			\Dec(\mu_a) \circ
			\left(
				\Dec(\mu_a) \otimes \Dec(\mu_b) \otimes \Unit_{a-2}
			\right)\circ
			B_{b,c,a-2}
		\end{align*}	
		where $\mathcal{s}$ is the sign prefactor appearing in the statement and 
		$$ 
			S = |\mu_c| (b+a) + |\mu_b|(a-1) + b(c-1)
		$$
		comes from applying lemma \ref{Lemma:CalcoloDecalageGerstenhaber} two times.			In particular it is implied that $\mathcal{s} = S \mod{2}$.
	\end{proof}

	\subsubsection{The \RN product on $M^{\sym}(V)$ and $M^{\skew}(V)$}
	Let us now focus again on the spaces of degree $k$ graded (skew)-symmetric $n$-multilinear homogeneous maps, previously denoted as 
	\begin{displaymath}
		M_{n,k}^{\sym}(V,W):= \underline{\Hom}^k (V^{\odot n},W)
		~,\quad
		M_{n,k}^{\skew}(V,W):= \underline{\Hom}^k (V^{\wedge n},W)
		~.
	\end{displaymath}
	Recall that	It follows from the splitting sequence \eqref{eq:symskewsplitting} that $M_{n,k}(V,W)=M_{n,k}^{\sym}(V,W)\oplus M_{n,k}^{\skew}(V,W)$.
	 As before, when $W=V$ we will lighten the notation omitting the second entry.
	 Recall that the operator $\dec$ determines an isomorphism $V^{\wedge n}[n]\cong V[1]^{\odot n}$.
	 This last property is also suitably transferred to multibrackets.
	 \\ 
	 	\begin{remark}[D\'ecalage of symmetric multibrackets]\label{Rem:DecalageofMultiBrackets}
	 Note that a relation similar to the one on the restriction of the d\'ecalage isomorphism to (skew)symmetric tensor spaces (see lemma \ref{Lemma:DecalageRestrictToSymTens}) also holds when dealing with multilinear maps. Namely, the following diagram commutes
	 \begin{displaymath}
		\begin{tikzcd}[column sep = large]
			M^{\sym}_{n, k}(V) \ar[r,"\Dec\big|_{M^{\sym}_{n, k}(V)}"] \ar[d,hook]
			&[3em]
			M^{\skew}_{n, k+n-1}(V[1])\ar[d,hook]
			\\
			M_{n, k}(V) \ar[r,"\Dec"] 
			&
			M_{n, k+n-1}(V[1])
			\\
			M^{\skew}_{n, k}(V) \ar[r,"\Dec\big|_{M^{\skew}_{n, k}(V)}"'] \ar[u,hook]&
			M^{\sym}_{n, k+n-1}(V[1]) \ar[u,hook]
		\end{tikzcd}
	 \end{displaymath}
		where the hooked arrows denote the standard inclusion of symmetric and skew-symmetric maps as a subset of all multilinear maps.
	\end{remark}
	
	\begin{notation}\label{Notation:FixingDec}
	 To avoid confusion,  from now on we will assume that
	 \begin{displaymath}
	 	\Dec: M^{\skew}_{n, k}(V) \to M^{\sym}_{n, k+n-1}(V[1])
	 \end{displaymath}
	 and
	 \begin{displaymath}
	 	\Dec^{-1}: M^{\sym}_{n, k}(V[1]) \to	M^{\skew}_{n, k+1-n}(V)
	 	~.
	 \end{displaymath}
	\end{notation}
 	
 	We are now interested in replicating the reasoning of remark \ref{rem:DecasGradedMorph} to the case of (skew)-symmetric multibrackets.
	Namely, we want to consider suitable combinations of arities and degrees in order to produce a honest $\ZZ$-graded vector space of (skew)-symmetric  multibrackets.
 	\\ 
 	Conventionally, we introduce the \emph{graded vector spaces of graded symmetric and graded skew-symmetric multilinear maps} from $V$ to $V$ as vector subspaces of ${{M(V,W)}}$ and $\overline{\overline{M(V,W)}}$ respectively
 	\begin{definition}[The $\ZZ$-graded vector space of (skew)-symmetric multilinear brackets]\label{Def:VSpacesRN}
 		In order to read $\Dec$ as a graded morphism between skew-symmetric and symmetric multilinear maps, it is necessary to equip their corresponding spaces with two different gradings.
 		Let  $V$ and $W$ be two graded vector spaces.
 		We call \emph{($\ZZ$)-graded vector space of graded symmetric multilinear maps} the graded vector subspace $M^{\sym}(V,W)\hookrightarrow M(V,W)$ given by:
 		\begin{displaymath}
 			M^{sym}(V,W):= \left( k \mapsto \bigoplus_{n} M^{\sym}_{n, k}(V,W)\right)
 			~.
 		\end{displaymath}
 		We call \emph{($\ZZ$)-graded vector space of graded skew-symmetric multilinear maps} the graded vector subspace $M^{\skew}(V,W)\hookrightarrow \overline{\overline{M(V,W)}}$ given by:
 		\begin{displaymath}
 			M^{skew}(V,W):= \left( k \mapsto \bigoplus_{n+i=k+1} M^{\skew}_{n, i}(V,W)\right)
 		~.
 		\end{displaymath}
 		Consider a multilinear maps $\mu\in M_{a,k}$ which is also graded skew-symmetric, we denote $|\mu|=k$ its degree as an homogeneous map and we denote by $\parallel\mu\parallel = a+k-1$ its degree as an element in $\overline{\overline{M(V,W)}}$.
 	\end{definition}
 	This unusual choice for the grading of $M^{\skew}(V,W)$ is due to enforce the commutativity of the following diagram in the category of graded vector spaces
 	\begin{equation}\label{eq:diagramMultibracketsSpaces}
 		\begin{tikzcd}
 			M^{skew}(V) \ar[r,"\Dec"] \ar[d,hook]&
 			M^{sym}(V[1])\ar[d,hook]
 			\\
 			\overline{\overline{M(V)}} \ar[r,"\Dec"] & M(V[1])
 		\end{tikzcd}
 		~.
 	\end{equation}
In other words, this convention assures that the d\'ecalage operator of multilinear maps restricts to a well-defined isomorphism $M^{sym}(V[1])\cong M^{\skew}(V)$.
	The preliminary work of subsection \ref{subsec:RNMV} allows us to conclude that the above diagram can be replicated inside the category of graded algebras:
	\begin{theorem}\label{Thm:ManettiFactorizationOnGerstenhaberAlgebras}
		Let $V$ be a graded vector spaces. Consider the graded algebras of multilinear maps (see definition \ref{Def:VSpacesRN} and remark \ref{rem:DecasGradedMorph}). The following diagram commutes in the category of graded  vector spaces:
		\begin{displaymath}
			\begin{tikzcd}
				(\overline{\overline{M(V)}},\overline{\overline{\gerst}})
				\ar[r,"\Dec"',"\sim"]\ar[d,two heads,"-\circ N_a"] & 
				(M(V),\gerst)\ar[d,two heads,"-\circ N_s"]
				\\
				(M^{\skew}(V),\ca)\ar[d,hook,"-\circ \pi_a"]
				\ar[r,"\Dec"',"\sim"] & 
				(M^{\sym}(V),\cs)\ar[d,hook,"-\circ \pi_s"]
				\\
				(\overline{\overline{M(V)}},\ca)
				\ar[r,"\Dec"',"\sim"] & 
				(M(V),\cs)								
			\end{tikzcd}
		\end{displaymath}	
		where:
		\begin{itemize}
			\item $\gerst$ is the Gerstenhaber product (definition \ref{Def:FullGerstenhaberProduct});
			\item $\overline{\overline{\gerst}}$ is its pullback along $\Dec$ (see remark \ref{rem:gerstStrano});
			\item $\ca$ and $\cs$ are the \RN products (see definitions \ref{Def:SymGerstProd} and \ref{Def:SkewGerstProd});
			\item $-\circ F$ denotes precomposition with $F$ and $N_a,N_s,\pi_s,\pi_a$ are the maps introduced in remark \ref{Remark:ManettiNotation} and equation \eqref{eq:symskewsplitting}.
	\end{itemize}				
	\end{theorem}
	\begin{proof}
		Observe that, according to Lemma \ref{Lemma:GerstenhaberSymmetrization}, the monoid structures $\symgerst$ and $\skewgerst$ restrict correctly to a monoid structure on $M^{\sym}(V)$ and $M^{\skew}(V)$ respectively.
		Hence, $\symgerst$ (resp. $\skewgerst$) takes values in $M^{\sym}(V)$ (resp. $M^{\skew}(V)$) when computed on a pair of symmetric (resp. skew-symmetric) multilinear maps.
		\\
		Moreover, the symmetric (resp. skew-symmetric) Gerstenhaber product between a generic multilinear map $\mu_n \in M_{n,k}(V)$ with $n\leq 2$ and a symmetric (resp. skew-symmetric) multibracket $\mu_m$ yields a symmetric multibracket $\mu_n \symgerst \mu_m \in M^{\sym}(V)$ (resp. $\mu_n \skewgerst \mu_m \in M^{\skew}(V)$).
		Resuming the notation given in remark \ref{Remark:ManettiNotation}, lemma \ref{Lemma:GerstenhaberSymmetrization} implies that 
		\begin{displaymath}
			\ell_n \gerst \ell_m \circ \frac{N \circ \pi}{\cancel{(n+m-1)!}}
			=
			\cancel{\left(\frac{n! m!}{(n+m-1)!}\right)}
			\left(\ell_n \circ \frac{N \circ \pi}{\cancel{n!}}\right)
			\symgerst
			\left(\ell_m \circ \frac{N \circ \pi}{\cancel{m!}} \right)
			~,
		\end{displaymath}		
		hence, the following diagram commutes
		\begin{displaymath}
			\begin{tikzcd}[column sep = huge]
				(M(V),\symgerst) & 
				(M(V),\gerst) \ar[l,"-\circ N \circ \pi"'] \ar[d,two heads,"-\circ N"]
				\\
				& (M^{\sym}(V),\symgerst) \ar[ul,hook,"-\circ \pi"]
			\end{tikzcd}
			~.
		\end{displaymath}
		The same reasoning can be replicated in the skew-symmetric case replacing $N$ with $N_a$ and $\pi$ with $\pi_a$.
		Summing up, the whole situation is subsumed by the claimed diagram.
	\end{proof}

	The pairs $(M^{\sym}(V),\cs)$ and $(M^{\skew}(V))$ are also known as the \RN algebras of symmetric and skew-symmetric multilinear maps.
	These algebras are right pre-Lie:

	\begin{proposition}\label{prop:RNExplictPreLie}
		Consider three graded symmetric multilinear operators $\nu_i$ and three skew-symmetric multilinear operators $\mu_i$.
		The corresponding associators satisfy the following symmetry properties:
		\begin{align}
				\alpha(\symgerst; \nu_\ell,\nu_m,\nu_n) 
				=&
				(-)^{|\nu_m||\nu_n|}
				\alpha(\symgerst; \nu_\ell,\nu_n,\nu_m)
				~,
				\label{Equation:SymGerstAssociatorSymmetry-ripetizione}			
				\\
				\alpha(\skewgerst; \mu_\ell,\mu_m,\mu_n) 
				=&
				(-)^{(|\mu_n| + n - 1)(|\mu_m| + m - 1)}
				\alpha(\skewgerst; \mu_\ell,\mu_n,\mu_m)
				\label{Equation:SkewGerstAssociatorSymmetry-ripetizione}
				~.
		\end{align}
	\end{proposition}		
	\begin{proof}
		A complete proof of these statements can be obtained employing the combinatorial properties of the unshuffles. 
		We refer to proposition \ref{Prop:SymmetricGerstenhaberAssociators-AppendixC} in appendix \ref{Section:AppendixProofPreLie} for a complete argument.
		We only observe here that the second claim follows from the first since, for any $\nu_i= \Dec(\mu_i)$ one has
		\begin{displaymath}
			\begin{aligned}
				\alpha(\ca;\mu_a,\mu_b,\mu_c)
				=&~ \Dec^{-1}(\alpha(\cs;\nu_a,\nu_b,\nu_c) ) =
				\\
				=&~
				(-)^{|\nu_b||\nu_c|}\Dec^{-1}(\alpha(\cs;\nu_a,\nu_c,\nu_b) )
				\\
				=&~
				(-)^{|\nu_b||\nu_c|}\alpha(\ca;\mu_a,\mu_c,\mu_b)
				\\
				=&~
				(-)^{\parallel\mu_b\parallel\parallel\mu_c\parallel}\alpha(\ca;\mu_a,\mu_c,\mu_b)~.
			\end{aligned}
		\end{displaymath}	
	\end{proof}

	\begin{remark}[Pre-Lie structure on $M^{\sym}(V)$ and $M^{\skew}(V)$]
	The last two results of proposition \ref{Prop:SymmetricGerstenhaberAssociators} can be read as a pre-Lie property of $\symgerst$ and $\skewgerst$ with respect to a certain choice of (single) grading on the bigraded vector space $M(V)$.
	Let us denote this new single grading as $\lVert\cdot \rVert$ in order to not confuse it with $|\cdot|$ that gives the degree as an homogeneous map.
	\\
	Equation \eqref{Equation:SymGerstAssociatorSymmetry-ripetizione} implies that on the graded vector space
	\begin{displaymath}
		M(V) := \big( k \mapsto \oplus_{a \geq 1} M_{a,k}(V) \big)
		~,
	\end{displaymath}
	\ie when considering the grading given uniquely by the degree as an homogeneous map (meaing that $\lVert\cdot \rVert = |\cdot|$), the \RN product gives a pre-Lie structure.
	In particular, $(M^{\sym}(V),\symgerst)$, implicitly considered with the former grading, forms a pre-Lie algebra and
	\begin{displaymath}
		\morphism{{[\cdot,\cdot]}}
		{M^{\sym}(V)\otimes M^{\sym}(V)}
		{M^{\sym}(V)}
		{\mu_m\otimes \mu_n}
		{\mu_m\symgerst \mu_n - (-)^{|\mu_m||\mu_n|}\mu_n\symgerst\mu_m}
	\end{displaymath}
	satisfies the \emph{graded Jacobi equation}.
	
	Similarly, Equation \eqref{Equation:SkewGerstAssociatorSymmetry-ripetizione} implies a pre-Lie property for $\skewgerst$ when considering a grading interweaving the arity of the multilinear map with the degree, namely
	\begin{equation}\label{Skew-symmetric grading}
		\lVert \mu_n \rVert = |\mu_n| + n -1 \qquad \forall \mu_n \in M_{n,|\mu_n|}(V)
		~.
	\end{equation}
	In other words, $\skewgerst$ is pre-Lie on the graded vector space
	\begin{displaymath}
		\tot(M_{\bullet \bullet}(V))[1] :=
		\left( k \mapsto \bigoplus_{a+d=k +1 } M_{a,d}(V) = \bigoplus_{a\geq 1 } M_{a,k+1-a}(V)   \right)
	\end{displaymath}	
	where $\tot$ denotes the total graded vector space of a bigraded vector space as defined in
	\ref{Section:CategoricalGradedMonoidalStructure}.
	\\
	In particular, $(M^{\skew}(V),\skewgerst)$, implicitly considered with the former grading, forms a pre-Lie algebra and
	\begin{displaymath}
		\morphism{{[\cdot,\cdot]}}
		{M^{\skew}(V)\otimes M^{\sym}(V)}
		{M^{\skew}(V)}
		{\mu_m\otimes \mu_n}
		{\mu_m\skewgerst \mu_n - (-)^{\lVert\mu_m\rVert  \lVert\mu_n\rVert}\mu_n\skewgerst\mu_m}
	\end{displaymath}
	satisfies the Jacobi equation.
	\\
	\note{ MZ: 	I guess another way to say this is: replace V by V[1], to absorb the arity in the homogeneous degree.
	
		Questa è di sicuro una ripetizione, non riesco a correggere prima della consegna.	
	}
	\end{remark}

\section{Algebra of coderivations}
Let us now return again to the notions of graded algebras and coalgebras introduced in section \ref{section:GradedAlgebras}.
Another relevant class of graded linear maps between graded (co)-algebras is given by (co)-derivations.
In this section, we recall the abstract framework of $F$-coderivations and in particular prove how they give rise to a certain algebra isomorphic to the algebras of multilinear maps introduced in the previous section.

\subsection{ $F$-(co)derivations on an arbitrary (co)algebra}\label{sec:abtractFcoderivations}
In this subsection, we recall the abstract framework of $F$-coderivations. 
The property that will be central to us is the existence of a lift procedure from homogeneous maps to coderivations (see subsection \ref{subsec:UniPropeLifts}).

Consider a graded algebra $(A,\mu)$ and a graded coalgebra $(C,\Delta)$, one can define: 

	\begin{paracol}{2}
		\begin{leftcolumn}
			\begin{adjustwidth}{-2em}{0em}	
	%
	\begin{definition}[$F$-Derivation]\label{Def:Fderivation}
		Given two graded algebras $(A,\mu)$ and $(A',\mu')$, consider a morphism of graded algebras $F\in\Hom_{\text{Alg}}(A,A')$, we call \emph{$F$-derivation from $A$ to $A'$} any homogeneous map $\theta \in \underline{\Hom}(A,A')$  such that the following diagram commutes
		\begin{displaymath}
			\begin{tikzcd}[ampersand replacement=\&]
				A\otimes A \ar[rr,"(\theta\otimes F + F\otimes\theta)"] \ar[d,"\mu"']
				\& \&
				A'\otimes A' \ar[d,"\mu'"]
				\\
				A \ar[rr,"\theta"] \& \& A'
			\end{tikzcd}
		\end{displaymath}
	\end{definition}
	\begin{remark}		
		From the Koszul convention (Remark \ref{Rem:TensorHomogeneousMaps}) 	
		one has for any homogeneous vectors $a,b\in A$ that
		\begin{displaymath}
			\theta(a\cdot b) = \theta(a)\cdot F(b) + (-)^{|\theta||a|}F(a)\cdot \theta(b)
			~,
		\end{displaymath}
		where we denoted by $\cdot$ the action of the multiplication $\mu$.		
	\end{remark}
	
			\end{adjustwidth}
		\end{leftcolumn}
		\begin{rightcolumn}
			\begin{adjustwidth}{0em}{-2em}	

	%
	\begin{definition}[$F$-coDerivation]\label{Def:FcoDerivations}
		Given two graded coalgebras $(C,\Gamma)$ and $(C',\Gamma')$, consider a morphism of graded coalgebras $F\in\Hom_{\text{coAlg}}(C,C')$, 
		 we call \emph{$F$-coderivation from $C$ to $C'$} any homogeneous map $Q \in \underline{\Hom}(C,C')$  such that the following diagram commutes
		\begin{displaymath}
			\begin{tikzcd}[ampersand replacement=\&]
				C \ar[rr,"Q"] \ar[d,"\Gamma"']
				\& \&
				C' \ar[d,"\Gamma'"]
				\\
				C\otimes C \ar[rr,"(Q \otimes F + F\otimes Q)"] \& \& C'\otimes C'
			\end{tikzcd}
			~.
		\end{displaymath}
	\end{definition}
			\end{adjustwidth}
		\end{rightcolumn}
	\end{paracol}
\sidebyside{
	\begin{notation}
		Fixed a graded coalgebra morphism, the set of all $F$-derivations in degree $p$ forms a graded vector space denoted by $\Der^p(A,A';F)$,.
		Accordingly, we denote by $\Der(A,A';F)$ the space of coderivations in any degree.
		\\
		In the special case where $A=A'$ and $F= \id_A$ we simply talk about \emph{derivations over $A$}  and we denote as
		\begin{displaymath}
			 \Der^p(A):= \Der^p(A,A;\id_A)
		\end{displaymath}				
		the corresponding graded vector space.		
	\end{notation}

}{
	\begin{notation}
		Fixed a graded coalgebra morphism, the set of all $F$-coderivations in degree $p$ form a graded vector space denoted by $\coDer^p(C,C';F)$.
		Accordingly, we denote by $\coDer(C,C';F)$ the space of coderivations in any degree.
		\\
		In the special case where $C=C'$ and $F= \id_C$ we simply talk about \emph{coderivations over $C$}  and we denote as
		\begin{displaymath}
			 \coDer^p(C):= \coDer^p(C,C;\id_C)
		\end{displaymath}				
		the corresponding graded vector space.		
	\end{notation}
	\begin{lemma}
		Given $F:(C,\Gamma)\to(C',\Gamma')$ a coalgebra morphism and $Q,Q'$ coderivations on $C$ and $C'$ respectively, then $F\circ Q$ and $Q'\circ F$ are $F$-coderivations from $C$ to $C'$.
	\end{lemma}
}

\begin{remark}
	In case that the considered (co)algebra is (co)unital, we also require a (co)derivation to be  zero at the (co)unit.
\end{remark}	
The upshot is that the notion of graded coderivation is dual to the notion of graded derivation. 
For our purposes, we will focus only on coderivations. 
Similar results can also be enunciated in the case of derivations and possibly further generalized from vector spaces to modules (\cite{Manetti-website-coalgebras,Doubek2007}).

\begin{remark}[Composition of coderivations is not a coderivation]\label{Rem:CompositionofCoderivation}
	Although $\coDer(C)$ constitutes a subspace of $\underline{Hom}(C,C)$, it does not form a subalgebra with respect to the composition of linear operator.
	Namely, if one consider two coderivations $M$ and $N$, the following diagram commutes
	\begin{displaymath}
		\begin{tikzcd}[column sep=huge]
			C \ar[rr,"\Gamma"]\ar[d,"M"'] 
			&& C\otimes C \ar[d,"M\otimes\mathbb{1} + \mathbb{1}\otimes M"]
			\\
			C \ar[rr,"\Gamma"]\ar[d,"N"'] 
			&& C\otimes C  \ar[d,"N\otimes\mathbb{1} + \mathbb{1}\otimes N"]
			\\
			C \ar[rr,"\Gamma"] 
			&& C\otimes C					
		\end{tikzcd}
		~.
	\end{displaymath}
	Nevertheless, the two rightmost vertical arrows composed to 
	$$(N\circ M) \otimes \mathbb{1} + \mathbb{1}\otimes (N\circ M) + (-)^{|M||N|} M\otimes N + N\otimes M$$
	which in general differs to $(N\circ M) \otimes \mathbb{1} + \mathbb{1}\otimes (N\circ M)$.
	On the other hand, the graded vector spaces of coderivations is a Lie subalgebra with respect to the Lie bracket given by the commutator:
	\begin{align*}
		\Gamma&\circ [N,M]_{\circ} =
		\\
		=&~
		\Gamma\circ (N\circ M - (-)^{|N||M|}M\circ N)
		\\
		=&~
			\big(
			N\circ M \otimes \mathbb{1} + \mathbb{1}\otimes N\circ M +{ (-)^{|M||N|} M\otimes N + N\otimes M}
			+
		\\
		&\phantom{\big(}-(-)^{|M||N|}(M\circ N \otimes \mathbb{1} + \mathbb{1}\otimes M\circ N + {M\otimes N + (-)^{|M||N|}N\otimes M)\big)\circ \Gamma} =
		\\
		=&~
		\big(
		N\circ M \otimes \mathbb{1} + \mathbb{1}\otimes N\circ M -(-)^{|M||N|}(M\circ N \otimes \mathbb{1} + \mathbb{1}\otimes M\circ N \big)\circ \Gamma	
		\\			
		=&~
		\left([N,M]_\circ \otimes \Unit + \Unit \otimes 	[N,M]_\circ \right)\circ \Gamma
		~.
	\end{align*}
	The latter property suggests to look for a pre-Lie structure on $\coDer(V)$. 
	This will be achieved in the case of tensor coalgebras (see sections \ref{Subsection:CoderivationsAlgebra} and \ref{SubSection:CoderivationonSymmetricTensor}).
\end{remark}

An interesting feature of degree zero coderivations is that they can be used to build coalgebra morphisms:
		\begin{proposition}
			Given $\alpha\in \coDer^0(C)$ a degree $0$ nilpotent coderivation , \ie $\exists k \in \mathbb{N}$ such that $\alpha^{k'} = 0 ~\forall k' > k$, 
			then
			\begin{displaymath}
				e^\alpha := \sum_{n\geq 0} \frac{\alpha^n}{n!}:C\to C
			\end{displaymath}
			is a coalgebra automorphism with inverse given by $e^{-\alpha}$.
		\end{proposition}
		\begin{proof}
			Note first that the nilpotency condition guarantees that the series defining the exponential of $\alpha$ converge to a linear operator $C\to C$. Requiring also $\alpha$ to be in degree $0$ guarantees that the exponential will be homogeneous in degree $0$.
			Compatibility with the comultiplication $\Gamma$ follows from the \emph{Baker-Campbell-Hausdorff formula}\footnote{Precisely, we use the BCH formula for the two commutating elements $\alpha\otimes\Unit$ and $\Unit\otimes\alpha$ in the Lie algebra of linear endomorphisms on $C\otimes C$ with the Lie bracket given by the usual commutator.}
			(BCH):
			\begin{displaymath}
				\begin{aligned}
				\Gamma e^\alpha 
				\equal{}& 
				\sum_{n\geq 0 } \frac{1}{n!} \Gamma \alpha^n 
				=  \sum_{n\geq 0 }\frac{1}{n!}(\alpha\otimes \Unit + \Unit \otimes \alpha ) \Gamma 
				= (e^{\alpha\otimes \Unit + \Unit \otimes \alpha})~\Gamma =
				\\
				\equal{BCH}&
				(e^{\alpha\otimes \Unit})~\circ (e^{\Unit \otimes \alpha})~\Gamma =
				(e^{\alpha}\otimes \Unit)~\circ (\Unit \otimes e^{\alpha})~ \Gamma =
				(e^{\alpha}\otimes e^{\alpha})~ \Gamma	
				~.
				\end{aligned}
			\end{displaymath}
			Invertibility is again a consequence of the BCH formula.
		\end{proof}
		\begin{remark}
			Taking $C= T(V)$ one see that the first component of the exponential of a given coderivation $\alpha \in \coDer(T(V))$ is given by
			\begin{displaymath}
				\pr \circ e^\alpha \eval_V = 
				\sum_{n\geq 0} \frac{1}{n!} \pr \circ \alpha^n \eval_V = 
				\sum_{n \geq 0} \frac{1}{n!} \alpha_1^n = 
				e^{\alpha_1}
				~,
			\end{displaymath}
			where $\pr:T(V)\twoheadrightarrow V$ denotes the standard projection.
			Claerly, not every coalgebra automorphism can be build as an exponential of  some coderivation.
		\end{remark}
		
		\note{Falsa Pista:\\ non è vero che ogni polinomio $P(\alpha) = \sum_i c_i \alpha^i$ costruito a partire da una coderivazione costituisce un isomorfismo. La formula BCH vale solo per l'esponenziale.}
		
		\note{
			Given $\alpha \in \coDer^0(\overline{T(V)},\overline{T(V)})$ we know from the universal property that there exists a morphism $a \in \underline{\Hom}(\overline{T(V)},V)$ such that $\widetilde{L}(a)=\alpha$ hence one can construct two coalgebra morphism $L(a)$ and $e^{\widetilde{L}(a)}$ in $\Hom_{\text{coAlg}}(\overline{T(V)},\overline{T(V)})$.
			
			TODO: che relazione c'è fra i due? coincidono? esiste $\varphi\in \Iso_{\text{coAlg}}(T(V),T(V))$ such that $L(a) = \varphi \circ e^{\widetilde{L}(a)} \circ \varphi^-1$ ?		
		}

\subsubsection{Universal Properties and Lifts}\label{subsec:UniPropeLifts}
In section \ref{Section:TensorCoalgebras} we have seen how a graded morphism $f\in {\Hom}(C,V)$ can be uniquely lifted to a coalgebra morphism $F \in {\Hom_{\text{coAlg}}(C,\overline{T(V)}}$. 
In the case of a homogeneous map  $f\in \underline{\Hom}(C,V)$, not necessarily in degree $0$, an analogous lift construction holds. 
In the latter case, one talk about \emph{lift to a coderivation}.
\\
From now on, we will assume that all considered coalgebras will be coassociative and conilpotent.

\begin{proposition}[Universal property of homogeneous maps]\label{Prop:UniversalPropertyCoderivation}
	Consider a coalgebra $(C,\Gamma)$ a graded morphism $f:C \to V$ and denote by $F=L(f) \in \Hom_{\text{coAlg}}(C,\overline{T(V)})$ the unique lift of $f$ to a coalgebra morphism. 
	The following property are satisfied:
	\begin{enumerate}
		\item for any homogeneous map $q \in \underline{\Hom}^k(C,V)$ in degree $k$ and graded morphism $f\in\Hom(C,V)$, there exists a unique $F$-coderivation 
		$ Q \in \coDer(C,\overline{T(V)};F)$
	such that $F$ equals the lift to a coalgebra morphism $L(f)$ and the following diagram commutes in the graded vector spaces category:
	\begin{displaymath}
		\begin{tikzcd}
			C \ar[r,"\exists!~Q"] \ar[dr,"\forall~q"']
			&
			\overline{T(V)}[k] \ar[d,"{\pr[k]}"]
			\\
			& 
			V[k]
		\end{tikzcd}
	\end{displaymath}
		
		\item Explicitly, $Q$, the unique lift of $q$ to a coderivation, is given by the following commutative diagram in $\GVect$:
			\begin{displaymath}
				\begin{tikzcd}
					\overline{T(C)}\ar[rrd,"R(q)"]
					\\
					C \ar[u,"\sum_{n>0}\Gamma^n"] \ar[rr,"Q"'] & & \overline{T(V)	}	[k]				
				\end{tikzcd}
			\end{displaymath}
			where 
			\begin{equation}\label{Eq:operatorR}
				\morphism{R}
				{\underline{\Hom}(C,V)}
				{\coDer\big(\overline{T(C)},\overline{T(V)};T(f)\big)}
				{q}
				{\displaystyle R(q)=\sum_{n\geq 1} \left[\sum_{i=0}^{n-1} f^{\otimes i} \otimes q \otimes f^{\otimes(n-1-i)}\right]}
			\end{equation}
			and hence 
			\begin{equation}\label{Equation:explicitLiftToCoderivation}
				Q =  \sum_{n \geq 1} \left[
					\sum_{i=0}^{n-1} (f^{\otimes i} \otimes q \otimes f^{\otimes (n-1-i)}) \circ \Gamma^{n-1}
				\right] ~.
			\end{equation}
	\end{enumerate}
\end{proposition}
\begin{proof}
	One has only to check that			
			\begin{displaymath}
				R(q) = \sum_{n \geq 1} \left[
					\sum_{i=0}^{n-1} f^{\otimes i} \otimes q \otimes f^{\otimes (n-1-i)}
				\right]
			\end{displaymath}
			is a $T(f)$-coderivation in the sense of definition \ref{Def:FcoDerivations} and remark \ref{Remark:TensorFunctors} (see \cite[Prop. 1.3.5]{Manetti-website-coalgebras} for the complete computation).
			If the latter is true, the precomposition of $R(q)$ with the canonical coalgebra morphism $\sum_{n\geq 1}\Gamma^n \in \Hom_{\text{coAlg}}(C,\overline{T(C)})$
			happen to be a $L(f)$-coderivation, where $L$ denotes the lift to a coalgebra morphism (see equation \eqref{Eq:LiftMorphismMapping}) because 
			\begin{displaymath}
				(T(f))\circ \sum_{n\geq 1}\Gamma^n = L(f)
				~.
			\end{displaymath}
\end{proof}

\begin{notation}[Lift to a coderivation operator]\label{Notation:CoderLiftMapping}
	We denote by
	\begin{displaymath}
		\morphism{\widetilde{L}_{(f)}}
		{\underline{\Hom}^k(C,V)}
		{\coDer^k(C,\overline{T(V)};F)}
		{q}
		{Q}
	\end{displaymath}
	the mapping giving the unique \emph{lift to a $F$-coderivation} (with $F=L(f)$) prescribed by proposition \ref{Prop:UniversalPropertyCoderivation}.
	We decorated this operator with a tilde to avoid possible confusion when dealing with degree zero map $C\to V$ since, according to remark \ref{Rem:LiftToMorphism}, they also admit a \emph{lift to a coalgebra morphism}.
	We we usually omit the subscript $(f)$ when it is clear from the context.
\end{notation}

The similarity with the case of lifts to coalgebra morphisms, see section \ref{Section:TensorCoalgebras}, continue by noticing that this universal property can be read as an isomorphism of graded vector spaces (\cf  with theorem \ref{Theorem:HomCoAlgISO}):
	\begin{theorem}[Universal property as an isomorphism between hom-spaces]\label{Theorem:HomCoAlgISO-Coder}
		Let $(C,\Gamma)$ be a graded coalgebra.
		The following graded vector spaces are isomorphic
		\begin{displaymath}
			\coDer^k(C,\overline{T(V)};F) \cong \underline{\Hom^k}(C,V)
			~,
		\end{displaymath}
		where $F=L(f)$ is the unique lift to a coalgebra morphism of an arbitrary $f\in \Hom(C,V)$ as defined in section \ref{Section:TensorCoalgebras}.
		The above isomorphism is induced from the corestriction according to the following diagram:
		\begin{displaymath}
			\begin{tikzcd}
				\underline{\Hom}^k(C,\overline{T(V)}) \ar[r,"P"] &
				\underline{\Hom}^k(C,V) \ar[d,equal]
				\\
				\coDer^k(C,\overline{T(V)};F) \ar[u,hook]&
				\underline{\Hom}^k(C,V) \ar[l,"\widetilde{L}"]
			\end{tikzcd}~,
		\end{displaymath}
		where $P$ denotes post-composition with the standard projection $\pr:T(V) \to V$.
		In particular, one also have that
		\begin{displaymath}
			\coDer^k(C,\overline{T(V)};F) \cong \coDer^k(C,\overline{T(V)};G) 
		\end{displaymath}
		for any choice of $F,G\in \Hom_{\text{coAlg}}(C,\overline{T(V)})$.

	\end{theorem}
	\begin{proof}
		For any $q \in \underline{\Hom}^k(C,V)$ one has that
		\begin{displaymath}
			P(\widetilde{L}(q)) = \pr		\sum_{n \geq 1} \left[
					\sum_{i=0}^{n-1} (f^{\otimes i} \otimes q \otimes f^{\otimes (n-1-i)}) \circ \Gamma^{n-1}
				\right] = q
		\end{displaymath}
		since the only non-zero summand is the one with $n=1$ and $i=0$.
		Hence $P\circ L = \id_{\underline{\Hom}^k(C,V)}$.
		\\
		For any $Q \in \coDer^k(C,\overline{T(V)})$ one has
		\begin{align*}
		 \widetilde{L} (P (Q)) \equal{\quad}&
		 \widetilde{L}(\pr\circ Q) = 
		 \\
		 \equal{}& \sum_{n\geq 1} \left( \sum_{i=0}^{n-1} f^{\otimes i}\otimes (\pr\circ Q) \otimes f^{n-1-i}\right) \circ \Gamma^{n-1} =
		 \\
		 \equal{Rem. \ref{Rem:LiftToMorphism}} &
		 \sum_{n\geq 1} \left( \sum_{i=0}^{n-1} (\pr\circ F)^{\otimes i}\otimes (\pr\circ Q) \otimes (\pr \circ F)^{n-1-i}\right) \circ \Gamma^{n-1} =
		 \\
		 \equal{\quad}&
		 \sum_{n\geq 1} \pr^{\otimes n}\circ\left( \sum_{i=0}^{n-1} F^{\otimes i}\otimes Q \otimes F^{n-1-i}\right) \circ \Gamma^{n-1} =			 
		 \\
		  \equal{ $F$-coderivation}&\quad
		   \sum_{n\geq 1} \pr^{\otimes n} \circ \Delta^{n-1} \circ Q=
		  \\
		  \equal{}&		 
		   \sum_{n\geq 1} \pr\eval_{V^{\otimes n}} \circ Q =
		   \\
		   \equal{}& Q
		\end{align*}
		where $\pr\eval_{V^{\otimes n}}: T(V) \to V^{\otimes n}$ is the standard projection.
	\end{proof}

	In the particular case where the coalgebra $C$ is cocommutative, \eg when $C=\overline{S(V)}$, a similar universal property holds.
	The difference is that now the lift maps into the symmetric tensor coalgebra.
\begin{proposition}[Universal property of homogeneous maps (cocommutative case)]\label{Prop:UniversalPropertyCoderivationCocommutative}
	Consider a cocommutative coalgebra $(C,\Gamma)$ a graded morphism $f:C \to V$ and denote by $F=L_{\sym}(f) \in \Hom_{\text{coAlg}}(C,\overline{S(V)})$ the unique lift of $f$ to a coalgebra morphism. 
	The following properties are satisfied:
	\begin{enumerate}
		\item for any homogeneous map $q \in \underline{\Hom}^k(C,V)$ in degree $k$ and $ f\in\Hom(C,V)$, there exists a unique $F$-coderivation 
		$ Q \in \coDer(C,\overline{S(V)};F)$,
	with $F=L_{\sym}(f)$ is the lift to a cocommutative coalgebra morphism, 
	such that the following diagram commutes in the graded vector space category:
	\begin{displaymath}
		\begin{tikzcd}
			C \ar[r,"\exists!~Q"] \ar[dr,"\forall~q"']
			&
			\overline{S(V)}[k] \ar[d,"{\pr[k]}"]
			\\
			& 
			V[k]
		\end{tikzcd}
	\end{displaymath}
		
		\item Explicitly, the unique lift is given by the following commutative diagram in $\GVect$:
			\begin{displaymath}
				\begin{tikzcd}[column sep = huge]
					\overline{T(C)}\ar[rr,"R(q)"] & 
					& \overline{T(V)}[k]\ar[dd,"{\pi[k]}"]
					\\
					\overline{S(C)} \ar[u,hook,"N"]
					\\
					C \ar[u,"\sum\limits_{n>0}\frac{\pi}{n!}\Gamma^n"'] \ar[uu,bend left=60, "\small\sum\limits_{n>0}\Gamma^n"] \ar[rr,"Q"'] \ar[rruu,"\widetilde{L}(q)",sloped]& 
					& \overline{S(V)	}[k]						
				\end{tikzcd}
			\end{displaymath}
			where $R$ is the operator defined in equation \eqref{Eq:operatorR} and $N$ is the injection of $S(V)\hookrightarrow T(V)$ introduced in remark \ref{Remark:ManettiNotation}.
	\end{enumerate}
\end{proposition}
\begin{proof}
	The proof runs in two steps. 
	The first is to prove that when $C$ is cocommutative the image of $\sum_{n>0} \Gamma^n$ is $S_n$-invariant, \ie the operator factors trough $N$. 
	The second is to prove that the image of $R(q)$ restricted to $\overline{S(V)}$, \ie $\pi \circ \widetilde{L}(q)$, is a coderivation with respect to $\pi \circ T(f) \circ N = S(f)$ (see for instance \cite[\S 1.5]{Manetti-website-coalgebras}).
\end{proof}
\begin{notation}\label{Notation:CoderLiftMappingCocommutative}
	We denote by
	\begin{displaymath}
		\morphism{\widetilde{L}_{\sym}= \pi \circ \widetilde{L}}
		{\underline{\Hom}^k(C,V)}
		{\coDer^k(C,\overline{S(V)};F)}
		{q}
		{Q}
	\end{displaymath}
	the mapping giving the unique lift to a coderivation prescribed by proposition \ref{Prop:UniversalPropertyCoderivationCocommutative}.
\end{notation}

In the cocommutative case, it is possible to give an a alternative explicit expression of the lift employing the sum on all unshuffles.
	\begin{proposition}[Unshuffles notation for the symmetric lift]\label{Prop:SymLiftUnshufflesNotation}
		Given $(C,\Gamma)$ a cocommutative coalgebra with coproduct $\Gamma$, for any $q \in \underline{\Hom}(C,V)$ and $f\in\Hom(C,V)$, the unique lift is given by
		$$\widetilde{L}_{\sym}(q) = H(q) \sum_{n>0}\Gamma^n$$ where
		\begin{displaymath}
			\morphism{H}
			{\underline{\Hom}(C,V)}
			{\coDer\big(\overline{S(C)},\overline{S(V)};S(f)\big)}
			{q}
			{\displaystyle H(q)=\sum_{n\geq 1} \pi \circ
			\left( q \otimes f^{\otimes n-1} \circ B_{1,n-1} \right) 
			\circ N}
		\end{displaymath}
	\end{proposition}
	\begin{proof}
		The claim follows from the commutativity of the following diagram
		\begin{displaymath}
			\begin{tikzcd}[column sep = huge]
				\overline{S(C)} \ar[hook,d,"N"] \ar[drr,"H(q)"]&&
				\\
				\overline{T(C)}\ar[drr,"R(q)"] && \overline{S(V)}\ar[d,hook,"N"]
				\\
				C \ar[rr,"\widetilde{L}(q)"] \ar[uu,bend left=60,"\sum_{n>0}\frac{\pi}{n!}\Gamma^n"] \ar[u,"\sum_{n>0}\Gamma^n"']\ar[drr,"q"]
				&& \overline{T(V)} \ar[d]
				\\
				&& V
			\end{tikzcd}
		\end{displaymath}
		(compare with diagram in theorem \ref{Prop:UniversalPropertyCocommutativeGradedCoalgebras} )
		where the leftmost triangle is the factorization of $\sum_{n>0}\Gamma^n$ when $C$ is cocomutative, and the uppermost square is consequence of the following computations (involving the operator $\mathcal{C}_{(n)}$ of cyclic permutation, see appendix \ref{App:UnshuffleAtors}):
		\begin{align*}
			R_n(b) &\mathcal{S}_{(n+1)}
			=
			\\
			\equal{}&
			\left(\sum_{i+j=n} f^{\otimes i}\otimes b \otimes f^{\otimes j} \right) \circ \mathcal{S}_{(n+1)}
			=
			\\				
			\equal{}&
			\frac{1}{(n+1)!} \sum_{i+j=n} \sum_{\sigma \in S_{(n+1)}}\left(\sum_{i+j=n} f^{\otimes i}\otimes b \otimes f^{\otimes j} \right)\circ B_{\sigma}
			=
			\\
			\equal{Lem: \ref{Lemma:CyclicCommutation}}&
			\frac{1}{(n+1)!} \sum_{i+j=n} \sum_{\sigma \in S_{(n+1)}}
			\left((\mathcal{C}_{i+1}\otimes \Unit_j)
			\circ(b \otimes f^{\otimes n})\circ \mathcal{C}^{-1}_{(i+1)} \right)\circ B_{\sigma}
			=
			\\
			\equal{}&
			\frac{1}{(n+1)!} \sum_{i+j=n} \sum_{\sigma \in S_{(n+1)}}
			\left((\mathcal{C}_{i+1}\otimes \Unit_j)
			\circ(b \otimes f^{\otimes n}) \right)\circ B_{\sigma}
			=
			\\
			\equal{Eq: \ref{Eq:DecompositionofSymmetrizator}}&
			\frac{n!}{(n+1)!} \sum_{i+j=n}
			\left((\mathcal{C}_{i+1}\otimes \Unit_j)
			\circ(b \otimes f^{\otimes n}\mathcal{S}_{(n)}) \right)\circ B_{1,n}				
			=
			\\
			\equal{Prop: \ref{Proposition:PermutingHomogeneousMaps}}&
			\left[
			\frac{1}{n+1} \sum_{i=0}^n
			(\mathcal{C}_{i+1}\otimes \Unit_{n-i}) \circ (\Unit\otimes \mathcal{S}_{(n)})
			\right]
			\circ(b \otimes f^{\otimes n}) \circ B_{1,n}				
			=
			\\
			\equal{}& \mathcal{S}_{(n+1)} \circ (b \otimes f^{\otimes n}) \circ B_{1,n}
			~.
		\end{align*}
	The latter implies that
		\begin{align*}
				R(b) \circ N 
				=&~ 
				R(b) \circ \mathcal{S} \circ N =
				\\
				=&~
				N \circ \left[ \pi \circ \left(\sum_{n>0} b \otimes f^{\otimes n} \circ B_{1,n} \right) \circ N \right] =
				\\
				=&~ N \circ H(b)~.
		\end{align*}
	\end{proof}

	Similarly to theorem \ref{Theorem:HomCoAlgISO-Coder}, this universal property can be also read as an isomorphism of graded vector spaces:
	\begin{theorem}[Universal property as an isomorphism between hom-spaces]\label{Theorem:HomCoAlgISOCoCommutative-Coder}
		Given a garded cocommutative coalgebra $C$, the following graded vector spaces are isomorphic
		\begin{displaymath}
			\coDer^k(C,\overline{S(V)};F) \cong \underline{\Hom^k}(C,V)
		\end{displaymath}
		where $F= L_{\sym}(f)$ is the unique lift to a coalgebra morphism of an arbitrary $f\in \Hom(C,V)$ as defined in section \ref{Section:TensorCoalgebras}.
		The isomorphism is induced from the corestriction according to the following diagram:
		\begin{displaymath}
			\begin{tikzcd}
				\underline{\Hom}^k(C,\overline{S(V)}) \ar[r,"N\circ-"] &[-0.5em]
				\underline{\Hom}^k(C,\overline{T(V)}) \ar[r,"P"] &[-0.5em]
				\underline{\Hom}^k(C,V) \ar[d,equal]
				\\
				\coDer^k(C,\overline{S(V)};F) \ar[u,hook] \ar[r,leftrightarrow,"\cong"]&
				\coDer^k(C,\overline{T(V)};N \circ F) \ar[u,hook]&
				\underline{\Hom}^k(C,V) \ar[l,"\widetilde{L}"]
			\end{tikzcd}~,
		\end{displaymath}
		where $P$ denotes post-composition with the standard projection $T(V) \to V$ and $N\circ-$ denotes precomposition with the injection $N:S(V)\hookrightarrow T(V)$.
	\end{theorem}
	\begin{proof}
		We already showed that $\widetilde{L}$ factorizes through $\underline{\Hom}^k(C,\overline{S(V)})$ since the image of $\widetilde{L}(q)$ is invariant under permutations for any $q \in \underline{\Hom}^k(C,V)$.
		It remains to exhibit the bottom left isomorphism.
		Note first that for any $Q\in \coDer^k(C,\overline{S(V)};F)$ one has
		\begin{displaymath}
			\begin{aligned}
			\Delta \circ N \circ Q =&~ N\otimes N \circ \Xi \circ Q =
			\\
			=&~ N \otimes N \circ ( Q \otimes F + F \otimes Q ) \circ \Gamma
			\end{aligned}
		\end{displaymath}
		where $\Delta$, $\Xi$ and $\Gamma$ denote the coproduct in $T(V)$, $S(V)$ and $C$ respectively.
		Therefore, post-composition with $N$ yields a morphism 
		$$N\circ - : \coDer^k(C,\overline{S(V)};F) \to \coDer^k(C,\overline{T(V)};N \circ F) ~.$$
		This morphism is invertible, with inverse given by $\pi$ since any $B\in \coDer^k(C,\overline{T(V)};N \circ F)$ lies in the image of $\widetilde{L}$ which is contained in the space of permutation invariants.		
	\end{proof}
	In other words, any coderivation $\theta\in\coDer^k(S(V))$ is fully determined by its corestriction \ie by the projection onto $V$.

\subsection{Algebra of coderivations on tensor spaces}\label{Subsection:CoderivationsAlgebra}
	We now focus our attention to the space of coderivations on a given (reduced) tensor coalgebra.
	Namely, we are going to specialize the construction of the previous subsection to the case of $F$-coderivations on $C=\overline{T(V)}$ with
		\begin{displaymath}
			F = L(\pr) = \id_{\overline{T(V)}}	~,		
		\end{displaymath}	
	where $\pr: \overline{T(V)}\to V$ denotes the standard projection.
	Recall that we denote the graded vector spaces of such coderivations as
		\begin{displaymath}
		\mathclap{
		\coDer(\overline{T(V)}) :=
		\left \lbrace
			Q \in \underline{\Hom}(\overline{T(V)},\overline{T(V)}) ~ \left\vert \quad
			\Gamma\circ Q = (Q\otimes \Unit + \Unit \otimes Q) \circ  \Delta
		\right.			
		\right \rbrace
		~,}
	\end{displaymath}
	where $\Gamma$ denotes the \emph{deconcatenation coproduct} acting on decomposable elements as follows:
	\begin{displaymath}
		\morphism{\Gamma}
		{\overline{T(V)}}{\overline{T(V)}\otimes \overline{T(V)}}
		{x_1\otimes\dots \otimes x_n}
		{\displaystyle\sum_{i=1}^{n-1}(x_1\otimes\dots\otimes x_i)\otimes (x_{i+1}\otimes\dots\otimes x_n)}
		~.
	\end{displaymath}	
	
\note{
		\begin{displaymath}
		\mathclap{
		\coDer(\overline{T(V)},\overline{T(W)};F) :=
		\left \lbrace
			Q \in \underline{\Hom}(\overline{T(V)},\overline{T(W)}) ~ \left\vert \quad
			\Delta\circ Q = (Q\otimes F + F \otimes Q) \circ  \Delta
		\right.			
		\right \rbrace
		~.}
	\end{displaymath}
}	
The interest for this case, stems from the fact that $\coDer(\overline{T(V)})$ can be proved to be isomorphic to the spaces of the multilinear operators  and in particular they provide an alternative language for dealing with multibrackets:

\begin{corollary}[(of prop. \ref{Prop:UniversalPropertyCoderivation})]
	Let  $C = \overline{T(V)}$ be the reduced tensor coalgebra of a graded vector space $V$. Consider the graded morphism $f=\pr: \overline{T(V)}\to V$, the unique lift prescribed by proposition \ref{Prop:UniversalPropertyCoderivation} is given by:
	\begin{displaymath}
		\morphism{\widetilde{L}}
		{\underline{\Hom}^k(\overline{T(V)},V)}
		{\coDer^k(\overline{T(V)})}
		{q}
		{Q=\displaystyle \sum_{n\geq1} \sum_{s=1}^{n} \left( \sum_{i=0}^{s-1} \Unit_i\otimes q_{n-s+1}\otimes\Unit_{s-1-i} \right)}
	\end{displaymath}
	Where $q_i = \pr \circ q \eval_{V^{\otimes i}}$ denote, as usual, the $i$-the corestriction of $q$.
	On any given n-tuple of vectors $v_1,\dots,v_n$, one has explicitly that
	\begin{displaymath}
		Q(v_1\otimes \dots \otimes v_n)= \sum_{i,\ell} (-)^{k(|v_1|+\dots+|v_i|)} v_1\otimes \dots \otimes v_i \otimes q(v_{i+1}\otimes\dots\otimes v_{i+\ell})\otimes \dots \otimes v_n ~.
	\end{displaymath}
\end{corollary}

For any graded vector space the naturally associated graded space of coderivations $\coDer(\overline{T(V)})$ (recall that the grading is given by the degree as homogeneous map) can be equipped with a (non-associative) algebra structure:
\begin{definition}[Gerstenhaber product of coderivations]\label{Def:GerstProdOfCoderivations}
 We call \emph{Gerstenhaber product} of coderivations the (non-associative) product given by
 \begin{displaymath}
 	\morphism{-\gerst-}
 	{\coDer(\overline{T(V)})\otimes \coDer(\overline{T(V)}) }
 	{\coDer(\overline{T(V)})}
 	{Q\otimes R}
 	{\widetilde{L}(\pr \circ Q\circ R)}
 \end{displaymath}
 where $\circ$ denotes the usual composition of graded linear maps.
\end{definition}
In other words $Q\gerst R$ is the unique coderivation with corestriction given by $p \circ Q\circ R$. 
Recall that $Q\circ R$ is not a coderivation in general (see remark \ref{Rem:CompositionofCoderivation}).

The graded map $\widetilde{L}$ defined in remark \ref{Notation:CoderLiftMapping} yields not only an isomorphism at the level of graded vector space but also one at the level of algebras:
\begin{lemma}[Gerstenhaber algebras isomorphism]\label{Lemma:GerstenAlgebrasIso}
	The graded vector spaces isomorphism
	\begin{displaymath}
		\widetilde{L}~:~ M(V) \equiv (\underline{\Hom}(\overline{T(V)},V), \gerst) ~\to~ (\coDer(\overline{T(V)}),\gerst)
	\end{displaymath}
	is an isomorphism in the category of graded algebras.
\end{lemma}
\begin{proof}
	Being $\widetilde{L}$ invertible, is sufficient to prove that the inverse $\pr$ preserve the product. For any fixed $n>0$ one has
	\begin{displaymath}
		\begin{aligned}
		p(Q\gerst R) \vert_{V^{\otimes n}} 
		\equal{\quad}&
		(q \circ R ) \vert_{V^{\otimes n}} =
		\\
		\equal{\quad}&		
		\sum_{s=1}^n q_s \circ (\sum_{i=0}^{s-1} \Unit_i \otimes r_{n-s+1} \otimes \Unit_{s-1-i} )=
		\\
		\equal{Def: \ref{Def:ithGerstenhaberProduct}}& \sum_{s=1}^n \sum_{i=0}^{s-1} q_s \gerst_i r_{n-s+1} =
		\\
		\equal{\quad}&	
		\sum_{s=1}^n q_s \gerst r_{n-s+1} = 
		\\
		\equal{\quad}&			
		(q \gerst r)_n = 
		\\
		\equal{\quad}&			
		(\pr\circ Q \gerst \pr\circ R)\vert_{V^{\otimes n}}
		~.
		\end{aligned}
	\end{displaymath}
\end{proof}
\begin{remark}
	We ought to notice that, in the literature (see \eg \cite{Manetti-website-coalgebras}), the definition of the Gerstenhaber product of multibrackets (see definition  \ref{Def:FullGerstenhaberProduct}) is often introduced via the previous lemma.
	In other words, the operator $\gerst$ on $M(V)$ is introduced by pulling back the product $\gerst$ on $\coDer(\overline{T(V)})$ along $\widetilde{L}$.
\end{remark}

\begin{remark}[$\gerst$ is pre-Lie]\label{Remark:GerstPreLie}
	Note that for any given $Q,R \in \coDer(\overline{T(V)})$ one has
	\begin{displaymath}
		\begin{aligned}
		[Q,R]_{\gerst} 
		=&~ 
		Q \gerst R - (-)^{|Q||R|} R \gerst Q =
		\\
		=&~
		\widetilde{L} \circ  p \circ \left( Q\circ R - (-)^{|Q||R|} R \circ Q \right) = 
		\\
		=&~\widetilde{L} \circ P \circ ([Q,R]_{\circ})
		\end{aligned}
	\end{displaymath}
	where the subscripts $\gerst$ and $\circ$ denotes the commutator bracket with respect to the Gerstenhaber product and  the composition of graded linear maps respectively.
	In particular, one has 
	\begin{displaymath}
		J_{\gerst}(Q,R,S) = \widetilde{L} \circ P (J_{\circ}(Q,R,S)) = L \circ P (0) = 0
	\end{displaymath}
	hence the Gerstenhaber product defined on the space of coderivation is pre-Lie.
	This observation, together with lemma \ref{Lemma:GerstenAlgebrasIso}, provide a more conceptual proof of the pre-Lie property of the Gerstenhaber product $\gerst$ on $M(V)$.
\end{remark}

\subsection{Algebra of coderivations on symmetric tensor spaces}\label{SubSection:CoderivationonSymmetricTensor}
Consider now a graded vector space $V$ and let us focus on the associated reduced symmetric tensor coalgebra $\overline{S(V)}= \oplus_{k\geq 1} V^{\odot k}$.
Recall that $\overline{S(V)}$ carries a natural coassociative coalgebra structure with coproduct given (see section \ref{Section:TensorCoalgebras}) by the unshuffled deconcatenation $\Xi$ 
	\begin{displaymath}
			\morphism{\Xi}
		{\overline{S(V)}}
		{\overline{S(V)}\otimes \overline{S(V)}}
		{x_1\odot\dots \odot x_n}
		{\displaystyle\sum_{i=1}^{n-1}
		\mkern-60mu\sum_{\mkern80mu\sigma\in \ush{i,n-i}}\mkern-60mu
		(x_1\odot\dots\odot x_i)\otimes (x_{i+1}\odot\dots\odot x_n)}
	\end{displaymath}
	where $\ush{i,j}\subset S_{i+j}$ denotes the subgroup of unshuffles permutations (see appendix \ref{App:UnshuffleAtors}).

Specializing theorem \ref{Theorem:HomCoAlgISOCoCommutative-Coder} to this specific cocommutative coalgebra yields a graded isomorphism between the space of coderivations and symmetric multilinear maps:
\begin{corollary}[of theorem \ref{Theorem:HomCoAlgISOCoCommutative-Coder}]
	Considering $C = \overline{S(V)}$ and $f=\pr\circ N=\id_V$, the unique lift prescribed by proposition \ref{Prop:UniversalPropertyCoderivation} is given by:
	\begin{displaymath}
		\mathclap{
		\morphism{\widetilde{L}_{\sym}~}
		{\underline{\Hom}^k(\overline{S(V)},V)}
		{\coDer^k(\overline{S(V)})}
		{q}
		{Q=\displaystyle \sum_{n\geq1} \sum_{s=1}^n ~ \pi\circ \left(
		q_{n-s+1}\otimes \Unit_{s-1}\right) \circ B_{n-s+1,s-1}  \circ N}
		}
	\end{displaymath}
	where $\pi$ and $N$ are the operators introduced in remark \ref{Remark:ManettiNotation}.
	On any given n-tuple of vectors $v_1,\dots,v_n$, one has explicitly that
	\begin{displaymath}
		Q(v_1\odot \dots \odot v_n)=
		\sum_{i=1}^n 
				\mkern-60mu\sum_{\mkern80mu\sigma\in \ush{i,n-i}}\mkern-60mu
		 \epsilon(\sigma) 
		q_i ( v_{\sigma_1}\odot\dots\odot v_{\sigma_k}) \odot v_{\sigma_{k+1}}\odot \dots \odot v_{\sigma_n}
		~,
	\end{displaymath}
	where $\ush{i,j}$ denotes as usual the $(i,j)$-unshuffles and $\epsilon(\sigma)$ is the Koszul sign.
\end{corollary}
\begin{proof}
 Applying the result of proposition \ref{Prop:SymLiftUnshufflesNotation} to this special case we get
 \begin{align*}
 		\widetilde{L}(q) =& \sum_{n> 0} H_n(q)\Big\vert_{\overline{S(V)}^{\odot n}} \circ  \Xi^n
 		=
 		\\
 		=&
 		\sum_{n> 0} \pi \circ (q \otimes \Unit_{n-1}) \circ B_{1,n-1} \circ N \circ \Xi^n
 		=
 		\\
 		=&
 		\sum_{n> 0}\sum_{m\geq n} \pi \circ (q \otimes \Unit_{n-1}) \circ B_{1,n-1} \circ N \circ \Xi^n \Big\vert_{V^{\odot m}}
 		=
 		\\
 		=&
 		 \sum_{m> 0} \pi \circ 
		 \left(
			 \sum_{s=1}^m
			 \left(
			 	q_{m-s+1}\otimes\Unit_{s-1}
			 \right)
			 \circ B_{m-s+1, s-1}
		 \right)
		 \circ N
		 ~,
 \end{align*}
 where $q_k$ denotes the corestriction of $q$ to $V^{\odot k}$, 
 $\Unit_k$ and $B_{1,k-1}$ are seen as operators on $\overline{S(V)}^{\odot k}$ in the penultimate equation and as operators on $V^{\odot k}$ in the last one.
\end{proof}

Similarly to what has been done in section \ref{Subsection:CoderivationsAlgebra}, we can equip the graded vector space $\coDer(\overline{S(V)})$ with an algebraic product:
\begin{definition}[\RN product of coderivations]\label{Def:SymmetricGerstProdOfCoderivations}
 We call \emph{\RN product} of coderivations the (non-associative) product given by
 \begin{displaymath}
 	\morphism{-\symgerst-}
 	{\coDer(\overline{S(V)})\otimes \coDer(\overline{S(V)}) }
 	{\coDer(\overline{S(V)})}
 	{A\otimes B}
 	{\widetilde{L}_{\sym}(\pr \circ A \circ B)}
 	~.
 \end{displaymath}
\end{definition}

Again as before, the graded map $\widetilde{L}_{\sym}$ defined in remark \ref{Notation:CoderLiftMappingCocommutative} yields not only an isomorphism at the level of graded vector space but also one at the level of non-associative algebras.
\begin{lemma}[Gerstenhaber algebras isomorphism]\label{Lemma:GerstenAlgebrasIsoSymmetric}
	The graded vector spaces isomorphism
	\begin{displaymath}
		\widetilde{L}_{\sym}~:~ M^{\sym}(V) \equiv (\underline{\Hom}(\overline{S(V)},V), \symgerst) ~\to~ (\coDer(\overline{S(V)}),\symgerst)
	\end{displaymath}
	is an isomorphism in the category of graded algebras.
\end{lemma}
\begin{proof}
	Being $\widetilde{L}_{\sym}$ invertible, it is sufficient to prove that the inverse $p$ preserves the product.
	Consider $R=\widetilde{L}_{\sym}(r)$ and $Q=\widetilde{L}_{\sym}(q)$ in $\coDer(\overline{S(V)})$.
	 For any fixed $n>0$ one has
	\begin{displaymath}
		\begin{aligned}
		\pr (R\symgerst Q) \vert_{V^{\otimes n}} \equal{Def: \ref{Def:SymGerstProd}}&
		(r \circ Q ) \vert_{V^{\otimes n}} =
		\\
		\equal{}&
		\sum_{s=1}^n r_s \circ \pi \left[(q_{n-s+1}\otimes \Unit_{s-1}) \circ B_{n-s+1,s-1}\right] \circ N =
		\\
		\equal{Def: \ref{Def:ithGerstenhaberProduct}}& 
		\sum_{s=1}^n r_s \symgerst q_{n-s+1} = 
		\\
		\equal{}&		
		(r \symgerst q)_n = 
		\\
		\equal{}&		
		\left((\pr\circ R) \symgerst (\pr\circ Q)\right)\vert_{V^{\otimes n}}
		~.
		\end{aligned}
	\end{displaymath}
\end{proof}
	Therefore, both $\coDer(S(V)$ and $M^{\sym(V)}= \bigoplus_{n\geq 0} \underline{\Hom}(V^{\odot n},V)$, when endowed with the \RN product $\cs$, form a right pre-Lie algebra.
	In particular, one has that the graded Lie brackets associated to any two $R,Q \in (\coDer(S(V)),\cs)$, \ie
	\begin{displaymath}
	[Q,R]:= Q\cs R - (-)^{|Q||R|} R \cs Q \in \coDer(S(V))
	~,
	\end{displaymath}
	coincides with the standard commutator of linear operators.

\section{The threefold nature of the Nijenhuis$-$Richardson algebra}\label{section:AppBConclusions}
	In this last section, we want to draw the conclusions of this long algebraic excursus.
	The previous sections' core was to identify an algebraic framework useful for formalizing expressions involving compositions of multilinear maps.  Synthetically, the core of this whole chapter can be subsumed by the following theorem:

 \begin{theorem}[Algebras of Multibrackets] \label{Theorem:RecapGerstenhaber}
	Let $V$ be a graded vector space, denote by $M^{\sym}(V)$ and $M^{\skew}(V)$ the graded vector spaces of symmetric and skew-symmetric multilinear maps from $V$ into itself.
	The following diagram commutes in the category of  graded right pre-Lie algebras
	  \begin{equation}\label{eq:diagrammone}
		 \begin{tikzcd}[row sep =large]
		 	(\overline{\overline{M(V[-1])}},\overline{\overline{\gerst}}) \ar[r,"\Dec","\cong"']
		 	\ar[d,two heads,"- \circ N_a"]
		 	& ({M(V)} ,\gerst) \ar[r,"\widetilde{L} ","\cong"'] 
		 	\ar[d,two heads,"- \circ N_s"]
		 	& (\coDer({T(V)}),\gerst)	 
		 	\ar[d,two heads,"\pi_s\circ - \circ N_a"]
		 	\\
		 	(M^{\skew}(V[-1]),\skewgerst) \ar[r,"\Dec","\cong"']
		 	\ar[d,hook,"- \circ \pi_a"]
		 	&  	(M^{\sym}(V),\symgerst) \ar[r,"\widetilde{L}_{\sym}","\cong"']
		 	\ar[d,hook,"-\circ \pi_s"]
		 	& (\coDer({S(V)}),\symgerst) 
		 	\ar[d,hook,"N_s \circ - \circ \pi_s"]
		 	\\
		 	(\overline{\overline{{M(V[-1])}}},\skewgerst)\ar[r,"\Dec","\cong"'] 
		 	&  	(M(V) ,\symgerst) \ar[r,"\widetilde{L} ","\cong"'] 
		 	& (\coDer(\overline{T(V)}),\symgerst)
		 	\\
		 \end{tikzcd} 
	 \end{equation}
	where:
	\begin{itemize}
		\item[-] $(\overline{M(V)} ,\gerst)$ denotes the \emph{Gerstenhaber algebra of multibrackets} (see remark \ref{rem:BigradedSpaceofMultilinearMaps} and definition \ref{Def:FullGerstenhaberProduct});
		\item[-]  $(\overline{\overline{M(V)}} ,\overline{\overline{\gerst}})$ denotes the algebra of graded multilinear maps taken with a different grading which intertwines arity and degree of the maps (see remark \ref{rem:DecasGradedMorph});
		\item[-] $\cs$ and $\ca$ denote respectively the symmetric and skew-symmetric \RN product (see definitions \ref{Def:SymGerstProd} and \ref{Def:SkewGerstProd} for the multibrackets case and definition \ref{Def:SymmetricGerstProdOfCoderivations} for the coderivations case);
		\item[-]  $\symgerst$ on $\coDer(\overline{T(V)}$ is defined by the following equation\footnote{We will not make any particular use of this algebra. Note that we are only saying that the algebra structure here is obtained by pulling back the product $\cs$ from $M(V)$.}
	 	\begin{displaymath}
	 	 Q \symgerst R = \widetilde{L}\left((\pr \circ Q) \symgerst(\pr \circ R)\right) ~;
		\end{displaymath}	
		\item[-] $\coDer(T(V))$ and $\coDer(S(V))$ denote the space of coderivations on the Tensor coalgebra and the symmetric Tensor coalgebra respectively;
		\item[-] $\pi_s: T(V)\to S(V)$ and $\pi_a:T(V)\to \Lambda(V)$ denote the canonical projections of $T(V)=S(V)\oplus\Lambda(V)$;
		\item[-] $N_s: S(V) \to T(V)$ and $N_a: \Lambda(V) \to T(V)$ denote the canonical injections;
		\item[-] $f\circ -$ denotes postcomposition with the function $f$ and $- \circ g$ denotes precomposition with function $g$;
		\item[-] $\Dec$ is the \emph{d\'ecalage} operator pertaining to multilinear maps (see remark \ref{Notation:FixingDec});
		\item[-] $\widetilde{L}$ and $\widetilde{L}_{\sym}$ are the lift operators to the tensor and symmetric tensor coalgebra respectively.
	\end{itemize}	 
 \end{theorem}
 	\begin{proof}
 		The horizontal isomorphism have been proved in proposition \ref{Prop:DecalageAsCoalgebrasISO} and lemmas \ref{Lemma:GerstenAlgebrasIso} and \ref{Lemma:GerstenAlgebrasIsoSymmetric}.
 		The central vertical line has been justified in theorem \ref{Thm:ManettiFactorizationOnGerstenhaberAlgebras}.
 		The two rightmost squares simply follows by observing that
 		\begin{displaymath}
 			\begin{aligned}
 			 \pr \circ (\pi\circ Q \circ N) =& p\circ Q \circ N \qquad \forall Q \in \coDer(\overline{T(V)})
 			 \\
 			 \pr \circ (N\circ R \circ \pi) =& p\circ R \circ \pi \qquad \forall R \in \coDer(\overline{S(V)})
 		 	\end{aligned}
 		\end{displaymath}
 		denoting with the same letter $\pr$ the canonical projections from $S(V)$ and $T(V)$ to $V$.
 	\end{proof}
	
	\begin{remark}[Coderivations on the anti-commutative tensor coalgebra]
		For the sake of completeness, we mention that one could add another column to the left of diagram \eqref{eq:diagrammone} introducing the algebra of coderivations on $\Lambda(V)$. 
		That case can be recovered, mutatis mutandis, from the analogous symmetric construction given in section \ref{SubSection:CoderivationonSymmetricTensor}.
	\end{remark}

	\begin{remark}\label{Remark:HomSVSVEmbedding}
		Observe that the reasoning of remark \ref{Remark:MSymEmbedding} can be easily extended to more general homogeneous map, namely one gets a an embedding diagram for any homogeneous maps between tensor spaces
		\begin{displaymath}
			\begin{tikzcd}[column sep = huge]
				\underline{\Hom}(\overline{T(V)},\overline{T(V)}) 
				&& \underline{\Hom}(\overline{T(V)},\overline{T(V)})
				 \ar[ll,"\mathcal{S\circ-\circ S}"'] \ar[d,two heads,"\pi-\circ-N"]
				\\
				\underline{\Hom}(\overline{S(V)},\overline{S(V)})
				\ar[u,hook,"N-\circ-\pi"]  \ar[urr,hook,"N-\circ-\pi"]
				&& 				\underline{\Hom}(\overline{S(V)},\overline{S(V)}) 				
				\ar[ll,"\sum_{n,m\geq 0} \frac{1}{n!m!} \pr_{V^{\odot n}} \circ - \pr_{V^{\odot n}}"]
			\end{tikzcd}
		\end{displaymath}
		where $\pr_{V^{\odot n}}$ denotes the canonical projector $S(V) \to V^{\odot n}$.	
		Will be clearer from theorem \ref{Theorem:RecapGerstenhaber} that this arrows factors also through coderivation seen as subspace of the space of homogeneous maps.
	\end{remark}

	\begin{remark}[THE \RN algebra]
		Since the horizontal lines in diagram \eqref{eq:diagrammone} are isomorphisms, 
		it would be completely legitimate to regard each line as a single object in the category of right pre-Lie algebras.
		In particular, the three isomorphic algebras contained in the central line can be thought as \emph{the} \RN algebra pertaining to the graded vector space $V[-1]$.
		\\
		One should think at the different descriptions of the \RN algebra (in term of symmetric multibrackets, skew-symmetric multibrackets or coderivation) as different presentation of the same algebraic objects.
	\end{remark}

\ifstandalone
	\bibliographystyle{../../hep} 
	\bibliography{../../mypapers,../../websites,../../biblio-tidy}
\fi

\cleardoublepage


%% file: chapters/permutators/permutators.tex
\chapter{Graded permutators and Unshuffles}\label{App:UnshuffleAtors}
In this appendix, we deal with the combinatorics involved when working with $L_\infty$-algebras.
Specifically, we will talk about unshuffle permutations and how they act on the tensor algebra $T(V)$ of a given graded vector space $V$.
The key results of this appendix are given in section \ref{Section:AppendixProofPreLie} which contains direct computations for the associator of the \RN algebras.

\section{Unshuffles}\label{subsection:UnshufflesAbstract}
	The term "shuffle" evokes the idea of shuffling a deck of cards. 
	Suppose we have a pack of $p$ cards and a pack of $q$ cards and we build a pack of $p+q$ cards, whilst retaining the order on the two "sub-packs".
	The result is a $(p,q)$-shuffle.
	A real-life shuffling consists of several repetitions of the previous operation.
	$(p,q)$-shuffles form a subgroup of the group $S_{p+q}$ of permutations of $p+q$ elements.
	The reverse of shuffling two decks of cards is what is called an "unshuffle".
	\begin{definition}[Unshuffles]\label{Def:Unshuffles}
		Let be $k_1,..k_\ell$ positive integers summing up to the integer $n$.
		Consider the group of permutation of $n$-elements $S_n$.
		The group of $(k_1,\dots,k_\ell)$-unshuffles is the subgroup of $S_n$ defined as follows:
		\begin{equation}\label{Eq:UnshufflesSet}
			\ush{k_1,k_2,\dots,k_\ell} : =
				\left\lbrace
					\sigma \in S_{k_1+\dots+k_\ell} \,\left|~
					\begin{aligned}
						&\sigma_i < \sigma_{i+1}\\
						&\forall i~\neq~ k_1,~ (k_1+k_2), \dots,~(k_1+\dots+k_{\ell})
					\end{aligned}
			\right.\right\rbrace
			~.
	\end{equation}	
\end{definition}

	\begin{reminder}[Cauchy representation and product of permutations]\label{rem:cauchyRep}
		Recall at first that any permutation $\sigma\in S_{k+\ell}$ can be denoted by its corresponding canonical Cauchy two-lines representation
	$$
	\sigma=
	\begin{pmatrix}
		1 & \dots & k & k+1 & \dots & k+\ell \\
		\sigma_1 & \dots & \sigma_k & \sigma_{k+1} & \dots & \sigma_{k+\ell}
	\end{pmatrix};
	$$
		(permutations are in a one-to-one correspondence with bijections on a finite set).
		Recall also that the direct product of two permutations can be denoted by the cartesian product of their canonical representation. 
		In other terms, given
	$$
	\mu=
	\begin{pmatrix}
		1 & \dots & k \\
		\mu_1 & \dots & \mu_k
	\end{pmatrix}
	\qquad
	\nu=
	\begin{pmatrix}
		1 & \dots & \ell \\
		\nu_1 & \dots & \nu_\ell
	\end{pmatrix}
	$$
	one has
	$$
	\mu\times \nu=
	\begin{pmatrix}
		1 & \dots & k & k+1 & \dots & k+\ell \\
		\mu_1 & \dots & \mu_k & \nu_{1} & \dots & \nu_{\ell}
	\end{pmatrix}~.
	$$
	\end{reminder}


	\begin{remark}\label{Rem:UnshufflesAsCoset}
		Observe that the group of unshuffles $\ush{p,q}$ is a set of representatives for the left cosets of the canonical embedding $S_p\times S_q \hookrightarrow S_{p+q}$. 
		In other words, for any $\eta \in S_{p+q}$ exists an unique decomposition $\eta = \sigma \circ \tau$ with $\sigma \in \ush{p,q}$ and $\tau \in S_p\times S_q$.
		\\
		In simpler terms, an unshuffle $\sigma\in\ush{p,q}$ rearrange the ordered list $(1,2,\dots,p+q)$ into a list that is separately ordered in the first $p$ and the second $q$ indexes.
	\end{remark}

\begin{remark}\label{Rem:UnshufflesCardinality}
	Recall that:
	$$
		\# \ush{k_1,\dots,k_\ell} = 
		\dfrac{(k_1+k_2+\dots+k_\ell)!}{k_1! k_2! \dots k_\ell!}
	$$
	since $S_{(k_1+\dots+k_\ell)}= S_{k_1}\times S_{k_2}\times \dots S_{k_\ell} \circ \ush{k_1,\dots,k_\ell}$.
	In particular
	\begin{equation}\label{Eq:UnshufflesCardinality}
		\# \ush{k,\ell} = \frac{(k+\ell)!}{k!\ell!} =
	\binom{k+\ell}{k} 
	\end{equation}
	where $\binom{\cdot}{\cdot}$ denotes the Newton binomial coefficient.
\end{remark}
\begin{remark}\label{Rem:UnshuffleSign}
	Observe that the generic element $\sigma \in \ush{2,n-2}$ can be written as 
	$$
		\sigma=(i,j;1,\dots,\hat{i},\dots,\hat{j},\dots,n) \qquad \text{for}~ 1\leq i < j \leq n~,
	$$
	therefore $|\sigma| = (-)^{i+j+1}$.
	Generalizing, any element $\sigma\in\ush{k,n-k}$ can be written as
	$$
		\sigma=(i_1,\dots,i_k;1,\dots,\hat{i_1},\dots,\hat{i_k},\dots,n) \qquad \text{for}~ 1\leq i_1 < i_2 <\dots< i_k \leq n~,
	$$
	with
	\begin{displaymath}
		\begin{aligned}
			|\sigma| =&~ (-)^{i_1-1}(-)^{i_2-2}\dots(-)^{i_k-k}
			\\
			=&~ (-)^{\sum_{j=1}^k i_j}(-)^{-\sum_{j=1}^k 1}=
			\\
			=&~ (-)^{\sum_{j=1}^k i_j}(-)^{\frac{k(k+1)}{2}}
			~.
		\end{aligned}
	\end{displaymath}
\end{remark}

	In defining the notion of $L_\infty$-morphism in the multibrackets presentation, one has also to deal with the following subgroup of the unshuffles:
	\begin{definition}[Ordered unshuffles]\label{Def:OrderedUnshufflesAbstract}
		Let be $k_1 \leq k_2 \dots \leq k_\ell$ non-zero integers summing up to the integer $n$.
		The group of $(k_1,\dots,k_\ell)$ \emph{ordered unshuffles} is the subgroup of $\ush{k_1,\dots,k_\ell}$ defined as follows:
		
	\begin{displaymath}\label{Eq:OrderedUnshufflesSet}
		\mathclap{
		\ush{k_1,\dots,k_\ell}^<
		=
		\left\lbrace
		\sigma \in \ush{k_1,k_2,\dots,k_\ell} \,\left|~
			\begin{aligned}
				&\sigma(k_1+\dots+k_{j-1}+1)<\sigma(k_1+\dots+k_{j}+1)
				\\[-1em] 
			 	&\qquad \forall j \text{~s.t.~} k_{j-1}=k_j
			\end{aligned}
		\right.\right\rbrace		~.
		}
	\end{displaymath}	
	\end{definition}

	When working with unshuffles is useful to introduce also the so-called \emph{Cyclic permutations}
	\begin{definition}[Cyclic permutation]\label{def:CyclycPermutation}
		We call \emph{cyclic permutation} (of $n$ elements) the permutation $\varsigma_{(n)}\in S_n$ given by the following Cauchy representation
		\begin{displaymath}
			\varsigma_{(n)}=
				\begin{pmatrix}
				1 & 2 &\dots & n-1 & n \\
				2 & 3 &\dots & n & 1 
			\end{pmatrix}~.
		\end{displaymath}
			We omit the subscript $(n)$ when there is no ambiguity on the number of elements that are cyclically permuted.
			We denote with $\varsigma^k$ the consecutive application of the cyclic permutation $k$-times.
	\end{definition}
		\begin{remark}\label{ref:cyclingUnhsuffles}
			Observe that an unshuffle $\sigma\in \ush{k,\ell}$ is a permutation that preserves the order of the first $k$ and the last $\ell$ elements, 
		then $\varsigma^k \circ \sigma$ is a permutation which preserves the order of the first $\ell$ and the last $k$ elements.
		In other words:
		\begin{displaymath}
			\Big\lbrace\varsigma^k \circ \sigma \,\Big|~ \sigma \in \ush{k,\ell}
			\Big\rbrace = \ush{\ell,k}	
			~.
		\end{displaymath} 
\end{remark}		

\begin{remark}[Decomposition of unshuffles]\label{rem:DecompositionOfUnshuffles}
		Recall that when a generic $(k,\ell)$-unshuffle acts on a list of elements $(1,2,\dots,k+\ell)$, in the resulting list the element $1$ can only appear in the first entry or at entry $(k+1)$. 
	This is due to the fact that index $1$ is the minimum in the total order $(1,\dots, k+\ell)$	 then, according to the definition of the group $\ush{k,\ell}$, it can only occur as the first element in the first batch of $k$ elements or the first in the second batch of $\ell$ elements.
	\\
	According to that, one can see that the group of $(k,\ell)$-unshuffles is the union of two subsets $\ush{k,\ell} = X \cup Y$ where 
	\begin{displaymath}
		X = \lbrace \mathbb{1}\times \tilde{\sigma} \quad\vert \quad \tilde{\sigma}\in \ush{k-1,\ell} \rbrace
	\end{displaymath}
	is the subgroup that keeps fixed the first index (\ie $\sigma_1 = 1$), and
	\begin{displaymath}
		Y = \left\lbrace (\varsigma_{(k+1)}\times \mathbb{1}_{\ell-1}) 
		\circ (\mathbb{1}\times \tilde{\sigma})
		\quad\left\vert \quad 
		\tilde{\sigma}\in \ush{k,\ell-1} \right.\right\rbrace
	\end{displaymath}
	is the subset that keeps the first index in the $k+1$ entry (i.e $\sigma_{k+1}=1$).
	Note that $\varsigma_{(k+1)}$ is the cyclic permutation of the first $k+1$ elements that takes in account the final position of the first index.
	Alternatively, elements in $Y$ can be expressed as
	\begin{displaymath}
		(\varsigma_{(k+1)}\times \mathbb{1}_{\ell-1}) \circ 
		(\mathbb{1}\times \tilde{\sigma}) = 
		(\mathbb{1}_k\times \varsigma_{(\ell)}^{-1})\circ
		(\tilde{\sigma}\times \mathbb{1}) \circ 
		(\varsigma_{(k+\ell)})
	\end{displaymath}
	where $\tilde{\sigma}\in \ush{k,\ell-1}$. 
	Here the first cyclic permutation is responsible for putting index $1$ in the last position and the second one is an inverse cyclic permutation of the last $\ell$ indices and is responsible of putting $1$ in position $k+1$.
\end{remark}

\section{Permutators on the unshuffles}

In section \ref{Section:ActionsonTensorSpaces} has been introduced the action of the permutation group on graded tensor spaces (tensor product of a given vector space with itself several times).
\\
The main subtlety with respect to ordinary vector space comes from the Koszul signs convention embedded in the graded braiding operator.

More precisely, fixed a graded vector space $L$ and given $\sigma \in S_n$, one can define two linear operator:
the even representation of $\sigma$:
\begin{displaymath}
	\morphism{B_\sigma}
	{L^{\otimes n}}
	{L^{\otimes n}}
	{v_1\otimes v_2\otimes \dots\otimes v_n}
	{\epsilon(\sigma,v_i)
	 v_{\sigma_1}\otimes v_{\sigma_2}\otimes \dots\otimes v_{\sigma_n}}
	~,
\end{displaymath}
where $	\epsilon(\sigma,v_i)$ is the Koszul sign (produced by the 
permutation of homogeneous graded elements),
and the odd representation of $\sigma$
\begin{displaymath}
	\morphism{P_\sigma}
	{L^{\otimes n}}
	{L^{\otimes n}}
	{v_1\otimes v_2\otimes \dots\otimes v_n}
	{\chi(\sigma,v_i)
	 v_{\sigma_1}\otimes v_{\sigma_2}\otimes \dots\otimes v_{\sigma_n}}
	~,
\end{displaymath}
where $\chi(\sigma,v_i)=(-)^\sigma	 \epsilon(\sigma,v_i)$.
Clearly one has that $P_\sigma = (-)^\sigma B_\sigma$.

Being the previous mappings two bona fide actions of the permutations group on $L^{\otimes n}$, they also satisfy the following property:
\begin{displaymath}
	B_\sigma \circ B_\sigma' = B_{\sigma\circ \sigma'}
\end{displaymath}

Considering a subset of the permutation group $\Omega \subset S_n$, one can define the permutator operators:
\begin{definition}[(even,odd) Permutator of $\Omega\subset S_n$]
	We call even (odd) permutator of $\Omega\subset S_n$ the linear operator
	$$B_{\Omega} = \sum_{\sigma \in \Omega} B_\sigma \qquad (P_{\Omega} = \sum_{\sigma \in \Omega}(-)^\sigma B_\sigma)$$
\end{definition}
The even (resp. odd) permutator corresponding to the complete cyclic group $\Omega =S_n$ takes the special name of \emph{simmetrizator}, denoted as $\mathcal{S}_n$ in definition \ref{Def:Symmetrizator} (resp. \emph{skew-simmetrizator} and denoted by $\mathcal{A}_n$).
The reason for such a naming comes from the fact that precomposing a multilinear operator with a (skew)-symmetrizator yields a graded (skew)-symmetric operator.
Note also that, for any given multilinear operator $H$ and any permutation $\sigma$, one has
\begin{displaymath}
	H \circ B_{\sigma} =
	\begin{cases}
		H & \quad \text{if $H$ is graded symmetric}\\
		(-)^{\sigma}	 H & \quad \text{if $H$ is graded skew-symmetric}
	\end{cases}
\end{displaymath}
and in particular, given any subset $\Omega < S_n$ and $H$ graded symmetric, one has $H \circ B_{\Omega} = (\#\Omega) H$.

In what concerns $L_\infty$-algebras, one has mainly to deal with the following three permutators (corresponding to the permutation subgroups introduced in subsection \ref{subsection:UnshufflesAbstract}):
\begin{definition}[Unshuffles permutator (Unshuffleator)]\label{Def:Unshuffleator}
	We call  even (odd) \emph{unshuffleator} the permutator
	\begin{displaymath}
			B_{k_1,k_2,\dots,k_\ell} := B_{\ush{k_1,k_2,\dots,k_\ell}}
			\qquad
			(P_{k_1,k_2,\dots,k_\ell} = P_{\ush{k_1,k_2,\dots,k_\ell}})
	\end{displaymath}	 
	pertaining to the subgroup of unshuffles (see definition \ref{Def:Unshuffles}).
\end{definition}
\begin{remark}
	Remark \ref{Rem:UnshufflesAsCoset} can be recast in terms of unshuffleators into the following equation
		\begin{equation}\label{Eq:DecompositionofSymmetrizator}
			\mathcal{S}_n = 
			\frac{\ell!(n-\ell)!}{n!}  
			\left(			
			\mathcal{S}_\ell \otimes \mathcal{S}_{n-\ell} 
			\right)
			\circ B_{n,n-\ell}
			~.
		\end{equation}
\end{remark}

\note{
	\Eg $B_{k,\ell}$ contains all the permutations that preserve the order of the first $k$ elements and of the last $\ell$ elements.
	\\
	A similar definition can be stated for the odd representation.
}

\begin{definition}[Ordered unshuffleator]\label{Def:OrderedUnshuffleator}
	We call  even (odd)  ordered \emph{unshuffleator} the permutator
	\begin{displaymath}
			B_{k_1,k_2,\dots,k_\ell}^< := B_{\ush{k_1,k_2,\dots,k_\ell}^<}
			\qquad
			(P_{k_1,k_2,\dots,k_\ell}^< = P_{\ush{k_1,k_2,\dots,k_\ell}^<})
	\end{displaymath}	 
	pertaining to the subgroup of ordered unshuffles (see definition \ref{Def:OrderedUnshufflesAbstract}).
	Recall in particular that $k_1 \leq k_2 \dots \leq k_\ell$.
\end{definition}

\begin{remark}
	Observe that when $\ell=2$ and $k_1=k_2=k$ one has
	$$ B_{k,k}^< = \mathbb{1}\otimes B_{k-1,k}$$
	and in particular
	$$ 
		B_{\underbrace{1,\dots,1}_{j~\text{times}},k}^<
		= B_{j,k} ~.
	$$	 
\end{remark}
%
%
\begin{definition}[Cyclic permutator]\label{Notation:Cyclic permutator}
	We call even \emph{cyclic permutator} the permutator
	\begin{displaymath}
			\cycPermutator_{(n)} = B_{\varsigma(n)}
	\end{displaymath}	 
	pertaining to the cyclic permutation of $n$ elements (see definition \ref{def:CyclycPermutation}).
	We denote as $\widetilde{\cycPermutator}_{(n)} = P_{\varsigma(n)}$ the corresponding odd permutator.
\end{definition}
\begin{remark}[The Koszul sign of the cyclyc permutator]\label{rem:KoszulSignCyclic}
	The cyclic permutation can be obtained by the consecutive application of $n-1$ contiguous swaps.
	Namely
	\begin{displaymath}
		\cycPermutator_{(n)} = 
		B_{n\leftrightarrow n-1}\circ 
		B_{n-1\leftrightarrow n-2}
		\dots B_{2\leftrightarrow 1}
		~,
	\end{displaymath}
	where $B_{n+1\leftrightarrow n}$ denotes the permutation swapping the $n$-the element with its immediate successor.
	Therefore  $(-)^{|\varsigma|} = (-)^{n-1}$.
	When evaluated, $\cycPermutator_n$ yields the following expressions:
	\begin{align*}
		\cycPermutator_{(n)} (v_1,\dots,v_n) &= (-)^{|v_1|(\sum_{i=2}^n|v_i|)}
		v_2\otimes\dots v_{n}\otimes v_1
		~,
		\\
		\cycPermutator_{(n)}^{-1} (v_1,\dots,v_n) &= (-)^{|v_n|(\sum_{i=1}^{n-1}|v_i|)}
		v_n \otimes v_1 \otimes \dots \otimes v_{n-1}
		~.
	\end{align*}
	Hence 
	\begin{align*}
		\epsilon(\varsigma,v_i)=&~(-)^{|v_1|(\sum_{i=2}^n|v_i|)} ~,
		\\
		\chi(\varsigma,v_i) =&~ (-)^{n-1}(-)^{|v_1|(\sum_{i=2}^n|v_i|)} ~.
	\end{align*}
	Note also that from the very definition of cyclic permutation one has that
	$\cycPermutator^{-\ell}_{(k+\ell)} = 	\cycPermutator^k_{(k+\ell)}  $. 
	\end{remark}
	
	\begin{lemma}[Properties of unshuffleators]
		Consider the even representation of the group of permutations, the following equations holds
		\begin{align}
			B_{k,\ell,m} 
			~=&~
			\left(B_{k,\ell}\otimes \mathbb{1}_m\right) \circ B_{k+\ell,m}
			\label{Eq:compUnsh}\\
			B_{\ell,k} 
			~=&~ 
			\cycPermutator^k_{(k+\ell)} \circ B_{k,\ell} 
			\label{Eq:CyclicUnshuffle}\\
			B_{k,\ell} 
			~=&~
			\mathbb{1}\otimes B_{k-1,\ell} +
			(\mathbb{1}_k\otimes\cycPermutator_{(\ell)}^{-1}) \circ 
			(B_{k,\ell-1}\otimes\mathbb{1}) \circ 
			(\cycPermutator_{(k+\ell)})
			\label{Eq:decompUnsh}
		\end{align}
		The same equations hold for the odd representation by replacing $B_\sigma $ with $P_\sigma$ and $\cycPermutator_{(n)}$ with $\widetilde{\cycPermutator}_{(n)}$.
	\end{lemma}
	\begin{proof}
		The statement follows from the properties of the group of unshuffles.
	\begin{itemize}
		\item 	Equation \eqref{Eq:compUnsh} follows from the observation that a permutation $\sigma \in \ush{k,\ell,m}$, being a permutation which preserves the order of the first $k$, the second $\ell$ and the last $m$ elements, can be decomposed as 
		$\sigma = (\alpha \otimes \mathbb{1}_m) \circ \beta$ with $\alpha \in \ush{k,\ell}$ and $\beta \in \ush{k+\ell,m}$.
		
		\item 	
		Equation \eqref{Eq:CyclicUnshuffle} follows from remark \ref{ref:cyclingUnhsuffles}

		\item Equation \eqref{Eq:decompUnsh} follows from remark \ref{rem:DecompositionOfUnshuffles}.

	\end{itemize}
	\end{proof}
	\begin{example}
		Equation \eqref{Eq:decompUnsh} could be better understood through examples. Inspecting the action of the unshuffleator on vectors $v_1,v_2,\dots$, denoted simply as $1,2,\dots$, one gets
		\begin{align*}
			B_{2,2} (1,2,3,4) 
			=&
			[ + (1,2;3,4) - (1,3;2,4) + (1,4;2,3) ] +
			\\
			&	[ + (3,4;1,2) - (2,4;1,3) + (2,3;1,4) ] 	
			=
			\\
			=& (1)\otimes P_{1,2} (2,3,4) 
			- (\mathbb{1}_2\otimes\cycPermutator_{2})
			\circ \left(P_{2,1}(2,3,4)\otimes (1)\right)
			\\
			=& \mathbb{1}\otimes P_{1,2} (1,2,3,4) 
			- (\mathbb{1}_2\otimes\cycPermutator_{2})
			\left(\circ P_{2,1}\otimes \mathbb{1} \right)
			(2,3,4,1)
			\\
			=& \mathbb{1}\otimes P_{1,2} (1,2,3,4) 
			+ (\mathbb{1}_2\otimes\cycPermutator_{(2)})\circ 
			\left(P_{2,1}\otimes \mathbb{1} \right)
			\circ\cycPermutator_{(4)}	
			(1,2,3,4)
			~;
		\end{align*}
		\begin{align*}
			B_{2,3} (1,2,3,4,5) 
			=&
			[ + (1,2;3,4,5) - (1,3;2,4,5) + (1,4;2,3,5) - (1,5;2,3,4) ] +
			\\
			&	[ + (2,3;1,4,5) - (2,4;1,3,5) + (2,5;1,3,4) +
			\\
			&	\phantom{[} + (3,4;1,2,5) - (3,5;1,2,4) + (4,5;1,2,3)] 	
						~.
		\end{align*}
	\end{example}

\section{Explicit expressions for the associators of the Nijenhuis$-$Richardson products}\label{Section:AppendixProofPreLie}
	We focus now on the actions of the previous permutations on the tensor products of homogeneous multilinear maps.
	In this section we deliver an explicit proof of the pre-Lie property of \RN products (see propositions \ref{Prop:SymmetricGerstenhaberAssociators} and \ref{prop:RNExplictPreLie}).

	\begin{lemma}\label{Lemma:CyclicCommutation}
	Given any homogeneous multilinear map $\mu_k: L^{\otimes k}\to L$, one has
	\begin{displaymath}
		\cycPermutator_{(\ell+1)}^{-1} \circ 
		\left(\mathbb{1}_\ell \otimes \mu_k \right) \circ 
		\cycPermutator_{(\ell+k)}^k 
		= 
		\mu_k \otimes \mathbb{1}_\ell
	\end{displaymath}
	and
	\begin{displaymath}
		\widetilde{\cycPermutator}_{(\ell+1)}^{-1} \circ 
		\left(\mathbb{1}_\ell \otimes \mu_k \right)\circ 
		\widetilde{\cycPermutator}_{(\ell+k)}^k 
		= (-)^{\ell(k+1)}\mu_k \otimes \mathbb{1}_\ell
	\end{displaymath}
\end{lemma}
\begin{proof}
	By inspection, let us denote as $1,\dots,i,\dots$ the graded vectors $v_1,\dots,v_i,\dots$. 
	From the definitions of cyclic permutator, one has that
	\begin{displaymath}
	\begin{aligned}
		\cycPermutator_{(\ell+k)} &(1,\dots,k,k+1,\dots,k+\ell)
		=
		\\
		=&~
		(-)^{|1|\cdot(|2|+\dots+|k+\ell|)}
		(2,\dots,k,k+1,\dots,k+\ell,1)
		~.
	\end{aligned}
	\end{displaymath}
	Iterating, one gets
	\begin{displaymath}
	\begin{aligned}
		\cycPermutator_{(\ell+k)}^2 &(1,\dots,k,k+1,\dots,k+\ell)
		=
		\\
		=&~
		(-)^{|1|\cdot(|2|+\dots+|k+\ell|)}
		\\
		&~(-)^{|2|\cdot(|1|+|3|+\dots+|k+\ell|)}
		(3,\dots,k,k+1,\dots,k+\ell,1,2)
	\end{aligned}
	\end{displaymath}
	and
	\begin{align*}
		\cycPermutator_{(\ell+k)}^k &(1,\dots,k,k+1,\dots,k+\ell)
		=
		\\
		=&~
		(-)^{|1|\cdot(|2|+\dots+|k+\ell|)}
		\\
		&~(-)^{|2|\cdot(|1|+|3|+\dots+|k+\ell|)}
		\\
		&~\vdots
		\\
		&~(-)^{|k|\cdot(|k+1|+\dots+|k+\ell| + |1| + \dots + |k-1|)}
		(k+1,\dots,k+\ell,1,\dots,k)
		=
		\\
		=&~
		(-)^{|1|\cdot(|k+1|+\dots+|k+\ell|)}
		(-)^{|1|\cdot(|2|+\dots+|k|)}
		\\
		&~(-)^{|2|\cdot(|k+1|+\dots+|k+\ell|)}
		(-)^{|2|\cdot(|1|+|3|+\dots+|k|)}		
		\\
		&~\vdots
		\\
		&~(-)^{|k|\cdot(|k+1|+\dots+|k+\ell|)}
		(-)^{|k|\cdot(|1|+\dots+|k-1|)}	
		(k+1,\dots,k+\ell,1,\dots,k)
		=
		\\
		=&~
		(-)^{(|k+1|+\dots+|k+\ell|)(|1|+\dots+|k|)}
		(k+1,\dots,k+\ell,1,\dots,k)	
		~.
	\end{align*}	
	Recall that the Koszul sign convention implies
	\begin{align*}
		\left(\mathbb{1}_\ell \otimes \mu_k \right)&(k+1,\dots,k+\ell,1,\dots,k)
		=
		\\
		=&~
		(-)^{|\mu_k|(|k+1|+\dots+|k+\ell|)}
		(k+1,\dots,k+\ell,\mu_k(1,\dots,k))	
		~,
	\end{align*}		
	therefore
	\begin{align*}
		\cycPermutator_{\ell+1}^{-1}& (k+1,\dots,k+\ell, \mu_k(1,\dots,k)) =
		\\
		=&~
		(-)^{|\mu_k(1,\dots,k)|\cdot(|k+1|+\dots+|k+\ell|)} (\mu_k(1,\dots,k),k+1,\dots,k+\ell) =
		\\
		=&~
		(-)^{(|\mu_k| +|1| +\dots+|k|)\cdot(|k+1|+\dots+|k+\ell|)} 
		\cdot \left(\mu_k \otimes \mathbb{1}_\ell\right) (1,\dots,k,k+1,\dots,k+\ell)		
		~.
	\end{align*}
	Composing all the previous steps we get
	\begin{displaymath}
	\begin{aligned}
		\left(\cycPermutator^{-1}_{\ell+1}\circ (\mathbb{1}_\ell \otimes \mu_k) \circ \cycPermutator^k_{\ell+k}\right)  &(1,\dots,k+\ell)
		=
		\\
		=&~
		\cancel{(-)^{|\mu_k|(|k+1|+\dots+|k+\ell|)}}
		\\
		&~
		\cancel{(-)^{(|1|+\dots+|k|)(|k+1|+\dots+|k+\ell|)}}
		\\
		&~ \cancel{ (-)^{|\mu_k|(|k+1|+\dots+|k+\ell|)}} 
		\\
		&~
		\cancel{(-)^{(|1|+\dots+|k|)(|k+1|+\dots+|k+\ell|)}} 
		\\
		& \mu_k \otimes \mathbb{1}_\ell (1,\dots,k,k+1,\dots,k+\ell)		 
	\end{aligned}
	\end{displaymath}
	which yields the first equation.
	The extra sign in the second equation follows from noting that $\widetilde{\cycPermutator}_n = (-)^{n-1} \cycPermutator_n$ (see remark \ref{rem:KoszulSignCyclic}).
\end{proof}
	Out of the cyclic permutator one can express a commutation rule involving the tensor product of multilinear maps:	
	\begin{corollary}\label{Cor:CyclicCommutationMultiBracket}
		Given any two homogeneous multilinear maps $\mu_k: L^{\otimes k}\to L$
		and $\eta_\ell:L^{\otimes\ell}\to L$, one has
		\begin{displaymath}
			\begin{aligned}
				\cycPermutator_{(2)} \circ \left(\mu_k\otimes\eta_\ell\right) \circ \cycPermutator^k_{(k+\ell)} 
				&=
				(-)^{|\eta_\ell||\mu_k|} \left(\eta_\ell \otimes \mu_k\right)
				\\
				\widetilde{\cycPermutator}_{(2)} \circ \left(\mu_k\otimes\eta_\ell \right)\circ \widetilde{\cycPermutator}^k_{(k+\ell)} 
				&=
				(-)^{k\ell +1}(-)^{|\eta_\ell||\mu_k|} \left(\eta_\ell \otimes \mu_k\right)
			\end{aligned}
			~.
		\end{displaymath}		
	\end{corollary}
	\begin{proof}
		Employing the previous lemma twice, one gets the first equation easily
		\begin{displaymath}
			\begin{aligned}
				\cycPermutator_{(2)} \circ 
				( \mu_k\otimes \eta_\ell) 
				\circ \cycPermutator^\ell_{(\ell+k)}
				=&~
				\cycPermutator_{(2)} \circ 
				(\mu_k\otimes\mathbb{1}) \circ 
				(\mathbb{1}_k\otimes \eta_\ell)
				\circ \cycPermutator^\ell_{(\ell+k)}=
				\\
				=&~	
				(\mathbb{1}\otimes \mu_k) \circ
				\cancel{
					\cycPermutator^k_{(k+1)} \circ
					\cycPermutator_{(k+1)} \circ
				}
				(\eta_\ell \otimes \mathbb{1}_k)
				=					
				\\
				=&~	
				(-)^{|\mu_k||\eta_\ell|} \eta_\ell \otimes \mu_k
				~,
			\end{aligned}
		\end{displaymath}
		where the last sign comes again from the Koszul convention on homogeneous maps, see equation \ref{Eq:TensorHomogeneousMaps}.
		\\
		Regarding the second equation, one has only to note that
		\begin{displaymath}
			\widetilde{\cycPermutator}_{(2)} \circ \left(\mu_k\otimes\eta_\ell\right) \circ \widetilde{\cycPermutator}^k_{(k+\ell)}
			=
			-(-)^{(k+\ell+1)k}
			\cycPermutator_{(2)} \circ \left(\mu_k\otimes\eta_\ell \right)\circ \cycPermutator^k_{(k+\ell)}
			~.
		\end{displaymath}
	\end{proof}

\begin{lemma}\label{Lemma:UnshuffleSandwich}
\begin{align*}
	(B_{m,k}) \circ & (\mu_n \otimes \mathbb{1}_{m+k-1}) \circ (B_{n,m+k-1}) 
	=
	\\
	=&~
	(\mu_n \otimes \mathbb{1}_{m+k-1}) \circ
	(B_{n,m-1,k}) +
		\left(\mathbb{1}_m\otimes
		\mu_n\otimes \mathbb{1}_{k-1}
		\right) \circ
		(B_{m,n,k-1})
		~;
	\\[.2em]
	(P_{m,k}) \circ& (\mu_n \otimes \mathbb{1}_{m+k-1}) \circ (P_{n,m+k-1}) 
	=
	\\
	=&~ 
	(\mu_n \otimes \mathbb{1}_{m+k-1}) \circ
			(P_{n,m-1,k}) +
		(-)^{m(n+1)}~
		\left(\mathbb{1}_m\otimes
		\mu_n\otimes \mathbb{1}_{k-1}
		\right) \circ
		(P_{m,n,k-1})	
		~.
\end{align*}

\end{lemma}
\begin{proof}
	Applying equation \eqref{Eq:decompUnsh} to $B_{m,k}$, one can write
	\begin{align*}
			(B_{m,k}) \circ &(\mu_n \otimes \mathbb{1}_{m+k-1}) \circ (B_{n,m+k-1})
			=
			\\
			=&~
			+(\mu_n\otimes B_{m-1,k})\circ (B_{n,m+k-1}) +
			\\
			&~	
			+\left[
			(\mathbb{1}_m\otimes\cycPermutator_{(k)}^{-1}) \circ 
			(B_{m,k-1}\otimes \mathbb{1}) \circ \cycPermutator_{(m+k)}
			\right]
			\circ \left(\mu_n \otimes \mathbb{1}_{m+k-1}\right)
			\circ (B_{n,m+k-1}) =
			\\
			=&~ X + Y
			~.
	\end{align*}
	From equation \eqref{Eq:compUnsh}, the first summand yields:
	\begin{displaymath}
		X = (\mu_n \otimes \mathbb{1}_{m+k-1}) \circ
			(B_{n,m-1,k}) ~;
	\end{displaymath}
	the second summand results:
	
	\begin{displaymath}
		\mathclap{
	\begin{aligned}
		Y 
		\equal{asso.}&
		(\mathbb{1}_m\otimes\cycPermutator_{(k)}^{-1}) \circ 
		(B_{m,k-1}\otimes \mathbb{1}) \circ 
		\left[
			\cycPermutator_{(m+k)}
			\circ \left(\mu_n \otimes \mathbb{1}_{m+k-1}\right)
		\right]
		\circ (B_{n,m+k-1}) =
		\\		
		\equal{Lem. \ref{Lemma:CyclicCommutation}}&
		(\mathbb{1}_m\otimes\cycPermutator_{(k)}^{-1}) \circ 
		(B_{m,k-1}\otimes \mathbb{1}) \circ 
		\left[
			\left(\mathbb{1}_{m+k-1}\otimes\mu_n \right) \circ
			\cycPermutator_{(n+m+k-1)}^n
		\right]
		\circ (B_{n,m+k-1}) =
		\\
		\equal{asso.}&
		(\mathbb{1}_m\otimes\cycPermutator_{(k)}^{-1}) \circ 
		(\mathbb{1}_{m+k-1}\otimes\mu_n) \circ
		(B_{m,k-1}\otimes \mathbb{1}) \circ 
		\left[
			\cycPermutator_{(n+m+k-1)}^n \circ
			(B_{n,m+k-1}) 
		\right]
		=
		\\
		\equal{Eq. \eqref{Eq:CyclicUnshuffle}}&
		(\mathbb{1}_m\otimes\cycPermutator_{(k)}^{-1}) \circ 
		(\mathbb{1}_{m+k-1}\otimes\mu_n) \circ
		(B_{m,k-1}\otimes \mathbb{1}) \circ 
		(B_{m+k-1,n})
		=
		\\
		\equal{Eq. \eqref{Eq:compUnsh}}&
		(\mathbb{1}_m\otimes\cycPermutator_{(k)}^{-1}) \circ 
		(\mathbb{1}_{m+k-1}\otimes\mu_n) \circ
		(B_{m,k-1,n})
		=
		\\
		\equal{asso.}&
		\left[\mathbb{1}_m\otimes
			\left(\cycPermutator_{(k)}^{-1} \circ \left(\mathbb{1}_{k-1}\otimes\mu_n\right)\right)
		\right] \circ 
		(B_{m,k-1,n})
		=
		\\
		\equal{Lem. \ref{Lemma:CyclicCommutation}}&
		\left(\mathbb{1}_m\otimes
		\mu_n\otimes \mathbb{1}_{k-1}
		\right) \circ
		(\mathbb{1}_m\otimes \cycPermutator^{-n}_{(n+k-1)})\circ 
		(B_{m,k-1,n})
		=
		\\
		\equal{Eq. \eqref{Eq:CyclicUnshuffle}}&
		\left(\mathbb{1}_m\otimes
		\mu_n\otimes \mathbb{1}_{k-1}
		\right) \circ
		(B_{m,n,k-1})	
		~.
	\end{aligned}
	}
	\end{displaymath}

	Similarly, applying equation \eqref{Eq:decompUnsh} to $P_{m,k}$, one can write
	\begin{align*}
		(P_{m,k}) \circ& (\mu_n \otimes \mathbb{1}_{m+k-1}) \circ (P_{n,m+k-1})
		=
		\\
		=&
		+(\mu_n\otimes P_{m-1,k})\circ (P_{n,m+k-1}) +
		\\
		&	
		+\left[
		(\mathbb{1}_m\otimes\widetilde{\cycPermutator}_{(k)}^{-1}) \circ 
		(P_{m,k-1}\otimes \mathbb{1}) \circ \widetilde{\cycPermutator}_{(m+k)}
		\right]
		\circ \left(\mu_n \otimes \mathbb{1}_{m+k-1}\right)
		\circ (P_{n,m+k-1}) =
		\\
		=&~ X + Y
		~.
	\end{align*}
	From equation \eqref{Eq:compUnsh}, the first summand yields:
	\begin{displaymath}
		X = (\mu_n \otimes \mathbb{1}_{m+k-1}) \circ
			(P_{n,m-1,k}) ~;
	\end{displaymath}
	the second summand results:
	\begin{displaymath}
	\begin{aligned}
		Y 
		\equal{asso.}&
		(\mathbb{1}_m\otimes\widetilde{\cycPermutator}_{(k)}^{-1}) \circ 
		(P_{m,k-1}\otimes \mathbb{1}) \circ 
		\left[
			\widetilde{\cycPermutator}_{(m+k)}
			\circ \left(\mu_n \otimes \mathbb{1}_{m+k-1}\right)
		\right]
		\circ (P_{n,m+k-1}) =
		\\		
		\equal{Lem. \ref{Lemma:CyclicCommutation}}&
		(\mathbb{1}_m\otimes\widetilde{\cycPermutator}_{(k)}^{-1}) \circ 
		(P_{m,k-1}\otimes \mathbb{1}) \circ 
		\\
		\phantom{\equal{}}&
		\left[
			(-)^{(m+k-1)(n+1)}\circ
			\left(
			\mathbb{1}_{m+k-1}\otimes\mu_n \right)
			\circ
			\widetilde{\cycPermutator}_{(n+m+k-1)}^n
		\right]\circ
		(P_{n,m+k-1}) =
		\\
		\equal{asso.}&
		(-)^{(m+k-1)(n+1)}~
		(\mathbb{1}_m\otimes\widetilde{\cycPermutator}_{(k)}^{-1}) \circ 
		(\mathbb{1}_{m+k-1}\otimes\mu_n) \circ
		(P_{m,k-1}\otimes \mathbb{1}) \circ 
		\\
		\phantom{\equal{}}&\circ
		\left[
			\widetilde{\cycPermutator}_{(n+m+k-1)}^n \circ
			(P_{n,m+k-1}) 
		\right]
		=
		\\
		\equal{Eq. \eqref{Eq:CyclicUnshuffle}}&
		(-)^{(m+k-1)(n+1)}~
		(\mathbb{1}_m\otimes\widetilde{\cycPermutator}_{(k)}^{-1}) \circ 
		(\mathbb{1}_{m+k-1}\otimes\mu_n) \circ
		(P_{m,k-1}\otimes \mathbb{1}) \circ 
		\\
		\phantom{\equal{}}&\circ		
		(P_{m+k-1,n})
		=
		\\
		\equal{Eq. \eqref{Eq:compUnsh}}&
		(-)^{(m+k-1)(n+1)}~
		(\mathbb{1}_m\otimes\widetilde{\cycPermutator}_{(k)}^{-1}) \circ 
		(\mathbb{1}_{m+k-1}\otimes\mu_n) \circ
		(P_{m,k-1,n})
		=
		\\
		\equal{asso.}&
		(-)^{(m+k-1)(n+1)}~
		\left[\mathbb{1}_m\otimes
			\left(\widetilde{\cycPermutator}_{(k)}^{-1} \circ \left(\mathbb{1}_{k-1}\otimes\mu_n\right)\right)
		\right] \circ 
		(P_{m,k-1,n})
		=
		\\
		\equal{Lem. \ref{Lemma:CyclicCommutation}}&
		(-)^{(m+k-1)(n+1)+(k-1)(n+1)}~
		\left(\mathbb{1}_m\otimes
		\mu_n\otimes \mathbb{1}_{k-1}
		\right) \circ
		\\
		\phantom{\equal{}}&\circ
		(\mathbb{1}_m\otimes \widetilde{\cycPermutator}^{-n}_{(n+k-1)})\circ 
		(P_{m,k-1,n})
		=
		\\
		\equal{Eq. \eqref{Eq:CyclicUnshuffle}}&
		(-)^{m(n+1)}~
		\left(\mathbb{1}_m\otimes
		\mu_n\otimes \mathbb{1}_{k-1}
		\right) \circ
		(P_{m,n,k-1})	
		~.
	\end{aligned}
	\end{displaymath}


\end{proof}	
We are now ready to compute the associator of the products $\symgerst$ and $\skewgerst$ defined in section \ref{Section:MultibracketsAlgebra}.
	\begin{proposition}[Explicit associator of $\symgerst$ and $\skewgerst$]\label{Prop:ExplicitAssociators}
		Given any three multilinear operators $\mu_\ell,\mu_m,\mu_n \in M(V)$ the corresponding associators with respect to $\symgerst$ and $\skewgerst$ result:
		\begin{equation}
		\label{Eq:ExplicitAssociatorSym}
			\alpha(\symgerst;\mu_\ell,\mu_m,\mu_n) 
			=
			\mu_\ell \circ \left[\left(\left(\mu_m\otimes\mu_n \right)\circ B_{m,n}\right)\otimes \Unit_{\ell-2}\right] \circ B_{m+n,\ell-2}
		\end{equation}
		\begin{equation}\label{Eq:ExplicitAssociatorSkew}
			\alpha(\skewgerst; \mu_\ell,\mu_m,\mu_n) 
			=
				(-)^{\mathcal{s}}\mu_\ell \circ \left[\left(\left(\mu_m\otimes\mu_n \right)\circ P_{m,n}\right)\otimes \Unit_{\ell-2}\right] \circ P_{m+n,\ell-2}~,
		\end{equation}
		where the sign prefactor is given by:
		\begin{displaymath}
			(-)^{\mathcal{s}}=(-)^{{|\mu_n|(m+\ell)} +{|\mu_m|(\ell-1)} +m(n+1)}
		\end{displaymath}
	\end{proposition}
\begin{proof}
	By definition, the associator reads as follows:
	\begin{displaymath}
		\alpha(\symgerst;\mu_\ell,\mu_m,\mu_n)
		\colon=
		(\mu_\ell \symgerst \mu_m) \symgerst \mu_n -
		\mu_\ell \symgerst ( \mu_m \symgerst \mu_n) 
		~.
	\end{displaymath}
	The second term on the right-hand side results:
	\begin{align*}
		\mu_\ell \symgerst& ( \mu_m \symgerst \mu_n) 
		=~
		\mu_\ell \symgerst \left( \mu_m \circ \left(\mu_n \otimes \mathbb{1}_{m-1}\right) \circ B_{n,m-1}\right)
		=\\
		=&~
		\mu_\ell \circ \left[\left( \mu_m \circ \left( \mu_n \otimes \mathbb{1}_{m-1}\right) \circ B_{n,m-1}\right)
		\otimes \mathbb{1}_{\ell-1}\right] \circ B_{m+n-1,\ell-1}
		=\\
		=&~			
		(\mu_\ell) \circ 
		( \mu_m\otimes \mathbb{1}_{\ell-1}) \circ 
		(\mu_n \otimes \mathbb{1}_{m+\ell-2}) \circ
		(B_{n,m-1} \otimes \mathbb{1}_{\ell-1}) \circ 
		B_{m+n-1,\ell-1}
		=\\
		=&~			
		(\mu_\ell) \circ 
		( \mu_m\otimes \mathbb{1}_{\ell-1}) \circ 
		(\mu_n \otimes \mathbb{1}_{m+\ell-2}) \circ
		B_{n,m-1,\ell-1}
		~.
	\end{align*}
	The first term results
	\begin{align*}
		(\mu_\ell \symgerst  \mu_m)& \symgerst \mu_n 
		=~
		\left(\mu_\ell \circ \left(\mu_m \otimes \mathbb{1}_{\ell-1}\right) \circ B_{m,\ell-1} \right) 
		\symgerst \mu_n	
		=\\
		\equal{}&~
		\left(\mu_\ell \circ \left(\mu_m \otimes \mathbb{1}_{\ell-1} \right)\circ B_{m,\ell-1} \right) 
		\circ \left(\mu_n \otimes \mathbb{1}_{m+\ell-2} \right) 
		\circ B_{n,m+\ell-2}	
		=\\
		\equal{}&~			
		(\mu_\ell) \circ 
		(\mu_m \otimes \mathbb{1}_{\ell-1}) \circ 
		\big[
		(B_{m,\ell-1} ) \circ 
		(\mu_n \otimes \mathbb{1}_{m+\ell-2}) \circ 
		(B_{n,m+\ell-2})
		\big]
		=
		\\
		\equal{Lem. \ref{Lemma:UnshuffleSandwich}}&
		+(\mu_\ell) \circ 
		(\mu_m \otimes \mathbb{1}_{\ell-1}) \circ 
		(\mu_n\otimes\mathbb{1}_{m+\ell-2}) \circ
		B_{n,m-1,\ell-1}
		+
		\\
		& + s\cdot 
		(\mu_\ell) \circ 
		(\mu_m \otimes \mathbb{1}_{\ell-1}) \circ 
		(\mathbb{1}_m\otimes \mu_n \otimes \mathbb{1}_{\ell-2})	\circ
		(B_{m,n,\ell-2})
		=
		\\		
		\equal{}&
		+ \mu_\ell \symgerst ( \mu_m \symgerst \mu_n) 
		+
		\\
		&+ s\cdot 
		(\mu_\ell) \circ 
		(\mu_m \otimes \mathbb{1}_{\ell-1}) \circ 
		(\mathbb{1}_m\otimes \mu_n \otimes \mathbb{1}_{\ell-2})	\circ
		(B_{m,n,\ell-2})
		~,
	\end{align*}
	where $s$ is, according to lemma \ref{Lemma:UnshuffleSandwich}, a numerical constant equal to one.
	Therefore
	\begin{align*}
			\alpha(\symgerst;\mu_\ell,&\mu_m,\mu_n) 
			=~
			\\
			\equal{}&
			s \cdot(\mu_\ell) \circ  
			(\mu_m \otimes \mathbb{1}_{\ell-1}) \circ 
			(\mathbb{1}_m\otimes \mu_n \otimes \mathbb{1}_{\ell-2})	\circ
			(B_{m,n,\ell-2})
			=
			\\		
			\equal{Eq. \eqref{Eq:compUnsh}}&
			s \cdot(\mu_\ell) \circ 
			(\mu_m \otimes \mu_n \otimes \mathbb{1}_{\ell-2})	\circ
			(B_{m,n} \otimes \mathbb{1}_{\ell-2})	\circ
			(B_{m+n,\ell-2})
			=
			\\		
			\equal{}&
			s \cdot(\mu_\ell) \circ 
			\left[\left(\left(\mu_m \otimes \mu_n\right) \circ
			B_{m,n}\right)
			 \otimes \mathbb{1}_{\ell-2} \right]	\circ
			(B_{m+n,\ell-2})
			~.
	\end{align*}
		The odd case is simply obtained substituting $B_{k,\ell}$ with $P_{k,\ell}$, $s$ with the sign $(-)^{m(n+1)}$ given by lemma \ref{Lemma:UnshuffleSandwich} and remembering that the definition of skewsymmetric Gerstenhaber takes an extra sign. Namely:
		\begin{displaymath}
			(\mu_\ell \skewgerst \mu_m) \skewgerst \mu_n = 
			(-)^\mathcal{s}~
			\mu_\ell \circ \left[ \mu_m \circ \left(\mu_n \otimes \mathbb{1}_{m-1}\right) \circ P_{n,m-1}\right]
			\otimes \mathbb{1}_{\ell-1}) \circ P_{m+n-1,\ell-1}
			~,
		\end{displaymath}
		where
		\begin{displaymath}
			\mathcal{s}=
			{|\mu_n|(m+\ell-2)} +{|\mu_m|(\ell-1)}
			~.
		\end{displaymath}
	\end{proof}	
	\begin{proposition}[Pre-Lie property of $\symgerst$ and $\skewgerst$]
	\label{Prop:SymmetricGerstenhaberAssociators-AppendixC}
		Given any three multilinear operators $\mu_\ell,\mu_m,\mu_n \in M(V)$, the corresponding associators satisfy the following symmetry properties in the last two entries:
		\begin{equation}\label{Equation:SymGerstAssociatorSymmetry}
			\alpha(\symgerst; \mu_\ell,\mu_m,\mu_n) 
			=
			(-)^{|\mu_m||\mu_n|}
			\alpha(\symgerst; \mu_\ell,\mu_n,\mu_m)
			~;
		\end{equation}
		\begin{equation}\label{Equation:SkewGerstAssociatorSymmetry}
			\alpha(\skewgerst; \mu_\ell,\mu_m,\mu_n) 
			=
			(-)^{(|\mu_n| + n - 1)(|\mu_m| + m - 1)}
			\alpha(\skewgerst; \mu_\ell,\mu_n,\mu_m)
			~.
		\end{equation}	
	\end{proposition}	
	\begin{proof}
		Plugging the result of corollary \ref{Cor:CyclicCommutationMultiBracket} inside equation \eqref{Eq:ExplicitAssociatorSym} yields
		\begin{align*}
			\alpha(\mu_\ell & ,\mu_m,\mu_n) =
			\\				
			\equal{}&
			\mu_\ell \circ \left[\left(\left( \mu_m \otimes \mu_n \right) \circ B_{m,n}\right)
			\otimes \mathbb{1}_{\ell-2}\right]
			\circ P_{m+n,\ell-2}
			=
			\\
			\equal{Cor. \ref{Cor:CyclicCommutationMultiBracket}}&
			(-)^{|\mu_m||\mu_n|}
			\mu_\ell \circ 
			\left[\left(\cycPermutator_{(2)} \circ \left(\mu_m \otimes \mu_n\right) \circ \cycPermutator_{(m+n)}^n 
			\circ B_{m,n}\right)
			\otimes \mathbb{1}_{\ell-2}\right]
			\circ P_{m+n,\ell-2}
			=
			\\
			\equal{}&
			(-)^{|\mu_m||\mu_n|}			
			\mu_\ell \circ \left[\left(\left( \mu_n \otimes \mu_m \right)\circ B_{n,m}\right)
			\otimes \mathbb{1}_{\ell-2}\right]
			\circ P_{m+n,\ell-2}
			=
			\\
			\equal{}&	
			(-)^{|\mu_m||\mu_n|}
			\alpha(\mu_\ell,\mu_n,\mu_m)
			~.
		\end{align*}
		Plugging the same corollary in equation \eqref{Eq:ExplicitAssociatorSkew} yields extra signs:
		\begin{align*}
			\alpha(\mu_\ell & ,\mu_m,\mu_n) =
			\\
			\equal{}&
			(-)^s ~
			\mu_\ell \circ \left[\left(\left( \mu_m \otimes \mu_n \right)\circ P_{m,n}\right)
			\otimes \mathbb{1}_{\ell-2}\right]
			\circ P_{m+n,\ell-2}
			=
			\\
			\equal{Cor. \ref{Cor:CyclicCommutationMultiBracket}}&
			(-)^{s + t} ~
			\mu_\ell \circ
			\left[\left(\tilde\cycPermutator_{(2)} \circ( \mu_n \otimes \mu_m) \circ \tilde\cycPermutator_{(m+n)}^n\circ P_{m,n}\right)
			\otimes \mathbb{1}_{\ell-2}\right]
			\circ P_{m+n,\ell-2}
			=
			\\
			\equal{}&
			(-)^{s + t} ~		
			\mu_\ell \circ \left[\left(\left( \mu_n \otimes \mu_m \right)\circ P_{n,m}\right)
			\otimes \mathbb{1}_{\ell-2}\right]
			\circ P_{m+n,\ell-2}
			=
			\\
			\equal{}&	
			(-)^{s + t + s'} ~
			\alpha(\mu_\ell,\mu_n,\mu_m)~,
		\end{align*}			
		where $s$ and $s'$ are the signs contained in the explicit expression of the associators (see proposition \ref{Prop:ExplicitAssociators}):
		\begin{align*}
			s=&
			{{|\mu_n|(m+\ell)} +{|\mu_m|(\ell-1)} +m(n+1)} ~,
			\\
			s'=&
			{{|\mu_m|(n+\ell)} +{|\mu_n|(\ell-1)} +n(m+1)} ~,
		\end{align*}
		and $t$ is the sign coming from corollary \ref{Cor:CyclicCommutationMultiBracket}:
		\begin{displaymath}
			t = 	mn+1 +|\mu_m||\mu_n|
			~.
		\end{displaymath}
		Computing $s+s'+t \mod 2 $ gives the exponent contained in equation 
		 \eqref{Equation:SkewGerstAssociatorSymmetry}.
	\end{proof}

\ifstandalone
	\bibliographystyle{../../hep} 
	\bibliography{../../mypapers,../../websites,../../biblio-tidy}
\fi

\cleardoublepage


%% file: chapters/gradedprelie/gradedprelie.tex
\chapter{Graded pre-Lie algebras}\label{App:PreLie}
This appendix contains some definitions and basic properties of other algebraic structures mentioned in this work.
In particular, we talk about \emph{differential graded Lie algebras} (DGLA) and \emph{pre-Lie algebras}.
We do not aim to give a comprehensive account; the general idea is to lay the foundations for a later reference in the thesis's body.

\section{Differential graded Lie algebras}
	$L_\infty$-algebras are a generalization of \emph{differential graded Lie algebra} (DGLA). Let us record some basic definitions:
	%
	\begin{definition}[Lie algebra]
		A \emph{Lie algebra} is a vector space $\mathfrak{g}\in \Vect$ equipped with a bilinear skew-symmetric map  $[\cdot,\cdot]:\mathfrak{g}\wedge\mathfrak{g}\to \mathfrak{g}$ which satisfies the Jacobi identity:
		\begin{displaymath}
			[x,[y,z]]+[z,[x,y]]+[y,[z,x]]=0 \qquad \forall x,y,z\in \mathfrak{g}~.
		\end{displaymath}
	\end{definition}
	The whole definition could be subsumed by the following diagram\footnote{The notion of Lie algebra may be formulated internally to any symmetric monoidal linear category. We do not insist here in this direction, see \cite{nlab:lie_algebra} for the general idea.}
	 in the category of vector spaces:
	\begin{equation}\label{eq:LieAlgebraDiagram}
		\begin{tikzcd}[column sep=7em]
			& \wedge^3 \mathfrak{g} \ar[ddl,"{[\cdot,\cdot]\ca [\cdot,\cdot]}",sloped]
			 \ar[d,"{\left([\cdot,\cdot]\otimes \Unit\right) \circ P_{2,1}}"]
			\\
			& \wedge^2 \mathfrak{g}
				\ar[d,"{[\cdot,\cdot]}"]
			\\
			0 \ar[r,hook]
			& \mathfrak{g}
		\end{tikzcd}
	\end{equation}
	Reading the above diagram in the category $\GVect$ of graded vector spaces one gets the notion of \emph{graded Lie algebra}:
	\begin{definition}[Graded Lie algebra (GLA)]
		A \emph{graded Lie algebra} is a graded vector space $\mathfrak{L}\in \GVect$ equipped with a bilinear graded skew-symmetric degree $0$ homogeneous map  $[\cdot,\cdot]:\mathfrak{L}\wedge\mathfrak{L}\to \mathfrak{L}$ which satisfies the (graded) Jacobi identity:
		\begin{equation}\label{eq:GradedJacobiDGLA}
			(-)^{|x||z|}[x,[y,z]]+(-)^{|z||y|}[z,[x,y]]+(-)^{|y||x|}[y,[z,x]]=0 \qquad \forall x,y,z\in \mathfrak{L}~.
		\end{equation}	
	\end{definition}
	\begin{remark}[Any associative algebra determines a Lie algebra]\label{rem:anyassociativegiveLie}
		Given a graded associative algebra $(A,\bullet)$. 
		Consider the graded commutator  $[\cdot,\cdot]_\bullet$ defined on homogeneous elements as
		\begin{displaymath}
			[x,y]_{\bullet}= x\bullet y - (-)^{|x||y|} y \bullet x	
		\end{displaymath}				
		(see definition \ref{def:gradedCommutator} below).
		The pair $(A,[\cdot,\cdot]_\bullet)$ is a GLA.
		\\
		Observe furthermore that the commutator of an associative graded algebra acts as a graded derivation (see definition \ref{Def:Fderivation}). 
		Namely, for any fixed homogeneous elements $x,y,z\in A$,  the following equation holds
		\begin{displaymath}
			[x,y\bullet z]_\bullet = [x,y]_\bullet \bullet z + (-)^{|x||y|} y \bullet [x,z]_\bullet
			~.
		\end{displaymath}
		In other terms, for any given $x\in A$, the unary operator $[x,\cdot]_\bullet$ is a degree $|x|$ derivation from the graded algebra $A$ into itself.
	\end{remark}

	If one reads diagram \eqref{eq:LieAlgebraDiagram} in the category of graded differential vector spaces, \ie cochain complexes (see section \ref{sec:HomologicalAlgebrasConventions}), one gets the notion of \emph{differential graded Lie algebra}:
	\begin{definition}[Differential Graded Lie Algebra (DGLA)]\label{def:DGLA}
		A \emph{differential graded Lie algebra} (DGLA) is a graded vector space $\mathfrak{L}\in\GVect$ equipped with a bilinear graded skew-symmetric degree $0$ homogeneous map  $[\cdot,\cdot]:\mathfrak{L}\wedge\mathfrak{L}\to \mathfrak{L}$	 and with a unary degree $1$ homogeneous map $d$ such that
		\begin{itemize}
			\item $d$ is a coboundary, \ie $d\circ d=0$;
			\item $[\cdot,\cdot]$ satisfies the graded Jacobi equation \eqref{eq:GradedJacobiDGLA};
			\item $\d$ and $[\cdot,\cdot]$ are compatible in the sense of the Liebniz rule\footnote{$d$ is a degree $1$ derivation from the non-associative anticommutative graded algebra $(A,[\cdot,\cdot])$ into itself.}
			\begin{displaymath}
				d[x,y]= [d x,y] + (-1)^{|x|}[x, d y]
				~.		
			\end{displaymath}
		\end{itemize}				
	\end{definition}
	The three conditions defining a DGLA are subsumed by the following commutative diagram in the category of graded vector spaces
	\begin{equation}\label{eq:DGLADiagram}
		\begin{tikzcd}[column sep=7em, row sep = large]
			& \wedge^3 \mathfrak{g} \ar[ddl,"{[\cdot,\cdot]\ca [\cdot,\cdot]}",sloped]
			 \ar[d,"{\left([\cdot,\cdot]\otimes \Unit\right) \circ P_{2,1}}"] & &
			\\
			& \wedge^2 \mathfrak{g}
				\ar[d,"{[\cdot,\cdot]}"]
				\ar[r,"d\otimes\Unit + \Unit\otimes d"] 
			& \wedge^2 \mathfrak{g}	[1] \ar[d,"{[\cdot,\cdot]}_{[1]}"]
			\\
			0 \ar[r,hook]
			& \mathfrak{g} \ar[r,"d"] \ar[rr,bend right=20,"0"']& \mathfrak{g}[1] \ar[r,"{d[1]}"] 
			& \mathfrak{g}[2]
		\end{tikzcd}
	\end{equation}

	\begin{remark}[Any GLA determines (several) DGLAs]\label{rem:DglaFromGla}
		Consider a GLA $(\mathfrak{L},[\cdot,\cdot])$.
		Let be $x\in \mathfrak{L}^1$ a degree $1$ homogeneous element that commutes with itself:
		\begin{displaymath}
			[x,x]=0 
		\end{displaymath}
		(notice that we are in the graded setting, hence the above equation is not automatic on elements in odd degree).
		Then, the graded Jacobi identity implies that the degree $1$ homogeneous map $[x,\cdot]$  is a derivation on $(\mathfrak{L},[\cdot,\cdot])$ and is $2$-nilpotent.
		In other terms, $[x,\cdot]$ is a coboundary operator which satisfies the Liebniz equation, hence $(\mathfrak{L},[\cdot,\cdot],[x,\cdot])$ is a DGLA.
		\\
		This idea is one of the fundamental concepts underlying the so-called \emph{deformation theory} (see \cite{Doubek2007} for some introductory lecture notes).
	\end{remark}
	
	The special elements introduced in remark \ref{rem:DglaFromGla} are subcases of a more general class:
	\begin{definition}[Maurer-Cartan elements]\label{def:MCelements}
		Let be $\mathfrak{L}=(\mathfrak{L},d,[\cdot,\cdot])$ a DGLA. We call \emph{Maurer-Cartan elements} of $\mathfrak{L}$ the set of degree $1$ elements satisfying the Maurer-Cartan equation:
		\begin{displaymath}
			MC(\mathfrak{L}):=
			\left\{
				x \in \mathfrak{L}^1 ~\left\vert~ 
					d x + \frac{1}{2}[x,x] =0
			\right\}\right.
		\end{displaymath}
	\end{definition}
	\begin{remark}[Maurer-Cartan elements in a $L_\infty$-algebra]
		Notice that the Maurer-Cartan equation can be formally extended to any (curved) $L_\infty$-algebra  $(L,\{\mu_k\}_{k\geq 0})$ as follows:
		\begin{displaymath}
			\sum_{k=0}^\infty \dfrac{1}{k!}\mu_k(x,\dots x) = 0
			~.
		\end{displaymath}
		In this case, however, it is essential to pay extra attention to the convergence conditions of the infinite series.
	\end{remark}

\section{Graded pre-Lie algebras}
Given any associative algebra, it is a standard result (see remark \ref{rem:anyassociativegiveLie}) that the corresponding commutator yields a Lie algebra structure.
Roughly speaking, a pre-Lie algebra is a non-associative algebra determining a Lie structure in a similar fashion. 
One of the early adopters of this concept was Gerstenhaber (see \cite{Gerstenhaber1963a}) see also \cite{Manchon2011a} for a survey.

Consider a non-associative (non-commutative) graded $\mathbb{R}$-algebra $(X,\MBComp)$, we introduce the following auxiliary operators:
\begin{definition}[Associator]\label{def:gradedAssociator}
	Given a (not necessarily associative) graded algebra $(X,\ca)$, we call \emph{associator} the graded tri-linear graded morphism
	\begin{displaymath}
		\alpha(\ca\,;\cdot,\cdot,\cdot): X^{\otimes 3} \to X
	\end{displaymath}
	defined on arbitrary homogeneous elements $x,y,z\in X$ as 
	\begin{displaymath}
		\alpha(\ca\,;x,y,z) := (x\ca y)\ca z - x \ca (y \ca z)~.
	\end{displaymath}	
\end{definition}
One can construct a graded skew-symmetric binary bracket out of any graded algebra:
\begin{definition}[(Graded) Commutator]\label{def:gradedCommutator}
	Given a (not necessarily associative) graded algebra $(X,\ca)$, 
	we call \emph{commutator} the graded skew-symmetric bi-linear graded morphism
	\begin{displaymath}
		[\cdot,\cdot]_{\ca}: X^{\wedge 2} \to X
	\end{displaymath}
	defined on arbitrary homogeneous elements $x,y\in X$ as	
	\begin{displaymath}
		[x,y]_{\ca} := (x\ca y)- (-)^{|x||y|} (y \ca x)~.
	\end{displaymath}	
\end{definition}

The role of the associator is to measure the failing of the associativity of the product $\ca$. When $(X,\ca)$ is associative $\alpha$ is automatically zero.
In the case that $\alpha$ satisfy certain symmetry properties one talks about \emph{pre-Lie algebras}:

\sidebyside{
\begin{definition}[Left pre-Lie Algebra]
	A non-associative graded algebra $(X,\ca)$ is said to be a \emph{left pre-Lie algebra} if the corresponding associator is graded symmetric in the two leftmost entries;
	\\ 
	\ie, for any  $x,y,z\in X$, one has:
	\begin{displaymath}
		\alpha(\ca;x,y,z) = (-)^{|x||y|}\alpha(\ca;y,x,z)~.
	\end{displaymath}	
\end{definition}

}{
\begin{definition}[Right pre-Lie Algebra]
	A non-associative graded algebra $(X,\ca)$ is said to be a \emph{right pre-Lie algebra} if the corresponding associator is graded symmetric in the two rightmost entries;  
	\\
	\ie, for any  $x,y,z\in X$, one has:
	\begin{displaymath}
		\alpha(\ca;x,y,z) = (-)^{|y||z|}\alpha(\ca;x,z,y)
		~.
	\end{displaymath}	
\end{definition}
}

\smallskip
\begin{notation}[Short-hand notation]
In the following we will focus on right pre-Lie algebras and we will omit the symbol $\ca$ when it is possible. 
The action of the product $\ca$ will be denoted by juxtaposition. 
\end{notation}

The next proposition justifies why these algebras are called \emph{pre-}Lie:
\begin{proposition}
	Given a right pre-Lie $(X,\ca)$, then the corresponding graded commutator satisfies the 
	\emph{(right) graded Jacobi} equation
	\begin{displaymath}
	 J(x,y,z) := (-)^ {|x||z|}[x,[y,z]] + (-)^ {|x||y|}[y,[z,x]] + (-)^ {|y||z|}[z,[x,y]]
	 =0~.
	\end{displaymath}
\end{proposition}
\begin{proof}
	Expanding the definition of commutator one gets that:
	\begin{displaymath}
		\begin{aligned}
		[x,&[y,z]] =\\
		=&~x(yz) - (-)^{|y||z|} x(zy) - (-)^{(|y|+|z|)|x|}(yz)x + (-)^{(|y|+|z|)|x|}(-)^{|y||z|}(zy)x
		~;
		\end{aligned}
	\end{displaymath}
	therefore:
	\begin{align*}
		&(-)^{|x||z|}[x,[y,z]]
		\\ 
		&=~ (-)^{|x||z|}x(yz) - (-)^{(|x|+|y|)|z|} x(zy) 
		 - (-)^{|y||x|}(yz)x + (-)^{|y|(|z|+|x|)}(zy)x	
		 ~;
		\\[1em]
		&(-)^{|y||x|}[y,[z,x]]
		\\
		&=~ (-)^{|y||x|}y(zx) - (-)^{(|y|+|z|)|x|} y(xz) 
		 - (-)^{|z||y|}(zx)y + (-)^{|z|(|x|+|y|)}(xz)y	
		~;
		\\[1em]
		&(-)^{|z||y|}[z,[x,y]]
		\\ 
		&=~ (-)^{|z||y|}z(xy) - (-)^{(|z|+|x|)|y|} z(yx) 
		 - (-)^{|x||z|}(xy)z + (-)^{|x|(|y|+|z|)}(yx)z		
		 ~.
	\end{align*}
	Summing up all the previous terms, one gets
	\begin{align*}
		J(x,y,z)
		=&+ (-)^{|x||z|} ( -\alpha(x,y,z) + (-)^{|z||y|} \alpha(x,z,y))+
		\\
		 &+ (-)^{|y||x|} ( -\alpha(y,z,x) + (-)^{|x||z|} \alpha(y,x,z))+
		\\
		 &+ (-)^{|z||y|} ( -\alpha(z,x,y) + (-)^{|y||x|} \alpha(z,y,z))=
		\\
		=& 0~.
	\end{align*}
\end{proof}

In a pre-Lie algebra the associator measure the failure of the latter property:
%
%
%
\begin{proposition}
	Let be $(X,\ca)$ right pre-Lie. Consider three homogeneous elements $x,y,z\in X$. 
	Then
	\begin{displaymath}
	\begin{aligned}
		[x,yz] =& [x,y]z + (-)^{|x||y|}y[x,z] + \alpha(x,y,z)
		~;
		\\
		[yz,x] =& y[z,x] + (-)^{|x||z|}[y,x]z - (-)^{|x|(|y|+|z|)}\alpha(x,y,z)		
		~.
	\end{aligned}
	\end{displaymath}

\end{proposition}
\begin{proof}
	The first claim follows by the simple expansion of the terms:
	\begin{displaymath}
	\begin{aligned}
		[x,yz] =& x(yz) -(-)^{|x|(|y|+|z|)}(yz)x 
		~;
		\\[1em]
		y[x,z] =& y(xz) -(-)^{|x||z|} y(zx) = 
		\\
		=&~ y(xz) -(-)^{|x||z|} (yz)x + (-)^{|x||z|} \alpha(y,z,x)
		~;
		\\[1em]
		[x,y]z =& (xy)z -(-)^{|x||y|}(yx)z = 
		\\
		=&~ x(yz) + \alpha(x,y,z) 
		-(-)^{|x||y|}y(xz) -(-)^{|x||y|}\alpha(y,x,z)
		~.
	\end{aligned}
	\end{displaymath}
	The second claim follows from the first and from the graded commutativity of the commutator:
	\begin{displaymath}
	\begin{aligned}
		[yz,x] =& -(-)^{|x|(|y|+|z|)}[x,yz] 
		=
		\\
		=&  -(-)^{|x|(|y|+|z|)}([x,y]z + (-)^{|x||y|}y[x,z] + \alpha(x,y,z))
		=
		\\
		=& (-)^{|x||z|}[y,x]z + y[z,x] - (-)^{|x|(|y|+|z|)}\alpha(x,y,z)
		~.
	\end{aligned}
	\end{displaymath}

\end{proof}

\ifstandalone
	\bibliographystyle{../../hep} 
	\bibliography{../../mypapers,../../websites,../../biblio-tidy}
\fi

\cleardoublepage


%% file: thesis-arxiv.bbl
\newcommand{\etalchar}[1]{$^{#1}$}
\begin{thebibliography}{CMDNY13}

\bibitem[AB84]{MR721448}
M.~F. Atiyah and R.~Bott, \textsl{ The moment map and equivariant cohomology},
\newblock Topology \textbf{ 23}(1), 1--28 (1984).

\bibitem[AK98]{Arn-Khe}
V.~I. Arnold and B.~A. Khesin,
\newblock \textsl{ Topological methods in hydrodynamics}, volume 125,
\newblock New York, NY: Springer, 1998.

\bibitem[All10]{Allocca2010}
M.~P. Allocca,
\newblock \textsl{ L $\infty$ Algebra Representation Theory},
\newblock PhD thesis, North Carolina State University, 2010.

\bibitem[AMM78]{Abraham1978}
R.~Abraham, J.~E. Marsden and J.~E. Marsden,
\newblock \textsl{ Foundations of mechanics}, volume~36,
\newblock Benjamin/Cummings Publishing Company Reading, Massachusetts, ii
  edition, 1978.

\bibitem[Arn66]{Arnold66}
V.~Arnold,
\newblock Sur la g{\'e}om{\'e}trie diff{\'e}rentielle des groupes de Lie de
  dimension infinie et ses applications {\`a} l'hydrodynamique des fluides
  parfaits,
\newblock in \textsl{ Annales de l'institut Fourier}, volume~16, pages
  319--361, 1966.

\bibitem[Ati57]{Atiyah1957}
M.~F. Atiyah, \textsl{ Complex analytic connections in fibre bundles},
\newblock Transactions of the American Mathematical Society \textbf{ 85}(1),
  181--181 (January 1957).

\bibitem[Av80]{Kuznetsov-Mikhailov80}
K.~E.~M. Av, \textsl{ On The Topological Meaning Of Canonical Clebsch
  Variables},
\newblock Phys. Lett. Sect. A; Issn 0375-9601; Nld; Da. 1980; Vol. 77; No 1;
  Pp. 37-38; Bibl. 6 Ref.  (1980).

\bibitem[Awa92]{Awane1992}
A.~Awane, \textsl{ k ‐symplectic structures},
\newblock J. Math. Phys. \textbf{ 33}(12), 4046--4052 (December 1992).

\bibitem[Ban16]{Bandiera2016}
R.~Bandiera, \textsl{ Homotopy abelian $L_\infty$ algebras and splitting
  property},
\newblock Rend. di Mat. e delle Sue Appl. \textbf{ 37}(1-2), 105--122 (2016).

\bibitem[BC03]{Baez2003}
J.~C. Baez and A.~S. Crans, \textsl{ Higher-Dimensional Algebra VI: Lie
  2-Algebras},
\newblock Theory Appl. Categ. \textbf{ 12}(1), 492--538 (July 2003),
  {arXiv:0307263}.

\bibitem[Bes08]{MR2371700}
A.~L. Besse,
\newblock \textsl{ Einstein manifolds},
\newblock Classics in Mathematics, Springer-Verlag, Berlin, 2008.

\bibitem[BFLS98]{Barnich1998}
G.~Barnich, R.~Fulp, T.~Lada and J.~Stasheff, \textsl{ The sh {L}ie structure
  of {P}oisson brackets in field theory},
\newblock Communications in Mathematical Physics \textbf{ 191}(3), 585--601
  (February 1998), {9702176}.

\bibitem[BHR10]{Baez2010}
J.~C. Baez, A.~E. Hoffnung and C.~L. Rogers, \textsl{ Categorified Symplectic
  Geometry and the Classical String},
\newblock Communications in Mathematical Physics \textbf{ 293}(3), 701--725
  (February 2010), {arXiv:0808.0246}.

\bibitem[Bla19]{Blacker2019}
C.~Blacker, \textsl{ Polysymplectic reduction and the moduli space of flat
  connections},
\newblock J. Phys. A Math. Theor. \textbf{ 52}(33), 335201 (July 2019),
  {arXiv:1810.04924}.

\bibitem[BM13]{Buijs2013}
U.~Buijs and A.~Murillo, \textsl{ Algebraic models of non-connected spaces and
  homotopy theory of $L\infty$ algebras},
\newblock Adv. Math. (N. Y). \textbf{ 236}, 60--91 (2013).

\bibitem[Bor49]{MR0029915}
A.~Borel, \textsl{ Some remarks about {L}ie groups transitive on spheres and
  tori},
\newblock Bull. Amer. Math. Soc. \textbf{ 55}, 580--587 (1949).

\bibitem[Bor50]{MR0034768}
A.~Borel, \textsl{ Le plan projectif des octaves et les sph\`eres comme espaces
  homog\`enes},
\newblock C. R. Acad. Sci. Paris \textbf{ 230}, 1378--1380 (1950).

\bibitem[Bor55]{Borel1955}
A.~Borel, \textsl{ Topology of Lie groups and characteristic classes},
\newblock Bulletin of the American Mathematical Society \textbf{ 61}(5),
  397--433 (September 1955).

\bibitem[Bor60]{Borel1960}
A.~Borel,
\newblock \textsl{ Seminar on transformation groups. With contributions by G.
  Bredon, E. E. Floyd, D. Montgomery and R. Palais}, volume~46,
\newblock Princeton University Press, Princeton, NJ, 1960.

\bibitem[Bou17]{Bourles2017}
H.~Bourl{\`{e}}s,
\newblock \textsl{ Fundamentals of advanced mathematics 1. Categories,
  algebraic structures, linear and homological algebra},
\newblock Amsterdam: Elsevier/ISTE Press, 2017.

\bibitem[BPS09]{Bursztyn2008}
H.~Bursztyn, D.~I. Ponte and P.~Severa, \textsl{ Courant morphisms and moment
  maps},
\newblock Mathematical Research Letters \textbf{ 16}(2), 215--232 (January
  2009), {0801.1663}.

\bibitem[Bry93]{Bry}
J.-L. Brylinski,
\newblock \textsl{ Loop Spaces, Characteristic Classes and Geometric
  Quantization},
\newblock Birkh{\"{a}}user Boston, Boston, MA, reprint of edition, 1993.

\bibitem[Bry06]{MR2253159}
R.~L. Bryant, \textsl{ On the geometry of almost complex 6-manifolds},
\newblock Asian J. Math. \textbf{ 10}(3), 561--605 (2006).

\bibitem[BS06]{BeSpe06}
A.~Besana and M.~Spera, \textsl{ On some symplectic aspects of knot framings},
\newblock Journal of Knot Theory and Its Ramifications \textbf{ 15}(07),
  883--912 (2006).

\bibitem[BS11]{Bi2011a}
Y.~Bi and Y.~Sheng, \textsl{ On higher analogues of Courant algebroids},
\newblock Sci. China Math. \textbf{ 54}(3), 437--447 (March 2011),
  {arXiv:1003.1350v2}.

\bibitem[BT13]{Bott-Tu82}
R.~Bott and L.~W. Tu,
\newblock \textsl{ Differential forms in algebraic topology}, volume~82 of
  \textsl{ Graduate Texts in Mathematics},
\newblock Springer Science \& Business Media, 2013.

\bibitem[Car18]{Carosso2018}
A.~Carosso, \textsl{ Geometric Quantization},
\newblock (January 2018), {arXiv:1801.02307}.

\bibitem[CCI91]{Carinena1991b}
J.~F. Cari{\~{n}}ena, M.~Crampin and L.~A. Ibort, \textsl{ On the
  multisymplectic formalism for first order field theories},
\newblock Differ. Geom. its Appl. \textbf{ 1}(4), 345--374 (1991).

\bibitem[CdS01]{CannasdaSilva2001}
A.~Cannas~da Silva,
\newblock \textsl{ Lectures on symplectic geometry}, volume 1764 of \textsl{
  Lecture Notes in Mathematics},
\newblock Springer-Verlag, Berlin, 2001.

\bibitem[CFRZ16]{Callies2016}
M.~Callies, Y.~Fr{\'{e}}gier, C.~L. Rogers and M.~Zambon, \textsl{ Homotopy
  moment maps},
\newblock Adv. Math. (N. Y). \textbf{ 303}, 954--1043 (November 2016),
  {arXiv:1304.2051}.

\bibitem[{Che}77]{Chen}
K.-T. {Chen}, \textsl{ Iterated path integrals},
\newblock Bull. Am. Math. Soc. \textbf{ 83}, 831--879 (1977).

\bibitem[CID99]{CatIbort}
F.~Cantrijn, A.~Ibort and M.~{De Le{\'{o}}n}, \textsl{ On the geometry of
  multisymplectic manifolds},
\newblock J. Aust. Math. Soc. \textbf{ 66}(3), 303--330 (June 1999).

\bibitem[CMDNY13]{Cappelletti2013}
B.~Cappelletti-Montano, A.~De~Nicola and I.~Yudin, \textsl{ A Survey On
  Cosymplectic Geometry},
\newblock Rev. Math. Phys. \textbf{ 25}(10), 1343002 (November 2013),
  {arXiv:1305.3704}.

\bibitem[Cou90]{Courant1990}
T.~J. Courant, \textsl{ Dirac manifolds},
\newblock Trans. Am. Math. Soc. \textbf{ 319}(2), 631--661 (February 1990).

\bibitem[Crn88]{Crnkovic}
C.~Crnkovic, \textsl{ Symplectic geometry of the convariant phase space},
\newblock Classical and Quantum Gravity \textbf{ 5}(12), 1557 (1988).

\bibitem[CS11]{CATTANEO2011}
A.~S. Cattaneo and F.~Sch{\"{a}}tz, \textsl{ Introduction To Supergeometry},
\newblock Rev. Math. Phys. \textbf{ 23}(06), 669--690 (July 2011),
  {arXiv:1011.3401v2}.

\bibitem[CT01]{Chen4}
K.-T. Chen and P.~Tondeur,
\newblock \textsl{ Collected Papers of KT Chen},
\newblock Springer Science \& Business Media, 2001.

\bibitem[{de }35]{DeDonder35}
T.~{de Donder},
\newblock Th\'eorie invariantive du calcul des variations. Nouvelle edit,
\newblock Paris: Gauthier-Villars \& Cie. XXI, 230 p. (1935)., 1935.

\bibitem[Del87]{Deligne}
P.~Deligne,
\newblock A letter to Millson,
\newblock \url{https://publications.ias.edu/sites/default/files/millson.pdf}do,
  1987.

\bibitem[Del15]{Delgado2015}
N.~L. Delgado,
\newblock $A_\infty$- and $L_\infty$-algebras from the point of view of formal
  geometry,
\newblock in \textsl{ Supergeometry, Operads and L$\infty$-algebras}, Rabat,
  2015.

\bibitem[Del18a]{Delgado2018}
N.~L. Delgado,
\newblock \textsl{ Lagrangian field theories: ind/pro-approach and L-infinity
  algebra of local observables},
\newblock PhD thesis, Rheinischen Friedrich-Wilhelms-Universit{\"{a}}t Bonn,
  2018.

\bibitem[Del18b]{Delgado2018b}
N.~L. Delgado, \textsl{ Multisymplectic Structures and Higher Momentum Maps},
\newblock (November 2018), {arXiv:1811.01415}.

\bibitem[DK00]{Duistermaat2000}
J.~J. Duistermaat and J.~A.~C. Kolk,
\newblock \textsl{ Lie Groups},
\newblock Universitext, Springer Berlin Heidelberg, Berlin, Heidelberg, 2000.

\bibitem[DMS88]{DeLeon1988}
M.~{De Le{\'{o}}n}, I.~M{\'{e}}ndez and M.~Salgado, \textsl{ p-Almost tangent
  structures},
\newblock Rend. del Circ. Mat. di Palermo \textbf{ 37}(2), 282--294 (May 1988).

\bibitem[DMZ07]{Doubek2007}
M.~Doubek, M.~Markl and P.~Zima, \textsl{ Deformation theory (lecture notes)},
\newblock Archivum Mathematicum \textbf{ 43}(5), 333--371 (05 2007),
  {arXiv:0705.3719v3}.

\bibitem[Dol07]{Dolgushev2007}
V.~A. Dolgushev, \textsl{ Erratum to: "A Proof of Tsygan's Formality Conjecture
  for an Arbitrary Smooth Manifold"},
\newblock (March 2007), {arXiv:0703113}.

\bibitem[DP16]{Dotsenko2016}
V.~Dotsenko and N.~Poncin, \textsl{ A Tale of Three Homotopies},
\newblock Appl. Categ. Struct. \textbf{ 24}(6), 845--873 (December 2016),
  {arXiv:1208.4695}.

\bibitem[dR55]{dR}
G.~de~Rham,
\newblock \textsl{ Vari{\'e}t{\'e}s diff{\'e}rentiables},
\newblock Publications de l'Institut de Math\'ematique de l'Universit\'e de
  Nancago III. Actualit\'es scientifiques et industrielles, Paris: Hermann.,
  1955.

\bibitem[Dro73]{Dror1973}
E.~Dror, \textsl{ Homology spheres},
\newblock Israel Journal of Mathematics \textbf{ 15}(2), 115--129 (1973).

\bibitem[DS18]{Deser2018b}
A.~Deser and C.~S{\"{a}}mann, \textsl{ Extended Riemannian Geometry I: Local
  Double Field Theory},
\newblock Ann. Henri Poincar{\'{e}} \textbf{ 19}(8), 2297--2346 (August 2018),
  {arXiv:1711.03308}.

\bibitem[EM70]{Eb-Mar}
D.~G. Ebin and J.~Marsden, \textsl{ Groups of diffeomorphisms and the motion of
  an incompressible fluid},
\newblock Annals of Mathematics , 102--163 (1970).

\bibitem[EPS15]{Enciso2015}
A.~Enciso and D.~Peralta-Salas, \textsl{ Existence of knotted vortex tubes in
  steady Euler flows},
\newblock Acta Math. \textbf{ 214}(1), 61--134 (October 2015).

\bibitem[FF83]{Fenn}
R.~Fenn and R.~A. Fenn,
\newblock \textsl{ Techniques of geometric topology}, volume~57,
\newblock CUP Archive, 1983.

\bibitem[FH13]{fulton2013representation}
W.~Fulton and J.~Harris,
\newblock \textsl{ Representation theory: a first course}, volume 129,
\newblock Springer Science \& Business Media, 2013.

\bibitem[FLGZ15]{Fregier2015}
Y.~Fr{\'{e}}gier, C.~Laurent-Gengoux and M.~Zambon, \textsl{ A cohomological
  framework for homotopy moment maps},
\newblock J. Geom. Phys. \textbf{ 97}, 119--132 (September 2015),
  {arXiv:1409.3142}.

\bibitem[FM07]{Fiorenza2006}
D.~Fiorenza and M.~Manetti, \textsl{ $L_\infty$ structures on mapping cones},
\newblock Algebr. Number Theory \textbf{ 1}(3), 301--330 (January 2007),
  {arXiv:0601312}.

\bibitem[FPRR03]{Forger2003}
M.~Forger, C.~Paufler, H.~Roemer and H.~R{\"{o}}mer, \textsl{ The Poisson
  bracket for Poisson forms in multisymplectic field theory},
\newblock Rev. Math. Phys. \textbf{ 15}(7), 705--743 (September 2003),
  {arXiv:0202043v1}.

\bibitem[FR05]{Forger2005}
M.~Forger and S.~V. Romero, \textsl{ Covariant Poisson Brackets in Geometric
  Field Theory},
\newblock Communications in Mathematical Physics \textbf{ 256}(2), 375--410
  (June 2005), {arXiv:0408008}.

\bibitem[FRS14]{Fiorenza2014a}
D.~Fiorenza, C.~L. Rogers and U.~Schreiber, \textsl{ $L_\infty$-algebras of
  local observables from higher prequantum bundles},
\newblock Homology Homotopy Appl. \textbf{ 16}(2), 107--142 (2014).

\bibitem[FSB95]{Freyd1995}
P.~Freyd, R.~Street and M.~Barr, \textsl{ On the size of categories},
\newblock Theory and Applications of Categories \textbf{ 1}(9), 18--25 (1995).

\bibitem[FYH{\etalchar{+}}85]{Freyd-etal}
P.~{Freyd}, D.~{Yetter}, J.~{Hoste}, W.~B.~R. {Lickorish}, K.~{Millett} and
  A.~{Ocneanu}, \textsl{ A new polynomial invariant of knots and links},
\newblock Bull. Am. Math. Soc., New Ser. \textbf{ 12}, 239--246 (1985).

\bibitem[FZ15]{Fregier2015b}
Y.~Fr{\'{e}}gier and M.~Zambon, \textsl{ Simultaneous deformations of algebras
  and morphisms via derived brackets},
\newblock J. Pure Appl. Algebr. \textbf{ 219}(12), 5344--5362 (December 2015),
  {arXiv:1301.4864}.

\bibitem[Ger63]{Gerstenhaber1963a}
M.~Gerstenhaber, \textsl{ The Cohomology Structure of an Associative Ring},
\newblock Annals of Mathematics \textbf{ 78}(2), 267 (September 1963).

\bibitem[Ger64]{Gerstenhaber1964}
M.~Gerstenhaber, \textsl{ On the Deformation of Rings and Algebras},
\newblock Annals of Mathematics \textbf{ 79}(1), 59--103 (1964).

\bibitem[Get10]{Getzler1991}
E.~Getzler, \textsl{ Higher derived brackets},
\newblock (October 2010), {arXiv:1010.5859}.

\bibitem[GHV73]{MR0336651}
W.~Greub, S.~Halperin and R.~Vanstone,
\newblock \textsl{ Connections, curvature, and cohomology. {V}ol. {II}: {L}ie
  groups, principal bundles, and characteristic classes},
\newblock Academic Press [A subsidiary of Harcourt Brace Jovanovich,
  Publishers], New York-London, 1973.

\bibitem[GIM{\etalchar{+}}98]{Gimmsy1}
M.~Gotay, J.~Isenberg, J.~Marsden, R.~Montgomery, J.~{\'S}niatycki and
  P.~Yasskin, \textsl{ Momentum Maps and Classical Fields. Part I: Covariant
  Field Theory (1998)},
\newblock arXiv preprint physics/9801019 (August 2004), 68 (January 1998),
  {9801019}.

\bibitem[Gol71]{Goldin71}
G.~A. Goldin, \textsl{ Nonrelativistic current algebras as unitary
  representations of groups},
\newblock Journal of Mathematical Physics \textbf{ 12}(3), 462--487 (1971).

\bibitem[Gol88]{Goldin87}
G.~A. Goldin,
\newblock Parastatistics, $\theta$-Statistics, and Topological Quantum
  Mechanics from Unitary Representations of Diffeomorphism Groups,
\newblock in \textsl{ Proceedings of the XV International Conference on
  Differential Geometric Methods in Physics, ed. H.-D. Doebner and JD Hennig
  (World Scientific, Singapore, 1987)}, page 197, World Scientific, 1988.

\bibitem[Gol12]{Goldin12}
G.~A. Goldin,
\newblock Diffeomorphism groups and nonlinear quantum mechanics,
\newblock in \textsl{ Journal of Physics: Conference Series}, volume 343, page
  012006, IOP Publishing, 2012.

\bibitem[{Got}91]{Gotay91}
M.~J. {Gotay},
\newblock A multisymplectic framework for classical field theory and the
  calculus of variations. I: Covariant Hamiltonian formalism,
\newblock Mechanics, analysis and geometry: 200 years after Lagrange, 203-235
  (1991)., 1991.

\bibitem[GS77]{Gui-Ste}
V.~{Guillemin} and S.~{Sternberg},
\newblock \textsl{ Geometric asymptotics}, volume~14,
\newblock American Mathematical Society (AMS), Providence, RI, 1977.

\bibitem[GS84]{GS84}
V.~Guillemin and S.~Sternberg,
\newblock \textsl{ Symplectic techniques in physics},
\newblock Cambridge university press, 1984.

\bibitem[GS15]{Grutzmann2015}
M.~Gr{\"{u}}tzmann and T.~Strobl, \textsl{ General Yang–Mills type gauge
  theories for p-form gauge fields: From physics-based ideas to a mathematical
  framework or From Bianchi identities to twisted Courant algebroids},
\newblock Int. J. Geom. Methods Mod. Phys. \textbf{ 12}(01), 1550009 (January
  2015), {arXiv:1407.6759}.

\bibitem[GSB99]{Guillemin1999}
V.~W. Guillemin, S.~Sternberg and J.~Br{\"{u}}ning,
\newblock \textsl{ Supersymmetry and Equivariant de Rham Theory},
\newblock Springer Berlin Heidelberg, Berlin, Heidelberg, 1999.

\bibitem[{Gua}93]{Gua}
E.~{Guadagnini},
\newblock \textsl{ The link invariants of the Chern-Simons field theory. New
  developments in topological quantum field theory}, volume~10,
\newblock Berlin: Walter de Gruyter, 1993.

\bibitem[Gua03]{Gualtieri2004}
M.~Gualtieri,
\newblock \textsl{ Generalized complex geometry},
\newblock PhD thesis, University of Oxford, 2003.

\bibitem[Gua11]{Gu2}
M.~Gualtieri, \textsl{ Generalized complex geometry},
\newblock Ann. of Math. (2) \textbf{ 174}(1), 75--123 (2011).

\bibitem[G{\"{u}}n87]{Gunther1987a}
C.~G{\"{u}}nther, \textsl{ The polysymplectic Hamiltonian formalism in field
  theory and calculus of variations. I. The local case},
\newblock J. Differ. Geom. \textbf{ 25}(1), 23--53 (1987).

\bibitem[H\"71]{Hor}
L.~H\"ormander, \textsl{ Fourier integral operators. I},
\newblock Acta Math. \textbf{ 127}, 79--183 (1971).

\bibitem[Hai87]{Hain}
R.~M. Hain,
\newblock The geometry of the mixed Hodge structure on the fundamental group,
\newblock in \textsl{ Proc. Symp. Pure Math}, volume~46, pages 247--282, 1987.

\bibitem[Hat02]{MR1867354}
A.~Hatcher,
\newblock \textsl{ Algebraic topology},
\newblock Cambridge University Press, Cambridge, 2002.

\bibitem[H{\'{e}}l11]{Helein2011c}
F.~H{\'{e}}lein,
\newblock Multisymplectic formalism and the covariant phase space,
\newblock in \textsl{ Var. Probl. Differ. Geom.}, edited by R.~Bielawski,
  K.~Houston and J.~Speight, pages 94--126, Cambridge University Press,
  Cambridge, June 2011.

\bibitem[Her18a]{Herman2018}
J.~Herman, \textsl{ Existence and uniqueness of weak homotopy moment maps},
\newblock Journal of Geometry and Physics \textbf{ 131}, 52--65 (September
  2018).

\bibitem[Her18b]{Herman2017}
J.~Herman, \textsl{ Noether's theorem in multisymplectic geometry},
\newblock Differ. Geom. its Appl. \textbf{ 56}, 260--294 (February 2018),
  {arXiv:1705.05818}.

\bibitem[Hsi75]{Hsiang1975}
W.~Y. Hsiang,
\newblock \textsl{ Cohomology theory of topological transformation groups},
  volume~85,
\newblock Springer-Verlag, Berlin, 1975.

\bibitem[HT12]{Hebda-Tsau12}
J.~Hebda and C.~Tsau, \textsl{ An approach to higher order linking invariants
  through holonomy and curvature},
\newblock Transactions of the American Mathematical Society \textbf{ 364}(8),
  4283--4301 (2012).

\bibitem[Ibo01]{MR1939543}
A.~Ibort,
\newblock Multisymplectic geometry: generic and exceptional,
\newblock in \textsl{ Proceedings of the {IX} {F}all {W}orkshop on {G}eometry
  and {P}hysics ({V}ilanova i la {G}eltr\'{u}, 2000)}, volume~3 of \textsl{
  Publ. R. Soc. Mat. Esp.}, pages 79--88, R. Soc. Mat. Esp., Madrid, 2001.

\bibitem[JRSW19]{Jurco2018}
B.~Jur{\v{c}}o, L.~Raspollini, C.~S{\"{a}}mann and M.~Wolf, \textsl{
  $L_\infty$‐Algebras of Classical Field Theories and the
  Batalin–Vilkovisky Formalism},
\newblock Fortschritte der Phys. \textbf{ 67}(7), 1900025 (July 2019),
  {arXiv:1809.09899}.

\bibitem[{Kau}01]{Kauffman}
L.~H. {Kauffman},
\newblock \textsl{ Knots and physics. 3rd ed}, volume~1,
\newblock Singapore: World Scientific, 3rd ed. edition, 2001.

\bibitem[Khe13]{Khe}
B.~Khesin, \textsl{ The vortex filament equation in any dimension},
\newblock Procedia IUTAM \textbf{ 7}, 135--140 (2013).

\bibitem[KI13]{Kleckner2013}
D.~Kleckner and W.~T.~M. Irvine, \textsl{ Creation and dynamics of knotted
  vortices},
\newblock Nat. Phys. \textbf{ 9}(4), 253--258 (April 2013).

\bibitem[Kij73]{Kijowski1973}
J.~Kijowski, \textsl{ A Finite-dimensional Canonical Formalism in the Classical
  Field Theory},
\newblock  \textbf{ 30}, 99--128 (1973).

\bibitem[Kir01]{Kirillov01}
A.~A. Kirillov,
\newblock Geometric quantization,
\newblock in \textsl{ Dynamical systems IV}, pages 139--176, Springer, 2001.

\bibitem[KM97]{Kri-Mich}
A.~{Kriegl} and P.~W. {Michor},
\newblock \textsl{ The convenient setting of global analysis}, volume~53,
\newblock Providence, RI: American Mathematical Society, 1997.

\bibitem[KMM20]{Khesin2020a}
B.~Khesin, G.~Misiolek and K.~Modin, \textsl{ Geometric Hydrodynamics of
  Compressible Fluids},
\newblock (January 2020), {arXiv:2001.01143}.

\bibitem[KN96]{Kobayashi1996}
S.~Kobayashi and K.~Nomizu,
\newblock \textsl{ Foundations of Differential Geometry Volume I}, volume~1,
\newblock Wiley-Interscience, 1996.

\bibitem[Kna13]{knapp2013lie}
A.~W. Knapp,
\newblock \textsl{ Lie groups beyond an introduction}, volume 140,
\newblock Springer Science \& Business Media, 2013.

\bibitem[Kon03]{Kontsevich2003}
M.~Kontsevich, \textsl{ Deformation Quantization of Poisson Manifolds},
\newblock Lett. Math. Phys. \textbf{ 66}(3), 157--216 (December 2003),
  {9709040}.

\bibitem[Kos70]{Kostant70}
B.~Kostant,
\newblock Quantization and unitary representations. {I}. {P}requantization,
\newblock in \textsl{ Lectures in modern analysis and applications, III}, pages
  87--208. Lecture Notes in Math., Vol. 170, Springer, Berlin, 1970.

\bibitem[KS76]{KS}
J.~Kijowski and W.~Szczyrba, \textsl{ A canonical structure for classical field
  theories},
\newblock Communications in Mathematical Physics \textbf{ 46}(2), 183--206
  (1976).

\bibitem[KS06a]{Kajiura2006b}
H.~Kajiura and J.~Stasheff, \textsl{ Homotopy Algebras Inspired by Classical
  Open-Closed String Field Theory},
\newblock Communications in Mathematical Physics \textbf{ 263}(3), 553--581
  (May 2006).

\bibitem[KS06b]{Kajiura2006}
H.~Kajiura and J.~Stasheff, \textsl{ Open-closed homotopy algebra in
  mathematical physics},
\newblock J. Math. Phys. \textbf{ 47}(2), 023506 (February 2006).

\bibitem[KSM93]{Kolar1993}
I.~Kol{\'{a}}ř, J.~Slov{\'{a}}k and P.~W. Michor,
\newblock \textsl{ Natural Operations in Differential Geometry},
\newblock Springer Berlin Heidelberg, Berlin, Heidelberg, 1993.

\bibitem[KSV95]{Kimura1995}
T.~Kimura, J.~Stasheff and A.~A. Voronov, \textsl{ On operad structures of
  moduli spaces and string theory},
\newblock Communications in Mathematical Physics \textbf{ 171}(1), 1--25 (July
  1995).

\bibitem[KV20]{krepski2020multiplicative}
D.~Krepski and J.~Vaughan, \textsl{ Multiplicative vector fields on bundle
  gerbes},
\newblock (2020), {arXiv:2003.12874}.

\bibitem[Lap72]{Laplaza1972}
M.~L. Laplaza,
\newblock Coherence for distributivity,
\newblock in \textsl{ Coherence in Categories, Lect. Notes Math. 281}, pages
  29--65, 1972.

\bibitem[LBM09]{Li-Bland2009}
D.~Li-Bland and E.~Meinrenken, \textsl{ Courant Algebroids and Poisson
  Geometry},
\newblock Int. Math. Res. Not.  (April 2009).

\bibitem[Lic97]{Lickorish1997}
W.~B.~R. Lickorish,
\newblock \textsl{ An Introduction to Knot Theory},
\newblock Springer New York, New York, NY, 1997.

\bibitem[LM95]{LadaMarkl}
T.~Lada and M.~Markl, \textsl{ Strongly homotopy {L}ie algebras},
\newblock Comm. Algebra \textbf{ 23}(6), 2147--2161 (1995), {9406095}.

\bibitem[LMS92]{Lecomte1992}
P.~Lecomte, P.~Michor and H.~Schicketanz, \textsl{ The multigraded
  Nijenhuis–Richardson algebra, its universal property and applications},
\newblock J. Pure Appl. Algebr. \textbf{ 77}(1), 87--102 (February 1992).

\bibitem[LR12]{Liu-Ricca12}
X.~{Liu} and R.~L. {Ricca}, \textsl{ The Jones polynomial for fluid knots from
  helicity},
\newblock J. Phys. A, Math. Theor. \textbf{ 45}(20), 14 (2012),
\newblock Id/No 205501.

\bibitem[LR15]{Liu-Ricca15}
X.~{Liu} and R.~L. {Ricca}, \textsl{ On the derivation of the HOMFLYPT
  polynomial invariant for fluid knots},
\newblock J. Fluid Mech. \textbf{ 773}, 34--48 (2015).

\bibitem[LS93]{LadaStasheff}
T.~Lada and J.~Stasheff, \textsl{ Introduction to {SH} {L}ie algebras for
  physicists},
\newblock Internat. J. Theoret. Phys. \textbf{ 32}(7), 1087--1103 (July 1993),
  {9209099}.

\bibitem[LV12]{Loday}
J.-L. Loday and B.~Vallette,
\newblock \textsl{ Algebraic Operads}, volume 346 of \textsl{ Grundlehren der
  mathematischen Wissenschaften},
\newblock Springer Berlin Heidelberg, Berlin, Heidelberg, 2012.

\bibitem[{Mac}78]{MacLane1978}
S.~{Mac Lane},
\newblock \textsl{ Categories for the Working Mathematician}, volume~5 of
  \textsl{ Graduate Texts in Mathematics},
\newblock Springer New York, New York, NY, 1978.

\bibitem[Mam20]{Mammadova2020a}
L.~Mammadova,
\newblock \textsl{ Moment maps in multisymplectic geometry},
\newblock Doctoral thesis, KU Leuven, 2020.

\bibitem[Man11a]{Manchon2011a}
D.~Manchon,
\newblock A short survey on pre-Lie algebras,
\newblock in \textsl{ Noncommutative Geom. Phys. Renorm. Motiv. Index Theory},
  pages 89--102, European Mathematical Society Publishing House, Zuerich,
  Switzerland, 2011.

\bibitem[Man11b]{Manetti-website}
M.~Manetti,
\newblock Deformation Theory (lecture notes),
\newblock \url{https://www1.mat.uniroma1.it/people/manetti/DT2011/DT.html},
  2011.

\bibitem[Man11c]{Manetti-website-coalgebras}
M.~Manetti,
\newblock Graded coalgebras (lecture notes),
\newblock
  \url{https://www1.mat.uniroma1.it/people/manetti/DT2011/gradedcoalgebras.pdf},
  2011.

\bibitem[Mar98]{Markl1998}
M.~Markl, \textsl{ Homotopy Algebras via Resolutions of Operads},
\newblock Circ. Mat. di Palermo (63), 157--164 (1998).

\bibitem[{Mas}72]{Mas}
V.~P. {Maslov},
\newblock Th\'eorie des perturbations et m\'ethodes asymptotiques. Suivi de
  deux notes complementaires de V. I. Arnol'd et V. C. Bouslaev. Traduit par J.
  Lascoux et R. Seneor,
\newblock Etudes mathematiques. Paris: Dunod; Paris: Gauthier-Villars. XVI,384
  p. F 170.00 (1972)., 1972.

\bibitem[Mei04]{Meinrenken2004}
E.~Meinrenken,
\newblock Equivariant Cohomology and the Cartan Model,
\newblock in \textsl{ Encycl. Math. Phys. Five-Volume Set}, pages 242--250,
  2004.

\bibitem[Mit15]{Miti2015}
A.~M. Miti,
\newblock \textsl{ Algebraic quantization of Jacobi fields and geometric
  approach to Peierls brackets},
\newblock Master's thesis, Universit{\`{a}} Degli Studi di Milano, 2015.

\bibitem[Mol07]{Molitor2007}
M.~Molitor,
\newblock \textsl{ Grassmanniennes non-lin$\backslash$'eaires, groupes de
  diff$\backslash$'eomorphismes unimodulaires et quelques
  $\backslash$'equations hamiltoniennes en dimension infinie},
\newblock PhD thesis, Universit$\backslash$'e Paul Verlaine-Metz, 2007.

\bibitem[MPSW01]{MPSW}
J.~E. Marsden, S.~Pekarsky, S.~Shkoller and M.~West, \textsl{ Variational
  methods, multisymplectic geometry and continuum mechanics},
\newblock Journal of Geometry and Physics \textbf{ 38}(3-4), 253--284 (2001).

\bibitem[MR92]{Moffatt-Ricca92}
H.~K. {Moffatt} and R.~L. {Ricca}, \textsl{ {Helicity and the
  C\u{a}lug\u{a}reanu invariant}},
\newblock Proc. R. Soc. Lond., Ser. A \textbf{ 439}(1906), 411--429 (1992).

\bibitem[MR20a]{Mammadova2020}
L.~Mammadova and L.~Ryvkin, \textsl{ On the extension problem for weak moment
  maps},
\newblock (January 2020), {arXiv:2001.00264}.

\bibitem[MR20b]{Miti2019}
A.~M. Miti and L.~Ryvkin, \textsl{ Multisymplectic actions of compact Lie
  groups on spheres},
\newblock Journal of Symplectic Geometry \textbf{ 18}(6), 1751--1785 (2020),
  {1906.08790}.

\bibitem[MS43]{MR0008817}
D.~Montgomery and H.~Samelson, \textsl{ Transformation groups of spheres},
\newblock Ann. of Math. (2) \textbf{ 44}, 454--470 (1943).

\bibitem[MS98]{McD-Sal}
D.~{McDuff} and D.~{Salamon},
\newblock \textsl{ Introduction to symplectic topology. 2nd ed},
\newblock New York, NY: Oxford University Press, 2nd ed. edition, 1998.

\bibitem[MS12a]{Madsen2013}
T.~B. Madsen and A.~Swann, \textsl{ Closed forms and multi-moment maps},
\newblock Geometriae Dedicata \textbf{ 165}(1), 25--52 (September 2012).

\bibitem[MS12b]{Madsen2012}
T.~B. Madsen and A.~Swann, \textsl{ Multi-moment maps},
\newblock Adv. Math. (N. Y). \textbf{ 229}(4), 2287--2309 (2012).

\bibitem[MS20]{Miti2019a}
A.~M. Miti and M.~Spera, \textsl{ Derivation of the HOMFLYPT knot polynomial
  via helicity and geometric quantization},
\newblock Bollettino dell'Unione Matematica Italiana  (August 2020),
  {1910.13400}.

\bibitem[MS21]{Miti2018}
A.~M. Miti and M.~Spera, \textsl{ A Hydrodynamical Homotopy Co-momentum Map And
  A Multisymplectic Interpretation Of Higher-order Linking Numbers},
\newblock Journal of the Australian Mathematical Society , 1--20 (February
  2021), {1805.01696}.

\bibitem[{Mug}10]{Muger2010}
M.~{Muger}, \textsl{ Tensor categories: a selective guided tour},
\newblock Rev. Uni\'on Mat. Argent. \textbf{ 51}(1), 95--163 (2010).

\bibitem[MW83]{MW83}
J.~Marsden and A.~Weinstein, \textsl{ Coadjoint orbits, vortices, and Clebsch
  variables for incompressible fluids},
\newblock Physica D: Nonlinear Phenomena \textbf{ 7}(1-3), 305--323 (1983).

\bibitem[MZ20]{Mammadova2019}
L.~Mammadova and M.~Zambon, \textsl{ Lie 2-algebra moment maps in
  multisymplectic geometry},
\newblock Differential Geometry and its Applications \textbf{ 70}, 101631 (June
  2020), {arXiv:1901.10842}.

\bibitem[MZ21]{Miti2020}
A.~M. Miti and M.~Zambon,
\newblock Observables on multisymplectic manifolds and higher Courant
  algebroids,
\newblock (expected) Spring 2021.

\bibitem[na20a]{nlab:algebroid}
nLab authors,
\newblock algebroid,
\newblock \url{http://ncatlab.org/nlab/show/algebroid}, March 2020,
\newblock \href{http://ncatlab.org/nlab/revision/algebroid/14}{Revision 14}.

\bibitem[na20b]{nlab:cartesian_closed_category}
nLab authors,
\newblock cartesian closed category,
\newblock \url{http://ncatlab.org/nlab/show/cartesian%20closed%20category},
  February 2020,
\newblock
  \href{http://ncatlab.org/nlab/revision/cartesian%20closed%20category/32}{Revision
  32}.

\bibitem[na20c]{nlab:closed_monoidal_category}
nLab authors,
\newblock closed monoidal category,
\newblock \url{http://ncatlab.org/nlab/show/closed%20monoidal%20category},
  March 2020,
\newblock
  \href{http://ncatlab.org/nlab/revision/closed%20monoidal%20category/41}{Revision
  41}.

\bibitem[na20d]{nlab:currying}
nLab authors,
\newblock currying,
\newblock \url{http://ncatlab.org/nlab/show/currying}, October 2020,
\newblock \href{http://ncatlab.org/nlab/revision/currying/4}{Revision 4}.

\bibitem[na20e]{nlab:distributive_monoidal_category}
nLab authors,
\newblock distributive monoidal category,
\newblock
  \url{http://ncatlab.org/nlab/show/distributive%20monoidal%20category},
  February 2020,
\newblock
  \href{http://ncatlab.org/nlab/revision/distributive%20monoidal%20category/7}{Revision
  7}.

\bibitem[na20f]{nlab:double_complex}
nLab authors,
\newblock double complex,
\newblock \url{http://ncatlab.org/nlab/show/double%20complex}, October 2020,
\newblock \href{http://ncatlab.org/nlab/revision/double%20complex/20}{Revision
  20}.

\bibitem[na20g]{nlab:framed_link}
nLab authors,
\newblock framed link,
\newblock \url{http://ncatlab.org/nlab/show/framed%20link}, November 2020,
\newblock \href{http://ncatlab.org/nlab/revision/framed%20link/3}{Revision 3}.

\bibitem[na20h]{nlab:functor_category}
nLab authors,
\newblock functor category,
\newblock \url{http://ncatlab.org/nlab/show/functor%20category}, February 2020,
\newblock
  \href{http://ncatlab.org/nlab/revision/functor%20category/20}{Revision 20}.

\bibitem[na20i]{nlab:gelfand_duality}
nLab authors,
\newblock {{G}}elfand duality,
\newblock \url{http://ncatlab.org/nlab/show/Gelfand%20duality}, October 2020,
\newblock \href{http://ncatlab.org/nlab/revision/Gelfand%20duality/29}{Revision
  29}.

\bibitem[na20j]{nlab:higher_symplectic_geometry}
nLab authors,
\newblock higher symplectic geometry,
\newblock \url{http://ncatlab.org/nlab/show/higher%20symplectic%20geometry},
  October 2020,
\newblock
  \href{http://ncatlab.org/nlab/revision/higher%20symplectic%20geometry/7}{Revision
  7}.

\bibitem[na20k]{nlab:horizontal_categorification}
nLab authors,
\newblock horizontal categorification,
\newblock \url{http://ncatlab.org/nlab/show/horizontal%20categorification},
  November 2020,
\newblock
  \href{http://ncatlab.org/nlab/revision/horizontal%20categorification/20}{Revision
  20}.

\bibitem[na20l]{nlab:l-infinity-algebra}
nLab authors,
\newblock {{L}}-infinity-algebra,
\newblock \url{http://ncatlab.org/nlab/show/L-infinity-algebra}, November 2020,
\newblock
  \href{http://ncatlab.org/nlab/revision/L-infinity-algebra/97}{Revision 97}.

\bibitem[na20m]{nlab:lie_algebra}
nLab authors,
\newblock {{L}}ie algebra,
\newblock \url{http://ncatlab.org/nlab/show/Lie%20algebra}, November 2020,
\newblock \href{http://ncatlab.org/nlab/revision/Lie%20algebra/49}{Revision
  49}.

\bibitem[na20n]{nlab:monoidal_category}
nLab authors,
\newblock monoidal category,
\newblock \url{http://ncatlab.org/nlab/show/monoidal%20category}, November
  2020,
\newblock
  \href{http://ncatlab.org/nlab/revision/monoidal%20category/126}{Revision
  126}.

\bibitem[na20o]{nlab:rig_category}
nLab authors,
\newblock rig category,
\newblock \url{http://ncatlab.org/nlab/show/rig%20category}, February 2020,
\newblock \href{http://ncatlab.org/nlab/revision/rig%20category/10}{Revision
  10}.

\bibitem[na20p]{nlab:split_exact_sequence}
nLab authors,
\newblock split exact sequence,
\newblock \url{http://ncatlab.org/nlab/show/split%20exact%20sequence}, November
  2020,
\newblock
  \href{http://ncatlab.org/nlab/revision/split%20exact%20sequence/17}{Revision
  17}.

\bibitem[na21]{nlab:internal_hom}
nLab authors,
\newblock internal hom,
\newblock \url{http://ncatlab.org/nlab/show/internal%20hom}, February 2021,
\newblock \href{http://ncatlab.org/nlab/revision/internal%20hom/67}{Revision
  67}.

\bibitem[Nes03]{Nestruev2010}
J.~Nestruev,
\newblock \textsl{ Smooth Manifolds and Observables}, volume 220 of \textsl{
  Graduate Texts in Mathematics},
\newblock Springer-Verlag, New York, 2003.

\bibitem[NR67]{Nijenhuis1967}
A.~Nijenhuis and R.~{Richardson, Jr.}, \textsl{ Deformations of Lie Algebra
  Structures},
\newblock Indiana University Mathematics Journal \textbf{ 17}(1), 89--105
  (1967).

\bibitem[Pri10]{Pridham2010}
J.~Pridham, \textsl{ Unifying derived deformation theories},
\newblock Adv. Math. (N. Y). \textbf{ 224}(3), 772--826 (June 2010),
  {arXiv:0705.0344v7}.

\bibitem[PS89]{Pe-Spe89}
V.~Penna and M.~Spera, \textsl{ A geometric approach to quantum vortices},
\newblock Journal of mathematical physics \textbf{ 30}(12), 2778--2784 (1989).

\bibitem[PS92]{Pe-Spe92}
V.~Penna and M.~Spera, \textsl{ On coadjoint orbits of rotational perfect
  fluids},
\newblock Journal of mathematical physics \textbf{ 33}(3), 901--909 (1992).

\bibitem[PS00]{Pe-Spe00}
V.~Penna and M.~Spera, \textsl{ String limit of vortex current algebra},
\newblock Physical Review B \textbf{ 62}(21), 14547 (2000).

\bibitem[PS02]{Pe-Spe02}
V.~Penna and M.~Spera, \textsl{ Higher order linking numbers, curvature and
  holonomy},
\newblock Journal of Knot Theory and Its Ramifications \textbf{ 11}(05),
  701--723 (2002).

\bibitem[PT87]{PT}
J.~H. {Przytycki} and P.~{Traczyk}, \textsl{ Conway algebras and skein
  equivalence of links},
\newblock Proc. Am. Math. Soc. \textbf{ 100}, 744--748 (1987).

\bibitem[Rei19]{Reinhold2019}
B.~Reinhold, \textsl{ $L_\infty$-algebras and their cohomology},
\newblock Emergent Sci. \textbf{ 3}, 4 (June 2019).

\bibitem[RGK94]{RAMADEVI1994}
P.~Ramadevi, T.~Govindarajan and R.~Kaul, \textsl{ Chirality Of Knots $9_{42}$
  And $10_{71}$ And Chern-simons Theory},
\newblock Mod. Phys. Lett. A \textbf{ 09}(34), 3205--3217 (November 1994).

\bibitem[{Rie}16]{Riehl2016}
E.~{Riehl},
\newblock \textsl{ Category theory in context},
\newblock Mineola, NY: Dover Publications, 2016.

\bibitem[RN11]{Ricca2011}
R.~L. Ricca and B.~Nipoti, \textsl{ Gauss' linking Number Revisited},
\newblock Journal of Knot Theory and Its Ramifications \textbf{ 20}(10),
  1325--1343 (October 2011).

\bibitem[Rog01]{Roger2012}
C.~Roger,
\newblock Unimodular vector fields and deformation quantization,
\newblock in \textsl{ Deform. Quantization}, pages 135--148, De Gruyter,
  Berlin, Boston, 2001.

\bibitem[Rog05]{Roger2005}
C.~Roger,
\newblock The group of volume preserving diffeomorphisms and the Lie algebra of
  unimodular vector fields: survey of some classical and not-so-classical
  results,
\newblock in \textsl{ Twenty years Bialowieza A Math. Anthol. Asp. Differ.
  Geom. methods Phys.}, pages 79--98, Hackensack, NJ: World Scientific, 2005.

\bibitem[Rog11]{Rogers2011}
C.~L. Rogers,
\newblock \textsl{ Higher Symplectic Geometry},
\newblock Phd thesis, University of California Riverside, June 2011.

\bibitem[Rog12]{Rogers2010}
C.~L. Rogers, \textsl{ $L_\infty$-Algebras from Multisymplectic Geometry},
\newblock Letters in Mathematical Physics \textbf{ 100}(1), 29--50 (April
  2012), {arXiv:1005.2230}.

\bibitem[Rog13]{Rogers2013}
C.~L. Rogers, \textsl{ 2-Plectic geometry, Courant algebroids, and categorified
  prequantization},
\newblock J. Symplectic Geom. \textbf{ 11}(1), 53--91 (September 2013),
  {arXiv:1009.2975}.

\bibitem[Rol03]{Rolfsen}
D.~Rolfsen,
\newblock \textsl{ Knots and links}, volume 346,
\newblock American Mathematical Soc., 2003.

\bibitem[RR75]{Rasetti-Regge75}
M.~Rasetti and T.~Regge, \textsl{ Vortices in He II, current algebras and
  quantum knots},
\newblock Physica A: Statistical Mechanics and its Applications \textbf{
  80}(3), 217--233 (1975).

\bibitem[RS15]{Ritter2015a}
P.~Ritter and C.~Saemann, \textsl{ Automorphisms of Strong Homotopy Lie
  Algebras of Local Observables},
\newblock arXiv:1507.00972  (07 2015), {arXiv:1507.00972}.

\bibitem[RW98]{Roytenberg1998}
D.~Roytenberg and A.~Weinstein, \textsl{ Courant Algebroids and Strongly
  Homotopy Lie Algebras},
\newblock Lett. Math. Phys. \textbf{ 46}(1), 81--93 (February 1998), {9802118}.

\bibitem[RW15]{zbMATH06448534}
L.~Ryvkin and T.~Wurzbacher, \textsl{ Existence and unicity of co-moments in
  multisymplectic geometry.},
\newblock Differ. Geom. Appl. \textbf{ 41}, 1--11 (2015).

\bibitem[RW19]{Ryvkin2018}
L.~Ryvkin and T.~Wurzbacher, \textsl{ An invitation to multisymplectic
  geometry},
\newblock Journal of Geometry and Physics \textbf{ 142}, 9--36 (August 2019).

\bibitem[RWZ20]{Ryvkin2016}
L.~Ryvkin, T.~Wurzbacher and M.~Zambon, \textsl{ Conserved quantities on
  multisymplectic manifolds},
\newblock Journal of the Australian Mathematical Society \textbf{ 108}(1),
  120--144 (February 2020), {arXiv:1610.05592}.

\bibitem[Ryv16]{Ryvkin2016a}
L.~Ryvkin,
\newblock \textsl{ Observables and Symmetries of $n$-Plectic Manifolds},
\newblock Springer Fachmedien Wiesbaden, Wiesbaden, 2016.

\bibitem[Sau89]{Saunders1989}
D.~J. Saunders,
\newblock \textsl{ The Geometry of Jet Bundles},
\newblock Cambridge University Press, March 1989.

\bibitem[Sch09]{Schatz2009}
F.~Sch\"atz,
\newblock \textsl{ Coisotropic submanifolds and the BFV-complex},
\newblock PhD thesis, Zurich, 2009.

\bibitem[Soc20]{encyclopedia:Isotropy}
E.~M. Society,
\newblock Isotropy group,
\newblock \url{https://www.encyclopediaofmath.org/}, February 2020.

\bibitem[{Sou}66]{Souriau66}
J.-M. {Souriau}, \textsl{ Quantification g\'eom\'etrique},
\newblock Communications in Mathematical Physics \textbf{ 1}, 374--398 (1966).

\bibitem[Sou70]{Souriau70}
J.-M. Souriau,
\newblock \textsl{ Structure des syst\'emes dynamiques: ma\^itrises de
  math\`ematiques},
\newblock 1970.

\bibitem[Spe06]{Spe06}
M.~Spera, \textsl{ A survey on the differential and symplectic geometry of
  linking numbers},
\newblock Milan J. Math. \textbf{ 74}(1), 139--197 (2006).

\bibitem[{Spe}11]{Spe11}
M.~{Spera}, \textsl{ A note on \(n\)-gerbes and transgressions},
\newblock Port. Math. (N.S.) \textbf{ 68}(4), 381--387 (2011).

\bibitem[Spe16]{Spera16}
M.~Spera, \textsl{ Moment map and gauge geometric aspects of the
  Schr{\"{o}}dinger and Pauli equations},
\newblock Int. J. Geom. Methods Mod. Phys. \textbf{ 13}(04), 1630004 (2016).

\bibitem[SSS12]{Sati2012a}
H.~Sati, U.~Schreiber and J.~Stasheff, \textsl{ Twisted Differential String and
  Fivebrane Structures},
\newblock Communications in Mathematical Physics \textbf{ 315}(1), 169--213
  (October 2012), {arXiv:0910.4001v2}.

\bibitem[Sta08]{Stacey2008}
A.~Stacey, \textsl{ Comparative Smootheology},
\newblock Theory Appl. Categ. \textbf{ 25}, 64--117 (February 2008),
  {arXiv:0802.2225}.

\bibitem[{Sta}18]{stacks-project}
T.~{Stacks Project Authors},
\newblock \textit{Stacks Project},
\newblock \url{https://stacks.math.columbia.edu}, 2018.

\bibitem[Sta19]{Stasheff2019}
J.~Stasheff, \textsl{ L-infinity and A-infinity structures},
\newblock High. Struct. \textbf{ 3}(1), 292--326 (2019).

\bibitem[SW20]{SevestreWurzbacherPreq}
G.~Sevestre and T.~Wurzbacher, \textsl{ On the prequantisation map for
  2-plectic manifolds},
\newblock (08 2020), {arXiv:2008.05184}.

\bibitem[Swa91]{Swann1991}
A.~Swann, \textsl{ HyperK$\backslash$:aehler and quaternionic
  K$\backslash$:ahler geometry},
\newblock Math. Ann. \textbf{ 289}(1), 421--450 (March 1991).

\bibitem[SZ15]{Shahbazi2016}
C.~S. Shahbazi and M.~Zambon, \textsl{ Products of multisymplectic manifolds
  and homotopy moment maps},
\newblock J. Lie Theory \textbf{ 26}(4), 1037--1067 (April 2015),
  {arXiv:1504.08194}.

\bibitem[{Tak}09]{Agoh}
{Takashi {Agoh} and Karl {Dilcher}}, \textsl{ Shortened recurrence relations
  for Bernoulli numbers},
\newblock Discrete Math. \textbf{ 309}(4), 887--898 (2009).

\bibitem[Tav94]{Tavares}
J.~Tavares, \textsl{ Chen integrals, generalized loops and loop calculus},
\newblock International Journal of Modern Physics A \textbf{ 9}(26), 4511--4548
  (1994).

\bibitem[Tu11]{Tu2011}
L.~W. Tu, \textsl{ What is ... equivariant cohomology?},
\newblock Notices Am. Math. Soc. \textbf{ 58}(3), 423--426 (2011),
  {arXiv:1305.4293}.

\bibitem[Val14]{Vallette2014}
B.~Vallette, \textsl{ Algebra + Homotopy = Operad},
\newblock Symplectic, Poisson Noncommutative Geom. \textbf{ 62}, 229--290
  (February 2014), {arXiv:1202.3245v1}.

\bibitem[Vau14]{VaughanJ}
J.~Vaughan, \textsl{ Metaplectic-c Quantomorphisms},
\newblock Symmetry, Integrability and Geometry: Methods and Applications
  \textbf{ 11} (10 2014).

\bibitem[Vin90]{zbMATH04172022}
A.~M. Vinogradov, \textsl{ Combination of Schouten and Nijenhuis brackets,
  cohomologies and superdifferential operators},
\newblock Mat. Zametki \textbf{ 47}(6), 138--140 (1990).

\bibitem[Vit14]{Vitagliano2013}
L.~Vitagliano, \textsl{ On the strong homotopy Lie–Rinehart algebra of a
  foliation},
\newblock Commun. Contemp. Math. \textbf{ 16}(06), 1450007 (December 2014),
  {arXiv:1204.2467}.

\bibitem[Vor04]{Voronov2004}
T.~T. Voronov, \textsl{ Higher Derived Brackets for Arbitrary Derivations},
\newblock Trav. Math{\'{e}}matiques, XVI (2005), 163–186, , 24 (2004),
  {arXiv:0412202}.

\bibitem[Vor05]{Voronov2005}
T.~Voronov, \textsl{ Higher derived brackets and homotopy algebras},
\newblock J. Pure Appl. Algebr. \textbf{ 202}(1-3), 133--153 (November 2005).

\bibitem[War83]{Warner}
F.~W. Warner,
\newblock \textsl{ Foundations of Differentiable Manifolds and Lie Groups},
  volume~94 of \textsl{ Graduate Texts in Mathematics},
\newblock Springer New York, New York, NY, 1983.

\bibitem[Wei]{Weisstein}
E.~W. Weisstein,
\newblock Bernoulli Number. {From MathWorld---A Wolfram Web Resource},
\newblock
  \href{https://mathworld.wolfram.com/BernoulliNumber.html}{https://mathworld.wolfram.com/BernoulliNumber.html},
\newblock Accessed: 2020-09-23.

\bibitem[Wei77]{weinstein1977lectures}
A.~Weinstein,
\newblock \textsl{ Lectures on symplectic manifolds},
\newblock Number~29, American Mathematical Soc., 1977.

\bibitem[Wei94]{Weibel}
C.~A. Weibel,
\newblock \textsl{ An Introduction to Homological Algebra},
\newblock Cambridge University Press, Cambridge, April 1994.

\bibitem[{Wey}35]{Weyl35}
H.~{Weyl}, \textsl{ Geodesic fields in the calculus of variation for multiple
  integrals},
\newblock Ann. Math. (2) \textbf{ 36}, 607--629 (1935).

\bibitem[{Wit}89]{Witten89}
E.~{Witten}, \textsl{ Quantum field theory and the Jones polynomial},
\newblock Communications in Mathematical Physics \textbf{ 121}(3), 351--399
  (1989).

\bibitem[Woo97]{Woodhouse97}
N.~Woodhouse,
\newblock \textsl{ Geometric quantization},
\newblock Oxford: Clarendon Press, 2nd ed. edition, 1997.

\bibitem[WZ05]{Weinstein2005a}
A.~Weinstein and M.~Zambon,
\newblock Variations on prequantization,
\newblock in \textsl{ Proc. 4th Conf. Poisson Geom. Luxemb. June 7--11, 2004},
  volume~16, pages 187--219, Luxembourg: Universit{\'{e}} du Luxembourg, 2005.

\bibitem[Zam12]{Zambon2012}
M.~Zambon, \textsl{ $L_\infty$-algebras and higher analogues of Dirac
  structures and Courant algebroids},
\newblock Journal of Symplectic Geometry \textbf{ 10}(4), 563--599 (2012),
  {arXiv:1003.1004}.

\bibitem[Zuc87]{Zuckerman87}
G.~J. Zuckerman,
\newblock Action Principles and Global Geometry,
\newblock in \textsl{ Math. Asp. String Theory}, pages 259--284, World
  Scientific, September 1987.

\end{thebibliography}
